\documentclass{amsbook}

\usepackage{graphicx}

\newtheorem{theorem}{Theorem}[chapter]
\newtheorem{lemma}[theorem]{Lemma}
\newtheorem{corollary}[theorem]{Corollary}

\theoremstyle{definition}
\newtheorem{definition}[theorem]{Definition}
\newtheorem{example}[theorem]{Example}

\theoremstyle{remark}
\newtheorem{remark}[theorem]{Remark}

\numberwithin{section}{chapter}
\numberwithin{theorem}{section}
\numberwithin{equation}{section}

\begin{document}
\frontmatter
\title[Spatial Complex Variables and Applications]{Spatial Complex Variables and Applications}
%\title{Space Complex Variables \\
%Theory and Applications}

\author{Shan-guang Tan}
%\address{Zhejiang University, Hangzhou, 310027, CHINA}
\curraddr{Zhejiang University, Hangzhou, 310027, CHINA}
\email{tanshanguang@163.com}

%\date{June 6, 2026}

\begin{abstract}
We introduce a new kind of numbers $s=e_{xy}(x+iy)+z$ called Spatial Complex Numbers and establish a new theory of Spatial Complex Variables here where $e_{xy}$ denotes the $xy$ coordinate plane.

Studies of the theory of spatial complex variables include: (1) Definition and Algebras; (2) Analytic functions, Cauchy-Riemann equations, and spatial harmonic functions; (3) Elementary functions; (4) Spatial contour integrals; (5) Surface integrals; (6) Taylor and Laurent series; (7) Residue theorem; (8) Conformal mapping; (9) Applications of three-dimensional fluid flow and Navier-Stokes equations; (10) The Schwarz-Christoffel Transformation; (11) Spatial circle and spherical surface integral formulas of the Poisson type; (12) Applications of Residues in two-dimensional improper integrals; (13) Joukowski Transformation and Applications for airfoils.

Results show that spatial complex numbers are with good properties of fields just as the two-dimensional complex numbers. Spatial complex numbers follow the commutative and associative laws of multiplication of numbers, the distributive law of addition and multiplication of numbers, and other laws of numbers. The theory of spatial complex variables can be used to solve many three-dimensional problems in sciences and technologies, such as expressed by the three-dimensional Laplace equation. So numbers are developed from the real and two-dimensional complex numbers to three-dimensional complex numbers, and the number system is enlarged.
\end{abstract}

\maketitle

\setcounter{page}{4}
\tableofcontents

%-----------------------------------------------------------------------------
% Beginning of preface.tex
%-----------------------------------------------------------------------------
%
% AMS-LaTeX 1.2 sample file for a monograph, based on amsbook.cls.
% This is a data file input by chapter.tex.
%%%%%%%%%%%%%%%%%%%%%%%%%%%%%%%%%%%%%%%%%%%%%%%%%%%%%%%%%%%%%%%%%%%%%%%%

\chapter*{Preface}

It's well known that the numbers are composed of the real and two-dimensional complex numbers. A two-dimensional complex number is with a pair of reals and forms a complex plane.

Researches on developing the three-dimensional complex numbers have been made since years ago such as the Hyper-complex Numbers given by Hamilton, Gelasiman or others~\cite{F,G}, but few important results were derived.
Normally, it's well known that there are no commutative fields for multiplications of numbers, which are three dimensional vector spaces over the reals in the rectangle coordinate systems with the form of $x+iy+jz$ where $i=\sqrt{-1}$ and $j=\sqrt{-1}$, other than the complex numbers and if one allows non-commutative fields, the only extra possibility is the 4-dimensional algebra of quaternion numbers in the rectangle coordinate systems.

Our purpose here is to introduce a new kind of numbers called Spatial Complex Numbers and to establish a new theory of spatial complex variables.

In the first chapter, the spatial complex numbers and their operations of addition and multiplication are defined.

First, the basic definition of spatial complex numbers is given out as $s=e_{xy}(x+iy)+z$ where $e_{xy}$ denotes the $xy$ coordinate plane and $i=\sqrt{-1}$.

Next, the algebras of spatial complex numbers are defined, particularly, a key definition that makes spatial complex numbers be a number field is the product of two spatial complex numbers $s_{1}=e_{xy}(x_{1}+iy_{1})+z_{1}$ and $s_{2}=e_{xy}(x_{2}+iy_{2})+z_{2}$, which is defined as follows:
\[
s_{1}s_{2}=e_{xy}[(x_{1}+iy_{1})(x_{2}+iy_{2})+(x_{1}+iy_{1})z_{2}+(x_{2}+iy_{2})z_{1}]+z_{1}z_{2}.
\]
Note that $(x_{1}+iy_{1})(x_{2}+iy_{2})$ is the product of two plane complex numbers $(x_{1}+iy_{1})$ and $(x_{2}+iy_{2})$ in the $xy$ coordinate plane, which is also a plane complex number in the $xy$ coordinate plane; $(x_{1}+iy_{1})z_{2}$ and $(x_{2}+iy_{2})z_{1}$ are two products of a real and a plane complex number in the $xy$ coordinate plane, which are also two plane complex numbers in the $xy$ coordinate plane; so that the expression
\[
(x_{1}+iy_{1})(x_{2}+iy_{2})+(x_{1}+iy_{1})z_{2}+(x_{2}+iy_{2})z_{1}
\]
is the sum of three plane complex numbers in the $xy$ coordinate plane, which is also a plane complex number in the $xy$ coordinate plane. For the sake of simplicity, we can take the symbol $e_{xy}$ as an operator in the $xy$ coordinate plane, then there can be $e_{xy}^{n}=e_{xy}$ for any integer $n\neq0$, and for a term $e_{xy}x^{j}y^{k}$, the operator $e_{xy}$ becomes $e_{xy}^{m}$ where $m=|j|+|k|$.

Then, it is proved that the operations of addition and multiplication of spatial complex numbers follow the commutative law and the associative law of numbers, and follow the distributive law of addition and multiplication of numbers. Other laws of numbers are also followed by spatial complex numbers.

Third, the geometry of spatial complex numbers is investigated, such as the moduli, complex conjugates, and exponential form of spatial complex numbers.

Finally, some special properties, such as products and quotients in exponential form, the roots of spatial complex numbers, and some spatial planes and regions in the complex space are discussed.

Hence the spatial complex numbers form a number field. The spatial complex numbers are with good properties just as the real numbers and the classic plane complex numbers. So the normal complex numbers are developed from one in complex plane with two dimensions to another in complex space with three dimensions, and the number system is enlarged also.

In the second chapter, analysis of spatial complex functions of a spatial complex variable is made. First, the basic format of spatial complex functions of a spatial complex variable is defined as
\[
f(s)=e_{xy}[u(x,y,z)+iv(x,y,z)]+w(z)
\]
where $u$ and $v$ are real functions of three real variables, and $w$ is a real function of one real variable, then some properties of the spatial complex functions are analyzed and discussed, such as the limit, continuity, derivative, analytic functions, Cauchy-Riemann equations in complex space, three-dimensional complex harmonic functions of three real variables $x$, $y$, and $z$, and the three-dimensional orthogonal property of the first-order partial derivatives of a pair of complex harmonic conjugate functions.

In the third chapter, some elementary functions of a spatial complex variable, such as the polynomial, exponential, logarithmic, trigonometric, hyperbolic, inverse trigonometric, and inverse hyperbolic functions, are defined. Their power series are also given out. Most of elementary functions $f(s)$ are analytic in a domain $D$ and can be expressed as
\[
f(s)=e_{xy}[f(x+iy+z)-f(z)]+f(z)
\]
and $w(z)=f(z)$,
\[
u(x,y,z)+w(z)=Re\texttt{ }f(x+iy+z),\texttt{ }
v(x,y,z)=Im\texttt{ }f(x+iy+z).
\]

In the fourth chapter, the spatial contour integral of spatial complex functions is analyzed. First, the spatial
contour integral of spatial complex functions are defined and Cauchy spatial contour integral theorem of the spatial complex functions is put forward and proved. Next, extensions of Cauchy spatial contour integral theorem are made and proved. Then, the Cauchy-formula of spatial contour integral of spatial complex functions is put forward and proved, and the infinite differential properties of the analytic function are discussed. Finally, some theorems of the spatial contour integral of spatial complex functions are put forward and proved, which are Morera theorem, Liouville's Theorem, and the fundamental theorem of algebra, and the maximum modulus principle involving maximum values of the moduli of analytic functions is analyzed and proved.

In the fifth chapter, the surface integral of spatial complex functions is analyzed.

First, the surface integrals and Cauchy-theorem of surface integrals of the second kind and first kind are put forward and proved. When a spatial complex function $f(s)=\varpi=(U(x,y,z),V(x,y,z),W(z))$ is analytic, its surface integral of the second kind satisfies the Gauss-Ostrogradskii Formula, and for a closed surface $S$ there is
\[
\iint_{S}X(x,y,z)dydz+Y(x,y,z)dzdx+Z(x,y,z)dxdy=0
\]
where
\[
X=\frac{\partial W}{\partial y}-\frac{\partial V}{\partial z},\texttt{ }
Y=\frac{\partial U}{\partial z}-\frac{\partial W}{\partial x},\texttt{ }
Z=\frac{\partial V}{\partial x}-\frac{\partial U}{\partial y};
\]
and its boundary contour integrals of the second kind satisfies the Stokes' Formula, and for a closed surface $S$ or a closed boundary contour $C$ there is
\[
\oint_{C}U(x,y,z)dx+V(x,y,z)dy+W(z)dz=0.
\]
The Cauchy theorem of surface integrals of the first kind is also proved: Suppose that a spatial complex function $f(s)$ is analytic in a simply connected complex space $D$, and $S$ is any simple closed surface, then, there is
\[
\iint_{S}f(s)d\sigma=0.
\]

Next, a theorem called Boundary-Integration of the spatial complex functions is put forward and proved, which means that the surface integral of the spatial complex functions could be expressed by the integral of spatial complex functions on the boundary curve of the surface. Then, the Cauchy-formula of surface integral of spatial complex functions is put forward and proved, and the infinite differential properties of the analytic function are discussed.

Finally, the Morera theorem of a spatial complex variable is put forward and proved, and the maximum modulus principle and Schwarz lemma of a spatial complex variable are analyzed and proved.

In the sixth chapter, the series of spatial complex functions of a spatial complex variable is analyzed. First, Taylor series and its properties of spatial complex functions are defined and proved. The coefficients of Taylor series can be expressed by spatial circle integrals in spatial planes and spherical surface integrals in space, respectively. Then the power series and its convergence sphere is discussed. Next, Laurent series in a spherical shell are defined, and its properties are discussed and proved. The coefficients of Laurent series can also be expressed by spatial circle integrals in spatial planes and spherical surface integrals in space, respectively. Finally the properties of Taylor series and Laurent series of spatial complex functions of a spatial complex variable in spherical domain is put forward and proved.

In the seventh chapter, the residue theorem of spatial complex functions is analyzed with contour integrals in spatial planes and surface integrals in space, respectively. First, the residue of spatial complex functions is defined. Next, the residue theorem of spatial complex functions is proved. Finally, the types and properties of residues are discussed, including isolated singular points, poles and zeros.

In the eighth chapter, mapping by elementary functions of a spatial complex variable is made. First, fractional linear functions are analyzed and discussed. It is proved that the spheres in the extended complex space are conformally mapped onto spheres by fractional linear functions. Next, some special functions are discussed also here, including the transformation $\varpi$=l/s, linear fractional transformations, mappings of the upper half plane, the transformation $\varpi=\sin s$, mappings by $s^{2}$ and branches of $s^{1/2}$, and square roots of polynomials.

In the ninth chapter, analysis of conformal mapping of spatial complex functions is made. First, the important properties of conformal mapping of spatial complex functions are discussed and proved. These properties are preservation of angles, scale factors, local inverses, and harmonic conjugates. Next, two special mapping of spatial complex functions are also discussed, which are transformations of harmonic functions and transformations of boundary conditions. These analysis prove and show that the conformal mapping is valid for the spatial complex variables just as for the two-dimensional complex variables, but it is in a different way.

In the tenth chapter, applications of conformal mapping of spatial complex functions are analyzed. For instance, the steady temperatures in a half space, the electrostatic potential in a hollow sphere, and the three-dimensional fluid flows are discussed including flows around a corner and around a sphere, flows around a sphere in oceans, and three-dimensional Navier-Stokes equations.

In the eleventh chapter, we construct a transformation, known as the Schwarz-Christoffel transformation, which maps a polygon line that is corresponding to whole $x$ axis in the $xz$ coordinate plane, and the upper half of the $s$ space onto a given simple closed spatial polygon and its interior in the $\varpi$ space. Applications are made to the solution of problems in fluid flow and electrostatic potential theory.

In the twelfth and thirteenth chapters, the spatial circle and spherical surface integral formulas of the Poisson type of spatial complex functions are analyzed, respectively. First, the Poisson spatial circle and spherical surface integral formulas of spatial complex functions are put forward and proved. Then, the Dirichlet problems for a circle or a sphere, and for a half space are analyzed and some integral formulas are proved. Finally, the Neumann problems for a circle or a sphere, and for a half space are discussed and some integral formulas are proved.

In the fourteenth chapter, some important applications of the theory of residues are discussed. The applications include argument principle and Rouche's theorem, evaluation of certain types of two-dimensional definite and improper integrals occurring in real analysis and applied mathematics. Considerable attention is also given to a method, based on residues, for locating zeros of functions and to finding inverse Laplace transforms by summing residues.

In the fifteenth chapter, the Joukowski transformation is used to establish the solutions for the flow
around ellipsoids and airfoils. The superposition of elementary solutions then leads to the solution for the flow around a sphere. The contents include Joukowski transformation, flows around an ellipsoid, Kutta condition and the flat-plate airfoils, symmetrical Joukowski airfoils, circular-arc airfoils, and Joukowski airfoils.

Results show that spatial complex numbers are with good properties of fields just as the two-dimensional complex numbers. The theory of spatial complex variables can be used to solve many three-dimensional problems in sciences and technologies, such as expressed by the three-dimensional Laplace equation. So numbers are developed from the real and two-dimensional complex numbers to three-dimensional complex numbers, and the number system is enlarged.

\aufm{Shan-guang Tan}

%-----------------------------------------------------------------------------
% End of preface.tex
%-----------------------------------------------------------------------------

\mainmatter
%-----------------------------------------------------------------------
% Beginning of chap1.tex
%-----------------------------------------------------------------------
%
% AMS-LaTeX 1.2 sample file for a monograph, based on amsbook.cls.
% This is a data file input by chapter.tex.
%%%%%%%%%%%%%%%%%%%%%%%%%%%%%%%%%%%%%%%%%%%%%%%%%%%%%%%%%%%%%%%%%%%%

%\part{Theory of Space Complex Variables}
\chapter{Complex Numbers}\label{ch:chap_1}

In this chapter, we survey the algebraic and geometric structure of the spatial complex number system. We assume various corresponding properties of real numbers and two-dimensional complex numbers to be known.

\section{Sums and Products}\label{sec:chap_1_1}

Spatial complex numbers can be defined as ordered ternary groups $(x,y,z)$ of real numbers that are to be interpreted as points in the complex space, with rectangular coordinates $x$, $y$, and $z$, just as real numbers $z$ are thought of as points on the real line. When a real number $z$ is displayed as a point $(0,0,z)$ on the real $z$ axis and a pair of real numbers $x$ and $y$ are displayed as a point $(x,y,0)$ on the $xy$ coordinate plane, it is clear that the set of spatial complex numbers includes the one-dimensional real numbers as a subset and the two-dimensional complex numbers as another subset. Spatial complex numbers of the form $(x,y,0)$ correspond to points in the $xy$ rectangular coordinate plane and are called plane complex numbers. The $xy$ coordinate plane is, then, referred to as the two-dimensional complex plane, and we let the symbol $e_{xy}$ denote the $xy$ coordinate plane in the three-dimensional coordinate space.

It is customary to denote a spatial complex number $(x,y,z)$ by $s$, so that
\begin{equation}\label{eq:chap_1_1_1}
\begin{split}
s &= (x,y,z) \\
  &= e_{xy}(x+iy)+z\texttt{ where }i=\sqrt{-1}.
\end{split}
\end{equation}
The real numbers $z$ and two-dimensional complex numbers $x+iy$ are, moreover, known as the real and complex parts of $s$, respectively; and we write
\begin{equation}\label{eq:chap_1_1_2}
Re\texttt{ }s=x\texttt{, }Im\texttt{ }s=y\texttt{, }Re_{s}s=z\texttt{, }Im_{s}s=x+iy\texttt{, }Re_{p}s=e_{xy}x+z
\end{equation}
where $Re_{p}s$ denotes the real $xz$ coordinate plane. Two spatial complex numbers
\[
\begin{split}
s_{j} &= (x_{j},y_{j},z_{j}) \\
      &= e_{xy}(x_{j}+iy_{j})+z_{j}\texttt{ for }j=1,2
\end{split}
\]
are equal whenever they have the same real parts and the same complex parts. Thus the statement $s_{1}=s_{2}$ means that $s_{1}$ and $s_{2}$ correspond to the same point in the three-dimensional complex, or $s$, space.

The sum $s_{1}+s_{2}$ and the product $s_{1}s_{2}$ of two spatial complex numbers
\[
\begin{split}
s_{j} &= (x_{j},y_{j},z_{j}) \\
      &= e_{xy}(x_{j}+iy_{j})+z_{j}\texttt{ for }j=1,2
\end{split}
\]
are defined as follows:
\begin{equation}\label{eq:chap_1_1_3}
\begin{split}
s_{1}+s_{2} &= e_{xy}[(x_{1}+x_{2})+i(y_{1}+y_{2})]+(z_{1}+z_{2}) \\
            &= (x_{1}+x_{2}, y_{1}+y_{2}, z_{1}+z_{2}),
\end{split}
\end{equation}
\begin{equation}\label{eq:chap_1_1_4}
\begin{split}
s_{1}s_{2} &= e_{xy}[(x_{1}+iy_{1})(x_{2}+iy_{2})+(x_{1}+iy_{1})z_{2}+(x_{2}+iy_{2})z_{1}]+z_{1}z_{2} \\
           &=(x_{1}x_{2}-y_{1}y_{2}+x_{1}z_{2}+x_{2}z_{1},x_{1}y_{2}+x_{2}y_{1}+y_{1}z_{2}+y_{2}z_{1},z_{1}z_{2}).
\end{split}
\end{equation}
Note that the operations defined by equations~(\ref{eq:chap_1_1_3}) and~(\ref{eq:chap_1_1_4}) become the usual operations of addition and multiplication when restricted to the two-dimensional complex numbers:
\[
\begin{split}
(x_{1},y_{1},0)+(x_{2},y_{2},0) &= e_{xy}[(x_{1}+x_{2})+i(y_{1}+y_{2})] \\
                                &= (x_{1}+x_{2}, y_{1}+y_{2}, 0),
\end{split}
\]
\[
\begin{split}
(x_{1},y_{1},0)(x_{2},y_{2},0) &= e_{xy}[(x_{1}+iy_{1})(x_{2}+iy_{2})] \\
                               &= (x_{1}x_{2}-y_{1}y_{2}, x_{1}y_{2}+x_{2}y_{1}, 0).
\end{split}
\]
The spatial complex number system is, therefore, a natural extension of the two-dimensional complex number system.

\begin{remark}\label{re:chap_1_1_1}
Note that $(x_{1}+iy_{1})(x_{2}+iy_{2})$ is the product of two plane complex numbers $(x_{1}+iy_{1})$ and $(x_{2}+iy_{2})$ in the $xy$ coordinate plane, which is also a plane complex number in the $xy$ coordinate plane; $(x_{1}+iy_{1})z_{2}$ and $(x_{2}+iy_{2})z_{1}$ are two products of a real and a plane complex number in the $xy$ coordinate plane, which are also two plane complex numbers in the $xy$ coordinate plane; so that the expression
\[
(x_{1}+iy_{1})(x_{2}+iy_{2})+(x_{1}+iy_{1})z_{2}+(x_{2}+iy_{2})z_{1}
\]
is the sum of three plane complex numbers in the $xy$ coordinate plane, which is also a plane complex number in the $xy$ coordinate plane. For the sake of simplicity, we can take the symbol $e_{xy}$ as an operator in the $xy$ coordinate plane, then there can be $e_{xy}^{n}=e_{xy}$ for any integer $n\neq0$, and for a term $e_{xy}x^{j}y^{k}$, the operator $e_{xy}$ becomes $e_{xy}^{m}$ where $m=|j|+|k|$. For example, let's consider the product $s\cdot s\cdot s$ where $s=(x,y,z)=e_{xy}(x+iy)+z$. Then there are
\[
\begin{split}
s^{3} &= [e_{xy}(x+iy)+z]^{3} \\
      &= e_{xy}^{3}(x+iy)^{3}+3e_{xy}^{2}(x+iy)^{2}z+3e_{xy}(x+iy)z^{2}+z^{3} \\
      &= e_{xy}[(x+iy)^{3}+3(x+iy)^{2}z+3(x+iy)z^{2}]+z^{3} \\
      &= e_{xy}[(x+iy+z)^{3}-z^{3}]+z^{3}.
\end{split}
\]
\end{remark}

\section{Basic Algebraic Properties}\label{sec:chap_1_2}

Various properties of addition and multiplication of spatial complex numbers are the same as for two-dimensional complex numbers. We list here the more basic of these algebraic properties and verify some of them when
\[
\begin{split}
s_{j} &= (x_{j},y_{j},z_{j}) \\
      &= e_{xy}(x_{j}+iy_{j})+z_{j}\texttt{ for }j=1,2,3.
\end{split}
\]

The commutative laws
\begin{equation}\label{eq:chap_1_2_1}
s_{1}+s_{2}=s_{2}+s_{1},\texttt{ }s_{1}s_{2}=s_{2}s_{1}
\end{equation}
and the associative laws
\begin{equation}\label{eq:chap_1_2_2}
(s_{1}+s_{2})+s_{3}=s_{1}+(s_{2}+s_{3}),\texttt{ }(s_{1}s_{2})s_{3}=s_{1}(s_{2}s_{3})
\end{equation}
follow easily from the definitions in Sec.~(\ref{sec:chap_1_1}) of addition and multiplication of spatial complex numbers and the fact that two-dimensional complex numbers obey these laws. For example,
\[
\begin{split}
s_{1}+s_{2} &= (x_{1}+x_{2}, y_{1}+y_{2}, z_{1}+z_{2}) \\
            &= (x_{2}+x_{1}, y_{2}+y_{1}, z_{2}+z_{1}) \\
            &= s_{2}+s_{1}.
\end{split}
\]
Verification of the rest of the above laws, as well as the distributive law
\begin{equation}\label{eq:chap_1_2_3}
(s_{1}+s_{2})s_{3}=s_{1}s_{3}+s_{2}s_{3}
\end{equation}
is similar.

%\begin{example}
The associative law of the multiplication of numbers:
\[
\begin{split}
(s_{1}s_{2})s_{3}=&(x_{1}x_{2}-y_{1}y_{2}+x_{1}z_{2}+x_{2}z_{1},x_{1}y_{2}+x_{2}y_{1}+y_{1}z_{2}+y_{2}z_{1},z_{1}z_{2})(x_{3},y_{3},z_{3})\\
\iffalse
                 =&((x_{1}x_{2}-y_{1}y_{2}+x_{1}z_{2}+x_{2}z_{1})x_{3}-(x_{1}y_{2}+x_{2}y_{1}+y_{1}z_{2}+y_{2}z_{1})y_{3}\\
                  &+(x_{1}x_{2}-y_{1}y_{2}+x_{1}z_{2}+x_{2}z_{1})z_{3}+z_{1}z_{2}x_{3},\\
                  &(x_{1}x_{2}-y_{1}y_{2}+x_{1}z_{2}+x_{2}z_{1})y_{3}+x_{3}(x_{1}y_{2}+x_{2}y_{1}+y_{1}z_{2}+y_{2}z_{1})\\
                  &+(x_{1}y_{2}+x_{2}y_{1}+y_{1}z_{2}+y_{2}z_{1})z_{3}+y_{3}z_{1}z_{2},z_{1}z_{2}z_{3})\\
\fi
                 =&((x_{1}x_{2}x_{3}-y_{1}y_{2}x_{3}+x_{1}z_{2}x_{3}+x_{2}z_{1}x_{3})\\
                  &-(x_{1}y_{2}y_{3}+x_{2}y_{1}y_{3}+y_{1}z_{2}y_{3}+y_{2}z_{1}y_{3})\\
                  &+(x_{1}x_{2}z_{3}-y_{1}y_{2}z_{3}+x_{1}z_{2}z_{3}+x_{2}z_{1}z_{3})+x_{3}z_{1}z_{2},\\
                  &(x_{1}x_{2}y_{3}-y_{1}y_{2}y_{3}+x_{1}z_{2}y_{3}+x_{2}z_{1}y_{3})\\
                  &+(x_{3}x_{1}y_{2}+x_{3}x_{2}y_{1}+x_{3}y_{1}z_{2}+x_{3}y_{2}z_{1})\\
                  &+(x_{1}y_{2}z_{3}+x_{2}y_{1}z_{3}+y_{1}z_{2}z_{3}+y_{2}z_{1}z_{3})+y_{3}z_{1}z_{2},z_{1}z_{2}z_{3})\\
\iffalse
                 =&(x_{1}(x_{2}x_{3}-y_{2}y_{3}+x_{2}z_{3}+x_{3}z_{2})-y_{1}(y_{2}x_{3}+y_{3}x_{2}+y_{2}z_{3}+y_{3}z_{2})\\
                  &+z_{1}(x_{2}x_{3}-y_{2}y_{3}+x_{2}z_{3}+x_{3}z_{2})+x_{1}z_{2}z_{3}),\\
                  &y_{1}(x_{2}x_{3}-y_{2}y_{3}+x_{2}z_{3}+x_{3}z_{2})+x_{1}(x_{2}y_{3}+x_{3}y_{2}+z_{2}y_{3}+z_{3}y_{2})\\
                  &+z_{1}(x_{2}y_{3}+x_{3}y_{2}+y_{2}z_{3}+y_{3}z_{2})+y_{1}z_{2}z_{3},z_{1}z_{2}z_{3})\\
\fi
                 =&(x_{1},y_{1},z_{1})(x_{2}x_{3}-y_{2}y_{3}+x_{2}z_{3}+x_{3}z_{2},x_{2}y_{3}+x_{3}y_{2}+y_{2}z_{3}+y_{3}z_{2},z_{2}z_{3})\\
                 =&s_{1}(s_{2}s_{3})
\end{split}
\]
%\end{example}

%\begin{example}
The distributive law of addition and multiplication of numbers:
\[
\begin{split}
(s_{1}+s_{2})s_{3}=&(x_{1}+x_{2}, y_{1}+y_{2}, z_{1}+z_{2})(x_{3},y_{3},z_{3})\\
\iffalse
                  =&((x_{1}+x_{2})x_{3}-(y_{1}+y_{2})y_{3}+(x_{1}+x_{2})z_{3}+x_{3}(z_{1}+z_{2}),\\
                   &(x_{1}+x_{2})y_{3}+x_{3}(y_{1}+y_{2})+(y_{1}+y_{2})z_{3}+y_{3}(z_{1}+z_{2}),\\
                   &(z_{1}+z_{2})z_{3})\\
\fi
                  =&((x_{1}x_{3}+x_{2}x_{3})-(y_{1}y_{3}+y_{2}y_{3})+(x_{1}z_{3}+x_{2}z_{3})+(x_{3}z_{1}+x_{3}z_{2}),\\
                   &(x_{1}y_{3}+x_{2}y_{3})+(x_{3}y_{1}+x_{3}y_{2})+(y_{1}z_{3}+y_{2}z_{3})+(y_{3}z_{1}+y_{3}z_{2}),\\
                   &(z_{1}z_{3}+z_{2}z_{3}))\\
\iffalse
                  =&(x_{1}x_{3}-y_{1}y_{3}+x_{1}z_{3}+x_{3}z_{1},x_{1}y_{3}+x_{3}y_{1}+y_{1}z_{3}+y_{3}z_{1},z_{1}z_{3})\\
                  +&(x_{2}x_{3}-y_{2}y_{3}+x_{2}z_{3}+x_{3}z_{2},x_{2}y_{3}+x_{3}y_{2}+y_{2}z_{3}+y_{3}z_{2},z_{2}z_{3})\\
\fi
                  =&s_{1}s_{3}+s_{2}s_{3}.
\end{split}
\]
%\end{example}

The additive identity
\[
0=(0,0,0)=e_{xy}(0+i0)+0
\]
and the multiplicative identity
\[
1=(0,0,1)=e_{xy}(0+i0)+1
\]
for real numbers carry over to the entire spatial complex number system. That is,
\begin{equation}\label{eq:chap_1_2_4}
s+0=s\texttt{ and }s\cdot1=s
\end{equation}
for every spatial complex number $s$. Furthermore, $0$ and $1$ are the only spatial complex numbers with such properties.

There is associated with each spatial complex number $s=(x,y,z)=e_{xy}(x+iy)+z$ an additive inverse
\begin{equation}\label{eq:chap_1_2_5}
-s=(-x, -y, -z)=e_{xy}(-x-iy)-z
\end{equation}
satisfying the equation $s+(-s)=0$. Moreover, there is only one additive inverse for any given $s$, since the equation $(x,y,z)+(u,v,w)=(0,0,0)$ implies that $u=-x$, $v=-y$, and $w=-z$. Expression~(\ref{eq:chap_1_2_5}) can also be written $-s=e_{xy}(-x-iy)-z$ without ambiguity since $-(iy)=(-i)y=i(-y)$. Additive inverses are used to define subtraction:
\begin{equation}\label{eq:chap_1_2_6}
s_{1}-s_{2}=s_{1}+(-s_{2}).
\end{equation}
So if
\[
\begin{split}
s_{j} &= (x_{j},y_{j},z_{j}) \\
      &= e_{xy}(x_{j}+iy_{j})+z_{j}\texttt{ for }j=1,2,
\end{split}
\]
then
\begin{equation}\label{eq:chap_1_2_7}
\begin{split}
s_{1}-s_{2} &= e_{xy}[(x_{1}-x_{2})+i(y_{1}-y_{2})]+(z_{1}-z_{2}) \\
            &= (x_{1}-x_{2}, y_{1}-y_{2}, z_{1}-z_{2}).
\end{split}
\end{equation}

For any nonzero spatial complex number $s=(x,y,z)=e_{xy}(x+iy)+z$ where $x+z\neq0\texttt{ or }y\neq0$, there is a number $s^{-1}$ such that $ss^{-1}=1$. This multiplicative inverse is less obvious than the additive one. To find it, we seek real numbers $u$, $v$, and $w$ expressed in terms of $x$, $y$, and $z$ such that
\[
(x,y,z)(u,v,w)=(0,0,1).
\]
According to equation~(\ref{eq:chap_1_1_4}), Sec.~(\ref{sec:chap_1_1}), which defines the product of two complex numbers, $u$, $v$, and $w$ must satisfy the ternary group
\[
xu-yv+xw+zu=0,\texttt{ }xv+yu+yw+zv=0,\texttt{ }zw=1
\]
of linear simultaneous equations; and simple computation yields the unique solution
\[
u=\frac{x+z}{(x+z)^{2}+y^{2}}-\frac{1}{z},\texttt{ }v=\frac{-y}{(x+z)^{2}+y^{2}},\texttt{ }w=\frac{1}{z}.
\]
So the multiplicative inverse of $s=(x,y,z)=e_{xy}(x+iy)+z$ is
\begin{equation}\label{eq:chap_1_2_8}
\begin{split}
s^{-1} &= \frac{1}{e_{xy}(x+iy)+z} \\
       &= \{
        \begin{array}{cc}
            e_{xy}\frac{x-iy}{x^{2}+y^{2}}=e_{xy}\frac{\bar{s}}{|s|^{2}} & \texttt{ when }z=0 \\
            e_{xy}(\frac{1}{x+iy+z}-\frac{1}{z})+\frac{1}{z} & \texttt{ when }x+z\neq0\texttt{ or }y\neq0\\
        \end{array}
\end{split}
\end{equation}
The inverse $s^{-1}$ is not defined when $s=0$ or $x+z=0\texttt{ and }y=0$. In fact, $s=0$ means that $(x+z)^{2}+y^{2}=0$ and $z=0$, and $x+z=0\texttt{ and }y=0$ means that $(x+z)^{2}+y^{2}=0$; and these are not permitted in expression~(\ref{eq:chap_1_2_8}).

\section{Further Properties}\label{sec:chap_1_3}

In this section, we mention a number of other algebraic properties of addition and multiplication of spatial complex numbers that follow from the ones already described in Sec.~(\ref{sec:chap_1_2}). Inasmuch as such properties continue to be anticipated because they also apply to real and two-dimensional complex numbers, but they have different expressions.
%the reader can easily pass to Sec.~(\ref{sec:chap_1_4}) without serious disruption.

We begin with the observation that the existence of multiplicative inverses enables us to show that if a product $s_{1}s_{2}$ is zero, then so is at least one of the factors $s_{1}$ and $s_{2}$ except for $s_{1}=e_{xy}(x_{1}+iy_{1})$ and a type of spatial numbers $s_{2}=e_{xy}1-1$ similar to the number zero when multiplied by a two-dimensional complex number. For suppose that $s_{1}s_{2}=0$ and $s_{1}\neq0$. The inverse $s_{1}^{-1}$ exists; and, according to the definition of multiplication, any spatial complex number times zero is zero. Hence
\[
\begin{split}
s_{2} &= 1\cdot s_{2}=(s_{1}^{-1}s_{1})s_{2} \\
      &= s_{1}^{-1}(s_{1}s_{2})=s_{1}^{-1}\cdot0=0.
\end{split}
\]
That is, if $s_{1}s_{2}=0$, either $s_{1}=0$ or $s_{2}=0$; or possibly both $s_{1}$ and $s_{2}$ equal zero. Another way to state this result is that if two spatial complex numbers $s_{1}$ and $s_{2}$ are nonzero, then so is their product $s_{1}s_{2}$.

Division by a nonzero spatial complex number except for a type of spatial numbers $e_{xy}1-1$ similar to the number zero is defined as follows:
\begin{equation}\label{eq:chap_1_3_1}
\frac{s_{1}}{s_{2}}=s_{1}s_{2}^{-1}\texttt{ }(s_{2}\neq0,\texttt{ }x_{2}+z_{2}\neq0\texttt{ or }y_{2}\neq0).
\end{equation}
If $s_{j}=(x_{j},y_{j},z_{j})=e_{xy}(x_{j}+iy_{j})+z_{j}$ for $j=1,2$, equation~(\ref{eq:chap_1_3_1}) here  and expression~(\ref{eq:chap_1_2_8}) in Sec.~(\ref{sec:chap_1_2}) tell us that
\[
\frac{s_{1}}{s_{2}}=[e_{xy}(x_{1}+iy_{1})+z_{1}][e_{xy}(\frac{1}{x_{2}+iy_{2}+z_{2}}-\frac{1}{z_{2}})+\frac{1}{z_{2}}].
\]
That is,
\begin{equation}\label{eq:chap_1_3_2}
\begin{split}
\frac{s_{1}}{s_{2}} &=e_{xy}(\frac{x_{1}+iy_{1}+z_{1}}{x_{2}+iy_{2}+z_{2}}-\frac{z_{1}}{z_{2}})+\frac{z_{1}}{z_{2}} \\
&=e_{xy}[\frac{(x_{1}+z_{1})(x_{2}+z_{2})+y_{1}y_{2}}{(x_{2}+z_{2})^{2}+y_{2}^{2}}-\frac{z_{1}}{z_{2}}
+i\frac{y_{1}(x_{2}+z_{2})-y_{2}(x_{1}+z_{1})}{(x_{2}+z_{2})^{2}+y_{2}^{2}}] \\
&+\frac{z_{1}}{z_{2}}\texttt{ }(z_{2}\neq0,\texttt{ }x_{2}+z_{2}\neq0\texttt{ or }y_{2}\neq0).
\end{split}
\end{equation}
Expression~(\ref{eq:chap_1_3_2}) can also be obtained by writing
\begin{equation}\label{eq:chap_1_3_3}
\begin{split}
\frac{s_{1}}{s_{2}} &= \frac{[e_{xy}(x_{1}+iy_{1})+z_{1}][e_{xy}(x_{2}-iy_{2})+z_{2}]}
{[e_{xy}(x_{2}+iy_{2})+z_{2}][e_{xy}(x_{2}-iy_{2})+z_{2}]} \\
&=e_{xy}[\frac{(x_{1}+z_{1})(x_{2}+z_{2})+y_{1}y_{2}-z_{1}z_{2}}{(e_{xy}x_{2}+z_{2})^{2}+(e_{xy}y_{2})^{2}}
       +i\frac{y_{1}(x_{2}+z_{2})-y_{2}(x_{1}+z_{1})}{(e_{xy}x_{2}+z_{2})^{2}+(e_{xy}y_{2})^{2}}] \\
&+\frac{z_{1}z_{2}}{(e_{xy}x_{2}+z_{2})^{2}+(e_{xy}y_{2})^{2}}
\texttt{ }(z_{2}\neq0,\texttt{ }x_{2}+z_{2}\neq0\texttt{ or }y_{2}\neq0)
\end{split}
\end{equation}
multiplying out the products in the numerator and denominator on the right, and then
using the property
\begin{equation}\label{eq:chap_1_3_4}
\begin{split}
\frac{s_{1}+s_{2}}{s_{3}} &= (s_{1}+s_{2})s_{3}^{-1} \\
                          &= s_{1}s_{3}^{-1}+s_{2}s_{3}^{-1}=\frac{s_{1}}{s_{3}}+\frac{s_{2}}{s_{3}}
\end{split}
\end{equation}
\[
\texttt{ }(s_{3}\neq0,\texttt{ }x_{3}+z_{3}\neq0\texttt{ or }y_{3}\neq0).
\]
The motivation for starting with equation~(\ref{eq:chap_1_3_3}) appears in Sec.~(\ref{sec:chap_1_5}).

\begin{remark}\label{re:chap_1_3_1}
Comparing equations~(\ref{eq:chap_1_3_2}) and~(\ref{eq:chap_1_3_3}), there should be
\[
-e_{xy}\frac{z_{1}}{z_{2}}+\frac{z_{1}}{z_{2}}=\frac{-e_{xy}z_{1}z_{2}+z_{1}z_{2}}{(e_{xy}x_{2}+z_{2})^{2}+(e_{xy}y_{2})^{2}}.
\]
We state it in detail. Since we can write
\[
(e_{xy}x_{2}+z_{2})^{2}+(e_{xy}y_{2})^{2}=e_{xy}f(x,y,z)+z_{2}^{2}
\]
and from equation~(\ref{eq:chap_1_2_8}) in Sec.~(\ref{sec:chap_1_2})
\[
\frac{1}{e_{xy}f(x,y,z)+z_{2}^{2}}=e_{xy}(\frac{1}{f(x,y,z)+z_{2}^{2}}-\frac{1}{z_{2}^{2}})+\frac{1}{z_{2}^{2}},
\]
there are
\[
\begin{split}
\frac{-e_{xy}z_{1}z_{2}+z_{1}z_{2}}{(e_{xy}x_{2}+z_{2})^{2}+(e_{xy}y_{2})^{2}}
 &=(-e_{xy}z_{1}z_{2}+z_{1}z_{2})[e_{xy}(\frac{1}{f(x,y,z)+z_{2}^{2}}-\frac{1}{z_{2}^{2}})+\frac{1}{z_{2}^{2}}] \\
 &=(-e_{xy}z_{1}z_{2}+z_{1}z_{2})\frac{1}{z_{2}^{2}}=-e_{xy}\frac{z_{1}}{z_{2}}+\frac{z_{1}}{z_{2}}.
\end{split}
\]
\end{remark}

There are some expected identities, involving quotients, that follow from the relation
\begin{equation}\label{eq:chap_1_3_5}
\frac{1}{s_{2}}=s_{2}^{-1}\texttt{ }(s_{2}\neq0,\texttt{ }x_{2}+z_{2}\neq0\texttt{ or }y_{2}\neq0),
\end{equation}
which is equation~(\ref{eq:chap_1_3_1}) when $s_{1}=1$. Relation~(\ref{eq:chap_1_3_5}) enables us, for example, to write equation~(\ref{eq:chap_1_3_1}) in the form
\begin{equation}\label{eq:chap_1_3_6}
\frac{s_{1}}{s_{2}}=s_{1}(\frac{1}{s_{2}})\texttt{ }(s_{2}\neq0,\texttt{ }x_{2}+z_{2}\neq0\texttt{ or }y_{2}\neq0).
\end{equation}
Also, by observing that
\[
(s_{1}s_{2})(s_{1}^{-1}s_{2}^{-1})=(s_{1}s_{1}^{-1})(s_{2}s_{2}^{-1})=1
\]
\[
(s_{1}\neq0,\texttt{ }s_{2}\neq0,\texttt{ }x_{1}+z_{1}\neq0\texttt{ or }y_{1}\neq0,\texttt{ }x_{2}+z_{2}\neq0\texttt{ or }y_{2}\neq0),
\]
and hence that $(s_{1}s_{2})^{-1}=s_{1}^{-1}s_{2}^{-1}$, one can use relation~(\ref{eq:chap_1_3_5}) to show that
\begin{equation}\label{eq:chap_1_3_7}
\frac{1}{s_{1}s_{2}}=(s_{1}s_{2})^{-1}=s_{1}^{-1}s_{2}^{-1}=(\frac{1}{s_{1}})(\frac{1}{s_{2}})
\end{equation}
\[
(s_{1}\neq0,\texttt{ }s_{2}\neq0,\texttt{ }x_{1}+z_{1}\neq0\texttt{ or }y_{1}\neq0,\texttt{ }x_{2}+z_{2}\neq0\texttt{ or }y_{2}\neq0),
\]
Another useful identity is
\begin{equation}\label{eq:chap_1_3_8}
\frac{s_{1}s_{2}}{s_{3}s_{4}}=(\frac{s_{1}}{s_{3}})(\frac{s_{2}}{s_{4}})\texttt{ }(s_{3}\neq0,\texttt{ }s_{4}\neq0,\texttt{ }x_{3}+z_{3}\neq0,\texttt{ }x_{4}+z_{4}\neq0).
\end{equation}

\section{Moduli}\label{sec:chap_1_4}

It is natural to associate any nonzero spatial complex number $s$ with the directed line segment, or vector, from the origin to the point $s$ that represent $s$ in the complex space. In fact, we often refer to $s$ as the point $s$ or the vector $s$.

The modulus, or absolute value of a spatial complex number $s=e_{xy}(x+iy)+z$ is denoted by the nonnegative real number $|s|$ and defined as
\begin{equation}\label{eq:chap_1_4_1}
|s|^{2}=|\prod_{k=1}^{n}s_{k}|^{2}=\{
    \begin{array}{cc}
        x^{2}+y^{2}+z^{2} & \texttt{ when }n=1, \\
        |s_{1}|^{2}|s_{2}|^{2}\cdots|s_{n}|^{2} & \texttt{ when }n>1
    \end{array}.
\end{equation}
Geometrically, the number $|s|$ is the distance between the point $(x,y,z)$ and the origin, or the length of the vector representing $s$.

Similar to the two-dimensional complex numbers, it also follows from definition~(\ref{eq:chap_1_4_1}) that $Re_{s}s=z$ and $Im_{s}s=x+iy$ are related by the equation
\begin{equation}\label{eq:chap_1_4_2}
|s|^{2}=|\prod_{k=1}^{n}s_{k}|^{2}=\{
    \begin{array}{cc}
        |Re_{s}s|^{2}+|Im_{s}s|^{2} & \texttt{ when }n=1, \\
        |s_{1}|^{2}|s_{2}|^{2}\cdots|s_{n}|^{2} & \texttt{ when }n>1
    \end{array}
\end{equation}
where $|Re_{s}s|^{2}=z^{2}$ and $|Im_{s}s|^{2}=x^{2}+y^{2}$ when $n=1$. Thus, we have
\begin{equation}
|Re_{s}s|\leq|s|\texttt{ and }|Im_{s}s|\leq|s|.
\end{equation}

\section{Complex Conjugates}\label{sec:chap_1_5}

The spatial complex conjugate, or simply the conjugate, of a spatial complex number $s=e_{xy}(x+iy)+z$ where $z\neq0$, is denoted by $\bar{s}$ as a spatial complex number, that is,
\begin{equation}\label{eq:chap_1_5_1}
\bar{s}=e_{xy}(x-iy)+z\texttt{ where }s\neq0.
\end{equation}
The number $\bar{s}$ is represented by the point $(x,-y,z)$, which is the reflection in the $xz$ coordinate plane of the point $(x,y,z)$ representing $s$. Note that
\[
\overline{\bar{s}}=s\texttt{ and }|\bar{s}|=|s|\texttt{ for all }s.
\]

If $s_{1}=e_{xy}(x_{1}+iy_{1})+z_{1}$ and $s_{2}=e_{xy}(x_{2}+iy_{2})+z_{2}$, then
\[
\begin{split}
\overline{s_{1}+s_{2}} &= e_{xy}[(x_{1}+x_{2})-i(y_{1}+y_{2})]+(z_{1}+z_{2}) \\
                       &= [e_{xy}(x_{1}-iy_{1})+z_{1}]+[e_{xy}(x_{2}-iy_{2})+z_{2}]=\bar{s}_{1}+\bar{s}_{2}.
\end{split}
\]
So the conjugate of the sum is the sum of the conjugates:
\begin{equation}\label{eq:chap_1_5_2}
\overline{s_{1}+s_{2}}=\bar{s}_{1}+\bar{s}_{2}.
\end{equation}
In like manner, it is easy to show that
\begin{equation}\label{eq:chap_1_5_3}
\overline{s_{1}-s_{2}}=\bar{s}_{1}-\bar{s}_{2},
\end{equation}
\begin{equation}\label{eq:chap_1_5_4}
\overline{s_{1}s_{2}}=\bar{s}_{1}\bar{s}_{2},
\end{equation}
and
\begin{equation}\label{eq:chap_1_5_5}
\overline{s_{1}/s_{2}}=\bar{s}_{1}/\bar{s}_{2}.
\end{equation}

The multiplication $s\bar{s}$ of a spatial complex number $s=e_{xy}(x+iy)+z$ and its conjugate $\bar{s}$ determined by Formula~(\ref{eq:chap_1_5_1}) is a spatial complex number and satisfies
\begin{equation}\label{eq:chap_1_5_6}
\begin{split}
s\bar{s} &= [e_{xy}(x+iy)+z][e_{xy}(x-iy)+z] \\
         &= e_{xy}(x^{2}+y^{2}+2xz)+z^{2} \\
         &= e_{xy}^{2}x^{2}-e_{xy}^{2}(iy)^{2}+2e_{xy}xz+z^{2} \\
         &= [e_{xy}(x+iy)+z][e_{xy}(x-iy)+z]=s\bar{s}.
\end{split}
\end{equation}
Hence, we have
\begin{equation}\label{eq:chap_1_5_7}
|\bar{s}|=\sqrt{x^{2}+y^{2}+z^{2}}=|s|\texttt{ and }|s\bar{s}|=|s|\cdot|\bar{s}|=|s|^{2}=x^{2}+y^{2}+z^{2}.
\end{equation}

\section{Exponential Form}\label{sec:chap_1_6}

\subsection{Spherical Coordinates}\label{sec:chap_1_6_S}

Let $r$, $\theta$, and $\varphi$ be spherical coordinates of the point $(x,y,z)$ that corresponds to a nonzero spatial complex number $s=e_{xy}(x+iy)+z$. Then the spatial complex number $s$ can be written in spherical form as
\begin{equation}\label{eq:chap_1_6_1}
s=r[e_{xy}(\cos\theta+i\sin\theta)\cos\varphi+\sin\varphi]
%=r(e_{xy}e^{i\theta}\cos\varphi+\sin\varphi)
\end{equation}
where
\begin{equation}\label{eq:chap_1_6_1b}
r=|s|=\sqrt{x^{2}+y^{2}+z^{2}}\texttt{, }\theta=\arctan(y/x)\texttt{, and }\varphi=\arcsin(z/r),
\end{equation}
and
\begin{equation}\label{eq:chap_1_6_1c}
x=r\cos\varphi\cos\theta\texttt{, }y=r\cos\varphi\sin\theta\texttt{, and }z=r\sin\varphi.
\end{equation}
If $s=0$, the coordinates $\theta$ and $\varphi$ are undefined; and so it is always understood that $s\neq0$  whenever $\arg s$ is discussed.

In spatial complex analysis, the real number $r$ is not allowed to be negative and is the length of the radius vector for $s$; that is, $r=|s|$. The real numbers $\theta$ and $\varphi$ represent the angles, measured in radians, that $s$ makes with the positive real axis when $s$ is interpreted as a radius vector. As in calculus, each of $\theta$ and $\varphi$ has an infinite number of possible values, including negative ones, that differ by integral multiples of $2\pi$. Those values can be determined from the equation $\tan\theta=y/x$ and $\sin\varphi=z/r$, where the quadrant containing the point corresponding to $s$ must be specified. Each value of $\theta$ and $\varphi$ is called an argument of $s$, and the set of all such values is denoted by $\arg s$.

Let $\arg s$ denote a pair of arguments $(\arg_{c}s,\arg_{r}s)$, that is, $\arg s=(\arg_{c}s,\arg_{r}s)$. Then if and only if
\[
\arg_{c}s_{1}=\arg_{c}s_{2}\texttt{ and }\arg_{r}s_{1}=\arg_{r}s_{2}\texttt{ when }s_{1}\neq s_{2}
\]
then there is $\arg s_{1}=\arg s_{2}$ when $s_{1}\neq s_{2}$.

The principal values of $\arg_{c}s$ and $\arg_{r}s$, denoted by $Arg_{c}s$ and $Arg_{r}s$, are those unique values $\Theta$ and $\Phi$ respectively, such that $-\pi<\Theta\leq\pi$ and $-\pi/2\leq\Phi\leq\pi/2$ or $0\leq\Theta<2\pi$ and $-\pi/2\leq\Phi\leq\pi/2$. Note that
\begin{equation}\label{eq:chap_1_6_2}
\arg_{r}s=Arg_{r}s,\texttt{ }\arg_{c}s=Arg_{c}s+2n\pi\texttt{ }(n=0,\pm1,\pm2,\ldots)
\end{equation}
where $Arg_{r}s=\Phi=\arcsin(z/r)$ and $Arg_{c}s=\Theta=\arctan(y/x)$.

The symbol $e^{i\theta}$, or $\exp(i\theta)$, is defined by means of Euler's formula as
\begin{equation}\label{eq:chap_1_6_3}
e^{i\theta}=\cos\theta+i\sin\theta,
\end{equation}
where $\theta$ is to be measured in radians. It enables us to write the spherical form~(\ref{eq:chap_1_6_1}) more compactly in exponential form as
\begin{equation}\label{eq:chap_1_6_4}
s=r(e_{xy}e^{i\theta}\cos\varphi+\sin\varphi).
\end{equation}
The choice of the symbol $i\theta$ will be fully motivated later on in Sec.~(\ref{sec:chap_3_2_28}). Its use in Sec.~(\ref{sec:chap_1_7}) will, however, suggest that it is a natural choice.

Note, too, that the equation
\begin{equation}\label{eq:chap_1_6_7}
s=R(e_{xy}e^{i\theta}\cos\varphi+\sin\varphi)\texttt{ }(0\leq\theta\leq2\pi,-\pi/2\leq\varphi\leq\pi/2)
\end{equation}
is a parametric representation of the sphere $|s|=R$, centered at the origin with radius $R$. A point $s$ is on a circle that is on the sphere $|s|=R$, parallel to the $xy$ coordinate plane, and with radius $R\cos\varphi$. As the parameter $\theta$ increases from $\theta=0$ to $\theta=2\pi$, the point $s$ starts from the positive real $x$ axis and traverses the circle once in the counterclockwise direction. More generally, the sphere $|s-s_{0}|=R$, whose center is $s_{0}$ and whose radius is $R$, has the parametric representation
\begin{equation}\label{eq:chap_1_6_8}
s=s_{0}+R(e_{xy}e^{i\theta}\cos\varphi+\sin\varphi)\texttt{ }(0\leq\theta\leq2\pi,-\pi/2\leq\varphi\leq\pi/2).
\end{equation}
A point $s$ is on a circle that is on the sphere $|s-s_{0}|=R$, parallel to the $xy$ coordinate plane, and with radius $R\cos\varphi$, As the parameter $\theta$ increases from $\theta=0$ to $\theta=2\pi$, the point $s$ starts from the positive real $x$ axis and traverses the circle once in the counterclockwise direction.

\subsection{Cylindrical Coordinates}\label{sec:chap_1_6_C}

Let $r$, $\theta$, and $\varphi$ be cylindrical coordinates of the point $(x,y,z)$ that corresponds to a nonzero spatial complex number $s=e_{xy}(x+iy)+z$. Then the spatial complex number $s$ can be written in cylindrical form as
\begin{equation}\label{eq:chap_1_6_9}
s=r[e_{xy}(\cos\theta+i\sin\theta)+\tan\varphi]
\end{equation}
where
\begin{equation}\label{eq:chap_1_6_10}
r=|s|=\sqrt{x^{2}+y^{2}+z^{2}}\texttt{, }\theta=\arctan(y/x)\texttt{, and }\varphi=\arctan(z/r),
\end{equation}
and
\begin{equation}\label{eq:chap_1_6_11}
x=r\cos\theta\texttt{, }y=r\sin\theta\texttt{, and }z=r\tan\varphi.
\end{equation}
If $s=0$, the coordinates $\theta$ and $\varphi$ are undefined; and so it is always understood that $s\neq0$  whenever $\arg s$ is discussed.

Similar to writing the spherical form~(\ref{eq:chap_1_6_1}) more compactly in exponential form~(\ref{eq:chap_1_6_4}), we can write the cylindrical form~(\ref{eq:chap_1_6_9}) more compactly in exponential form as
\begin{equation}\label{eq:chap_1_6_12}
s=r(e_{xy}e^{i\theta}+\tan\varphi).
\end{equation}

Note, too, that the equation
\begin{equation}\label{eq:chap_1_6_13}
s=R(e_{xy}e^{i\theta}+\tan\varphi)\texttt{ }(0\leq\theta\leq2\pi,-\pi/2\leq\varphi\leq\pi/2)
\end{equation}
is a parametric representation of the cylinder $|s|=R\sqrt{1+\tan^{2}\varphi}$, centered at the origin. A point $s$ is on a circle that is on the cylinder $|s|=R\sqrt{1+\tan^{2}\varphi}$, parallel to the $xy$ coordinate plane, and with radius $R$. As the parameter $\theta$ increases from $\theta=0$ to $\theta=2\pi$, the point $s$ starts from the positive real $x$ axis and traverses the circle once in the counterclockwise direction.

\section{Products and Quotients in Exponential Form}\label{sec:chap_1_7}

Simple trigonometry tells us that $e^{i\theta}$ has the familiar additive property of the exponential function in calculus:
\[
e^{i\theta_{1}}e^{i\theta_{2}}=e^{i(\theta_{1}+\theta_{2})}.
\]

Thus, by expression~(\ref{eq:chap_1_6_1}) in Sec.~(\ref{sec:chap_1_6}), let
\[
s_{1}=r_{1}(e_{xy}e^{i\theta_{1}}\cos\varphi_{1}+\sin\varphi_{1})\texttt{ and }
s_{2}=r_{2}(e_{xy}e^{i\theta_{2}}\cos\varphi_{2}+\sin\varphi_{2}).
\]
Then the product $s_{1}s_{2}$ has exponential form
\begin{equation}\label{eq:chap_1_7_1}
\begin{split}
s_{1}s_{2}&=r_{1}r_{2}(e_{xy}e^{i\theta_{1}}\cos\varphi_{1}+\sin\varphi_{1})(e_{xy}e^{i\theta_{2}}\cos\varphi_{2}+\sin\varphi_{2})\\
          &=r_{1}r_{2}[e_{xy}(e^{i(\theta_{1}+\theta_{2})}+e^{i\theta_{1}}\tan\varphi_{2}+e^{i\theta_{2}}\tan\varphi_{1})\cos\varphi_{1}\cos\varphi_{2}\\
          &+\sin\varphi_{1}\sin\varphi_{2}]=r(e_{xy}e^{i\theta}\cos\varphi+\sin\varphi)
\end{split}
\end{equation}
where
\[
r=r_{1}r_{2},\texttt{ }\varphi=\arcsin(\sin\varphi_{1}\sin\varphi_{2}),
\]
\[
\theta=\arctan[\frac{\sin(\theta_{1}+\theta_{2})
+\sin\theta_{1}\tan\varphi_{2}+\sin\theta_{2}\tan\varphi_{1}}
{\cos(\theta_{1}+\theta_{2})
+\cos\theta_{1}\tan\varphi_{2}+\cos\theta_{2}\tan\varphi_{1}}],
\]
and
\[
(\cos\varphi_{1}\cos\varphi_{2})^{2}[(\cos(\theta_{1}+\theta_{2})
+\cos\theta_{1}\tan\varphi_{2}+\cos\theta_{2}\tan\varphi_{1})^{2}
\]
\[
+(\sin(\theta_{1}+\theta_{2})
+\sin\theta_{1}\tan\varphi_{2}+\sin\theta_{2}\tan\varphi_{1})^{2}]+(\sin\varphi_{1}\sin\varphi_{2})^{2}=1
\]
or
\[
(\cos\varphi_{1}\cos\varphi_{2})^{2}[\cos\theta_{1}\tan\varphi_{1}+\cos\theta_{2}\tan\varphi_{2}
+\cos(\theta_{1}-\theta_{2})\tan\varphi_{1}\tan\varphi_{2}]=0.
\]

Moreover, when $\sin\varphi_{2}\neq0$,
\begin{equation}\label{eq:chap_1_7_2}
\begin{split}
\frac{s_{1}}{s_{2}}&=\frac{r_{1}(e_{xy}e^{i\theta_{1}}\cos\varphi_{1}+\sin\varphi_{1})}{r_{2}(e_{xy}e^{i\theta_{2}}\cos\varphi_{2}+\sin\varphi_{2})}\\
                   &=\frac{r_{1}}{r_{2}}[e_{xy}(\frac{e^{i\theta_{1}}\cos\varphi_{1}+\sin\varphi_{1}}{e^{i\theta_{2}}\cos\varphi_{2}+\sin\varphi_{2}}
                    -\frac{\sin\varphi_{1}}{\sin\varphi_{2}})+\frac{\sin\varphi_{1}}{\sin\varphi_{2}}] \\
                   &=r(e_{xy}e^{i\theta}\cos\varphi+\sin\varphi)
\end{split}
\end{equation}
where
\[
r=\frac{r_{1}}{r_{2}}, \varphi=\arcsin(\frac{\sin\varphi_{1}}{\sin\varphi_{2}}),
\]
\[
\theta=\arctan[\frac{\sin(\theta_{1}-\theta_{2})+\sin\theta_{1}\tan\varphi_{2}-\sin\theta_{2}\tan\varphi_{1}}
{\cos(\theta_{1}-\theta_{2})+\cos\theta_{1}\tan\varphi_{2}-\cos\theta_{2}\tan\varphi_{1}-\tan\varphi_{1}\cot\varphi_{2}}],
\]
and
\[
[(\cos(\theta_{1}-\theta_{2})+\cos\theta_{1}\tan\varphi_{2}-\cos\theta_{2}\tan\varphi_{1}-\tan\varphi_{1}\cot\varphi_{2})^{2}
\]
\[
+(\sin(\theta_{1}-\theta_{2})+\sin\theta_{1}\tan\varphi_{2}-\sin\theta_{2}\tan\varphi_{1})^{2}]
(\frac{\cos\varphi_{1}\cos\varphi_{2}}{1+\cos\theta_{2}\sin2\varphi_{2}})^{2}
+(\frac{\sin\varphi_{1}}{\sin\varphi_{2}})^{2}=1.
\]

Because $1=1\cdot[e_{xy}e^{i0}\cos(\pi/2)+\sin(\pi/2)]$, it follows from expression~(\ref{eq:chap_1_7_2}) that the inverse of any nonzero spatial complex number $s$ with $\cos\theta\sin2\varphi\neq-1$ is
\begin{equation}\label{eq:chap_1_7_3}
\begin{split}
\frac{1}{s}&=\frac{1}{r}[e_{xy}(\frac{1}{e^{i\theta}\cos\varphi+\sin\varphi}-\frac{1}{\sin\varphi})+\frac{1}{\sin\varphi}]\\
           &=\frac{1}{r}\{e_{xy}[\frac{e^{-i\theta}\cos\varphi}{1+\cos\theta\sin2\varphi}+\kappa_{x}(\theta,\varphi)]+\kappa_{z}(\varphi)\}\\
\end{split}
\end{equation}
where if $\sin\varphi=0$ then $\kappa_{x}(\theta,\varphi)=0$ and $\kappa_{z}(\varphi)=0$ else
\[
\kappa_{x}(\theta,\varphi)=\frac{\sin\varphi}{1+\cos\theta\sin2\varphi}-\frac{1}{\sin\varphi}\texttt{ and }\kappa_{z}(\varphi)=\frac{1}{\sin\varphi}.
\]
Expressions~(\ref{eq:chap_1_7_1}),~(\ref{eq:chap_1_7_2}), and~(\ref{eq:chap_1_7_3}) are, of course, easily remembered by applying the usual algebraic rules for real numbers and $e^{x}$.

Now we consider the product $s_{1}s_{2}\cdots s_{n}$. It has exponential form
\begin{equation}\label{eq:chap_1_7_5}
\begin{split}
s &= \prod_{k=1}^{n}s_{k}=\prod_{k=1}^{n}r_{k}(e_{xy}e^{i\theta_{k}}\cos\varphi_{k}+\sin\varphi_{k}) \\
&=\prod_{k=1}^{n}r_{k}\{e_{xy}[\prod_{k=1}^{n}(e^{i\theta_{k}}\cos\varphi_{k}+\sin\varphi_{k})-\prod_{k=1}^{n}\sin\varphi_{k}]+\prod_{k=1}^{n}\sin\varphi_{k}\}\\
&=r(e_{xy}e^{i\theta}\cos\varphi+\sin\varphi)\texttt{ for }n=2,3,\ldots
\end{split}
\end{equation}
where
\[
r=\prod_{k=1}^{n}r_{k},\texttt{ }\varphi=\arcsin(\prod_{k=1}^{n}\sin\varphi_{k}),
\]
\[
\theta=\arctan\frac{Im\prod_{k=1}^{n}(e^{i\theta_{k}}+\tan\varphi_{k})}{
Re\prod_{k=1}^{n}(e^{i\theta_{k}}+\tan\varphi_{k})-\prod_{k=1}^{n}\tan\varphi_{k}}.
\]
It is easily verified by mathematical induction. We first note that expressions~(\ref{eq:chap_1_7_1}) show it true when $n=2$. Next, we assume that it is valid when $n=m\geq2$. In view of expression~(\ref{eq:chap_1_7_1}) for the product of two nonzero spatial complex numbers in exponential form, it is then valid for $n=m+1$:
\[
\begin{split}
s &= \prod_{k=1}^{m+1}s_{k}=r_{m+1}(e_{xy}e^{i\theta_{m+1}}\cos\varphi_{m+1}+\sin\varphi_{m+1}) \\
  &\times\prod_{k=1}^{m}r_{k}\{e_{xy}[\prod_{k=1}^{m}(e^{i\theta_{k}}\cos\varphi_{k}+\sin\varphi_{k})-\prod_{k=1}^{m}\sin\varphi_{k}]+\prod_{k=1}^{m}\sin\varphi_{k}\}\\
  &=\prod_{k=1}^{m+1}r_{k}\{e_{xy}[\prod_{k=1}^{m+1}(e^{i\theta_{k}}\cos\varphi_{k}+\sin\varphi_{k})-\prod_{k=1}^{m+1}\sin\varphi_{k}]+\prod_{k=1}^{m+1}\sin\varphi_{k}\}.
\end{split}
\]

Also we consider the product $s_{1}^{-1}s_{2}^{-1}\cdots s_{n}^{-1}$. It has exponential form
\begin{equation}\label{eq:chap_1_7_6}
\begin{split}
s&=\prod_{k=1}^{n}s_{k}^{-1}=\prod_{k=1}^{n}\frac{1}{r_{k}}[e_{xy}(\frac{1}{e^{i\theta_{k}}\cos\varphi_{k}+\sin\varphi_{k}}-\frac{1}{\sin\varphi_{k}})+\frac{1}{\sin\varphi_{k}}]\\
&=\prod_{k=1}^{n}\frac{1}{r_{k}}[e_{xy}(\prod_{k=1}^{n}\frac{1}{e^{i\theta_{k}}\cos\varphi_{k}+\sin\varphi_{k}}-\prod_{k=1}^{n}\frac{1}{\sin\varphi_{k}})+\prod_{k=1}^{n}\frac{1}{\sin\varphi_{k}}]\\
&=r(e_{xy}e^{i\theta}\cos\varphi+\sin\varphi)\texttt{ for }n=1,2,\ldots
\end{split}
\end{equation}
where
\[
r=\prod_{k=1}^{n}\frac{1}{r_{k}},\texttt{ }\varphi=\arcsin(\prod_{k=1}^{n}\frac{1}{\sin\varphi_{k}}),
\]
\[
\theta=\arctan\frac{Im\prod_{k=1}^{n}(e^{-i\theta_{k}}\cos\varphi_{k}+\sin\varphi_{k})}{
Re\prod_{k=1}^{n}(e^{-i\theta_{k}}\cos\varphi_{k}+\sin\varphi_{k})-\prod_{k=1}^{n}(1+\cos\theta_{k}\sin2\varphi_{k})/\sin\varphi_{k}}.
\]
It is also easily verified by mathematical induction. To be specific, we first note that it is true when $n=1$. Next, we assume that it is valid when $n=m$, where $m$ is any positive integer. In view of expression~(\ref{eq:chap_1_7_1}) for the product of two nonzero spatial complex numbers in exponential form, it is then valid for $n=m+1$:
\[
\begin{split}
s&=\prod_{k=1}^{m+1}s_{k}^{-1}=\frac{1}{r_{m+1}}[e_{xy}(\frac{1}{e^{i\theta_{m+1}}\cos\varphi_{m+1}+\sin\varphi_{m+1}}-\frac{1}{\sin\varphi_{m+1}})+\frac{1}{\sin\varphi_{m+1}}]\\
&\times\prod_{k=1}^{m}\frac{1}{r_{k}}[e_{xy}(\prod_{k=1}^{m}\frac{1}{e^{i\theta_{k}}\cos\varphi_{k}+\sin\varphi_{k}}-\prod_{k=1}^{m}\frac{1}{\sin\varphi_{k}})+\prod_{k=1}^{m}\frac{1}{\sin\varphi_{k}}]\\
&=\prod_{k=1}^{m+1}\frac{1}{r_{k}}[e_{xy}(\prod_{k=1}^{m+1}\frac{1}{e^{i\theta_{k}}\cos\varphi_{k}+\sin\varphi_{k}}-\prod_{k=1}^{m+1}\frac{1}{\sin\varphi_{k}})+\prod_{k=1}^{m+1}\frac{1}{\sin\varphi_{k}}].
\end{split}
\]

Another important result that can be obtained formally by applying rules for real numbers to $s=r(e_{xy}e^{i\theta}\cos\varphi+\sin\varphi)$ is
\begin{equation}\label{eq:chap_1_7_7}
\begin{split}
s^{n} &=r^{n}(e_{xy}e^{i\theta}\cos\varphi+\sin\varphi)^{n} \\
      &=r^{n}\{e_{xy}[(e^{i\theta}\cos\varphi+\sin\varphi)^{n}-\sin^{n}\varphi]+\mathrm{\mathrm{\mathrm{\sin^{n}\varphi}}}\}
\end{split}
\end{equation}
for $n=0,\pm1,\pm2,\ldots$. It is easily verified for positive values of $n$ by mathematical induction. To be specific, we first note that it becomes $s=r(e_{xy}e^{i\theta}\cos\varphi+\sin\varphi)$ when $n=1$. Next, we assume that it is valid when $n=m$, where $m$ is any positive integer. In view of expression~(\ref{eq:chap_1_7_1}) for the product of two nonzero spatial complex numbers in exponential form, it is then valid for $n=m+1$:
\[
\begin{split}
s^{m+1}&=ss^{m}=r(e_{xy}e^{i\theta}\cos\varphi+\sin\varphi)r^{m}(e_{xy}e^{i\theta}\cos\varphi+\sin\varphi)^{m}\\
&=r^{m+1}(e_{xy}e^{i\theta}\cos\varphi+\sin\varphi)\{e_{xy}[(e^{i\theta}\cos\varphi+\sin\varphi)^{m}-\sin^{m}\varphi]+\sin^{m}\varphi\}\\
&=r^{m+1}\{e_{xy}[(e^{i\theta}\cos\varphi+\sin\varphi)^{m+1}-\sin^{m+1}\varphi]+\sin^{m+1}\varphi\}.
\end{split}
\]
Expression~(\ref{eq:chap_1_7_7}) is thus verified when $n$ is a positive integer. It also holds when $n=0$,
with the convention that $s^{0}=1$. If $n=-1,-2,\ldots$, on the other hand, we define $s^{n}$ in terms of  the multiplicative inverse of $s$ by writing
\[
s^{n}=(s^{-1})^{m}\texttt{ where }m=-n=1,2,\ldots.
\]
Then, since expression~(\ref{eq:chap_1_7_7}) is valid for positive integral powers, it follows from the
exponential form~(\ref{eq:chap_1_7_3}) of $s^{-1}$ that
\[
\begin{split}
s^{n}&=(\frac{1}{r}\frac{1}{e_{xy}e^{i\theta}\cos\varphi+\sin\varphi})^{m} \\
&=(\frac{1}{r})^{m}\{e_{xy}[(\frac{1}{e^{i\theta}\cos\varphi+\sin\varphi})^{m}-(\frac{1}{\sin\varphi})^{m}]+(\frac{1}{\sin\varphi})^{m}\}\\
&=(\frac{1}{r})^{-n}\{e_{xy}[(\frac{1}{e^{i\theta}\cos\varphi+\sin\varphi})^{-n}-(\frac{1}{\sin\varphi})^{-n}]+(\frac{1}{\sin\varphi})^{-n}\}\\
&=r^{n}\{e_{xy}[(e^{i\theta}\cos\varphi+\sin\varphi)^{n}-\sin^{n}\varphi]+\sin^{n}\varphi\}\texttt{ }(n=-1,-2,\ldots).
\end{split}
\]
Expression~(\ref{eq:chap_1_7_7}) is now established for all integral powers.

%Observe that if $r=1$, expression~(\ref{eq:chap_1_7_7}) becomes

From expression~(\ref{eq:chap_1_7_7}), let $s_{n}=s^{n}$ for $n=0,1,2,\ldots$. Then we get
\begin{equation}\label{eq:chap_1_7_8}
s_{n}=r_{n}(e_{xy}e^{i\theta_{n}}\cos\varphi_{n}+\sin\varphi_{n})\texttt{ }(n=0,1,2,\ldots)
\end{equation}
where
\[
r_{n}=r^{n},\texttt{ }\varphi_{n}=\arcsin(\sin^{n}\varphi),\texttt{ }
\theta_{n}=\arctan(\frac{\sum_{k=0}^{n-1}C_{n}^{k}\sin[(n-k)\theta]\tan^{k}\varphi}
{\sum_{k=0}^{n-1}C_{n}^{k}\cos[(n-k)\theta]\tan^{k}\varphi})
\]
because of
\[
(e^{i\theta}\cos\varphi+\sin\varphi)^{n}-\sin^{n}\varphi=\sum_{k=0}^{n}C_{n}^{k}(e^{i\theta}\cos\varphi)^{n-k}\sin^{k}\varphi-\sin^{n}\varphi
\]
\[
=\sum_{k=0}^{n-1}C_{n}^{k}e^{i(n-k)\theta}\cos^{n-k}\varphi\sin^{k}\varphi\texttt{ where }C_{n}^{k}=\frac{n!}{(n-k)!k!}.
\]
Also from expression~(\ref{eq:chap_1_7_7}), let $s_{m}=s^{-m}$ for $m=1,2,\ldots$. Then we get
\begin{equation}\label{eq:chap_1_7_9}
s_{m}=r_{m}(e_{xy}e^{i\theta_{m}}\cos\varphi_{m}+\sin\varphi_{m})\texttt{ }(m=1,2,\ldots)
\end{equation}
where $r_{m}=r^{-m}$, $\varphi_{m}=\arcsin(\sin^{-m}\varphi)$,
\[
\theta_{m}=\arctan(\frac{\sum_{k=0}^{m-1}C_{m}^{k}\sin[(m-k)\theta]\tan^{k}\varphi}
{\sum_{k=0}^{m-1}C_{m}^{k}\{(\cot\varphi+2\cos\theta)^{m-k}-\cos[(m-k)\theta]\}\tan^{k}\varphi})
\]
because of
\[
(\frac{1}{e^{i\theta}\cos\varphi+\sin\varphi})^{m}-(\frac{1}{\sin\varphi})^{m}=\frac{\sum_{k=0}^{m-1}C_{m}^{k}e^{-i(m-k)\theta}\cos^{m-k}\varphi\sin^{k}\varphi}{(1+\cos\theta\sin2\varphi)^{m}}
\]
\[
-\frac{\sum_{k=0}^{m-1}C_{m}^{k}(\cot\varphi+2\cos\theta)^{m-k}\cos^{m-k}\varphi\sin^{k}\varphi}{(1+\cos\theta\sin2\varphi)^{m}}\texttt{ where }C_{m}^{k}=\frac{m!}{(m-k)!k!}.
\]

Expressions~(\ref{eq:chap_1_7_7}-\ref{eq:chap_1_7_9}) can be useful in finding powers of spatial complex numbers even when they are given in rectangular form and the result is desired in that form.

\section{Roots of Spatial Complex Numbers}\label{sec:chap_1_8}

Consider now a point
\[
s=r(e_{xy}e^{i\theta}\cos\varphi+\sin\varphi),
\]
lying on a circle centered at the origin with radius $r$ in a spatial plane. As $\theta$ is increased, $s$ moves around the circle in the counterclockwise direction. In particular, when $\theta$ is increased by $2\pi$, we arrive at the original point, but while $\theta$ is increasing, $\varphi$ varies back and forth between $-\varphi_{max}$ and $\varphi_{max}$ where $0\leq\varphi_{max}\leq\pi/2$; and the same is true when $\theta$ is decreased by $2\pi$. It is, therefore, evident that \textit{two nonzero complex numbers}
\[
s_{1}=r_{1}(e_{xy}e^{i\theta_{1}}\cos\varphi_{1}+\sin\varphi_{1})
\]
\textit{and}
\[
s_{2}=r_{2}(e_{xy}e^{i\theta_{2}}\cos\varphi_{2}+\sin\varphi_{2})
\]
\textit{are equal if and only if}
\[
r_{1}=r_{2},\texttt{ }\varphi_{1}=\varphi_{2},\texttt{ \textit{and} }\theta_{1}=\theta_{2}+2k\pi,
\]
\textit{where $k$ is some integer $(k=0,\pm1,\pm2,\ldots)$}.

This observation, together with the expression
\[
s^{n}=r^{n}(e_{xy}e^{i\theta}\cos\varphi+\sin\varphi)^{n}
\]
in Sec.~(\ref{sec:chap_1_7}) for integral powers of spatial complex numbers $s$, is useful in finding the $n$th roots of any nonzero spatial complex number
\[
s_{0}=r_{0}(e_{xy}e^{i\theta_{0}}\cos\varphi_{0}+\sin\varphi_{0}),
\]
where $n$ has one of the values $n=2,3,\ldots$. The method starts with the fact that an $n$th root of $s_{0}$ is a nonzero number $s$ such that $s^{n}=s_{0}$, or
\[
r^{n}(e_{xy}e^{i\theta}\cos\varphi+\sin\varphi)^{n}=r_{0}(e_{xy}e^{i\theta_{0}}\cos\varphi_{0}+\sin\varphi_{0}).
\]

From expression~(\ref{eq:chap_1_4_2}) in Sec.~(\ref{sec:chap_1_4}) we get
\[
|e_{xy}e^{i\theta}\cos\varphi+\sin\varphi|^{n}=1\texttt{ }(n=2,3,\ldots)
\]
because of
\[
\begin{split}
|e_{xy}e^{i\theta}\cos\varphi+\sin\varphi|^{2} &= |e^{i\theta}\cos\varphi|^{2}+|\sin\varphi|^{2} \\
                                               &= \cos^{2}\varphi+\sin^{2}\varphi=1.
\end{split}
\]

According to the statement above, then, we define the roots of spatial complex numbers as follows:

\begin{definition}\label{def:chap_1_8_1}
Given a positive integer $n$, roots of spatial complex numbers $s_{0}$ with
\[
z_{0}\geq0,\texttt{ }r_{0}=|s_{0}|>0,\texttt{ }\theta_{0}=\arctan(y_{0}/x_{0}),\texttt{ and }\varphi_{0}=\arcsin(z_{0}/r_{0})
\]
are defined as
\begin{equation}\label{eq:chap_1_8_0}
\begin{split}
s &= \sqrt[n]{s_{0}} \\
  &= e_{xy}(\sqrt[n]{x_{0}+iy_{0}+z_{0}}-\sqrt[n]{z_{0}})+\sqrt[n]{z_{0}} \\
  &= r(e_{xy}e^{i\theta}\cos\varphi+\sin\varphi)\texttt{ }(n=2,3,\ldots)
\end{split}
\end{equation}
where
\[
r=\sqrt[n]{r_{0}},\texttt{ }\theta=\theta_{b}(k)
\]
\[
\sin\varphi=\sqrt[n]{\sin\varphi_{0}},\texttt{ }\cos\varphi=\sqrt[n]{r_{a}}r_{b}(k),
\]
\[
r_{a}^{2}=1+\cos\theta_{0}\sin2\varphi_{0},\texttt{ }\theta_{a}=\arctan\frac{\sin\theta_{0}}{\cos\theta_{0}+\tan\varphi_{0}},
\]
\[
r_{b}^{2}(k)=1-2\cos\frac{\theta_{a}+2k\pi}{n}\sqrt[n]{\frac{\sin\varphi_{0}}{r_{a}}}+(\sqrt[n]{\frac{\sin\varphi_{0}}{r_{a}}})^{2},
\]
and
\[
\theta_{b}(k)=\arctan\frac{\sin(\theta_{a}/n+2k\pi/n)}{\cos(\theta_{a}/n+2k\pi/n)-\sqrt[n]{(\sin\varphi_{0})/r_{a}}}\texttt{ }(k=0,\pm1,\pm2,\ldots).
\]
\end{definition}

First, from the following equation
\[
\begin{split}
x_{0}+iy_{0}+z_{0} &=r_{0}r_{a}(\frac{\cos\theta_{0}\cos\varphi_{0}+\sin\varphi_{0}}{r_{a}}+i\frac{\sin\theta_{0}\cos\varphi_{0}}{r_{a}})\\
&=r_{0}r_{a}\exp(i\theta_{a}),
\end{split}
\]
we get
\[
\begin{split}
r_{a}^{2} &= (\cos\theta_{0}\cos\varphi_{0}+\sin\varphi_{0})^{2}+(\sin\theta_{0}\cos\varphi_{0})^{2} \\
          &= 1+\cos\theta_{0}\sin2\varphi_{0}
\end{split}
\]
and
\[
\begin{split}
\theta_{a} &= \arctan\frac{\sin\theta_{0}\cos\varphi_{0}}{\cos\theta_{0}\cos\varphi_{0}+\sin\varphi_{0}} \\
           &= \arctan\frac{\sin\theta_{0}}{\cos\theta_{0}+\tan\varphi_{0}}.
\end{split}
\]

Next, from the following equation
\[
\begin{split}
\sqrt[n]{x_{0}+iy_{0}+z_{0}}-\sqrt[n]{z_{0}}
&=r\sqrt[n]{r_{a}}[\exp(i\frac{\theta_{a}+2k\pi}{n})-\sqrt[n]{\frac{\sin\varphi_{0}}{r_{a}}}] \\
%=r\sqrt[n]{r_{a}}r_{b}(k)[\frac{\cos(\theta_{a}/n+2k\pi/n)-\sqrt[n]{(\sin\varphi_{0})/r_{a}}}{r_{b}(k)}+i\frac{\sin(\theta_{a}/n+2k\pi/n)}{r_{b}(k)}]
&=r\sqrt[n]{r_{a}}r_{b}(k)\exp[i\theta_{b}(k)]\texttt{ }(k=0,\pm1,\pm2,\ldots),
\end{split}
\]
we get
\[
r_{b}^{2}(k)=1-2\cos\frac{\theta_{a}+2k\pi}{n}\sqrt[n]{\frac{\sin\varphi_{0}}{r_{a}}}+(\sqrt[n]{\frac{\sin\varphi_{0}}{r_{a}}})^{2}
\]
and
\[
\theta_{b}(k)=\arctan\frac{\sin(\theta_{a}/n+2k\pi/n)}{\cos(\theta_{a}/n+2k\pi/n)-\sqrt[n]{(\sin\varphi_{0})/r_{a}}}\texttt{ }(k=0,\pm1,\pm2,\ldots).
\]

Finally, from the following equation
\[
\begin{split}
|s|^{2} &= |\sqrt[n]{x_{0}+iy_{0}+z_{0}}-\sqrt[n]{z_{0}}|^{2}+(\sqrt[n]{z_{0}})^{2} \\
        &= r^{2}[(\sqrt[n]{r_{a}}r_{b}(k))^{2}+(\sqrt[n]{\sin\varphi_{0}})^{2}],
\end{split}
\]
we get
\[
(\sqrt[n]{r_{a}}r_{b}(k))^{2}+(\sqrt[n]{\sin\varphi_{0}})^{2}=1\texttt{ and }\sqrt[n]{r_{a}}r_{b}(k)=\cos\varphi.
\]

Thus equation~(\ref{eq:chap_1_8_0}) holds. So $r=\sqrt[n]{r_{0}}$, where this radical denotes
the unique positive $n$th root of the positive real number $r_{0}$, $\varphi=\arcsin(\sqrt[n]{\sin\varphi_{0}})$, and
\[
\begin{split}
\theta &= \theta_{b}(k) \\
       &= \arctan\frac{\sin(\theta_{a}/n+2k\pi/n)}{\cos(\theta_{a}/n+2k\pi/n)-\sqrt[n]{(\sin\varphi_{0})/r_{a}}}\texttt{ }(k=0,\pm1,\pm2,\ldots)
\end{split}
\]
where
\[
\theta_{a}=\arctan\frac{\sin\theta_{0}}{\cos\theta_{0}+\tan\varphi_{0}}.
\]

Consequently, the spatial complex numbers $s$ expressed by equation~(\ref{eq:chap_1_8_0}) are the $n$th roots  of $s_{0}$. We are able to see immediately from this exponential form of the roots that they all lie on the circle $|s|=\sqrt[n]{r_{0}}$ about the origin and are unequally spaced radians, starting with argument $\theta(0)$ where
\begin{equation}\label{eq:chap_1_8_00}
\theta(k)=\arctan\frac{\sin(\theta_{a}/n+2k\pi/n)}{\cos(\theta_{a}/n+2k\pi/n)-\sqrt[n]{(\sin\varphi_{0})/r_{a}}}\texttt{ }(k=0,1,2,\ldots,n-1).
\end{equation}
Evidently, then, all of the distinct roots are obtained when $k=0,1,2,\ldots,n-1$, and no further roots arise with other values of $k$. We let $c_{k}$ $(k=0,1,2,\ldots,n-1)$ denote these distinct roots and write
\begin{equation}\label{eq:chap_1_8_1}
c_{k}=\sqrt[n]{r_{0}}(e_{xy}e^{i\theta(k)}\cos\varphi+\sin\varphi)\texttt{ }(k=0,1,2,\ldots,n-1).
\end{equation}

The number $\sqrt[n]{r_{0}}$ is the length of each of the radius vectors representing the $n$ roots. The first root $c_{0}$ has argument $\theta(0)$; and the $n$ roots lie at the vertices of an irregular polygon of $n$ sides inscribed in that circle.

We shall let $s_{0}^{1/n}$ denote the set of $n$th roots of $s_{0}$. If, in particular, $s_{0}$ is a positive real number $r_{0}$, the symbol $r_{0}^{1/n}$ denotes the entire set of roots; and the symbol $\sqrt[n]{r_{0}}$ in expression~(\ref{eq:chap_1_8_1}) is reserved for the one positive root. When the value of $\theta_{0}$ that is used in expression~(\ref{eq:chap_1_8_1}) is the principal value of $\arg_{c}s_{0}$ $(-\pi<\theta_{0}\leq\pi)$, the number $c_{0}$ is referred to as the principal root. Thus when $s_{0}$ is a positive real number $r_{0}$, its principal root is $\sqrt[n]{r_{0}}$.

Finally, a convenient way to remember expression~(\ref{eq:chap_1_8_1}) is to write $s_{0}$ in its most
general exponential form
\begin{equation}\label{eq:chap_1_8_2}
s_{0}=r_{0}(e_{xy}e^{i(\theta_{0}+2k\pi)}\cos\varphi_{0}+\sin\varphi_{0})\texttt{ }(k=0,\pm1,\pm2,\ldots).
\end{equation}
and to formally apply laws of fractional exponents involving real numbers, keeping in mind that there are precisely $n$ roots:
\[
\begin{split}
s_{0}^{1/n} &= [r_{0}(e_{xy}e^{i(\theta_{0}+2k\pi)}\cos\varphi_{0}+\sin\varphi_{0})]^{1/n} \\
            &= \sqrt[n]{r_{0}}(e_{xy}e^{i\theta(k)}\cos\varphi+\sin\varphi)
\end{split}
\]
where $\theta(k)$ is determined by expression~(\ref{eq:chap_1_8_00}) for $k=0,1,2,n-1$.

\section{Spatial Planes in the Complex Space}\label{sec:chap_1_9}

\subsection{Spatial planes determined by an argument $\varphi$}\label{sec:chap_1_9_1}

The expression of a spatial complex number $s$ in a three-dimensional spherical coordinates system is
\[
s=r(e_{xy}e^{i\theta}\cos\varphi+\sin\varphi)\texttt{ }(r=|s|,\texttt{ }0\leq\theta\leq2\pi,\texttt{ }-\pi/2\leq\varphi\leq\pi/2).
\]
Given an argument $\varphi$, a spatial plane is determined by the fixed argument $\varphi$ and parallel to the $xy$ coordinate plane.

Let $S_{\rho}$ denote the sphere $|s|=\rho$, and let $C_{\varphi}$ denote a positively oriented circle on the sphere $S_{\rho}$, which is in the spatial plane determined by the argument $\varphi$ and parallel to the $xy$ coordinate plane. Then there is the parametric representation
\begin{equation}\label{eq:chap_1_9_1}
s=\rho(e_{xy}e^{i\theta}\cos\varphi+\sin\varphi)\texttt{ for }C_{\varphi}\texttt{ }(0\leq\theta\leq2\pi, \varphi\texttt{ is fixed}).
\end{equation}

\subsection{Spatial planes determined by an argument $\theta$}\label{sec:chap_1_9_2}

The expression of a spatial complex number $s$ in a three-dimensional spherical coordinates system is
\[
s=r(e_{xy}e^{i\theta}\cos\varphi+\sin\varphi)\texttt{ }(r=|s|,\texttt{ }0\leq\theta\leq2\pi,\texttt{ }-\pi/2\leq\varphi\leq\pi/2).
\]
Given an argument $\theta$, a spatial plane is determined by the fixed argument $\theta$ and perpendicular to the $xy$ coordinate plane.

Let $S_{\rho}$ denote the sphere $|s|=\rho$, and let $C_{\theta}$ denote a positively oriented circle on the sphere $S_{\rho}$, which is in the spatial plane determined by the argument $\theta$. Then there is the parametric representation
\begin{equation}\label{eq:chap_1_9_2}
s=\rho(e_{xy}e^{i\theta}\cos\varphi+\sin\varphi)\texttt{ for }C_{\theta}\texttt{ }(\theta\texttt{ is fixed},\texttt{ }-\pi/2\leq\varphi\leq\pi/2).
\end{equation}

\subsection{Spatial planes determined by an argument $\varphi$ and the $x$ axis}\label{sec:chap_1_9_3}

The expression of a spatial complex number $s$ in three-dimensional rectangular and spherical coordinates systems is
\[
\begin{split}
s &= e_{xy}(x+iy)+z \\
  &= r(e_{xy}e^{i\theta}\cos\varphi+\sin\varphi)\texttt{ }(r=|s|,\texttt{ }0\leq\theta\leq2\pi,\texttt{ }-\pi/2\leq\varphi\leq\pi/2).
\end{split}
\]
Given an argument $\varphi$, a spatial plane is determined by the argument $\varphi$ and the $x$ axis, which is the result that the $xy$ coordinate plane rotates around the $x$ axis for an angle $\varphi$.

Let the origin be unchanged and $x^{*}$, $y^{*}$, and $z^{*}$ denote three axes of the transformed three-dimensional rectangular coordinates system. Then the $x$ axis becomes the $x^{*}$ axis, points $s$ become $s^{*}$ in the transformed coordinates system, and there are
\[
\begin{split}
s^{*} &= e_{xy}(x^{*}+iy^{*}) \\
      &= e_{xy}r^{*}e^{i\theta^{*}}\texttt{ }(r^{*}=|s^{*}|,\texttt{ }0\leq\theta^{*}\leq2\pi,\texttt{ }\varphi^{*}=0).
\end{split}
\]

Since the origin is unchanged, there are
\[
r^{*}=|s^{*}|=|s|=r,\texttt{ }x^{*}=r\cos\theta^{*},\texttt{ }y^{*}=r\sin\theta^{*},\texttt{ }z^{*}=0.
\]

Let $S_{\rho}$ denote the sphere $|s|=\rho$, and let $C_{\varphi}$ denote a positively oriented circle on the sphere $S_{\rho}$, which is in the spatial plane determined by the argument $\varphi$ and the $x$ axis. Then there is the parametric representation
\begin{equation}\label{eq:chap_1_9_3}
s^{*}=\rho e^{i\theta^{*}}\texttt{ for }C_{\varphi}\texttt{ }(0\leq\theta^{*}\leq2\pi,\varphi^{*}=0).
\end{equation}

\subsection{Spatial planes determined by an argument $\theta$ and the $z$ axis}\label{sec:chap_1_9_4}

The expression of a spatial complex number $s$ in three-dimensional rectangular and spherical coordinates systems is
\[
\begin{split}
s &= e_{xy}(x+iy)+z \\
  &= r(e_{xy}e^{i\theta}\cos\varphi+\sin\varphi)\texttt{ }(r=|s|,\texttt{ }0\leq\theta\leq2\pi,\texttt{ }-\pi/2\leq\varphi\leq\pi/2).
\end{split}
\]
Given an argument $\theta$, a spatial plane is determined by the argument $\theta$ and the $z$ axis, which is the result that the $xy$ coordinate plane rotates around the $z$ axis for an angle $\theta$ and derives the new $x^{*}$ axis, then rotates around the $x^{*}$ axis for an angle $\pi/2$.

Let the origin be unchanged and $x^{*}$, $y^{*}$, and $z^{*}$ denote three axes of the transformed three-dimensional rectangular coordinates system. Then points $s$ become $s^{*}$ in the transformed coordinates system and there are
\[
\begin{split}
s^{*} &= (x^{*},y^{*},0) \\
      &= e_{xy}(x^{*}+iy^{*})=e_{xy}r^{*}e^{i\theta^{*}}
\end{split}
\]
where $r^{*}=|s^{*}|$ and $\varphi^{*}=0$.

Since the origin is unchanged, there are
\[
r^{*}=|s^{*}|=|s|=r,\texttt{ }x^{*}=r\cos\theta^{*},\texttt{ }y^{*}=r\sin\theta^{*},\texttt{ }z^{*}=0.
\]

Let $S_{\rho}$ denote the sphere $|s|=\rho$, and let $C_{\rho}$ denote a positively oriented circle on the sphere $S_{\rho}$, which is in the spatial plane determined by the argument $\theta$ and the $z$ axis. Then there is the parametric representation
\begin{equation}\label{eq:chap_1_9_4}
s^{*}=\rho e^{i\theta^{*}}\texttt{ for }C_{\rho}\texttt{ }(0\leq\theta^{*}\leq2\pi,\varphi^{*}=0).
\end{equation}

\subsection{Spatial planes determined by points $s_{0}$, $s_{1}$, and the origin}\label{sec:chap_1_9_5}

Given three points $s_{0}$, $s_{1}$, and the origin in three-dimensional rectangular and spherical coordinates systems where
\[
\begin{split}
s_{0} &= (x_{0},y_{0},z_{0})=e_{xy}(x_{0}+iy_{0})+z_{0} \\
      &= r_{0}(e_{xy}e^{i\theta_{0}}\cos\varphi_{0}+\sin\varphi_{0}),\texttt{ }r_{0}=|s_{0}|
\end{split}
\]
and
\[
\begin{split}
s_{1} &= (x_{1},y_{1},z_{1})=e_{xy}(x_{1}+iy_{1})+z_{1} \\
      &= r_{1}(e_{xy}e^{i\theta_{1}}\cos\varphi_{1}+\sin\varphi_{1}),\texttt{ }r_{1}=|s_{1}|,
\end{split}
\]
a spatial plane is determined by the three points which are in the plane. Let the origin be unchanged and $x^{*}$, $y^{*}$, and $z^{*}$ denote three axes of the transformed three-dimensional rectangular coordinates system. Then two points $s_{0}$ and $s_{1}$ become $s_{0}^{*}$ and $s_{1}^{*}$ in the transformed coordinates system and there are
\[
\begin{split}
s_{0}^{*} &= (x_{0}^{*},y_{0}^{*},0) \\
          &= e_{xy}(x_{0}^{*}+iy_{0}^{*})=e_{xy}r_{0}^{*}e^{i\theta_{0}^{*}}\texttt{ where }r_{0}^{*}=|s_{0}^{*}|\texttt{ and }\varphi_{0}^{*}=0,
\end{split}
\]
and
\[
\begin{split}
s_{1}^{*} &= (x_{1}^{*},y_{1}^{*},0) \\
          &= e_{xy}(x_{1}^{*}+iy_{1}^{*})=e_{xy}r_{1}^{*}e^{i\theta_{1}^{*}}\texttt{ where }r_{1}^{*}=|s_{1}^{*}|\texttt{ and }\varphi_{1}^{*}=0.
\end{split}
\]

Since the origin is unchanged, there are
\[
r_{0}^{*}=|s_{0}^{*}|=|s_{0}|=r_{0}
\]
and
\[
r_{1}^{*}=|s_{1}^{*}|=|s_{1}|=r_{1}.
\]
Thus we get
\[
|s_{1}^{*}-s_{0}^{*}|=|s_{1}-s_{0}|,
\]
\[
(x_{0}^{*})^{2}+(y_{0}^{*})^{2}=x_{0}^{2}+y_{0}^{2}+z_{0}^{2},
\]
and
\[
(x_{1}^{*})^{2}+(y_{1}^{*})^{2}=x_{1}^{2}+y_{1}^{2}+z_{1}^{2}.
\]

Let $\rho$ be the distance between $s_{1}$ and $s_{0}$, and let $S_{\rho}$ denote the sphere $|s-s_{0}|=\rho$ centered at $s_{0}$, and let $C_{\rho}$ denote a positively oriented circle on the sphere $S_{\rho}$, which passes through $s_{1}$ and is in the spatial plane determined by the three points $s_{0}$, $s_{1}$, and the origin. Then there is the parametric representation
\begin{equation}\label{eq:chap_1_9_6}
s^{*}=s_{0}^{*}+\rho e^{i\theta^{*}}\texttt{ for }C_{\rho}\texttt{ }(0\leq\theta^{*}\leq2\pi,\varphi^{*}=0).
\end{equation}

\section{Regions in the Complex Space}\label{sec:chap_1_10}

In this section, we are concern with sets of spatial complex numbers, or points in the $s$ space, and their closeness to one another. Our basic tool is the concept of an $\varepsilon$ \texttt{neighborhood}
\begin{equation}\label{eq:chap_1_10_1}
|s-s_{0}|<\varepsilon
\end{equation}
of a given point $s_{0}$. It consists of all points $s$ lying inside but not on a sphere centered at $s_{0}$ and with a specified positive radius $\varepsilon$. When the value of $\varepsilon$ is understood or immaterial in discussion, the set~(\ref{eq:chap_1_10_1}) is often referred to as just a neighborhood. Occasionally, it is convenient to speak of a \texttt{deleted neighborhood}, or spherical shell,
\begin{equation}\label{eq:chap_1_10_2}
0<|s-s_{0}|<\varepsilon
\end{equation}
consisting all points $s$ an $\varepsilon$ neighborhood of $s_{0}$ except for the point $s_{0}$ itself.

A point $s_{0}$ is said to be an interior point of a set $S$ whenever there is some neighborhood of $s_{0}$ that contains only points of $S$; it is called an exterior point of $S$ when there exists a neighborhood of it containing no points of $S$. If $s_{0}$ is neither of these, it is a boundary point of $S$. A boundary point is, therefore, a point all of whose neighborhoods contain points in $S$ and points not in $S$. The totality of all boundary points is called the boundary of $S$. The sphere $|s|=1$, for instance, is the boundary of each of the sets
\begin{equation}\label{eq:chap_1_10_3}
|s|<1\texttt{ and }|s|\leq1.
\end{equation}

A set is open if it contains none of its boundary points. if and only if each of its points is an interior point. A set is closed if it contains all of its boundary points; and the closure of a set $S$ is the closed set consisting of all points in $S$ together with the boundary of $S$. Note that the first of the sets~(\ref{eq:chap_1_10_3}) is open and that the second is its closure.

Some sets are, of course, neither open nor closed. For a set to be not open, there must be a boundary point that is contained in the set; and if a set is not closed, there exists a boundary point not contained in the set. Observe that the deleted sphere $0<|s|\leq1$ is neither open nor closed. The set of all spatial complex numbers is, on the other hand, both open and closed since it has no boundary points.

An open set $S$ is connected if each pair of points $s_{1}$ and $s_{2}$ in it can be joined by a polygonal line, consisting of a finite number of line segments joined end to end, that lies entirely in $S$. The open set $|s|<1$ is connected. The annulus $1<|s|<2$ is, of course, open and it is also connected. An open set that is connected is called a domain. Note that any neighborhood is a domain. A domain together with some, none, or all of its boundary points is referred to as a region.

A set $S$ is bounded if every point of $S$ lies inside some sphere $|s|=R$; otherwise, it is unbounded. Both  of the sets~(\ref{eq:chap_1_10_3}) are bounded regions, and the half space $Re$ $s\geq0$ is unbounded.

A point $s_{0}$ is said to be an accumulation point of a set $S$ if each deleted neighborhood of $s_{0}$ contains at least one point of $S$. It follows that if a set $S$ is closed, then it contains each of its accumulation points. For if an accumulation point $s_{0}$ were not in $S$, it would be a boundary point of $S$; but this contradicts the fact that a closed set contains all of its boundary points.
%It is left as an exercise to show that the converse is, in fact, true.
Thus, a set is closed if and only if it contains all of its accumulation points.

Evidently, a point $s_{0}$ is not an accumulation point of a set $S$ whenever there exists
some deleted neighborhood of $s_{0}$ that does not contain points of $S$. Note that the origin
is the only accumulation point of the set $s_{n}=i/n$ $(n=1,2,\ldots)$.

%-----------------------------------------------------------------------
% Beginning of chap2.tex
%-----------------------------------------------------------------------
%
% AMS-LaTeX 1.2 sample file for a monograph, based on amsbook.cls.
% This is a data file input by chapter.tex.
%%%%%%%%%%%%%%%%%%%%%%%%%%%%%%%%%%%%%%%%%%%%%%%%%%%%%%%%%%%%%%%%%%%%

%\part{This is a Part Title Sample}

\chapter{Analytic Functions}\label{ch:chap_2}

We now consider functions of a spatial complex variable and develop a theory of differentiation for them. The main goal of the chapter is to introduce spatial analytic functions, which play a central role in spatial complex analysis.

\section{Functions of a Spatial Complex Variable}\label{sec:chap_2_1_11}

Let $S$ be a set of spatial complex numbers. A function $f$ defined on $S$ is a rule that assigns to each $s$ in $S$ a spatial complex number $\varpi$. The number $\varpi$ is called the value of $f$ at $s$ and is denoted by $f(s)$; that is, $\varpi=f(s)$. The set $S$ is called the domain of definition of $f$.

It must be emphasized that both a domain of definition and a rule are needed in order for a function to be well defined. When the domain of definition is not mentioned, we agree that the largest possible set is to be taken. Also, it is not always convenient to use notation that distinguishes between a given function and its values.

\begin{example}\label{ex:chap_2_1_1}
If $f$ is defined on the set $s\neq0$ by means of the equation $\varpi=1/s$, it may be referred to only as the function $\varpi=1/s$, or simply the function $1/s$.

Suppose that $\varpi=e_{xy}(u+iv)+w$ is the value of a function $f$ at $s=e_{xy}(x+iy)+z$, so that
\[
e_{xy}(u+iv)+w=f[e_{xy}(x+iy)+z].
\]
Each of the real numbers $u$, $v$, and $w$ depends on the real variables $x$, $y$, and $z$, and it follows that $f(s)$ can be expressed in terms of a ternary group of real-valued functions of the real variables $x$, $y$, and $z$:
\begin{equation}\label{eq:chap_2_1_1}
f(s)=e_{xy}[u(x,y,z)+iv(x,y,z)]+w(z).
\end{equation}
If the spherical coordinates $r$, $\theta$, and $\varphi$, instead of $x$, $y$ and $z$ are used, then
\[
e_{xy}(u+iv)+w=f[r(e_{xy}e^{i\theta}\cos\varphi+\sin\varphi)],
\]
where $\varpi=e_{xy}(u+iv)+w$ and $s=r(e_{xy}e^{i\theta}\cos\varphi+\sin\varphi)$. In that case, we may write
\begin{equation}\label{eq:chap_2_1_2}
f(s)=e_{xy}[u(r,\theta,\varphi)+iv(r,\theta,\varphi)]+w(r,\varphi).
\end{equation}
\end{example}

\begin{example}\label{ex:chap_2_1_2}
If $f(s)=s^{2}$, then
\[
\begin{split}
f[e_{xy}(x+iy)+z] &= [e_{xy}(x+iy)+z]^{2} \\
                  &= e_{xy}[(x^{2}-y^{2}+2xz)+i2(x+z)y]+z^{2}.
\end{split}
\]
Hence
\[
u(x,y,z)=x^{2}-y^{2}+2xz,\texttt{ }v(x,y,z)=2(x+z)y,\texttt{ and }w(z)=z^{2}.
\]
When spherical coordinates are used,
\[
\begin{split}
f[r(e_{xy}e^{i\theta}\cos\varphi+\sin\varphi)] &= r^{2}(e_{xy}e^{i\theta}\cos\varphi+\sin\varphi)^{2} \\
&= r^{2}[e_{xy}(e^{i2\theta}\cos^{2}\varphi+e^{i\theta}\sin2\varphi)+\sin^{2}\varphi].
\end{split}
\]
Consequently,
\[
u(r,\theta,\varphi)=r^{2}(\cos2\theta\cos^{2}\varphi+\cos\theta\sin2\varphi),
\]
\[
v(r,\theta,\varphi)=r^{2}(\sin2\theta\cos^{2}\varphi+\sin\theta\sin2\varphi),
\]
and
\[
w(r,\varphi)=r^{2}\sin^{2}\varphi.
\]
\end{example}

If, in either of equations~(\ref{eq:chap_2_1_1}) and~(\ref{eq:chap_2_1_2}), the functions $u$ and $v$ always have value zero, then the value of $f$ is always real. That is, $f$ is a real-valued function of a spatial complex variable. Also if, in either of equations~(\ref{eq:chap_2_1_1}) and~(\ref{eq:chap_2_1_2}), the function $w$ always has value zero, then the value of $f$ is always a two-dimensional complex number. That is, $f$ is a two-dimensional complex-valued function of a spatial complex variable.

If $n$ is zero or a positive integer and if $a_{0}$, $a_{1}$, $a_{2}$, $\ldots$, $a_{n}$ are spatial complex constants, where $a_{n}\neq0$, the function
\[
P(s)=a_{0}+a_{1}s+a_{2}s^{2}+\ldots+a_{n}s^{n}
\]
is a polynomial of degree $n$. Note that the sum here has a finite number of terms and that the domain of definition is the entire $s$ space. Quotients $P(s)/Q(s)$ of polynomials are called rational functions and are defined at each point $s$ where $Q(s)\neq0$. Polynomials and rational functions constitute elementary, but important, classes of functions of a spatial complex variable.

A generalization of the concept of function is a rule that assigns more than one value to a point $S$ in the domain  of definition. These multiple-valued functions occur in the theory of functions of a spatial complex variable, just as they do in the case of real or two-dimensional complex variables. When multiple-valued functions are studied, usually just one of the possible values assigned to each point is taken, in a systematic manner, and a (single-valued) function is constructed from the multiple-valued function.

\begin{example}\label{ex:chap_2_1_4}
Let $s$ denote any nonzero spatial complex number. We know from Sec.~(\ref{sec:chap_1_8})
that $s^{1/2}$ has the two values
\[
s^{1/2}=\sqrt{r}(e_{xy}e^{i\theta_{b}(k)}\sqrt{r_{a}}r_{b}(k)+\sqrt{\sin\varphi})\texttt{ }(k=0,1)
\]
where $r=|s|$ and $\Theta$ $(-\pi<\Theta\leq\pi)$ is the principal value of $\arg_{c}s$
\[
r_{a}^{2}=1+\cos\Theta\sin2\varphi,\texttt{ }\theta_{a}=\arctan\frac{\sin\Theta}{\cos\Theta+\tan\varphi},
\]
\[
r_{b}^{2}(k)=1-2\cos(\theta_{a}/2+k\pi)\sqrt{\frac{\sin\varphi}{r_{a}}}+(\sqrt{\frac{\sin\varphi}{r_{a}}})^{2},
\]
and
\[
\theta_{b}(k)=\arctan\frac{\sin(\theta_{a}/2+k\pi)}{\cos(\theta_{a}/2+k\pi)-\sqrt{(\sin\varphi)/r_{a}}}\texttt{ }(k=0,1).
\]
But, if we choose only the positive value of $\pm\sqrt{r}$ when $k=0$ and write
\begin{equation}\label{eq:chap_2_1_3}
f(s)=s^{1/2}=\sqrt{r}(e_{xy}e^{i\theta_{b}(0)}\sqrt{r_{a}}r_{b}(0)+\sqrt{\sin\varphi})\texttt{ }(r>0,\texttt{ }-\pi<\Theta\leq\pi).
\end{equation}
the (single-valued) function~(\ref{eq:chap_2_1_3}) is well defined on the set of nonzero numbers in the $s$ space.  Since zero is the only square root of zero, we also write $f(0)=0$. The function $f$ is then well defined on the entire space.
\end{example}

\section{Mappings}\label{sec:chap_2_2_12}

Properties of a real-valued function of a real variable are often exhibited by the graph of the function. But when $\varpi=f(s)$, where $s$ and $\varpi$ are spatial complex, no such convenient graphical representation of the function $f$ is available because each of the numbers $s$ and $\varpi$ is located in a space rather than on a line. One can, however, display some information about the function by indicating pairs of corresponding points $s=(x,y,z)$ and $\varpi=(u,v,w)$. To do this, it is generally simpler to draw the $s$ and $\varpi$ spaces separately.

When a function $f$ is thought of in this way, it is often referred to as a mapping, or transformation. The image  of a point $s$ in the domain of definition $S$ is the point $\varpi=f(s)$, and the set of images of all points in a set $T$ that is contained in $S$ is called the image of $T$. The image of the entire domain of definition $S$ is called the range of $f$. The inverse image of a point $\varpi$ is the set of all points $s$ in the domain of definition of $f$ that have $\varpi$ as their image. The inverse image of a point may contain just one point, many points, or none at all. The last case occurs, of course, when $s$ is not in the range of $f$.

Terms such as translation, rotation, and reflection are used to convey dominant geometric characteristics of certain mappings. In such cases, it is sometimes convenient to consider the $s$ and $\varpi$ spaces to be the same. %For example, the mapping

More information is usually exhibited by sketching images of curves and regions than by simply indicating images  of individual points. In the following examples, we illustrate this with the transformation $\varpi=s^{2}$.

We begin by finding the images of some curves in the $s$ space.

\begin{example}\label{ex:chap_2_2_1}
According to Example~(\ref{ex:chap_2_1_2}) in Sec.~(\ref{sec:chap_2_1_11}), the mapping $\varpi=s^{2}$ can be thought of as the transformation
\begin{equation}\label{eq:chap_2_2_1}
u(x,y,z)=(x+z)^{2}-y^{2}-z^{2},\texttt{ }v(x,y,z)=2(x+z)y,\texttt{ and }w(z)=z^{2}
\end{equation}
from the $xyz$ space to the $uvw$ space. This form of the mapping is especially useful in finding the images of certain hyperbolas.

First, the third of equations~(\ref{eq:chap_2_2_1}) that $w=z^{2}$ maps the $z$ plane parallel to the $xy$ coordinate plane with a distance $z$ in $s$ space onto the $w$ plane parallel to the $uv$ coordinate plane with a distance $w$ in $\varpi$ space. We note that $z^{2}=c_{w}$. Then $w=c_{w}$ implies a $w=c_{w}$ plane in $\varpi$ space.

Second, it is easy to show, for instance, that each branch of a spatial hyperbola
\begin{equation}\label{eq:chap_2_2_2}
(x+z)^{2}-y^{2}-z^{2}=c_{u}\texttt{ }(c_{u}>0)
\end{equation}
is mapped in a one to one manner onto the vertical $u=c_{u}$ plane parallel to the $vw$ coordinate plane with a distance $u=c_{u}$ in $\varpi$ space. We start by noting from the first of equations~(\ref{eq:chap_2_2_1}) that $u=c_{u}$ when $(x,y,z)$ is a point lying on either branch. When, in particular, it lies on the right-hand branch, the second of equations~(\ref{eq:chap_2_2_1}) tells us that $v=2y\sqrt{y^{2}+c_{u}+c_{w}}$. Thus the image of the right-hand branch can be expressed parametrically as
\[
u=c_{u},\texttt{ }v=2y\sqrt{y^{2}+c_{u}+c_{w}},\texttt{ }w=c_{w}\texttt{ }(-\infty<y<\infty);
\]
and it is evident that the image of a point $(x, y, z)$ on that branch moves upward along the entire spatial curve as $(x, y, z)$ traces out the branch in the upward direction. Likewise, since the component equations of $\varpi=s^{2}$
\[
u=c_{u},\texttt{ }v=-2y\sqrt{y^{2}+c_{u}+c_{w}},\texttt{ }w=c_{w}\texttt{ }(-\infty<y<\infty)
\]
furnishes a parametric representation for the image of the left-hand branch of the hyperbola, the image of a point going downward along the entire left-hand branch is seen to move up the entire curve in the $u=c_{u}$ plane.

Finally, on the other hand, each branch of a hyperbola
\begin{equation}\label{eq:chap_2_2_3}
2y(x+z)=c_{v}\texttt{ }(c_{v}>0)
\end{equation}
is transformed into the $v=c_{v}$ plane parallel to the $uw$ coordinate plane with a distance $v=c_{v}$ in $\varpi$ space. To verify this, we note from the second of equations~(\ref{eq:chap_2_2_1}) that $v=c_{v}$ when $(x, y, z)$ is a point on either branch. Suppose that it lies on the branch lying in the first quadrant. Then, since $y=c_{v}/(2x+2z)$, the first of equations~(\ref{eq:chap_2_2_1}) reveals that the branch's image has parametric representation
\[
u=(x+z)^{2}-\frac{c_{v}}{4(x+z)^{2}}-c_{w},\texttt{ }v=c_{v},\texttt{ }w=c_{w}\texttt{ }(0<x+z<\infty)
\]
Observe that
\[
\lim_{x+z \to 0,\texttt{ }x>0}u=-\infty\texttt{ and }\lim_{x \to \infty}u=\infty,
\]
Since $u$ depends continuously on $x+z$, then, it is clear that as $(x, y, z)$ travels down the entire upper branch of hyperbola~(\ref{eq:chap_2_2_3}), its image moves to the right along the entire horizontal curve in the $v=c_{v}$ plane. Inasmuch as the image of the lower branch has parametric representation
\[
u=\frac{c_{v}}{4y^{2}}-y^{2}-c_{w},\texttt{ }v=c_{v},\texttt{ }w=c_{w}\texttt{ }(-\infty<y<0)
\]
and since
\[
\lim_{y \to -\infty}u=-\infty\texttt{ and }\lim_{y \to 0,\texttt{ }y<0}u=\infty,
\]
it follows that the image of a point moving upward along the entire lower branch also travels to the right along the entire curve in the $v=c_{v}$ plane.
\end{example}

We shall now use Example~(\ref{ex:chap_2_2_1}) to find the image of a certain region.

\begin{example}\label{ex:chap_2_2_2}
The domain $x>0$, $y>0$, $z>0$, and $(x+z)y<1$ consists of all points lying on the upper branches of hyperbolas from the family $2(x+z)y=c$, where $0<c<2$. We know from Example~(\ref{ex:chap_2_2_1}) that as a point travels downward along the entirety of one of these branches, its image under the transformation $\varpi=s^{2}$ moves to the right along the entire curve in the $v=c$ plane. Since, for all values of $c$ between $0$ and $2$, the branches fill out the domain $x>0$, $y>0$, $z>0$, and $(x+z)y<1$, that domain is mapped onto the horizontal strip $0<v<2$.

In view of equations~(\ref{eq:chap_2_2_1}), the image of a point $(0, y, 0)$ in the $s$ space is $(-y^{2}, 0, 0)$. Hence as $(0, y, 0)$ travels downward to the origin along the $y$ axis, its image moves to the right along the negative $u$ axis and reaches the origin in the $\varpi$ space. Then, since the image of a point $(x, 0, z)$ is $(x^{2}+2xz, 0, z^{2})$, that image moves to the right from the origin along the $u$ and $w$ axes as $(x, 0, z)$ moves to the right from the origin along the $x$ and $z$ axes, respectively. The image of the upper branch of the hyperbola $(x+z)y=1$ is, of course, the horizontal line $v=2$. Evidently, then, the closed region $x\geq0$, $\geq0$, $z\geq0$, and $(x+z)y\leq1$ is mapped onto the closed strip $0\leq v\leq2$.
\end{example}

Our last example here illustrates how spherical coordinates can be useful in analyzing certain mappings.

\begin{example}\label{ex:chap_2_2_3}
The mapping $\varpi=s^{2}$ becomes
\begin{equation}\label{eq:chap_2_2_4}
\varpi=r^{2}[e_{xy}(e^{i2\theta}+2e^{i\theta}\tan\varphi)\cos^{2}\varphi+\sin^{2}\varphi]
\end{equation}
when $s=r(e_{xy}e^{i\theta}\cos\varphi+\sin\varphi)$. Hence if $\varpi=\rho(e_{xy}e^{i\phi}\cos\psi+\sin\psi)$, we have
\[
\rho(e_{xy}e^{i\phi}\cos\psi+\sin\psi)=r^{2}[e_{xy}(e^{i2\theta}+2e^{i\theta}\tan\varphi)\cos^{2}\varphi+\sin^{2}\varphi];
\]
and the statement in italics near the beginning of Sec.~(\ref{sec:chap_1_8}) tells us that
\begin{equation}\label{eq:chap_2_2_5}
\rho=r^{2},\texttt{ }\psi=\arcsin(\sin^{2}\varphi),\texttt{ }\phi=\arctan(\frac{\sin2\theta+2\sin\theta\tan\varphi}{\cos2\theta+2\cos\theta\tan\varphi})+2k\pi,
\end{equation}
where $k$ has one of the values $k=0,\pm1,\pm2,\ldots$. Evidently, then, the image of any nonzero point $s$ is found by squaring the modulus of $s$ and calculating a pair of arguments $\arg\varpi=(\arg_{c}\varpi,\arg_{r}\varpi)$ (Sec.~(\ref{sec:chap_1_6})) from expressions~(\ref{eq:chap_2_2_5}).

Observe that points $s=r_{0}(e_{xy}e^{i\theta}\cos\varphi+\sin\varphi)$ on a sphere $r=r_{0}$ are transformed into points $\varpi=r_{0}^{2}(e_{xy}e^{i\phi}\cos\psi+\sin\psi)$ on the sphere $\rho=r_{0}^{2}$. As a point on the first sphere moves counterclockwise from the positive real axis to the positive imaginary axis in the $xy$ plane, its image on the second sphere moves counterclockwise from the positive real axis to the negative real axis in the $uv$ plane. So, as all possible positive values of $r_{0}$ are chosen, the corresponding arcs in the $s$ and $\varpi$ spaces fill out the first quadrant and the upper half space, respectively. The transformation $\varpi=s^{2}$ is, then, a one to one mapping of the first quadrant $r\geq0$, $0\leq\theta\leq\pi/2$, and $-\pi/2\leq\varphi\leq\pi/2$ in the $s$ space onto the upper half $\rho\geq0$, $0\leq\phi\leq\pi$, and $-\pi/2\leq\psi\leq\pi/2$ of the $\varpi$ space. The point $s=0$ is, of course, mapped onto the point $\varpi=0$.

The transformation $\varpi=s^{2}$ also maps the upper half space $r\geq0$, $0\leq\theta\leq\pi$, and $-\pi/2\leq\varphi\leq\pi/2$ onto the entire $\varpi$ space. However, in this case, the transformation is not one to one since both the positive and negative real axes in the $s$ space are mapped onto the positive real axis in the $\varpi$ space.

When $n$ is a positive integer greater than $2$, various mapping properties of the
transformation $\varpi=s^{n}$, or $\rho=r^{n}$, $\psi=\arcsin(\sin^{n}\varphi)$, and
\[
\phi=\arctan(\frac{\sum_{j=0}^{n-1}C_{n}^{j}\sin[(n-j)\theta]\tan^{j}\varphi}
{\sum_{j=0}^{n-1}C_{n}^{j}\cos[(n-j)\theta]\tan^{j}\varphi})+2k\pi\texttt{ }(C_{n}^{j}=\frac{n!}{(n-j)!j!})
\]
are similar to those of $\varpi=s^{2}$ where $k$ has one of the values $k=0,\pm1,\pm2,\ldots$. Such a transformation maps the entire $s$ space onto the entire $\varpi$ space, where each nonzero point in the $\varpi$ space is the image of $n$ distinct points in the $s$ space. The sphere $r=r_{0}$ is mapped onto the sphere $\rho=r_{0}^{n}$; and the sector $r\leq r_{0}$, $0\leq\theta\leq2\pi/n$, and $-\pi/2\leq\varphi\leq\pi/2$ is mapped onto the sphere $\rho\leq r_{0}^{n}$, but not in a one to one manner.
\end{example}

\section{Mappings by the Exponential Function}\label{sec:chap_2_3_13}

In Chap.~(\ref{ch:chap_3}) we shall introduce and develop properties of a number of elementary functions which do not involve polynomials. That chapter will contain the exponential function
\begin{equation}\label{eq:chap_2_3_1}
e^{s}=e^{e_{xy}(x+iy)}e^{z}\texttt{ }(s=e_{xy}(x+iy)+z),
\end{equation}
the three factors $e^{e_{xy}x}$, $e^{iy}$, and $e^{z}$ being well defined at this time (see Sec.~(\ref{sec:chap_1_6})). Note that definition~(\ref{eq:chap_2_3_1}), which can also be written
\[
e^{e_{xy}(x+iy)+z}=e^{z}e^{e_{xy}(x+iy)},
\]
is suggested by the familiar property
\[
e^{x_{1}+x_{2}}=e^{x_{1}}e^{x_{2}},
\]
of the exponential function in calculus.

The object of this section is to use the function $e^{s}$ to provide the reader with additional examples of mappings that continue to be reasonably simple. We begin by examining the images of vertical and horizontal lines.

\begin{example}\label{ex:chap_2_3_1}
The transformation
\begin{equation}\label{eq:chap_2_3_2}
\varpi=e^{s}
\end{equation}
can be written $\rho(e_{xy}e^{i\phi}\cos\psi+\sin\psi)$, where
\[
s=e_{xy}(x+iy)+z\texttt{ and }\varpi=\rho(e_{xy}e^{i\phi}\cos\psi+\sin\psi).
\]
Thus, transformation~(\ref{eq:chap_2_3_2}) can be expressed in the form
\begin{equation}\label{eq:chap_2_3_3}
    \begin{array}{cc}
        \rho=(e^{2x}+e^{2z})^{1/2},\texttt{ }\phi=y+2n\pi,\texttt{ }\psi=\arctan e^{z-x} & \texttt{ when }z\neq0 \\
        \rho=e^{x},\texttt{ }\phi=y+2n\pi,\texttt{ }\psi=0 & \texttt{ when }z=0
    \end{array}
\end{equation}
where $n$ is some integer (see Sec.~(\ref{sec:chap_1_8})).

The image of a typical point $s=(c_{x},y,c_{z})$ on a vertical line $x=c_{x}$ and $z=c_{z}$ has spherical coordinates
\[
\rho=\sqrt{\exp(2c_{x})+\exp(2c_{z})},\texttt{ }\phi=y,\texttt{ and }\psi=\arctan e^{c_{z}-c_{x}}
\]
when $c_{z}\neq0$ or $\rho=\exp(c_{x})$, $\phi=y$, and $\psi=0$ when $c_{z}=0$ in the $\varpi$ space. That image moves counterclockwise around the circle as $s$ moves up the line. The image of the line is evidently the entire circle; and each point on the circle is the image of an infinite number of points, spaced $2\pi$ units apart, along the line.

A horizontal line $y=c_{y}$ and $z=c_{z}$ is mapped in a one to one manner onto the ray $\phi=c_{y}$. To see that this is so, we note that the image of a point $s=(x,c_{y},c_{z})$ has spherical coordinates
\[
\rho=\sqrt{\exp(2x)+\exp(2c_{z})},\texttt{ }\phi=c_{y},\texttt{ and }\psi=\arctan e^{c_{z}-x}
\]
when $c_{z}\neq0$ or $\rho=\exp(x)$, $\phi=c_{y}$, and $\psi=0$ when $c_{z}=0$ in the $\varpi$ space. Evidently, then, as that point $s$ moves along the entire line from left to right, its image moves outward along the entire ray $\phi=c_{y}$.
\end{example}

Vertical and horizontal line \textit{segments} are mapped onto portions of circles and rays, respectively, and images of various regions are readily obtained from observations made in Example~(\ref{ex:chap_2_3_1}). This is illustrated in the following example.

\begin{example}\label{ex:chap_2_3_2}
Let us show that the transformation $\varpi=e^{s}$ maps the rectangular region $a\leq x\leq b$, $c\leq y\leq d$ in a plane $z=c_{z}$ where $c_{z}\neq0$ onto the region
\[
\sqrt{e^{2a}+e^{2c_{z}}}\leq\rho\leq\sqrt{e^{2b}+e^{2c_{z}}},\texttt{ }c\leq\phi\leq d,\texttt{ and }-\pi/2\leq\psi\leq\pi/2.
\]
There are two regions and corresponding parts of their boundaries. The vertical line segment $AD$ is mapped onto the arc
\[
\rho=\sqrt{e^{2a}+e^{2c_{z}}},\texttt{ }c\leq\phi\leq d,\texttt{ and }\psi=\arctan e^{c_{z}-x}
\]
in the $\varpi$ space, which is labeled $A'D'$. The images of vertical line segments to the right of $AD$ and joining the horizontal parts of the boundary are larger arcs; eventually, the image of the line segment $BC$ is the arc
\[
\rho=\sqrt{e^{2b}+e^{2c_{z}}},\texttt{ }c\leq\phi\leq d,\texttt{ and }\psi=\arctan e^{c_{z}-x}
\]
labeled $B'C'$. The mapping is one to one if $d-c<2\pi$. In particular, if $c=0$ and $d=\pi$, then $0\leq\phi\leq\pi$; and the rectangular region is mapped onto half of a circular ring.
\end{example}

Our final example here uses the images of horizontal lines to find the image of a horizontal strip.

\begin{example}\label{ex:chap_2_3_3}
When $\varpi=e^{s}$, the image of the infinite strip $0\leq y\leq\pi$ in a plane $z=c_{z}$ where $c_{z}\neq0$ is the upper half $v\geq0$ of the plane $w=\exp(c_{z})$ in the $\varpi$ space. This is seen by recalling from Example~(\ref{ex:chap_2_3_1}) how a horizontal line $y=c$ in the plane $z=c_{z}$ is transformed into a ray $\phi=c$ from the origin. As the real number $c$ increases from $c=0$ to $c=\pi$, the $y$ intercepts of the lines increase from $0$ to $\pi$ and the angles of inclination of the rays increase from $\phi=0$ to $\phi=\pi$.
\end{example}

\section{Limits}\label{sec:chap_2_4_14}

Let a function $f$ be defined at all points $s$ in some deleted neighborhood (Sec.~(\ref{sec:chap_1_10})) of $s_{0}$. The statement that the limit of $f(s)$ as $s$ approaches $s_{0}$ is a number $\varpi_{0}$, or that
\begin{equation}\label{eq:chap_2_4_1}
\lim_{s \to s_{0}}f(s)=\varpi_{0}
\end{equation}
means that the point $\varpi=f(s)$ can be made arbitrarily close to $\varpi_{0}$ if we choose the point $s$ close enough to $s_{0}$ but distinct from it. We now express the definition of limit in a precise and usable form.

Statement~(\ref{eq:chap_2_4_1}) means that, for each positive number $\varepsilon$, there is a positive number $\delta$ such that
\begin{equation}\label{eq:chap_2_4_2}
|f(s)-\varpi_{0}|<\varepsilon\texttt{ whenever }0<|s-s_{0}|<\delta,
\end{equation}
Geometrically, this definition says that, for each $\varepsilon$ neighborhood $\varpi-\varpi_{0}<\varepsilon$ of $\varpi_{0}$, there is a deleted $\delta$ neighborhood $0<|s-s_{0}|<\delta$ of $s_{0}$ such that every point $s$ in it has an image $\varpi$ lying in the $\varepsilon$ neighborhood. Note that even though all points in the deleted neighborhood $0<|s-s_{0}|<\delta$ are to be considered, their images need not fill up the entire neighborhood $|\varpi-\varpi_{0}|<\varepsilon$. If $f$ has the constant value $\varpi_{0}$, for instance, the image of $s$ is always the center of that neighborhood. Note, too, that once a $\delta$ has been found, it can be replaced by any smaller positive number, such as $\delta/2$.

It is easy to show that when a limit of a function $f(s)$ exists at a point $s_{0}$, it is unique. To do this, we suppose that
\[
\lim_{s \to s_{0}}f(s)=\varpi_{0}\texttt{ and }\lim_{s \to s_{0}}f(s)=\varpi_{1}.
\]

Then, for any positive number $\varepsilon$, there are positive numbers $\delta_{0}$ and $\delta_{1}$ such that
\[
|f(s)-\varpi_{0}|<\varepsilon\texttt{ whenever }0<|s-s_{0}|<\delta_{0}
\]
and
\[
|f(s)-\varpi_{1}|<\varepsilon\texttt{ whenever }0<|s-s_{0}|<\delta_{1}.
\]
So if $0<|s-s_{0}|<\delta$, where $\delta$ denotes the smaller of the two numbers $\delta_{0}$ and $\delta_{1}$ we
find that
\[
\begin{split}
|\varpi_{1}-\varpi_{0}| &= |[f(s)-\varpi_{0}]-[f(s)-\varpi_{1}]| \\
                        &\leq|f(s)-\varpi_{0}|+|f(s)-\varpi_{1}| \\
                        &<\varepsilon+\varepsilon=2\varepsilon.
\end{split}
\]
But $|\varpi_{1}-\varpi_{0}|$ is a nonnegative constant, and $\varepsilon$ can be chosen arbitrarily small. Hence
\[
\varpi_{1}-\varpi_{0}=0\texttt{ or }\varpi_{1}=\varpi_{0}.
\]

Definition~(\ref{eq:chap_2_4_2}) requires that $f$ be defined at all points in some deleted neighborhood of $s_{0}$. Such a deleted neighborhood, of course, always exists when $s_{0}$ is an interior point of a space on which $f$ is defined. We can extend the definition of limit to the case in which $s_{0}$ is a boundary point of the space by agreeing that the first of inequalities~(\ref{eq:chap_2_4_2}) need be satisfied by only those points $s$ that lie in both the space and the deleted neighborhood.

If $s_{0}$ is an interior point of the domain of definition of $f$, and limit~(\ref{eq:chap_2_4_1}) is to exist, the first of inequalities~(\ref{eq:chap_2_4_2}) must hold for all points in the deleted neighborhood $0<|s-s_{0}|<\delta$. Thus the symbol $s \to s_{0}$ implies that $s$ is allowed to approach $s_{0}$ in an arbitrary manner, not just from some particular direction.

While definition~(\ref{eq:chap_2_4_2}) provides a means of testing whether a given point $\varpi_{0}$ is a limit, it does not directly provide a method for determining that limit. Theorems on limits, presented in the next section, will enable us to actually find many limits.

\section{Theorems on Limits}\label{sec:chap_2_5_15}

We can expedite our treatment of limits by establishing a connection between limits of functions of a spatial complex variable and limits of real-valued functions of three real variables. Since limits of the latter type are studied in calculus, we use their definition and properties freely.

\begin{theorem}\label{th:chap_2_5_l}
Suppose that
\[
f(s)=e_{xy}[u(x,y,z)+iv(x,y,z)]+w(z),
\]
\[
s_{0}=e_{xy}(x_{0}+iy_{0})+z_{0},\texttt{ and }\varpi_{0}=e_{xy}(u_{0}+iv_{0})+w_{0}.
\]
Then
\begin{equation}\label{eq:chap_2_5_l}
\lim_{s \to s_{0}}f(s)=\varpi_{0}
\end{equation}
if and only if
\begin{equation}\label{eq:chap_2_5_2}
\begin{split}
        &\lim_{(x,y,z) \to (x_{0},y_{0},z_{0})}\texttt{ }u(x,y,z)=u_{0}, \\
        &\lim_{(x,y,z) \to (x_{0},y_{0},z_{0})}\texttt{ }v(x,y,z)=v_{0}, \\
        &\lim_{(x,y,z) \to (x_{0},y_{0},z_{0})}\texttt{ }w(z)=w_{0}.
\end{split}
\end{equation}
\end{theorem}

To prove the theorem, we first assume that limits~(\ref{eq:chap_2_5_2}) hold and obtain limit~(\ref{eq:chap_2_5_l}). Limits~(\ref{eq:chap_2_5_2}) tell us that, for each positive number $\varepsilon$, there exist positive numbers $\delta_{1}$, $\delta_{2}$, and $\delta_{3}$ such that
\begin{equation}\label{eq:chap_2_5_3}
|u-u_{0}|<\varepsilon/3\texttt{ whenever }0<\sqrt{(x-x_{0})^{2}+(y-y_{0})^{2}+(z-z_{0})^{2}}<\delta_{1},
\end{equation}
\begin{equation}\label{eq:chap_2_5_4}
|v-v_{0}|<\varepsilon/3\texttt{ whenever }0<\sqrt{(x-x_{0})^{2}+(y-y_{0})^{2}+(z-z_{0})^{2}}<\delta_{2},
\end{equation}
and
\begin{equation}\label{eq:chap_2_5_5}
|w-w_{0}|<\varepsilon/3\texttt{ whenever }0<\sqrt{(x-x_{0})^{2}+(y-y_{0})^{2}+(z-z_{0})^{2}}<\delta_{3}.
\end{equation}
Let $\delta$ denote the smaller of the three numbers $\delta_{1}$, $\delta_{2}$, and $\delta_{3}$. Since
\[
\begin{split}
|[e_{xy}(u+iv)+w]-[e_{xy}(u_{0}+iv_{0})+w_{0}]| &= |e_{xy}[(u-u_{0})+i(v-v_{0})]+(w-w_{0})| \\
                                                &\leq|u-u_{0}|+|v-v_{0}|+|w-w_{0}|
\end{split}
\]
and
\[
\begin{split}
\sqrt{(x-x_{0})^{2}+(y-y_{0})^{2}+(z-z_{0})^{2}} &= |e_{xy}[(x-x_{0})+i(y-y_{0})]+(z-z_{0})| \\
                                                 &= |[e_{xy}(x+iy)+z]-[e_{xy}(x_{0}+iy_{0})+z_{0}]|,
\end{split}
\]
it follows from statements~(\ref{eq:chap_2_5_3}),~(\ref{eq:chap_2_5_4}), and~(\ref{eq:chap_2_5_5}) that
\[
|[e_{xy}(u+iv)+w]-[e_{xy}(u_{0}+iv_{0})+w_{0}]|<\varepsilon/3+\varepsilon/3+\varepsilon/3=\varepsilon
\]
whenever
\[
0<|[e_{xy}(x+iy)+z]-[e_{xy}(x_{0}+iy_{0})+z_{0}]|<\delta.
\]
That is, limit~(\ref{eq:chap_2_5_l}) holds.

Let us now start with the assumption that limit~(\ref{eq:chap_2_5_l}) holds. With that assumption, we know that, for each positive numbers $\varepsilon$, there is a positive number $\delta$ such that
\begin{equation}\label{eq:chap_2_5_6}
|[e_{xy}(u+iv)+w]-[e_{xy}(u_{0}+iv_{0})+w_{0}]|<\varepsilon
\end{equation}
whenever
\begin{equation}\label{eq:chap_2_5_7}
0<|[e_{xy}(x+iy)+z]-[e_{xy}(x_{0}+iy_{0})+z_{0}]|<\delta.
\end{equation}
But
\[
\begin{split}
|u-u_{0}| &\leq|e_{xy}[(u-u_{0})+i(v-v_{0})]+(w-w_{0})| \\
          &=|[e_{xy}(u+iv)+w]-[e_{xy}(u_{0}+iv_{0})+w_{0}]|,
\end{split}
\]
\[
\begin{split}
|v-v_{0}| &\leq|e_{xy}[(u-u_{0})+i(v-v_{0})]+(w-w_{0})| \\
          &=|[e_{xy}(u+iv)+w]-[e_{xy}(u_{0}+iv_{0})+w_{0}]|,
\end{split}
\]
\[
\begin{split}
|w-w_{0}| &\leq|e_{xy}[(u-u_{0})+i(v-v_{0})]+(w-w_{0})| \\
          &=|[e_{xy}(u+iv)+w]-[e_{xy}(u_{0}+iv_{0})+w_{0}]|,
\end{split}
\]
and
\[
\begin{split}
|[e_{xy}(x+iy)+z]-[e_{xy}(x_{0}+iy_{0})+z_{0}]| &= |e_{xy}[(x-x_{0})+i(y-y_{0})]+(z-z_{0})| \\
                                                &= \sqrt{(x-x_{0})^{2}+(y-y_{0})^{2}+(z-z_{0})^{2}}.
\end{split}
\]
Hence it follows from inequalities~(\ref{eq:chap_2_5_6}) and~(\ref{eq:chap_2_5_7}) that
\[
|u-u_{0}|<\varepsilon,\texttt{ }|v-v_{0}|<\varepsilon,\texttt{ and }|w-w_{0}|<\varepsilon
\]
whenever
\[
0<\sqrt{(x-x_{0})^{2}+(y-y_{0})^{2}+(z-z_{0})^{2}}<\delta.
\]
This establishes limits~(\ref{eq:chap_2_5_2}), and the proof of the theorem is complete.

\begin{theorem}\label{th:chap_2_5_2}
Suppose that
\begin{equation}\label{eq:chap_2_5_8}
\lim_{s \to s_{0}}f(s)=\varpi_{0}\texttt{ and }\lim_{s \to s_{0}}F(s)=\overline{W}_{0}.
\end{equation}
Then
\begin{equation}\label{eq:chap_2_5_9}
\lim_{s \to s_{0}}[f(s)+F(s)]=\varpi_{0}+\overline{W}_{0},
\end{equation}
\begin{equation}\label{eq:chap_2_5_10}
\lim_{s \to s_{0}}[f(s)F(s)]=\varpi_{0}\overline{W}_{0},
\end{equation}
and, if $\overline{W}_{0}\neq0$
\begin{equation}\label{eq:chap_2_5_11}
\lim_{s \to s_{0}}\frac{f(s)}{F(s)}=\frac{\varpi_{0}}{\overline{W}_{0}}.
\end{equation}
\end{theorem}

This important theorem can be proved directly by using the definition of the limit of a function of a spatial complex variable. But, with the aid of Theorem~(\ref{th:chap_2_5_l}), it follows almost immediately from theorems on limits of real-valued functions of three real variables.

To verify property~(\ref{eq:chap_2_5_10}), for example, we write
\[
f(s)=e_{xy}[u(x,y,z)+iv(x,y,z)]+w(z),
\]
\[
F(s)=e_{xy}[U(x,y,z)+iV(x,y,z)]+W(z),
\]
\[
s_{0}=e_{xy}(x_{0}+iy_{0})+z_{0},\texttt{ }\varpi_{0}=e_{xy}(u_{0}+iv_{0})+w_{0},\texttt{ }\overline{W}_{0}=e_{xy}(U_{0}+iV_{0})+W_{0}.
\]
Then, according to hypotheses~(\ref{eq:chap_2_5_8}) and Theorem~(\ref{th:chap_2_5_l}), the limits as $(x,y,z)$ approaches $(x_{0},y_{0},z_{0})$ of the functions $u$, $v$, $w$, $U$, $V$, and $W$ exist and have the values $u_{0}$, $v_{0}$, $w_{0}$, $U_{0}$, $V_{0}$, and $W_{0}$, respectively. So the real and two-dimensional complex components of the product
\[
\begin{split}
f(s)F(s) &= [e_{xy}(u+iv)+w][e_{xy}(U+iV)+W] \\
         &= e_{xy}[(u+iv)(U+iV)+(u+iv)W+(U+iV)w]+wW \\
         &= e_{xy}[(uU-vV+uW+Uw)+i(vU+uV+vW+Vw)]+wW,
\end{split}
\]
have the limits
\[
u_{0}U_{0}-v_{0}V_{0}+u_{0}W_{0}+U_{0}w_{0},\texttt{ }v_{0}U_{0}+u_{0}V_{0}+v_{0}W_{0}+V_{0}w_{0},\texttt{ and }w_{0}W_{0},
\]
respectively, as $(x,y,z)$ approaches
$(x_{0},y_{0},z_{0})$. Hence, by Theorem~(\ref{th:chap_2_5_l}) again, $f(s)F(s)$ has the limit
\[
e_{xy}[(u_{0}U_{0}-v_{0}V_{0}+u_{0}W_{0}+U_{0}w_{0})+i(v_{0}U_{0}+u_{0}V_{0}+v_{0}W_{0}+V_{0}w_{0})]+w_{0}W_{0}
\]
as $s$ approaches $s_{0}$: and this is equal to $\varpi_{0}\overline{W}_{0}$. Property~(\ref{eq:chap_2_5_10}) is thus established. Corresponding verifications of properties~(\ref{eq:chap_2_5_9}) and~(\ref{eq:chap_2_5_11}) can be given.

It is easy to see from definition~(\ref{eq:chap_2_4_2}), Sec.~(\ref{sec:chap_2_4_14}), of  limit that
\[
\lim_{s \to s_{0}}c=c\texttt{ and }\lim_{s \to s_{0}}s=s_{0},
\]
where $s_{0}$ and  $c$ are any spatial complex numbers; and, by property~(\ref{eq:chap_2_5_10})) and mathematical
induction, it follows that
\[
\lim_{s \to s_{0}}s^{n}=s_{0}^{n}\texttt{ }(n=1,2,\ldots).
\]
So, in view of properties~(\ref{eq:chap_2_5_9}) and~(\ref{eq:chap_2_5_10}), the limit of a polynomial
\[
P(s)=a_{0}+a_{1}s+a_{2}s^{2}+\cdots+a_{n}s^{n}
\]
as $s$ approaches a point $s_{0}$ is the value of the polynomial at that point:
\begin{equation}\label{eq:chap_2_5_12}
\lim_{s \to s_{0}}P(s)=P(s_{0}).
\end{equation}

\section{Limits Involving the Point at Infinity}\label{sec:chap_2_6_16}

It is sometimes convenient to include with the complex space the point at infinity, denoted by $\infty$, and to use limits involving it. The complex space together with this point is called the extended complex space.

Let us agree that, in referring to a point $s$, we mean a point in the finite space. Hereafter, when the point at infinity is to be considered, it will be specifically mentioned.

A meaning is now readily given to the statement
\[
\lim_{s \to s_{0}}f(s)=\varpi_{0}
\]
when either $s_{0}$ or $\varpi_{0}$, or possibly each of these numbers, is replaced by the point at infinity. In the definition of limit in Sec.~(\ref{sec:chap_2_4_14}), we simply replace the appropriate neighborhoods of $s_{0}$ and $\varpi_{0}$ by neighborhoods of $\infty$. The proof of the following theorem illustrates how this is done.

\begin{theorem}\label{th:chap_2_6_1}
If $s_{0}$ and $\varpi_{0}$ are points in the $s$ and $\varpi$ spaces, respectively, then
\begin{equation}\label{eq:chap_2_6_1}
\lim_{s \to s_{0}}f(s)=\infty\texttt{ if and only if }\lim_{s \to s_{0}}\frac{1}{f(s)}=0
\end{equation}
and
\begin{equation}\label{eq:chap_2_6_2}
\lim_{s \to \infty}f(s)=\varpi_{0}\texttt{ if and only if }\lim_{s \to 0}f(\frac{1}{s})=\varpi_{0}.
\end{equation}
Moreover,
\begin{equation}\label{eq:chap_2_6_3}
\lim_{s \to \infty}f(s)=\infty\texttt{ if and only if }\lim_{s \to 0}\frac{1}{f(1/s)}=0.
\end{equation}
\end{theorem}

We start the proof by noting that the first of limits~(\ref{eq:chap_2_6_1}) means that, for each positive
number $\varepsilon$, there is a positive number $\delta$ such that
\begin{equation}\label{eq:chap_2_6_4}
|f(s)|>1/\varepsilon\texttt{ whenever }0<|s-s_{0}|<\delta,
\end{equation}
That is, the point $\varpi=f(s)$ lies in the $\varepsilon$ neighborhood $|\varpi|>1/\varepsilon$ of $\infty$ whenever $s$ lies in the deleted neighborhood $0<|s-s_{0}|<\delta$ of $s_{0}$. Since statement~(\ref{eq:chap_2_6_4}) can be written
\[
|\frac{1}{f(s)}-0|<\varepsilon\texttt{ whenever }0<|s-s_{0}|<\delta,
\]
the second of limits~(\ref{eq:chap_2_6_1}) follows.

The first of limits~(\ref{eq:chap_2_6_2}) means that, for each positive number $\varepsilon$, a positive number
$\delta$ exists such that
\begin{equation}\label{eq:chap_2_6_5}
|f(s)-\varpi_{0}|<\varepsilon\texttt{ whenever }|s|>1/\delta.
\end{equation}
Replacing $s$ by $1/s$ in statement~(\ref{eq:chap_2_6_5}) and then writing the result as
\[
|f(\frac{1}{s})-\varpi_{0}|<\varepsilon\texttt{ whenever }0<|s-0|<\delta,
\]
we arrive at the second of limits~(\ref{eq:chap_2_6_2}).

Finally, the first of limits~(\ref{eq:chap_2_6_3}) is to be interpreted as saying that, for each positive
number $\varepsilon$, there is a positive number $\delta$ such that
\begin{equation}\label{eq:chap_2_6_6}
|f(s)|>\frac{1}{\varepsilon}\texttt{ whenever }|s|>\frac{1}{\delta}.
\end{equation}
When $s$ is replaced by $1/s$, this statement can be put in the form
\[
|\frac{1}{f(1/s)}-0|<\varepsilon\texttt{ whenever }0<|s-0|<\delta;
\]
and this gives us the second of limits~(\ref{eq:chap_2_6_3}).

\section{Continuity}\label{sec:chap_2_7_17}

A function $f$ is continuous at a point $s_{0}$ all three of the following conditions are
satisfied:
\begin{equation}\label{eq:chap_2_7_1}
\lim_{s \to s_{0}}f(s)\texttt{ exists},
\end{equation}
\begin{equation}\label{eq:chap_2_7_2}
f(s_{0})\texttt{ exists},
\end{equation}
\begin{equation}\label{eq:chap_2_7_3}
\lim_{s \to s_{0}}f(s)=f(s_{0}).
\end{equation}
Observe that statement~(\ref{eq:chap_2_7_3}) actually contains statements~(\ref{eq:chap_2_7_1}) and~(\ref{eq:chap_2_7_2}), since the existence of the quantity on each side of the equation there is implicit. Statement~(\ref{eq:chap_2_7_3}) says that, for each positive number $\varepsilon$, there is a positive number $\delta$ such that
\begin{equation}\label{eq:chap_2_7_4}
|f(s)-f(s_{0})|<\varepsilon\texttt{ whenever }|s-s_{0}|<\delta,
\end{equation}

A function of a complex variable is said to be continuous in a space $S$ if it is continuous at each point in $S$.

If two functions are continuous at a point, their sum and product are also continuous at that point; their quotient is continuous at any such point where the denominator is not zero. These observations are direct consequences of Theorem~(\ref{th:chap_2_5_2}), Sec.~(\ref{sec:chap_2_5_15}). Note, too, that a polynomial is continuous in the entire plane because of limit~(\ref{eq:chap_2_5_12}), Sec.~(\ref{sec:chap_2_5_15}).

We turn now to two expected properties of continuous functions whose verifications are not so immediate. Our proofs depend on definition~(\ref{eq:chap_2_7_4}), and we present the results as theorems.

\begin{theorem}\label{th:chap_2_7_l}
A composition of continuous functions is itself continuous.
\end{theorem}

A precise statement of this theorem is contained in the proof to follow. We let $\varpi=f(s)$ be a function that is defined for all $s$ in a neighborhood $|s-s_{0}|<\delta$ of a point $s_{0}$, and we let $\overline{W}=g(\varpi)$ be a function whose domain of definition contains the image (Sec.~(\ref{sec:chap_2_2_12})) of that neighborhood under $f$. The composition $\overline{W}=g[f(s)]$ is, then, defined for all $s$ in the neighborhood $|s-s_{0}|<\delta$. Suppose now that $f$ is continuous at $s_{0}$ and that $g$ is continuous at the point $f(s_{0})$ in the $\varpi$ space. In view of the continuity of $g$ at $f(s_{0})$, there is, for each positive number $\varepsilon$, a positive number $\gamma$ such that
\[
|g[f(s)]-g[f(s_{0})]|<\varepsilon\texttt{ whenever }|f(s)-f(s_{0})|<\gamma.
\]
But the continuity of $f$ at $s_{0}$ ensures that the neighborhood $|s-s_{0}|<\delta$ can be made small enough that the second of these inequalities holds. The continuity of the composition $g[f(s)]$ is, therefore, established.

\begin{theorem}\label{th:chap_2_7_2}
If a function $f(s)$ is continuous and nonzero at a point $s_{0}$, then $f(s)\neq0$ throughout some neighborhood  of that point.
\end{theorem}

Assuming that $f(s)$ is, in fact, continuous and nonzero at $s_{0}$, we can prove Theorem~(\ref{th:chap_2_7_2}) by assigning the positive value $|f(s_{0})|/2$ to the numbers in statement~(\ref{eq:chap_2_7_4}). This tells us that there is a positive number $\delta$ such that
\[
|f(s)-f(s_{0})|<|f(s_{0})|/2\texttt{ whenever }|s-s_{0}|<\delta.
\]
So if there is a point $s$ in the neighborhood $|s-s_{0}|<\delta$ at which $f(s)=0$, we have the contradiction
\[
|f(s_{0})|<|f(s_{0})|/2;
\]
and the theorem is proved.

The continuity of a function
\begin{equation}\label{eq:chap_2_7_5}
f(s)=e_{xy}[u(x,y,z)+iv(x,y,z)]+w(z)
\end{equation}
is closely related to the continuity of its component functions $u(x,y,z)$, $v(x,y,z)$, and $w(z)$. We note, for instance, how it follows from Theorem~(\ref{th:chap_2_5_l}) in Sec.~(\ref{sec:chap_2_5_15}) that the function~(\ref{eq:chap_2_7_5}) is continuous at a point $s_{0}=(x_{0},y_{0},z_{0})$ if and only if its component functions are continuous there. To illustrate the use of this statement, suppose that the function~(\ref{eq:chap_2_7_5}) is continuous in a space $S$ that is both closed and bounded (see Sec.~(\ref{sec:chap_1_10})). The function
\[
\sqrt{[u(x,y,z)]^{2}+[v(x,y,z)]^{2}+[w(z)]^{2}}
\]
is then continuous in $S$ and thus reaches a maximum value somewhere in that space. That is, $f$ is bounded on $S$ and $|f(s)|$ reaches a maximum value somewhere in $S$. More precisely, there exists a nonnegative real number $M$ such that
\begin{equation}\label{eq:chap_2_7_6}
|f(s)|\leq M\texttt{ for all $s$ in $S$}
\end{equation}
where equality holds for at least one such $s$.

\section{Derivatives}\label{sec:chap_2_8_18}

Let $f$ be a function whose domain of definition contains a neighborhood of a point $s_{0}$. The derivative of $f$ at $s_{0}$, written $f'(s_{0})$, is defined by the equation
\begin{equation}\label{eq:chap_2_8_1}
f'(s_{0})=\lim_{s \to s_{0}}\frac{f(s)-f(s_{0})}{s-s_{0}},
\end{equation}
provided this limit exists. The function $f$ is said to be differentiable at $s_{0}$ when its derivative at $s_{0}$ exists.

By expressing the variable $s$ in definition~(\ref{eq:chap_2_8_1}) in terms of the new spatial complex variable
\[
\Delta s=s-s_{0},
\]
we can write that definition as
\begin{equation}\label{eq:chap_2_8_2}
f'(s_{0})=\lim_{\Delta s \to 0}\frac{f(s_{0}+\Delta s)-f(s_{0})}{\Delta s}.
\end{equation}
Note that, because $f$ is defined throughout a neighborhood of $s_{0}$, the number
\[
f(s_{0}+\Delta s)
\]
is always defined for $|\Delta s|$ sufficiently small.

When taking form~(\ref{eq:chap_2_8_2}) of the definition of derivative, we often drop the subscript on $s_{0}$ and introduce the number
\[
\Delta\varpi=f(s+\Delta s)-f(s),
\]
which denotes the change in the value of $f$ corresponding to a change $\Delta s$ in the point at which $f$ is  evaluated. Then, if we write $d\varpi/ds$ for $f'(s)$, equation~(\ref{eq:chap_2_8_2}) becomes
\begin{equation}\label{eq:chap_2_8_3}
\frac{d\varpi}{ds}=\lim_{\Delta s \to 0}\frac{\Delta\varpi}{\Delta s}.
\end{equation}

\section{Differentiation Formulas}\label{sec:chap_2_9_19}

The definition of derivative in Sec.~(\ref{sec:chap_2_8_18}) is identical in form to that of the derivative of a real-valued function of a real variable. In fact, the basic differentiation formulas given below can be derived from that definition by essentially the same steps as the ones used in calculus. In these formulas, the derivative  of a function $f$ at a point $s$ is denoted by either
\[
\frac{d}{ds}f(s)\texttt{ or }f'(s)
\]
depending on which notation is more convenient.

Let $c$ be a spatial complex constant and let $f$ be a function whose derivative exists at a point $s$. It is easy to show that
\begin{equation}\label{eq:chap_2_9_1}
\frac{d}{ds}c=0\texttt{, }\frac{d}{ds}s=1\texttt{, }\frac{d}{ds}[cf(s)]=cf'(s).
\end{equation}
Also, if $n$ is a positive integer,
\begin{equation}\label{eq:chap_2_9_2}
\frac{d}{ds}s^{n}=ns^{n-1}.
\end{equation}
This formula remains valid when $n$ is a negative integer, provided that $s\neq0$.

If the derivatives of two functions $f$ and $F$ exist at a point $s$, then
\begin{equation}\label{eq:chap_2_9_3}
\frac{d}{ds}[f(s)+F(s)]=f'(s)+F'(s),
\end{equation}
\begin{equation}\label{eq:chap_2_9_4}
\frac{d}{ds}[f(s)F(s)]=f'(s)F(s)+f(s)F'(s);
\end{equation}
and, when $F(s)\neq0$,
\begin{equation}\label{eq:chap_2_9_5}
\frac{d}{ds}[\frac{f(s)}{F(s)}]=\frac{f'(s)F(s)-f(s)F'(s)}{[F(s)]^{2}}.
\end{equation}

Let us derive formula~(\ref{eq:chap_2_9_4}). To do this, we write the following expression for the change in the product $\varpi=f(s)F(s)$:
\[
\begin{split}
\Delta\varpi &= f(s+\Delta s)F(s+\Delta s)-f(s)F(s) \\
             &= [f(s+\Delta s)-f(s)]F(s+\Delta s)+f(s)[F(s+\Delta s)-F(s)].
\end{split}
\]
Thus
\[
\frac{\Delta\varpi}{\Delta s}=\frac{f(s+\Delta s)-f(s)}{\Delta s}F(s+\Delta s)+f(s)\frac{F(s+\Delta s)-F(s)}{\Delta s};
\]
and, letting $\Delta s$ tend to zero, we arrive at the desired formula for the derivative of
$f(s)F(s)$. Here we have used the fact that $F$ is continuous at the point $s$, since $F'(s)$ exists; thus $F(s+\Delta s)$ tends to $F(s)$ as $\Delta s$ tends to zero.

There is also a chain rule for differentiating composite functions. Suppose that $f$
has a derivative at $s_{0}$ and that $g$ has a derivative at the point $f(s_{0})$. Then the function
$F(s)=g[f(s)]$ has a derivative at $s_{0}$, and
\begin{equation}\label{eq:chap_2_9_6}
F'(s_{0})=g'[f(s_{0})]f'(s_{0}).
\end{equation}

If we write $\varpi=f(s)$ and $\overline{W}=g(\varpi)$, so that $\overline{W}=F(s)$, the chain rule becomes
\[
\frac{d\overline{W}}{ds}=\frac{d\overline{W}}{d\varpi}\frac{d\varpi}{ds}.
\]

To start the proof of formula~(\ref{eq:chap_2_9_6}), choose a specific point $s_{0}$ at which $f'(s_{0})$
exists. Write $\varpi_{0}=f(s_{0})$ and also assume that $g'(\varpi_{0})$ exists. There is, then, some
$\varepsilon$ neighborhood $|\varpi-\varpi_{0}|<\varepsilon$ of $\varpi_{0}$ such that, for all points $\varpi$ in  that neighborhood, we can define a function $\Phi$ which has the values $\Phi(\varpi_{0})=0$ and
\begin{equation}\label{eq:chap_2_9_7}
\Phi(\varpi)=\frac{g(\varpi)-g(\varpi_{0})}{\varpi-\varpi_{0}}-g'(\varpi_{0})\texttt{ when }\varpi\neq\varpi_{0}.
\end{equation}
Note that, in view of the definition of derivative,
\begin{equation}\label{eq:chap_2_9_8}
\lim_{\varpi \to \varpi_{0}}\Phi(\varpi)=0.
\end{equation}
Hence $\Phi$ is continuous at $\varpi_{0}$.

Now expression~(\ref{eq:chap_2_9_7}) can be put in the form
\begin{equation}\label{eq:chap_2_9_9}
g(\varpi)-g(\varpi_{0})=[g'(\varpi_{0})+\Phi(\varpi)](\varpi-\varpi_{0})\texttt{ }(|\varpi-\varpi_{0}|<\varepsilon),
\end{equation}
which is valid even when $\varpi=\varpi_{0}$; and, since $f'(s_{0})$ exists and $f$ is, therefore, continuous at $s_{0}$, we can choose a positive number $\delta$ such that the point $f(s)$ lies in the $\varepsilon$ neighborhood  $|\varpi-\varpi_{0}|<\varepsilon$ of $\varpi_{0}$ if $s$ lies in the $\delta$ neighborhood $|s-s_{0}|<\delta$ of $s_{0}$. Thus it is legitimate to replace the variable $\varpi$ in equation~(\ref{eq:chap_2_9_9}) by $f(s)$ when $s$ is any point in the neighborhood $|s-s_{0}|<\delta$. With that substitution, and with $\varpi_{0}=f(s_{0})$, equation~(\ref{eq:chap_2_9_9}) becomes
\begin{equation}\label{eq:chap_2_9_10}
\frac{g(\varpi)-g(\varpi_{0})}{s-s_{0}}=[g'(\varpi_{0})+\Phi(\varpi)]\frac{f(s)-f(s_{0})}{s-s_{0}}\texttt{ }(0<|s-s_{0}|<\delta),
\end{equation}
where we must stipulate that $s\neq s_{0}$ so that we are not dividing by zero. As already noted, $f$ is continuous at $s_{0}$ and $\Phi$ is continuous at the point $\varpi_{0}=f(s_{0})$. Thus the composition $\Phi(f(s))$ is continuous at $s_{0}$; and, since $\Phi(\varpi_{0})=0$,
\[
\lim_{s \to s_{0}}\Phi(f(s))=0.
\]
So equation~(\ref{eq:chap_2_9_10}) becomes equation~(\ref{eq:chap_2_9_6}) in the limit as $s$ approaches $s_{0}$.

\section{Cauchy-Riemann Equation}\label{sec:chap_2_10_20}

In this section, we obtain a pair of equations that the first-order partial derivatives of the component functions $u$ and $v$ of a spatial complex function
\begin{equation}\label{eq:chap_2_10_1}
f(s)=e_{xy}[u(x,y,z)+iv(x,y,z)]+w(z)
\end{equation}
must satisfy at a point $s_{0}=(x_{0},y_{0},z_{0})$ when the derivative of $f$ exists there. We also show how to express $f'(s_{0})$ in terms of those partial derivatives.

We start by writing $s_{0}=e_{xy}(x_{0}+iy_{0})+z_{0}$, $\Delta s=e_{xy}(\Delta x+i\Delta y)+\Delta z$, and
\[
\begin{split}
\Delta\varpi &= f(s_{0}+\Delta s)-f(s_{0}) \\
             &= e_{xy}\{[u(x_{0}+\Delta x,y_{0}+\Delta y,z_{0}+\Delta z)-u(x_{0},y_{0},z_{0})] \\
             &+ i[v(x_{0}+\Delta x,y_{0}+\Delta y,z_{0}+\Delta z)-v(x_{0},y_{0},z_{0})]\} \\
             &+ [w(z_{0}+\Delta z)-w(z_{0})] \\
             &= e_{xy}(\Delta u+i\Delta v)+\Delta w.
\end{split}
\]
Thus we have
\[
\begin{split}
\frac{\Delta\varpi}{\Delta s} &= \frac{e_{xy}(\Delta u+i\Delta v)+\Delta w}{e_{xy}(\Delta x+i\Delta y)+\Delta z}\\
                              &= e_{xy}(\frac{\Delta u+i\Delta v+\Delta w}{\Delta x+i\Delta y+\Delta z}-\frac{\Delta w}{\Delta z})+\frac{\Delta w}{\Delta z}.
\end{split}
\]

Assuming that the derivative
\begin{equation}\label{eq:chap_2_10_2}
f'(s_{0})=\lim_{\Delta s \to 0}\frac{\Delta\varpi}{\Delta s}
\end{equation}
exists, we know from Theorem~(\ref{th:chap_2_5_l}) in Sec.~(\ref{sec:chap_2_5_15}) that
\begin{equation}\label{eq:chap_2_10_3}
f'(s_{0})=e_{xy}\lim_{(\Delta x,\Delta y,\Delta z) \to (0,0,0)}Im_{s}\frac{\Delta\varpi}{\Delta s}
+\lim_{(\Delta x,\Delta y,\Delta z) \to (0,0,0)}Re_{s}\frac{\Delta\varpi}{\Delta s}
\end{equation}
where
\[
Im_{s}\frac{\Delta\varpi}{\Delta s}
=\frac{\Delta u+i\Delta v+\Delta w}{\Delta x+i\Delta y+\Delta z}-\frac{\Delta w}{\Delta z}\texttt{ and }
Re_{s}\frac{\Delta\varpi}{\Delta s}=\frac{\Delta w(z)}{\Delta z}.
\]

Now it is important to keep in mind that expression~(\ref{eq:chap_2_10_3}) is valid as $(\Delta x,\Delta y,\Delta z)$ tends to $(0,0,0)$ in any manner that we may choose just as in the case of the two-dimensional complex plane. In particular, we first let $(\Delta x,\Delta y,\Delta z)$ tend to $(0,0,0)$ through the point $(\Delta x,0,0)$ along the $x$ axis. Inasmuch as $\Delta y=\Delta z=0$, the quotient $\Delta\varpi/\Delta s$ becomes
\[
\frac{\Delta\varpi}{\Delta s}=e_{xy}[\frac{u(x_{0}+\Delta x,y_{0},z_{0})-u(x_{0},y_{0},z_{0})}{\Delta x}
+i\frac{v(x_{0}+\Delta x,y_{0},z_{0})-v(x_{0},y_{0},z_{0})}{\Delta x}].
\]
Thus
\[
\begin{split}
\lim_{(\Delta x,\Delta y,\Delta z) \to (0,0,0)}Re\frac{\Delta\varpi}{\Delta s}
&=\lim_{\Delta x \to 0}e_{xy}\frac{u(x_{0}+\Delta x,y_{0},z_{0})-u(x_{0},y_{0},z_{0})}{\Delta x} \\
&=e_{xy}u_{x}(x_{0},y_{0},z_{0})
\end{split}
\]
and
\[
\begin{split}
\lim_{(\Delta x,\Delta y,\Delta z) \to (0,0,0)}Im\frac{\Delta\varpi}{\Delta s}
&=\lim_{\Delta x \to 0}e_{xy}\frac{v(x_{0}+\Delta x,y_{0},z_{0})-v(x_{0},y_{0},z_{0})}{\Delta x} \\
&=e_{xy}v_{x}(x_{0},y_{0},z_{0}).
\end{split}
\]
where $u_{x}(x_{0},y_{0},z_{0})$, $v_{x}(x_{0},y_{0},z_{0})$, and $w_{x}(x_{0},y_{0},z_{0})$ denote the first-order partial derivatives with respect to $x$ of the functions $u$, $v$, and $w$, respectively, at $(x_{0},y_{0},z_{0})$. Substitution of these limits into expression~(\ref{eq:chap_2_10_3}) tells us that
\begin{equation}\label{eq:chap_2_10_4}
f'(s_{0})=e_{xy}[u_{x}(x_{0},y_{0},z_{0})+iv_{x}(x_{0},y_{0},z_{0})].
\end{equation}

Next, we might have let $(\Delta x,\Delta y,\Delta z)$ tend to $(0,0,0)$ through the point $(0,\Delta y,0)$ along the $y$ axis. In that case, $\Delta x=\Delta z=0$ and
\[
\begin{split}
\frac{\Delta\varpi}{\Delta s} &= e_{xy}[\frac{u(x_{0},y_{0}+\Delta y,z_{0})-u(x_{0},y_{0},z_{0})}{i\Delta y}
+i\frac{v(x_{0},y_{0}+\Delta y,z_{0})-v(x_{0},y_{0},z_{0})}{i\Delta y}] \\
                              &= e_{xy}[\frac{v(x_{0},y_{0}+\Delta y,z_{0})-v(x_{0},y_{0},z_{0})}{\Delta y}
-i\frac{u(x_{0},y_{0}+\Delta y,z_{0})-u(x_{0},y_{0},z_{0})}{\Delta y}].
\end{split}
\]
Evidently, then,
\[
\begin{split}
\lim_{(\Delta x,\Delta y,\Delta z) \to (0,0,0)}Re\frac{\Delta\varpi}{\Delta s}
&=\lim_{\Delta y \to 0}e_{xy}\frac{v(x_{0},y_{0}+\Delta y,z_{0})-v(x_{0},y_{0},z_{0})}{\Delta y} \\
&=e_{xy}v_{y}(x_{0},y_{0},z_{0})
\end{split}
\]
and
\[
\begin{split}
\lim_{(\Delta x,\Delta y,\Delta z) \to (0,0,0)}Im\frac{\Delta\varpi}{\Delta s}
&=-\lim_{\Delta y \to 0}e_{xy}\frac{u(x_{0},y_{0}+\Delta y,z_{0})-u(x_{0},y_{0},z_{0})}{\Delta y} \\
&=-e_{xy}u_{y}(x_{0},y_{0},z_{0}).
\end{split}
\]
Hence it follows from expression~(\ref{eq:chap_2_10_3}) that
\begin{equation}\label{eq:chap_2_10_5}
f'(s_{0})=e_{xy}[v_{y}(x_{0},y_{0},z_{0})-iu_{y}(x_{0},y_{0},z_{0})]
\end{equation}
where the partial derivatives of the spatial functions $u$ and $v$ are, this time, with respect to $y$. Note that
equation~(\ref{eq:chap_2_10_5}) can also be written in the form
\[
f'(s_{0})=-ie_{xy}[u_{y}(x_{0},y_{0},z_{0})+iv_{y}(x_{0},y_{0},z_{0})].
\]

Finally, we might have let $(\Delta x,\Delta y,\Delta z)$ tend to $(0,0,0)$ through the point $(0,0,\Delta z)$ along the $z$ axis. In that case, $\Delta x=\Delta y=0$ and
\[
\frac{\Delta\varpi}{\Delta s}=\frac{w(z_{0}+\Delta z)-w(z_{0})}{\Delta z}
\]
\[
+e_{xy}[\frac{u(x_{0},y_{0},z_{0}+\Delta z)-u(x_{0},y_{0},z_{0})}{\Delta z}
+i\frac{v(x_{0},y_{0},z_{0}+\Delta z)-v(x_{0},y_{0},z_{0})}{\Delta z}].
\]
Evidently, then,
\[
\begin{split}
\lim_{(\Delta x,\Delta y,\Delta z) \to (0,0,0)}Re\frac{\Delta\varpi}{\Delta s}
&=\lim_{\Delta z \to 0}e_{xy}\frac{u(x_{0},y_{0},z_{0}+\Delta z)-u(x_{0},y_{0},z_{0})}{\Delta z} \\
&=e_{xy}u_{z}(x_{0},y_{0},z_{0}),
\end{split}
\]
\[
\begin{split}
\lim_{(\Delta x,\Delta y,\Delta z) \to (0,0,0)}Im\frac{\Delta\varpi}{\Delta s}
&=\lim_{\Delta z \to 0}e_{xy}\frac{v(x_{0},y_{0},z_{0}+\Delta z)-v(x_{0},y_{0},z_{0})}{\Delta z} \\
&=e_{xy}v_{z}(x_{0},y_{0},z_{0}),
\end{split}
\]
and
\[
\lim_{(\Delta x,\Delta y,\Delta z) \to (0,0,0)}Re_{s}\frac{\Delta\varpi}{\Delta s}
=\lim_{\Delta z \to 0}\frac{w(z_{0}+\Delta z)-w(z_{0})}{\Delta z}=w_{z}(z_{0})
\]
Hence it follows from expression~(\ref{eq:chap_2_10_3}) that
\begin{equation}\label{eq:chap_2_10_6}
f'(s_{0})=e_{xy}[u_{z}(x_{0},y_{0},z_{0})+iv_{z}(x_{0},y_{0},z_{0})]+w_{z}(z_{0})
\end{equation}
where the partial derivatives of the spatial functions $u$, $v$, and $w$ are, this time, with respect to $z$.

Equations~(\ref{eq:chap_2_10_4}),~(\ref{eq:chap_2_10_5}), and~(\ref{eq:chap_2_10_6}) not only give $f'(s_{0})$ in terms of partial derivatives of the component functions $u$, $v$, and $w$, but they also provide necessary conditions for the existence of $f'(s_{0})$. For, on equating the real and imaginary parts on the right-hand sides of these equations, we see that the existence of $f'(s_{0})$ requires that
\begin{equation}\label{eq:chap_2_10_7}
\begin{split}
& e_{xy}u_{x}(x_{0},y_{0},z_{0})=e_{xy}v_{y}(x_{0},y_{0},z_{0})=e_{xy}u_{z}(x_{0},y_{0},z_{0})+w_{z}(z_{0}) \\
& e_{xy}v_{x}(x_{0},y_{0},z_{0})=e_{xy}v_{z}(x_{0},y_{0},z_{0})=-e_{xy}u_{y}(x_{0},y_{0},z_{0})
\end{split}
\end{equation}
Equations~(\ref{eq:chap_2_10_7}) are the Cauchy-Riemann equations of spatial complex functions.

We summarize the above results as follows.

\begin{theorem}\label{th:chap_2_10_l}
Suppose that
\[
f(s)=e_{xy}[u(x,y,z)+iv(x,y,z)]+w(z)
\]
and that $f'(s)$ exists at a point $s_{0}=e_{xy}(x_{0}+iy_{0})+z_{0}$. Then the first-order partial derivatives
of $u$, $v$ and $w$ must exist at $(x_{0},y_{0},z_{0})$, and they must satisfy the Cauchy-Riemann equations of spatial complex functions
\begin{equation}\label{eq:chap_2_10_8}
e_{xy}u_{x}=e_{xy}v_{y}=e_{xy}u_{z}+w_{z},\texttt{ }
e_{xy}v_{x}=e_{xy}v_{z}=-e_{xy}u_{y}.
\end{equation}
Also, $f'(s_{0})$ can be written
\begin{equation}\label{eq:chap_2_10_9}
\begin{split}
f'(s_{0}) &= e_{xy}(u_{x}+iv_{x}) \\
          &= e_{xy}(v_{y}-iu_{y}) \\
          &= e_{xy}(u_{z}+iv_{z})+w_{z}
\end{split}
\end{equation}
and these partial derivatives are to be evaluated at $(x_{0},y_{0},z_{0})$.
\end{theorem}

\begin{example}
Let
\[
f(s)=s^{2}=e_{xy}[u(x,y,z)+iv(x,y,z)]+w(z)
\]
where
\[
u(x,y,z)=x^{2}-y^{2}+2xz,\texttt{ }
v(x,y,z)=2xy+2yz,\texttt{ }
w(z)=z^{2}.
\]
Then $f'(s)$ exists at a point $s_{0}=e_{xy}(x_{0}+iy_{0})+z_{0}$. The first-order partial derivatives
of $u$, $v$ and $w$ exist at $(x_{0},y_{0},z_{0})$ and satisfy the Cauchy-Riemann equations
\[
e_{xy}u_{x}=e_{xy}v_{y}=e_{xy}u_{z}+w_{z}=2(e_{xy}x+z),
\]
and
\[
e_{xy}v_{x}=e_{xy}v_{z}=-e_{xy}u_{y}=2e_{xy}y
\]
where from remark~(\ref{re:chap_1_1_1}) for a term $e_{xy}x^{j}y^{k}$, the operator $e_{xy}$ becomes $e_{xy}^{m}$ where $m=|j|+|k|$, and
\[
\begin{split}
f'(s_{0}) &= e_{xy}(u_{x}+iv_{x}) \\
          &= e_{xy}(v_{y}-iu_{y}) \\
          &= e_{xy}(u_{z}+iv_{z})+w_{z} \\
          &= 2[e_{xy}(x_{0}+iy_{0})+z_{0}]=2s_{0}.
\end{split}
\]
\end{example}

\section{Sufficient Conditions for Differentiability}\label{sec:chap_2_11_21}

Satisfaction of the Cauchy-Riemann equations at a point $s_{0}=(x_{0},y_{0},z_{0})$ is not sufficient to ensure the existence of the derivative of a function $f(s)$ at that point. But, with certain continuity conditions, we have the following useful theorem.

\begin{theorem}\label{th:chap_2_11_l}
Let the spatial complex function
\[
f(s)=e_{xy}[u(x,y,z)+iv(x,y,z)]+w(z)
\]
be defined throughout some $\varepsilon$ neighborhood of a point $s_{0}=e_{xy}(x_{0}+iy_{0})+z_{0}$, and suppose
that the first-order partial derivatives of the functions $u$, $v$ and $w$ with respect to $x$, $y$, and $z$
exist everywhere in that neighborhood. If those partial derivatives are continuous at $(x_{0},y_{0},z_{0})$ and satisfy the Cauchy-Riemann equations
\[
e_{xy}u_{x}=e_{xy}v_{y}=e_{xy}u_{z}+w_{z},\texttt{ }
e_{xy}v_{x}=e_{xy}v_{z}=-e_{xy}u_{y},
\]
at $(x_{0},y_{0},z_{0})$ then $f'(s_{0})$ exists.
\end{theorem}

To start the proof, we write $\Delta s=e_{xy}(\Delta x+i\Delta y)+\Delta z$ where $0<|\Delta s|<\varepsilon$, and
\begin{equation}\label{eq:chap_2_11_1}
\Delta\varpi=f(s_{0}+\Delta s)-f(s_{0})=e_{xy}(\Delta u+i\Delta v)+\Delta w
\end{equation}
where
\[
\Delta u=u(x_{0}+\Delta x,y_{0}+\Delta y,z_{0}+\Delta z)-u(x_{0},y_{0},z_{0}),
\]
\[
\Delta v=v(x_{0}+\Delta x,y_{0}+\Delta y,z_{0}+\Delta z)-v(x_{0},y_{0},z_{0}),
\]
and
\[
\Delta w=w(z_{0}+\Delta z)-w(z_{0}).
\]
The assumption that the first-order partial derivatives of $u$, $v$, and $w$ are continuous at the point $(x_{0},y_{0},z_{0})$ enables us to write
\begin{equation}\label{eq:chap_2_11_2}
\begin{split}
\Delta u &= u_{x}(x_{0},y_{0},z_{0})\Delta x+u_{y}(x_{0},y_{0},z_{0})\Delta y+u_{z}(x_{0},y_{0},z_{0})\Delta z \\
         &+ \varepsilon_{1}\sqrt{(\Delta x)^{2}+(\Delta y)^{2}+(\Delta z)^{2}},
\end{split}
\end{equation}
\begin{equation}\label{eq:chap_2_11_3}
\begin{split}
\Delta v &= v_{x}(x_{0},y_{0},z_{0})\Delta x+v_{y}(x_{0},y_{0},z_{0})\Delta y+v_{z}(x_{0},y_{0},z_{0})\Delta z \\
         &+ \varepsilon_{2}\sqrt{(\Delta x)^{2}+(\Delta y)^{2}+(\Delta z)^{2}},
\end{split}
\end{equation}
and
\begin{equation}\label{eq:chap_2_11_3w}
\Delta w=w_{z}(z_{0})\Delta z+\varepsilon_{3}|\Delta z|
\end{equation}
where $\varepsilon_{1}$, $\varepsilon_{2}$, and $\varepsilon_{3}$ tend to $0$ as $(\Delta x,\Delta y,\Delta z)$ approaches $(0,0,0)$ in the $\Delta s$ space. Substitution of expressions~(\ref{eq:chap_2_11_2}),~(\ref{eq:chap_2_11_3}), and~(\ref{eq:chap_2_11_3w}) into equation~(\ref{eq:chap_2_11_1}) now tells us that
\begin{equation}\label{eq:chap_2_11_4}
\begin{split}
\Delta\varpi &= e_{xy}\{[u_{x}(x_{0},y_{0},z_{0})\Delta x+u_{y}(x_{0},y_{0},z_{0})\Delta y+u_{z}(x_{0},y_{0},z_{0})\Delta z] \\
&+ i[v_{x}(x_{0},y_{0},z_{0})\Delta x+v_{y}(x_{0},y_{0},z_{0})\Delta y+v_{z}(x_{0},y_{0},z_{0})\Delta z]\} \\
&+ w_{z}(z_{0})\Delta z+e_{xy}[(\varepsilon_{1}+i\varepsilon_{2})\sqrt{(\Delta x)^{2}+(\Delta y)^{2}+(\Delta z)^{2}}]
+\varepsilon_{3}|\Delta z|.
\end{split}
\end{equation}

Assuming that the Cauchy-Riemann equations are satisfied at $(x_{0},y_{0},z_{0})$, we can
replace $u_{y}(x_{0},y_{0},z_{0})$ by $-v_{x}(x_{0},y_{0},z_{0})$ and $v_{y}(x_{0},y_{0},z_{0})$ by $u_{x}(x_{0},y_{0},z_{0})$ in equation~(\ref{eq:chap_2_11_4}) and then divide through by $\Delta s$ to get
\begin{equation}\label{eq:chap_2_11_5}
\begin{split}
\frac{\Delta\varpi}{\Delta s} &= e_{xy}[u_{x}(x_{0},y_{0},z_{0})+iv_{x}(x_{0},y_{0},z_{0})]\frac{\Delta x+i\Delta y}{\Delta s} \\
&+ \{e_{xy}[u_{z}(x_{0},y_{0},z_{0})+iv_{z}(x_{0},y_{0},z_{0})]+w_{z}(z_{0})\}\frac{\Delta z}{\Delta s} \\
&+ e_{xy}(\varepsilon_{1}+i\varepsilon_{2})\frac{\sqrt{(\Delta x)^{2}+(\Delta y)^{2}+(\Delta z)^{2}}}{\Delta s}
+\varepsilon_{3}\frac{|\Delta z|}{\Delta s}.
\end{split}
\end{equation}
But $|\Delta x+i\Delta y|+|\Delta z|=|\Delta s|$, and $\sqrt{(\Delta x)^{2}+(\Delta y)^{2}+(\Delta z)^{2}}=|\Delta s|$, and so
\[
|\frac{\Delta z}{\Delta s}|\leq1,\texttt{ }
|\frac{\Delta x+i\Delta y}{\Delta s}|+
|\frac{\Delta z}{\Delta s}|=1,\texttt{ and }
|\frac{\sqrt{(\Delta x)^{2}+(\Delta y)^{2}+(\Delta z)^{2}}}{\Delta s}|=1.
\]
Also, $\varepsilon_{1}+i\varepsilon_{2}$ and $\varepsilon_{3}$ tend to $0$ as $(\Delta x,\Delta y,\Delta z)$ approaches $(0,0,0)$. So the last two terms on the right in equation~(\ref{eq:chap_2_11_5}) tend to $0$ as the variable $\Delta s=e_{xy}(\Delta x+i\Delta y)+\Delta z$ tends to $0$. This means that the limit of the left-hand side of equation~(\ref{eq:chap_2_11_5}) exists and there are
\begin{equation}\label{eq:chap_2_11_6}
\begin{split}
f'(s_{0}) &= e_{xy}(u_{x}+iv_{x}) \\
          &= e_{xy}(v_{y}-iu_{y}) \\
          &= e_{xy}(u_{z}+iv_{z})+w_{z},
\end{split}
\end{equation}
and these partial derivatives are to be evaluated at $(x_{0},y_{0},z_{0})$. So Theorem~(\ref{th:chap_2_11_l}) holds and the proof of Theorem~(\ref{th:chap_2_11_l}) is completed.

\section{Spherical Coordinates}\label{sec:chap_2_12_22}

Assuming that $s_{0}\neq0$, we shall in this section use the coordinate transformation
\begin{equation}\label{eq:chap_2_12_1}
x=r\cos\varphi\cos\theta\texttt{, }y=r\cos\varphi\sin\theta\texttt{, and }z=r\sin\varphi.
\end{equation}
to restate the theorem in Sec.~(\ref{sec:chap_2_11_21}) in spherical coordinates.

Depending on whether we write
\[
s=e_{xy}(x+iy)+z\texttt{ or }s=r(e_{xy}e^{i\theta}\cos\varphi+\sin\varphi)
\]
when $\varpi=f(s)$, the complex and real parts of $\varpi=e_{xy}(u+iv)+w$ are expressed in terms of either the variables $x$, $y$, and $z$ or $r$, $\theta$, and $\varphi$. Suppose that the first-order partial derivatives of $u$, $v$, and $w$ with respect to $x$, $y$, and $z$ exist everywhere in some neighborhood of a given nonzero point  $s_{0}$ and are continuous at that point. The first-order partial derivatives with respect to $r$, $\theta$, and $\varphi$ also have these properties, and the chain rule for differentiating real-valued functions of three real variables can be used to write them in terms of the ones with respect to $x$, $y$, and $z$. More precisely, since
\[
\begin{split}
& \frac{\partial u}{\partial r}=\frac{\partial u}{\partial x}\frac{\partial x}{\partial r}+\frac{\partial u}{\partial y}\frac{\partial y}{\partial r}+\frac{\partial u}{\partial z}\frac{\partial z}{\partial r}, \\
& \frac{\partial u}{\partial\theta}=\frac{\partial u}{\partial x}\frac{\partial x}{\partial\theta}+\frac{\partial u}{\partial y}\frac{\partial y}{\partial\theta}+\frac{\partial u}{\partial z}\frac{\partial z}{\partial\theta}, \\
& \frac{\partial u}{\partial\varphi}=\frac{\partial u}{\partial x}\frac{\partial x}{\partial\varphi}+\frac{\partial u}{\partial y}\frac{\partial y}{\partial\varphi}+\frac{\partial u}{\partial z}\frac{\partial z}{\partial\varphi},
\end{split}
\]
and
\[
\begin{split}
& \frac{\partial v}{\partial r}=\frac{\partial v}{\partial x}\frac{\partial x}{\partial r}+\frac{\partial v}{\partial y}\frac{\partial y}{\partial r}+\frac{\partial v}{\partial z}\frac{\partial z}{\partial r}, \\
& \frac{\partial v}{\partial\theta}=\frac{\partial v}{\partial x}\frac{\partial x}{\partial\theta}+\frac{\partial v}{\partial y}\frac{\partial y}{\partial\theta}+\frac{\partial v}{\partial z}\frac{\partial z}{\partial\theta}, \\
& \frac{\partial v}{\partial\varphi}=\frac{\partial v}{\partial x}\frac{\partial x}{\partial\varphi}+\frac{\partial v}{\partial y}\frac{\partial y}{\partial\varphi}+\frac{\partial v}{\partial z}\frac{\partial z}{\partial\varphi},
\end{split}
\]
and
\[
\begin{split}
& \frac{\partial w}{\partial r}=\frac{\partial w}{\partial z}\frac{\partial z}{\partial r}, \\
& \frac{\partial w}{\partial\varphi}=\frac{\partial w}{\partial z}\frac{\partial z}{\partial\varphi},
\end{split}
\]
one can write
\begin{equation}\label{eq:chap_2_12_2}
\begin{split}
u_{r} &= u_{x}\cos\varphi\cos\theta+u_{y}\cos\varphi\sin\theta+u_{z}\sin\varphi, \\
u_{\theta} &= r(-u_{x}\cos\varphi\sin\theta+u_{y}\cos\varphi\cos\theta), \\
u_{\varphi} &= -r(u_{x}\sin\varphi\cos\theta+u_{y}\sin\varphi\sin\theta)+u_{z}r\cos\varphi, \\
\end{split}
\end{equation}
\begin{equation}\label{eq:chap_2_12_3}
\begin{split}
v_{r} &= v_{x}\cos\varphi\cos\theta+v_{y}\cos\varphi\sin\theta+v_{z}\sin\varphi, \\
v_{\theta} &= r(-v_{x}\cos\varphi\sin\theta+v_{y}\cos\varphi\cos\theta), \\
v_{\varphi} &= -r(v_{x}\sin\varphi\cos\theta+v_{y}\sin\varphi\sin\theta)+v_{z}r\cos\varphi, \\
\end{split}
\end{equation}
\begin{equation}\label{eq:chap_2_12_3w}
w_{r}=w_{z}\sin\varphi,\texttt{ and }w_{\varphi}=w_{z}r\cos\varphi.
\end{equation}

If the partial derivatives with respect to $x$, $y$, and $z$ also satisfy the Cauchy-Riemann equations
\begin{equation}\label{eq:chap_2_12_4}
e_{xy}u_{x}=e_{xy}v_{y}=e_{xy}u_{z}+w_{z},\texttt{ }
e_{xy}v_{x}=e_{xy}v_{z}=-e_{xy}u_{y}
\end{equation}
at $s_{0}$, then from equations~(\ref{eq:chap_2_12_2}) and~(\ref{eq:chap_2_12_3}) we get
\begin{equation}\label{eq:chap_2_12_5}
\begin{split}
      u_{r} &= u_{x}\cos\varphi\cos\theta+u_{y}\cos\varphi\sin\theta+u_{z}\sin\varphi, \\
 u_{\theta} &= r(-u_{x}\cos\varphi\sin\theta+u_{y}\cos\varphi\cos\theta), \\
u_{\varphi} &= -r(u_{x}\sin\varphi\cos\theta+u_{y}\sin\varphi\sin\theta)+u_{z}r\cos\varphi, \\
\end{split}
\end{equation}
\begin{equation}\label{eq:chap_2_12_6}
\begin{split}
      v_{r} &= -u_{y}\cos\varphi\cos\theta+u_{x}\cos\varphi\sin\theta+v_{z}\sin\varphi, \\
 v_{\theta} &= r(u_{y}\cos\varphi\sin\theta+u_{x}\cos\varphi\cos\theta), \\
v_{\varphi} &= r(u_{y}\sin\varphi\cos\theta-u_{x}\sin\varphi\sin\theta)+v_{z}r\cos\varphi, \\
\end{split}
\end{equation}
at that point. It is then clear from equations~(\ref{eq:chap_2_12_5}) and~(\ref{eq:chap_2_12_6}) that
\[
ru_{r}= v_{\theta}+u_{z}r\sin\varphi\texttt{ and }
rv_{r}=-u_{\theta}+v_{z}r\sin\varphi
\]
at the point $s_{0}$.

Because $r=(x^{2}+y^{2}+z^{2})^{1/2}$, $\theta=\arctan(y/x)$, and $\varphi=\arcsin(z/r)$, there are
\[
\begin{split}
u_{x}&=u_{r}\frac{\partial r}{\partial x}+u_{\theta}\frac{\partial\theta}{\partial x}+u_{\varphi}\frac{\partial\varphi}{\partial x}=u_{r}\frac{x}{r}-u_{\theta}\frac{y}{x^{2}+y^{2}}-u_{\varphi}\frac{xz/r^{3}}{[1-(z/r)^{2}]^{1/2}} \\
&=u_{r}\cos\varphi\cos\theta-\frac{u_{\theta}}{r}\frac{\sin\theta}{\cos\varphi}-\frac{u_{\varphi}}{r}\sin\varphi\cos\theta, \\
\end{split}
\]
\[
\begin{split}
u_{y}&=u_{r}\frac{\partial r}{\partial y}+u_{\theta}\frac{\partial\theta}{\partial y}+u_{\varphi}\frac{\partial\varphi}{\partial y}=u_{r}\frac{y}{r}-u_{\theta}\frac{x}{x^{2}+y^{2}}-u_{\varphi}\frac{yz/r^{3}}{[1-(z/r)^{2}]^{1/2}} \\
&=u_{r}\cos\varphi\sin\theta-\frac{u_{\theta}}{r}\frac{\cos\theta}{\cos\varphi}-\frac{u_{\varphi}}{r}\sin\varphi\sin\theta, \\
\end{split}
\]
\[
\begin{split}
u_{z}&=u_{r}\frac{\partial r}{\partial z}+u_{\varphi}\frac{\partial\varphi}{\partial z}=u_{r}\frac{z}{r}+\frac{u_{\varphi}\cos^{2}\varphi}{r[1-(z/r)^{2}]^{1/2}} \\
&=u_{r}\sin\varphi+\frac{u_{\varphi}}{r}\cos\varphi.
\end{split}
\]

Similarly, there are
\[
\begin{split}
v_{x}&=v_{r}\cos\varphi\cos\theta-\frac{v_{\theta}}{r}\frac{\sin\theta}{\cos\varphi}-\frac{v_{\varphi}}{r}\sin\varphi\cos\theta,\\ v_{y}&=v_{r}\cos\varphi\sin\theta-\frac{v_{\theta}}{r}\frac{\cos\theta}{\cos\varphi}-\frac{v_{\varphi}}{r}\sin\varphi\sin\theta,\\ v_{z}&=v_{r}\sin\varphi+\frac{v_{\varphi}}{r}\cos\varphi.
\end{split}
\]

Thus we get
\begin{equation}\label{eq:chap_2_12_7}
ru_{r}\cos^{2}\varphi= v_{\theta}+u_{\varphi}\sin\varphi\cos\varphi\texttt{ and }
rv_{r}\cos^{2}\varphi=-u_{\theta}+v_{\varphi}\sin\varphi\cos\varphi
\end{equation}
at the point $s_{0}$.

If, on the other hand, equations~(\ref{eq:chap_2_12_7}) are known to hold at $s_{0}$, it is straightforward
to know that equations~(\ref{eq:chap_2_12_4}) must hold there. Equations~(\ref{eq:chap_2_12_7}) are, therefore,
an alternative form of the Cauchy-Riemann equations~(\ref{eq:chap_2_12_4}).

We can now restate the theorem in Sec.~(\ref{sec:chap_2_11_21}) using spherical coordinates.

\begin{theorem}\label{th:chap_2_12_l}
Let the function
\[
f(s)=e_{xy}[u(r,\theta,\varphi)+iv(r,\theta,\varphi)]+w(r,\varphi)
\]
be defined throughout some $\varepsilon$ neighborhood of a nonzero point
\[
s_{0}=r_{0}[e_{xy}\exp(i\theta_{0})\cos\varphi_{0}+\sin\varphi_{0}],
\]
and suppose that the first-order partial derivatives of the functions $u$, $v$, and $w$ with respect to $r$, $\theta$, and $\varphi$ exist everywhere in that neighborhood. If those partial derivatives are continuous at $(r_{0},\theta_{0},\varphi_{0})$ and satisfy the spherical form
\[
ru_{r}\cos^{2}\varphi= v_{\theta}+u_{\varphi}\sin\varphi\cos\varphi\texttt{ and }
rv_{r}\cos^{2}\varphi=-u_{\theta}+v_{\varphi}\sin\varphi\cos\varphi
\]
of the Cauchy-Riemann equations at $(r_{0},\theta_{0},\varphi_{0})$, then $f'(s_{0})$ exists.
\end{theorem}

The derivative $f'(s_{0})$ here can be written
\begin{equation}\label{eq:chap_2_12_8}
f'(s_{0})=f_{r}(s_{0})=f_{\varphi}(s_{0})=f_{\theta}(s_{0})
\end{equation}
where
\[
\begin{split}
&f_{r}(s_{0})=\frac{e_{xy}(u_{r}+iv_{r})+w_{r}}{e^{i\theta}\cos\varphi+\sin\varphi}, \\
&f_{\varphi}(s_{0})=\frac{e_{xy}(u_{\varphi}+iv_{\varphi})+w_{\varphi}}{r(-e^{i\theta}\sin\varphi+\cos\varphi)},\\
&f_{\theta}(s_{0})=\frac{e_{xy}(v_{\theta}-iu_{\theta})}{re^{i\theta}\cos\varphi},
\end{split}
\]
and the right-hand sides are to be evaluated at $(r_{0},\theta_{0},\varphi_{0})$.

From equations~(\ref{eq:chap_2_10_9}) in Sec.~(\ref{sec:chap_2_10_20}), there are
\[
\begin{split}
f_{r}(s_{0}) &= \frac{e_{xy}(u_{r}+iv_{r})+w_{r}}{e^{i\theta}\cos\varphi+\sin\varphi} \\
&=\frac{e_{xy}(u_{x}+iv_{x})e^{i\theta}\cos\varphi}{e^{i\theta}\cos\varphi+\sin\varphi}
+\frac{[e_{xy}(u_{z}+iv_{z})+w_{z}]\sin\varphi}{e^{i\theta}\cos\varphi+\sin\varphi} \\
&=\frac{f'(s_{0})e^{i\theta}\cos\varphi}{e^{i\theta}\cos\varphi+\sin\varphi}
+\frac{f'(s_{0})\sin\varphi}{e^{i\theta}\cos\varphi+\sin\varphi}
=f'(s_{0}),
\end{split}
\]
\[
\begin{split}
f_{\varphi}(s_{0}) &=\frac{e_{xy}(u_{\varphi}+iv_{\varphi})+w_{\varphi}}{r(-e^{i\theta}\sin\varphi+\cos\varphi)}\\
&=\frac{-e_{xy}(u_{x}+iv_{x})e^{i\theta}\sin\varphi}{-e^{i\theta}\sin\varphi+\cos\varphi}
+\frac{[e_{xy}(u_{z}+iv_{z})+w_{z}]\cos\varphi}{-e^{i\theta}\sin\varphi+\cos\varphi} \\
&=\frac{-f'(s_{0})e^{i\theta}\sin\varphi}{-e^{i\theta}\sin\varphi+\cos\varphi}
+\frac{f'(s_{0})\cos\varphi}{-e^{i\theta}\sin\varphi+\cos\varphi}
=f'(s_{0}),
\end{split}
\]
and
\[
\begin{split}
f_{\theta}(s_{0}) &= \frac{e_{xy}(v_{\theta}-iu_{\theta})}{re^{i\theta}\cos\varphi} \\
                  &= e_{xy}(v_{y}-iu_{y})=f'(s_{0}).
\end{split}
\]
So Theorem~(\ref{th:chap_2_12_l}) holds and the proof of Theorem~(\ref{th:chap_2_12_l}) is completed.

\section{Analytic Functions}\label{sec:chap_2_13_23}

We are now ready to introduce the concept of an analytic function. A function $f$ of the complex variable $s$ is  analytic in an open set if it has a derivative at each point in that set. If we should speak of a function $f$ that is analytic in a set $S$ which is not open, it is to be understood that $f $ is analytic in an open set containing $S$. In particular, $f$ is analytic at a point $s_{0}$ if it is analytic throughout some neighborhood of $s_{0}$.

We note, for instance, that the function $f(s)=1/s$ is analytic at each nonzero point in the finite space. But the function $f(s)=|s|^{2}$ is not analytic at any point since its derivative exists only at $s=0$ and not throughout any neighborhood.

An entire function is a function that is analytic at each point in the entire finite space. Since the derivative  of a polynomial exists everywhere, it follows that every polynomial is an entire function.

If a function $f$ fails to be analytic at a point $s_{0}$ but is analytic at some point in every neighborhood of $s_{0}$, then $s_{0}$ is called a singular point, or singularity, of $f$.

The point $s=0$ is evidently a singular point of the function $f(s)=1/s$. The function $f(z)=|s|^{2}$, on the other hand, has no singular points since it is nowhere analytic.

A necessary, but by no means sufficient, condition for a function $f$ to be analytic in a domain $D$ is clearly the continuity of $f$ throughout $D$. Satisfaction of the Cauchy-Riemann equations is also necessary, but not sufficient. Sufficient conditions for analyticity in $D$ are provided by the theorems in Secs.~(\ref{sec:chap_2_11_21}) and~(\ref{sec:chap_2_12_22}).

Other useful sufficient conditions are obtained from the differentiation formulas
in Sec.~(\ref{sec:chap_2_9_19}). The derivatives of the sum and product of two functions exist wherever the functions themselves have derivatives. Thus, if two functions are analytic in a domain $D$, their sum and their product are both analytic in $D$. Similarly, their quotient is analytic in $D$ provided the function in the denominator does not vanish at any point in $D$. In particular, the quotient $P(s)/Q(s)$ of two polynomials is analytic in any domain throughout which $Q(s)\neq0$.

From the chain rule for the derivative of a composite function, we find that a composition of two analytic functions is analytic. More precisely, suppose that a function $f(s)$ is analytic in a domain $D$ and that the image (Sec.~(\ref{sec:chap_2_2_12})) of $D$ under the transformation $\varpi=f(s)$ is contained in the domain of definition of a function $g(\varpi)$.

Then the composition $g[f(s)]$ is analytic in $D$, with derivative
\[
\frac{d}{ds}g[f(s)]=g'[f(s)]f'(s).
\]

The following theorem is especially useful, in addition to being expected.

\begin{theorem}\label{th:chap_2_13_l}
If $f'(s)=0$ everywhere in a domain $D$, then $f(s)$ must be constant throughout $D$.
\end{theorem}

We start the proof by writing $f(s)=e_{xy}[u(x,y,z)+iv(x,y,z)]+w(z)$. Assuming that $f'(s)=0$ in $D$, then we note that $e_{xy}(u_{x}+iv_{x})=0$; and, in view of the Cauchy-Riemann equations, $v_{y}-iu_{y}=0$, and from equations~(\ref{eq:chap_2_11_6}), Sec.~(\ref{sec:chap_2_11_21}) $e_{xy}(u_{z}+iv_{z})+w_{z}=0$. Consequently,
\[
u_{x}=u_{y}=u_{z}=v_{x}=v_{y}=v_{z}=w_{z}=0
\]
at each point in $D$.

Next, we show that $u(x,y,z)$ is constant along any line segment $L$ extending from a point $P$ to a point $P'$ and lying entirely in $D$. We let $\varsigma$ denote the distance along $L$ from the point $P$ and let $U$ denote the unit vector along $L$ in the direction of increasing $\varsigma$. We know from calculus that the directional derivative $du/d\varsigma$ can be written as the dot product
\begin{equation}\label{eq:chap_2_13_1}
\frac{du}{d\varsigma}=(grad\texttt{ }u)\cdot U,
\end{equation}
where grad $u$ is the gradient vector
\begin{equation}\label{eq:chap_2_13_2}
grad\texttt{ }u=u_{x}\texttt{i}+u_{y}\texttt{j}+u_{z}\texttt{k},
\end{equation}

Because $u_{x}$, $u_{y}$ and $u_{z}$ are zero everywhere in $D$, then, grad  $u$ is the zero vector at all
points on $L$. Hence it fo1lows from equation~(\ref{eq:chap_2_13_1}) that the derivative $du/d\varsigma$ is zero along $L$; and this means that $u$ is constant on $L$.

Finally, since there is always a finite number of such line segments, joined end to end, connecting any two points  $P$ and $Q$ in $D$ (Sec.~(\ref{sec:chap_1_10})), the values of $u$ at $P$ and $Q$ must be the same. We may conclude, then, that there is a real constant $a$ such that $u(x,y,z)=a$ throughout $D$. Similarly, $v(x,y,z)=b$ and $w(z)=c$; and we find that $f(s)=e_{xy}(a+bi)+c$ at each point in $D$.

\begin{theorem}\label{th:chap_2_13_2}
If a spatial complex function $f(s)$ is analytic everywhere in a domain $D$, then its first-order partial derivatives are orthogonal in the domain $D$, i.e., for $f(s)=e_{xy}[u(x,y,z)+iv(x,y,z)]+w(z)$ there is
\begin{equation}\label{eq:chap_2_13_3}
u_{x}(x,y,z)v_{x}(x,y,z)+u_{y}(x,y,z)v_{y}(x,y,z)=0.
\end{equation}
\end{theorem}

\begin{proof}
When a spatial complex function $f(s)$ determined by Definition~(\ref{eq:chap_2_1_1}) is analytic, by the Cauchy-Riemann equations~(\ref{eq:chap_2_10_7}), Sec.~(\ref{sec:chap_2_10_20}), there are
\[
u_{x}(x,y,z)=v_{y}(x,y,z)\texttt{ and }
v_{x}(x,y,z)=-u_{y}(x,y,z).
\]
Thus Equation~(\ref{eq:chap_2_13_3}) holds and the proof of the theorem is completed.
\end{proof}

\begin{theorem}\label{th:chap_2_13_3}
Suppose that a function $f$ is analytic throughout an open sphere $|s|<R$ with radius $R$. Then $f(s)$ has the expression
\begin{equation}\label{eq:chap_2_13_4}
\begin{split}
f(s) &= e_{xy}[f(x+iy+z)-f(z)]+f(z) \\
     &= e_{xy}[u(x,y,z)+iv(x,y,z)]+w(z).
\end{split}
\end{equation}
where
\[
\begin{split}
u(x,y,z) &= Re\texttt{ }f(x+iy+z)-f(z), \\
v(x,y,z) &= Im\texttt{ }f(x+iy+z),\texttt{ and }w(z)=f(z).
\end{split}
\]
\end{theorem}

\begin{proof}
Let $\sigma=x+iy$. Then $\sigma$ and $\sigma+z$ are two-dimensional complex numbers, and $s=e_{xy}\sigma+z$.

Because $f$ is analytic, from the Taylor series theorems in Calculus and the two-dimensional complex variables theory, the real function $f(z)$ has Taylor series
\[
f(z)=\sum_{n=0}^{\infty}\frac{f^{(n)}(0)}{n!}z^{n}
\]
and the plane complex function $f(\sigma+z)$ is also analytic and has Taylor series
\[
f(\sigma+z)=\sum_{n=0}^{\infty}\frac{f^{(n)}(0)}{n!}(\sigma+z)^{n}.
\]

For each index $n$ in the Taylor series of $f(z)$ and $f(\sigma+z)$, there are $z^{n}$ and
\[
(\sigma+z)^{n}=\sigma^{n}+n\sigma^{n-1}z+\cdots+n\sigma z^{n-1}+z^{n}.
\]

Since both $\sigma+z$ and its powers are points in the two-dimensional complex plane and from remark~(\ref{re:chap_1_1_1}) for a term $e_{xy}x^{j}y^{k}$, the operator $e_{xy}$ becomes $e_{xy}^{m}$ where $m=|j|+|k|$, we get
\[
\begin{split}
e_{xy}[(\sigma+z)^{n}-z^{n}]+z^{n}
&= (e_{xy}\sigma)^{n}+n(e_{xy}\sigma)^{n-1}z+\cdots+n(e_{xy}\sigma)z^{n-1}+z^{n} \\
&= (e_{xy}\sigma+z)^{n}=s^{n}\texttt{ for }n=1,2,\ldots.
\end{split}
\]

Thus for $|s|<\infty$ and $z\neq0$ there are $|\sigma+z|<\infty$ and
\[
\begin{split}
&e_{xy}[f(\sigma+z)-f(z)]+f(z) \\
=&e_{xy}[\sum_{n=0}^{\infty}\frac{f^{(n)}(0)}{n!}(\sigma+z)^{n}-\sum_{n=0}^{\infty}\frac{f^{(n)}(0)}{n!}z^{n}]
+\sum_{n=0}^{\infty}\frac{f^{(n)}(0)}{n!}z^{n} \\
=&e_{xy}\sum_{n=1}^{\infty}\frac{f^{(n)}(0)}{n!}[(\sigma+z)^{n}-z^{n}]+\sum_{n=0}^{\infty}\frac{f^{(n)}(0)}{n!}z^{n} \\
=&\sum_{n=0}^{\infty}\frac{f^{(n)}(0)}{n!}(e_{xy}\sigma+z)^{n}
=\sum_{n=0}^{\infty}\frac{f^{(n)}(0)}{n!}s^{n}
=f(s).
\end{split}
\]
So the Taylor series of the spatial complex analytic function $f(s)$ exits and is
\begin{equation}\label{eq:chap_2_13_5}
f(s)=\sum_{n=0}^{\infty}\frac{f^{(n)}(0)}{n!}s^{n}\texttt{ }(|s|<\infty).
\end{equation}

Hence expression~(\ref{eq:chap_2_13_4}) holds and the proof of the theorem is completed.
\end{proof}

\section{Cauchy-Riemann Equation of $f(x+iy+z)$}\label{sec:chap_2_14_24}

In this section, we obtain a pair of equations that the first-order partial derivatives of the component functions  $u+w$ and $v$ of a spatial complex function
\begin{equation}\label{eq:chap_2_14_1}
\begin{split}
f(s) &= e_{xy}[f(x+iy+z)-f(z)]+f(z) \\
     &= e_{xy}[u(x,y,z)+iv(x,y,z)]+w(z)
\end{split}
\end{equation}
must satisfy at a point $s_{0}=(x_{0},y_{0},z_{0})$ when the derivative of $f$ exists there.
%We also show how to express $f'(s_{0})$ in terms of those partial derivatives.

Let $\sigma=x+iy+z$ and
\[
f(\sigma)=f(x+iy+z)=\phi(x,y,z)+i\psi(x,y,z).
\]
Then there are
\[
\phi(x,y,z)=u(x,y,z)+w(z)\texttt{ and }\psi(x,y,z)=v(x,y,z).
\]

We start by writing
\[
\sigma_{0}=x_{0}+z_{0}+iy_{0},\texttt{ }
\Delta\sigma=\Delta x+\Delta z+i\Delta y,
\]
and
\[
\begin{split}
\Delta\omega &= f(\sigma_{0}+\Delta\sigma)-f(\sigma_{0}) \\
             &= \phi(x_{0}+\Delta x,y_{0}+\Delta y,z_{0}+\Delta z)-\phi(x_{0},y_{0},z_{0}) \\
             &+ i[\psi(x_{0}+\Delta x,y_{0}+\Delta y,z_{0}+\Delta z)-\psi(x_{0},y_{0},z_{0})]=\Delta\phi+i\Delta\psi.
\end{split}
\]
Thus we have
\[
\frac{\Delta\omega}{\Delta\sigma}=\frac{\Delta\phi+i\Delta\psi}{\Delta x+\Delta z+i\Delta y}.
\]

Assuming that the derivative
\begin{equation}\label{eq:chap_2_14_2}
f'(\sigma_{0})=\lim_{\Delta\sigma \to 0}\frac{\Delta\omega}{\Delta\sigma}
\end{equation}
exists, we know from two-dimensional complex theory that
\begin{equation}\label{eq:chap_2_14_3}
f'(\sigma_{0})=\lim_{(\Delta x+\Delta z,\Delta y) \to (0,0)}Re\frac{\Delta\omega}{\Delta\sigma}
+i\lim_{(\Delta x+\Delta z,\Delta y) \to (0,0)}Im\frac{\Delta\omega}{\Delta\sigma}.
\end{equation}

Now it is important to keep in mind that $\Delta\sigma$ and $\Delta\omega/\Delta\sigma$ are in the $e_{xy}$ coordinate plane or in the two-dimensional complex plane, so that expression~(\ref{eq:chap_2_14_3}) is valid as $(\Delta x+\Delta z,\Delta y)$ tends to $(0,0)$ in any manner that we may choose just as in the case of the two-dimensional complex plane. In particular, we let $(\Delta x+\Delta z,\Delta y)$ tend to $(0,0)$ through the point $(\Delta x+\Delta z,0)$ along the $x$ axis. Inasmuch as $\Delta y=0$,  the quotient $\Delta\omega/\Delta\sigma$ becomes
\[
\begin{split}
\frac{\Delta\omega}{\Delta\sigma}&=\frac{\phi(x_{0}+\Delta x,y_{0},z_{0}+\Delta z)-\phi(x_{0},y_{0},z_{0})}{\Delta x+\Delta z} \\
&+i\frac{\psi(x_{0}+\Delta x,y_{0},z_{0}+\Delta z)-\psi(x_{0},y_{0},z_{0})}{\Delta x+\Delta z}.
\end{split}
\]
Thus
\[
\begin{split}
\lim_{(\Delta x+\Delta z,\Delta y) \to (0,0)}Re\frac{\Delta\omega}{\Delta\sigma}
&=\lim_{\Delta x+\Delta z \to 0}\frac{\phi(x_{0}+\Delta x,y_{0},z_{0}+\Delta z)-\phi(x_{0},y_{0},z_{0})}{\Delta x+\Delta z} \\
&=\phi_{x+z}(x_{0},y_{0},z_{0})
\end{split}
\]
and
\[
\begin{split}
\lim_{(\Delta x+\Delta z,\Delta y) \to (0,0)}Im\frac{\Delta\omega}{\Delta\sigma}
&=\lim_{\Delta x+\Delta z \to 0}\frac{\psi(x_{0}+\Delta x,y_{0},z_{0}+\Delta z)-\psi(x_{0},y_{0},z_{0})}{\Delta x+\Delta z} \\
&=\psi_{x+z}(x_{0},y_{0},z_{0})
\end{split}
\]
where $\phi_{x+z}(x_{0},y_{0},z_{0})$ and $\psi_{x+z}(x_{0},y_{0},z_{0})$ denote the first-order partial derivatives with respect to the sum of $x$ and $z$ of the functions $\phi$ and $\psi$ respectively, at $(x_{0},y_{0},z_{0})$. Furthermore, if $\Delta z=0$ or $\Delta\sigma=\Delta x$ then there are
\[
\begin{split}
\lim_{(\Delta x+\Delta z,\Delta y) \to (0,0)}Re\frac{\Delta\omega}{\Delta\sigma}
&=\lim_{\Delta x \to 0}\frac{\phi(x_{0}+\Delta x,y_{0},z_{0})-\phi(x_{0},y_{0},z_{0})}{\Delta x} \\
&=\phi_{x}(x_{0},y_{0},z_{0})
\end{split}
\]
and
\[
\begin{split}
\lim_{(\Delta x+\Delta z,\Delta y) \to (0,0)}Im\frac{\Delta\omega}{\Delta\sigma}
&=\lim_{\Delta x \to 0}\frac{\psi(x_{0}+\Delta x,y_{0},z_{0})-\psi(x_{0},y_{0},z_{0})}{\Delta x+\Delta z} \\
&=\psi_{x}(x_{0},y_{0},z_{0})
\end{split}
\]
where $\phi_{x}(x_{0},y_{0},z_{0})$ and $\psi_{x}(x_{0},y_{0},z_{0})$ denote the first-order partial derivatives with respect to $x$ of the functions $\phi$ and $\psi$ respectively, at $(x_{0},y_{0},z_{0})$; or if $\Delta x=0$ or $\Delta\sigma=\Delta z$ then there are
\[
\begin{split}
\lim_{(\Delta x+\Delta z,\Delta y) \to (0,0)}Re\frac{\Delta\omega}{\Delta\sigma}
&=\lim_{\Delta z \to 0}\frac{\phi(x_{0},y_{0},z_{0}+\Delta z)-\phi(x_{0},y_{0},z_{0})}{\Delta z} \\
&=\phi_{z}(x_{0},y_{0},z_{0})
\end{split}
\]
and
\[
\begin{split}
\lim_{(\Delta x+\Delta z,\Delta y) \to (0,0)}Im\frac{\Delta\omega}{\Delta\sigma}
&=\lim_{\Delta z \to 0}\frac{\psi(x_{0},y_{0},z_{0}+\Delta z)-\psi(x_{0},y_{0},z_{0})}{\Delta z} \\
&=\psi_{z}(x_{0},y_{0},z_{0})
\end{split}
\]
where $\phi_{z}(x_{0},y_{0},z_{0})$ and $\psi_{z}(x_{0},y_{0},z_{0})$ denote the first-order partial derivatives with respect to $z$ of the functions $\phi$ and $\psi$ respectively, at $(x_{0},y_{0},z_{0})$. Substitution of these limits into expression~(\ref{eq:chap_2_14_3}) tells us that
\begin{equation}\label{eq:chap_2_14_4}
\begin{split}
\frac{\partial}{\partial\sigma}f(\sigma_{0}) &= \phi_{x+z}(x,y,z)+i\psi_{x+z}(x,y,z) \\
                                             &= \phi_{x}(x,y,z)+i\psi_{x}(x,y,z) \\
                                             &= \phi_{z}(x,y,z)+i\psi_{z}(x,y,z).
\end{split}
\end{equation}

We might have let $(\Delta x+\Delta z,\Delta y)$ tend to $(0,0)$ through the point $(0,\Delta y)$ along the $y$ axis. In that case, $\Delta x+\Delta z=0$ and
\[
\begin{split}
\frac{\Delta\omega}{\Delta\sigma}&=\frac{\phi(x_{0},y_{0}+\Delta y,z_{0})-\phi(x_{0},y_{0},z_{0})}{i\Delta y}
+i\frac{\psi(x_{0},y_{0}+\Delta y,z_{0})-\psi(x_{0},y_{0},z_{0})}{i\Delta y} \\
&=\frac{\psi(x_{0},y_{0}+\Delta y,z_{0})-\psi(x_{0},y_{0},z_{0})}{\Delta y}
-i\frac{\phi(x_{0},y_{0}+\Delta y,z_{0})-\phi(x_{0},y_{0},z_{0})}{\Delta y}.
\end{split}
\]
Evidently, then,
\[
\begin{split}
\lim_{(\Delta x+\Delta z,\Delta y) \to (0,0)}Re\frac{\Delta\omega}{\Delta\sigma}
&=\lim_{\Delta y \to 0}\frac{\psi(x_{0},y_{0}+\Delta y,z_{0})-\psi(x_{0},y_{0},z_{0})}{\Delta y} \\
&=\psi_{y}(x_{0},y_{0},z_{0})
\end{split}
\]
and
\[
\begin{split}
\lim_{(\Delta x,\Delta y) \to (0,0)}Im\frac{\Delta\omega}{\Delta\sigma}
&=-\lim_{\Delta y \to 0}\frac{\phi(x_{0},y_{0}+\Delta y,z_{0})-\phi(x_{0},y_{0},z_{0})}{\Delta y} \\
&=-\phi_{y}(x_{0},y_{0},z_{0}).
\end{split}
\]
Hence it follows from expression~(\ref{eq:chap_2_14_3}) that
\begin{equation}\label{eq:chap_2_14_5}
\frac{\partial}{\partial\sigma}f(\sigma_{0})=\psi_{y}(x_{0},y_{0},z_{0})-i\phi_{y}(x_{0},y_{0},z_{0})
\end{equation}
where the partial derivatives of the spatial functions $\phi$ and $\psi$ are, this time, with respect to $y$. Note that equation~(\ref{eq:chap_2_14_5}) can also be written in the form
\[
\frac{\partial}{\partial\sigma}f(\sigma_{0})=-i[\phi_{y}(x_{0},y_{0},z_{0})+i\psi_{y}(x_{0},y_{0},z_{0})].
\]

Equations~(\ref{eq:chap_2_14_4}) and~(\ref{eq:chap_2_14_5}) not only give $\frac{\partial}{\partial\sigma}f(\sigma_{0})$ in terms of partial derivatives of the component functions $\phi$ and $\psi$, but they also provide necessary conditions for the existence of $f'(\sigma_{0})$. For, on equating the real and imaginary parts on the right-hand sides of these equations, we see that the existence of $f'(\sigma_{0})$ requires that
\begin{equation}\label{eq:chap_2_14_6}
\phi_{x+z}=\phi_{x}=\phi_{z}=\psi_{y},\texttt{ }\psi_{x+z}=\psi_{x}=\psi_{z}=-\phi_{y}.
\end{equation}
Equations~(\ref{eq:chap_2_14_6}) are the Cauchy-Riemann equations of a spatial complex function $f(x+iy+z)$.

We summarize the above results as follows.

\begin{theorem}\label{th:chap_2_14_l}
Suppose that the derivative of a spatial complex function
\[
\begin{split}
f(s) &= e_{xy}[f(x+iy+z)-f(z)]+f(z) \\
     &= e_{xy}[u(x,y,z)+iv(x,y,z)]+w(z)
\end{split}
\]
exists at a point $s_{0}=(x_{0},y_{0},z_{0})$, so that also the derivative of the plane complex function of a plane complex variable $\sigma=x+z+iy$ composed of three real variables $x$, $y$, and $z$
\[
f(\sigma)=f(x+iy+z)=\phi(x,y,z)+i\psi(x,y,z)
\]
exists at a point $\sigma_{0}=x_{0}+z_{0}+iy_{0}$ in $e_{xy}$ complex plane where
\[
\phi(x,y,z)=u(x,y,z)+w(z)\texttt{ and }\psi(x,y,z)=v(x,y,z).
\]
Then the first-order partial derivatives of $\phi$ and $\psi$ must exist at $(x_{0}+z_{0},y_{0})$, and they must satisfy the Cauchy-Riemann equations of a two-dimensional complex function $f(x+iy+z)$
\begin{equation}\label{eq:chap_2_14_7}
\phi_{x+z}=\phi_{x}=\phi_{z}=\psi_{y},\texttt{ }\psi_{x+z}=\psi_{x}=\psi_{z}=-\phi_{y}
\end{equation}
there. Also, $f'(\sigma_{0})$ can be written
\begin{equation}\label{eq:chap_2_14_8}
f'(\sigma_{0})=f_{x+z}(\sigma_{0})=f_{x}(\sigma_{0})=f_{z}(\sigma_{0})=f_{y}(\sigma_{0})
\end{equation}
where
\[
\begin{split}
&f_{x+z}(\sigma)=\phi_{x+z}+i\psi_{x+z}, \\
&f_{x}(\sigma)=\phi_{x}+i\psi_{x}, \\
&f_{z}(\sigma)=\phi_{z}+i\psi_{z}, \\
&f_{y}(\sigma)=\psi_{y}-i\phi_{y},
\end{split}
\]
and these partial derivatives are to be evaluated at $(x_{0},y_{0},z_{0})$.
\end{theorem}

\section{Harmonic Functions}\label{sec:chap_2_15_25}

A real-valued function $H$ of three real variables $x$, $y$, and $z$ is said to be harmonic in a given domain of  the $xyz$ space if throughout that domain, it has continuous partial derivatives of the first and second order and satisfies the partial differential equation
\begin{equation}\label{eq:chap_2_15_0}
H_{xx}(x,y,z)+H_{yy}(x,y,z)+H_{zz}(x,y,z)=0,
\end{equation}
known as Laplace's equation.

\begin{theorem}\label{th:chap_2_15_1}
If a spatial complex function
\[
f(s)=e_{xy}[u(x,y,z)+iv(x,y,z)]+w(z)
\]
is analytic in a domain $D$, then its complex component functions $u$ and $v$ are two-dimensional harmonic in $D$, i.e., there are
\begin{equation}\label{eq:chap_2_15_1}
\Delta u(x,y,z)=0\texttt{ and }\Delta v(x,y,z)=0\texttt{ for }\Delta=\frac{\partial^{2}}{(\partial x)^{2}}+\frac{\partial^{2}}{(\partial y)^{2}}
\end{equation}
with boundary conditions in the $xyz$ coordinate space.
\end{theorem}

To prove this, we need a result that is proved in the two-dimensional complex variable theory. Namely, if a function of a complex variable is analytic at a point, then its components have continuous partial derivatives of all orders at that point.

Assuming that $f$ is analytic in $D$, we start with the observation that the first-order partial derivatives of its complex component functions $u$ and $v$ must satisfy the Cauchy-Riemann equations throughout $D$:
\begin{equation}\label{eq:chap_2_15_2}
u_{x}=v_{y}\texttt{, }u_{y}=-v_{x}.
\end{equation}
Differentiating both sides of these equations with respect to $x$, we have
\begin{equation}\label{eq:chap_2_15_3}
u_{xx}=v_{yx}\texttt{, }u_{yy}=-v_{xy}.
\end{equation}
Likewise, differentiation with respect to $y$ yields
\begin{equation}\label{eq:chap_2_15_4}
u_{xy}=v_{yy}\texttt{, }u_{yx}=-v_{xx}.
\end{equation}
Now, by a theorem in advanced calculus, the continuity of the partial derivatives of $u$ and $v$ ensures that $u_{yx}=u_{xy}$ and $v_{yx}=v_{xy}$. It then follows from equations~(\ref{eq:chap_2_15_3}) and~(\ref{eq:chap_2_15_4}) that
\[
u_{xx}+u_{yy}=0\texttt{ and }v_{xx}+v_{yy}=0.
\]
That is, $u$ and $v$ are two-dimensional harmonic in $D$.

Thus Theorem~(\ref{th:chap_2_15_1}) holds and the proof of the theorem is completed.

If two given functions $u$ and $v$ are harmonic in a domain $D$ and their first-order partial derivatives satisfy the Cauchy-Riemann equations~(\ref{eq:chap_2_15_2}) throughout $D$, $v$ is said to be a harmonic conjugate of $u$. The meaning of the word conjugate here is, of course, different from that in Sec.~(\ref{sec:chap_1_5}), where $\bar{s}$ is defined.

\begin{theorem}\label{th:chap_2_15_1b}
A spatial complex function
\[
f(s)=e_{xy}[u(x,y,z)+iv(x,y,z)]+w(z)
\]
is analytic in a domain $D$ if and only if $v$ is a harmonic conjugate of $u$.
\end{theorem}

The proof is easy. If $v$ is a harmonic conjugate of $u$ in $D$, the theorem in Sec.~(\ref{sec:chap_2_11_21}) tells us that $f$ is analytic in $D$. Conversely, if $f$ is analytic in $D$, we know from Theorem~(\ref{th:chap_2_15_1}) above that $u$ and $v$ are harmonic in $D$; and, in view of the theorem in Sec.~(\ref{sec:chap_2_10_20}) the Cauchy-Riemann equations are satisfied.

It is evident that if $v$ is a harmonic conjugate of $u$ in some domain, it is not, in general, true that $u$ is a harmonic conjugate of $v$ there.

In Chap.~(\ref{ch:chap_9}) (Sec.~(\ref{sec:chap_9_4_97})) we shall show that a function $u$ which is harmonic in a
domain of a certain type always has a harmonic conjugate. Thus, in such domains, every harmonic function is the real part of an analytic function. It is also true that a harmonic conjugate, when it exists, is unique except for an additive constant.

\begin{theorem}\label{th:chap_2_15_1s}
If a spatial complex function
\[
f(s)=e_{xy}[u(r,\theta,\varphi)+iv(r,\theta,\varphi)]+w(r,\varphi)
\]
is analytic in a domain $D$, then its complex component functions $u$ and $v$ satisfy
\[
r\frac{\partial}{\partial r}(ru_{r}\cos^{2}\varphi)+\frac{u_{\theta\theta}}{\cos^{2}\varphi}
+\tan\varphi\frac{\partial}{\partial\varphi}(u_{\varphi}\sin\varphi\cos\varphi)
-2\sin\varphi\frac{\partial}{\partial\varphi}(ru_{r}\cos\varphi)=0.
\]
and
\[
r\frac{\partial}{\partial r}(rv_{r}\cos^{2}\varphi)+\frac{v_{\theta\theta}}{\cos^{2}\varphi}
+\tan\varphi\frac{\partial}{\partial\varphi}(v_{\varphi}\sin\varphi\cos\varphi)
-2\sin\varphi\frac{\partial}{\partial\varphi}(rv_{r}\cos\varphi)=0.
\]
with boundary conditions in the $xyz$ coordinate space.
\end{theorem}

To prove this, we need a result that is proved in the two-dimensional complex variable theory. Namely, if a function of a complex variable is analytic at a point, then its components have continuous partial derivatives of all orders at that point.

According to the Cauchy-Riemann equations in Theorem~(\ref{th:chap_2_12_l})
\[
ru_{r}\cos^{2}\varphi= v_{\theta}+u_{\varphi}\sin\varphi\cos\varphi\texttt{ and }
rv_{r}\cos^{2}\varphi=-u_{\theta}+v_{\varphi}\sin\varphi\cos\varphi
\]
there are
\[
\frac{\partial}{\partial r}(ru_{r}\cos^{2}\varphi)=v_{r\theta}+u_{r\varphi}\sin\varphi\cos\varphi\texttt{ and }
rv_{r\theta}\cos^{2}\varphi=-u_{\theta\theta}+v_{\theta\varphi}\sin\varphi\cos\varphi.
\]
That is
\[
r\cos^{2}\varphi\frac{\partial}{\partial r}(ru_{r}\cos^{2}\varphi)+u_{\theta\theta}=(ru_{r\varphi}\cos^{2}\varphi+v_{\theta\varphi})\sin\varphi\cos\varphi.
\]

Because $v_{\theta}=ru_{r}\cos^{2}\varphi-u_{\varphi}\sin\varphi\cos\varphi$, there are
\[
v_{\theta\varphi}=\frac{\partial}{\partial\varphi}(ru_{r}\cos^{2}\varphi-u_{\varphi}\sin\varphi\cos\varphi)
\]
and
\[
r\frac{\partial}{\partial r}(ru_{r}\cos^{2}\varphi)+\frac{u_{\theta\theta}}{\cos^{2}\varphi}
+\tan\varphi\frac{\partial}{\partial\varphi}(u_{\varphi}\sin\varphi\cos\varphi)
-2\sin\varphi\frac{\partial}{\partial\varphi}(ru_{r}\cos\varphi)=0.
\]

Likewise, there are
\[
\frac{\partial}{\partial r}(rv_{r}\cos^{2}\varphi)=-u_{r\theta}+v_{r\varphi}\sin\varphi\cos\varphi\texttt{ and }
ru_{r\theta}\cos^{2}\varphi=v_{\theta\theta}+u_{\theta\varphi}\sin\varphi\cos\varphi.
\]
That is
\[
r\cos^{2}\varphi\frac{\partial}{\partial r}(rv_{r}\cos^{2}\varphi)+v_{\theta\theta}=(rv_{r\varphi}\cos^{2}\varphi-u_{\theta\varphi})\sin\varphi\cos\varphi.
\]

Because $u_{\theta}=v_{\varphi}\sin\varphi\cos\varphi-rv_{r}\cos^{2}\varphi$, there are
\[
u_{\theta\varphi}=\frac{\partial}{\partial\varphi}(v_{\varphi}\sin\varphi\cos\varphi-rv_{r}\cos^{2}\varphi)
\]
and
\[
r\frac{\partial}{\partial r}(rv_{r}\cos^{2}\varphi)+\frac{v_{\theta\theta}}{\cos^{2}\varphi}
+\tan\varphi\frac{\partial}{\partial\varphi}(v_{\varphi}\sin\varphi\cos\varphi)
-2\sin\varphi\frac{\partial}{\partial\varphi}(rv_{r}\cos\varphi)=0.
\]

Hence Theorem~(\ref{th:chap_2_15_1s}) holds and the proof of the theorem is completed.

\begin{theorem}\label{th:chap_2_15_2}
If a spatial complex function
\[
\begin{split}
f(s) &= e_{xy}[f(x+iy+z)-f(z)]+f(z) \\
     &= e_{xy}[u(x,y,z)+iv(x,y,z)]+w(z)
\end{split}
\]
is analytic in a domain $D$ where $f(z)=w(z)$,
\[
\begin{split}
f(x+iy+z) &= \phi(x,y,z)+i\psi(x,y,z) \\
          &= u(x,y,z)+w(z)+iv(x,y,z),
\end{split}
\]
then $\phi(x,y,z)$ and $\psi(x,y,z)$ are two-dimensional harmonic in $D$, i.e., there are
\begin{equation}\label{eq:chap_2_15_1b}
\Delta\phi(x,y,z)=0\texttt{ and }\Delta\psi(x,y,z)=0\texttt{ for }\Delta=\frac{\partial^{2}}{(\partial x+\partial z)^{2}}+\frac{\partial^{2}}{(\partial y)^{2}}
\end{equation}
with boundary conditions in the $xyz$ coordinate space.
\end{theorem}

To prove this, we need a result that is proved in the two-dimensional complex variable theory. Namely, if a function of a complex variable is analytic at a point, then its components have continuous partial derivatives of all orders at that point.

Assuming that $f(s)$ is analytic in $D$ so that $f(x+iy+z)$ is also analytic in $D$, we start with Theorem~(\ref{th:chap_2_14_l}) that the first-order partial derivatives of its complex component functions $u$, $v$, and $w$ must satisfy the Cauchy-Riemann equations throughout $D$ for $(u+w+iv)$:
\begin{equation}\label{eq:chap_2_15_2b}
\phi_{x+z}=\psi_{y}\texttt{, }\phi_{y}=-\psi_{x+z}.
\end{equation}
Differentiating both sides of these equations with respect to $x+z$, we have
\begin{equation}\label{eq:chap_2_15_3b}
\phi_{(x+z)(x+z)}=\psi_{y(x+z)}\texttt{, }\phi_{y(x+z)}=-\psi_{(x+z)(x+z)}.
\end{equation}
Likewise, differentiation with respect to $y$ yields
\begin{equation}\label{eq:chap_2_15_4b}
\phi_{(x+z)y}=\psi_{yy}\texttt{, }\phi_{yy}=-\psi_{(x+z)y}.
\end{equation}
Now, by a theorem in advanced calculus, the continuity of the partial derivatives of $\phi$ and $\psi$ ensures that $\phi_{y(x+z)}=\phi_{(x+z)y}$ and $\psi_{y(x+z)}=\psi_{(x+z)y}$. It then follows from equations~(\ref{eq:chap_2_15_3b}) and~(\ref{eq:chap_2_15_4b}) that
\[
\phi_{(x+z)(x+z)}+\phi_{yy}=0\texttt{ and }\psi_{(x+z)(x+z)}+\psi_{yy}=0.
\]
That is, $\phi$ and $\psi$ are two-dimensional harmonic in $D$.

Thus Theorem~(\ref{th:chap_2_15_2}) holds and the proof of the theorem is completed.

\begin{theorem}\label{th:chap_2_15_3}
If a spatial complex function $f(s)$ is analytic in a domain $D$ where
\[
\begin{split}
f(s) &= e_{xy}[f(x+iy+z)-f(z)]+f(z) \\
     &= e_{xy}[u(x,y,z)+iv(x,y,z)]+w(z)
\end{split}
\]
and
\[
\begin{split}
f(x+iy+z) &= \phi(x,y,z)+i\psi(x,y,z), \\
u(x,y,z) &= Re\texttt{ }f(x+iy+z)-f(z), \\
v(x,y,z) &= Im\texttt{ }f(x+iy+z)\texttt{ and }w(z)=f(z),
\end{split}
\]
let $\alpha$, $\beta$, and $\gamma$ be positive reals and satisfy $\alpha^{2}+\gamma^{2}=\beta^{2}$, and let $x$, $y$, and $z$ be replaced by $\alpha x$, $\beta y$, and $\gamma z$ in the function $f(x+iy+z)$, respectively, that is,
\[
f(\alpha x+i\beta y+\gamma z)=\phi(\alpha x, \beta y, \gamma z)+i\psi(\alpha x, \beta y, \gamma z),
\]
then $\phi$ and $\psi$ are three-dimensional harmonic in $D$, i.e., there are
\begin{equation}\label{eq:chap_2_15_1c}
\Delta\phi=0\texttt{ and }\Delta\psi=0\texttt{ for }\Delta=\frac{\partial^{2}}{(\partial x)^{2}}+\frac{\partial^{2}}{(\partial y)^{2}}+\frac{\partial^{2}}{(\partial z)^{2}}
\end{equation}
with boundary conditions in the $xyz$ coordinate space.
\end{theorem}

To prove this, we need a result that is proved in the two-dimensional complex variable theory. Namely, if a function of a complex variable is analytic at a point, then its components have continuous partial derivatives of all orders at that point.

Assuming that $f(s)$ is analytic in $D$ so that $f(x+iy+z)$ is also analytic in $D$, we start with Theorem~(\ref{th:chap_2_14_l}) in Sec.~(\ref{sec:chap_2_14_24}) that the first-order and second-order partial derivatives of its complex component functions $\phi$ and $\psi$ exist throughout $D$.

For convenience, let the first-order and second-order partial derivatives with respect to $x$, $y$, and $z$ of the function $\phi(\alpha x, \beta y, \gamma z)$ be denoted by $\phi_{x}$, $\phi_{y}$, $\phi_{z}$, $\phi_{xx}$, $\phi_{yy}$, $\phi_{xx}$, and so on, respectively, and let the first-order and second-order partial derivatives with respect to $\alpha x$, $\beta y$, and $\gamma z$ of the function $\phi(\alpha x, \beta y, \gamma z)$ be denoted by $\phi_{\bar{x}}$, $\phi_{\bar{y}}$, $\phi_{\bar{z}}$, $\phi_{\overline{xx}}$, $\phi_{\overline{yy}}$, $\phi_{\overline{zz}}$, and so on, respectively. Then the relations between them are
\begin{equation}\label{eq:chap_2_15_2c}
\begin{array}{cccc}
\phi_{x}=\alpha\phi_{\alpha x}=\alpha\phi_{\bar{x}}, & \phi_{xx}=\alpha^{2}\phi_{\overline{xx}}, & \phi_{xy}=\alpha\beta\phi_{\overline{xy}}, & \phi_{xz}=\alpha\gamma\phi_{\overline{xz}}, \\
\phi_{y}=\beta\phi_{\beta y}=\beta\phi_{\bar{y}}, & \phi_{yy}=\beta^{2}\phi_{\overline{yy}}, & \phi_{yz}=\beta\gamma\phi_{\overline{yz}}, & \phi_{yx}=\alpha\beta\phi_{\overline{yx}}, \\
\phi_{z}=\gamma\phi_{\gamma z}=\gamma\phi_{\bar{z}}, & \phi_{zz}=\gamma^{2}\phi_{\overline{zz}}, & \phi_{zx}=\alpha\gamma\phi_{\overline{zx}}, & \phi_{zy}=\beta\gamma\phi_{\overline{zy}}.
\end{array}
\end{equation}
So does the function $\psi(\alpha x, \beta y, \gamma z)$.

First, from Cauchy-Riemann equations~(\ref{eq:chap_2_14_7}) in Sec.~(\ref{sec:chap_2_14_24}), we have
\[
\phi_{\overline{xx}}+\phi_{\overline{yy}}=0\texttt{ and }\phi_{\overline{zz}}+\phi_{\overline{yy}}=0.
\]
Because
\[
\phi_{\overline{xx}}=\phi_{xx}/\alpha^{2},\texttt{ }\phi_{\overline{yy}}=\phi_{yy}/\beta^{2},\texttt{ and }\phi_{\overline{zz}}=\phi_{zz}/\gamma^{2},
\]
there are
\[
\beta^{2}\phi_{xx}+\alpha^{2}\phi_{yy}=0\texttt{ and }\beta^{2}\phi_{zz}+\gamma^{2}\phi_{yy}=0.
\]
Because $\alpha^{2}+\gamma^{2}=\beta^{2}$, summing two equations above, we can derive
\[
\beta^{2}(\phi_{xx}+\phi_{zz})+(\alpha^{2}+\gamma^{2})\phi_{yy}
=\beta^{2}(\phi_{xx}+\phi_{zz}+\phi_{yy})=0.
\]
So does the function $\psi(\alpha x, \beta y, \gamma z)$. Thus we get
\[
\phi_{xx}+\phi_{yy}+\phi_{zz}=\Delta\phi=0\texttt{ and }\psi_{xx}+\psi_{yy}+\psi_{zz}=\Delta\psi=0.
\]

Next, let $\sigma=\alpha x+i\beta y+\gamma z$. Then there are $f(\alpha x+i\beta y+\gamma z)=f(\sigma)$, and
\[
\begin{split}
& \frac{\partial^{2}f(\sigma)}{(\partial x)^{2}}=\frac{\partial^{2}f(\sigma)}{(\partial\sigma)^{2}}\frac{(\partial\sigma)^{2}}{(\partial x)^{2}}=\frac{\partial^{2}f(\sigma)}{(\partial\sigma)^{2}}\alpha^{2}, \\
& \frac{\partial^{2}f(\sigma)}{(\partial y)^{2}}=\frac{\partial^{2}f(\sigma)}{(\partial\sigma)^{2}}\frac{(\partial\sigma)^{2}}{(\partial y)^{2}}=\frac{\partial^{2}f(\sigma)}{(\partial\sigma)^{2}}(i\beta)^{2}, \\
& \frac{\partial^{2}f(\sigma)}{(\partial z)^{2}}=\frac{\partial^{2}f(\sigma)}{(\partial\sigma)^{2}}\frac{(\partial\sigma)^{2}}{(\partial z)^{2}}=\frac{\partial^{2}f(\sigma)}{(\partial\sigma)^{2}}\gamma^{2}.
\end{split}
\]
Thus we get
\[
\begin{split}
\Delta f(\sigma) &=\frac{\partial^{2}f(\sigma)}{(\partial x)^{2}}+\frac{\partial^{2}f(\sigma)}{(\partial y)^{2}}+\frac{\partial^{2}f(\sigma)}{(\partial z)^{2}} \\
                 &=\alpha^{2}\frac{\partial^{2}f(\sigma)}{(\partial\sigma)^{2}}-\beta^{2}\frac{\partial^{2}f(\sigma)}{(\partial\sigma)^{2}}
+\gamma^{2}\frac{\partial^{2}f(\sigma)}{(\partial\sigma)^{2}} \\
                 &=(\alpha^{2}-\beta^{2}+\gamma^{2})\frac{\partial^{2}f(\sigma)}{(\partial\sigma)^{2}}=0.
\end{split}
\]
Because there are
\[
\Delta f(\sigma)=\Delta(\phi+i\psi)=\Delta\phi+i\Delta\psi,
\]
there must be
\[
\Delta\phi=0\texttt{ and }\Delta\psi=0.
\]

Hence, $\phi$ and $\psi$ are three-dimensional harmonic in $D$ and Theorem~(\ref{th:chap_2_15_3}) holds. The proof of the theorem is completed.

\begin{theorem}\label{th:chap_2_15_3s}
If a spatial complex function $f(s)$ is analytic in a domain $D$ where
\[
\begin{split}
f(s) &= e_{xy}[f(x+iy+z)-f(z)]+f(z) \\
     &= e_{xy}[u(x,y,z)+iv(x,y,z)]+w(z)
\end{split}
\]
and
\[
\begin{split}
f(x+iy+z) &= \phi(x,y,z)+i\psi(x,y,z), \\
u(x,y,z) &= Re\texttt{ }f(x+iy+z)-f(z), \\
v(x,y,z) &= Im\texttt{ }f(x+iy+z)\texttt{ and }w(z)=f(z),
\end{split}
\]
let $\alpha$, $\beta$, and $\gamma$ be positive reals and satisfy $\alpha^{2}+\gamma^{2}=\beta^{2}$, and let $x$, $y$, and $z$ be replaced by $\alpha x$, $\beta y$, and $\gamma z$ in the function $f(x+iy+z)$, respectively, that is,
\[
f(\alpha x+i\beta y+\gamma z)=\phi(\alpha x, \beta y, \gamma z)+i\psi(\alpha x, \beta y, \gamma z),
\]
then $\phi$ and $\psi$ are three-dimensional harmonic in $D$ in the spherical coordinates, i.e., there are
\[
\frac{\partial(r^{2}\phi_{r})}{r^{2}\partial r}+\frac{\partial\phi_{\theta}}{r^{2}\cos^{2}\varphi\partial\theta}+\frac{\partial(\phi_{\varphi}\cos\varphi)}{r^{2}\cos\varphi\partial\varphi}=0
\]
and
\[
\frac{\partial(r^{2}\psi_{r})}{r^{2}\partial r}+\frac{\partial\psi_{\theta}}{r^{2}\cos^{2}\varphi\partial\theta}+\frac{\partial(\psi_{\varphi}\cos\varphi)}{r^{2}\cos\varphi\partial\varphi}=0
\]
with boundary conditions in the spherical coordinate space.
\end{theorem}

From Theorem~(\ref{th:chap_2_15_3}) we get
\[
\Delta\phi=0\texttt{ and }\Delta\psi=0\texttt{ for }\Delta=\frac{\partial^{2}}{(\partial x)^{2}}+\frac{\partial^{2}}{(\partial y)^{2}}+\frac{\partial^{2}}{(\partial z)^{2}}.
\]
Then, in the spherical coordinates, these Laplace's equations above become
\[
\frac{\partial(r^{2}\phi_{r})}{r^{2}\partial r}+\frac{\partial\phi_{\theta}}{r^{2}\cos^{2}\varphi\partial\theta}+\frac{\partial(\phi_{\varphi}\cos\varphi)}{r^{2}\cos\varphi\partial\varphi}=0
\]
and
\[
\frac{\partial(r^{2}\psi_{r})}{r^{2}\partial r}+\frac{\partial\psi_{\theta}}{r^{2}\cos^{2}\varphi\partial\theta}+\frac{\partial(\psi_{\varphi}\cos\varphi)}{r^{2}\cos\varphi\partial\varphi}=0
\]
where because
\[
r=(x^{2}+y^{2}+z^{2})^{1/2},\texttt{ }\theta=\arctan(y/x),\texttt{ and }\varphi=\arctan(z/\sqrt{x^{2}+y^{2}}),
\]
there are
\[
\begin{split}
r_{x}&=\frac{x}{r},\texttt{ }r_{y}=\frac{y}{r},\texttt{ }r_{z}=\frac{z}{r},\\
\theta_{x}&=-\frac{y}{x^{2}+y^{2}},\texttt{ }
\theta_{y}=\frac{x}{x^{2}+y^{2}},\texttt{ }\theta_{z}=0,\\
\varphi_{x}&=-\frac{xz}{r^{2}(x^{2}+y^{2})^{1/2}},\texttt{ }
\varphi_{y} =-\frac{yz}{r^{2}(x^{2}+y^{2})^{1/2}},\texttt{ }
\varphi_{z} =\frac{\sqrt{x^{2}+y^{2}}}{r^{2}}.
\end{split}
\]
\[
\begin{split}
r_{xx}&=\frac{y^{2}+z^{2}}{r^{3}},\texttt{ }
r_{yy}=\frac{x^{2}+z^{2}}{r^{3}},\texttt{ }
r_{zz}=\frac{x^{2}+y^{2}}{r^{3}},\\
\theta_{xx}&=-\frac{2xy}{(x^{2}+y^{2})^{2}},\texttt{ }
\theta_{yy}=\frac{2xy}{(x^{2}+y^{2})^{2}},\texttt{ }
\theta_{zz}=0,\\
\varphi_{xx}&=-\frac{z}{r^{2}(x^{2}+y^{2})^{1/2}}+2\frac{x^{2}z}{r^{4}(x^{2}+y^{2})^{1/2}}+\frac{x^{2}z}{r^{2}(x^{2}+y^{2})^{3/2}},\\
\varphi_{yy}&=-\frac{z}{r^{2}(x^{2}+y^{2})^{1/2}}+2\frac{y^{2}z}{r^{4}(x^{2}+y^{2})^{1/2}}+\frac{y^{2}z}{r^{2}(x^{2}+y^{2})^{3/2}},\\
\varphi_{zz}&=-2\frac{z\sqrt{x^{2}+y^{2}}}{r^{4}},
\end{split}
\]
and
\[
\begin{split}
\phi_{x}&=\phi_{r}r_{x}+\phi_{\theta}\theta_{x}+\phi_{\varphi}\varphi_{x},\\ \phi_{y}&=\phi_{r}r_{y}+\phi_{\theta}\theta_{y}+\phi_{\varphi}\varphi_{y},\\
\phi_{z}&=\phi_{r}r_{z}+\phi_{\theta}\theta_{z}+\phi_{\varphi}\varphi_{z}.
\end{split}
\]
\[
\begin{split}
\phi_{xx}&=\phi_{r}r_{xx}+\phi_{rr}r_{x}^{2}+\phi_{r\theta}r_{x}\theta_{x}+\phi_{r\varphi}r_{x}\varphi_{x} \\
&+\phi_{\theta}\theta_{xx}+\phi_{r\theta}r_{x}\theta_{x}+\phi_{\theta\theta}\theta_{x}^{2}+\phi_{\theta\varphi}\theta_{x}\varphi_{x}\\
&+\phi_{\varphi}\varphi_{xx}+\phi_{r\varphi}r_{x}\varphi_{x}+\phi_{\theta\varphi}\theta_{x}\varphi_{x}+\phi_{\varphi\varphi}\varphi_{x}^{2},\\
\phi_{yy}&=\phi_{r}r_{yy}+\phi_{rr}r_{y}^{2}+\phi_{r\theta}r_{y}\theta_{y}+\phi_{r\varphi}r_{y}\varphi_{y} \\
&+\phi_{\theta}\theta_{yy}+\phi_{r\theta}r_{y}\theta_{y}+\phi_{\theta\theta}\theta_{y}^{2}+\phi_{\theta\varphi}\theta_{y}\varphi_{y}\\
&+\phi_{\varphi}\varphi_{yy}+\phi_{r\varphi}r_{y}\varphi_{y}+\phi_{\theta\varphi}\theta_{y}\varphi_{y}+\phi_{\varphi\varphi}\varphi_{y}^{2},\\
\phi_{zz}&=\phi_{r}r_{zz}+\phi_{rr}r_{z}^{2}+\phi_{r\theta}r_{z}\theta_{z}+\phi_{r\varphi}r_{z}\varphi_{z} \\
&+\phi_{\theta}\theta_{zz}+\phi_{r\theta}r_{z}\theta_{z}+\phi_{\theta\theta}\theta_{z}^{2}+\phi_{\theta\varphi}\theta_{z}\varphi_{z}\\
&+\phi_{\varphi}\varphi_{zz}+\phi_{r\varphi}r_{z}\varphi_{z}+\phi_{\theta\varphi}\theta_{z}\varphi_{z}+\phi_{\varphi\varphi}\varphi_{z}^{2},\\
\end{split}
\]
so that there are
\[
\phi_{xx}+\phi_{yy}+\phi_{zz}=\frac{\partial(r^{2}\phi_{r})}{r^{2}\partial r}+\frac{\partial\phi_{\theta}}{r^{2}\cos^{2}\varphi\partial\theta}+\frac{\partial(\phi_{\varphi}\cos\varphi)}{r^{2}\cos\varphi\partial\varphi}
\]
and similarly
\[
\psi_{xx}+\psi_{yy}+\psi_{zz}=\frac{\partial(r^{2}\psi_{r})}{r^{2}\partial r}+\frac{\partial\psi_{\theta}}{r^{2}\cos^{2}\varphi\partial\theta}+\frac{\partial(\psi_{\varphi}\cos\varphi)}{r^{2}\cos\varphi\partial\varphi}.
\]

Hence, $\phi$ and $\psi$ are three-dimensional harmonic in $D$ in the spherical coordinates, and Theorem~(\ref{th:chap_2_15_3s}) holds. The proof of the theorem is completed.

\begin{theorem}\label{th:chap_2_15_4}
If a spatial complex function $f(s)$ is analytic in a domain $D$ where
\[
\begin{split}
f(s) &= e_{xy}[f(x+iy+z)-f(z)]+f(z) \\
     &= e_{xy}[u(x,y,z)+iv(x,y,z)]+w(z)
\end{split}
\]
and
\[
\begin{split}
f(x+iy+z) &= \phi(x,y,z)+i\psi(x,y,z), \\
u(x,y,z) &= Re\texttt{ }f(x+iy+z)-f(z), \\
v(x,y,z) &= Im\texttt{ }f(x+iy+z)\texttt{ and }w(z)=f(z),
\end{split}
\]
let $\alpha$, $\beta$, and $\gamma$ be positive reals and satisfy $\alpha^{2}+\gamma^{2}=\beta^{2}$, and let $x$, $y$, and $z$ be replaced by $\alpha x$, $\beta y$, and $\gamma z$ in the function $f(x+iy+z)$, respectively, that is,
\[
f(\alpha x+i\beta y+\gamma z)=\phi(\alpha x, \beta y, \gamma z)+i\psi(\alpha x, \beta y, \gamma z),
\]
then the first-order partial derivatives of $\phi$ and $\psi$ satisfy the three-dimensional orthogonal equation
\begin{equation}\label{eq:chap_2_15_1d}
\phi_{x}\psi_{x}+\phi_{y}\psi_{y}+\phi_{z}\psi_{z}=0.
\end{equation}
\end{theorem}

To prove this, we use the result of Theorem~(\ref{th:chap_2_13_l}) in Sec.~(\ref{sec:chap_2_13_23}), and expressions~(\ref{eq:chap_2_15_2c}) in Sec.~(\ref{sec:chap_2_15_25}), we have
\[
\phi_{\bar{x}}\psi_{\bar{x}}+\phi_{\bar{y}}\psi_{\bar{y}}=0\texttt{ and }
\phi_{\bar{z}}\psi_{\bar{z}}+\phi_{\bar{y}}\psi_{\bar{y}}=0.
\]
Because
\[
\phi_{\bar{x}}\psi_{\bar{x}}=\phi_{x}\psi_{x}/\alpha^{2},\texttt{ }
\phi_{\bar{y}}\psi_{\bar{y}}=\phi_{y}\psi_{y}/\beta^{2},\texttt{ and }
\phi_{\bar{z}}\psi_{\bar{z}}=\phi_{z}\psi_{z}/\gamma^{^{2}},
\]
there are
\[
\beta^{2}\phi_{x}\psi_{x}+\alpha^{2}\phi_{y}\psi_{y}=0\texttt{ and }
\beta^{2}\phi_{z}\psi_{z}+\gamma^{2}\phi_{y}\psi_{y}=0.
\]
Because $\alpha^{2}+\gamma^{2}=\beta^{2}$, summing two equations above, we can derive
\[
\beta^{2}(\phi_{x}\psi_{x}+\phi_{z}\psi_{z})+(\alpha^{2}+\gamma^{2})\phi_{y}\psi_{y}
=\beta^{2}(\phi_{x}\psi_{x}+\phi_{z}\psi_{z}+\phi_{y}\psi_{y})=0.
\]
Thus we get
\[
\phi_{x}\psi_{x}+\phi_{y}\psi_{y}+\phi_{z}\psi_{z}=0.
\]

Hence, $\phi$ and $\psi$ are of the three-dimensional orthogonal property and Theorem~(\ref{th:chap_2_15_4}) holds. The proof of the theorem is completed.

\begin{theorem}\label{th:chap_2_15_5}
If a spatial complex function
\[
f(s)=e_{xy}[u(x,y,z)+iv(x,y,z)]+w(z)
\]
is analytic in a domain $D$ and there exist two real functions $p(x,y,z)$ and $q(x,y,z)$ which satisfy the Cauchy-Riemann equations throughout $D$, that is
\[
p_{x}(x,y,z)=q_{y}(x,y,z)\texttt{, }p_{y}(x,y,z)=-q_{x}(x,y,z),
\]
and there are
\[
\frac{\partial^{2}}{(\partial z)^{2}}(u(x,y,z)+p(x,y,z))=0\texttt{ and }\frac{\partial^{2}}{(\partial z)^{2}}(v(x,y,z)+q(x,y,z))=0,
\]
then let $\phi(x,y,z)=u(x,y,z)+p(x,y,z)$ and $\psi(x,y,z)=v(x,y,z)+q(x,y,z)$. Thus the real functions $\phi(x,y,z)$ and $\psi(x,y,z)$ satisfy the Cauchy-Riemann equations throughout $D$, that is
\begin{equation}\label{eq:chap_2_15_5}
\phi_{x}=\psi_{y}\texttt{, }\phi_{y}=-\psi_{x},
\end{equation}
and the real functions $\phi$ and $\psi$ are three-dimensional harmonic in $D$, that is
\begin{equation}\label{eq:chap_2_15_6}
\Delta\phi(x,y,z)=0\texttt{ and }\Delta\psi(x,y,z)=0\texttt{ for }\Delta=\frac{\partial^{2}}{(\partial x)^{2}}+\frac{\partial^{2}}{(\partial y)^{2}}+\frac{\partial^{2}}{(\partial z)^{2}}
\end{equation}
with boundary conditions in the $xyz$ coordinate space.
\end{theorem}

To prove this, we use the result of Theorem~(\ref{th:chap_2_15_1}). Assuming that $f$ is analytic in $D$, based on Theorem~(\ref{th:chap_2_15_1}), we get
\[
\Delta u(x,y,z)=0\texttt{ and }\Delta v(x,y,z)=0\texttt{ for }\Delta=\frac{\partial^{2}}{(\partial x)^{2}}+\frac{\partial^{2}}{(\partial y)^{2}}
\]
with boundary conditions in the $xyz$ coordinate space.

Let take a spatial complex polynomial function $p(s,2)$ and a two-dimensional complex polynomial function $p_{a}(z,2)$ from Sec.~(\ref{sec:chap_3_1_27}). Let $f(s)=p(s,2)$ and
\[
\begin{split}
\phi(x,y,z) &= Re[Im_{s}f(s)-p_{a}(z,2)]-Re_{s}p(s,2), \\
\psi(x,y,z) &= Im[Im_{s}f(s)-p_{a}(z,2)].
\end{split}
\]
Then lemma~(\ref{le:chap_2_15_1}) shows that the first-order partial derivatives of $\phi$ and $\psi$ satisfy the Cauchy-Riemann equations~(\ref{eq:chap_2_15_5}) and the second-order partial derivatives of $\phi$ and $\psi$ satisfy the Laplace's equations~(\ref{eq:chap_2_15_6}) with boundary conditions in the $xyz$ coordinate space.

Let take a spatial complex polynomial function $p(s,3)$ and two two-dimensional complex polynomial functions $p_{a}(z,3)$ and $p_{c}(z,3)$ from Sec.~(\ref{sec:chap_3_1_27}). Let $f(s)=p(s,3)$ and
\[
\begin{split}
\phi(x,y,z) &= Re[Im_{s}f(s)-(p_{a}(z,3)+p_{c}(z,3))]-Re_{s}p(s,3), \\
\psi(x,y,z) &= Im[Im_{s}f(s)-(p_{a}(z,3)+p_{c}(z,3))].
\end{split}
\]
Then lemma~(\ref{le:chap_2_15_2}) shows that the first-order partial derivatives of $\phi$ and $\psi$ satisfy the Cauchy-Riemann equations~(\ref{eq:chap_2_15_5}) and the second-order partial derivatives of $\phi$ and $\psi$ satisfy the Laplace's equations~(\ref{eq:chap_2_15_6}) with boundary conditions in the $xyz$ coordinate space.

Similarly, let take a spatial complex polynomial function $p(s,4)$ and two two-dimensional complex polynomial functions $p_{a}(z,4)$ and $p_{c}(z,4)$ from Sec.~(\ref{sec:chap_3_1_27}). Let $f(s)=p(s,4)$ and
\[
\begin{split}
\phi(x,y,z) &= Re[Im_{s}f(s)-(p_{a}(z,4)+p_{c}(z,4))]-Re_{s}p(s,4), \\
\psi(x,y,z) &= Im[Im_{s}f(s)-(p_{a}(z,4)+p_{c}(z,4))].
\end{split}
\]
Then lemma~(\ref{le:chap_2_15_3}) shows that the first-order partial derivatives of $\phi$ and $\psi$ satisfy the Cauchy-Riemann equations~(\ref{eq:chap_2_15_5}) and the second-order partial derivatives of $\phi$ and $\psi$ satisfy the Laplace's equations~(\ref{eq:chap_2_15_6}) with boundary conditions in the $xyz$ coordinate space.

Thus Theorem~(\ref{th:chap_2_15_5}) holds and the proof of the theorem is completed.

\begin{lemma}\label{le:chap_2_15_1}
Let take a spatial complex polynomial function $p(s,2)$ and a two-dimensional complex polynomial function $p_{a}(z,2)$ from Sec.~(\ref{sec:chap_3_1_27}). Let
\[
\begin{split}
\phi(x,y,z) &= Re(Im_{s}p(s,2)-p_{a}(z,2))-Re_{s}p(s,2), \\
\psi(x,y,z) &= Im(Im_{s}p(s,2)-p_{a}(z,2)).
\end{split}
\]
Then the first-order partial derivatives of $\phi$ and $\psi$ satisfy the Cauchy-Riemann equations~(\ref{eq:chap_2_15_5}) and the second-order partial derivatives of $\phi$ and $\psi$ satisfy the Laplace's equations~(\ref{eq:chap_2_15_6}) with boundary conditions in the $xyz$ coordinate space.
\end{lemma}

From Sec.~(\ref{sec:chap_3_1_27}), there are
\[
\begin{split}
    p(s,2) &= \sum_{n=1}^{2}a_{n}s^{n}=e_{xy}(Im_{s}p(s,2)-Re_{s}p(s,2))+Re_{s}p(s,2), \\
p_{a}(z,2) &= e_{xy}\sum_{n=1}^{2}(a_{n,x}+ia_{n,y})z^{n}
\end{split}
\]
where $a_{n}=e_{xy}(a_{n,x}+ia_{n,y})+a_{n,z}$ for $n=1,2$ are spatial complex constants and
\[
Re_{s}p(s,2)=\sum_{n=1}^{2}a_{n,z}z^{n},\texttt{ }
Im_{s}p(s,2)=\sum_{n=1}^{2}a_{n}(x+z+iy)^{n},
\]
and
\[
p_{a}(z,2)+Re_{s}p(s,2)=a_{1}z+a_{2}z^{2}.
\]

Let $a_{n,xz}$ denote $a_{n,x}+a_{n,z}$ for $n=1,2$. Then with definition $C_{n}^{k}=\frac{n!}{(n-k)!k!}$ there are
\[
\begin{split}
Im_{s}p(s,2)-Re_{s}p(s,2)-p_{a}(z,2) &= \sum_{n=1}^{2}a_{n}(x+z+iy)^{n}-\sum_{n=1}^{2}a_{n}z^{n} \\
                                     &= \sum_{n=1}^{2}a_{n}\sum_{k=0}^{n-1}C_{n}^{k}(x+iy)^{n-k}z^{k} \\
                                     &= a_{1}(x+iy)+a_{2}[(x+iy)^{2}+2(x+iy)z].
\end{split}
\]

Let $Im_{s}p(s,2)-Re_{s}p(s,2)-p_{a}(z,2)=\phi(x,y,z)+i\psi(x,y,z)$. Then we get
\begin{equation}\label{eq:chap_2_15_7}
\phi(x,y,z)=a_{1,xz}x-a_{1,y}y+a_{2,xz}(x^{2}+2xz-y^{2})-2a_{2,y}(x+z)y
\end{equation}
and
\begin{equation}\label{eq:chap_2_15_8}
\psi(x,y,z)=a_{1,xz}y+a_{1,y}x+2a_{2,xz}(x+z)y+a_{2,y}(x^{2}+2xz-y^{2}).
\end{equation}

First, we get the first-order partial derivatives of $\phi$ and $\psi$ which are
\begin{equation}\label{eq:chap_2_15_9}
\begin{split}
& \phi_{x}=a_{1,xz}+2a_{2,xz}(x+z)-2a_{2,y}y, \\
& \phi_{y}=-a_{1,y}-2a_{2,xz}y-2a_{2,y}(x+z), \\
& \phi_{z}=2a_{2,xz}x-2a_{2,y}y,
\end{split}
\end{equation}
and
\begin{equation}\label{eq:chap_2_15_10}
\begin{split}
& \psi_{x}=a_{1,y}+2a_{2,xz}y+2a_{2,y}(x+z), \\
& \psi_{y}=a_{1,xz}+2a_{2,xz}(x+z)-2a_{2,y}y, \\
& \psi_{z}=2a_{2,xz}y+2a_{2,y}x.
\end{split}
\end{equation}
These first-order partial derivatives of $\phi$ and $\psi$ satisfy the Cauchy-Riemann equations~(\ref{eq:chap_2_15_5})

Next, we get the second-order partial derivatives of $\phi$ and $\psi$ which are
\begin{equation}\label{eq:chap_2_15_11}
    \begin{array}{ccc}
        \phi_{xx}=2a_{2,xz}, & \phi_{xy}=-2a_{2,y}, & \phi_{xz}=2a_{2,xz}, \\
        \phi_{yx}=-2a_{2,y}, & \phi_{yy}=-2a_{2,xz}, & \phi_{yz}=-2a_{2,y}, \\
        \phi_{zx}=2a_{2,xz}, & \phi_{zy}=-2a_{2,y}, & \phi_{zz}=0,
    \end{array}
\end{equation}
and
\begin{equation}\label{eq:chap_2_15_12}
    \begin{array}{ccc}
        \psi_{xx}=2a_{2,y}, & \psi_{xy}=2a_{2,xz}, & \psi_{xz}=2a_{2,y}, \\
        \psi_{yx}=2a_{2,xz}, & \psi_{yy}=-2a_{2,y}, & \psi_{yz}=2a_{2,xz}, \\
        \psi_{zx}=2a_{2,y}, & \psi_{zy}=2a_{2,xz}, & \psi_{zz}=0.
    \end{array}
\end{equation}
These second-order partial derivatives of $\phi$ and $\psi$ satisfy the Laplace's equations~(\ref{eq:chap_2_15_6}) with boundary conditions in the $xyz$ coordinate space. The proof of lemma~(\ref{le:chap_2_15_1}) is completed.

\begin{lemma}\label{le:chap_2_15_2}
Let take a spatial complex polynomial function $p(s,3)$ and two two-dimensional complex polynomial functions $p_{a}(z,3)$ and $p_{c}(z,3)$ from Sec.~(\ref{sec:chap_3_1_27}). Let
\[
\begin{split}
\phi(x,y,z) &= Re[Im_{s}p(s,3)-(p_{a}(z,3)+p_{c}(z,3))]-Re_{s}p(s,3), \\
\psi(x,y,z) &= Im[Im_{s}p(s,3)-(p_{a}(z,3)+p_{c}(z,3))].
\end{split}
\]
Then the first-order partial derivatives of $\phi$ and $\psi$ satisfy the Cauchy-Riemann equations~(\ref{eq:chap_2_15_5}) and the second-order partial derivatives of $\phi$ and $\psi$ satisfy the Laplace's equations~(\ref{eq:chap_2_15_6}) with boundary conditions in the $xyz$ coordinate space.
\end{lemma}

From Sec.~(\ref{sec:chap_3_1_27}), there are
\[
p(s,3)=\sum_{n=1}^{3}a_{n}s^{n}=e_{xy}(Im_{s}p(s,3)-Re_{s}p(s,3))+Re_{s}p(s,3),
\]
\[
p_{a}(z,3)=e_{xy}\sum_{n=1}^{3}(a_{n,x}+ia_{n,y})z^{n},\texttt{ and }
p_{c}(z,3)=3e_{xy}(a_{3,xz}+ia_{3,y})(x+iy)z^{2}
\]
where $a_{n}=e_{xy}(a_{n,x}+ia_{n,y})+a_{n,z}$ for $n=1,2,3$ are spatial complex constants and
\[
Re_{s}p(s,3)=\sum_{n=1}^{3}a_{n,z}z^{n},\texttt{ }
Im_{s}p(s,3)=\sum_{n=1}^{3}a_{n}(x+z+iy)^{n},
\]
and
\[
p_{a}(z,3)+Re_{s}p(s,3)=a_{1}z+a_{2}z^{2}+a_{3}z^{3}.
\]

Let $a_{n,xz}$ denote $a_{n,x}+a_{n,z}$ for $n=1,2,3$. Then with definition $C_{n}^{k}=\frac{n!}{(n-k)!k!}$ there are
\[
Im_{s}p(s,3)-Re_{s}p(s,3)-(p_{a}(z,3)+p_{c}(z,3))
\]
\[
\begin{split}
&=\sum_{n=1}^{3}a_{n}(x+z+iy)^{n}
-\sum_{n=1}^{3}a_{n}z^{n}
-\sum_{n=3}^{3}a_{n}\sum_{k=2}^{n-1}C_{n}^{k}(x+iy)^{n-k}z^{k} \\
&=\sum_{n=1}^{3}a_{n}\sum_{k=0}^{n-1}C_{n}^{k}(x+iy)^{n-k}z^{k}
-\sum_{n=3}^{3}a_{n}\sum_{k=2}^{n-1}C_{n}^{k}(x+iy)^{n-k}z^{k} \\
&=\sum_{n=1}^{2}a_{n}\sum_{k=0}^{n-1}C_{n}^{k}(x+iy)^{n-k}z^{k}
+\sum_{n=3}^{3}a_{n}\sum_{k=0}^{1}C_{n}^{k}(x+iy)^{n-k}z^{k} \\
&=a_{1}(x+iy)+a_{2}[(x+iy)^{2}+2(x+iy)z]+a_{3}[(x+iy)^{3}+3(x+iy)^{2}z].
\end{split}
\]

Let $Im_{s}p(s,3)-Re_{s}p(s,3)-(p_{a}(z,3)+p_{c}(z,3))=\phi(x,y,z)+i\psi(x,y,z)$. Then we get
\begin{equation}\label{eq:chap_2_15_13}
\begin{split}
\phi(x,y,z) &= a_{1,xz}x-a_{1,y}y+a_{2,xz}(x^{2}+2xz-y^{2})-2a_{2,y}(x+z)y \\
            &+ a_{3,xz}[x^{3}-3xy^{2}+3(x^{2}-y^{2})z]-a_{3,y}[3(x^{2}+2xz)y-y^{3}],
\end{split}
\end{equation}
\begin{equation}\label{eq:chap_2_15_14}
\begin{split}
\psi(x,y,z) &= a_{1,xz}y+a_{1,y}x+2a_{2,xz}(x+z)y+a_{2,y}(x^{2}+2xz-y^{2}) \\
            &+ a_{3,y}[x^{3}-3xy^{2}+3(x^{2}-y^{2})z]+a_{3,xz}[3(x^{2}+2xz)y-y^{3}].
\end{split}
\end{equation}

First, we get the first-order partial derivatives of $\phi$ and $\psi$ which are
\begin{equation}\label{eq:chap_2_15_15}
\begin{split}
& \phi_{x}=a_{1,xz}+2a_{2,xz}(x+z)-2a_{2,y}y+3a_{3,xz}(x^{2}+2xz-y^{2}) \\
& \texttt{ }-6a_{3,y}(x+z)y, \\
& \phi_{y}=-a_{1,y}-2a_{2,xz}y-2a_{2,y}(x+z)-6a_{3,xz}(x+z)y \\
& \texttt{ }-3a_{3,y}(x^{2}+2xz-y^{2}), \\
& \phi_{z}=2a_{2,xz}x-2a_{2,y}y+3a_{3,xz}(x^{2}-y^{2})-6a_{3,y}xy,
\end{split}
\end{equation}
\begin{equation}\label{eq:chap_2_15_16}
\begin{split}
& \psi_{x}=a_{1,y}+2a_{2,xz}y+2a_{2,y}(x+z)+6a_{3,xz}(x+z)y \\
& \texttt{ }+3a_{3,y}(x^{2}+2xz-y^{2}), \\
& \psi_{y}=a_{1,xz}+2a_{2,xz}(x+z)-2a_{2,y}y+3a_{3,xz}(x^{2}+2xz-y^{2}) \\
& \texttt{ }-6a_{3,y}(x+z)y, \\
& \psi_{z}=2a_{2,xz}y+2a_{2,y}x+6a_{3,xz}xy+3a_{3,y}(x^{2}-y^{2}).
\end{split}
\end{equation}
These first-order partial derivatives of $\phi$ and $\psi$ satisfy the Cauchy-Riemann equations~(\ref{eq:chap_2_15_5})

Next, we get the second-order partial derivatives of $\phi$ and $\psi$ which are
\begin{equation}\label{eq:chap_2_15_17}
\begin{split}
& \phi_{xx}=2a_{2,xz}+6a_{3,xz}(x+z)-6a_{3,y}y, \\
& \phi_{xy}=-2a_{2,y}-6a_{3,xz}y-6a_{3,y}(x+z), \\
& \phi_{xz}=2a_{2,xz}+6a_{3,xz}x-6a_{3,y}y, \\
& \phi_{yx}=-2a_{2,y}-6a_{3,xz}y-6a_{3,y}(x+z), \\
& \phi_{yy}=-2a_{2,xz}-6a_{3,xz}(x+z)+6a_{3,y}y, \\
& \phi_{yz}=-2a_{2,y}-6a_{3,xz}y-6a_{3,y}x, \\
& \phi_{zx}=2a_{2,xz}+6a_{3,xz}x-6a_{3,y}y, \\
& \phi_{zy}=-2a_{2,y}-6a_{3,xz}y-6a_{3,y}x, \\
& \phi_{zz}=0,
\end{split}
\end{equation}
\begin{equation}\label{eq:chap_2_15_18}
\begin{split}
& \psi_{xx}=2a_{2,y}+6a_{3,xz}y+6a_{3,y}(x+z), \\
& \psi_{xy}=2a_{2,xz}+6a_{3,xz}(x+z)-6a_{3,y}y, \\
& \psi_{xz}=2a_{2,y}+6a_{3,xz}y+6a_{3,y}x, \\
& \psi_{yx}=2a_{2,xz}+6a_{3,xz}(x+z)-6a_{3,y}y, \\
& \psi_{yy}=-2a_{2,y}-6a_{3,xz}y-6a_{3,y}(x+z), \\
& \psi_{yz}=2a_{2,xz}+6a_{3,xz}x-6a_{3,y}y, \\
& \psi_{zx}=2a_{2,y}+6a_{3,xz}y+6a_{3,y}x, \\
& \psi_{zy}=2a_{2,xz}+6a_{3,xz}x-6a_{3,y}y, \\
& \psi_{zz}=0.
\end{split}
\end{equation}
These second-order partial derivatives of $\phi$ and $\psi$ satisfy the Laplace's equations~(\ref{eq:chap_2_15_6}) with boundary conditions in the $xyz$ coordinate space. The proof of lemma~(\ref{le:chap_2_15_2}) is completed.

\begin{lemma}\label{le:chap_2_15_3}
Let take a spatial complex polynomial function $p(s,4)$ and two two-dimensional complex polynomial functions $p_{a}(z,4)$ and $p_{c}(z,4)$ from Sec.~(\ref{sec:chap_3_1_27}). Let
\[
\begin{split}
\phi(x,y,z) &= Re[Im_{s}p(s,4)-(p_{a}(z,4)+p_{c}(z,4))]-Re_{s}p(s,4), \\
\psi(x,y,z) &= Im[Im_{s}p(s,4)-(p_{a}(z,4)+p_{c}(z,4))].
\end{split}
\]
Then the first-order partial derivatives of $\phi$ and $\psi$ satisfy the Cauchy-Riemann equations~(\ref{eq:chap_2_15_5}) and the second-order partial derivatives of $\phi$ and $\psi$ satisfy the Laplace's equations~(\ref{eq:chap_2_15_6}) with boundary conditions in the $xyz$ coordinate space.
\end{lemma}

From Sec.~(\ref{sec:chap_3_1_27}), there are
\[
\begin{split}
\texttt{ }p(s,4)&=\sum_{n=1}^{4}a_{n}s^{n}=e_{xy}(Im_{s}p(s,4)-Re_{s}p(s,4))+Re_{s}p(s,4),\\
p_{a}(z,4)&=e_{xy}\sum_{n=1}^{4}(a_{n,x}+ia_{n,y})z^{n},\\
p_{c}(z,4)&=e_{xy}\sum_{n=3}^{4}a_{n}\sum_{k=2}^{n-1}C_{n}^{k}(x+iy)^{n-k}z^{k}\texttt{ }(C_{n}^{k}=\frac{n!}{(n-k)!k!})
\end{split}
\]
where $a_{n}=e_{xy}(a_{n,x}+ia_{n,y})+a_{n,z}$ for $n=1,2,3,4$ are spatial complex constants and
\[
Re_{s}p(s,4)=\sum_{n=1}^{4}a_{n,z}z^{n},\texttt{ }
Im_{s}p(s,4)=\sum_{n=1}^{4}a_{n}(x+z+iy)^{n},
\]
and
\[
p_{a}(z,4)+Re_{s}p(s,4)=a_{1}z+a_{2}z^{2}+a_{3}z^{3}+a_{4}z^{4}.
\]

Let $a_{n,xz}$ denote $a_{n,x}+a_{n,z}$ for $n=1,2,3,4$. Then there are
\[
Im_{s}p(s,4)-Re_{s}p(s,4)-(p_{a}(z,4)+p_{c}(z,4))
\]
\[
\begin{split}
&=\sum_{n=1}^{4}a_{n}(x+z+iy)^{n}
-\sum_{n=1}^{4}a_{n}z^{n}
-\sum_{n=3}^{4}a_{n}\sum_{k=2}^{n-1}C_{n}^{k}(x+iy)^{n-k}z^{k} \\
&=\sum_{n=1}^{4}a_{n}\sum_{k=0}^{n-1}C_{n}^{k}(x+iy)^{n-k}z^{k}
-\sum_{n=3}^{4}a_{n}\sum_{k=2}^{n-1}C_{n}^{k}(x+iy)^{n-k}z^{k} \\
\end{split}
\]
\[
\begin{split}
&=\sum_{n=1}^{2}a_{n}\sum_{k=0}^{n-1}C_{n}^{k}(x+iy)^{n-k}z^{k}
+\sum_{n=3}^{4}a_{n}\sum_{k=0}^{1}C_{n}^{k}(x+iy)^{n-k}z^{k} \\
&=a_{1}(x+iy)+a_{2}[(x+iy)^{2}+2(x+iy)z]+a_{3}[(x+iy)^{3}+3(x+iy)^{2}z] \\
&+a_{4}[(x+iy)^{4}+4(x+iy)^{3}z].
\end{split}
\]

Let $Im_{s}p(s,4)-Re_{s}p(s,4)-(p_{a}(z,4)+p_{c}(z,4))=\phi(x,y,z)+i\psi(x,y,z)$. Then we get
\begin{equation}\label{eq:chap_2_15_19}
\begin{split}
\phi(x,y,z) &= a_{1,xz}x-a_{1,y}y+a_{2,xz}(x^{2}+2xz-y^{2})-2a_{2,y}(x+z)y \\
            &+ a_{3,xz}[x^{3}-3xy^{2}+3(x^{2}-y^{2})z]-a_{3,y}(3x^{2}y-y^{3}+6xyz) \\
            &+ a_{4,xz}[x^{4}-6x^{2}y^{2}+y^{4}+4(x^{3}-3xy^{2})z] \\
            &- 4a_{4,y}[x^{3}y-xy^{3}+(3x^{2}y-y^{3})z]
\end{split}
\end{equation}
and
\begin{equation}\label{eq:chap_2_15_20}
\begin{split}
\psi(x,y,z) &= a_{1,xz}y+a_{1,y}x+2a_{2,xz}(x+z)y+a_{2,y}(x^{2}+2xz-y^{2}) \\
            &+ a_{3,y}[x^{3}-3xy^{2}+3(x^{2}-y^{2})z]+a_{3,xz}(3x^{2}y-y^{3}+6xyz) \\
            &+ a_{4,y}[x^{4}-6x^{2}y^{2}+y^{4}+4(x^{3}-3xy^{2})z] \\
            &+ 4a_{4,xz}[x^{3}y-xy^{3}+(3x^{2}y-y^{3})z].
\end{split}
\end{equation}

First, we get the first-order partial derivatives of $\phi$ and $\psi$ which are
\begin{equation}\label{eq:chap_2_15_21}
\begin{split}
\phi_{x} &= a_{1,xz}+2a_{2,xz}(x+z)-2a_{2,y}y+3a_{3,xz}(x^{2}+2xz-y^{2}) \\
         &- 6a_{3,y}(x+z)y-4a_{4,y}[3x^{2}y-y^{3}+6xyz] \\
         &+ 4a_{4,xz}[x^{3}-3xy^{2}+3(x^{2}-y^{2})z], \\
\phi_{y} &= -a_{1,y}-2a_{2,xz}y-2a_{2,y}(x+z)-3a_{3,y}(x^{2}+2xz-y^{2}) \\
         &- 6a_{3,xz}(x+z)y-4a_{4,xz}(3x^{2}y-y^{3}+6xyz) \\
         &- 4a_{4,y}[x^{3}-3xy^{2}+3(x^{2}-y^{2})z], \\
\phi_{z} &= 2a_{2,xz}x-2a_{2,y}y+3a_{3,xz}(x^{2}-y^{2})-6a_{3,y}xy \\
         &+ 4a_{4,xz}(x^{3}-3xy^{2})-4a_{4,y}(3x^{2}y-y^{3})
\end{split}
\end{equation}
and
\begin{equation}\label{eq:chap_2_15_22}
\begin{split}
\psi_{x} &= a_{1,y}+2a_{2,xz}y+2a_{2,y}(x+z)+3a_{3,y}(x^{2}+2xz-y^{2}) \\
         &+ 6a_{3,xz}(x+z)y+4a_{4,xz}(3x^{2}y-y^{3}+6xyz) \\
         &+ 4a_{4,y}[x^{3}-3xy^{2}+3(x^{2}-y^{2})z], \\
\psi_{y} &= a_{1,xz}+2a_{2,xz}(x+z)-2a_{2,y}y+3a_{3,xz}(x^{2}+2xz-y^{2}) \\
         &- 6a_{3,y}(x+z)y-4a_{4,y}(3x^{2}y-y^{3}+6xyz) \\
         &+ 4a_{4,xz}[x^{3}-3xy^{2}+3(x^{2}-y^{2})z], \\
\psi_{z} &= 2a_{2,xz}y+2a_{2,y}x+3a_{3,y}(x^{2}-y^{2})+6a_{3,xz}xy \\
         &+ 4a_{4,y}(x^{3}-3xy^{2})+4a_{4,xz}(3x^{2}y-y^{3}).
\end{split}
\end{equation}
These first-order partial derivatives of $\phi$ and $\psi$ satisfy the Cauchy-Riemann equations~(\ref{eq:chap_2_15_5})

Next, we get the second-order partial derivatives of $\phi$ and $\psi$ which are
\begin{equation}\label{eq:chap_2_15_23}
\begin{split}
\phi_{xx} &= 2a_{2,xz}+6a_{3,xz}(x+z)-6a_{3,y}y+12a_{4,xz}(x^{2}+2xz-y^{2}) \\
          &- 24a_{4,y}(x+z)y, \\
\phi_{xy} &= -2a_{2,y}-6a_{3,xz}y-6a_{3,y}(x+z)-12a_{4,y}(x^{2}+2xz-y^{2}) \\
          &- 24a_{4,xz}(x+z)y, \\
\phi_{xz} &= 2a_{2,xz}+6a_{3,xz}x-6a_{3,y}y+12a_{4,xz}(x^{2}-y^{2})-24a_{4,y}xy, \\
\end{split}
\end{equation}
\begin{equation}\label{eq:chap_2_15_24}
\begin{split}
\phi_{yx} &= -2a_{2,y}-6a_{3,xz}y-6a_{3,y}(x+z)-12a_{4,y}[x^{2}+2xz-y^{2}) \\
          &- 24a_{4,xz}(x+z)y, \\
\phi_{yy} &= -2a_{2,xz}-6a_{3,xz}(x+z)+6a_{3,y}y-12a_{4,xz}(x^{2}+2xz-y^{2}) \\
          &+ 24a_{4,y}(x+z)y, \\
\phi_{yz} &= -2a_{2,y}-6a_{3,xz}y-6a_{3,y}x-12a_{4,y}(x^{2}-y^{2})-24a_{4,xz}xy, \\
\end{split}
\end{equation}
\begin{equation}\label{eq:chap_2_15_25}
\begin{split}
\phi_{zx} &= 2a_{2,xz}+6a_{3,xz}x-6a_{3,y}y+12a_{4,xz}(x^{2}-y^{2})-24a_{4,y}xy, \\
\phi_{zy} &= -2a_{2,y}-6a_{3,xz}y-6a_{3,y}x-12a_{4,y}(x^{2}-y^{2})-24a_{4,xz}xy, \\
\phi_{zz} &= 0,
\end{split}
\end{equation}
and
\begin{equation}\label{eq:chap_2_15_26}
\begin{split}
\psi_{xx} &= 2a_{2,y}+6a_{3,xz}y+6a_{3,y}(x+z)+12a_{4,y}(x^{2}+2xz-y^{2}) \\
          &+ 24a_{4,xz}(x+z)y, \\
\psi_{xy} &= 2a_{2,xz}+6a_{3,xz}(x+z)-6a_{3,y}y+12a_{4,xz}(x^{2}+2xz-y^{2}) \\
          &- 24a_{4,y}(x+z)y, \\
\psi_{xz} &= 2a_{2,y}+6a_{3,xz}y+6a_{3,y}x+12a_{4,y}(x^{2}-y^{2})+24a_{4,xz}xy, \\
\end{split}
\end{equation}
\begin{equation}\label{eq:chap_2_15_27}
\begin{split}
\psi_{yx} &= 2a_{2,xz}+6a_{3,xz}(x+z)-6a_{3,y}y+12a_{4,xz}(x^{2}+2xz-y^{2}) \\
          &- 24a_{4,y}(x+z)y, \\
\psi_{yy} &= -2a_{2,y}-6a_{3,xz}y-6a_{3,y}(x+z)-12a_{4,y}(x^{2}+2xz-y^{2}) \\
          &- 24a_{4,xz}(x+z)y, \\
\psi_{yz} &= 2a_{2,xz}+6a_{3,xz}x-6a_{3,y}y+12a_{4,xz}(x^{2}-y^{2})-24a_{4,y}xy, \\
\end{split}
\end{equation}
\begin{equation}\label{eq:chap_2_15_28}
\begin{split}
\psi_{zx} &= 2a_{2,y}+6a_{3,xz}y+6a_{3,y}x+12a_{4,y}(x^{2}-y^{2})+24a_{4,xz}xy, \\
\psi_{zy} &= 2a_{2,xz}+6a_{3,xz}x-6a_{3,y}y+12a_{4,xz}(x^{2}-y^{2})-24a_{4,y}xy, \\
\psi_{zz} &= 0.
\end{split}
\end{equation}
These second-order partial derivatives of $\phi$ and $\psi$ satisfy the Laplace's equations~(\ref{eq:chap_2_15_6}) with boundary conditions in the $xyz$ coordinate space. The proof of lemma~(\ref{le:chap_2_15_3}) is completed.

\section{Uniquely Determined Analytic Functions}\label{sec:chap_2_16_26}

We conclude this chapter with two sections dealing with how the values of an analytic function in a domain $D$ are affected by its values in a subdomain or on a line segment lying in $D$. While these sections are of considerable theoretical interest, they are not central to our development of analytic functions in later chapters. The reader may pass directly to Chap.~(\ref{ch:chap_3}) at this time and refer back when necessary.

\begin{lemma}\label{le:chap_2_16_1}
Suppose that

(i) a function $f$ is analytic throughout a domain $D$;

(ii) $f(s)=0$ at each point $s$ of a domain or line segment contained in $D$.

Then $f(s)\equiv0$ in $D$; that is, $f(s)$ is identically equal to zero throughout $D$.
\end{lemma}

To prove this lemma, we let $f$ be as stated in its hypothesis and let $s_{0}$ be any point of the subdomain or line segment at each point of which $f(s)=0$. Since $D$ is a connected open set (Sec.~(\ref{sec:chap_1_10})), there is a polygonal line $L$, consisting of a finite number of line segments joined end to end and lying entirely in $D$, that extends from $s_{0}$ to any other point $P$ in $D$. We let $d$ be the shortest distance from points on $L$ to the boundary of $D$, unless $D$ is the entire space; in that case, $d$ may be any positive number. We then form a finite sequence of points
\[
s_{0},s_{1},s_{2},\ldots,s_{n-1},s_{n}
\]
along $L$, where the point $s_{n}$ coincides with $P$ and where each point is sufficiently close to the adjacent ones that
\[
|s_{k}-s_{k-1}|<d\texttt{ }(k=1,2,\ldots,n).
\]
Finally, we construct a finite sequence of neighborhoods
\[
N_{0},N_{1},N_{2},\ldots,N_{n-1},N_{n}
\]
where each neighborhood $N_{k}$ is centered at $s_{k}$ and has radius $d$. Note that these
neighborhoods are all contained in $D$ and that the center $s_{k}$ of any neighborhood $N_{k}$ $(k=1,2,\ldots,n)$ lies in the preceding neighborhood $N_{k-l}$.

At this point, we need to use a result that is proved later on in Chap.~(\ref{ch:chap_7}). Namely, Theorem~(\ref{th:chap_7_7_3}) in Sec.~(\ref{sec:chap_7_7_68}) tells us that since $f$ is analytic in the domain $N_{0}$ and since $f(s)=0$ in a domain or on a line segment containing $s_{0}$, then $f(s)\equiv0$ in $N_{0}$. But the point $s_{1}$ lies in the domain $N_{0}$. Hence a second application of the same theorem reveals that $f(s)\equiv0$ in $N_{1}$; and, by continuing in this manner, we arrive at the fact that $f(s)\equiv0$ in $N_{n}$. Since $N_{n}$ is centered at the point $P$ and since $P$ was arbitrarily selected in $D$, we may conclude that $f(s)\equiv0$ in $D$. This completes the proof of the lemma.

Suppose now that two functions $f$ and $g$ are analytic in the same domain $D$ and that $f(s)=g(s)$ at each point  $s$ of some domain or line segment contained in $D$. The difference
\[
h(s)=f(s)-g(s)
\]
is also analytic in $D$, and $h(s)=0$ throughout the subdomain or along the line segment. According to the above lemma, then, $h(s)=0$ throughout $D$; that is, $f(s)=g(s)$ at each point $s$ in $D$. We thus arrive at the following important theorem.

\begin{theorem}\label{th:chap_2_16_1}
A function that is analytic in a domain $D$ is uniquely determined over $D$ by its values in a domain, or along a line segment, contained in $D$.
\end{theorem}

This theorem is useful in studying the question of extending the domain of definition of an analytic function. More precisely, given two domains $D_{1}$ and $D_{2}$, consider the intersection $D_{1}\bigcap D_{2}$, consisting of all points that lie in both $D_{1}$ and $D_{2}$. If $D_{1}$ and $D_{2}$ have points in common and a function $f_{1}$ is analytic in $D_{1}$, there may exist a function $f_{2}$, which is analytic in $D_{2}$, such that $f_{2}(s)=f_{1}(s)$ for each $s$ in the intersection $D_{1}\bigcap D_{2}$. If so, we call $f_{2}$ an analytic continuation of $f_{1}(s)$ into the second domain $D_{2}$.

Whenever that analytic continuation exists, it is unique, according to the theorem just proved. That is, not more than one function can be analytic in $D_{2}$ and assume the value $f_{1}(s)$ at each point $s$ of the domain $D_{1}\bigcap D_{2}$ interior to $D_{2}$. However, if there is an analytic continuation $f_{3}$ of $f_{2}$ from $D_{2}$ into a domain $D_{3}$ which intersects $D_{1}$, it is not necessarily true that $f_{3}(s)=f_{1}(s)$ for each $s$ in $D_{1}\bigcap D_{3}$.

If $f_{2}$ is the analytic continuation of $f_{1}$ from a domain $D_{1}$ into a domain $D_{2}$, then the function  $F$ defined by the equations
\[
F(s)=\{
    \begin{array}{cc}
        f_{1}(s) & \texttt{ when }s\texttt{ is in }D_{1}, \\
        f_{2}(s) & \texttt{ when }s\texttt{ is in }D_{2}
    \end{array}
\]
is analytic in the union $D_{1}\bigcup D_{2}$, which is the domain consisting of all points that lie
in either $D_{1}$ or $D_{2}$. The function $F$ is the analytic continuation into $D_{1}\bigcup D_{2}$ of either
$f_{1}$ or $f_{2}$ and $f_{1}$ and $f_{2}$ are called elements of $F$.

\section{Reflection Principle}\label{sec:chap_2_17_27}

The theorem in this section concerns the fact that some analytic functions possess the property that $\overline{f(s)}=f(\overline{s})$ for all points $s$ in certain domains, while others do not. We note, for example, that $s+1$ and $s^{2}$ have that property when $D$ is the entire finite space; but the same is not true of $s+i$ and $iz^{2}$. The theorem, which is known as the reflection principle, provides a way of predicting when $\overline{/(s)}=/(\overline{s})$.

\begin{theorem}\label{th:chap_2_17_l}
Suppose that a function $f$ is analytic in some domain $D$ which contains a segment of the $x$ axis and whose lower  half is the reflection of the upper half with respect to that axis in the $xy$ coordinate plane for each $z$. Then
\begin{equation}\label{eq:chap_2_17_1}
\overline{f(s)}=f(\overline{s})
\end{equation}
for each point $s$ in the domain $if$ and only if $f(e_{xy}(x+i0)+0)$ is real for each point $x$ on the segment.
\end{theorem}

We start the proof by assuming that $f(e_{xy}(x+i0)+0)$ is real at each point $x$ on the segment. Once we show that the function
\begin{equation}\label{eq:chap_2_17_2}
F(s)=\overline{f(\overline{s})}.
\end{equation}
is analytic in $D$, we shall use it to obtain equation~(\ref{eq:chap_2_17_1}). To establish the analyticity of
$F(s)$, we write
\[
f(s)=e_{xy}[u(x,y,z)+iv(x,y,z)]+w(z),
\]
\[
F(s)=e_{xy}[U(x,y,z)+iV(x,y,z)]+W(z)
\]
and observe how it follows from equation~(\ref{eq:chap_2_17_2}) that, since
\begin{equation}\label{eq:chap_2_17_3}
\overline{f(\overline{s})}=e_{xy}[u(x,-y,z)-iv(x,-y,z)]+w(z),
\end{equation}
the components of $F(s)$ and $f(s)$ are related by the equations
\begin{equation}\label{eq:chap_2_17_4}
U(x,y,z)=u(x,t,z)\texttt{, }V(x,y,z)=-v(x,t,z)\texttt{, }W(z)=w(z),
\end{equation}
where $t=-y$. Now, because $f(e_{xy}(x+it)+z)$ is an analytic function of $e_{xy}(x+it)+z$, the first-order partial derivatives of the functions $u(x,t,z)$, $v(x,t,z)$, and $w(z)$ are continuous throughout $D$ and satisfy the Cauchy-Riemann equations
\begin{equation}\label{eq:chap_2_17_5}
e_{xy}u_{x}(x,t,z)=e_{xy}v_{t}(x,t,z)\texttt{, }e_{xy}u_{t}(x,t,z)=-e_{xy}v_{x}(x,t,z).
\end{equation}
Furthermore, in view of equations~(\ref{eq:chap_2_17_4}),
\[
U_{x}=u_{x}\texttt{, }V_{y}=-v_{t}\frac{dt}{dy}=v_{t}\texttt{, }W_{z}=w_{z},
\]
and it follows from these and the first of equations~(\ref{eq:chap_2_17_5}) that $U_{x}=V_{y}$. Similarly,
\[
U_{y}=u_{t}\frac{dt}{dy}=-u_{t}\texttt{, }V_{x}=-v_{x};
\]
and the second of equations~(\ref{eq:chap_2_17_5}) tells us that $U_{y}=-V_{x}$. Inasmuch as the first-order
partial derivatives of $U(x,y,z)$ and $V(x,y,z)$ are now shown to satisfy the Cauchy-Riemann equations and since those derivatives are continuous, we find that the function $F(s)$ is analytic in $D$. Moreover, since $f(e_{xy}(x+i0)+0)$ is real on the segment of the real axis in the $xy$ coordinate plane lying in $D$, $v(x,0,0)=0$ and $w(0)=0$ on that segment; and, in view of equations~(\ref{eq:chap_2_17_4}), this means that
\[
F(e_{xy}(x+i0)+0)=e_{xy}[U(x,0,0)+iV(x,0,0)]+W(0)=e_{xy}U(x,0,0)
\]
That is,
\begin{equation}\label{eq:chap_2_17_6}
F(s)=f(s).
\end{equation}
at each point on the segment. We now refer to the theorem in Sec.~(\ref{sec:chap_2_16_26}), which tells us that an analytic function defined on a domain $D$ is uniquely determined by its values along any line segment lying in $D$. Thus equation~(\ref{eq:chap_2_17_6}) actually holds throughout $D$. Because of definition~(\ref{eq:chap_2_17_2}) of  the function $F(s)$, then,
\begin{equation}\label{eq:chap_2_17_7}
\overline{f(\overline{s})}=f(s).
\end{equation}
and this is the same as equation~(\ref{eq:chap_2_17_1}).

To prove the converse of the theorem, we assume that equation~(\ref{eq:chap_2_17_1}) holds and note
that, in view of expression~(\ref{eq:chap_2_17_3}), the form~(\ref{eq:chap_2_17_7}) of equation~(\ref{eq:chap_2_17_1}) can be written
\[
e_{xy}[u(x,-y,z)-iv(x,-y,z)]+w(z)=e_{xy}[u(x,y,z)+iv(x,y,z)]+w(z).
\]
In particular, if $(x,0,0)$ is a point on the segment of the real axis that lies in $D$,
\[
e_{xy}[u(x,0,0)-iv(x,0,0)]+w(0)=e_{xy}[u(x,0,0)+iv(x,0,0)]+w(0);
\]
and, by equating imaginary parts here, we see that $v(x,0,0)=0$. Hence $f(e_{xy}(x+i0)+0)$ is real on
the segment of the real axis in the $xy$ coordinate plane lying in $D$.

%-----------------------------------------------------------------------
% Beginning of chap3.tex
%-----------------------------------------------------------------------
%
% AMS-LaTeX 1.2 sample file for a monograph, based on amsbook.cls.
% This is a data file input by chapter.tex.
%%%%%%%%%%%%%%%%%%%%%%%%%%%%%%%%%%%%%%%%%%%%%%%%%%%%%%%%%%%%%%%%%%%%

%\part{This is a Part Title Sample}

\chapter{Elementary Functions}\label{ch:chap_3}

We consider here various elementary functions studied in calculus and two-dimensional complex numbers, and define corresponding spatial complex functions of a spatial complex variable, which are easily proved true by definitions~(\ref{eq:chap_1_1_1}) and~(\ref{eq:chap_1_1_2}), and the expansion of each elementary function into a power series. To be specific, we define analytic functions of a spatial complex variable $s$ that reduce to the elementary functions in calculus when $s=e_{xy}(0+i0)+z$. We start by defining the spatial complex exponential function and then use it to develop the others.

Most of elementary functions $f(s)$ defined here can be expressed as follows:
\[%begin{equation}\label{eq:chap_3_1_0}
\begin{split}
f(s) &= e_{xy}[f(x+iy+z)-f(z)]+f(z) \\
     &= e_{xy}[u(x,y,z)+iv(x,y,z)]+w(z)
\end{split}
\]%end{equation}
where
\[
\begin{split}
u(x,y,z) &= Re\texttt{ }f(x+iy+z)-f(z), \\
v(x,y,z) &= Im\texttt{ }f(x+iy+z),\texttt{ and }w(z)=f(z).
\end{split}
\]
For $|s|=1$, there are $x+z=\cos\theta\cos\varphi+\sin\varphi$ and $-\sqrt{2}\leq x+z\leq\sqrt{2}$.

\section{Polynomial Functions}\label{sec:chap_3_1_27}

\subsection{Polynomial Functions $p(s,m)$}\label{sec:chap_3_1_27_1}

The spatial complex polynomial function $p(s,m)$ of a spatial complex variable $s=e_{xy}(x+iy)+z$ is defined as
\begin{equation}\label{eq:chap_3_1_1}
\begin{split}
p(s,m) &= \sum_{n=1}^{m}a_{n}s^{n}=\sum_{n=0}^{m}a_{n}[e_{xy}(x+iy)+z]^{n}\texttt{ }(|s|<+\infty) \\
       &= e_{xy}(p(x+iy+z,m)-p(z,m))+p(z,m) \\
       &= e_{xy}Im_{s}p(s,m)+Re_{s}p(s,m)
\end{split}
\end{equation}
where spatial complex constants $a_{n}=e_{xy}(a_{n,x}+ia_{n,y})+a_{n,z}$ for $n=1,2,\ldots,m$, and
\[
\begin{split}
Re_{s}p(s,m) &= \sum_{n=1}^{m}a_{n,z}z^{n}\texttt{ and } \\
Im_{s}p(s,m) &= \sum_{n=1}^{m}a_{n}(x+iy+z)^{n}-Re_{s}p(s,m).
\end{split}
\]

Let $a_{n,xz}$ denote $a_{n,x}+a_{n,z}$. Then there is $e_{xy}a_{n}=e_{xy}(a_{n,xz}+ia_{n,y})$. For $n=1,2,\ldots,m$ with definition $C_{n}^{k}=\frac{n!}{(n-k)!k!}$ there are
\[
(x+z+iy)^{n}=\sum_{k=0}^{n}C_{n}^{k}(x+z)^{n-k}(iy)^{k}
\]
\[
=\sum_{k=0}^{n/2}(-1)^{k}C_{n}^{2k}(x+z)^{n-2k}y^{2k}
+i\sum_{k=1}^{(n+1)/2}(-1)^{k-1}C_{n}^{2k-1}(x+z)^{n-2k+1}y^{2k-1}.
\]

Thus we get
\[
\begin{split}
Im_{s}p(s,m)+Re_{s}p(s,m) &= \sum_{n=1}^{m}a_{n}(x+z+iy)^{n} \\
                          &= \sum_{n=1}^{m}(a_{n,xz}+ia_{n,y})(x+z+iy)^{n} \\
\end{split}
\]
\[
=\sum_{n=1}^{m}[a_{n,xz}\sum_{k=0}^{n/2}(-1)^{k}C_{n}^{2k}(x+z)^{n-2k}y^{2k}-a_{n,y}\sum_{k=1}^{(n+1)/2}(-1)^{k-1}C_{n}^{2k-1}(x+z)^{n-2k+1}y^{2k-1}]
\]
\[
+i\sum_{n=1}^{m}[a_{n,xz}\sum_{k=1}^{(n+1)/2}(-1)^{k-1}C_{n}^{2k-1}(x+z)^{n-2k+1}y^{2k-1}+a_{n,y}\sum_{k=0}^{n/2}(-1)^{k}C_{n}^{2k}(x+z)^{n-2k}y^{2k}].
\]

That is
\begin{equation}\label{eq:chap_3_1_2}
Re\texttt{ }Im_{s}p(s,m)=\sum_{n=1}^{m}[a_{n,xz}\sum_{k=0}^{n/2}(-1)^{k}C_{n}^{2k}(x+z)^{n-2k}y^{2k}
\end{equation}
\[
-a_{n,y}\sum_{k=1}^{(n+1)/2}(-1)^{k-1}C_{n}^{2k-1}(x+z)^{n-2k+1}y^{2k-1}]-Re_{s}p(s,m),
\]
and
\begin{equation}\label{eq:chap_3_1_3}
Im\texttt{ }Im_{s}p(s,m)=\sum_{n=1}^{m}[a_{n,xz}\sum_{k=1}^{(n+1)/2}(-1)^{k-1}C_{n}^{2k-1}(x+z)^{n-2k+1}y^{2k-1}
\end{equation}
\[
+a_{n,y}\sum_{k=0}^{n/2}(-1)^{k}C_{n}^{2k}(x+z)^{n-2k}y^{2k}].
\]

\subsection{Polynomial Functions $q(s,m)$}\label{sec:chap_3_1_27_2}

The spatial complex polynomial function $q(s,m)$ of a spatial complex variable $s=e_{xy}(x+iy)+z$ is defined as
\begin{equation}\label{eq:chap_3_1_4}
\begin{split}
q(s,m) &= \sum_{n=1}^{m}\frac{b_{n}}{s^{n}}=\sum_{n=1}^{m}\frac{b_{n}}{[e_{xy}(x+iy)+z]^{n}}\texttt{ }(s\neq0) \\
       &= e_{xy}(q(x+iy+z,m)-q(z,m))+q(z,m) \\
       &= e_{xy}Im_{s}q(s,m)+Re_{s}q(s,m)
\end{split}
\end{equation}
where spatial complex constants $b_{n}=e_{xy}(b_{n,x}+ib_{n,y})+b_{n,z}$ for $n=1,2,\ldots,m$, and
\[
\begin{split}
Re_{s}q(s,m) &= \sum_{n=1}^{m}\frac{b_{n,z}}{z^{n}}\texttt{ and } \\
Im_{s}q(s,m) &= \sum_{n=1}^{m}\frac{b_{n}}{(x+iy+z)^{n}}-Re_{s}q(s,m).
\end{split}
\]

Let $b_{n,xz}$ denote $b_{n,x}+b_{n,z}$. Then there is $e_{xy}b_{n}=e_{xy}(b_{n,xz}+ib_{n,y})$. For $n=1,2,\ldots,m$ with definition $C_{n}^{k}=\frac{n!}{(n-k)!k!}$ there are
\[
\begin{split}
\frac{1}{(x+z+iy)^{n}} &= (\frac{x+z-iy}{(x+z)^{2}+y^{2}})^{n} \\
                       &= (\frac{1}{(x+z)^{2}+y^{2}})^{n}\sum_{k=0}^{n}C_{n}^{k}(x+z)^{n-k}(-iy)^{k} \\
                       &= (\frac{1}{(x+z)^{2}+y^{2}})^{n}[\sum_{k=0}^{n/2}(-1)^{k}C_{n}^{2k}(x+z)^{n-2k}y^{2k} \\
                       &+ i\sum_{k=1}^{(n+1)/2}(-1)^{k}C_{n}^{2k-1}(x+z)^{n-2k+1}y^{2k-1}].
\end{split}
\]

Thus we get
\[
Im_{s}q(s,m)+Re_{s}q(s,m)=\sum_{n=1}^{m}\frac{b_{n}}{(x+z+iy)^{n}}=\sum_{n=1}^{m}\frac{b_{n,xz}+ib_{n,y}}{(x+z+iy)^{n}}
\]
\[
=\sum_{n=1}^{m}(\frac{1}{(x+z)^{2}+y^{2}})^{n}[b_{n,xz}\sum_{k=0}^{n/2}(-1)^{k}C_{n}^{2k}(x+z)^{n-2k}y^{2k}
\]
\[
-b_{n,y}\sum_{k=1}^{(n+1)/2}(-1)^{k}C_{n}^{2k-1}(x+z)^{n-2k+1}y^{2k-1}]
\]
\[
+i\sum_{n=1}^{m}(\frac{1}{(x+z)^{2}+y^{2}})^{n}[b_{n,xz}\sum_{k=1}^{(n+1)/2}(-1)^{k}C_{n}^{2k-1}(x+z)^{n-2k+1}y^{2k-1}
\]
\[
+\frac{b_{n,y}}{(x+z+iy)^{n}}\sum_{k=0}^{n/2}(-1)^{k}C_{n}^{2k}(x+z)^{n-2k}y^{2k}].
\]

That is
\begin{equation}\label{eq:chap_3_1_5}
\begin{split}
Re\texttt{ }Im_{s}q(s,m) &= \sum_{n=1}^{m}(\frac{1}{(x+z)^{2}+y^{2}})^{n}[b_{n,xz}\sum_{k=0}^{n/2}(-1)^{k}C_{n}^{2k}(x+z)^{n-2k}y^{2k} \\
&-b_{n,y}\sum_{k=1}^{(n+1)/2}(-1)^{k}C_{n}^{2k-1}(x+z)^{n-2k+1}y^{2k-1}]-Re_{s}q(s,m)
\end{split}
\end{equation}
and
\begin{equation}\label{eq:chap_3_1_6}
\begin{split}
Im\texttt{ }Im_{s}q(s,m) &= \sum_{n=1}^{m}(\frac{1}{(x+z)^{2}+y^{2}})^{n}[
b_{n,y}\sum_{k=0}^{n/2}(-1)^{k}C_{n}^{2k}(x+z)^{n-2k}y^{2k}] \\
&+b_{n,xz}\sum_{k=1}^{(n+1)/2}(-1)^{k-1}C_{n}^{2k-1}(x+z)^{n-2k+1}y^{2k-1}.
\end{split}
\end{equation}

\subsection{Polynomial Functions $p_{a}(z,m)$ and $q_{a}(z,m)$}\label{sec:chap_3_1_27_3}

Let spatial complex constants
\[
\begin{split}
a_{n} &= e_{xy}(a_{n,x}+ia_{n,y})+a_{n,z}\texttt{ and } \\
b_{n} &= e_{xy}(b_{n,x}+ib_{n,y})+b_{n,z}\texttt{ }(n=1,2,\ldots,m).
\end{split}
\]

The two-dimensional complex polynomial function $p_{a}(z,m)$ of a real variable $z$ is defined as
\begin{equation}\label{eq:chap_3_1_7}
p_{a}(z,m)=e_{xy}\sum_{n=1}^{m}(a_{n,x}+ia_{n,y})z^{n}\texttt{ }(|z|<+\infty)
\end{equation}
where
\begin{equation}\label{eq:chap_3_1_8}
p_{a}(z,m)+Re_{s}p(s,m)=\sum_{n=1}^{m}a_{n}z^{n}\texttt{ }(|z|<+\infty).
\end{equation}

The two-dimensional complex polynomial function $q_{a}(z,m)$ of a real variable $z$ is defined as
\begin{equation}\label{eq:chap_3_1_9}
q_{a}(z,m)=e_{xy}\sum_{n=1}^{m}\frac{b_{n,x}+ib_{n,y}}{z^{n}}\texttt{ }(z\neq0)
\end{equation}
where
\begin{equation}\label{eq:chap_3_1_10}
q_{a}(z,m)+Re_{s}q(s,m)=\sum_{n=1}^{m}\frac{b_{n}}{z^{n}}\texttt{ }(z\neq0).
\end{equation}

\subsection{Polynomial Functions $p_{c}(z,m)$ and $q_{c}(z,m)$}\label{sec:chap_3_1_27_4}

Let a spatial complex constant $a_{n}=e_{xy}(a_{n,x}+ia_{n,y})+a_{n,z}$ and $a_{n,xz}$ denote $a_{n,x}+a_{n,z}$ for $n=1,2,\ldots,m$.

The two-dimensional complex polynomial function $p_{c}(z,m)$ of a two-dimensional complex variable $x+iy$ and a real variable $z$ is defined as
\begin{equation}\label{eq:chap_3_1_11}
p_{c}(z,m)=e_{xy}\sum_{n=3}^{m}a_{n}\sum_{k=2}^{n-1}C_{n}^{k}(x+iy)^{n-k}z^{k}\texttt{ }(|(x+iy)z|<+\infty)
\end{equation}

Let a spatial complex constant $b_{n}=e_{xy}(b_{n,x}+ib_{n,y})+b_{n,z}$ and $b_{n,xz}$ denote $b_{n,x}+b_{n,z}$ for $n=1,2,\ldots,m$.

The two-dimensional complex polynomial function $q_{c}(z,m)$ of a real variable $z$ is defined as
\begin{equation}\label{eq:chap_3_1_12}
q_{c}(z,m)=e_{xy}\sum_{n=3}^{m}b_{n}\sum_{k=2}^{n-1}C_{n}^{k}(\frac{1}{x+iy})^{n-k}z^{k}\texttt{ }(|x+iy|>0,\texttt{ }|z|<+\infty).
%q_{c}(z,m)=e_{xy}\sum_{n=1}^{m}\frac{b_{n,x}+ib_{n,y}}{z^{n}}\texttt{ }(z\neq0)
\end{equation}

\section{The Exponential Function}\label{sec:chap_3_2_28}

\begin{definition}\label{def:chap_3_2_1}
The spatial complex exponential function $e^{s}$ is defined as
\begin{equation}\label{eq:chap_3_2_1}
\begin{split}
e^{s} &= e^{e_{xy}(x+iy)+z} \\
      &= \{
    \begin{array}{cc}
        e_{xy}(e^{x+iy+z}-e^{z})+e^{z} & \texttt{ when }z\neq0 \\
        e_{xy}e^{x+iy} & \texttt{ when }z=0
    \end{array}\texttt{ }(|s|<+\infty).
\end{split}
\end{equation}
\end{definition}

Now, we shall prove equation~(\ref{eq:chap_3_2_1}).

Let $\sigma=x+iy$. Then $\sigma$ and $\sigma+z$ are two-dimensional complex numbers. The Taylor series of the real exponential function $e^{z}$ and two-dimensional complex exponential function $e^{\sigma+z}$ are
\[
e^{z}=\sum_{n=0}^{\infty}\frac{z^{n}}{n!}\texttt{ and }
e^{\sigma+z}=\sum_{n=0}^{\infty}\frac{(\sigma+z)^{n}}{n!}\texttt{ }(|\sigma+z|<\infty).
\]

For an index $n$ in the Taylor series of $e^{z}$ and $e^{\sigma+z}$, there are $z^{n}$ and
\[
(\sigma+z)^{n}=\sigma^{n}+n\sigma^{n-1}z+\cdots+n\sigma z^{n-1}+z^{n}.
\]

Since both $\sigma+z$ and its powers are points in the two-dimensional complex plane and from remark~(\ref{re:chap_1_1_1}) for a term $e_{xy}x^{j}y^{k}$, the operator $e_{xy}$ becomes $e_{xy}^{m}$ where $m=|j|+|k|$, we get
\[
\begin{split}
e_{xy}[(\sigma+z)^{n}-z^{n}]+z^{n}
&= (e_{xy}\sigma)^{n}+n(e_{xy}\sigma)^{n-1}z+\cdots+n(e_{xy}\sigma)z^{n-1}+z^{n} \\
&= (e_{xy}\sigma+z)^{n}=s^{n}\texttt{ for }n=1,2,\ldots.
\end{split}
\]

Thus for $|s|<\infty$ and $z\neq0$ there are $|\sigma+z|<\infty$ and
\[
\begin{split}
e_{xy}(e^{\sigma+z}-e^{z})+e^{z}
&=e_{xy}[\sum_{n=0}^{\infty}\frac{(\sigma+z)^{n}}{n!}-\sum_{n=0}^{\infty}\frac{z^{n}}{n!}]
+\sum_{n=0}^{\infty}\frac{z^{n}}{n!} \\
&=e_{xy}\sum_{n=1}^{\infty}\frac{(\sigma+z)^{n}-z^{n}}{n!}+\sum_{n=0}^{\infty}\frac{z^{n}}{n!} \\
&=\sum_{n=0}^{\infty}\frac{(e_{xy}\sigma+z)^{n}}{n!}
=\sum_{n=0}^{\infty}\frac{s^{n}}{n!}
=e^{s}.
\end{split}
\]
So the Taylor series of the spatial complex exponential function $e^{s}$ exits and is
\begin{equation}\label{eq:chap_3_2_2}
e^{s}=\sum_{n=0}^{\infty}\frac{s^{n}}{n!}\texttt{ }(|s|<\infty).
\end{equation}

According to definition~(\ref{eq:chap_3_2_1}), the property's extension,
\begin{equation}\label{eq:chap_3_2_3}
e^{s_{1}}e^{s_{2}}=e_{xy}(e^{x_{1}+x_{2}+i(y_{1}+y_{2})+z_{1}+z_{2}}-e^{z_{1}+z_{2}})+e^{z_{1}+z_{2}}\texttt{ }(|s_{1}+s_{2}|<+\infty),
\end{equation}
to spatial complex analysis is easy to prove. To do this, we write
\[
s_{1}=e_{xy}(x_{1}+iy_{1})+z_{1}\texttt{ and }s_{2}=e_{xy}(x_{2}+iy_{2})+z_{2}.
\]
Then
\[
\begin{split}
e^{s_{1}}e^{s_{2}} &= e^{e_{xy}(x_{1}+x_{2}+i(y_{1}+y_{2}))+z_{1}+z_{2}} \\
                   &= e_{xy}(e^{x_{1}+x_{2}+i(y_{1}+y_{2})+z_{1}+z_{2}}-e^{z_{1}+z_{2}})+e^{z_{1}+z_{2}}.
\end{split}
\]

Observe how property~(\ref{eq:chap_3_2_3}) enables us to write $e^{s_{1}-s_{2}}e^{s_{2}}=e^{s_{1}}$, or
\begin{equation}\label{eq:chap_3_2_4}
\frac{e^{s_{1}}}{e^{s_{2}}}=e_{xy}(e^{x_{1}-x_{2}+i(y_{1}-y_{2})+z_{1}-z_{2}}-e^{z_{1}-z_{2}})+e^{z_{1}-z_{2}}\texttt{ }(|s_{1}-s_{2}|<+\infty).
\end{equation}

From this and the fact that $e^{0}=1$, it follows that $\frac{1}{e^{s}}=e^{-s}$.

There are a number of other important properties of $e^{s}$ that are expected. For instance,
\begin{equation}\label{eq:chap_3_2_5}
\frac{d}{ds}e^{s}=e^{s}
\end{equation}
everywhere in the $s$ space. Note that the differentiability of $e^{s}$ for all $s$ tells us that $e^{s}$ is  entire (Sec.~(\ref{sec:chap_2_13_23})). It is also true that
\begin{equation}\label{eq:chap_3_2_6}
e^{s}\neq0\texttt{ for any spatial complex number }s.
\end{equation}

This is evident upon writing definition~(\ref{eq:chap_3_2_1}) in the form
\begin{equation}\label{eq:chap_3_2_7}
\begin{split}
e^{s} &= e^{z}e^{e_{xy}(x+iy)} \\
      &= \rho(e_{xy}e^{i\phi}\cos\psi+\sin\psi) \\
      &= e_{xy}[u(x,y,z)+iv(x,y,z)]+w(z)\texttt{ }(|s|<+\infty)
\end{split}
\end{equation}
where
\[
    \begin{array}{cc}
        \rho=(e^{2x}+e^{2z})^{1/2},\texttt{ }\phi=y+2n\pi,\texttt{ }\psi=\arctan e^{z-x} & \texttt{ when }z\neq0 \\
        \rho=e^{x},\texttt{ }\phi=y+2n\pi,\texttt{ }\psi=0 & \texttt{ when }z=0
    \end{array}
\]
where $n$ is some integer (see Sec.~(\ref{sec:chap_1_8})); and
\[
    \begin{array}{cc}
        u(x,y,z)=e^{x+z}\cos y-e^{z},\texttt{ }v(x,y,z)=e^{x+z}\sin y,\texttt{ }w(z)=e^{z} & \texttt{ when }z\neq0 \\
        u(x,y,z)=e^{x}\cos y,\texttt{ }v(x,y,z)=e^{x}\sin y,\texttt{ }w(z)=0 & \texttt{ when }z=0
    \end{array},
\]
which tells us that
\begin{equation}\label{eq:chap_3_2_8}
|e^{s}|=(e^{2x}+e^{2z})^{1/2},\texttt{ }\arg_{c}e^{s}=y+2n\pi\texttt{ }(n=0,\pm1,\pm2,\ldots)
\end{equation}
and
\[
\arg_{r}e^{s}=\{
    \begin{array}{cc}
        \arctan e^{z-x} & \texttt{ when }z\neq0 \\
        0 & \texttt{ when }z=0
    \end{array}.
\]

Some properties of $e^{s}$ are, however, \textit{not expected}. For example, since
\[
e^{s+2\pi i}=e^{s}e^{2\pi i}\texttt{ and }e^{2\pi i}=1,
\]
we find that $e^{s}$ is periodic, with a pure imaginary period $2\pi i$:
\begin{equation}\label{eq:chap_3_2_9}
e^{s+2\pi i}=e^{s}.
\end{equation}

\section{The Logarithmic Function}\label{sec:chap_3_3_29}

Our motivation for the definition of the logarithmic function is based on solving the equation
\begin{equation}\label{eq:chap_3_3_1}
e^{\varpi}=s
\end{equation}
for $\varpi$, where $s$ is any nonzero spatial complex number. To do this, we note that when $s$ and $\varpi$ are  written
\[
s=r(e_{xy}e^{i\theta}\cos\varphi+\sin\varphi)\texttt{ }(-\pi<\theta<\pi,\texttt{ }0\leq\varphi\leq\pi/2)
\]
where $r=|s|$, and $\varpi=e_{xy}(u+iv)+w$, equation~(\ref{eq:chap_3_3_1}) becomes
\begin{equation}\label{eq:chap_3_3_1b}
e^{\varpi}=e^{w}[e_{xy}(e^{u}e^{iv}-1)+1].
\end{equation}
Then, in view of the statement in italics in Sec.~(\ref{sec:chap_1_8}) regarding the equality of two spatial complex numbers expressed in exponential form,
\begin{equation}\label{eq:chap_3_3_1c}
\begin{split}
&u=\frac{1}{2}\ln[1+\cos\theta\sin2\varphi]-\ln\sin\varphi=\frac{1}{2}\ln[(x+z)^{2}+y^{2}]-\ln z, \\
&v=\arctan\frac{\sin\theta}{\cos\theta+\tan\varphi}+2n\pi=\arctan\frac{y}{x+z}+2n\pi, \\
&w=\{
    \begin{array}{cc}
        \ln(r\sin\varphi) & \texttt{ when }\sin\varphi>0 \\
        \ln r & \texttt{ when }\sin\varphi=0
    \end{array}, \\
&\texttt{ }(-\pi<\theta<\pi,\texttt{ }0\leq\varphi\leq\pi/2)
\end{split}
\end{equation}
where $n$ is any integer, and there are
\[
u_{x}=v_{y}=\frac{x+z}{(x+z)^{2}+y^{2}}\texttt{ and }u_{y}=-v_{x}=\frac{y}{(x+z)^{2}+y^{2}}.
\]
Since the equation $e^{w}=r\sin\varphi$ is the same as $w=\ln(r\sin\varphi)$, it follows that equation~(\ref{eq:chap_3_3_1}) is satisfied if and only if $\varpi$ has one of the values
\[
\varpi=\ln|s|+\{
    \begin{array}{cc}
        \log(e_{xy}e^{i\theta}\cos\varphi+\sin\varphi) & \texttt{ when }\sin\varphi>0 \\
        e_{xy}i\theta & \texttt{ when }\sin\varphi=0
    \end{array}
\]
\[
\texttt{ }(-\pi<\theta<\pi,\texttt{ }0\leq\varphi\leq\pi/2)
\]
Thus, if we write
\begin{equation}\label{eq:chap_3_3_2}
\log s=\ln|s|+\{
    \begin{array}{cc}
        \log(e_{xy}e^{i\theta}\cos\varphi+\sin\varphi) & \texttt{ when }\sin\varphi>0 \\
        e_{xy}i\theta & \texttt{ when }\sin\varphi=0
    \end{array}
\end{equation}
\[
\texttt{ }(|s|<+\infty\texttt{ and }0\leq\varphi\leq\pi/2),
\]
and from lemma~(\ref{le:chap_3_3_1}) for $|s|<+\infty$ and $0<\varphi\leq\pi/2$
\begin{equation}\label{eq:chap_3_3_2b}
\log s=\ln|s|
\end{equation}
\[
+\{
    \begin{array}{cc}
        \ln\sin\varphi+e_{xy}[\log(e^{i\theta}\cos\varphi+\sin\varphi)-\ln\sin\varphi] & \texttt{ when }\sin\varphi>0 \\
        e_{xy}i\theta & \texttt{ when }\sin\varphi=0
    \end{array}
\]
\[
=\ln|s|+\{
    \begin{array}{cc}
        \ln\sin\varphi+e_{xy}[\ln r_{L}(\theta,\varphi)+i\theta_{L}(\theta,\varphi)-\ln\sin\varphi] & \texttt{ when }\sin\varphi>0 \\
        e_{xy}i\theta & \texttt{ when }\sin\varphi=0
    \end{array}
\]
for $0\leq\varphi\leq\pi/2$ and $\theta\neq(2n+1)\pi$ where $n=0,1,2,\ldots$
\[
\begin{split}
r_{L}(\theta,\varphi) &= \sqrt{1+\cos\theta\sin2\varphi},\texttt{ } \\
\theta_{L}(\theta,\varphi) &= \arctan\frac{\sin\theta}{\cos\theta+\tan\varphi}+2n\pi,
\end{split}
\]
and $0<\cos\varphi+\sin\varphi\leq\sqrt{2}$ when $0\leq\varphi\leq\pi/2$, then we have the simple relation
\begin{equation}\label{eq:chap_3_3_3}
e^{\log s}=s\texttt{ }(s\neq0),
\end{equation}
which serves to motivate expression~(\ref{eq:chap_3_3_2}) as the definition of the (multiple-valued)
logarithmic function of a nonzero spatial complex variable $s=r(e_{xy}e^{i\theta}\cos\varphi+\sin\varphi)$.

It should be emphasized that it is not true that the left-hand side of equation~(\ref{eq:chap_3_3_3}) with the order of the exponential and logarithmic functions reversed reduces to just $s$. More precisely, since expression~(\ref{eq:chap_3_3_2}) can be written
\[
\log s=\ln|s|+\{
    \begin{array}{cc}
        \log(e_{xy}e^{i\theta}\cos\varphi+\sin\varphi) & \texttt{ when }\sin\varphi>0 \\
        e_{xy}i\theta & \texttt{ when }\sin\varphi=0
    \end{array}
\]
\[
\texttt{ }(-\pi<\theta<\pi,\texttt{ }0\leq\varphi\leq\pi/2)
\]
and since (Sec.~(\ref{sec:chap_1_6})) let $\arg s$ denote a pair of arguments $(\arg_{c}s,\arg_{r}s)$, that is, $\arg s=(\arg_{c}s,\arg_{r}s)$ where $\arg_{c}s=\theta$ and $\arg_{r}s=\varphi$. The principal values of $\arg_{c}s$ and $\arg_{r}s$, denoted by $Arg_{c}s$ and $Arg_{r}s$, are those unique values $\Theta$ and $\Phi$ respectively, such that $-\pi<\Theta<\pi$ and $0\leq\Phi\leq\pi/2$. Note that
\begin{equation}\label{eq:chap_3_3_4}
\arg_{r}s=Arg_{r}s,\texttt{ }\arg_{c}s=Arg_{c}s+2n\pi\texttt{ }(n=0,\pm1,\pm2,\ldots)
\end{equation}
where $Arg_{r}s=\Phi=\arcsin(z/r)$ and $Arg_{c}s=\Theta=\arctan(y/x)$.

The principal value of $\log s$ is the value obtained from equation~(\ref{eq:chap_3_3_2}) when $n=0$ there and is denoted by Log $s$. Thus
\begin{equation}\label{eq:chap_3_3_5}
Log\texttt{ }s=\ln|s|+\{
    \begin{array}{cc}
        \log(e_{xy}e^{i\Theta}\cos\Phi+\sin\Phi) & \texttt{ when }\sin\Phi>0 \\
        e_{xy}i\Theta & \texttt{ when }\sin\Phi=0
    \end{array}
\end{equation}
\[
\texttt{ }(0<\Theta<\pi,\texttt{ }0\leq\Phi\leq\pi/2).
\]
Note that Log $s$ is well defined and single-valued when $s\neq0$ and that generally equation~(\ref{eq:chap_3_3_2}) is used.

It reduces to the usual logarithm in two-dimensional complex variables when $s$ is a two-dimensional complex number $s=e_{xy}(x+iy)$.  To see this, one need only write $s=re_{xy}e^{i\theta}$, in which case equation~(\ref{eq:chap_3_3_5}) becomes
\[
\begin{split}
Log\texttt{ }s &= \ln|s|+\log(e_{xy}e^{i\Theta}) \\
               &= \ln|s|+e_{xy}i\Theta\texttt{ }(0<\Theta<\pi)\texttt{ for }\Phi=0.
\end{split}
\]

\begin{lemma}\label{le:chap_3_3_1}
A spatial complex function $\log(e_{xy}e^{i\theta}\cos\varphi+\sin\varphi)$ is analytic and
\begin{equation}\label{eq:chap_3_3_6}
\log(e_{xy}e^{i\theta}\cos\varphi+\sin\varphi)=e_{xy}[\log(e^{i\theta}\cos\varphi+\sin\varphi)-\log\sin\varphi]+\log\sin\varphi
\end{equation}
\[
=e_{xy}[\log r_{L}(\theta,\varphi)+i\theta_{L}(\theta,\varphi)-\log\sin\varphi]+\log\sin\varphi
\]
for $0<\varphi\leq\pi/2$ and $\theta\neq(2n+1)\pi$ where $n=0,1,2,\ldots$
\[
\begin{split}
r_{L}(\theta,\varphi) &= \sqrt{1+\cos\theta\sin2\varphi},\texttt{ } \\
\theta_{L}(\theta,\varphi) &= \arctan\frac{\sin\theta}{\cos\theta+\tan\varphi}.
\end{split}
\]
\end{lemma}

\begin{proof}
From expression~(\ref{eq:chap_1_3_2}) in Sec.~(\ref{sec:chap_1_3}), there is
\[
\frac{s_{1}}{s_{2}}=e_{xy}(\frac{x_{1}+iy_{1}+z_{1}}{x_{2}+iy_{2}+z_{2}}-\frac{z_{1}}{z_{2}})+\frac{z_{1}}{z_{2}}\texttt{ }(z_{2}\neq0,\texttt{ }x_{2}+z_{2}\neq0\texttt{ or }y_{2}\neq0).
\]

Let $s_{1}=e_{xy}e^{i\theta}\cos\varphi+0$ and $s_{2}=e_{xy}e^{i\theta}\cos\varphi+2\sin\varphi$. Then there are
\[
\begin{split}
\frac{s_{1}}{s_{2}} &= \frac{e_{xy}e^{i\theta}\cos\varphi+0}{e_{xy}e^{i\theta}\cos\varphi+2\sin\varphi} \\
                    &= e_{xy}\frac{e^{i\theta}\cos\varphi}{e^{i\theta}\cos\varphi+2\sin\varphi},
\end{split}
\]

We know that the power series of a real function $\log(a+x)$ $(a>0,a+x>0)$ is
\[
\log(a+x)=\log a+2\sum_{n=0}^{\infty}\frac{1}{2n+1}(\frac{x}{2a+x})^{2n+1}\texttt{ }(a>0,-a<x<\infty).
\]

Let $a=\sin\varphi$ and $x=e_{xy}e^{i\theta}\cos\varphi$. Then when $0<\varphi\leq\pi/2$ and $\theta\neq(2n+1)\pi$ for $n=0,1,2,\ldots$, there are $a>0$ and $a+x\neq0$. Thus formally we have
\[
\begin{split}
\log(\sin\varphi+e_{xy}e^{i\theta}\cos\varphi)&=\log\sin\varphi+\sum_{n=0}^{\infty}\frac{2}{2n+1}(\frac{e_{xy}e^{i\theta}\cos\varphi}{2\sin\varphi+e_{xy}e^{i\theta}\cos\varphi})^{2n+1}\\
&=\log\sin\varphi+e_{xy}\sum_{n=0}^{\infty}\frac{2}{2n+1}(\frac{e^{i\theta}\cos\varphi}{2\sin\varphi+e^{i\theta}\cos\varphi})^{2n+1}\\
&=\log\sin\varphi+e_{xy}[\log(\sin\varphi+e^{i\theta}\cos\varphi)-\log\sin\varphi].
\end{split}
\]

When $0<\varphi\leq\pi/2$ and $\theta\neq(2n+1)\pi$ for $n=0,1,2,\ldots$, each term on the right-hand side of the equalities above is analytic, so the term on the left-hand side of the equalities above is also analytic.

This completes the proof of the lemma.
\end{proof}

\section{Branches and Derivatives of Logarithms}\label{sec:chap_3_4_30}

If $s=r(e_{xy}e^{i\theta}\cos\varphi+\sin\varphi)$ is a nonzero spatial complex number, each of the arguments $\theta$ and $\varphi$ has any one of the values $\theta=\Theta+2n\pi$ or $\varphi=\Phi+2n\pi$ $(n=0,\pm1,\pm2,\ldots)$, where $\Theta=Arg_{c}$ $s$ and $\Phi=Arg_{r}$ $s$. Hence the definition
\[
\log s=\log|s|+\log[e_{xy}e^{i(\Theta+2n\pi)}\cos(\Phi+2n\pi)+\sin(\Phi+2n\pi)]
\]
\[
\texttt{ }(n=0,\pm1,\pm2,\ldots)
\]
of the multiple-valued logarithmic function in Sec.~(\ref{sec:chap_3_3_29}) can be written
\begin{equation}\label{eq:chap_3_4_1}
\log s=\log|s|+\log(e_{xy}e^{i\theta}\cos\varphi+\sin\varphi).
\end{equation}

If we let $\alpha$ denote any real number and restrict the value of $\theta$ in expression~(\ref{eq:chap_3_4_1}) so that $\alpha<\theta<\alpha+2n\pi$ and $0<\varphi\leq\pi/2$, the function
\begin{equation}\label{eq:chap_3_4_2}
\begin{split}
\log s &= \log|s|+\log(e_{xy}e^{i\theta}\cos\varphi+\sin\varphi) \\
       &= \log(r\sin\varphi)+\log(e_{xy}e^{i\theta}\cot\varphi+1)
\end{split}
\end{equation}
where $r=|s|>0$, $\alpha<\theta<\alpha+2n\pi$, and $0<\varphi\leq\pi/2$, with components
%\begin{equation}\label{eq:chap_3_4_3}
\[
w(r,\varphi)=\log(r\sin\varphi)\texttt{ and }
e_{xy}[u(\theta,\varphi)+iv(\theta,\varphi)]=\log(e_{xy}e^{i\theta}\cot\varphi+1),
\]
%\end{equation}
is single-valued and continuous in the stated domain. Note that if the function~(\ref{eq:chap_3_4_2}) were to be defined on the ray $\theta=\alpha$, it would not be continuous there. For, if $s$ is a point on that ray, there are points arbitrarily close to $s$ at which the values of $v$ are near $\alpha$ and also points such that the values of $v$ are near $\alpha+2\pi$.

The function~(\ref{eq:chap_3_4_2}) is not only continuous but also analytic in the domain $r>0$,
$\alpha<\theta<\alpha+2n\pi$, and $0<\varphi\leq\pi/2$ since the first-order partial derivatives of $u$, $v$, and $w$ are continuous there and satisfy the spherical form (Sec.~(\ref{sec:chap_2_12_22}))
\[
ru_{r}=v_{\theta}+u_{z}r\sin\varphi\texttt{ and }u_{\theta}=-rv_{r}+v_{z}r\sin\varphi
\]
of the Cauchy-Riemann equations. Furthermore, according to Sec.~(\ref{sec:chap_2_12_22}),
\[
\begin{split}
\frac{d}{ds}\log s &= \frac{e_{xy}(u_{r}+iv_{r})+w_{r}}{e_{xy}e^{i\theta}\cos\varphi+\sin\varphi}
                    = \frac{e_{xy}(0+i0)+1/r}{e_{xy}e^{i\theta}\cos\varphi+\sin\varphi} \\
                   &= \frac{1}{r(e_{xy}e^{i\theta}\cos\varphi+\sin\varphi)}=\frac{1}{s},
\end{split}
\]
that is,
\begin{equation}\label{eq:chap_3_4_4}
\frac{d}{ds}\log s=\frac{1}{s}\texttt{ }(|s|>0,\alpha<\theta<\alpha+2n\pi,0<\varphi\leq\pi/2).
\end{equation}
In particular,
\begin{equation}\label{eq:chap_3_4_5}
\frac{d}{ds}Log\texttt{ }s=\frac{1}{s}\texttt{ }(|s|>0,-\pi<\Theta<\pi,0<\Phi\leq\pi/2).
\end{equation}
A branch of a multiple-valued function $f$ is any single-valued function $F$ that is analytic in some domain at  each point $s$ of which the value $F(s)$ is one of the values $f(s)$. The requirement of analyticity, of course, prevents $F$ from taking on a random selection of the values of $f$. Observe that, for each fixed $\alpha$, the single-valued function~(\ref{eq:chap_3_4_2}) is a branch of the multiple-valued function~(\ref{eq:chap_3_4_1}). The function
\begin{equation}\label{eq:chap_3_4_6}
Log\texttt{ }s=\log|s|+\log(e_{xy}e^{i\Theta}\cos\Phi+\sin\Phi)\texttt{ }(0<\Theta<\pi,\texttt{ }0<\Phi\leq\pi/2)
\end{equation}
is called the principal branch.

A branch cut is a portion of a line or curve that is introduced in order to define a branch $F$ of a multiple-valued function $f$. Points on the branch cut for $F$ are singular points (Sec.~(\ref{sec:chap_2_13_23}))  of $F$, and any point that is common to all branch cuts of $f$ is called a branch point. The origin and the ray $\theta=\alpha$ make up the branch cut for the branch~(\ref{eq:chap_3_4_2}) of the logarithmic function. The branch cut for the principal branch~(\ref{eq:chap_3_4_6}) consists of the origin and the rays $\Theta=\pi$ and $\Phi=\pi$. The origin is evidently a branch point for branches of the multiple-valued logarithmic function.

\section{Some Identities Involving Logarithms}\label{sec:chap_3_5_31}

As suggested by relations~(\ref{eq:chap_3_3_3}) and~(\ref{eq:chap_3_3_4}) in. Sec.~(\ref{sec:chap_3_3_29}), some identities involving logarithms in calculus carry over to spatial complex analysis and others do not. In this section, we derive a few that do carry over, sometimes with qualifications as to how they are to be interpreted. A reader who wishes to pass to Sec.~(\ref{sec:chap_3_6_32}) can simply refer to results here when needed.

If $s_{1}$ and $s_{2}$ denote any two nonzero spatial complex numbers, it is straightforward to show that
\begin{equation}\label{eq:chap_3_5_1}
\log(s_{1}s_{2})=\log s_{1}+\log s_{2}.
\end{equation}
This statement, involving a multiple-valued function, is to be interpreted in the same way that the statement
\begin{equation}\label{eq:chap_3_5_2}
\arg(s_{1}s_{2})=\arg s=(\arg_{c}s,\arg_{r}s)
\end{equation}
where $\arg_{r}s=\varphi=\arcsin(\sin\varphi_{1}\sin\varphi_{2})$ and
\[
\arg_{c}s=\theta=\arctan\frac{\sin(\theta_{1}+\theta_{2})
+\sin\theta_{1}\tan\varphi_{2}+\sin\theta_{2}\tan\varphi_{1}}
{\cos(\theta_{1}+\theta_{2})
+\cos\theta_{1}\tan\varphi_{2}+\cos\theta_{2}\tan\varphi_{1}}
\]
were in Sec.~(\ref{sec:chap_1_7}). That is, if values of two of the three logarithms are specified, then there is a value of the third logarithm such that equation~(\ref{eq:chap_3_5_1}) holds.

The proof of statement~(\ref{eq:chap_3_5_1}) can be based on statement~(\ref{eq:chap_3_5_2}) in the following way. Since $|s_{1}s_{2}|=|s_{1}||s_{2}|$ and since these moduli are all positive real numbers, we know from experience with logarithms of such numbers in calculus that
\[
\ln|s_{1}s_{2}|=\ln|s_{1}|+\ln|s_{2}|.
\]
So it follows from this and equation~(\ref{eq:chap_3_5_2}) that
\begin{equation}\label{eq:chap_3_5_3}
\log s_{1}s_{2}=\log|s_{1}s_{2}|+\log(e_{xy}e^{i\theta}\cos\varphi+\sin\varphi)
\end{equation}
where $|s_{1}s_{2}|<+\infty$ and $0<\varphi\leq\pi/2$. Finally, because of the way in which equations~(\ref{eq:chap_3_5_1}) and~(\ref{eq:chap_3_5_2}) are to be interpreted, equation~(\ref{eq:chap_3_5_3}) is the same as equation~(\ref{eq:chap_3_5_1}).

Verification of the statement
\begin{equation}\label{eq:chap_3_5_4}
\log(s_{1}/s_{2})=\log s_{1}-\log s_{2}
\end{equation}
is to be interpreted in the same way as statement~(\ref{eq:chap_3_5_1}) that the statement
\begin{equation}\label{eq:chap_3_5_4b}
\arg(s_{1}/s_{2})=\arg s=(\arg_{c}s,\arg_{r}s)
\end{equation}
where $\arg_{r}s=\varphi=\arcsin(\sin\varphi_{1}/\sin\varphi_{2})$ and
\[
\arg_{c}s=\theta=\arctan\frac{\sin(\theta_{1}-\theta_{2})+\sin\theta_{1}\tan\varphi_{2}-\sin\theta_{2}\tan\varphi_{1}}
{\cos(\theta_{1}-\theta_{2})+\cos\theta_{1}\tan\varphi_{2}-\cos\theta_{2}\tan\varphi_{1}-\tan\varphi_{1}\cot\varphi_{2}}.
\]
were in Sec.~(\ref{sec:chap_1_7}). That is, if values of two of the three logarithms are specified, then there is a value of the third logarithm such that equation~(\ref{eq:chap_3_5_4}) holds.

We include here two other properties of $\log s$ that will be of special interest in
Sec.~(\ref{sec:chap_3_6_32}). If $s$ is a nonzero spatial complex number, then
\begin{equation}\label{eq:chap_3_5_5}
s^{n}=e^{n\log s}\texttt{ }(n=0,\pm1,\pm2,\ldots)
\end{equation}
for any value of $\log s$ that is taken. When $n=1$, this reduces, of course, to relation~(\ref{eq:chap_3_3_3}),
Sec.~(\ref{sec:chap_3_3_29}). Equation~(\ref{eq:chap_3_5_5}) is readily verified by writing $s=r(e_{xy}e^{i\theta}\cos\varphi+\sin\varphi)$ and noting that each side becomes $r^{n}(e_{xy}e^{i\theta}\cos\varphi+\sin\varphi)^{n}$.

It is also true that when $s\neq0$,
\begin{equation}\label{eq:chap_3_5_6}
s^{1/n}=\exp(\frac{1}{n}\log s)\texttt{ }(n=0,1,2,\ldots).
\end{equation}
That is, the term on the right here has $n$ distinct values, and those values are the $n$th
roots of $s$. To prove this, we write $s=r(e_{xy}e^{i(\Theta+2k\pi)/n}\cos\varphi+\sin\varphi)$ where $\Theta$ is the principal value of
$\arg_{c}s$. Then, in view of definition~(\ref{eq:chap_3_3_2}), Sec.~(\ref{sec:chap_3_3_29}), of $\log s$,
\[
\exp(\frac{1}{n}\log s)=\exp(\frac{1}{n}\log|s|+\frac{1}{n}\log(e_{xy}e^{i\theta}\cos\varphi+\sin\varphi))
\]
where $k=0,\pm1,\pm2,\ldots$. Thus
\begin{equation}\label{eq:chap_3_5_7}
\exp(\frac{1}{n}\log s)=\sqrt[n]{r}(e_{xy}e^{i\theta_{b}(k)}\cos\varphi_{n}+\sin\varphi_{n})\texttt{ }(n=2,3,\ldots)
\end{equation}
where $\varphi_{n}=\arcsin(\sqrt[n]{\sin\varphi})$,
\[
\begin{split}
r_{a}^{2} &= 1+\cos\theta\sin2\varphi,\texttt{ }\theta_{a}=\arctan\frac{\sin\theta}{\cos\theta+tan\varphi}, \\
r_{b}^{2}(k) &= 1-2\cos\frac{\theta_{a}+2k\pi}{n}\sqrt[n]{\frac{\sin\varphi}{r_{a}}}+(\sqrt[n]{\frac{\sin\varphi}{r_{a}}})^{2}, \\
\theta_{b}(k) &= \arctan\frac{\sin(\theta_{a}/n+2k\pi/n)}{\cos(\theta_{a}/n+2k\pi/n)-\sqrt[n]{(\sin\varphi)/r_{a}}}\texttt{ }(k=0,\pm1,\pm2,\ldots).
\end{split}
\]

Because $\exp(i2k\pi/n)$ has distinct values only when $k=0,1,2,\ldots,n-1$, the right-hand side of equation~(\ref{eq:chap_3_5_7}) has only $n$ values. That right-hand side is, in fact, an expression for the $n$th  roots of $s$ (Sec.~(\ref{sec:chap_1_8})), and so it can be written $s^{1/n}$. This establishes
property~(\ref{eq:chap_3_5_6}), which is actually valid when $n$ is a negative integer too.

\section{Complex Exponent}\label{sec:chap_3_6_32}

When $s\neq0$ and the exponent $c$ is any spatial complex number, the function $s^{c}$ is defined by means of the equation
\begin{equation}\label{eq:chap_3_6_1}
s^{c}=e^{c\log s}=e_{xy}[(1+e^{i\theta}\cot\varphi)^{c}-1]+1+(|s|\sin\varphi)^{c}
\end{equation}
where $|s|<+\infty$, $0<\varphi\leq\pi/2$, and $\log$ $s$ denotes the multiple-valued logarithmic function. Equation~(\ref{eq:chap_3_6_1}) provides a consistent definition of $s^{c}$ in the sense that it is already known  to be valid (see Sec.~(\ref{sec:chap_3_5_31})) when $c=n$ $(n=0,\pm1,\pm2,\ldots)$ and $c=1/n$ $(n=0,\pm1,\pm2,\ldots)$. Definition~(\ref{eq:chap_3_6_1}) is, in fact, suggested by those particular choices of $c$.

Since the exponential function has the property $1/e^{s}=e^{-s}$, one can see that
\[
\frac{1}{s^{c}}=\frac{1}{\exp(c\log s)}=\exp(-c\log s)=s^{-c}.
\]

If $s=r(e_{xy}e^{i\theta}\cos\varphi+\sin\varphi)$, and $\alpha$ is any real number, the branch
\[
\log s=\log r+\log\sin\varphi+e_{xy}[\log(e^{i\theta}\cos\varphi+\sin\varphi)-\log\sin\varphi]
\]
\[
(r>0,\alpha<\theta<\alpha+2n\pi,0<\varphi\leq\pi/2)
\]
of the logarithmic function is single-valued and analytic in the indicated domain (Sec.~(\ref{sec:chap_3_4_30})).  When that branch is used, it follows that the function $s^{c}=\exp(c\log s)$ is single-valued and analytic in the same domain. The derivative of such a branch of $s^{c}$ is found by first using the chain rule to write
\[
\frac{d}{ds}s^{c}=\frac{d}{ds}\exp(c\log s)=\frac{c}{s}\exp(c\log s)
\]
and then recalling (Sec.~(\ref{sec:chap_3_3_29})) the identity $s=\exp(\log s)$. That yields the result
\[
\frac{d}{ds}s^{c}=c\frac{\exp(c\log s)}{\exp(\log s)}=c\exp[(c-1)\log s],
\]
or
\begin{equation}\label{eq:chap_3_6_4}
\frac{d}{ds}s^{c}=cs^{c-1}\texttt{ }(|s|>0,\alpha<\theta<\alpha+2n\pi,0<\varphi\leq\pi/2).
\end{equation}

The principal value of $s^{c}$ occurs when $\log s$ is replaced by $Log$ $s$ in definition~(\ref{eq:chap_3_6_1}):
\begin{equation}\label{eq:chap_3_6_5}
P.V.s^{c}=e^{cLog\texttt{ }s}.
\end{equation}
Equation~(\ref{eq:chap_3_6_5}) also serves to define the principal branch of the function $s^{c}$ on the domain
$|s|>0$, $-\pi<Arg_{c}s\leq\pi$, and $0<Arg_{r}s\leq\pi/2$.

According to definition~(\ref{eq:chap_3_6_1}), the exponential function with base $c$, where $c$ is any nonzero spatial complex constant, is written
\begin{equation}\label{eq:chap_3_6_8}
c^{s}=e^{s\log c}.
\end{equation}
Note that although $e^{s}$ is, in general, multiple-valued according to definition~(\ref{eq:chap_3_6_8}), the
usual interpretation of $e^{s}$ occurs when the principal value of the logarithm is taken. For the principal value  of $\log e$ is unity.

When a value of $\log c$ is specified, $c^{s}$ is an entire function of $s$. In fact,
\[
\frac{d}{ds}c^{s}=\frac{d}{ds}e^{s\log c}=e^{s\log c}\log c;
\]
and this shows that
\[
\frac{d}{ds}c^{s}=c^{s}\log c.
\]

\section{Trigonometric Functions}\label{sec:chap_3_7_33}

\begin{definition}\label{def:chap_3_7_1}
The spatial complex sine function $\sin s$ is defined as
\begin{equation}\label{eq:chap_3_7_01}
\sin s=e_{xy}[\sin(x+iy+z)-\sin z]+\sin z\texttt{ }(|s|<+\infty).
\end{equation}
\end{definition}

\begin{definition}\label{def:chap_3_7_2}
The spatial complex cosine function $\cos s$ is defined as
\begin{equation}\label{eq:chap_3_7_02}
\cos s=\{
    \begin{array}{cc}
        e_{xy}[\cos(x+iy+z)-\cos z]+\cos z & \texttt{ when }z\neq0 \\
        e_{xy}\cos(x+iy) & \texttt{ when }z=0.
    \end{array}
\end{equation}
\end{definition}

\begin{definition}\label{def:chap_3_7_3}
The spatial complex tan function $\tan s$ is defined as
\begin{equation}\label{eq:chap_3_7_03}
\tan s=e_{xy}[\tan(x+iy+z)-\tan z]+\tan z\texttt{ }(|s|<\pi/2).
\end{equation}
\end{definition}

First, we shall prove equation~(\ref{eq:chap_3_7_01}).

Let $\sigma=x+iy$. Then $\sigma$ and $\sigma+z$ are two-dimensional complex numbers. The Taylor series of the real sine function $\sin z$ and two-dimensional complex sine function $\sin(\sigma+z)$ are
\[
\sin z=\sum_{n=0}^{\infty}\frac{(-1)^{n}z^{2n+1}}{(2n+1)!}\texttt{ and }
\sin(\sigma+z)=\sum_{n=0}^{\infty}\frac{(-1)^{n}(\sigma+z)^{2n+1}}{(2n+1)!}\texttt{ }(|\sigma+z|<\infty).
\]

For an index $n$ in the Taylor series of $\sin z$ and $\sin(\sigma+z)$, there are $z^{2n+1}$ and
\[
(\sigma+z)^{2n+1}=\sigma^{2n+1}+(2n+1)\sigma^{2n}z+\cdots+(2n+1)\sigma z^{2n}+z^{2n+1}.
\]

Since both $\sigma+z$ and its powers are points in the two-dimensional complex plane and from remark~(\ref{re:chap_1_1_1}) for a term $e_{xy}x^{j}y^{k}$, the operator $e_{xy}$ becomes $e_{xy}^{m}$ where $m=|j|+|k|$, we get
\[
e_{xy}[(\sigma+z)^{2n+1}-z^{2n+1}]+z^{2n+1}
\]
\[
\begin{split}
&=(e_{xy}\sigma)^{2n+1}+(2n+1)(e_{xy}\sigma)^{2n}z+\cdots+(2n+1)(e_{xy}\sigma)z^{2n}+z^{2n+1} \\
&=(e_{xy}\sigma+z)^{2n+1}=s^{2n+1}.
\end{split}
\]

Thus there are
\[
e_{xy}[\sin(\sigma+z)-\sin z]+\sin z
\]
\[
\begin{split}
&=e_{xy}[\sum_{n=0}^{\infty}\frac{(-1)^{n}(\sigma+z)^{2n+1}}{(2n+1)!}
-\sum_{n=0}^{\infty}\frac{(-1)^{n}z^{2n+1}}{(2n+1)!}]
+\sum_{n=0}^{\infty}\frac{(-1)^{n}z^{2n+1}}{(2n+1)!} \\
&=\sum_{n=0}^{\infty}\frac{(-1)^{n}(e_{xy}\sigma+z)^{2n+1}}{(2n+1)!}
=\sum_{n=0}^{\infty}\frac{(-1)^{n}s^{2n+1}}{(2n+1)!}
=\sin s.
\end{split}
\]

Next, other equations~(\ref{eq:chap_3_7_02}) and~(\ref{eq:chap_3_7_03}) can also be proved in this way. So the Taylor series of spatial complex trigonometric functions $\sin s$, $\cos s$, and $\tan s$ are
\begin{equation}\label{eq:chap_3_7_04}
\begin{split}
\sin s &= \sum_{n=0}^{\infty}\frac{(-1)^{n}s^{2n+1}}{(2n+1)!}\texttt{ }(|s|<\infty), \\
\cos s &= \sum_{n=0}^{\infty}\frac{(-1)^{n}s^{2n}}{(2n)!}\texttt{ }(|s|<\infty),
\end{split}
\end{equation}
and
\begin{equation}\label{eq:chap_3_7_05}
\tan s=\sum_{n=1}^{\infty}\frac{4^{n}(4^{n}-1)^{n}B_{n}s^{2n-1}}{(2n)!}\texttt{ }(|s|<\pi/2)
\end{equation}
where $B_{n}$ is Bernoulli number.

Euler's formula (Sec.~(\ref{sec:chap_1_6})) tells us that
\[
e^{ix}=\cos x+i\sin x\texttt{ and }e^{-ix}=\cos x-i\sin x
\]
for every real number $x$, and it follows from these equations that
\[
e^{ix}-e^{-ix}=2i\sin x\texttt{ and }e^{ix}+e^{-ix}=2\cos x.
\]
That is,
\[
\sin x=\frac{e^{ix}-e^{-ix}}{2i}\texttt{ and }\cos x=\frac{e^{ix}+e^{-ix}}{2}.
\]
It is, therefore, natural to define the sine and cosine functions of a spatial complex variable $s$ as follows:
\begin{equation}\label{eq:chap_3_7_1}
\sin s=\frac{e^{is}-e^{-is}}{2i}\texttt{ and }\cos s=\frac{e^{is}+e^{-is}}{2}.
\end{equation}

These functions are entire since they are linear combinations of the entire functions $s^{c}$. Knowing the derivatives of those exponential functions, we find from equations~(\ref{eq:chap_3_7_01}) and~(\ref{eq:chap_3_7_02}) that
\begin{equation}\label{eq:chap_3_7_2}
\frac{d}{ds}\sin s=\cos s,\texttt{ }\frac{d}{ds}\cos s=-\sin s.
\end{equation}
It is easy to see from definitions~(\ref{eq:chap_3_7_01}) that
\begin{equation}\label{eq:chap_3_7_3}
\sin(-s)=-\sin s,\texttt{ }\cos(-s)=\cos s;
\end{equation}
and a variety of other identities from trigonometry are valid with spatial complex variables.

\begin{example}\label{ex:chap_3_7_1}
In order to show that
\begin{equation}\label{eq:chap_3_7_4}
2\sin s_{1}\cos s_{2}=\sin(s_{1}+s_{2})+\sin(s_{1}-s_{2}),
\end{equation}
using definitions~(\ref{eq:chap_3_7_1}) and properties of the exponential function, we first write
\[
\begin{split}
2\sin s_{1}\cos s_{2} &=2(\frac{e^{is_{1}}-e^{-is_{1}}}{2i})(\frac{e^{is_{2}}+e^{-is_{2}}}{2}) \\
                      &=\frac{e^{i(s_{1}+s_{2})}-e^{-i(s_{1}+s_{2})}}{2i}+\frac{e^{i(s_{1}-s_{2})}-e^{-i(s_{1}-s_{2})}}{2i},
\end{split}
\]
or
\[
\sin(s_{1}+s_{2})+\sin(s_{1}-s_{2});
\]
and identity~(\ref{eq:chap_3_7_4}) is  established.
\end{example}

Identity~(\ref{eq:chap_3_7_4}) leads to the identities
\begin{equation}\label{eq:chap_3_7_5}
\sin(s_{1}+s_{2})=\sin s_{1}\cos s_{2}+\cos s_{1}\sin s_{2},
\end{equation}
\begin{equation}\label{eq:chap_3_7_6}
\cos(s_{1}+s_{2})=\cos s_{1}\cos s_{2}-\sin s_{1}\sin s_{2},
\end{equation}
and from these it follows that
\begin{equation}\label{eq:chap_3_7_7}
\sin^{2}s+\cos^{2}s=1,
\end{equation}
\begin{equation}\label{eq:chap_3_7_8}
\sin2s=2\sin s\cos s,\texttt{ }\cos2s=\cos^{2}s-\sin^{2}s,
\end{equation}
\begin{equation}\label{eq:chap_3_7_9}
\sin(s+\pi/2)=\cos s,\texttt{ }\sin(s-\pi/2)=-\cos s.
\end{equation}

When $y$ is any real number, one can use definitions~(\ref{eq:chap_3_7_1}) and the hyperbolic functions
\begin{equation}\label{eq:chap_3_7_06}
\sinh y=\frac{e^{y}-e^{-y}}{2}\texttt{ and }\cosh y=\frac{e^{y}+e^{-y}}{2}
\end{equation}
from calculus to write
\begin{equation}\label{eq:chap_3_7_10}
\sin(iy)=i\sinh y\texttt{ and }\cos(iy)=\cosh y.
\end{equation}
The real and imaginary parts of $\sin\sigma$ and $\cos\sigma$ are then readily displayed by writing $\sigma_{1}=x$ and $\sigma_{2}=iy$ in identities~(\ref{eq:chap_3_7_06}) and~(\ref{eq:chap_3_7_10}):
\begin{equation}\label{eq:chap_3_7_11}
\sin\sigma=\sin x\cosh y+i\cos x\sinh y\texttt{ }(|\sigma|<+\infty),
\end{equation}
\begin{equation}\label{eq:chap_3_7_12}
\cos\sigma=\cos x\cosh y-i\sin x\sinh y\texttt{ }(|\sigma|<+\infty).
\end{equation}
where $\sigma=x+iy$.

So definitions~(\ref{eq:chap_3_7_01}) and~(\ref{eq:chap_3_7_02}) become
\begin{equation}\label{eq:chap_3_7_07}
\sin s=e_{xy}[\sin(x+z)\cosh y+i\cos(x+z)\sinh y-\sin z]+\sin z\texttt{ }(|s|<+\infty),
\end{equation}
and if $z=0$ then $\cos s=e_{xy}(\cos x\cosh y-i\sin x\sinh y)$ else
\begin{equation}\label{eq:chap_3_7_08}
\cos s=e_{xy}[\cos(x+z)\cosh y-i\sin(x+z)\sinh y-\cos z]+\cos z\texttt{ }(|s|<+\infty).
\end{equation}

Particularly, there are
\[
\begin{split}
\sin s &= (-1)^{n}\sin[e_{xy}(x+iy)] \\
       &= (-1)^{n}\sin(e_{xy}\sigma)\texttt{ when }z=n\pi\texttt{ }(n=0,1,2,\ldots), \\
\cos s &= (-1)^{n}\cos[e_{xy}(x+iy)] \\
       &= (-1)^{n}\cos(e_{xy}\sigma)\texttt{ when }z=\pi/2+n\pi\texttt{ }(n=0,1,2,\ldots).
\end{split}
\]
Then definitions~(\ref{eq:chap_3_7_07}) and~(\ref{eq:chap_3_7_08}) become
\[
\sin(e_{xy}\sigma)=e_{xy}\sin\sigma=e_{xy}(\sin x\cosh y+i\cos x\sinh y)\texttt{ }(|\sigma|<+\infty),
\]
\[
\cos(e_{xy}\sigma)=e_{xy}\cos\sigma=e_{xy}(\cos x\cosh y-i\sin x\sinh y)\texttt{ }(|\sigma|<+\infty).
\]

A number of important properties of $\sin s$ and $\cos s$ follow immediately from
expressions~(\ref{eq:chap_3_7_07}) and~(\ref{eq:chap_3_7_08}). The periodic character of these functions, for example, is evident:
\begin{equation}\label{eq:chap_3_7_13}
\sin(s+2\pi)=\sin s,\texttt{ }\sin(s+\pi)=-\sin s\texttt{ }(|s|<+\infty),
\end{equation}
\begin{equation}\label{eq:chap_3_7_14}
\cos(s+2\pi)=\cos s,\texttt{ }\cos(s+\pi)=-\cos s\texttt{ }(|s|<+\infty).
\end{equation}
Also according to equations~(\ref{eq:chap_3_7_07}) and~(\ref{eq:chap_3_7_08}), there are
\[
\begin{split}
|\sin s|^{2} &= [\sin(x+z)\cosh y-\sin z]^{2}+[\cos(x+z)\sinh y]^{2}+\sin^{2}z \\
             &= \sin^{2}(x+z)+\sinh^{2}y-2[\sin(x+z)\cosh y-\sin z]\sin z \\
             &= \sin^{2}(x+z)+\sinh^{2}y-2Re\sin s\sin z,
\end{split}
\]
and if $z=0$ then $|\cos s|^{2}=\cos^{2}x+\sinh^{2}y$ else
\[
\begin{split}
|\cos s|^{2} &= [\cos(x+z)\cosh y-\cos z]^{2}+[\sin(x+z)\sinh y]^{2}+\cos^{2}z \\
             &= \cos^{2}(x+z)+\sinh^{2}y-2(\cos(x+z)\cosh y-\cos z)\cos z \\
             &= \cos^{2}(x+z)+\sinh^{2}y-2Re\cos s\cos z.
\end{split}
\]
So we get
\begin{equation}\label{eq:chap_3_7_15}
|\sin s|^{2}=\sin^{2}(x+z)+\sinh^{2}y-2Re\sin s\sin z,
\end{equation}
\begin{equation}\label{eq:chap_3_7_16}
|\cos s|^{2}=\{
    \begin{array}{cc}
        \cos^{2}(x+z)+\sinh^{2}y-2Re\cos s\cos z & \texttt{ when }z\neq0 \\
        \cos^{2}x+\sinh^{2}y & \texttt{ when }z=0
    \end{array}.
\end{equation}
Inasmuch as $\sinh y$ tends to infinity as $y$ tends to infinity, it is clear from these two equations that $\sin s$ and $\cos s$ are not bounded on the complex space, whereas the absolute values of $\sin(x+z)$ and $\cos(x+z)$ are less than or equal to unity for all values of $x+z$. (See the definition of boundedness at the end of Sec.~(\ref{sec:chap_2_7_17}).)

A zero of a given function $f(s)$ is a number $s_{0}$ such that $f(s_{0})=0$. Since $\sin s$ becomes the usual sine function in calculus when $s$ is real, we know that the real $xz$ coordinate plane numbers $s=e_{xy}m\pi+n\pi$ $(m, n=0, \pm1, \pm2, \ldots)$ are all zeros of $\sin s$. To show that there are no other zeros, we assume that $\sin s=0$ so that there are $Re\sin s=0$ and $Re_{s}\sin s=\sin z=0$, and note how it follows from equation~(\ref{eq:chap_3_7_15}) that
\[
\sin^{2}(x+z)+\sinh^{2}y=2Re\sin s\sin z=0.
\]
Thus
\[
\sin(x+z)=0,\texttt{ }\sinh y=0,\texttt{ and }\sin z=0.
\]
Evidently, then, $x=m\pi$, $y=0$, and $z=n\pi$ $(m, n=0, \pm1, \pm2, \ldots)$; that is, if and only if
\begin{equation}\label{eq:chap_3_7_17}
s=e_{xy}m\pi+n\pi\texttt{ }(m, n=0, \pm1, \pm2, \ldots)\texttt{ then }\sin s=0.
\end{equation}

Since $\cos s=-\sin(s-\pi/2)$, according to the second of identities~(\ref{eq:chap_3_7_9}), if and only if
\begin{equation}\label{eq:chap_3_7_18}
s=e_{xy}(\pi/2+m\pi)+\pi/2+n\pi\texttt{ }(m. n=0, \pm1, \pm2, \ldots)\texttt{ then }\cos s=0.
\end{equation}
So, as was the case with $\sin s$, the zeros of $\cos s$ are all pints $s=e_{xy}x+z=0$ on the $xz$ coordinate plane.

The other four trigonometric functions are defined in terms of the sine and cosine functions by the usual relations:
\begin{equation}\label{eq:chap_3_7_19}
\tan s=\frac{\sin s}{\cos s},\texttt{ }\cot s=\frac{\cos s}{\sin s},
\end{equation}
\begin{equation}\label{eq:chap_3_7_20}
\sec s=\frac{1}{\cos s},\texttt{ }\csc s=\frac{1}{\sin s}.
\end{equation}
Observe that the quotients $\tan s$ and $\sec s$ are analytic everywhere except at the
singularities (Sec.~(\ref{sec:chap_2_13_23}))
\[
s=e_{xy}(\pi/2+m\pi)+\pi/2+n\pi\texttt{ }(m, n=0,\pm1,\pm2,\ldots),
\]
which are the zeros of $\cos s$. Likewise, $\cot s$ and $\sec s$ have singularities at the zeros
of $\sin s$, namely
\[
s=e_{xy}m\pi+n\pi\texttt{ }(m, n=0,\pm1,\pm2,\ldots).
\]

By differentiating the right-hand sides of equations~(\ref{eq:chap_3_7_19}) and~(\ref{eq:chap_3_7_20}), we obtain the expected differentiation formulas
\begin{equation}\label{eq:chap_3_7_21}
\frac{d}{ds}\tan s=\sec^{2}s,\texttt{ }\frac{d}{ds}\cot s=\csc^{2}s,
\end{equation}
\begin{equation}\label{eq:chap_3_7_22}
\frac{d}{ds}\sec s=\sec s\tan s,\texttt{ }\frac{d}{ds}\csc s=-\tan s\cot s.
\end{equation}
The periodicity of each of the trigonometric functions defined by equations~(\ref{eq:chap_3_7_19}) and~(\ref{eq:chap_3_7_20}) follows readily from equations~(\ref{eq:chap_3_7_13}) and~(\ref{eq:chap_3_7_14}). For example,
\begin{equation}\label{eq:chap_3_7_23}
\tan(s+\pi)=\tan s.
\end{equation}

Mapping properties of the transformation $\varpi=\sin s$ are especially important in the applications later on. A  reader who wishes at this time to learn some of those properties is sufficiently prepared to read Sec.~(\ref{sec:chap_8_7_89}) (Chap.~(\ref{ch:chap_8})), where they are discussed.

\section{Hyperbolic Functions}\label{sec:chap_3_8_34}

\begin{definition}\label{def:chap_3_8_1}
The spatial complex hyperbolic sine function $\sinh s$ is defined as
\begin{equation}\label{eq:chap_3_8_1}
\sinh s=e_{xy}[\sinh(x+iy+z)-\sinh z]+\sinh z\texttt{ }(|s|<+\infty).
\end{equation}
\end{definition}

\begin{definition}\label{def:chap_3_8_2}
The spatial complex hyperbolic cosine function $\cosh s$ is defined as if $z=0$ then $\cosh s=e_{xy}\cosh(x+iy)$ $\texttt{ }(|s|<+\infty)$ else
\begin{equation}\label{eq:chap_3_8_2}
\cosh s=e_{xy}[\cosh(x+iy+z)-\cosh z]+\cosh z\texttt{ }(|s|<+\infty).
\end{equation}
\end{definition}

\begin{definition}\label{def:chap_3_8_3}
The spatial complex hyperbolic tan function $\tanh s$ is defined as
\begin{equation}\label{eq:chap_3_8_3}
\tanh s=e_{xy}[\tanh(x+iy+z)-\tanh z]+\tanh z\texttt{ }(|s|<\pi/2).
\end{equation}
\end{definition}

First, we shall prove equation~(\ref{eq:chap_3_8_1}).

Let $\sigma=x+iy$. Then $\sigma$ and $\sigma+z$ are two-dimensional complex numbers. The Taylor series of the real hyperbolic sine function $\sinh z$ and two-dimensional complex hyperbolic sine function $\sinh(\sigma+z)$ are
\[
\sinh z=\sum_{n=0}^{\infty}\frac{z^{2n+1}}{(2n+1)!}\texttt{ and }
\sinh(\sigma+z)=\sum_{n=0}^{\infty}\frac{(\sigma+z)^{2n+1}}{(2n+1)!}\texttt{ }(|\sigma+z|<\infty).
\]

For an index $n$ in the Taylor series of $\sinh z$ and $\sinh(\sigma+z)$, there are $z^{2n+1}$ and
\[
(\sigma+z)^{2n+1}=\sigma^{2n+1}+(2n+1)\sigma^{2n}z+\cdots+(2n+1)\sigma z^{2n}+z^{2n+1}.
\]

Since both $\sigma+z$ and its powers are points in the two-dimensional complex plane and from remark~(\ref{re:chap_1_1_1}) for a term $e_{xy}x^{j}y^{k}$, the operator $e_{xy}$ becomes $e_{xy}^{m}$ where $m=|j|+|k|$, we get
\[
e_{xy}[(\sigma+z)^{2n+1}-z^{2n+1}]+z^{2n+1}
\]
\[
\begin{split}
&=(e_{xy}\sigma)^{2n+1}+(2n+1)(e_{xy}\sigma)^{2n}z+\cdots+(2n+1)(e_{xy}\sigma)z^{2n}+z^{2n+1} \\
&=(e_{xy}\sigma+z)^{2n+1}=s^{2n+1}.
\end{split}
\]

Thus there are
\[
e_{xy}[\sinh(\sigma+z)-\sinh z]+\sinh z
\]
\[
\begin{split}
&=e_{xy}(\sum_{n=0}^{\infty}\frac{(\sigma+z)^{2n+1}}{(2n+1)!}
-\sum_{n=0}^{\infty}\frac{z^{2n+1}}{(2n+1)!})
+\sum_{n=0}^{\infty}\frac{z^{2n+1}}{(2n+1)!} \\
&=\sum_{n=0}^{\infty}\frac{(e_{xy}\sigma+z)^{2n+1}}{(2n+1)!}
=\sum_{n=0}^{\infty}\frac{s^{2n+1}}{(2n+1)!}
=\sinh s.
\end{split}
\]

Next, other equations~(\ref{eq:chap_3_8_2}) and~(\ref{eq:chap_3_8_3}) can also be proved in this way. So the Taylor series of spatial complex trigonometric functions $\sinh s$, $\cosh s$, and $\tanh s$ are
\begin{equation}\label{eq:chap_3_8_4}
\begin{split}
\sinh s &= \sum_{n=0}^{\infty}\frac{s^{2n+1}}{(2n+1)!}\texttt{ }(|s|<\infty), \\
\cosh s &= \sum_{n=0}^{\infty}\frac{s^{2n}}{(2n)!}\texttt{ }(|s|<\infty),
\end{split}
\end{equation}
and
\begin{equation}\label{eq:chap_3_8_5}
\tanh s=\sum_{n=1}^{\infty}\frac{(-1)^{n}4^{n}(4^{n}-1)^{n}B_{n}s^{2n-1}}{(2n)!}\texttt{ }(|s|<\pi/2)
\end{equation}
where $B_{n}$ is Bernoulli number.

On the other hand, the hyperbolic sine and the hyperbolic cosine of a spatial complex variable are defined as
they are with a real variable or a two-dimensional complex variable; that is,
\begin{equation}\label{eq:chap_3_8_6}
\sinh s=\frac{e^{s}-e^{-s}}{2},\texttt{ }\cosh s=\frac{e^{s}+e^{-s}}{2}.
\end{equation}
Since $e^{s}$ and $e^{-s}$ are entire, it follows from definitions~(\ref{eq:chap_3_8_6}) that $\sinh s$ and $\cosh  s$ are entire. Furthermore,
\begin{equation}\label{eq:chap_3_8_7}
\frac{d}{ds}\sinh s=\cosh s,\texttt{ }\frac{d}{ds}\cosh s=\sinh s.
\end{equation}

Because of the way in which the exponential function appears in definitions~(\ref{eq:chap_3_8_6}) and in the definitions (Sec.~(\ref{sec:chap_3_7_33}))
\[
\sin s=\frac{e^{is}-e^{-is}}{2i},\texttt{ }\cos s=\frac{e^{is}+e^{-is}}{2}
\]
of $\sin s$ and $\cos s$, the hyperbolic sine and cosine functions are closely related to those trigonometric functions:
\begin{equation}\label{eq:chap_3_8_8}
-i\sinh is=\sin s,\texttt{ }\cosh is=\cos s,
\end{equation}
\begin{equation}\label{eq:chap_3_8_9}
-i\sin is=\sinh s,\texttt{ }\cos is=\cosh s.
\end{equation}

Some of the most frequently used identities involving hyperbolic sine and cosine
functions are
\begin{equation}\label{eq:chap_3_8_10}
\sinh(-s)=-\sinh s,\texttt{ }\cosh(-s)=\cosh s,
\end{equation}
\begin{equation}\label{eq:chap_3_8_11}
\cosh^{2}s-\sinh^{2}s=1,
\end{equation}
\begin{equation}\label{eq:chap_3_8_12}
\sinh(s_{1}+s_{2})=\sinh s_{1}\cosh s_{2}+\cosh s_{1}\sinh s_{2},
\end{equation}
\begin{equation}\label{eq:chap_3_8_13}
\cosh(s_{1}+s_{2})=\cosh s_{1}\cosh s_{2}+\sinh s_{1}\sinh s_{2},
\end{equation}
and from equations~(\ref{eq:chap_3_8_1}) and~(\ref{eq:chap_3_8_2})
\begin{equation}\label{eq:chap_3_8_14}
\sinh s=e_{xy}[\sinh(x+z)\cos y+i\cosh(x+z)\sin y-\sinh z]+\sinh z,
\end{equation}
\begin{equation}\label{eq:chap_3_8_15}
\cosh s=e_{xy}[\cosh(x+z)\cos y+i\sinh(x+z)\sin y-\cosh z]+\cosh z,
\end{equation}
\begin{equation}\label{eq:chap_3_8_16}
|\sinh s|^{2}=\sinh^{2}(x+z)+\sin^{2}y-2Re\sinh s\sinh z,
\end{equation}
\begin{equation}\label{eq:chap_3_8_17}
|\cosh s|^{2}=\{
    \begin{array}{cc}
        \sinh^{2}(x+z)+\cos^{2}y-2Re\cosh s\cosh z & \texttt{ when }z\neq0 \\
        \sinh^{2}x+\cos^{2}y & \texttt{ when }z=0
    \end{array}
\end{equation}
where $s=e_{xy}(x+iy)+z$. While these identities follow directly from definitions~(\ref{def:chap_3_8_1}),~(\ref{def:chap_3_8_2}), and~(\ref{eq:chap_3_8_6}). they
are often more easily obtained from related trigonometric identities, with the aid of
relations~(\ref{eq:chap_3_8_8}) and~(\ref{eq:chap_3_8_9}).

\begin{example}\label{ex:chap_3_8_1}
To illustrate the method of proof just suggested, let us verify identity~(\ref{eq:chap_3_8_16}) and~(\ref{eq:chap_3_8_17}). According to equations~(\ref{eq:chap_3_8_14}) and~(\ref{eq:chap_3_8_15}), there are
\[
\begin{split}
|\sinh s|^{2} &= [\sinh(x+z)\cos y-\sinh z]^{2}+[\cosh(x+z)\sin y]^{2}+\sinh^{2}z \\
              &= \sinh^{2}(x+z)+\sin^{2}y-2[\sinh(x+z)\cos y-\sinh z]\sinh z \\
              &= \sinh^{2}(x+z)+\sin^{2}y-2Re\sinh s\sinh z,
\end{split}
\]
and if $z=0$ then $|\cosh s|^{2}=\sinh^{2}x+\cos^{2}y$ else
\[
\begin{split}
|\cosh s|^{2} &= [\cosh(x+z)\cos y-\cosh z]^{2}+[\sinh(x+z)\sin y]^{2}+\cosh^{2}z \\
              &= \sinh^{2}(x+z)+\cos^{2}y-2[\cosh(x+z)\cos y-\cosh z]\cosh z \\
              &= \sinh^{2}(x+z)+\cos^{2}y-2Re\cosh s\cosh z,
\end{split}
\]
where $s=e_{xy}(x+iy)+z$. So we get equations~(\ref{eq:chap_3_8_16}) and~(\ref{eq:chap_3_8_17}).
\end{example}

In view of the periodicity of $\sin s$ and $\cos s$, it follows immediately from relations~(\ref{eq:chap_3_8_9}) that $\sinh s$ and $\cosh s$ are periodic with period $2\pi i$. Relations~(\ref{eq:chap_3_8_9}) also reveal that
\begin{equation}\label{eq:chap_3_8_18}
\sinh s=0\texttt{ if and only if }s=n\pi i\texttt{ }(n=0,\pm1,\pm2,\ldots)
\end{equation}
and
\begin{equation}\label{eq:chap_3_8_18}
\cosh s=0\texttt{ if and only if }s=(\pi/2+n\pi)i\texttt{ }(n=0,\pm1,\pm2,\ldots).
\end{equation}

The hyperbolic tangent of $s$ is defined by the equation
\begin{equation}\label{eq:chap_3_8_19}
\tanh s=\frac{\sinh s}{\cosh s}
\end{equation}
and is analytic in every domain in which $\cosh s\neq0$. The functions $\coth s$, $sech$ $s$, and $csch$ $s$ are the reciprocals of $\tanh s$, $\cosh s$, and $\sinh s$, respectively. It is straightforward to verify the following differentiation formulas, which are the same as those established in calculus for the corresponding functions of a real variable and two-dimensional complex functions of a complex variable:
\begin{equation}\label{eq:chap_3_8_20}
\frac{d}{ds}\tanh s=sech^{2}s,\texttt{ }\frac{d}{ds}coth s=csch^{2}s,
\end{equation}
\begin{equation}\label{eq:chap_3_8_21}
\frac{d}{ds}sech\texttt{ }s=-sech\texttt{ }s\tanh s,\texttt{ }\frac{d}{ds}csch\texttt{ }s=csch\texttt{ }s\coth s.
\end{equation}

\section{Inverse Trigonometric and Hyperbolic Functions}\label{sec:chap_3_9_35}

%We use the common notation for all of these inverse functions is $\arcsin s$, etc. Also, we remark that the common alternative notation for all of these inverse functions is $\arcsin s$, etc.

\subsection{Inverse Trigonometric Functions}\label{sec:chap_3_9_35_1}

Inverses of the trigonometric functions can be defined as follows.

\begin{definition}\label{def:chap_3_9_1}
The spatial complex inverse sine function $\arcsin s$ is defined as
\begin{equation}\label{eq:chap_3_9_1}
\arcsin s=e_{xy}[\arcsin(x+iy+z)-\arcsin z]+\arcsin z\texttt{ }(|s|\leq1).
\end{equation}
\end{definition}

\begin{definition}\label{def:chap_3_9_2}
The spatial complex inverse cosine function $\arccos s$ is defined as
\begin{equation}\label{eq:chap_3_9_2}
\arccos s=e_{xy}[\arccos(x+iy+z)-\arccos z]+\arccos z\texttt{ }(|s|\leq1).
\end{equation}
where $\arccos z=\pi/2-\arcsin z$ and $\arccos(x+iy+z)=\pi/2-\arcsin(x+iy+z)$.
\end{definition}

\begin{definition}\label{def:chap_3_9_3}
The spatial complex inverse tan function $tan^{-1}s$ is defined as
\begin{equation}\label{eq:chap_3_9_3}
\arctan s=e_{xy}[\arctan(x+iy+z)-\arctan z]+\arctan z\texttt{ }(|s|\leq1).
\end{equation}
\end{definition}

First, we shall prove equation~(\ref{eq:chap_3_9_1}).

Let $\sigma=x+iy$. Then $\sigma$ and $\sigma+z$ are two-dimensional complex numbers. The Taylor series of the real inverse sine function $\arcsin z$ and two-dimensional complex inverse sine function $\arcsin(\sigma+z)$ are
\[
\arcsin z=z+\sum_{n=1}^{\infty}\frac{(2n-1)!!}{(2n)!!}\frac{z^{2n+1}}{2n+1}\texttt{ }(|z|\leq1)
\]
and
\[
\arcsin(\sigma+z)=(\sigma+z)+\sum_{n=1}^{\infty}\frac{(2n-1)!!}{(2n)!!}\frac{(\sigma+z)^{2n+1}}{2n+1}\texttt{ }(|\sigma+z|\leq1).
\]

For an index $n$ in the Taylor series of $\arcsin(\sigma+z)$ and $\arcsin z$, there are
\[
(\sigma+z)^{2n+1}=\sigma^{2n+1}+(2n+1)\sigma^{2n}z+\cdots+(2n+1)\sigma z^{2n}+z^{2n+1}
\]
and $z^{2n+1}$, respectively.

Since both $\sigma+z$ and its powers are points in the two-dimensional complex plane and from remark~(\ref{re:chap_1_1_1}) for a term $e_{xy}x^{j}y^{k}$, the operator $e_{xy}$ becomes $e_{xy}^{m}$ where $m=|j|+|k|$, we get
\[
e_{xy}[(\sigma+z)^{2n+1}-z^{2n+1}]+z^{2n+1}
\]
\[
\begin{split}
&=(e_{xy}\sigma)^{2n+1}+(2n+1)(e_{xy}\sigma)^{2n}z+\cdots+(2n+1)(e_{xy}\sigma)z^{2n}+z^{2n+1} \\
&=(e_{xy}\sigma+z)^{2n+1}=s^{2n+1}.
\end{split}
\]

Thus there are
\[
\arcsin z+e_{xy}[\arcsin(\sigma+z)-\arcsin z]
=(z+\sum_{n=1}^{\infty}\frac{(2n-1)!!}{(2n)!!}\frac{z^{2n+1}}{2n+1})
\]
\[
+e_{xy}[(\sigma+z+\sum_{n=1}^{\infty}\frac{(2n-1)!!}{(2n)!!}\frac{(\sigma+z)^{2n+1}}{2n+1})
-(z+\sum_{n=1}^{\infty}\frac{(2n-1)!!}{(2n)!!}\frac{z^{2n+1}}{2n+1})]
\]
\[
=(e_{xy}\sigma+z)+\sum_{n=1}^{\infty}\frac{(2n-1)!!}{(2n)!!}\frac{(e_{xy}\sigma+z)^{2n+1}}{2n+1}
=s+\sum_{n=1}^{\infty}\frac{(2n-1)!!}{(2n)!!}\frac{s^{2n+1}}{2n+1}
=\arcsin s.
\]

Next, other equation%s~(\ref{eq:chap_3_9_2}) and
~(\ref{eq:chap_3_9_3}) can also be proved in this way. So the Taylor series of spatial complex inverse trigonometric functions $\arcsin s$, %$\arccos s$,
and $\arctan s$ are
\begin{equation}\label{eq:chap_3_9_4}
\arcsin s=s+\sum_{n=1}^{\infty}\frac{(2n-1)!!}{(2n)!!}\frac{s^{2n+1}}{2n+1}\texttt{ }(|s|\leq1),
\end{equation}
%\begin{equation}\label{eq:chap_3_9_5}
%\arccos s=\frac{\pi}{2}-(s+\sum_{n=1}^{\infty}\frac{(2n-1)!!}{(2n)!!}\frac{s^{2n+1}}{2n+1})\texttt{ }(|s|\leq1),
%\end{equation}
and
\begin{equation}\label{eq:chap_3_9_6}
\arctan s=\sum_{n=0}^{\infty}\frac{(-1)^{n}s^{2n+1}}{2n+1}\texttt{ }(|s|\leq1).
\end{equation}

\subsection{Inverse Hyperbolic Functions}\label{sec:chap_3_9_35_2}

Inverse of the hyperbolic functions can be treated in a corresponding manner.

\begin{definition}\label{def:chap_3_9_4}
The spatial complex inverse hyperbolic sine function $arc\sinh s$ is defined as
\begin{equation}\label{eq:chap_3_9_7}
arc\sinh s=e_{xy}[arc\sinh(x+iy+z)-arc\sinh z]+arc\sinh z\texttt{ }(|s|<1).
\end{equation}
\end{definition}

%\begin{definition}\label{def:chap_3_9_5}
%The spatial complex inverse hyperbolic cosine function $arc\cosh s$ is defined as
%\begin{equation}\label{eq:chap_3_9_8}
%arc\cosh s=e_{xy}(arc\cosh(x+iy+z)-arc\cosh z)+arc\cosh z\texttt{ }(|s|<1).
%\end{equation}
%\end{definition}

\begin{definition}\label{def:chap_3_9_6}
The spatial complex inverse hyperbolic tan function $arc\tanh s$ is defined as
\begin{equation}\label{eq:chap_3_9_9}
arc\tanh s=e_{xy}[arc\tanh(x+iy+z)-arc\tanh z]+arc\tanh z\texttt{ }(|s|<1).
\end{equation}
\end{definition}

First, we shall prove equation~(\ref{eq:chap_3_9_7}).

Let $\sigma=x+iy$. Then $\sigma$ and $\sigma+z$ are two-dimensional complex numbers. The Taylor series of the real inverse hyperbolic sine function $arc\sinh z$ and two-dimensional complex inverse hyperbolic sine function $arc\sinh(\sigma+z)$ are
\[
arc\sinh z=z+\sum_{n=1}^{\infty}\frac{(-1)^{n}(2n-1)!!}{(2n)!!}\frac{z^{2n+1}}{2n+1}\texttt{ }(|z|<1)
\]
and
\[
arc\sinh(\sigma+z)=(\sigma+z)+\sum_{n=1}^{\infty}\frac{(-1)^{n}(2n-1)!!}{(2n)!!}\frac{(\sigma+z)^{2n+1}}{2n+1}\texttt{ }(|\sigma+z|<1).
\]

For an index $n$ in the Taylor series of $arc\sinh(\sigma+z)$ and $arc\sinh z$, there are
\[
(\sigma+z)^{2n+1}=\sigma^{2n+1}+(2n+1)\sigma^{2n}z+\cdots+(2n+1)\sigma z^{2n}+z^{2n+1}
\]
and $z^{2n+1}$, respectively.

Since both $\sigma+z$ and its powers are points in the two-dimensional complex plane and from remark~(\ref{re:chap_1_1_1}) for a term $e_{xy}x^{j}y^{k}$, the operator $e_{xy}$ becomes $e_{xy}^{m}$ where $m=|j|+|k|$, we get
\[
e_{xy}[(\sigma+z)^{2n+1}-z^{2n+1}]+z^{2n+1}
\]
\[
\begin{split}
&=(e_{xy}\sigma)^{2n+1}+(2n+1)(e_{xy}\sigma)^{2n}z+\cdots+(2n+1)(e_{xy}\sigma)z^{2n}+z^{2n+1} \\
&=(e_{xy}\sigma+z)^{2n+1}=s^{2n+1}.
\end{split}
\]

Thus there are
\[
arc\sinh z+e_{xy}[arc\sinh(\sigma+z)-arc\sinh z]
=(z+\sum_{n=1}^{\infty}\frac{(-1)^{n}(2n-1)!!}{(2n)!!}\frac{z^{2n+1}}{2n+1})
\]
\[
+e_{xy}[(\sigma+z+\sum_{n=1}^{\infty}\frac{(-1)^{n}(2n-1)!!}{(2n)!!}\frac{(\sigma+z)^{2n+1}}{2n+1})
-(z+\sum_{n=1}^{\infty}\frac{(-1)^{n}(2n-1)!!}{(2n)!!}\frac{z^{2n+1}}{2n+1})]
\]
\[
=(e_{xy}\sigma+z)+\sum_{n=1}^{\infty}\frac{(-1)^{n}(2n-1)!!}{(2n)!!}\frac{(e_{xy}\sigma+z)^{2n+1}}{2n+1}
\]
\[
=s+\sum_{n=1}^{\infty}\frac{(-1)^{n}(2n-1)!!}{(2n)!!}\frac{s^{2n+1}}{2n+1}
=arc\sinh s.
\]

Next, other equations~(\ref{eq:chap_3_9_7}) and~(\ref{eq:chap_3_9_9}) can also be proved in this way. So the Taylor series of spatial complex inverse hyperbolic trigonometric functions $arc\sinh s$, %$arc\cosh s$,
and $arc\tanh s$ are
\begin{equation}\label{eq:chap_3_9_10}
arc\sinh s=s+\sum_{n=1}^{\infty}\frac{(-1)^{n}(2n-1)!!}{(2n)!!}\frac{s^{2n+1}}{2n+1}\texttt{ }(|s|<1),
\end{equation}
%\begin{equation}\label{eq:chap_3_9_11}
%arc\cosh s=\frac{\pi}{2}-(s+\sum_{n=1}^{\infty}\frac{(-1)^{n}(2n-1)!!}{(2n)!!}\frac{s^{2n+1}}{2n+1})\texttt{}(|s|<1),
%\end{equation}
and
\begin{equation}\label{eq:chap_3_9_12}
arc\tanh s=\sum_{n=0}^{\infty}\frac{s^{2n+1}}{2n+1}\texttt{ }(|s|<1).
\end{equation}

\subsection{Inverse Trigonometric and Hyperbolic Functions Described in Terms of Logarithms}\label{sec:chap_3_9_35_3}

Inverses of the trigonometric and hyperbolic functions can be described in terms of logarithms.

In order to define the inverse sine function $\arcsin s$, we write
\[
\varpi=\arcsin s\texttt{ when }s=\sin\varpi.
\]
That is, $\varpi=\arcsin s$ when
\[
s=\frac{e^{i\varpi}-e^{-i\varpi}}{2i}.
\]
If we put this equation in the form
\[
(e^{i\varpi})^{2}-2ise^{i\varpi}-1=0,
\]
which is quadratic in $e^{i\varpi}$, and solve for $e^{i\varpi}$, we find that
\begin{equation}\label{eq:chap_3_9_13}
e^{i\varpi}=is+(1-s^{2})^{1/2},
\end{equation}
where $(1-s^{2})^{1/2}$ is, of course, a double-valued function of $s$. Taking logarithms of
each side of equation~(\ref{eq:chap_3_9_13}) and recalling that $\varpi=\arcsin s$, we arrive at the expression
\begin{equation}\label{eq:chap_3_9_14}
\arcsin s=-i\log[is+(1-s^{2})^{1/2}].
\end{equation}

One can apply the technique used to derive expression~(\ref{eq:chap_3_9_14}) for $\arcsin s$ to show that
\begin{equation}\label{eq:chap_3_9_15}
\arccos s=-i\log[s+i(1-s^{2})^{1/2}]
\end{equation}
and that
\begin{equation}\label{eq:chap_3_9_16}
\arctan s=\frac{i}{2}\log\frac{i+s}{i-s}.
\end{equation}
The functions $\arccos s$ and $\arctan s$ are also multiple-valued. When specific branches of the square root and logarithmic functions are used, all three inverse functions become single-valued and analytic because they are then compositions of analytic functions.

The derivatives of these three functions are readily obtained from the above expressions. The derivatives of the first two depend on the values chosen for the square roots:
\begin{equation}\label{eq:chap_3_9_17}
\frac{d}{ds}\arcsin s=\frac{1}{(1-s^{2})^{1/2}},
\end{equation}
\begin{equation}\label{eq:chap_3_9_18}
\frac{d}{ds}\arccos s=\frac{-1}{(1-s^{2})^{1/2}}.
\end{equation}
The derivative of the last one,
\begin{equation}\label{eq:chap_3_9_19}
\frac{d}{ds}\arctan s=\frac{1}{1+s^{2}},
\end{equation}
does not, however, depend on the manner in which the function is made single-valued.

Inverse hyperbolic functions can be treated in a corresponding manner. It turns out that
\begin{equation}\label{eq:chap_3_9_20}
arc\sinh s=\log[s+(s^{2}+1)^{1/2}],
\end{equation}
\begin{equation}\label{eq:chap_3_9_21}
arc\cosh s=\log[s+(s^{2}-1)^{1/2}],
\end{equation}
and
\begin{equation}\label{eq:chap_3_9_22}
arc\tanh s=\frac{1}{2}\log\frac{1+s}{1-s}.
\end{equation}

Finally, we compare expression~(\ref{eq:chap_3_9_1}) and expression~(\ref{eq:chap_3_9_14}).

Let $\sigma=x+iy$. Then $\sigma$ and $\sigma+z$ are two-dimensional complex numbers, and $s=e_{xy}\sigma+z$. Similarly, for a two-dimensional complex number $\sigma$ and a real $z$, there are
\[
\arcsin(\sigma+z)=-i\log[i(\sigma+z)+(1-(\sigma+z)^{2})^{1/2}]
\]
and
\[
\arcsin z=-i\log[iz+(1-z^{2})^{1/2}],
\]
respectively.

When specific branches of the square root and logarithmic functions are used, all three inverse functions become single-valued and analytic. So all these logarithmic functions are analytic. Then from Theorem~(\ref{th:chap_2_13_3}) in Sec.~(\ref{sec:chap_2_13_23}), the function $\log[is+(1-s^{2})^{1/2}]$ can be written as
\[
\log[is+(1-s^{2})^{1/2}]=e_{xy}\{\log[i(\sigma+z)+(1-(\sigma+z)^{2})^{1/2}]
\]
\[
-\log[iz+(1-z^{2})^{1/2}]\}+\log[iz+(1-z^{2})^{1/2}].
\]
Multiplying through this equation by $-i$ and then we obtain
\[
\arcsin s=e_{xy}[\arcsin(\sigma+z)-\arcsin z]+\arcsin z.
\]
That is expression~(\ref{eq:chap_3_9_1})
\[
\arctan s=e_{xy}[\arctan(x+iy+z)-\arctan z]+\arctan z\texttt{ }(|s|<1).
\]

%-----------------------------------------------------------------------
% Beginning of chap4.tex
%-----------------------------------------------------------------------
%
% AMS-LaTeX 1.2 sample file for a monograph, based on amsbook.cls.
% This is a data file input by chapter.tex.
%%%%%%%%%%%%%%%%%%%%%%%%%%%%%%%%%%%%%%%%%%%%%%%%%%%%%%%%%%%%%%%%%%%%

%\part{This is a Part Title Sample}

\chapter{Contour Integrals}\label{ch:chap_4}

\section{Derivatives of Functions $\varpi(t)$}\label{sec:chap_4_1_36}

We first consider derivatives of spatial complex-valued functions $\varpi(t)$ of a real variable $t$. Let us write
\begin{equation}\label{eq:chap_4_1_1}
\varpi(t)=e_{xy}[u(t)+iv(t)]+w(t)
\end{equation}
where the functions $u$, $v$, and $w$ are real-valued functions of $t$. The derivative $\varpi'(t)$, or $d[\varpi(t)]/dt$, of the function~(\ref{eq:chap_4_1_1}) at a point $t$ is defined as
\begin{equation}\label{eq:chap_4_1_2}
\varpi'(t)=e_{xy}[u'(t)+iv'(t)]+w'(t)
\end{equation}
provided each of the derivatives $u'$, $v'$, and $w'$ exists at $t$.

From definition~(\ref{eq:chap_4_1_2}), it follows that, for every spatial complex constant $s_{0}=e_{xy}(x_{0}+iy_{0})+z_{0}$,
\[
\frac{d}{dt}[s_{0}\varpi(t)]=\{[e_{xy}(x_{0}+iy_{0})+z_{0}][e_{xy}(u+iv)+w]\}'
\]
\[
=e_{xy}[(x_{0}u'-y_{0}v'+x_{0}w'+z_{0}u')+i(x_{0}v'+y_{0}u'+y_{0}w'+z_{0}v')]+z_{0}w'.
\]
But
\[
e_{xy}[(x_{0}u'-y_{0}v'+x_{0}w'+z_{0}u')+i(x_{0}v'+y_{0}u'+y_{0}w'+z_{0}v')]+z_{0}w'=s_{0}\varpi'(t),
\]
and so
\begin{equation}\label{eq:chap_4_1_3}
\frac{d}{dt}[s_{0}\varpi(t)]=s_{0}\varpi'(t).
\end{equation}

Another expected rule that we shall often use is
\begin{equation}\label{eq:chap_4_1_4}
\frac{d}{dt}e^{s_{0}t}=s_{0}e^{s_{0}t}.
\end{equation}

Various other rules learned in calculus and two-dimensional complex variables, such as the ones for differentiating sums
and products, apply just as they do for real-valued functions and two-dimensional complex-valued functions of $t$. As was the case with property~(\ref{eq:chap_4_1_3}) and formula~(\ref{eq:chap_4_1_4}), verifications may be based on corresponding rules in calculus and two-dimensional complex variables. It should be pointed out, however, that not every rule for derivatives in calculus and two-dimensional complex variables carries over to functions of type~(\ref{eq:chap_4_1_1}).

\section{Definite Integrals of Functions $\varpi(t)$}\label{sec:chap_4_2_37}

When $\varpi(t)$ is a spatial complex-valued function of a real variable $t$ and is written
\begin{equation}\label{eq:chap_4_2_1}
\varpi(t)=e_{xy}[u(t)+iv(t)]+w(t)
\end{equation}
where $u$, $v$, and $w$ are real-valued, the definite integral of $\varpi(t)$ over an interval $a\leq t\leq b$ is defined as
\begin{equation}\label{eq:chap_4_2_2}
\int_{a}^{b}\varpi(t)dt=e_{xy}[\int_{a}^{b}u(t)dt+i\int_{a}^{b}v(t)dt]+\int_{a}^{b}w(t)dt
\end{equation}
when the individual integrals on the right exist. Thus
\begin{equation}\label{eq:chap_4_2_3}
Re_{s}\int_{a}^{b}\varpi(t)dt=\int_{a}^{b}Re_{s}\texttt{ }\varpi(t)dt\texttt{ and }Im_{s}\int_{a}^{b}\varpi(t)dt=\int_{a}^{b}Im_{s}\varpi(t)dt.
\end{equation}

Improper integrals of $\varpi(t)$ over unbounded intervals are defined in a similar way.

The existence of the integrals of $u$, $v$, and $w$ in definition~(\ref{eq:chap_4_2_2}) is ensured if those functions are piecewise continuous on the interval $a\leq t\leq b$. Such a function is continuous everywhere in the stated interval except possibly for a finite number of points where, although discontinuous, it has one-sided limits. Of course, only the right-hand limit is required at $a$; and only the left-hand limit is required at $b$. When all $u$, $v$, and $w$ are piecewise continuous, the function $\varpi$ is said to have that property.

Anticipated rules for integrating a complex constant times a function $\varpi(t)$, for integrating sums of such functions, and for interchanging limits of integration are all valid. Those rules, as well as the property
\[
\int_{a}^{b}\varpi(t)dt=\int_{a}^{c}\varpi(t)dt+\int_{c}^{b}\varpi(t)dt,
\]
are easy to verify by recalling corresponding results in calculus.

The fundamental theorem of calculus, involving antiderivatives, can, moreover, be extended so as to apply to integrals of the type~(\ref{eq:chap_4_2_2}). To be specific, suppose that the functions
\[
\varpi(t)=e_{xy}[u(t)+iv(t)]+w(t)\texttt{ and }
\overline{W}(t)=e_{xy}[U(t)+iV(t)]+W(t)
\]
are continuous on the interval $a\leq t\leq b$. If $\overline{W}'(t)=\varpi(t)$ when $a\leq t\leq b$, then $U'(t)=u(t)$, $V'(t)=v(t)$, and $W'(t)=w(t)$. Hence, in view of definition~(\ref{eq:chap_4_2_2}),
\[
\begin{split}
\int_{a}^{b}\varpi(t)dt &= e_{xy}[U(t)|_{a}^{b}+iV(t)|_{a}^{b}]+W(t)|_{a}^{b} \\
                        &= e_{xy}[U(b)+iV(b)]+W(b)-e_{xy}[U(a)+iV(a)]-W(a).
\end{split}
\]
That is,
\begin{equation}\label{eq:chap_4_2_4}
|\int_{a}^{b}\varpi(t)dt|=\overline{W}(b)-\overline{W}(a)=\overline{W}(t)|_{a}^{b}.
\end{equation}

We finish here with an important property of moduli of integrals. Namely,
\begin{equation}\label{eq:chap_4_2_5}
|\int_{a}^{b}\varpi(t)dt|\leq\int_{a}^{b}|\varpi(t)|dt\texttt{ }(a\leq b).
\end{equation}
This inequality clearly holds when the value of the integral on the left is zero, in particular when $a=b$. Thus, in the verification, we may assume that its value is a nonzero spatial complex number. If $r_{0}$ is the modulus and $\theta_{0}$ and $\varphi_{0}$ are arguments of that constant, then
\[
\int_{a}^{b}\varpi(t)dt=r_{0}[e_{xy}e^{i\theta_{0}}\cos\varphi_{0}+\sin\varphi_{0}].
\]
Solving for $r_{0}$, we write
\begin{equation}\label{eq:chap_4_2_6}
r_{0}=\int_{a}^{b}\frac{\varpi(t)dt}{e_{xy}e^{i\theta_{0}}\cos\varphi_{0}+\sin\varphi_{0}}.
\end{equation}
Now the left-hand side of this equation is a real number, and so the right-hand side is too. Thus, using the fact that the real part of a real number is the number itself and referring to the first of properties~(\ref{eq:chap_4_2_3}), we see that the right-hand side of equation~(\ref{eq:chap_4_2_6}) can be rewritten in the following way:
\[
\begin{split}
\int_{a}^{b}\frac{\varpi(t)}{e_{xy}e^{i\theta_{0}}\cos\varphi_{0}+\sin\varphi_{0}}dt
&=Re_{s}\int_{a}^{b}\frac{\varpi(t)}{e_{xy}e^{i\theta_{0}}\cos\varphi_{0}+\sin\varphi_{0}}dt \\
&=\int_{a}^{b}Re_{s}\frac{\varpi(t)}{e_{xy}e^{i\theta_{0}}\cos\varphi_{0}+\sin\varphi_{0}}dt.
\end{split}
\]
Equation~(\ref{eq:chap_4_2_6}) then takes the form
\begin{equation}\label{eq:chap_4_2_7}
r_{0}=\int_{a}^{b}Re_{s}\frac{\varpi(t)}{e_{xy}e^{i\theta_{0}}\cos\varphi_{0}+\sin\varphi_{0}}dt.
\end{equation}
But
\[
\begin{split}
Re_{s}\frac{\varpi(t)}{e_{xy}e^{i\theta_{0}}\cos\varphi_{0}+\sin\varphi_{0}}
&\leq|\frac{\varpi(t)}{e_{xy}e^{i\theta_{0}}\cos\varphi_{0}+\sin\varphi_{0}}| \\
&=\frac{|\varpi(t)|}{|e_{xy}e^{i\theta_{0}}\cos\varphi_{0}+\sin\varphi_{0}|}=|\varpi(t)|;
\end{split}
\]
and so, according to equation~(\ref{eq:chap_4_2_7}),
\[
r_{0}\leq\int_{a}^{b}|\varpi(t)|dt.
\]
Because $r_{0}$ is, in fact, the left-hand side of inequality~(\ref{eq:chap_4_2_5}) when the value of the integral
there is nonzero, the verification is now complete.

With only minor modifications, the above discussion yields inequalities such as
\begin{equation}\label{eq:chap_4_2_8}
|\int_{0}^{\infty}\varpi(t)dt|\leq\int_{0}^{\infty}|\varpi(t)|dt,
\end{equation}
provided both improper integrals exist.

\section{Contours}\label{sec:chap_4_3_38}

Integrals of spatial complex-valued functions of a spatial complex variable are defined on curves in
the three-dimensional complex space, rather than on just intervals of the real line or in just areas of the two-dimensional complex plane. Classes of curves that are adequate for the study of such integrals are introduced in this section.

A set of points $s=(x,y,z)$ in the three-dimensional complex space is said to be an arc if
\begin{equation}\label{eq:chap_4_3_1}
x=x(t)\texttt{, }y=y(t)\texttt{, and }z=z(t)\texttt{ }(a\leq t\leq b),
\end{equation}
where $x(t)$, $y(t)$, and $z(t)$ are continuous functions of the real parameter $t$. This definition
establishes a continuous mapping of the interval $a\leq t\leq b$ into the $xyz$, or $s$, space; and
the image points are ordered according to increasing values of $t$. It is convenient to
describe the points of $C$ by means of the equation
\begin{equation}\label{eq:chap_4_3_2}
s=s(t)\texttt{ }(a\leq t\leq b)
\end{equation}
where
\begin{equation}\label{eq:chap_4_3_3}
s(t)=e_{xy}[x(t)+iy(t)]+z(t).
\end{equation}

The arc $C$ is a simple arc, or a Jordan arc, if it does not cross itself; that is, $C$ is
simple if $s(t_{1})\neq s(t_{2})$ when $t_{1}\neq t _{2}$. When the arc $C$ is simple except for the fact that
$s(b)=s(a)$, we say that $C$ is a simple closed curve, or a Jordan curve.

The geometric nature of a particular arc often suggests different notation for the parameter $t$ in equation~(\ref{eq:chap_4_3_2}).

The parametric representation used for any given arc $C$ is, of  course, not unique. It is, in fact, possible to change the interval over which the parameter ranges to any other interval. To be specific, suppose that
\begin{equation}\label{eq:chap_4_3_9}
t=\phi(\tau)\texttt{ }(\alpha\leq\tau\leq\beta),
\end{equation}
where $\phi$ is a real-valued function mapping an interval $\alpha\leq\tau\leq\beta$ onto the interval
$a\leq t\leq b$ in representation~(\ref{eq:chap_4_3_2}). We assume that $\phi$ is continuous with a continuous derivative. We also assume that $\phi'(\tau)>0$ for each $\tau$; this ensures that $t$ increases with $\tau$. Representation~(\ref{eq:chap_4_3_2}) is then transformed by equation~(\ref{eq:chap_4_3_9}) into
\begin{equation}\label{eq:chap_4_3_10}
s=S(\tau)\texttt{ }(\alpha\leq\tau\leq\beta)
\end{equation}
where
\begin{equation}\label{eq:chap_4_3_11}
S(\tau)=s(\phi(\tau)).
\end{equation}

Suppose now that the components $x'(t)$, $y'(t)$, and $z'(t)$ of the derivative (Sec.~(\ref{sec:chap_4_1_36}))
\begin{equation}\label{eq:chap_4_3_12}
s'(t)=e_{xy}[x'(t)+iy'(t)]+z'(t)
\end{equation}
of the function~(\ref{eq:chap_4_3_3}), used to represent $C$, are continuous on the entire interval $a\leq t\leq b$. The arc is then called a differentiable arc, and the real-valued function
\[
|s'(t)|=\sqrt{[x'(t)]^{2}+[y'(t)]^{2}+[z'(t)]^{2}}
\]
is integrable over the interval $a\leq t\leq b$. In fact, according to the definition of arc length
in calculus, the length of $C$ is the number
\begin{equation}\label{eq:chap_4_3_13}
L=\int_{a}^{b}|s'(t)|dt.
\end{equation}

The value of $L$ is invariant under certain changes in the representation for $C$ that is used, as one would expect. More precisely, with the change of variable indicated in equation~(\ref{eq:chap_4_3_9}), expression~(\ref{eq:chap_4_3_13}) takes the form
\[
L=\int_{\alpha}^{\beta}|s'(\phi(\tau))|\phi'(\tau)d\tau.
\]
So, if representation~(\ref{eq:chap_4_3_10}) is used for $C$, the derivative% (Exercise 4)
\begin{equation}\label{eq:chap_4_3_14}
S'(\tau)=s'(\phi(\tau))\phi'(\tau),
\end{equation}
enables us to write expression~(\ref{eq:chap_4_3_13}) as
\[
L=\int_{\alpha}^{\beta}|S'(\tau)|d\tau.
\]
Thus the same length of $C$ would be obtained if representation~(\ref{eq:chap_4_3_10}) were to be used.

If equation~(\ref{eq:chap_4_3_2}) represents a differentiable arc and if $s'(t)\neq0$ anywhere in the
interval $a<t<b$, then the unit tangent vector
\[
\mathbf{T}=\frac{s'(t)}{|s'(t)|}
\]
is well defined for all $t$ in that open interval, with angle of inclination $\arg s'(t)$. Also, when $\mathbf{T}$ turns, it does so continuously as the parameter $t$ varies over the entire interval $a<t<b$. This expression for $\mathbf{T}$ is the one learned in calculus when $s(t)$ is interpreted as a radius vector. Such an arc is said to be smooth. In referring to a smooth arc $s=s(t)$ $(a\leq t\leq b)$, then, we agree that the derivative $s'(t)$ is continuous on the closed interval $a\leq t\leq b$ and nonzero on the open interval $a<t<b$.

A contour, or piecewise smooth arc, is an arc consisting of a finite number of smooth arcs joined end to end. Hence  if equation~(\ref{eq:chap_4_3_2}) represents a contour, $s(t)$ is continuous, whereas its derivative $s'(t)$ is piecewise continuous. The polygonal line is, for example, a contour. When only the initial and final values of $s(t)$ are the same, a contour $C$ is called a simple closed contour.
%Examples are the circles (5) and (6), as well as the boundary of a triangle or a rectangle taken in a specific direction.
The length of a contour or a simple closed contour is the sum of the lengths of the smooth arcs that make up the contour.

Given a positively oriented simple closed contour $C$, let all points on the contour $C$ be in or on a minimum closed space $S$. So two positively oriented simple closed contours $C_{1}$ and $C_{2}$ correspond to two spaces $S_{1}$ and $S_{2}$, respectively. If all points on the contour $C_{2}$ are inside $S_{1}$ and all points on the contour $C_{1}$ are outside $S_{2}$, then we say that $C_{2}$ is interior to $C_{1}$.

%The points on any simple closed curve or simple closed contour $C$ are boundary points of two distinct domains, one of which is the interior of $C$ and is bounded. The other, which is the exterior of $C$, is unbounded. It will be convenient to accept this statement, known as the Jordan curve theorem, as geometrically evident; the proof is not easy.

\section{Contour Integrals}\label{sec:chap_4_4_39}

We turn now to integrals of spatial complex-valued functions $f$ of the spatial complex variable $s$. Such an integral is defined in terms of the values $f(s)$ along a given contour $C$, extending from a point $s=s_{1}$ to a point $s=s_{2}$ in the three-dimensional complex space. It is, therefore, a line integral; and its value depends, in general, on the contour $C$ as well as on the function $f$. It is written
\[
\int_{C}f(s)ds\texttt{ or }\int_{s_{1}}^{s_{2}}f(s)ds
\]
the latter notation often being used when the value of the integral is independent of
the choice of the contour taken between two fixed end points.

Suppose that the equation
\begin{equation}\label{eq:chap_4_4_1}
s=s(t)\texttt{ }(a\leq t\leq b)
\end{equation}
represents a contour $C$, extending from a point $s_{1}=s(a)$ to a point $s_{2}=s(b)$. Let the
function $f(s)$ be piecewise continuous on $C$; that is, $f[s(t)]$ is piecewise continuous
on the interval $a\leq t\leq b$. We define the line integral, or contour integral, of $f$ along
$C$ as follows:
\begin{equation}\label{eq:chap_4_4_2}
\int_{C}f(s)ds=\int_{s_{1}}^{s_{2}}f[s(t)]s'(t)dt.
\end{equation}
Note that, since $C$ is a contour, $s'(t)$ is also piecewise continuous on the interval
$a\leq t\leq b$; and so the existence of integral~(\ref{eq:chap_4_4_2}) is ensured.

The value of a contour integral is invariant under a change in the representation of its contour when the change is  of the type~(\ref{eq:chap_4_3_11}), Sec.~(\ref{sec:chap_4_3_38}). This can be seen by following the same general procedure that was used in Sec.~(\ref{sec:chap_4_3_38}) to show the invariance of arc length.

It follows immediately from definition~(\ref{eq:chap_4_4_2}) and properties of integrals of spatial complex-valued functions $\varpi(t)$ mentioned in Sec.~(\ref{sec:chap_4_2_37}) that
\begin{equation}\label{eq:chap_4_4_3}
\int_{C}s_{0}f(s)ds=s_{0}\int_{C}f(s)ds,
\end{equation}
for any complex constant $s_{0}$, and
\begin{equation}\label{eq:chap_4_4_4}
\int_{C}[f(s)+g(s)]ds=\int_{C}f(s)ds+\int_{C}g(s)ds.
\end{equation}

Associated with the contour $C$ used in integral~(\ref{eq:chap_4_4_2}) is the contour $-C$, consisting of the same set of points but with the order reversed so that the new contour extends from the point $s_{2}$ to the point $s_{1}$. The contour $-C$ has parametric representation
\[
s=s(-t)\texttt{ }(-b\leq t\leq-a);
\]
and so, %in view of  Exercise 1 (a), Sec. 3 7,
\[
\int_{-C}f(s)ds=\int_{-b}^{-a}f(s(-t))\frac{d}{dt}s(-t)dt=-\int_{-b}^{-a}f(s(-t))s'(-t)dt,
\]
where $s'(-t)$ denotes the derivative of $s(t)$ with respect to $t$, evaluated at $-t$. Making
the substitution $\tau=-t$ in this last integral,% and referring to Exercise  l(a),  Sec. 38,
we obtain the expression
\[
\int_{-C}f(s)ds=-\int_{a}^{b}f(s(\tau))s'(\tau)d\tau,
\]
which is the same as
\begin{equation}\label{eq:chap_4_4_5}
\int_{-C}f(s)ds=-\int_{C}f(s)ds,
\end{equation}

Consider now a path C, with representation~(\ref{eq:chap_4_4_1}), that consists of a contour $C_{1}$ from $s_{1}$ to $s_{2}$ followed by a contour $C_{2}$ from $s_{2}$ to $s_{3}$, the initial point of $C_{2}$ being the final
point of $C_{1}$. There is a value $c$ of $t$, where $a<c<b$, such that $s(c)=s_{2}$. Consequently, $C_{1}$ is represented by
\[
s=s(t)\texttt{ }(a\leq t\leq c)
\]
and $C_{2}$ is represented by
\[
s=s(t)\texttt{ }(c\leq t\leq b).
\]
Also, by a rule for integrals of functions $\varpi(t)$ that was noted in Sec.~(\ref{sec:chap_4_2_37}),
\[
\int_{a}^{b}f[s(t)]s'(t)dt=\int_{a}^{c}f[s(t)]s'(t)dt+\int_{c}^{b}f[s(t)]s'(t)dt.
\]
Evidently, then,
\begin{equation}\label{eq:chap_4_4_6}
\int_{C}f(s)ds=\int_{C_{1}}f(s)ds+\int_{C_{2}}f(s)ds.
\end{equation}
Sometimes the contour $C$ is called the sum of its legs $C_{1}$ and $C_{2}$ and is denoted by $C_{1}+C_{2}$. The sum of two contours $C_{1}$ and $-C_{2}$ is well defined when $C_{1}$ and $C_{2}$ have the same final points, and it is written $C_{1}-C_{2}$.

Definite integrals in calculus can be interpreted as areas, and they have other interpretations as well. Except in special cases, no corresponding helpful interpretation, geometric or physical, is available for integrals in the complex space.

\section{Examples}\label{sec:chap_4_5_40}

The purpose of this section is to provide examples of the definition in Sec.~(\ref{sec:chap_4_4_39}) of
contour integrals and to illustrate various properties that were mentioned there. We defer development of the concept of antiderivatives of the integrands $f(s)$ in contour integrals until Sec.~(\ref{sec:chap_4_7_42}).

\begin{example}\label{ex:chap_4_5_1}
Let us find the value of the integral
\begin{equation}\label{eq:chap_4_5_1}
I=\int_{C}\bar{s}ds
\end{equation}
when $C$ is the right-hand half
\[
s=r(e_{xy}e^{i\theta}\cos\varphi+\sin\varphi)\texttt{ }(-\pi/2\leq\theta\leq\pi/2)
\]
of the circle $|s|=r$ in a spatial plane, from $\theta=-\pi/2$ to $\theta=\pi/2$, where $r$ and $\varphi$ are constants. According to definition~(\ref{eq:chap_4_4_2}), Sec.~(\ref{sec:chap_4_4_39}),
\[
I=\int_{-\pi/2}^{\pi/2}\overline{r(e_{xy}e^{i\theta}\cos\varphi+\sin\varphi)}[r(e_{xy}e^{i\theta}\cos\varphi+\sin\varphi)]'d\theta;
\]
and, since
\[
\overline{r(e_{xy}e^{i\theta}\cos\varphi+\sin\varphi)}=r(e_{xy}e^{-i\theta}\cos\varphi+\sin\varphi)
\]
and
\[
[r(e_{xy}e^{i\theta}\cos\varphi+\sin\varphi)]'=re_{xy}ie^{i\theta}\cos\varphi,
\]
this means that
\[
\begin{split}
I&=\int_{-\pi/2}^{\pi/2}r(e_{xy}e^{-i\theta}\cos\varphi+\sin\varphi)re_{xy}ie^{i\theta}\cos\varphi d\theta\\
 &=r^{2}e_{xy}\int_{-\pi/2}^{\pi/2}(i\cos^{2}\varphi+\frac{1}{2}ie^{i\theta}\sin2\varphi)d\theta \\
 &=r^{2}e_{xy}(i\theta\cos^{2}\varphi+\frac{1}{2}e^{i\theta}\sin2\varphi)|_{-\pi/2}^{\pi/2} \\
 &=r^{2}e_{xy}i(\pi\cos^{2}\varphi+\sin2\varphi).
\end{split}
\]

Note that when a point $s$ is on the circle $|s|=r$, it follows that $s\bar{s}=r^{2}$, or $\bar{s}=r^{2}/s$. Hence the result $I=r^{2}e_{xy}i(\pi\cos^{2}\varphi+\sin2\varphi)$ can also be written
\begin{equation}\label{eq:chap_4_5_2}
I=\int_{C}\bar{s}ds=e_{xy}i(\pi\cos^{2}\varphi+\sin2\varphi).
\end{equation}
\end{example}

\begin{example}\label{ex:chap_4_5_2}
In this example, we first let $C_{1}$ denote the contour $OAB$ where $O$ denotes the origin, $A$ denotes the point $(0,1,0)$, and $B$ denotes the point $(1,1,2)$, and evaluate the integral
\begin{equation}\label{eq:chap_4_5_3}
\int_{C_{1}}f(s)ds=\int_{OA}f(s)ds+\int_{AB}f(s)ds,
\end{equation}
where
\[
f(s)=e_{xy}(y-x-i3x^{2})+z\texttt{ }(s=e_{xy}(x+iy)+z).
\]
The leg $OA$ may be represented parametrically as $s=e_{xy}(0+iy)+0$ $(0\leq y\leq1)$; and since $x=0$ at points on that leg, the values of $f$ there vary with the parameter $y$ according to the equation $f(s)=e_{xy}y$ $(0\leq y\leq1)$. Consequently,
\[
\int_{OA}f(s)ds=\int_{0}^{1}e_{xy}yidy=e_{xy}i\int_{0}^{1}ydy=e_{xy}\frac{i}{2}.
\]
On the leg $AB$, $s=e_{xy}(x+i)+z$ $(0\leq x\leq1,z=2x)$; and so
\[
\begin{split}
\int_{AB}f(s)ds &= \int_{0}^{1}[e_{xy}(1-x-i3x^{2})+x](e_{xy}1+2)dx \\
                &= (e_{xy}1+2)[e_{xy}\int_{0}^{1}(1-x-i3x^{2})dx+\int_{0}^{1}xdx] \\
                &= (e_{xy}1+2)[e_{xy}(x-\frac{x^{2}}{2}-ix^{3})+\frac{x^{2}}{2}]|_{0}^{1} \\
                &= (e_{xy}1+2)[e_{xy}(\frac{1}{2}-i)+\frac{1}{2}]=e_{xy}(2-3i)+1.
\end{split}
\]
In view of equation~(\ref{eq:chap_4_5_3}), we now see that
\begin{equation}\label{eq:chap_4_5_4}
\int_{C_{1}}f(s)ds=e_{xy}\frac{i}{2}+e_{xy}(2-3i)+1=e_{xy}(2-\frac{5}{2}i)+1.
\end{equation}

If $C_{2}$ denotes the segment $OB$ of the line $y=x$ and $z=2x$ with parametric representation $s=e_{xy}(x+iy)+z$ $(0\leq y=x\leq1,z=2x)$,
\begin{equation}\label{eq:chap_4_5_5}
\begin{split}
\int_{C_{2}}f(s)ds &= \int_{0}^{1}(-e_{xy}i3x^{2}+x)[e_{xy}(1+i)+2]dx \\
                   &= [e_{xy}(1+i)+2](-e_{xy}i\int_{0}^{1}3x^{2}dx+\int_{0}^{1}xdx) \\
                   &= [e_{xy}(1+i)+2](-e_{xy}ix^{3}+\frac{x^{2}}{2})|_{0}^{1} \\
                   &= [e_{xy}(1+i)+2](-e_{xy}i+1/2)=e_{xy}(\frac{3}{2}-\frac{5}{2}i)+1
\end{split}
\end{equation}
Evidently, then, the integrals of $f(s)$ along the two paths $C_{1}$ and $C_{2}$ have \textit{different values} even though those paths have the same initial and the same final points.

Observe how it follows that the integral of $f(s)$ over the simple closed contour
$OABO$, or $C_{1}-C_{2}$, has the \textit{nonzero} value
\[
\int_{C_{1}}f(s)ds-\int_{C_{2}}f(s)ds=[e_{xy}(2-\frac{5}{2}i)+1]-[e_{xy}(\frac{3}{2}-\frac{5}{2}i)+1]=e_{xy}\frac{1}{2}.
\]
\end{example}

\begin{example}\label{ex:chap_4_5_3}
We begin here by letting $C$ denote an arbitrary smooth arc
\[
s=s(t)\texttt{ }(a\leq t\leq b)
\]
from a fixed point $s_{1}$ to a fixed point $s_{2}$. In order to evaluate the integral
\[
I=\int_{C}sds=\int_{a}^{b}s(t)s'(t)dt.
\]
we note that,
\[
\frac{d}{dt}\frac{[s(t)]^{2}}{2}=s(t)s'(t).
\]
Thus
\[
I=\frac{[s(t)]^{2}}{2}|_{a}^{b}=\frac{[s(b)]^{2}-[s(a)]^{2}}{2}.
\]
But $s(b)=s_{2}$ and $s(a)=s_{1}$; and so $I=(s_{2}^{2}-s_{1}^{2})/2$. Inasmuch as the value of $I$
depends only on the end points of $C$, and is otherwise independent of the arc that is taken, we may write
\begin{equation}\label{eq:chap_4_5_6}
\int_{s_{1}}^{s_{2}}sds=\frac{s_{2}^{2}-s_{1}^{2}}{2}.
\end{equation}
(Compare Example~(\ref{ex:chap_4_5_2}), where the value of an integral from one fixed point to another depended on the path that was taken.)

Expression~(\ref{eq:chap_4_5_6}) is also valid when $C$ is a spatial contour in a space that is not necessarily smooth since a contour consists of a finite number of smooth arcs $C_{k}$ $(k=1, 2, \ldots, n)$, joined end to end. More precisely, suppose that each $C_{k}$ extends from $s_{k}$ to $s_{k+1}$. Then
\begin{equation}\label{eq:chap_4_5_7}
\int_{C}sds=\sum_{k=1}^{n}\int_{C_{k}}sds=\sum_{k=1}^{n}\frac{s_{k+1}^{2}-s_{k}^{2}}{2}=\frac{s_{n+1}^{2}-s_{1}^{2}}{2},
\end{equation}
$s_{1}$ being the initial point of $C$ and $s_{n+1}$ its final point.

It follows from expression~(\ref{eq:chap_4_5_7}) that the integral of the function $f(s)=s$ around each closed contour in the entire space has value zero. (Once again, compare Example~(\ref{ex:chap_4_5_2}), where the value of the integral of a given function around a certain closed path was \textit{not} zero.) The question of predicting when an integral around a closed contour has value zero will be discussed in Secs.~(\ref{sec:chap_4_7_42}),~(\ref{sec:chap_4_9_44}), and~(\ref{sec:chap_4_11_46}).
\end{example}

\begin{example}\label{ex:chap_4_5_4}
Let $C$ denote the semicircular path
\[
s=r(e_{xy}e^{i\theta}\cos\varphi+\sin\varphi)\texttt{ }(0\leq\theta\leq\pi)
\]
in a spatial plane about the origin from the point $s_{1}$ $(\theta=0)$ to the point $s_{2}$ $(\theta=\pi)$ where $r$ and $\varphi$ are constants. Although the branch (Sec.~(\ref{sec:chap_3_4_30}))
\begin{equation}\label{eq:chap_4_5_8}
f(s)=s^{1/2}=\sqrt{r}(e_{xy}e^{i\theta}\cos\varphi+\sin\varphi)^{1/2}\texttt{ }(r>0,0<\theta<2\pi)
\end{equation}
of the multiple-valued function $s^{1/2}$ is not defined at the initial point $s=s_{1}$ $(\theta=0)$ of the  contour $C$, the integral
\begin{equation}\label{eq:chap_4_5_9}
I=\int_{C}s^{1/2}ds
\end{equation}
of that branch nevertheless exists. For the integrand is piecewise continuous on $C$. To see that this is so, we observe that when $s(\theta)=r(e_{xy}e^{i\theta}\cos\varphi+\sin\varphi)$, the right-hand limits of the real and imaginary components of the function
\[
f[s(\theta)]=\sqrt{r}(e_{xy}\cos\theta\cos\varphi+\sin\varphi+ie_{xy}\sin\theta)^{1/2}\texttt{ }(0<\theta\leq\pi)
\]
at $\theta=0$ are $\sqrt{r}(e_{xy}\cos\varphi+\sin\varphi)^{1/2}$ and 0, respectively. Hence $f[s(\theta)]$ is continuous on the closed interval $0\leq\theta\leq\pi$ when its value at $\theta=0$ is defined as $\sqrt{r}(e_{xy}\cos\varphi+\sin\varphi)^{1/2}$. Consequently,
\[
\begin{split}
I &= \int_{0}^{\pi}\sqrt{r}(e_{xy}e^{i\theta}\cos\varphi+\sin\varphi)^{1/2}re_{xy}ie^{i\theta}\cos\varphi d\theta \\
  &= r^{3/2}\int_{0}^{\pi}(e_{xy}e^{i\theta}\cos\varphi+\sin\varphi)^{1/2}e_{xy}ie^{i\theta}\cos\varphi d\theta \\
  &= \frac{2}{3}r^{3/2}(e_{xy}e^{i\theta}\cos\varphi+\sin\varphi)^{3/2}|_{0}^{\pi} \\
  &=-\frac{2}{3}r^{3/2}[(e_{xy}\cos\varphi+\sin\varphi)^{3/2}+i(e_{xy}\cos\varphi-\sin\varphi)^{3/2}].
\end{split}
\]
Finally, then,
\[
I=-\frac{2}{3}r^{3/2}[(e_{xy}\cos\varphi+\sin\varphi)^{3/2}+i(e_{xy}\cos\varphi-\sin\varphi)^{3/2}].
\]
\end{example}

\section{Upper Bound for Moduli of Contour Integrals}\label{sec:chap_4_6_41}

When $C$ denotes a contour $s=s(t)(a\leq t\leq b)$, we know from definition~(\ref{eq:chap_4_4_2}), Sec.~(\ref{sec:chap_4_4_39}), and inequality~(\ref{eq:chap_4_2_4}) in Sec.~(\ref{sec:chap_4_2_37}) that
\[
|\int_{C}f(s)ds|=|\int_{s_{1}}^{s_{2}}f[s(t)]s'(t)dt|\leq\int_{s_{1}}^{s_{2}}|f[s(t)]||s'(t)|dt.
\]
So, for any nonnegative constant $M$ such that the values of $f$ on $C$ satisfy the inequality
$|f(s)|\leq M$,
\[
|\int_{C}f(s)ds|\leq M\int_{s_{1}}^{s_{2}}|s'(t)|dt.
\]
Since the integral on the right here represents the length $L$ of the contour (see Sec.~(\ref{sec:chap_4_3_38}),
it follows that the modulus of the value of the integral of $f$ along  $C$ does not exceed $ML$:
\begin{equation}\label{eq:chap_4_6_1}
|\int_{C}f(s)ds|\leq ML.
\end{equation}
This is, of course, a strict inequality when the values of $f$ on $C$ are such that $|f(s)|\leq M$.

Note that since all of the paths of integration to be considered here are contours
and the integrands are piecewise continuous functions defined on those contours, a
number $M$ such as the one appearing in inequality~(\ref{eq:chap_4_6_1}) will always exist. This is
because the real-valued function $|f[s(t)]|$ is continuous on the closed bounded interval
$a\leq t\leq b$ when $f$ is continuous on $C$; and such a function always reaches a maximum
value $M$ on that interval. Hence $|f(s)|$ has a maximum value on $C$ when $f$ is
continuous on it. It now follows immediately that the same is true when $f$ is piecewise
continuous on $C$.

\section{Antiderivatives}\label{sec:chap_4_7_42}

Although the value of a contour integral of a function $f(s)$ from a fixed point $s_{1}$ to a fixed point $s_{2}$ depends, in general, on the path that is taken, there are certain functions whose integrals from $s_{1}$ to $s_{2}$ have values that are independent of path. The theorem below is useful in determining when integration is independent of path and, moreover, when an integral around a closed path has value zero.

In proving the theorem, we shall discover an extension of the fundamental theorem of calculus that simplifies the evaluation of many contour integrals. That extension involves the concept of an antiderivative of a continuous function $f$ in a domain $D$, or a function $F$ such that $F'(s)=f(s)$ for all $s$ in $D$. Note that an antiderivative is, of necessity, an analytic function. Note, too, that an antiderivative of a given function $f$ is unique except for an additive complex constant. This is because the derivative of the difference $F(s)-G(s)$ of any two such antiderivatives $F(s)$ and $G(z)$ is zero; and, according to the theorem in Sec.~(\ref{sec:chap_2_13_23}), an analytic function is constant in a domain $D$ when its derivative is zero throughout $D$.

\begin{theorem}\label{th:chap_4_7_1}
Suppose that a function $f(s)$ is continuous on a domain $D$. If any one of the following statements is true, then so are the others:

(i) $f(s)$ has an antiderivative $F(s)$ in $D$;

(ii) the integrals of $f(s)$ along contours lying entirely in $D$ and extending from any
fixed point $s_{1}$ to any fixed point $s_{2}$ all have the same value;

(iii) the integrals of $f(s)$ around closed contours lying entirely in $D$ all have value
zero.
\end{theorem}

It should be emphasized that the theorem does not claim that any of these
statements is true for a given function $f$ and a given domain $D$. It says only that
all of them are true or that none of them is true. To prove the theorem, it is sufficient
to show that statement (i) implies statement (ii), that statement (ii) implies statement
(iii), and finally that statement (iii) implies statement (i).

Let us assume that statement (i) is true. If a contour $C$ from $s_{1}$ to $s_{2}$, lying in $D$, is
just a smooth arc, with parametric representation $s=s(t)(a\leq t\leq b)$, we know that
\[
\frac{d}{dt}F[s(t)]=F'[s(t)]s'(t)=f[s(t)]s'(t)\texttt{ }(a\leq t\leq b).
\]

Because the fundamental theorem of complex-valued functions of a real variable can be extended so as to apply to spatial complex-valued functions of a real variable (Sec.~(\ref{sec:chap_4_2_37})), it follows that.
\[
\int_{C}f(s)ds=\int_{a}^{b}f[s(t)]s'(t)dt=F[s(t)]|_{a}^{b}=F[s(b)]-F[s(a)].
\]
Since $s(b)=s_{2}$ and $s(a)=s_{1}$, the value of this contour integral is, then,
\[
F(s_{2})-F(s_{1});
\]
and that value is evidently independent of the contour $C$ as long as $C$ extends from $s_{1}$ to $s_{2}$ and lies entirely in $D$.  That is,
\begin{equation}\label{eq:chap_4_7_1}
\int_{a}^{b}f(s)ds=F(s_{2})-F(s_{1})=F(s)|_{s_{1}}^{s_{2}}
\end{equation}
when $C$ is smooth. Expression~(\ref{eq:chap_4_7_1}) is also valid when $C$ is any contour, not necessarily
a smooth one, that lies in $D$. For, if $C$ consists of a finite number of smooth arcs
$C_{k}(k=1,2,\ldots,n)$, each $C_{k}$ extending from a point $s_{k}$ to a point $s_{k+1}$ then
\[
\int_{C}f(s)ds=\sum_{k=1}^{n}\int_{Ck}f(s)ds=\sum_{k=1}^{n}[F(s_{k+1})-F(s_{k})]=F(s_{n+1})-F(s_{1}).
\]
The fact that statement (ii) follows from statement (i) is now established.

To see that statement (ii) implies statement (iii), we let $s_{1}$ and $s_{2}$ denote any two
points on a closed contour $C$ lying in $D$ and form two paths, each with initial point
$s_{1}$ and final point $s_{2}$, such that $C=C_{1}-C_{2}$. Assuming that statement (ii) is
true, one can write
\begin{equation}\label{eq:chap_4_7_2}
\int_{C_{1}}f(s)ds=\int_{C_{2}}f(s)ds,
\end{equation}
or
\begin{equation}\label{eq:chap_4_7_3}
\int_{C_{1}}f(s)ds+\int_{-C_{2}}f(s)ds=0.
\end{equation}
That is, the integral of $f(s)$ around the closed contour $C=C_{1}-C_{2}$ has value zero.

It remains to show that statement (iii) implies statement (i). We do this by assuming that statement (iii) is true, establishing the validity of statement (ii), and then arriving at statement (i). To see that statement (ii) is true, we let $C_{1}$ and $C_{2}$ denote any two contours, lying in $D$, from a point $s_{1}$ to a point $s_{2}$ and observe that, in view of statement (iii), equation~(\ref{eq:chap_4_7_3}) holds. Thus equation~(\ref{eq:chap_4_7_2}) holds. Integration is, therefore, independent of path in $D$; and we can define the function
\[
F(s)=\int_{s_{0}}^{s}f(s)ds
\]
on $D$. The proof of the theorem is complete once we show that $F'(s)=f(s)$ everywhere in $D$. We do this by letting $s+\Delta s$ be any point, distinct from $s$, lying in some neighborhood of $s$ that is small enough to be contained in $D$. Then
\[
F(s+\Delta s)-F(s)=\int_{s_{0}}^{s+\Delta s}f(s)ds-\int_{s_{0}}^{s}f(s)ds=\int_{s}^{s+\Delta s}f(s)ds,
\]
where the path of integration from $s$ to $s+\Delta s$ may be selected as a line segment. Since
\[
\int_{s}^{s+\Delta s}ds=\Delta s,
\]
we can write
\[
f(s)=\frac{1}{\Delta s}\int_{s}^{s+\Delta s}f(s)ds;
\]
and it follows that
\[
\frac{F(s+\Delta s)-F(s)}{\Delta s}-f(s)=\frac{1}{\Delta s}\int_{s}^{s+\Delta s}[f(\varsigma)-f(s)]d\varsigma.
\]
But $f$ is continuous at the point $s$. Hence, for each positive number $\varepsilon$, a positive number
$\delta$ exists such that
\[
|f(\varsigma)-f(s)|<\varepsilon\texttt{ whenever }|\varsigma-s|<\delta.
\]
Consequently, if the point $s+\Delta s$ is close enough to $s$ so that $|\Delta s|<\delta$, then
\[
|\frac{F(s+\Delta s)-F(s)}{\Delta s}-f(s)|<\frac{1}{|\Delta s|}\varepsilon|\Delta s|=\varepsilon;
\]
that is,
\[
\lim_{\Delta s \to 0}\frac{F(s+\Delta s)-F(s)}{\Delta s}=f(s)\texttt{ or }F'(s)=f(s).
\]

\section{Examples}\label{sec:chap_4_8_43}

The following examples illustrate the theorem in Sec.~(\ref{sec:chap_4_7_42}) and, in particular, the use of the extension~(\ref{th:chap_4_7_1}) of the fundamental theorem of calculus in that section.

\begin{example}\label{ex:chap_4_8_1}
The continuous function $f(s)=s^{2}$ has an antiderivative $F(s)=s^{3}/3$ throughout the $s$ space. Hence when $s_{1}=0$ and $s_{2}=e_{xy}(1+i)+1$
\[
\begin{split}
\int_{s_{1}}^{s_{2}}s^{2}ds &= \frac{s^{3}}{3}|_{s_{1}}^{s_{2}}=\frac{1}{3}[e_{xy}(1+i)+1]^{3} \\
                    &= \frac{1}{3}[e_{xy}^{3}(1+i)^{3}+3e_{xy}^{2}(1+i)^{2}+3e_{xy}(1+i)+1] \\
                    &= \frac{1}{3}\{e_{xy}[(1+i)^{3}+3(1+i)^{2}+3(1+i)]+1\} \\
                    &= \frac{1}{3}\{e_{xy}[(2+i)^{3}-1]+1\}
\end{split}
\]
for every contour from $s_{1}=0$ to $s_{2}=e_{xy}(1+i)+1$.
\end{example}

\begin{example}\label{ex:chap_4_8_2}
The function $f(s)=1/s^{2}$, which is continuous everywhere except at the origin, has an antiderivative $F(s)=-1/s$ in the domain $|s|>0$, consisting of the entire $s$ space with the origin deleted. Consequently,
\[
\int_{C}\frac{ds}{s^{2}}=0
\]
when $C$ is the positively oriented circle
\begin{equation}\label{eq:chap_4_8_1}
s=r(e_{xy}e^{i\theta}\cos\varphi+\sin\varphi)\texttt{ }(-\pi\leq\theta\leq\pi)
\end{equation}
in a spatial plane about the origin where $r$ and $\varphi$ are constants.

Note that the integral of the function $f(s)=1/s$ around the same circle cannot be evaluated in a similar way. For, although the derivative of any branch $F(s)$ of $\log s$ is $1/s$, $F(s)$ is not differentiable,  or even defined, along its branch cut. In particular, if a ray $\theta=\alpha$ from the origin is used to form the branch cut, $F'(s)$ fails to exist at the point where that ray intersects the circle $C$. So $C$  does not lie in a domain throughout which  $F'(s)=1/s$, and we cannot make direct use of an antiderivative. Example~(\ref{ex:chap_4_8_3}), just below, illustrates how a combination of two different antiderivatives can be used to evaluate  $f(s)=1/s$ around $C$.
\end{example}

\begin{example}\label{ex:chap_4_8_3}
Let $C_{1}$ denote the right half
\begin{equation}\label{eq:chap_4_8_2}
s=r(e_{xy}e^{i\theta}\cos\varphi+\sin\varphi)\texttt{ }(-\pi/2\leq\theta\leq\pi/2)
\end{equation}
of the circle $C$ in Example~(\ref{ex:chap_4_8_2}) where $r$ and $\varphi$ are constants. From equation~(\ref{eq:chap_3_3_5}) and lemma~(\ref{le:chap_3_3_1}) in Sec.~(\ref{sec:chap_3_3_29}), when $\sin\Phi>0$, the principal branch
\[
\begin{split}
Log\texttt{ }s &= \ln|s|+\log(e_{xy}e^{i\Theta}\cos\Phi+\sin\Phi)\texttt{ }(-\pi<\Theta<\pi) \\
               &= \ln|s|+e_{xy}[\log r_{L}(\Theta,\Phi)+i\theta_{L}(\Theta,\Phi)-\log\sin\Phi]+\log\sin\Phi
\end{split}
\]
where
\[
r_{L}(\Theta,\Phi)=\sqrt{1+\cos\Theta\sin2\Phi}\texttt{ and }
\theta_{L}(\Theta,\Phi)=\arctan\frac{\sin\Theta}{\cos\Theta+\tan\Phi}
\]
of the logarithmic function serves as an antiderivative of the function $1/s$ in the evaluation of the integral of $1/s$ along $C_{1}$:
\[
\begin{split}
\int_{C_{1}}\frac{ds}{s} &= \int_{\bar{s}_{0}}^{s_{0}}\frac{ds}{s}=Log\texttt{ }s|_{\bar{s}_{0}}^{s_{0}}=\ln|s|+\log(e_{xy}e^{i\Theta}\cos\Phi+\sin\Phi)|_{\bar{s}_{0}}^{s_{0}} \\
&= e_{xy}[\log r_{L}(\Theta,\Phi)+i\theta_{L}(\Theta,\Phi)]-e_{xy}[\log r_{L}(-\Theta,\Phi)+i\theta_{L}(-\Theta,\Phi)] \\
&= e_{xy}i2\arctan\frac{\sin\Theta}{\cos\Theta+\tan\Phi}.
\end{split}
\]
This integral was evaluated in another way in Example~(\ref{ex:chap_4_5_1}), Sec.~(\ref{sec:chap_4_5_40}), where representation~(\ref{eq:chap_4_8_2}) for the semicircle was used.

Next, let $C_{2}$ denote the left half
\begin{equation}\label{eq:chap_4_8_3}
s=r(e_{xy}e^{i\theta}\cos\varphi+\sin\varphi)\texttt{ }(\pi/2\leq\theta\leq3\pi/2)
\end{equation}
of the same circle $C$ and consider the branch
\[
\begin{split}
Log\texttt{ }s &= \ln|s|+\log(e_{xy}e^{i\Theta}\cos\Phi+\sin\Phi)\texttt{ }(r>0,0<\Theta<2\pi) \\
               &= \ln|s|+e_{xy}[\log r_{L}(\Theta,\Phi)+i\theta_{L}(\Theta,\Phi)-\log\sin\Phi]+\log\sin\Phi
\end{split}
\]
of the logarithmic function. One can write
\[
\begin{split}
\int_{C_{2}}\frac{ds}{s} &= \int_{\bar{s}}^{s}\frac{ds}{s}=Log\texttt{ }s|_{\bar{s}}^{s}=\ln|s|+\log(e_{xy}e^{i\Theta}\cos\Phi+\sin\Phi)|_{\bar{s}_{0}}^{s_{0}} \\
&= e_{xy}[\log r_{L}(\Theta,\Phi)+i\theta_{L}(\Theta,\Phi)]-e_{xy}[\log r_{L}(-\Theta,\Phi)+i\theta_{L}(-\Theta,\Phi)] \\
&= e_{xy}i2\arctan\frac{\sin\Theta}{\cos\Theta+\tan\Phi}.
\end{split}
\]

The value of the integral of $1/s$ around the entire circle $C=C_{1}+C_{2}$ is thus obtained:
\[
\begin{split}
\int_{C}\frac{ds}{s} &=\int_{C_{1}}\frac{ds}{s}+\int_{C_{2}}\frac{ds}{s} \\
                     &=e_{xy}i2\arctan\frac{\sin\Theta}{\cos\Theta+\tan\Phi}+e_{xy}i2\arctan\frac{\sin\Theta}{\cos\Theta+\tan\Phi}\\
                     &=e_{xy}i4\arctan\frac{\sin\Theta}{\cos\Theta+\tan\Phi}.
\end{split}
\]
\end{example}

\begin{example}\label{ex:chap_4_8_4}
Let us use an antiderivative to evaluate the integral
\begin{equation}\label{eq:chap_4_8_4}
\int_{C_{1}}s^{1/2}ds,
\end{equation}
where the integrand is the branch
\begin{equation}\label{eq:chap_4_8_5}
s^{1/2}=\sqrt{r}(e_{xy}e^{i\theta}\cos\varphi+\sin\varphi)^{1/2}\texttt{ }(r>0,0<\theta<2\pi)
\end{equation}
of the square root function where $r$ and $\varphi$ are constants, and where $C_{1}$ is any contour from $s_{1}=r(-e_{xy}\cos\varphi+\sin\varphi)$ to $s_{2}=r(e_{xy}\cos\varphi+\sin\varphi)$ that, except for its end points, lies above the $xz$ coordinate plane. Although the integrand is piecewise continuous on $C_{1}$, and the integral therefore exists, the branch~(\ref{eq:chap_4_8_5}) of $s^{1/2}$ is not defined on the ray $\theta=0$, in particular at the point $s=r(e_{xy}\cos\varphi+\sin\varphi)$. But another branch,
\[
f_{1}(s)=\sqrt{r}(e_{xy}e^{i\theta}\cos\varphi+\sin\varphi)^{1/2}\texttt{ }(r>0,-\pi/2<\theta<3\pi/2),
\]
is defined and continuous everywhere on $C_{1}$. The values of $f_{1}(s)$ at all points on $C_{1}$ except $s=r(e_{xy}\cos\varphi+\sin\varphi)$ coincide with those of our integrand~(\ref{eq:chap_4_8_5}); so the integrand can be replaced by $f_{1}(s)$. Since an antiderivative of $f_{1}(s)$ is the function
\[
F_{1}(s)=\frac{2}{3}s^{3/2}=\frac{2}{3}r\sqrt{r}(e_{xy}e^{i\theta}\cos\varphi+\sin\varphi)^{3/2}\texttt{ }(r>0,-\pi/2<\theta<3\pi/2),
\]
we can now write
\[
\begin{split}
\int_{C_{1}}s^{1/2}ds &=\int_{s_{1}}^{s_{2}}f_{1}(s)ds=F_{1}(s)|_{s_{1}}^{s_{2}}
=\frac{2}{3}r\sqrt{r}(e_{xy}e^{i\theta}\cos\varphi+\sin\varphi)^{3/2}|_{s_{1}}^{s_{2}} \\
&=\frac{2}{3}r\sqrt{r}[(e_{xy}e^{i\theta}\cos\varphi+\sin\varphi)^{3/2}-(-e_{xy}e^{i\theta}\cos\varphi+\sin\varphi)^{3/2}].
\end{split}
\]

The integral
\begin{equation}\label{eq:chap_4_8_6}
\int_{C_{2}}s^{1/2}ds
\end{equation}
of the function~(\ref{eq:chap_4_8_5}) over any contour $C_{2}$ that extends from $s_{1}=r(-e_{xy}\cos\varphi+\sin\varphi)$ to $s_{2}=r(e_{xy}\cos\varphi+\sin\varphi)$ below the $xz$ coordinate plane can be evaluated in a similar way. In this case, we can replace the integrand by the branch
\[
f_{2}(s)=\sqrt{r}(e_{xy}e^{i\theta}\cos\varphi+\sin\varphi)^{1/2}\texttt{ }(r>0,\pi/2<\theta<5\pi/2),
\]
whose values coincide with those of the integrand at $s=r(-e_{xy}\cos\varphi+\sin\varphi)$ and at all points on $C_{2}$ below the $xz$ coordinate plane. This enables us to use an antiderivative of $f_{2}(s)$ to evaluate integral~(\ref{eq:chap_4_8_6}).
\end{example}

\section{Cauchy-Goursat Theorem}\label{sec:chap_4_9_44}

In Sec.~(\ref{sec:chap_4_3_38}), we saw that when a continuous function $f$ has an antiderivative in a domain $D$, the integral of $f(s)$ around any given closed contour $C$ lying entirely in $D$ has value zero. In this section, we present a theorem giving other conditions on a function $f$, which ensure that the value of  the integral of $f(s)$ around a simple closed contour (Sec.~(\ref{sec:chap_4_3_38})) is zero. The theorem is central to the theory of functions of a spatial complex variable; and some extensions of it, involving certain special types of domains, will be given in Sec.~(\ref{sec:chap_4_11_46}).

We let $C$ denote a simple closed contour $s=s(t)$ $(a\leq t\leq b)$, described in the positive sense (counterclockwise), all points on the contour are in or on a sphere $S$ with a center point $s_{0}$ and a radius $R_{s}$, and we assume that $f$ is analytic at each point interior to and on $C$ and interior to and on $S$. According to Sec.~(\ref{sec:chap_4_4_39}),
\begin{equation}\label{eq:chap_4_9_1}
\int_{C}f(s)ds=\int_{a}^{b}f[s(t)]s'(t)dt;
\end{equation}
and if
\[
f(s)=e_{xy}[u(x,y,z)+iv(x,y,z)]+w(z)\texttt{ and }s(t)=e_{xy}[x(t)+iy(t)]+z(t),
\]
the integrand $f[s(t)]s'(t)$ in expression~(\ref{eq:chap_4_9_1}) is the product of the functions
\[
e_{xy}\{u[x(t),y(t),z(t)]+iv[x(t),y(t),z(t)]\}+w[z(t)]
\]
and
\[
e_{xy}[x'(t)+iy'(t)]+z'(t)
\]
of the real variable $t$. Thus
\begin{equation}\label{eq:chap_4_9_2}
\int_{C}f(s)ds=\int_{a}^{b}[e_{xy}(u+iv)+w][e_{xy}[x'(t)+iy'(t)]+z'(t)]dt
\end{equation}
\[
=e_{xy}\int_{a}^{b}[(ux'-vy')+i(vx'+uy')+(wx'+uz')+i(wy'+vz')]dt+\int_{a}^{b}wz'dt.
\]
In terms of line integrals of real-valued functions of three real variables, then,
\begin{equation}\label{eq:chap_4_9_3}
\begin{split}
\int_{C}f(s)ds &= e_{xy}(\int_{C}udx-vdy+i\int_{C}vdx+udy) \\
               &+ e_{xy}(\int_{C}wdx+udz+i\int_{C}wdy+vdz)+\int_{C}w dz.
\end{split}
\end{equation}
Observe that expression~(\ref{eq:chap_4_9_3}) can be obtained formally by replacing $f(s)$ and $ds$ on the
left with the trinomial
\[
e_{xy}(u+iv)+w\texttt{ and }e_{xy}(dx+idy)+dz,
\]
respectively, and expanding their product. Expression~(\ref{eq:chap_4_9_3}) is, of course, also valid when
$C$ is any contour, not necessarily a simple closed one, and $f[s(t)]$ is only piecewise
continuous on it.

We next recall a result from calculus that enables us to express the line integrals
\[
\int_{C}udx-vdy+i\int_{C}vdx+udy
\]
on the right in equation~(\ref{eq:chap_4_9_3}) as double integrals. Let the projection of the spatial contour on the $xy$ coordinate plane be a two-dimensional closed contour $C^{*}$ and a two-dimensional closed region $R$ consists of all points interior to and on the simple closed contour $C^{*}$. Then suppose that two real-valued functions $P(x,y)$ and $Q(x,y)$, together with their first-order partial derivatives of real variables $x$ and $y$, are continuous throughout the closed region $R$. According to Green's theorem,
\[
\int_{C}Pdx+Qdy=\iint_{R}(Q_{x}-P_{y})dA.
\]

Now $f$ is continuous in $R$, since it is analytic there. Hence the functions $u$ and $v$ are also continuous in $R$. Likewise, if the derivative $f'$ of $f$ is continuous in $R$, so are the first-order partial derivatives of $u$ and $v$. Green's theorem then enables us to rewrite equation~(\ref{eq:chap_4_9_3}) as
\begin{equation}\label{eq:chap_4_9_4}
\begin{split}
\int_{C}f(s)ds &= e_{xy}[\iint_{R}(-v_{x}-u_{y})dA+i\iint_{R}(u_{x}-v_{y})dA] \\
               &+ e_{xy}(\int_{C}wdx+udz+i\int_{C}wdy+vdz)+\int_{C}wdz.
\end{split}
\end{equation}

But, in view of the Cauchy-Riemann equations
\[
e_{xy}u_{x}=e_{xy}v_{y}\texttt{ and }e_{xy}u_{y}=-e_{xy}v_{x},
\]
the integrands of these two double integrals are zero throughout $S$ and $R$.

Moreover, let $z_{max}=\max\{z(t)|a\leq t\leq b\}$ and $z_{min}=\min\{z(t)|a\leq t\leq b\}$. Then there are
\[
\int_{C}w(z)dx=w(z)\int_{C}dx=0\texttt{, }\int_{C}w(z)dy=w(z)\int_{C}dy=0,
\]
and for $g=u,v,w$
\[
\begin{split}
\int_{C}gdz &= \int_{z_{min}}^{z_{max}}gdz+\int_{z_{max}}^{z_{min}}gdz \\
            &= \int_{z_{min}}^{z_{max}}gdz-\int_{z_{min}}^{z_{max}}gdz=0.
\end{split}
\]

So, when $f$ is analytic in $S$ and $f'$ is continuous there,
\begin{equation}\label{eq:chap_4_9_5}
\int_{C}f(s)ds=0.
\end{equation}

Note that, once it has been established that the value of this integral is zero, the orientation of $C$ is immaterial. That is, statement~(\ref{eq:chap_4_9_5}) is also true if $C$ is taken in the clockwise direction, since then
\[
\int_{C}f(s)ds=-\int_{-C}f(s)ds=0.
\]

Goursat was the first to prove that the condition of continuity on $f'$ can be omitted. Its removal is important and will allow us to show, for example, that the derivative $f'$ of an analytic function $f$ is analytic without having to assume the continuity of $f'$, which fo1lows as a consequence. We now state the revised form of Cauchy's result, known as the Cauchy-Goursat theorem.

\begin{theorem}\label{th:chap_4_9_1}
If a function $f$ is analytic at all points interior to a closed space $S$ and on a simple closed
contour $C$ where all points on the contour are in or on $S$, then
\begin{equation}\label{eq:cauchy_goursat}
\int_{C}f(s)ds=0.
\end{equation}
\end{theorem}

The proof is presented in the next section, where, to be specific, we assume that $C$ is positively oriented. The reader who wishes to accept this theorem without proof may pass directly to Sec.~(\ref{sec:chap_4_11_46}).

\section{Proof of the Cauchy-Goursat Theorem}\label{sec:chap_4_10_45}

We preface the proof of the Cauchy-Goursat theorem with a lemma. Let a positively oriented simple closed contour $C$ be entirely in a closed space $S$ consisting of all points on the contour $C$ together with all other boundary and interior points of $S$. We start by forming subsets of $S$. To do this, we draw equally spaced planes parallel to the $xy$, $yz$, and $zx$ coordinate planes such that the distance between adjacent planes which are parallel to each other is same. We thus form a finite number of closed cube subspaces, where each point of $S$ lies in at least one such subspace and each subspace contains points of $S$. We refer to these cube subspaces simply as cubes, always keeping in mind that by a cube we mean a boundary together with the points interior to it. If a particular cube contains points that are not in $S$, we remove those points and call what remains a partial cube. We thus cover or fill $S$ with a finite number of cubes and partial cubes, and our proof of the following lemma starts with this covering or filling.

\begin{lemma}\label{le:chap_4_10_1}
Let a positively oriented simple closed contour $C$ be entirely in a closed space $S$ consisting of all points on the contour $C$ together with all other boundary and interior points of $S$. Let $f$ be analytic throughout the closed space $S$. For any positive number $\varepsilon$, the closed space $S$ can be covered or filled with a finite number of cubes and partial cubes, indexed by $j=1, 2, \ldots, n$, such that in each one there is a fixed point $s_{j}$ for which the inequality
\begin{equation}\label{eq:chap_4_10_1}
|\frac{f(s)-f(s_{j})}{s-s_{j}}-f'(s_{j})|<\varepsilon\texttt{ }(s\neq s_{j})
\end{equation}
is satisfied by all other points in that cube or partial cube.
\end{lemma}

\begin{proof}
To start the proof, we consider the possibility that, in the covering or filling constructed just prior to the statement of the lemma, there is some cube or partial cube in which no point $s_{j}$ exists such that inequality~(\ref{eq:chap_4_10_1}) holds for all other points $s$ in it. If that subspace is a cube, we construct eight smaller cubes by drawing line segments joining the midpoints of its opposite sides. If the subspace is a partial cube, we treat the whole cube in the same manner and then let the portions that lie outside $S$ be discarded. If, in any one of these smaller subspaces, no point $s_{j}$ exists such that inequality~(\ref{eq:chap_4_10_1}) holds for all other points $s$ in it, we construct still smaller cubes and
partial cubes, etc. When this is done to each of the original subspaces that requires it, it turns out that, after a finite number of steps, the closed space $S$ can be covered or filled with a finite number of cubes and partial cubes such that the lemma is true.

To verify this, we suppose that the needed points $s_{j}$ do not exist after subdividing one of the original subspaces a finite number of times and reach a contradiction. We let $\sigma_{0}$ denote that subspace if it is a cube; if it is a partial cube, we let $\sigma_{0}$ denote the entire cube of which it is a part. After we subdivide $\sigma_{0}$, at least one of the eight smaller cubes, denoted by $\sigma_{1}$, must contain points of $S$ but no appropriate point $s_{j}$. We then subdivide $\sigma_{1}$ and continue in this manner. It may be that after a cube $\sigma_{k-1}$ $(k=1, 2, ... )$ has been subdivided, more than one of the eight smaller cubes constructed from it can be chosen. To make a specific choice, we take $\sigma_{k}$ to be the one lowest and then furthest to the left.

In view of the manner in which the nested infinite sequence
\begin{equation}\label{eq:chap_4_10_2}
\sigma_{0},\sigma_{1},\sigma_{2},\ldots,\sigma_{k-1},\sigma_{k},\ldots
\end{equation}
of cubes is constructed, it is easily known %(Exercise 9, Sec.~(\ref{sec:chap_4_11_46}))
that there is a point $s_{0}$ common to each $\sigma_{k}$; also, each of these cubes contains points of $S$ other than possibly $s_{0}$. Recall how the sizes of the cubes in the sequence are decreasing, and note that
any $\delta$ neighborhood $|s-s_{0}|<\delta$ of $s_{0}$ contains such cubes when their diagonals have lengths less than $\delta$. Every $\delta$ neighborhood $|s-s_{0}|<\delta$ therefore contains points of $S$ distinct from $s_{0}$, and this means that $s_{0}$ is an accumulation point of $S$. Since the space $S$ is a closed set, it follows that $s_{0}$ is a point in $S$.%  (See Sec. 10.)

Now the function $f$ is analytic throughout $S$ and, in particular, at $s_{0}$. Consequently, $f'(s_{0})$ exists, According to the definition of derivative (Sec.~(\ref{sec:chap_2_8_18})), there is, for each positive numbers $\varepsilon$, a $\delta$ neighborhood $|s-s_{0}|<\delta$ such that the inequality
\[
|\frac{f(s)-f(s_{0})}{s-s_{0}}-f'(s_{0})|<\varepsilon
\]
is satisfied by all points distinct from $s_{0}$ in that neighborhood. But the neighborhood $|s-s_{0}|<\delta$ contains a cube $\sigma_{K}$ when the integer $K$ is large enough that the length of a diagonal of that cube is less than $\delta$. Consequently, $s_{0}$ serves as the point $s_{j}$ in inequality~(\ref{eq:chap_4_10_1}) for the subspace consisting of the cube $\sigma_{K}$ or a part of $\sigma_{K}$. Contrary to the way in which the sequence~(\ref{eq:chap_4_10_2}) was formed, then, it is not necessary to subdivide $\sigma_{K}$. We thus arrive at a contradiction, and the proof of the lemma is complete.

%This completes the proof of the lemma.
\end{proof}

Continuing with a function $f$ which is analytic throughout a closed space $S$ consisting
of a positively oriented simple closed curve $C$ and points interior to $S$, we are now
to prove the Cauchy-Goursat theorem, namely that
\begin{equation}\label{eq:chap_4_10_3}
\int_{C}f(s)ds=0.
\end{equation}

Given an arbitrary positive number $\varepsilon$, we consider the covering or filling of $S$ in the statement of the lemma. Let us define on the $j$th cube or partial cube the following function, where $s_{j}$ is the fixed point in that subspace for which inequa1ity~(\ref{eq:cauchy_goursat}) holds:
\begin{equation}\label{eq:chap_4_10_4}
\delta_{j}(s)=\{
    \begin{array}{cc}
        \frac{f(s)-f(s_{j})}{s-s_{j}}-f'(s_{j}) & \texttt{ when }s\neq s_{j}, \\
        0 & \texttt{ when }s=s_{j}.
    \end{array}
\end{equation}
According to inequality~(\ref{eq:cauchy_goursat}),
\begin{equation}\label{eq:chap_4_10_5}
|\delta_{j}(s)|<\varepsilon
\end{equation}
at all points $s$ in the subspace on which $\delta_{j}(s)$ is defined. Also, the function $\delta_{j}(s)$ is
continuous throughout the subspace since $f(s)$ is continuous there and
\[
\lim_{s \to s_{j}}\delta_{j}(s)=f'(s_{j})-f'(s_{j})=0.
\]

Next, let $C_{j}$ $(j=1, 2, ... , n)$ denote the positively oriented boundaries of the above cubes or partial cubes covering or filling $S$. In view of definition~(\ref{eq:chap_4_10_4}), the value of $f$ at a point $s$ on any particular $C_{j}$ can be written
\[
f(s)=f(s_{j})-s_{j}f'(s_{j})+f'(s_{j})s+(s-s_{j})\delta_{j}(s);
\]
and this means that
\begin{equation}\label{eq:chap_4_10_6}
\int_{C_{j}}f(s)ds=[f(s_{j})-s_{j}f'(s_{j})]\int_{C_{j}}ds+f'(s_{j})\int_{C_{j}}sds+\int_{C_{j}}(s-s_{j})\delta_{j}(s)ds.
\end{equation}
But
\[
\int_{C_{j}}ds=0\texttt{ and }\int_{C_{j}}sds=0
\]
since the functions $1$ and $s$ possess antiderivatives everywhere in the finite space. So
equation~(\ref{eq:chap_4_10_6}) reduces to
\begin{equation}\label{eq:chap_4_10_7}
\int_{C_{j}}f(s)ds=\int_{C_{j}}(s-s_{j})\delta_{j}(s)ds\texttt{ }(j=1, 2, ... , n).
\end{equation}

The sum of all $n$ integrals on the left in equations~(\ref{eq:chap_4_10_7}) can be written
\[
\sum_{j=1}^{n}\int_{C_{j}}f(s)ds=\int_{C}f(s)ds.
\]
since the two integrals along the common boundary of every pair of adjacent subspaces cancel each other, the integral being taken in one sense along that line segment in one subspace and in the opposite sense in the other. Only the integrals along the arcs that are parts of $C$ remain. Thus, in view of equations~(\ref{eq:chap_4_10_7}),
\[
\int_{C}f(s)ds=\sum_{j=1}^{n}\int_{C_{j}}(s-s_{j})\delta_{j}(s)ds;
\]
and so
\begin{equation}\label{eq:chap_4_10_8}
|\int_{C}f(s)ds|\leq\sum_{j=1}^{n}|\int_{C_{j}}(s-s_{j})\delta_{j}(s)ds|.
\end{equation}

Let us now use property~(\ref{eq:chap_4_6_1}), Sec.~(\ref{sec:chap_4_6_41}) to find an upper bound for each absolute value on the right in inequality~(\ref{eq:chap_4_10_8}). To do this, we first recall that each $C_{j}$ coincides either entirely or partially with the boundary of a cube. In either case, we let $h_{j}$ denote the length of a side of the cube. Since, in the $j$th integral, both the variable $s$ and the point $s_{j}$ lie in that cube,
\[
|s-s_{j}|\leq\sqrt{3}h_{j}.
\]
In view of inequality~(\ref{eq:chap_4_10_5}), then, we know that each integrand on the right in inequality~(\ref{eq:chap_4_10_8}) satisfies the condition
\begin{equation}\label{eq:chap_4_10_9}
|(s-s_{j})\delta_{j}(s)|\leq\sqrt{3}h_{j}\varepsilon.
\end{equation}
As for the length of the path $C_{j}$, it is $8h_{j}$ if $C_{j}$ is the boundary of a cube. In that case,
we let $A_{j}$ denote the area of one side face of the cube and observe that
\begin{equation}\label{eq:chap_4_10_10}
|\int_{C_{j}}(s-s_{j})\delta_{j}(s)ds|<\sqrt{3}h_{j}\varepsilon8h_{j}=8\sqrt{3}A_{j}\varepsilon.
\end{equation}
If $C_{j}$ is the boundary of a partial cube, its length does not exceed $8h_{j}+L_{j}$ where
$L_{j}$ is the length of that part of $C_{j}$ which is also a part of $C$. Again letting $A_{j}$ denote
the area of one side face of the full cube, we find that
\begin{equation}\label{eq:chap_4_10_11}
|\int_{C_{j}}(s-s_{j})\delta_{j}(s)ds|<\sqrt{3}h_{j}\varepsilon(8h_{j}+L_{j})<8\sqrt{3}A_{j}\varepsilon+\sqrt{3}HL_{j}\varepsilon.
\end{equation}
where $H$ is the length of a side of some cube that encloses the entire contour $C$ as well as all of the cubes originally used in covering or filling $S$. Note that the sum of all the $A_{j}$ does not exceed $H^{2}$.

If $L$ denotes the length of $C$, it now follows from inequalities~(\ref{eq:chap_4_10_8}),~(\ref{eq:chap_4_10_10}), and~(\ref{eq:chap_4_10_11}) that
\begin{equation}
|\int_{C}f(s)ds|<(8\sqrt{3}H^{2}+\sqrt{3}HL)\varepsilon.
\end{equation}
Since the value of the positive numbers $\varepsilon$ is arbitrary, we can choose it so that the right-hand side of this last inequality is as small as we please. The left-hand side, which is independent of $s$, must therefore be equal to zero; and statement~(\ref{eq:chap_4_10_3}) fol1ows. This completes the proof of the Cauchy-Goursat theorem.

Thus the Cauchy-Goursat theorem of integrals of one spatial complex variable also hold.

\section{Simply and Multiply Connected Domains}\label{sec:chap_4_11_46}

A simply connected domain $D$ is a domain such that every simple closed spatial contour within
it encloses only points of $D$. The set of points interior to a simple closed spatial contour is an
example. The spherical shell domain between two concentric spheres is, however, not simply
connected. A domain that is not simply connected is said to be multiply connected.

The Cauchy-Goursat theorem can be extended in the following way, involving a simply connected domain.

\begin{theorem}\label{th:chap_4_11_1}
If a function $f$ is analytic throughout a simply connected domain $D$, then
\begin{equation}\label{eq:chap_4_11_1}
\int_{C}f(s)ds=0
\end{equation}
for every closed contour $C$ lying in $D$.
\end{theorem}

The proof is easy if $C$ is a simple closed contour or if it is a closed contour that intersects itself a finite number of times. For, if $C$ is simple and lies in $D$, the function $f$ is analytic at each point interior to and on $C$; and the Cauchy-Goursat theorem ensures that equation~(\ref{eq:chap_4_11_1}) holds. Furthermore, if $C$ is closed but intersects itself a finite number of times, it consists of a finite number of simple closed contours. Let the finite number be $n$. Then the simple closed curves $C_{k} (k=1, 2, \ldots, n)$ make up $C$. Since the value of the integral around each $C_{k}$ is zero, according to the Cauchy-Goursat theorem, it follows that
\[
\int_{C}f(s)ds=\sum_{k=1}^{n}\int_{C_{k}}f(s)ds=0.
\]

\begin{corollary}\label{co:chap_4_11_1}
A function $f$ that is analytic throughout a simply connected domain $D$
must have an antiderivative everywhere in $D$.
\end{corollary}

This corollary follows immediately from Theorem~(\ref{th:chap_4_11_1}) because of the theorem in Sec.~(\ref{sec:chap_4_7_42}), which tells us that a continuous function $f$ always has an antiderivative in a given domain when equation~(\ref{eq:chap_4_11_1}) holds for each closed contour $C$ in that domain. Note that, since the finite space is simply connected, Corollary~(\ref{co:chap_4_11_1}) tells us that entire functions always possess antiderivatives.

The Cauchy-Goursat theorem can also be extended in a way that involves integrals along the boundary of a multiply connected domain. The following theorem is such an extension.

\begin{theorem}\label{th:chap_4_11_2}
Suppose that

(i) $C$ is a simple closed contour, described in the counterclockwise direction;

(ii) $C_{k} (k =1, 2, ... , n)$ are simple closed contours interior to $C$, all described in the clockwise direction, that are disjoint and whose interiors have no points in common.

If a function $f$ is analytic on all of these contours and throughout the multiply connected domain consisting of  all points inside $C$ and exterior to each $C_{k}$, then
\begin{equation}\label{eq:chap_4_11_2}
\int_{C}f(s)ds=\sum_{k=1}^{n}\int_{C_{k}}f(s)ds=0.
\end{equation}
\end{theorem}

Note that, in equation~(\ref{eq:chap_4_11_2}), the direction of each path of integration is such that
the multiply connected domain lies to the left of that path.

To prove the theorem, we introduce a polygonal path $L_{1}$, consisting of a finite number of line segments joined end to end, to connect the outer contour $C$ to the inner contour $C_{1}$. We introduce another polygonal path $L_{2}$ which connects $C_{1}$ to $C_{2}$; and we continue in this manner, with $L_{n+l}$ connecting $C_{n}$ to $C$. Two simple closed contours $\Gamma_{1}$ and $\Gamma_{2}$ can be formed, each consisting of polygonal paths $L_{k}$ or $-L_{k}$ and pieces of $C$ and $C_{k}$ and each described in such a direction that the points enclosed by them lie to the left. The Cauchy-Goursat theorem can now be applied to $f$ on  $\Gamma_{1}$ and $\Gamma_{2}$, and the sum of the values of the integrals over those contours is found to be zero. Since the integrals in opposite directions along each path $L_{k}$ cancel, only the integrals along $C$ and $C_{k}$ remain; and we arrive at statement~(\ref{eq:chap_4_11_2}).

The following corollary is an especially important consequence of Theorem~(\ref{th:chap_4_11_2}).

Given a positively oriented simple closed contour $C$, let all points on the contour $C$ be in or on a minimum closed space $S$. So two positively oriented simple closed contours $C_{1}$ and $C_{2}$ correspond to two spaces $S_{1}$ and $S_{2}$, respectively. If all points on the contour $C_{2}$ are inside $S_{1}$ and all points on the contour $C_{1}$ are outside $S_{2}$, then we say that $C_{2}$ is interior to $C_{1}$.

\begin{corollary}\label{co:chap_4_11_2}
Let $C_{1}$ and $C_{2}$ denote positively oriented simple closed contours corresponding to two spaces $S_{1}$ and $S_{2}$, respectively, where $C_{2}$ is interior to $C_{1}$. If a function $f$ is analytic in the closed space consisting
of those contours and all points between them, then
\begin{equation}\label{eq:chap_4_11_3}
\int_{C_{1}}f(s)ds=\int_{C_{2}}f(s)ds.
\end{equation}
\end{corollary}

For a verification, we use Theorem~(\ref{th:chap_4_11_2}) to write
\[
\int_{C_{1}}f(s)ds+\int_{-C_{2}}f(s)ds=0
\]
and we note that this is just a different form of equation~(\ref{eq:chap_4_11_3}).

Corollary~(\ref{co:chap_4_11_2}) is known as the principle of deformation of paths since it tells us that if $C_{1}$ is continuously deformed into $C_{2}$, always passing through points at which $f$ is analytic, then the value of  the integral of $f$ over $C_{1}$ never changes.

\section{Cauchy Integral Formula}\label{sec:chap_4_12_47}

Another fundamental result will now be established.

\begin{theorem}\label{th:chap_4_12_0}
Let $f$ be analytic everywhere inside and on a simple closed contour $C$ that is in a plane parallel to the $xy$ coordinate plane, taken in the positive sense. If $s_{0}$ is any point interior to $C$, then
\begin{equation}\label{eq:chap_4_12_0}
f(s_{0})=\frac{1}{2\pi i}\int_{C}\frac{f(s)ds}{s-s_{0}}.
\end{equation}
\end{theorem}

Because $C$ is in a plane parallel to the $xy$ coordinate plane, according to the two-dimensional complex variable theory, formula~(\ref{th:chap_4_12_0}) is called the Cauchy integral formula for two-dimensional complex planes. It tells us that if a function $f$ is to be analytic within and on a simple closed contour $C$ in a plane parallel to the $xy$ coordinate plane, then the values of $f$ interior to $C$ are completely determined by the values of $f$ on $C$.

\begin{theorem}\label{th:chap_4_12_1}
Let $f$ be analytic everywhere inside and on a sphere $S$ and a simple closed contour $C$ in or/and on $S$ be in a plane passing through the origin, taken in the positive sense. If $s_{0}$ is any point in the plane and interior to $C$, then
\begin{equation}\label{eq:chap_4_12_1}
f(s_{0})=\frac{1}{2\pi i}\int_{C}\frac{f(s)ds}{s-s_{0}}.
\end{equation}
\end{theorem}

Formula~(\ref{th:chap_4_12_1}) is called the Cauchy integral formula for three-dimensional complex planes. It tells us that if a function $f$ is to be analytic within and on a simple closed contour $C$, then the values of $f$ interior to $C$ are completely determined by the values of $f$ on $C$.

When the Cauchy integral formula is written
\begin{equation}\label{eq:chap_4_12_2}
\int_{C}\frac{f(s)ds}{s-s_{0}}=2\pi i f(s_{0}),
\end{equation}
it can be used to evaluate certain integrals along simple closed contours.

We begin the proof of the theorem by letting $C_{\rho}$ denote a positively oriented circle $|s-s_{0}|=\rho$, which and $C$ are in a same plane, where $\rho$ is small enough that $C_{\rho}$ interior to $C$, i,e., all points on $C$ are outside $C_{\rho}$. Since the function $f(s)/(s-s_{0})$ is analytic between and on the contours $C$ and $C_{\rho}$, it follows from the principle of deformation of paths (Corollary~(\ref{co:chap_4_11_2}), Sec.~(\ref{sec:chap_4_11_46})) that
\[
\int_{C}\frac{f(s)ds}{s-s_{0}}=\int_{C_{\rho}}\frac{f(s)ds}{s-s_{0}}.
\]
This enables us to write
\begin{equation}\label{eq:chap_4_12_3}
\int_{C}\frac{f(s)ds}{s-s_{0}}-f(s_{0})\int_{C_{\rho}}\frac{ds}{s-s_{0}}=\int_{C_{\rho}}\frac{f(s)-f(s_{0})}{s-s_{0}}ds.
\end{equation}
According to the two-dimensional complex theory
\[
\int_{C_{\rho}}\frac{ds}{s-s_{0}}=2\pi i;
\]
and so equation~(\ref{eq:chap_4_12_3}) becomes
\begin{equation}\label{eq:chap_4_12_4}
\int_{C}\frac{f(s)ds}{s-s_{0}}-2\pi if(s_{0})=\int_{C_{\rho}}\frac{f(s)-f(s_{0})}{s-s_{0}}ds.
\end{equation}

Now the fact that $f$ is analytic, and therefore continuous, at $s_{0}$ ensures that, corresponding to each positive number $\varepsilon$, however small, there is a positive number $\delta$ such that
\begin{equation}\label{eq:chap_4_12_5}
|f(s)-f(s_{0}|<\varepsilon\texttt{ whenever }|s-s_{0}|<\delta.
\end{equation}
Let the radius $\rho$ of the circle $C_{\rho}$ be smaller than the number $\delta$ in the second of these inequalities. Since $|s-s_{0}|=\rho$ when $s$ is on $C_{\rho}$ it follows that the first of inequalities~(\ref{eq:chap_4_12_5}) holds when $s$ is such a point; and inequality~(\ref{eq:chap_4_6_1}), Sec.~(\ref{sec:chap_4_6_41}), giving upper bounds for the moduli of curve integrals, tells us that
\[
|\int_{C_{\rho}}\frac{f(s)-f(s_{0})}{s-s_{0}}ds|<\frac{\varepsilon}{\rho}2\pi\rho=2\pi\varepsilon.
\]
In view of equation~(\ref{eq:chap_4_12_4}), then,
\[
|\int_{C}\frac{f(s)-f(s_{0})}{s-s_{0}}ds-2\pi if(s_{0})|<2\pi\varepsilon.
\]
Since the left-hand side of this inequality is a nonnegative constant that is less than an arbitrarily small positive number, it must equal to zero. Hence equation~(\ref{eq:chap_4_12_2}) is valid,
and the theorem is proved.

\section{Derivatives of Analytic Functions}\label{sec:chap_4_13_48}

It follows from the Cauchy integral formula (Sec.~(\ref{sec:chap_4_12_47})) that if a function is analytic at a point, then its derivatives of all orders exist at that point and are themselves analytic there. To prove this, we start with a lemma that extends the Cauchy integral formula so as to apply to derivatives of the first and second order.

\begin{lemma}\label{le:chap_4_13_1}
Suppose that a function $f$ is analytic everywhere inside and on a simple closed contour $C$, taken in the positive sense. If $s$ is any point interior to $C$, then
\begin{equation}\label{eq:chap_4_13_1}
f'(s)=\frac{1}{2\pi i}\int_{C}\frac{f(\varsigma)d\varsigma}{(\varsigma-s)^{2}}\texttt{ and }
f''(s)=\frac{1}{\pi i}\int_{C}\frac{f(\varsigma)d\varsigma}{(\varsigma-s)^{3}}.
\end{equation}
\end{lemma}

Note that expressions~(\ref{eq:chap_4_13_1}) can be obtained formally, or without rigorous verification, by differentiating with respect to $s$ under the integral sign in the Cauchy integral formula
\begin{equation}\label{eq:chap_4_13_2}
f(s)=\frac{1}{2\pi i}\int_{C}\frac{f(\varsigma)d\varsigma}{\varsigma-s}
\end{equation}
where $s$ is interior to $C$ and $\varsigma$ denotes points on $C$.

To verify the first of expressions~(\ref{eq:chap_4_13_1}), we let $d$ denote the smallest distance from
$s$ to points on  $C$ and use formula~(\ref{eq:chap_4_13_2}) to write
\[
\begin{split}
\frac{f(s+\Delta s)-f(s)}{\Delta s} &= \frac{1}{2\pi i}\int_{C}(\frac{1}{\varsigma-s-\Delta s}-\frac{1}{\varsigma-s})\frac{f(\varsigma)}{\Delta s}d\varsigma \\
&= \frac{1}{2\pi i}\int_{C}\frac{f(\varsigma)d\varsigma}{(\varsigma-s-\Delta s)(\varsigma-s)}
\end{split}
\]
where $0<|\Delta s|<d$. Evidently, then,
\begin{equation}\label{eq:chap_4_13_3}
\frac{f(s+\Delta s)-f(s)}{\Delta s}-\frac{1}{2\pi i}\int_{C}\frac{f(\varsigma)d\varsigma}{(\varsigma-s)^{2}}=\frac{1}{2\pi i}\int_{C}\frac{\Delta s f(\varsigma)d\varsigma}{(\varsigma-s-\Delta s)(\varsigma-s)^{2}}
\end{equation}
Next, we let $M$ denote the maximum value of $|f(s)|$ on $C$ and observe that, since
$|\varsigma-s|\geq d$ and $|\Delta s|<d$,
\[
\begin{split}
|\varsigma-s-\Delta s| &= |(\varsigma-s)-\Delta s| \\
                       &\geq||\varsigma-s)|-|\Delta s||\geq d-|\Delta s|.
\end{split}
\]
Thus
\[
|\int_{C}\frac{\Delta s f(\varsigma)d\varsigma}{(\varsigma-s-\Delta s)(\varsigma-s)^{2}}|\leq\frac{|\Delta s|M}{(d-|\Delta s|)d^{2}}L,
\]
where $L$ is the length of $C$. Upon letting $|\Delta s|$ tend to zero, we find from this inequality
that the right-hand side of equation~(\ref{eq:chap_4_13_3}) also tends to zero. Consequently,
\[
\lim_{\Delta s \to 0}\frac{f(s+\Delta s)-f(s)}{\Delta s}-\frac{1}{2\pi i}\int_{C}\frac{f(\varsigma)d\varsigma}{(\varsigma-s)^{2}}=0;
\]
and the desired expression for $f'(s)$ is established.

The same technique can be used to verify the expression for $f''(s)$ in the statement of the lemma.

\begin{theorem}\label{th:chap_4_13_1}
If a function is analytic at a point, then its derivatives of all orders exist at that point. Those derivatives are, moreover, all analytic there.
\end{theorem}

To prove this remarkable theorem, we assume that a function $f$ is analytic at a point $s_{0}$. There must, then, be a neighborhood $|s-s_{0}|<\varepsilon$ of $s_{0}$ throughout which $f$ is analytic (see Sec.~(\ref{sec:chap_2_13_23})). Consequently, there is a positively oriented spatial circle $C_{0}$, centered at $s_{0}$ and with radius $\varepsilon/2$,  such that $f$ is analytic inside and on $C_{0}$. According to the above lemma,
\[
f''(s)=\frac{1}{\pi i}\int_{C}\frac{f(\varsigma)d\varsigma}{(\varsigma-s)^{3}}
\]
at each point $s$ interior to $C_{0}$, and the existence of $f''(s)$ throughout the neighborhood
$|s-s_{0}|<\varepsilon/2$ means that $f'$ is analytic at $s_{0}$. One can apply the same argument to the analytic function $f'$ to conclude that its derivative $f''$ is analytic, etc. Theorem~(\ref{th:chap_4_13_1}) is now established.

As a consequence, when a function
\[
f(s)=e_{xy}[u(x,y,z)+iv(x,y,z)]+w(z)
\]
is analytic at a point $s=e_{xy}(x+iy)+z$, the differentiability of $f'$ ensures the continuity of $f'$ there (Sec.~(\ref{sec:chap_2_8_18})). Then, since
\[
\begin{split}
f'(s) &= e_{xy}(u_{x}+iv_{x}) \\
      &= e_{xy}(v_{y}-iu_{y}) \\
      &= e_{xy}(u_{z}+iv_{z})+w_{z}
\end{split}
\]
(Sec.~(\ref{sec:chap_2_10_20})) we may conclude that the first-order partial derivatives of $u$, $v$, and $w$ are continuous at
that point. Furthermore, since $f''$ is analytic and continuous at $s$ and since
\[
\begin{split}
f''(s) &= e_{xy}(u_{xx}+iv_{xx}) \\
       &= e_{xy}(v_{yx}-iu_{yx}) \\
       &= e_{xy}(u_{zx}+iv_{zx})
\end{split}
\]
etc., we arrive at a corollary that was anticipated in Sec.~(\ref{sec:chap_2_15_25}), where harmonic functions were introduced.

\begin{corollary}\label{co:chap_4_13_1}
If a function
\[
f(s)=e_{xy}[u(x,y,z)+iv(x,y,z)]+w(z)
\]
is defined and analytic at a point
\[
s=e_{xy}(x+iy)+z
\]
then the component functions $u$, $v$, and $w$ have continuous partial derivatives of all orders at that point.
\end{corollary}

One can use mathematical induction to generalize formulas~(\ref{eq:chap_4_13_1}) to
\begin{equation}\label{eq:chap_4_13_4}
f^{(n)}(s)=\frac{n!}{2\pi i}\int_{C}\frac{f(\varsigma)d\varsigma}{(\varsigma-s)^{n+1}}\texttt{ }(n=1,2,\ldots).
\end{equation}
The verification is considerably more involved than for just $n=1$ and $n=2$. Note that, with the agreement that
\[
f^{(0)}(s)=f(s)\texttt{ and }0!=1,
\]
expression~(\ref{eq:chap_4_13_4}) is also valid when $n=0$, in which case it becomes the Cauchy integral
formula~(\ref{eq:chap_4_13_2}).

When written in the form
\begin{equation}\label{eq:chap_4_13_5}
\int_{C}\frac{f(s)ds}{(s-s_{0})^{n+1}}=\frac{2\pi i}{n!}f^{(n)}(s_{0})\texttt{ }(n=1,2,\ldots).
\end{equation}
expression~(\ref{eq:chap_4_13_4}) can be useful in evaluating certain integrals when $f$ is analytic inside
and on a simple closed contour $C$, taken in the positive sense, and $s_{0}$ is any point interior to $C$.
%It has already been illustrated in Sec.~(\ref{sec:chap_4_12_47}) when $n=0$.

We conclude this section with a theorem due to E. Morera (1856-1909). The proof here depends on the fact that the derivative of an analytic function is itself analytic, as stated in Theorem~(\ref{th:chap_4_13_1}).

\begin{theorem}\label{th:chap_4_13_2}
Let $f$ be continuous on a domain $D$. If
\begin{equation}\label{eq:chap_4_13_6}
\int_{C}f(s)ds=0
\end{equation}
for every closed contour $C$ lying in $D$, then $f$ is analytic throughout $D$.
\end{theorem}

In particular, when $D$ is simply connected, we have for the class of continuous functions on $D$ a converse of Theorem~(\ref{th:chap_4_11_1}) in Sec.~(\ref{sec:chap_4_11_46}), which is the extension of the Cauchy-Goursat theorem involving such domains.

To prove the theorem here, we observe that when its hypothesis is satisfied, the theorem in Sec.~(\ref{sec:chap_4_7_42}) ensures that $f$ has an antiderivative in $D$; that is, there exists an analytic function $F$ such that $F'(s)=f(s)$ at each point in $D$. Since $f$ is the derivative of $F$, it then follows from Theorem~(\ref{th:chap_4_13_1}) above that $f$ is analytic in $D$.

\section{Liouville's Theorem and the Fundamental Theorem of Algebra}\label{sec:chap_4_14_49}

This section is devoted to two important theorems that follow from the extension of the Cauchy integral formula in Sec.~(\ref{sec:chap_4_13_48}).

\begin{lemma}\label{le:chap_4_14_1}
Suppose that a function $f$ is analytic inside and on a sphere $S_{R}$ centered at $s_{0}$ and with radius $R$. If $M_{R}$ denotes the maximum
value of $|f(s)|$ on $S_{R}$, then
\begin{equation}\label{eq:chap_4_14_1}
|f^{(n)}(s_{0})|\leq\frac{n!M_{R}}{R^{n}}\texttt{ }(n=1,2,\ldots).
\end{equation}
\end{lemma}

Inequality~(\ref{eq:chap_4_14_1}) is called Cauchy's inequality and is an immediate consequence of the expression
\[
f^{(n)}(s_{0})=\frac{n!}{2\pi i}\int_{C_{R}}\frac{f(s)ds}{(s-s_{0})^{n+1}}\texttt{ }(n=1,2,\ldots)
\]
where $C_{R}$ is a positively oriented circle on the sphere $S_{R}$, which is a slightly different form of equation~(\ref{eq:chap_4_13_5}), Sec.~(\ref{sec:chap_4_13_48}). We need only apply inequality~(\ref{eq:chap_4_6_1}), Sec.~(\ref{sec:chap_4_6_41}), which gives upper bounds for the moduli of the values of contour integrals, to see that
\[
f^{(n)}(s_{0})\leq\frac{n!}{2\pi}\cdot\frac{M_{R}}{R^{n+1}}2\pi R\texttt{ }(n=1,2,\ldots).
\]
where $M_{_{R}}$ is as in the statement of the lemma. This inequality is, of  course, the same
as inequality~(\ref{eq:chap_4_14_1}) in the lemma.

The lemma can be used to show that no entire function except a constant is bounded in the complex space. Our first theorem here, which is known as Liouville's theorem,  states this result in a somewhat different way.

\begin{theorem}\label{th:chap_4_14_1}
If $f$ is entire and bounded in the complex space, then $f(s)$ is constant throughout the space.
\end{theorem}

To start the proof, we assume that $f$ is as stated in the theorem and note that, since $f$ is entire, Cauchy's inequality~(\ref{eq:chap_4_14_1}) with $n=1$ holds for any choices of $s_{0}$ and $R$:
\begin{equation}\label{eq:chap_4_14_2}
|f'(s_{0})|\leq\frac{M_{R}}{R}.
\end{equation}
Moreover, the boundedness condition in the statement of the theorem tells us that a nonnegative constant $M$ exists such that $|f(s)|<M$ for all $s$; and, because the constant
$M_{R}$ in inequality~(\ref{eq:chap_4_14_2}) is always less than or equal to $M$, it follows that
\begin{equation}\label{eq:chap_4_14_3}
|f'(s_{0})|\leq\frac{M}{R},
\end{equation}
where $s_{0}$ is any fixed point in the space and $R$ is arbitrarily large. Now the number $M$ in inequality~(\ref{eq:chap_4_14_3}) is independent of the value of $R$ that is taken. Hence that inequality can hold for arbitrarily large values of $R$ only if $f'(s_{0})=0$. Since the choice of $s_{0}$ was arbitrary, this means that $f'(s)=0$ everywhere in the complex space. Consequently, $f$ is a constant function, according to the theorem in Sec.~(\ref{sec:chap_2_13_23}).

The following theorem, known as the fundamental theorem of algebra, follows readily from Liouville's theorem.

\begin{theorem}\label{th:chap_4_14_2}
Any  polynomial
\[
P(s)=a_{0}+a_{1}s+a_{2}s^{2}+\cdots+a_{n}s^{n}\texttt{ }(a_{n}\neq0)
\]
of degree $n$ $(n>1)$ has at least one zero. That is, there exists at least one point $s_{0}$ such that $P(s_{0})=0$.
\end{theorem}

The proof here is by contradiction. Suppose that  $P(s)$ is not zero for any value of $s$. Then the reciprocal
\[
f(s)=\frac{1}{P(s)}
\]
is clearly entire, and it is also bounded in the complex space.

To show that it is bounded, we first write
\begin{equation}\label{eq:chap_4_14_4}
\varpi=\frac{a_{0}}{s^{n}}+\frac{a_{1}}{s^{n-1}}+\frac{a_{2}}{s^{n-2}}+\cdots+\frac{a_{n-1}}{s}\texttt{ }(a_{n}\neq0),
\end{equation}
so that $P(s)=(a_{n}+\varpi)s^{n}$. We then observe that a sufficiently large positive number
$R$ can be found such that the modulus of each of the quotients in expression~(\ref{eq:chap_4_14_3}) is less
than the number $|a_{n}|/(2n)$ when $|s|\geq R$. The generalized triangle inequality, applied
to $n$ complex numbers, thus shows that $|\varpi|<|a_{n}|/2$ for such values of $s$. Consequently,
when $|s|\geq R$,
\[
|a_{n}+\varpi|\geq||a_{n}|-|\varpi||>\frac{|a_{n}|}{2},
\]
and this enables us to write
\begin{equation}\label{eq:chap_4_14_5}
|P(s)|=|a_{n}+\varpi||s^{n}|>\frac{|a_{n}|}{2}|s|^{n}\geq\frac{|a_{n}|}{2}R^{n}\texttt{ whenever }(|s|\geq R).
\end{equation}
Evidently, then,
\[
|f(s)|=\frac{1}{|P(s)|}<\frac{2}{|a_{n}|R^{n}}\texttt{ whenever }(|s|\geq R).
\]
So $f$ is bounded in the region exterior to the sphere $|s|\leq R$. But $f$ is continuous in that closed sphere, and this means that $f$ is bounded there too. Hence $f$ is bounded in the entire space.

It now follows from Liouville's theorem that $f(s)$, and consequently $P(s)$, is constant. But $P(s)$ is not constant, and we have reached a contradiction.

The fundamental theorem tells us that any polynomial $P(s)$ of degree $n$ $(n>1)$ can be expressed as a product of  linear factors:
\begin{equation}\label{eq:chap_4_14_6}
P(s)=c(s-s_{1})(s-s_{2})\cdots(s-s_{n}),
\end{equation}
where $c$ and $s_{k}$ $(k=1,2,\ldots,n)$ are spatial complex constants. More precisely, the theorem
ensures that $P(s)$ has a zero $s_{1}$. Then
\[
P(s)=(s-s_{1})Q_{1}(s)
\]
where $Q_{1}(s)$ is a polynomial of degree $n-1$. The same argument, applied to $Q_{1}(s)$, reveals that there is a number $s_{2}$ such that
\[
P(s)=(s-s_{1})(s-s_{2})Q_{2}(s)
\]
where $Q_{2}(s)$ is a polynomial of degree $n-2$. Continuing in this way, we arrive at expression~(\ref{eq:chap_4_14_6}). Some of the constants $s_{k}$ in expression~(\ref{eq:chap_4_14_6}) may, of  course, appear more than once, and it is clear that $P(s)$ can have no more than $n$ distinct zeros.

\section{Maximum Modulus Principle}\label{sec:chap_4_15_50}

In this section, we derive an important result involving maximum values of the moduli of analytic functions. We begin with a needed lemma.

\begin{lemma}\label{le:chap_4_15_1}
Suppose that $|f(s)|\leq|f(s_{0})|$ at each point $s$ in some neighborhood $|s-s_{0}|<\varepsilon$ in which $f$ is analytic. Then $f(s)$ has the constant value $f(s_{0})$ throughout that neighborhood.
\end{lemma}

To prove this, we assume that $f$ satisfies the stated conditions and let $s_{1}$ be any point other than $s_{0}$ in the given neighborhood. We then let $\rho$ be the distance between $s_{1}$ and $s_{0}$. If $S_{\rho}$ denotes the sphere $|s-s_{0}|=\rho$ centered at $s_{0}$, and $C_{\rho}$ is a positively oriented circle on the sphere $S_{\rho}$, which passes through $s_{1}$ and is in the spatial plane determined by points $s_{0}$, $s_{1}$, and the origin (Subsec.~(\ref{sec:chap_1_9_5})), the Cauchy integral formula tells us that
\begin{equation}\label{eq:chap_4_15_1}
f(s_{0})=\frac{1}{2\pi i}\int_{C_{\rho}}\frac{f(s)ds}{s-s_{0}};
\end{equation}
and in the transformed coordinates system, the parametric representation
\[
s=s_{0}+\rho e^{i\theta}\texttt{ }(0\leq\theta\leq2\pi)
\]
for $C_{\rho}$ enables us to write equation~(\ref{eq:chap_4_15_1}) as
\begin{equation}\label{eq:chap_4_15_2}
f(s_{0})=\frac{1}{2\pi}\int_{0}^{2\pi}f(s_{0}+\rho e^{i\theta})d\theta.
\end{equation}
We note from expression~(\ref{eq:chap_4_15_2}) that when a function is analytic within and on a given
sphere, its value at the center is the arithmetic mean of its values on the circle. This result is called Gauss's mean value theorem.

From equation~(\ref{eq:chap_4_15_2}), we obtain the inequality
\begin{equation}\label{eq:chap_4_15_3}
|f(s_{0})|\leq\frac{1}{2\pi}\int_{0}^{2\pi}|f(s_{0}+\rho e^{i\theta})|d\theta.
\end{equation}
On the other hand, since
\begin{equation}\label{eq:chap_4_15_4}
|f(s_{0}+\rho e^{i\theta})|\leq|f(s_{0})|.
\end{equation}
we find that
\[
\int_{0}^{2\pi}|f(s_{0}+\rho e^{i\theta})|d\theta\leq\int_{0}^{2\pi}|f(s_{0})|d\theta=2\pi|f(s_{0})|.
\]
Thus
\begin{equation}\label{eq:chap_4_15_5}
|f(s_{0})|\geq\frac{1}{2\pi}\int_{0}^{2\pi}|f(s_{0}+\rho e^{i\theta})|d\theta.
\end{equation}
It is now evident from inequalities~(\ref{eq:chap_4_15_3}) and~(\ref{eq:chap_4_15_5}) that
\[
|f(s_{0})|=\frac{1}{2\pi}\int_{0}^{2\pi}|f(s_{0}+\rho e^{i\theta})|d\theta,
\]
or
\[
\int_{0}^{2\pi}(|f(s_{0})|-|f(s_{0}+\rho e^{i\theta})|)d\theta=0.
\]
The integrand in this last integral is continuous in the variable $\theta$; and, in view of condition~(\ref{eq:chap_4_15_4}), it is greater than or equal to zero on the entire interval $0\leq\theta\leq2\pi$. Because the value of the integral is zero, then, the integrand must be identically equal to zero. That is,
\begin{equation}\label{eq:chap_4_15_6}
|f(s_{0}+\rho e^{i\theta})|=|f(s_{0})|\texttt{ }0\leq\theta\leq2\pi.
\end{equation}
This shows that $|f(s)|=|f(s_{0})|$ for all points $s$ on the circle as well as the sphere $|s-s_{0}|=\rho$.

Finally, since $s_{1}$ is any point in the deleted neighborhood $0<|s-s_{0}|<\varepsilon$, we see that the equation $|f(s)|=|f(s_{0})|$ is, in fact, satisfied by all points $s$ lying on any sphere $|s-s_{0}|=\rho$, where $0<\rho<\varepsilon$. Consequently, $|f(s)|=|f(s_{0})|$ everywhere in the neighborhood $|s-s_{0}|<\varepsilon$. But we know that when the modulus of an analytic function is constant in a domain, the function itself is constant there. Thus $f(s)=f(s_{0})$ for each point $s$ in the neighborhood, and the proof of the lemma is complete.

This lemma can be used to prove the following theorem, which is known as the maximum modulus principle.

\begin{theorem}\label{th:chap_4_15_1}
If a function $f$ is analytic and not constant in a given domain $D$, then
$|f(s)|$ has no maximum value in $D$. That is, there is no point $s_{0}$ in the domain such
that $|f(s)|<|f(s_{0})|$ for all points $s$ in it.
\end{theorem}

Given that $f$ is analytic in $D$, we shall prove the theorem by assuming that $|f(s)|$ does have a maximum value at some point $s_{0}$ in $D$ and then showing that $f(s)$ must be constant throughout $D$.

The general approach here is similar to that taken in the proof of the lemma in Sec.~(\ref{sec:chap_2_16_26}).  We draw a polygonal line  $L$ lying in $D$ and extending from $s_{0}$ to any other point $P$ in $D$.  Also, $d$ represents the shortest distance from points on $L$ to the boundary of $D$. When $D$ is the entire plane, $d$ may have any positive value. Next, we observe that there is a finite sequence of points
\[
s_{0},s_{1},s_{2},\ldots,s_{n-1},s_{n}
\]
along $L$ such that $s_{n}$ coincides with the point $P$ and
\[
|s_{k}-s_{k-1}|<d\texttt{ }k=1,2,\ldots,n.
\]
On forming a finite sequence of neighborhoods
\[
N_{0},N_{1},N_{2},\ldots,N_{n-1},N_{n},
\]
where each $N_{k}$ has center $s_{k}$ and radius $d$, we see that $f$ is analytic in each of these neighborhoods, which are all contained in $D$, and that the center of each neighborhood
$N_{k}$ $(k=1,2,\ldots,n)$ lies in the neighborhood $N_{k-1}$.

Since $|f(s)|$ was assumed to have a maximum value in $D$ at $s_{0}$, it also has a maximum value in $N_{0}$ at that point. Hence, according to the preceding lemma, $f(s)$ has the constant value $f(s_{0})$ throughout $N_{0}$. In particular, $f(s_{1})=f(s_{0})$. This means that $|f(s)|\leq|f(s_{1})|$ for each point $s$ in $N_{1}$; and the lemma can be applied again, this time telling us that
\[
f(s)=f(s_{1})=f(s_{0})
\]
when $s$ is in $N_{1}$. Since $s_{2}$ is in $N_{1}$ then, $f(s_{2})=f(s_{0})$. Hence $|f(s)|\leq|f(s_{2})|$ when
$s$ is in $N_{2}$; and the lemma is once again applicable, showing that
\[
f(s)=f(s_{2})=f(s_{0})
\]
when $s$ is in $N_{2}$. Continuing in this manner, we eventually reach the neighborhood $N_{n}$ and arrive at the fact that $f(s_{n})=f(s_{0})$.

Recalling that $s_{n}$ coincides with the point $P$, which is any point other than $s_{0}$ in $D$, we may conclude that $f(s)=f(s_{0})$ for every point $s$ in $D$. Inasmuch as $f(s)$ has now been shown to be constant throughout  $D$, the theorem is proved.

If a function $f$ that is analytic at each point in the interior of a closed bounded region $R$ is also continuous throughout $R$, then the modulus $|f(s)|$ has a maximum value somewhere in $R$ (Sec.~(\ref{sec:chap_2_7_17})). That is, there exists a nonnegative constant $M$ such that $|f(s)|\leq M$ for all points $s$ in $R$, and equality holds for at least one such point. If $f$ is a constant function, then $|f(s)|=M$ for all $s$ in $R$. If, however, $f(s)$ is not constant, then, according to the maximum modulus principle,  $|f(s)|\neq M$ for any point $s$ in the interior of $R$. We thus arrive at an important corollary of the maximum modulus principle.

\begin{corollary}\label{co:chap_4_15_1}
Suppose that a function $f$ is continuous on a closed bounded region $R$ and that it is analytic and not constant in the interior of $R$. Then the maximum value of $|f(s)|$ in R, which is always reached, occurs somewhere on the boundary of $R$ and never in the interior.
\end{corollary}

When the function $f$ in the corollary is written
\[
f(s)=e_{xy}[u(x,y,z)+iv(x,y,z)]+w(z),
\]
the component function $u(x,y,z)$ also has a maximum value in $R$ which is assumed on the boundary of $R$ and never in the interior, where it is harmonic (Sec.~(\ref{sec:chap_2_15_25})).% For the composite function $g(s)=\exp(f(s))$ is continuous in $R$ and analytic and not constant in the interior. Consequently, its modulus $|g(s)|=\exp(u(x,y,z))$,  which is continuous in $R$, must assume its maximum value in $R$ on the boundary. Because of the increasing nature  of the exponential function, it follows that the maximum value of $u(x,y,z)$ also occurs on the boundary.

%-----------------------------------------------------------------------
% Beginning of chap5.tex
%-----------------------------------------------------------------------
%
% AMS-LaTeX 1.2 sample file for a monograph, based on amsbook.cls.
% This is a data file input by chapter.tex.
%%%%%%%%%%%%%%%%%%%%%%%%%%%%%%%%%%%%%%%%%%%%%%%%%%%%%%%%%%%%%%%%%%%%

%\part{This is a Part Title Sample}

\chapter{Surface Integrals}\label{ch:chap_5}

\section*{Introduction}

Integrals are extremely important in the study of a spatial complex variable. The theory of surface integration of spatial complex functions is to be developed in this chapter. A theorem of the surface integral of spatial complex functions will be put forward and proved. The surface integral is governed by this theorem of a spatial complex variable and may not be processed by the theory of the two-dimensional complex variables. This theorem may be called the Cauchy surface integral theorem. Another important theorem of the surface integral is put forward and proved, which is called Boundary theorem and expresses the relation of the integral of an open surface and the spatial contour integral of the boundary spatial contour of the open surface.

Representation formulas, and in particular integral representation formulas, play an important role in mathematics, since they allow us to discover a function on a large set from its behavior on a small set. Here, we will first prove a surface integral representation formula in a manner that it is independent of the theory of harmonic functions. Then we will prove other theorems of the surface integral of spatial complex functions.

These formula and theorems are similar to the Cauchy formula and theorems of the plane contour integral of the two-dimensional complex functions, such as the Mean value theorem, Cauchy inequalities, Liouville's theorem, Morera's theorem, maximum modulus principle and Schwarz's lemma. But they are of the surface integral which may not be processed by the theory of the two-dimensional complex variables.

\section{Cauchy-Theorem of Surface Integrals of the Second Kind}\label{sec:chap_5_1_51}

Theorems and formulas in next two subsections are referenced in this section.

\subsection{The Gauss-Ostrogradskii Formula}

Just as Green's formula connects the integral over the boundary of a plane domain with a corresponding integral over the domain itself, the Gauss-Ostrogradskii formula given below connects the integral over the boundary of a three-dimensional domain with an integral over the domain itself.

\begin{theorem}\label{th:chap_5_1_1}
Let $\mathbb{R}^{3}$ be three-dimensional space with a fixed coordinate system $x$,$y$,$z$ and $\overline{D}$ a compact domain in $\mathbb{R}^{3}$ bounded by piecewise-smooth surfaces $S$. Let $X$, $Y$, and $Z$ be smooth functions in the closed domain $\overline{D}$.

Then the following relation holds:
\begin{equation}\label{eq:chap_5_1_1}
\iint_{S}X(x,y,z)dydz+Y(x,y,z)dzdx+Z(x,y,z)dxdy
\end{equation}
\[
=\iiint_{\overline{D}}(\frac{\partial X}{\partial x}+\frac{\partial Y}{\partial y}+\frac{\partial Z}{\partial z})dxdydz
\texttt{ where $S$ is taken its outside.}
\]
\end{theorem}
Expression~(\ref{eq:chap_5_1_1}) is called Gauss-Ostrogradskii Formula.

\begin{theorem}\label{th:chap_5_1_2}
Suppose that $X(x,y,z)$, $Y(x,y,z)$, $Z(x,y,z)$, and their first-order partial derivatives are continuous on the compact domain $\overline{D}$. If any one of the following statements is true, then so are the others:

(i) the surface integrals
\[
\iint_{S}X(x,y,z)dydz+Y(x,y,z)dzdx+Z(x,y,z)dxdy
\]
has no relation with the form of the surface $S$, that is, it only has relations with the contour of $S$;

(ii) there exists a vector function
\[
B(x,y,z)=U(x,y,z)\mathbf{i}+V(x,y,z)\mathbf{j}+W(x,y,z)\mathbf{k}
\]
such that
\[
X=\frac{\partial W}{\partial y}-\frac{\partial V}{\partial z},\texttt{ }
Y=\frac{\partial U}{\partial z}-\frac{\partial W}{\partial x},\texttt{ }
Z=\frac{\partial V}{\partial x}-\frac{\partial U}{\partial y};
\]

(iii) $\forall (x,y,z)\in \overline{D}$ there is $\frac{\partial X}{\partial x}+\frac{\partial Y}{\partial y}+\frac{\partial Z}{\partial z}=0$;

(iv) for any $S$ in $\overline{D}$ there is
\[
\iiint_{\overline{D}}(\frac{\partial X}{\partial x}+\frac{\partial Y}{\partial y}+\frac{\partial Z}{\partial z})dxdydz=0
\texttt{ where $S$ is taken its outside.}
\]
\end{theorem}

\subsection{Stokes' Formula in $\mathbb{R}^{3}$}

\begin{theorem}\label{th:chap_5_1_3}
Let $S$ be an oriented piecewise-smooth compact two-dimensional surface with boundary $C$ embedded in a domain $D\subset\mathbb{R}^{3}$, in which a smooth 1-form $\omega=Udx+Vdy+Wdz$ is defined. Then the following relation holds:
\begin{equation}\label{eq:chap_5_1_2}
\oint_{C}Udx+Vdy+Wdz
\end{equation}
\[
=\iint_{S}(\frac{\partial W}{\partial y}-\frac{\partial V}{\partial z})dydz
+(\frac{\partial U}{\partial z}-\frac{\partial W}{\partial x})dzdx
+(\frac{\partial V}{\partial x}-\frac{\partial U}{\partial y})dxdy
\]
where the orientation of the boundary $C$ is chosen consistently with the orientation of the surface $S$.
\end{theorem}

\begin{theorem}\label{th:chap_5_1_4}
Suppose that $U(x,y,z)$, $V(x,y,z)$, $W(x,y,z)$, and their first-order partial derivatives are continuous on the compact domain $D$. If any one of the following statements is true, then so are the others:

(i) the curve integrals
\[
\oint_{C(A,B)}U(x,y,z)dx+V(x,y,z)dy+W(x,y,z)dz
\]
has no relation with the path of the curve $C(A,B)$, that is, it only has relations with the start point $A$ and end point $B$ of $C(A,B)$;

(ii) there exists a real function $\phi(x,y,z)$ and a vector function
\[
B(x,y,z)=U(x,y,z)\mathbf{i}+V(x,y,z)\mathbf{j}+W(x,y,z)\mathbf{k}
\]
with a vector $dr=dx\mathbf{i}+dy\mathbf{j}+dz\mathbf{k}$ such that
\[
X=\frac{\partial W}{\partial y}-\frac{\partial V}{\partial z},\texttt{ }
Y=\frac{\partial U}{\partial z}-\frac{\partial W}{\partial x},\texttt{ }
Z=\frac{\partial V}{\partial x}-\frac{\partial U}{\partial y},
\]
\[
d\phi=B(x,y,z)\cdot dr=U(x,y,z)dx+V(x,y,z)dy+W(x,y,z)dz;
\]

(iii) $\forall (x,y,z)\in D$ there are
\[
X=Y=Z=0\texttt{ or }\frac{\partial X}{\partial x}+\frac{\partial Y}{\partial y}+\frac{\partial Z}{\partial z}=0;
\]

(iv) for any oriented piecewise-smooth curve $\Gamma$ in $D$ there is
\[
\oint_{\Gamma}U(x,y,z)dx+V(x,y,z)dy+W(x,y,z)dz=0.
\]
\end{theorem}

\subsection{Cauchy-Theorem of Spatial Contour Integrals}

\begin{theorem}\label{th:chap_5_1_5}
Suppose that a spatial complex function
\[
f(s)=\varpi=(u(x,y,z),v(x,y,z),w(z))\texttt{ where }s=(x,y,z)
\]
is analytic in a simply connected closed complex space $D$, $S$ is any simple closed surface in the closed complex space $D$, and $C$ is any closed spatial contour on the surface $S$, then, there are
\begin{equation}\label{eq:chap_5_1_3}
\oint_{C}f(s)\cdot ds=\oint_{C}u(x,y,z)dx+v(x,y,z)dy+w(z)dz=0.
\end{equation}
\end{theorem}

\begin{proof}
First, from Theorem~(\ref{th:chap_5_1_3}) and Theorem~(\ref{th:chap_5_1_4}), there are
\[
\oint_{C}u(x,y,z)dx+v(x,y,z)dy+w(z)dz=\iint_{S}Xdydz+Ydzdx+Zdxdy
\]
where
\[
X=\frac{\partial w}{\partial y}-\frac{\partial v}{\partial z},\texttt{ }
Y=\frac{\partial u}{\partial z}-\frac{\partial w}{\partial x},\texttt{ }
Z=\frac{\partial v}{\partial x}-\frac{\partial u}{\partial y}.
\]

Next, from Theorem~(\ref{th:chap_5_1_1}) and Theorem~(\ref{th:chap_5_1_2}), there are
\[
\iint_{S}Xdydz+Ydzdx+Zdxdy
=\iiint_{D}(\frac{\partial X}{\partial x}+\frac{\partial Y}{\partial y}+\frac{\partial Z}{\partial z})dxdydz.
\]

The first-order partial derivatives of $X$, $Y$, and $Z$ are
\[
\frac{\partial X}{\partial x}=\frac{\partial^{2}w}{\partial y\partial x}-\frac{\partial^{2}v}{\partial z\partial x},\texttt{ }
\frac{\partial Y}{\partial y}=\frac{\partial^{2}u}{\partial z\partial y}-\frac{\partial^{2}w}{\partial x\partial y},\texttt{ }
\frac{\partial Z}{\partial z}=\frac{\partial^{2}v}{\partial x\partial z}-\frac{\partial^{2}u}{\partial y\partial z}.
\]

Now, by a theorem in advanced calculus, the continuity of the partial derivatives of $u$, $v$, and $w$ ensures that $u_{zy}=u_{yz}$, $v_{zx}=v_{xz}$, and $w_{yx}=w_{xy}$. It then follows from above expressions of the first-order partial derivatives of $X$, $Y$, and $Z$ that there are
\[
\frac{\partial X}{\partial x}+\frac{\partial Y}{\partial y}+\frac{\partial Z}{\partial z}
=\frac{\partial^{2}w}{\partial y\partial x}-\frac{\partial^{2}v}{\partial z\partial x}
+\frac{\partial^{2}u}{\partial z\partial y}-\frac{\partial^{2}w}{\partial x\partial y}
+\frac{\partial^{2}v}{\partial x\partial z}-\frac{\partial^{2}u}{\partial y\partial z}=0.
\]

Thus we get
\[
\oint_{C}u(x,y,z)dx+v(x,y,z)dy+w(z)dz
=\iiint_{D}(\frac{\partial X}{\partial x}+\frac{\partial Y}{\partial y}+\frac{\partial Z}{\partial z})dxdydz=0.
\]
The proof of Theorem~(\ref{th:chap_5_1_5}) is completed.
\end{proof}

\subsection{Cauchy-Theorem of Spatial Surface Integrals}

\begin{theorem}\label{th:chap_5_1_6}
Suppose that a spatial complex function
\[
f(s)=\varpi=(u(x,y,z),v(x,y,z),w(z))\texttt{ where }s=(x,y,z)
\]
is analytic in a simply connected closed complex space $D$, $S$ is any simple closed surface in the closed complex space $D$, then, there are
\begin{equation}\label{eq:chap_5_1_4}
\iint_{S}\mathbf{F}\cdot d\mathbf{S}=\iint_{S}X(x,y,z)dydz+Y(x,y,z)dzdx+Z(x,y,z)dxdy=0
\end{equation}
where
\[
\mathbf{F}(x,y,z)=X(x,y,z)\mathbf{i}+Y(x,y,z)\mathbf{j}+Z(x,y,z)\mathbf{k},\texttt{ }
d\mathbf{S}=dydz\mathbf{i}+dzdx\mathbf{j}+dxdy\mathbf{k},
\]
and
\[
X=\frac{\partial w}{\partial y}-\frac{\partial v}{\partial z},\texttt{ }
Y=\frac{\partial u}{\partial z}-\frac{\partial w}{\partial x},\texttt{ }
Z=\frac{\partial v}{\partial x}-\frac{\partial u}{\partial y}.
\]
\end{theorem}

\begin{proof}
From Theorem~(\ref{th:chap_5_1_1}) and Theorem~(\ref{th:chap_5_1_2}), there are
\[
\iint_{S}Xdydz+Ydzdx+Zdxdy
=\iiint_{D}(\frac{\partial X}{\partial x}+\frac{\partial Y}{\partial y}+\frac{\partial Z}{\partial z})dxdydz.
\]

The first-order partial derivatives of $X$, $Y$, and $Z$ are
\[
\frac{\partial X}{\partial x}=\frac{\partial^{2}w}{\partial y\partial x}-\frac{\partial^{2}v}{\partial z\partial x},\texttt{ }
\frac{\partial Y}{\partial y}=\frac{\partial^{2}u}{\partial z\partial y}-\frac{\partial^{2}w}{\partial x\partial y},\texttt{ }
\frac{\partial Z}{\partial z}=\frac{\partial^{2}v}{\partial x\partial z}-\frac{\partial^{2}u}{\partial y\partial z}.
\]

Now, by a theorem in advanced calculus, the continuity of the partial derivatives of $u$, $v$, and $w$ ensures that $u_{zy}=u_{yz}$, $v_{zx}=v_{xz}$, and $w_{yx}=w_{xy}$. It then follows from above expressions of the first-order partial derivatives of $X$, $Y$, and $Z$ that there are
\[
\frac{\partial X}{\partial x}+\frac{\partial Y}{\partial y}+\frac{\partial Z}{\partial z}
=\frac{\partial^{2}w}{\partial y\partial x}-\frac{\partial^{2}v}{\partial z\partial x}
+\frac{\partial^{2}u}{\partial z\partial y}-\frac{\partial^{2}w}{\partial x\partial y}
+\frac{\partial^{2}v}{\partial x\partial z}-\frac{\partial^{2}u}{\partial y\partial z}=0.
\]

Thus we get
\[
\iint_{S}X(x,y,z)dydz+Y(x,y,z)dzdx+Z(x,y,z)dxdy
=\iiint_{D}(\frac{\partial X}{\partial x}+\frac{\partial Y}{\partial y}+\frac{\partial Z}{\partial z})dxdydz=0.
\]
The proof of Theorem~(\ref{th:chap_5_1_6}) is completed.
\end{proof}

\section{Surface Integral of the First Kind}\label{sec:chap_5_2_52}

Assuming that in a complex space $D$, there is a simple surface $S$ connecting two points $s_{1}$ and $s_{2}$. The simple surface $S$ is smooth or smooth in separate pieces. Assume that $f(s)=e_{xy}[u(x,y,z)+iv(x,y,z)]+w(z)$ is a continuous spatial complex function on $S$, where $u(x,y,z)$, $v(x,y,z)$, and $w(z)$ are the components of $f(s)$. Let the surface $S$ is divided into $n$ smaller piecewise surfaces $\Delta S_{k}$ by divided points $s_{0},s_{1},\ldots,s_{n-1}$, where points $s_{k}(k=0,1,2,\ldots,n-1)$ are sorted from $s_{0}$ to $s_{n-1}$ in an order on the surface $S$. The area of a piecewise surfaces $\Delta S_{k}$ is approximately equal to $\Delta\sigma_{k}$. When $\varsigma_{k}$ is any point on $\Delta S_{k}$, then, the sum
\begin{equation}\label{eq:chap_5_2_1}
\sum^{n-1}_{k=0}f(\varsigma_{k})\Delta\sigma_{k}
\end{equation}
may be written as
\begin{equation}\label{eq:chap_5_2_2}
\sum^{n-1}_{k=0}f(\varsigma_{k})\Delta\sigma_{k}
=\sum^{n-1}_{k=0}\{e_{xy}[u(x_{k},y_{k},z_{k})+iv(x_{k},y_{k},z_{k})]+w(z_{k})\}\Delta\sigma_{k}
\end{equation}
\[
=e_{xy}[\sum^{n-1}_{k=0}u(x_{k},y_{k},z_{k})\Delta\sigma_{k} +i\sum^{n-1}_{k=0}v(x_{k},y_{k},z_{k})\Delta\sigma_{k}]
+\sum^{n-1}_{k=0}w(z_{k})\Delta\sigma_{k}.
\]

Hence, every sum on the right-hand side of equation~(\ref{eq:chap_5_2_2}) is the sum of real functions on the simple region of a spatial surface. According to the integral of the real functions, when the number $n$, or the sum of the divided points on the surface $S$ increases to infinite, and
\[
\max\{|\Delta\varsigma_{k}|,k=0,1,2,\ldots,n-1\}\rightarrow 0,
\]
all sums on the right-hand side of equation~(\ref{eq:chap_5_2_2}) are of
limits, separately,
\[
\lim_{n\rightarrow\infty}\sum^{n-1}_{k=0}u(x_{k},y_{k},z_{k})\Delta\sigma_{k}=\iint_{S}u(x,y,z)d\sigma
\]
\[
\lim_{n\rightarrow\infty}\sum^{n-1}_{k=0}v(x_{k},y_{k},z_{k})\Delta\sigma_{k}=\iint_{S}v(x,y,z)d\sigma,
\]
and
\[
\lim_{n\rightarrow\infty}\sum^{n-1}_{k=0}w(z_{k})\Delta\sigma_{k}=\iint_{S}w(z)d\sigma.
\]
Thus, the sum in equation~(\ref{eq:chap_5_2_1}) is of limit
\[
\lim_{n\rightarrow\infty}\sum^{n-1}_{k=0}f(\varsigma_{k})\Delta\sigma_{k}=\iint_{S}f(s)d\sigma.
\]
This is called the integral of the spatial complex function $f(s)$ on
the surface $S$, and denoted by $\iint_{S}f(s)d\sigma$. Then, there is
\begin{equation}\label{eq:chap_5_2_3}
\iint_{S}f(s)d\sigma=\iint_{S}\{e_{xy}[u(x,y,z)+iv(x,y,z)]+w(z)\}d\sigma.
\end{equation}
According to the definition of spatial complex numbers, we have
\[
\iint_{S}f(s)d\sigma=e_{xy}[U(x,y,z)+iV(x,y,z)]+W(z).
\]

Based on the definition of the surface integral of spatial complex functions, its basic properties may be derived, just as that in two-dimensional complex functions. Assuming that two spatial complex functions $f(s)$ and $g(s)$ are continuous on a simple smooth surface, then, there are

(1) $\iint_{S}\alpha f(s)d\sigma=\alpha\iint_{S}f(s)d\sigma$, where $\alpha$ is a spatial complex constant;

(2) $\iint_{S}[f(s)+g(s)]d\sigma=\iint_{S}f(s)d\sigma+\iint_{S}g(s)d\sigma$;

(3) $\iint_{S}f(s)d\sigma=\iint_{S_{1}}f(s)d\sigma+\iint_{S_{2}}f(s)d\sigma+\cdots+\iint_{S_{m}}f(s)d\sigma$,
where the surface $S$ is consist of some simple smooth surfaces $S_{1},S_{2},\ldots,S_{m}$.

(4) $\iint_{S^{-}}f(s)d\sigma=-\iint_{S}f(s)d\sigma$, where $S^{-}$ means that the integral is in an anti-direction on the surface $S$.

(5) When $|f(s)|\leq M$ on the surface $S$, there is $\iint_{S}f(s)d\sigma\leq M\Omega$, where $\Omega$ is the area of the surface $S$, $M$ and $\Omega$ are finite positive numbers.

The definition and properties of the surface integral of spatial complex functions are valid not only for a simple smooth surface but also for general smooth or smooth in piecewise surfaces.

\section{Cauchy-Theorem of Surface Integrals of the First Kind}\label{sec:chap_5_3_53}

\begin{theorem}[Cauchy Theorem of Surface Integrals]\label{th:chap_5_3_1}
Suppose that a spatial complex function $f(s)$ is analytic in a simply connected complex space $D$, and $S$ is any simple closed surface in $D$, then, there is
\begin{equation}\label{eq:chap_5_3_1}
\iint_{S}f(s)d\sigma=0.
\end{equation}
\end{theorem}

\begin{proof}
The proof is of three steps.

First step: $S$ is the bound of any tetrahedron in space $D$, Suppose
\[
|\iint_{S}f(s)d\sigma|=M,
\]
let prove that $M=0$.

Because every piecewise surface of a tetrahedron is a space triangle, let every edge of the tetrahedron is bisected, and every two neighborhood points of these bisected points are connected. In this way, the tetrahedron is divided into eight same smaller tetrahedrons, whose bounds are $S_{1},S_{2},\ldots,S_{8}$, respectively. Thus, there is
\begin{equation}\label{eq:chap_5_3_2}
\iint_{S}f(s)d\sigma=\iint_{S_{1}}f(s)d\sigma+\iint_{S_{2}}f(s)d\sigma+\cdots+\iint_{S_{8}}f(s)d\sigma
\end{equation}
On the bounds of these bisected smaller tetrahedron, which are not on the bounds of the original tetrahedron, the integral is just happen two times but in opposite direction, so, they are canceled out each other. According to equation~(\ref{eq:chap_5_3_2}), in the bound surfaces $S_{k}(k=1,2,\ldots,8)$, there should be at least one
surface on which the modulus of the surface integral is not less than $M/8$. For example, this surface is $S_{1}$:
\[
|\iint_{S_{1}}f(s)d\sigma|\geq \frac{M}{8}.
\]
Same as above, this tetrahedron is divided into eight smaller tetrahedrons by bisecting its edges. One of which is with bound surface $S^{(2)}$ and satisfies
\[
|\iint_{S^{(2)}}f(s)d\sigma|\geq \frac{M}{8^{2}}.
\]

Doing continuous in this way infinitely, a sequence of tetrahedrons with bound surfaces $S^{(1)}, S^{(2)},S^{(3)},\ldots,S^{(n)},\ldots$ may be obtained, in which, every one contains the next and there is
\begin{equation}\label{eq:chap_5_3_3}
|\iint_{S^{(n)}}f(s)d\sigma|\geq \frac{M}{8^{n}}(n=1,2,\ldots).
\end{equation}

Let $\Omega$ be the area of the surface $S$, so, the area of the surface $S^{(n)}$ is $\Omega/8^{n}(n=1,2,\ldots)$. Now, let estimate the modulus of $|\iint_{S^{(n)}}f(s)d\sigma|$. Because
every tetrahedron in the sequence contains all of its followings, and $\Omega/8^{n}\rightarrow 0(n\rightarrow\infty)$, hence, there should be at least one point $s_{0}$ belongs to all of these tetrahedrons in the sequence. Because that the function $f(s)$ is of its derivative $f'(s_{0})$ at point $s_{0}$, therefore, for an arbitrary $\varepsilon>0$, we may find a $\delta>0$, and it is true that when $s\in D$ and $0\leq|s-s_{0}|\leq\delta$, there is
\[
|\frac{f(s)-f(s_{0})}{s-s_{0}}-f'(s_{0})|<\varepsilon.
\]

Thus, when $s\in D$ and $0\leq|s-s_{0}|\leq\delta$, there is
\begin{equation}\label{eq:chap_5_3_4}
|f(s)-f(s_{0})-(s-s_{0})f'(s_{0})|<\varepsilon(s-s_{0}).
\end{equation}

It is obvious that when $n$ is sufficient large, $S^{(n)}$ is contained by a sphere determined by $|s-s_{0}|<\delta$. Hence, when $s\in S^{(n)}$, equation~(\ref{eq:chap_5_3_4}) is true, and
$|s-s_{0}|<\frac{\Omega}{8^{n}}$, thus,
\[
|f(s)-f(s_{0})-(s-s_{0})f'(s_{0})|<\frac{\Omega}{8^{n}}\varepsilon.
\]

Next, because $r\leq|s-s_{0}|\leq\delta$, $\iint_{S^{(n)}}f(s)d\sigma<4\pi\delta^{2}$, we have
\begin{equation}\label{eq:chap_5_3_5}
|\iint_{S^{(n)}}f(s)d\sigma|
=|\iint_{S^{(n)}}[f(s)-f(s_{0})-(s-s_{0})f'(s_{0})]d\sigma|<4\pi\delta^{2}\frac{\Omega}{8^{n}}\varepsilon
\end{equation}

By comparing equations~(\ref{eq:chap_5_3_3}) and~(\ref{eq:chap_5_3_5}), we obtain
\[
\frac{M}{8^{n}}<4\pi\delta^{2}\frac{\Omega}{8^{n}}\varepsilon
\]
that is,
\[
M<4\pi\varepsilon\delta^{2}\Omega.
\]
Because $\varepsilon$ is an arbitrary positive number, there is $M=0$.

Second step: $S$ is the bound surface $P$ of any simple space polyhedron.

Let the simple polyhedron with surface $P$ is divided into some simple tetrahedrons with its diagonal lines. Then, the integral on surface $P$ may be derived by the sum of integrals on the bound surfaces of these tetrahedrons. Same as the division of a tetrahedron, the integrals on the bound surfaces of these
tetrahedrons, which are not on the surface of the polyhedron, are done two times in opposite directions and canceled out each other. Hence, for the simple polyhedron, there is
\begin{equation}\label{eq:chap_5_3_6}
\iint_{P}f(s)d\sigma=0.
\end{equation}

Assuming that the surface $P$ is any closed polygonal planes, which may be intersect each other, according to the instance of the simple polyhedron stated above, equation~(\ref{eq:chap_5_3_6}) is also valid.

Third step: $S$ is any closed surface in space $D$.

(1) For an arbitrary positive number $\varepsilon > 0$, we may find a simple inscribed polyhedron to the surface $S$, and
\begin{equation}\label{eq:chap_5_3_7}
|\iint_{S}f(s)d\sigma-\iint_{P}f(s)d\sigma|<\varepsilon,
\end{equation}
that is, the integral $\iint_{S}f(s)d\sigma$ may be approximated to any accurate value by the integral on the bound surface of a simple inscribed polyhedron to the surface $S$ in space $D$.

To prove this fact, let consider a closed subspace $\overline{G}$ in space $D$ and let the entire surface $S$ is in $G$. Assume $\rho$ is the minimum distance between the bound of $G$ and the surface $S$. It is easy to know that $\rho>0$. Hence, a sphere with a center at any point on $S$ and a radius $\rho$, is entirely contained by
$\overline{G}$, and when a distance between any two points on $S$ is less than $\rho$, their connection line must be entirely in $\overline{G}$.

According to the assumption in the theorem, the function $f(s)$ is continuous on $\overline{G}$, thus, $f(s)$ is uniformly continuous on $\overline{G}$. So, for any small positive number
$\varepsilon>0$, there exists a positive number $\delta_{1}=\delta_{1}(\varepsilon)$, and when $s^{'}$ and $s^{''}$ are on $\overline{G}$ and satisfy $|s^{'}-s^{''}|<\delta_{1}$, the inequality
\[
|f(s^{'})-f(s^{''})|<\frac{\varepsilon}{2\sigma}
\]
is true, where $\sigma$ is the area of $S$.

It is obvious that $n$ points $(s_{1},s_{2},\ldots,s_{n})$ on $S$ may be taken according to the positive direction of the surface integral, to divide the surface $S$ into $n$ small piecewise
surfaces $(\sigma_{1},\sigma_{2},\ldots,\sigma_{n})$, and
\[
\max_{1\leq i\leq n}(\texttt{area of
}\sigma_{i})<\delta<\min(\delta_{1}, \rho).
\]

Hence, the simple polyhedron $P$ with points $s_{1},s_{2},\ldots,s_{n}$ as its vertexes is entirely contained by
$\overline{G}$(and entirely contained by $D$). Let $\tau_{k}$ denote the piecewise planes of $P$ where $k=1,2,\ldots,n$. Then the small piecewise planes $\tau_{1},\tau_{2},\ldots,\tau_{n}$ of $P$ are the inscribed small piecewise planes of the small piecewise surfaces $\sigma_{1},\sigma_{2},\ldots,\sigma_{n}$. Thus, there is
\[
\begin{split}
|\iint_{S}f(s)d\sigma-\iint_{P}f(s)d\sigma|
&=|\sum^{n}_{i=1}\iint_{\sigma_{i}}f(s)d\sigma-\sum^{n}_{i=1}\iint_{\tau_{i}}f(s)d\sigma| \\
&\leq\sum^{n}_{i=1}|\iint_{\sigma_{i}}f(s)d\sigma-\iint_{\tau_{i}}f(s)d\sigma|.
\end{split}
\]

Because of
\[
\iint_{\sigma_{i}}f(s_{i})d\sigma=f(s_{i})\iint_{\sigma_{i}}d\sigma=\iint_{\tau_{i}}f(s_{i})d\sigma
\]
there is
\[
|\iint_{\sigma_{i}}f(s)d\sigma-\iint_{\tau_{i}}f(s)d\sigma| \leq
|\iint_{\sigma_{i}}[f(s)-f(s_{i})]d\sigma|+|\iint_{\tau_{i}}[f(s)-f(s_{i})]d\sigma|
\]
\[
\leq\sup_{s\in\sigma_{i}}[f(s)-f(s_{i})](\texttt{area of
}\sigma_{i})+\sup_{s\in \tau_{i}}[f(s)-f(s_{i})](\texttt{area of
}\tau_{i})
\]

Because the distance between any two points on the small piecewise surfaces $\sigma_{i}$ and the small piecewise planes $\tau_{i}$ respectively is less than $\delta$, there is
\[
\sup_{s\in\sigma_{i}}[f(s)-f(s_{i})](\texttt{area of
}\sigma_{i})<\frac{\varepsilon}{2\sigma}
\]
\[
\sup_{s\in \tau_{i}}[f(s)-f(s_{i})](\texttt{area of
}\tau_{i})<\frac{\varepsilon}{2\sigma}.
\]
Hence we have
\[
|\iint_{\sigma_{i}}f(s)d\sigma-\iint_{\tau_{i}}f(s)d\sigma|
<\frac{\varepsilon}{2\sigma}(\texttt{area of }\sigma_{i}+\texttt{area of
}\tau_{i})<\frac{\varepsilon}{\sigma}(\texttt{area of }\sigma_{i}).
\]
Thus, there is
\[
|\iint_{S}f(s)d\sigma-\iint_{P}f(s)d\sigma|
<\frac{\varepsilon}{\sigma}(\texttt{area of }\sigma)=\varepsilon.
\]

(2) According to the results obtained in second step, for the simple polyhedron $P$ made in property (1) of last section, there is
\[
\iint_{P}f(s)d\sigma=0,
\]
and equation~(\ref{eq:chap_5_3_7}) becomes
\[
|\iint_{S}f(s)d\sigma|<\varepsilon.
\]
Because $\varepsilon$ is an arbitrary positive number, the must be
\[
|\iint_{S}f(s)d\sigma|=0.
\]
Hence, the Cauchy theorem of the surface integral of a spatial complex variable has been proved.
\end{proof}

\section{Extension of Cauchy-Theorem of Surface Integrals}\label{sec:chap_5_4_54}

Based on the above Cauchy-Theorem of the surface integral of a spatial complex variable, we may obtain extensive theorems of Cauchy-Theorem.

\begin{theorem}\label{th:chap_5_4_1}
Suppose that a spatial complex function $f(s)$ is analytic in a simply connected complex space $D$, and $S$ is any closed surface (which may be not simple) in $D$, then, there is
\begin{equation}\label{eq:chap_5_4_1}
\iint_{S}f(s)d\sigma=0.
\end{equation}
\end{theorem}

\begin{proof}
Because the surface $S$ may always be considered as one surface consisted of some finite surfaces in a space $D$ by connection, and according to the basic property (3) in section~(\ref{sec:chap_5_2_52}) of the surface integral of a spatial complex variable and the above Cauchy-Theorem, the theorem may be proved.
\end{proof}

\begin{theorem}\label{th:chap_5_4_2}
Suppose that a spatial complex function $f(s)$ is analytic in a simply connected complex space $D$, then, the surface integral of $f(s)$ in a space $D$ is independent of the shape of the integral surface $S$, and only dependent of the boundary contour. That is, for any closed spatial contour $C$ in $D$, the value of a surface integral
\begin{equation}\label{eq:chap_5_4_2}
\iint_{S}f(s)d\sigma=\frac{1}{2}\oint_{C}F(s)d\bar{s},
\end{equation}
is independent of the shape of the surface with a closed boundary spatial contour in space $D$, where $F(s)$ is the primitive function of $f(s)$ (See the proof in following section~(\ref{sec:chap_5_5_55})).
\end{theorem}

\begin{proof}
Suppose $S_{1}$ and $S_{2}$ are any two surfaces with a same closed boundary spatial contour in space $D$, then, surface $S_{1}$ in positive direction and surface $S_{1}$ in opposite direction are connected to form a
close surface $S$ in space $D$. Thus, according to theorem~(\ref{eq:chap_5_4_1}) and the basic property (3) in section~(\ref{sec:chap_5_2_52}) of the surface integral of a spatial complex variable, there is
\[
\iint_{S_{1}}f(s)d\sigma+\iint_{S^{-}_{2}}f(s)d\sigma=\iint_{S}f(s)d\sigma=0.
\]
Hence we have
\[
\iint_{S_{1}}f(s)d\sigma=\iint_{S_{2}}f(s)d\sigma.
\]
\end{proof}

Let prove that the following theorem is equivalent to the above Cauchy-Theorem.

\begin{theorem}\label{th:chap_5_4_3}
Suppose that surface $S$ is a closed surface, the space region $D$ is the internal of $S$, and a spatial complex function $f(s)$ is analytic in a closed space region $\overline{D}=D+S$, then
\[
\iint_{S}f(s)d\sigma=0.
\]
\end{theorem}

\begin{proof}
(1) From Theorem~(\ref{th:chap_5_3_1}) to Theorem~(\ref{th:chap_5_4_3}).

By the assumption of Theorem~(\ref{th:chap_5_3_1}), the spatial complex function $f(s)$ must be analytic in a simply connected complex space $G$ which contains $\overline{D}$. Hence, according to Theorem~(\ref{th:chap_5_3_1}), there should be $\iint_{S}f(s)d\sigma=0$.

(2) From Theorem~(\ref{th:chap_5_4_3}) to Theorem~(\ref{th:chap_5_3_1}).

By the assumption of Theorem~(\ref{th:chap_5_3_1}), the spatial complex function $f(s)$ is analytic in a simply connected complex space $D$ and $S$ is any closed surface in $D$. Assume that the space $G$ is in the internal of $S$, then, $f(s)$ must be analytic on the closed space region $\overline{D}=D+S$. Hence, according to
Theorem~(\ref{th:chap_5_4_3}), there should be $\iint_{S}f(s)d\sigma=0$.
\end{proof}

Then, let extend the Cauchy-Theorem of surface integration of spatial complex functions on the other hand, that is, extend the Cauchy-Theorem from a bounded simply connected space region with the bound of a simple closed surface to bounded multi-connected space regions with the bound of a complex surface consisted of multi-closed surfaces.

\begin{definition}\label{def:chap_5_4_1}
Considering $n+1$ simple closed surfaces $S_{0},S_{1},S_{2},\ldots,S_{n}$, where every closed surface of $S_{1},S_{2},\ldots,S_{n}$ is at the outsides of others, and all of them are in the internal of the surface $S_{0}$. These point sets in the internal of the surface $S_{0}$ and at the outsides of $S_{1},S_{2},\ldots,S_{n}$ form a bounded multi-connected space region $D$ which is of the bound consisted of simple closed surfaces $S_{0},S_{1}, S_{2},\ldots,S_{n}$. In this case, we define that the bound of
space region $D$ is a complex surface $S=S_{0}+S_{1}^{-}+S_{2}^{-}+\cdots+S_{n}^{-}$, which consists of
$S_{0}$ in positive direction and $S_{1},S_{2},\ldots,S_{n}$ in opposite direction.
\end{definition}

\begin{theorem}\label{th:chap_5_4_4}
Suppose that a bounded and multi-connected space region $D$ is bounded by a complex surface
$S=S_{0}+S_{1}^{-}+S_{2}^{-}+\cdots+S_{n}^{-}$, The spatial complex function $f(s)$ is analytic in $D$ and continuous in the closed space region $\overline{D}=D+S$, then
\[
\iint_{S}f(s)d\sigma=0,
\]
or rewritten as
\begin{equation}\label{eq:chap_5_4_3}
\iint_{S_{0}}f(s)d\sigma+\iint_{S_{1}^{-}}f(s)d\sigma+\iint_{S_{2}^{-}}f(s)d\sigma+\cdots+\iint_{S_{n}^{-}}f(s)d\sigma=0,
\end{equation}
or rewritten as
\begin{equation}\label{eq:chap_5_4_4}
\iint_{S_{0}}f(s)d\sigma=\iint_{S_{1}}f(s)d\sigma+\iint_{S_{2}}f(s)d\sigma+\cdots+\iint_{S_{n}}f(s)d\sigma.
\end{equation}
\end{theorem}

\begin{proof}
Take $n+1$ auxiliary simple smooth surfaces $L_{0},L_{1},L_{2},\ldots,L_{n}$, which are mutually disjoint each other and all in $D$ (except their bounds), as cutting surfaces. Connect them to $S_{0},S_{1},S_{2},\ldots,S_{n}$ in order. Assuming that the space region $D$ is cut along the cutting surfaces, then, $D$ is divided
into two simply connected space regions, each of which is bounded by a closed surface and denoted by $\Gamma_{1}$ and $\Gamma_{2}$. It is obvious that there are
\[
\iint_{\Gamma_{1}}f(s)d\sigma=0,\iint_{\Gamma_{2}}f(s)d\sigma=0.
\]
Adding the two equations and considering that the surface integral on the surfaces $L_{0},L_{1},L_{2},\ldots,L_{n}$ are done two times in opposite direction and canceled out each other, then, according
to the basic property (3) in section~(\ref{sec:chap_5_2_52}) of the surface integral of spatial complex
functions, we may obtain
\[
\iint_{S}f(s)d\sigma=0.
\]
Hence, there are equations~(\ref{eq:chap_5_4_3}) and~(\ref{eq:chap_5_4_4}).
\end{proof}

\section{Boundary Theorem of Surface Integrals}\label{sec:chap_5_5_55}

The Theorem~(\ref{th:chap_5_4_2}) points out that the surface integral of a spatial complex analytic function $f(s)$ in a simply connected space region $D$ is independent of the shape of the integral surface $S$ and only dependent of the bounded contour of the surface $S$. The following theorem determines the value of the surface integral or the bounded contour integral of the surface $S$.

\begin{theorem}[Boundary Integrals Theorem]\label{th:chap_5_5_1}
Suppose that a spatial complex function $f(s)$ is analytic in a simply connected complex space $D$, then, the surface integral of $f(s)$ in a space $D$ is independent of the shape of the integral surface $S$, and only dependent of the boundary contour $C$. The value of the surface integral of $f(s)$ may be expressed as
\begin{equation}\label{eq:chap_5_5_1}
\iint_{S}f(s)d\sigma=\frac{1}{2}\oint_{C}F(s)d\bar{s}-\int_{C}F(s)dz,
\end{equation}
where $F(s)$ is the primitive function of $f(s)$.
\end{theorem}

To prove the theorem, let give out definitions:
\begin{definition}\label{def:chap_5_5_1}
For a plane region $G$ bounded by a close contour $C$, assume that $x_{a}$ and $x_{b}$ are respectively the minimum and maximum values of the projections of contour $C$ on the $x$ axis. For every point between $x_{a}$ and $x_{b}$, make a line parallel to $y$ axis is intersected with the contour $C$, If the line is continuous in $G$, then $G$ may be called a \emph{continues region} to $x$, that is, any line parallel to $y$ axis and intersected with the contour $C$ is continuous.
\end{definition}

\begin{definition}\label{def:chap_5_5_2}
For a surface region $G$ bounded by a close contour C, assume that $x_{a}$ and $x_{b}$ are respectively the minimum and maximum values of the projections of contour $C$ on the $x$ axis. For every point between $x_{a}$ and $x_{b}$, make a cutting plane perpendicular to the $x$ axis and form a cutting contour $C_{p}$ on the surface $G$. If $C_{p}$ is continuous, then $G$ may be called a \emph{continuous surface} region to $x$, that is, any contour on the surface $G$ and on a plane perpendicular to the $x$ axis is continuous.
\end{definition}

\begin{proof}
(1) First, let consider that the boundary contour $C$ is on a plane, that is, all points of the contour $C$ are on a plane. Without loss of generality we assume that the plane is the $xy$ coordinate plane, otherwise, we may choose suitable coordinate system or by rotation mapping to map the plane to the $xy$ coordinate plane.

Because the surface integral of $f(s)$ in a space $D$ is independent of the shape of the integral surface $S$, hence when the integral surface $S$ is the plane region $G$ bounded by the contour $C$ on the $xy$ coordinate plane, the surface integral must be with the same value.

When $G$ is not a \emph{continues region} to $x$, we may always divide it into $m+1$ continues regions by using $m$ secant lines parallel to $x$ axis. When $G$ is a \emph{continues region} to $x$, the number $m$ may be considered as $m=0$. Hence we may write $G=G_{0}+G_{1}+G_{2}+\cdots+G_{m}$.

Now, we consider $G_{k}(k=0,1,2,\ldots,m)$. Let $S_{k}$ be the surface, which is a plane, of $G_{k}$ and $C_{k}$ be the bounded contour of $G_{k}$. Assume that $x_{a}$ and $x_{b}$ are respectively the minimum and maximum values of the projections of contour $C_{k}$ on $x$ axis. Corresponding to $x_{a}$ and $x_{b}$, there are $s_{a}$ and $s_{b}$. The contour $C_{k}$ is divided into $C_{k}=C_{k,1}+C_{k,2}$ by the two points.

Then, take $n$ points: $x_{a}=x_{1}<x_{2}<...<x_{n}=x_{b}$. For point $x_{i}$ along $x$ axis, let make a line $L_{i}$ parallel to $y$ axis. The line $L_{i}$ is intersected with the contour $C_{k}$. Because $G_{k}$ is a continues region, Line $L_{i}$ is continuous and there are two intersecting points $s^{'}$ and $s^{''}$ for line $L_{i}$ and contour $C_{k}$.

Hence for a small plane region $L_{i}\Delta x$, the integral $\Delta x\int_{L_{i}}f(s)ds$ is the surface integral on the small surface region. Because the integral is entirely on a complex plane, it is a two-dimensional complex integral on a complex plane and the spatial complex variable $s$ becomes a two-dimensional complex variable, according to the contour integral of the two-dimensional complex functions, we have
\[
\Delta x\int_{L_{i}}f(s)ds=[F(s^{'})-F(s^{''})]\Delta x.
\]
where $F(s)$ is the primitive function of $f(s)$.

Hence considering the sum for $x_{a}=x_{1}<x_{2}<...<x_{n}=x_{b}$, we have
\[
\sum^{n}_{i=1}\int_{L_{i}}f(s)ds\Delta x
=\sum^{n}_{i=1}F(s^{'})\Delta x -\sum^{n}_{i=1}F(s^{''})\Delta x.
\]

When $n\rightarrow\infty$ and $\Delta x\rightarrow0$, the sum above becomes following integrals:
\[
\iint_{S_{k}}f(s)d\sigma=\int_{C_{k,2}}F(s)dx-\int_{C^{-}_{k,1}}F(s)dx=\int_{C_{k}}F(s)dx.
\]
where $C_{k,1}$ and $C_{k,2}$ are the bound contours of the integrals.

For $G=G_{0}+G_{1}+G_{2}+\cdots+G_{m}$, we have
\[
\iint_{S}f(s)d\sigma=\sum^{m}_{k=0}\iint_{S_{k}}f(s)d\sigma=\sum^{m}_{k=0}\int_{C_{k}}F(s)dx.
\]
Because the integral on the $m$ secant lines is done two times in opposite direction and canceled out each other, we have
\begin{equation}\label{eq:chap_5_5_2}
\iint_{S}f(s)d\sigma=\sum^{m}_{k=0}\int_{C_{k}}F(s)dx=\int_{C}F(s)dx.
\end{equation}

According to the definition of spatial complex numbers (Sec.~(\ref{sec:chap_1_5})), there are
\[
x=\frac{Im_{s}s+Im_{s}\bar{s}}{2}=\frac{s+\bar{s}}{2}-z.
\]
Substituting the above expression into equation~(\ref{eq:chap_5_5_2}), we have
\begin{equation}\label{eq:chap_5_5_3}
\iint_{S}f(s)d\sigma=\frac{1}{2}\int_{C}F(s)ds+\frac{1}{2}\int_{C}F(s)d\bar{s}-\int_{C}F(s)dz.
\end{equation}
According to the definition of spatial complex functions and the Cauchy-Theorem of surface integration of a spatial complex variable, the first integral at the right-hand side of equation~(\ref{eq:chap_5_5_3}) is equal to zero. Hence the equation~(\ref{eq:chap_5_5_1}) is true.

(2) Next, let consider that the boundary contour $C$ is a spatial contour, that is, there is at least a point and all other points of the contour $C$ are not on a same plane. In this case, let take the projection of contour $C$ on the $xy$ coordinate plane. We may choose suitable coordinate system or by rotation mapping to make the projection be
the largest plane region.

Let point $s_{0}$ on $C$ be connected to all other points on $C$ to form a surface $G$. Same as the analysis above, the surface integral must be with the same value.

When $G$ is not a \emph{continues surface} to $x$, we may always divide it into $m+1$ continues surfaces by using $m$ cutting planes parallel to the $x$ axis and form $m$ cutting contours. When $G$ is a \emph{continues surface} to $x$, the number $m$ may be considered as $m=0$. Hence we may write $G=G_{0}+G_{1}+G_{2}+\cdots+G_{m}$.

Now, we consider $G_{k}(k=0,1,2,\ldots,m)$. Let $S_{k}$ be the surface of $G_{k}$ and $C_{k}$ be the bounded contour of $G_{k}$. Assume that $x_{a}$ and $x_{b}$ are respectively the minimum and maximum values of the projections of contour $C_{k}$ on the $x$ axis. Corresponding to $x_{a}$ and $x_{b}$, there are $s_{a}$ and $s_{b}$. The contour $C_{k}$ is divided into $C_{k}=C_{k,1}+C_{k,2}$ by the two points.

Then, take $n$ points: $x_{a}=x_{1}<x_{2}<...<x_{n}=x_{b}$. For point $x_{i}$ along $x$ axis, let make a space line $L_{i}$ parallel to $y$ axis. The line $L_{i}$ is intersected with the contour $C_{k}$. Because $G_{k}$ is a continues surface, Line $L_{i}$ is continuous and there are two intersecting points $s^{'}$ and $s^{''}$ for line $L_{i}$ and contour $C_{k}$.

Assume $P_{i}$ is the plane parallel to $y$ axis and passing line $L_{i}$. Let a contour $\Gamma_{i}$ is the intersecting contour of the plane $P_{i}$ and the contour $C_{k}$. The contour $\Gamma_{i}$ and line $L_{i}$ form a closed contour. By choose suitable local coordinate or rotation mapping, the plane $P_{i}$ may be considered
as a two-dimensional complex plane and the integral of a complex function on line $L_{i}$ is equal to the integral of the function on contour $\Gamma_{i}$ according to the theory of the two-dimensional complex functions.

Hence for a small plane region $L_{i}\Delta x$, the integral $\Delta x\int_{L_{i}}f(s)ds$ is equal to the surface integral on the small surface region. Because the integral is entirely on a complex plane, it is a two-dimensional complex integral on a complex plane and the spatial complex variable $s$ becomes a two-dimensional complex variable, according to the contour integral of the two-dimensional complex functions, we have
\[
\Delta x\int_{L_{i}}f(s)ds=[F(s^{'})-F(s^{''})]\Delta x.
\]
where F(s) is the primitive function of $f(s)$.

Hence considering the sum for $x_{a}=x_{1}<x_{2}<...<x_{n}=x_{b}$, we have
\[
\sum^{n}_{i=1}\int_{L_{i}}f(s)ds\Delta x
=\sum^{n}_{i=1}F(s^{'})\Delta x -\sum^{n}_{i=1}F(s^{''})\Delta x.
\]

When $n\rightarrow\infty$ and $\Delta x\rightarrow0$, the sum above becomes following integrals:
\[
\iint_{S_{k}}f(s)d\sigma=\int_{C_{k,2}}F(s)dx-\int_{C^{-}_{k,1}}F(s)dx=\int_{C_{k}}F(s)dx.
\]
where $C_{k,1}$ and $C_{k,2}$ are the bound contours of the integrals.

For $G=G_{0}+G_{1}+G_{2}+\cdots+G_{m}$, we have
\[
\iint_{S}f(s)d\sigma=\sum^{m}_{k=0}\iint_{S_{k}}f(s)d\sigma=\sum^{m}_{k=0}\int_{C_{k}}F(s)dx.
\]
Because the integral on the $m$ cutting contours is done two times in opposite direction and canceled out each other, we have
\begin{equation}\label{eq:chap_5_5_4}
\iint_{S}f(s)d\sigma=\sum^{m}_{k=0}\int_{C_{k}}F(s)dx=\int_{C}F(s)dx.
\end{equation}

According to the definition of spatial complex numbers (Sec.~(\ref{sec:chap_1_5})), there are
\[
x=\frac{Im_{s}s+Im_{s}\bar{s}}{2}=\frac{s+\bar{s}}{2}-z.
\]
Substituting the above expression into equation~(\ref{eq:chap_5_5_4}), we have
\begin{equation}\label{eq:chap_5_5_5}
\iint_{S}f(s)d\sigma=\frac{1}{2}\int_{C}F(s)ds+\frac{1}{2}\int_{C}F(s)d\bar{s}-\int_{C}F(s)dz.
\end{equation}
According to the definition of spatial complex functions and the Cauchy-Theorem of surface integration of a spatial complex variable, the first integral at the right-hand side of equation~(\ref{eq:chap_5_5_5}) is equal to zero. Hence the equation~(\ref{eq:chap_5_5_1}) is true.

The proof of Theorem~(\ref{th:chap_5_5_1}) is completed.
\end{proof}

\section{Cauchy's Surface Integral Formula}\label{sec:chap_5_6_56}

\subsection{Cauchy's formula}

Suppose the spatial complex function $f(s)$ is analytic in a closed sphere $D$ bounded by a closed sphere
$D_{0}:|s-s_{0}|\leq\rho_{0}(0<\rho_{0}<\infty)$. By Cauchy's surface integral Theorem~(\ref{th:chap_5_3_1}) of spatial complex functions, the surface integral of $f(s)$ on any simple closed surface $S$ in $D_{0}$ is
equal to zero. Consider the surface integral
\begin{equation}\label{eq:chap_5_6_1}
I=\iint_{S}\frac{f(s)}{(s-s_{0})^{2}}d\sigma.
\end{equation}
Because the integral function $\frac{f(s)}{(s-s_{0})^{2}}$ is continuous on $S$, the integral $I$ should exist. But $\frac{f(s)}{(s-s_{0})^{2}}$ is not analytic in the above closed sphere $D_{0}$, the value of the integral $I$ is not certainly equal to zero.

Now, consider that $f(s)$ is a general analytic function. Let make a spherical surface $S_{\rho}$ centered at point $s_{0}$ and with a radius $\rho(0<\rho<\rho_{0})$. By Cauchy's surface integral theorem of spatial complex functions, we have
\begin{equation}\label{eq:chap_5_6_2}
\iint_{S}\frac{f(s)}{(s-s_{0})^{2}}d\sigma=\iint_{S_{\rho}}\frac{f(s)}{(s-s_{0})^{2}}d\sigma.
\end{equation}
The equation~(\ref{eq:chap_5_6_2}) is true for any $\rho$ which satisfies $0<\rho<\rho_{0}$. It is shown that the value of the integral $I$ is only dependent of the values of $f(s)$ in the spherical neighborhood of point $s_{0}$.

Let $C_{\varphi}$ be a circle $e_{xy}\rho e^{i\theta}\cos\varphi$ in the $e_{xy}$ plane and on the spherical surface $S_{\rho}$ where $0\leq\theta\leq2\pi$, $-\pi/2\leq\varphi\leq\pi/2$, and $\varphi$ is the included angle of $s-s_{0}$ and the $xy$ coordinate plane.

Then the radius $r_{\varphi}$ of $C_{\varphi}$ is equal to $\rho\cos\varphi$ or the projection $|s-s_{0}|\cos\varphi$ of the vector $s-s_{0}$ on the $e_{xy}$ plane, and by the definition~(\ref{th:chap_5_3_1}) of the surface integral,
\[
\begin{split}
d\sigma &= (s-s_{0})\cos\varphi d\theta(s-s_{0})d\varphi \\
        &= (s-s_{0})^{2}\cos\varphi d\varphi d\theta.
\end{split}
\]

Hence
\begin{equation}%\label{eq:chap_5_6_2b}
\begin{split}
I=\iint_{S}\frac{f(s)}{(s-s_{0})^{2}}d\sigma &= \iint_{S}f[s_{0}+\rho(e_{xy}e^{i\theta}\cos\varphi+\sin\varphi)]\cos\varphi d\varphi d\theta \\
&=\int_{-\pi/2}^{\pi/2}\cos\varphi d\varphi\int_{0}^{2\pi}f[s_{0}+\rho(e_{xy}e^{i\theta}\cos\varphi+\sin\varphi)]d\theta.
\end{split}
\end{equation}
By equation~(\ref{eq:chap_5_6_2}), the value of the integral $I$ is independent of $\rho$. According to the continuity of $f(s)$ at point $s_{0}$, we may conjecture that $I=4\pi f(s_{0})$, that is,
\begin{equation}\label{eq:chap_5_6_3}
f(s_{0})=\frac{1}{4\pi}\iint_{S}\frac{f(s)}{(s-s_{0})^{2}}d\sigma.
\end{equation}
The equation~(\ref{eq:chap_5_6_3}) is true, and we will prove it as following:

\begin{proof}
When $\rho$ approximates to zero for the integral at the right-hand side of equation~(\ref{eq:chap_5_6_2}), we have
\begin{equation}\label{eq:chap_5_6_4}
\iint_{S}\frac{f(s)}{(s-s_{0})^{2}}d\sigma=f(s_{0})\iint_{S_{\rho}}\frac{1}{(s-s_{0})^{2}}d\sigma+\iint_{S_{\rho}}\frac{f(s)-f(s_{0})}{(s-s_{0})^{2}}d\sigma.
\end{equation}
Because $f(s)$ is continuous at $s_{0}$, for an arbitrary $\varepsilon>0$, we may find a $\delta>0(\delta\leq\rho_{0})$ to make $|f(s)-f(s_{0})|<\varepsilon$ when $0<\rho<\delta,s\in S_{\rho}$. Then
\[
|\iint_{S_{\rho}}\frac{f(s)-f(s_{0})}{(s-s_{0})^{2}}d\sigma|\leq\frac{\varepsilon}{\rho^{2}}4\pi\rho^{2}=4\pi\varepsilon.
\]
Hence, when $\rho$ approximates zero, the second integral at the right-hand side of equation~(\ref{eq:chap_5_6_4}) approximates to zero. By $d\sigma=(s-s_{0})^{2}\cos\varphi d\varphi d\theta$, as stated above, the first
integral at the right-hand side of equation~(\ref{eq:chap_5_6_4}) is
\[
\begin{split}
\iint_{S_{\rho}}\frac{1}{(s-s_{0})^{2}}d\sigma&=\int_{-\pi/2}^{\pi/2}\cos\varphi d\varphi\int_{0}^{2\pi}d\theta \\
&=\int_{-\pi/2}^{\pi/2}d\sin\varphi\int_{0}^{2\pi}d\theta=4\pi.
\end{split}
\]
By this equation and equation~(\ref{eq:chap_5_6_2}), we have equation~(\ref{eq:chap_5_6_3}).
\end{proof}

\begin{theorem}[Cauchy's formula]\label{th:chap_5_6_1}
Suppose $f(s)$ is analytic in an open set that contains the closure of a sphere $D$, If $S$ denotes the boundary surface of this sphere with the positive orientation, then
\begin{equation}\label{eq:chap_5_6_5}
f(s)=\frac{1}{4\pi}\iint_{S}\frac{f(\varsigma)}{(\varsigma-s)^{2}}d\sigma
 \texttt{ for any point s $\in$ D}.
\end{equation}
\end{theorem}

\begin{proof} Suppose $s \in D$. It is obvious that the spatial complex function $f(\varsigma)/(\varsigma-s)^{2}$ is analytic at the point $\varsigma$, which satisfy $s \in \overline{D}, \varsigma\neq s$. Let make a closed sphere centered at point $s$ and contained by $D$. Let $\rho$ and $S_{\rho}$ denote the radius and boundary spherical surface of the sphere, respectively.

When the sphere bounded by $S_{\rho}$ is removed from $D$, the set of remain points on $\overline{D}$ is a closed space region $\overline{D}_{\rho}$. The spatial complex functions $f(\varsigma)$ and $f(\varsigma)/(\varsigma-s)^{2}$ are analytic on $\overline{D}_{\rho}$. By Cauchy's surface integral Theorem~(\ref{th:chap_5_3_1}) of spatial complex functions, we have
\begin{equation}\label{eq:chap_5_6_6}
\iint_{S}\frac{f(\varsigma)}{(\varsigma-s)^{2}}d\sigma=\iint_{S_{\rho}}\frac{f(\varsigma)}{(\varsigma-s)^{2}}d\sigma
\end{equation}
where the surface integral on $S$ is along the positive direction and the surface integral on $S_{\rho}$ is along the opposite direction.

By equations~(\ref{eq:chap_5_6_6}) and~(\ref{eq:chap_5_6_3}), we have equation~(\ref{eq:chap_5_6_5}), that is, the theorem of Cauchy's formula is true.
\end{proof}

\subsection{Mean value theorem}

As a special case, we have the following \textbf{Mean value theorem} of the spatial complex analytic functions.

\begin{theorem}[Mean value theorem]\label{th:chap_5_6_2}
If a spatial complex function $f(s)$ is analytic in a sphere $|\varsigma-s_{0}|<R$ and continuous on a closed sphere $|\varsigma-s_{0}|\leq R$, then
\begin{equation}\label{eq:chap_5_6_7}
f(s_{0})=\frac{1}{4\pi}\iint_{S}f(s_{0}+R(e_{xy}e^{i\theta}\cos\varphi+\sin\varphi))\cos\varphi d\varphi d\theta,
\end{equation}
\end{theorem}
that is, the value of $f(s)$ at the spherical center $s_{0}$ is equal to the mean value of $f(s)$ on the spherical surface $S$.

\begin{proof}
Let $S$ denote a spherical surface $|\varsigma-s_{0}|=R$. Then
\[
\varsigma-s_{0}=R(e_{xy}e^{i\theta}\cos\varphi+\sin\varphi)\texttt{ or }\varsigma=s_{0}+R(e_{xy}e^{i\theta}\cos\varphi+\sin\varphi).
\]
where $0\leq\theta\leq2\pi$ and $-\pi/2\leq\varphi\leq\pi/2$. By the definition of the surface integral above and the Cauchy formula, we have
\[
\begin{split}
f(s_{0}) &=\frac{1}{4\pi}\iint_{S}\frac{f(\varsigma)}{(\varsigma-s_{0})^{2}}d\sigma \\
         &=\frac{1}{4\pi}\iint_{S}f(s_{0}+R(e_{xy}e^{i\theta}\cos\varphi+\sin\varphi))\cos\varphi d\varphi d\theta
\end{split}
\]
Hence Theorem~(\ref{th:chap_5_6_2}) is true.
\end{proof}

\section{Derivatives of Space Complex Analytic Functions}\label{sec:chap_5_7_57}

As a corollary to the Cauchy surface integral formula of spatial complex functions, we will obtain further integral formulas expressing the derivatives of the spatial complex analytic function $f(s)$ inside the sphere in terms of the values of $f(s)$ on the boundary surface.

\subsection{Derivatives of spatial complex analytic functions}

Suppose that a spatial complex function $f(s)$ is analytic everywhere inside and on a simple closed spherical surface $S$, taken in the positive sense. If $s$ is any point interior to $S$, then expressions
\begin{equation}\label{eq:chap_5_7_0}
f'(s)=\frac{1}{4\pi}\iint_{S}\frac{2!f(\varsigma)}{(\varsigma-s)^{3}}d\sigma\texttt{, }
f''(s)=\frac{1}{4\pi}\iint_{S}\frac{3!f(\varsigma)}{(\varsigma-s)^{4}}d\sigma
\end{equation}
may be obtained formally, or without rigorous verification, by differentiating with respect to $s$ under the integral sign in the Cauchy integral formula
\[
f(s)=\frac{1}{4\pi}\iint_{S}\frac{f(\varsigma)}{(\varsigma-s)^{2}}d\sigma
\]
where $s$ is interior to $S$ and $\varsigma$ denotes points on $S$.

We will prove that these formulas are true.

\begin{theorem}\label{th:chap_5_7_1}
If a spatial complex function $f(s)$ is analytic in an open set $\Omega$, then $f(s)$ has infinitely many spatial complex derivatives in $\Omega$. Moreover, if $D\subset\Omega$ is a closed sphere whose interior is also contained in $\Omega$, $S\subset\Omega$ is the boundary spherical surface of $D$, then
\begin{equation}\label{eq:chap_5_7_1}
f^{(n)}(s)=\frac{1}{4\pi}\iint_{S}\frac{(n+1)!f(\varsigma)}{(\varsigma-s)^{n+2}}d\sigma\texttt{ }(n=1,2,\ldots)\texttt{ for any point s $\in$ D}.
\end{equation}
\end{theorem}

\begin{proof}
First, we prove that equation~(\ref{eq:chap_5_7_1}) is true for $n=1$. By Theorem~(\ref{th:chap_5_6_1}) and for $\Delta s\neq 0$, we have
\[
\begin{split}
\frac{f(s+\Delta s)-f(s)}{\Delta s} &= \frac{1}{\Delta s}\frac{1}{4\pi}[\iint_{S}\frac{f(\varsigma)}{(\varsigma-s-\Delta
s)^{2}}d\sigma-\iint_{S}\frac{f(\varsigma)}{(\varsigma-s)^{2}}d\sigma] \\
                                    &= \frac{1}{4\pi}\iint_{S}\frac{[2(\varsigma-s)-\Delta
s]f(\varsigma)}{(\varsigma-s-\Delta s)^{2}(\varsigma-s)^{2}}d\sigma.
\end{split}
\]

We will prove that the difference
\begin{equation}\label{eq:chap_5_7_2}
|\frac{1}{4\pi}[\iint_{S}\frac{[2(\varsigma-s)-\Delta
s]f(\varsigma)}{(\varsigma-s-\Delta
s)^{2}(\varsigma-s)^{2}}d\sigma-\iint_{S}\frac{2f(\varsigma)}{(\varsigma-s)^{3}}d\sigma]|
\end{equation}
\[
=|\frac{1}{4\pi}\iint_{S}\frac{[3(\varsigma-s)-2\Delta s]\Delta
sf(\varsigma)}{(\varsigma-s-\Delta s)^{2}(\varsigma-s)^{3}}d\sigma|
\]
is not over an arbitrary positive $\varepsilon$ for a sufficiently small $|\Delta s|$.

Suppose $|f(\varsigma)|\leq M$ on $S$. Let $d$ denote the smallest distance from $s$ to points $\varsigma$ on $S$. Then, there is $|\varsigma-s|\geq d >0$ for $\varsigma\in S$.

Suppose $|\Delta s|<d/2$, then
\[
|\varsigma-s-\Delta s|\geq|\varsigma-s|-|\Delta s|>d/2\texttt{ and }|3(\varsigma-s)-2\Delta s|<3|\varsigma-s|.
\]
Hence the difference by equation~(\ref{eq:chap_5_7_2}) is not over
\[
\frac{1}{4\pi}\frac{3M\sigma}{(d/2)^{2}d^{2}}|\Delta s|
\]
where $\sigma$ is the area of $S$. To make the value of the above expression not over the positive number $\varepsilon$, we only need to take
\[
|\Delta s|<\delta=min(\frac{d}{2}, \frac{\pi d^{4}
\varepsilon}{3M\sigma}).
\]
Hence, equation~(\ref{eq:chap_5_7_1}) is proved.

Next, it needs only the use of the mathematical induction to complete the proof. Suppose that when $n=m$, equation~(\ref{eq:chap_5_7_1}) is true. Then we will prove that when $n=m+1$, equation~(\ref{eq:chap_5_7_1}) is also true. That is to prove that the following expression
\[
\frac{f^{(m)}(s+\Delta s)-f^{(m)}(s)}{\Delta s}=\frac{1}{\Delta
s}\frac{1}{4\pi}[\iint_{S}\frac{(m+1)!f(\varsigma)}{(\varsigma-s-\Delta
s)^{m+2}}d\sigma-\iint_{S}\frac{(m+1)!f(\varsigma)}{(\varsigma-s)^{m+2}}d\sigma]
\]
approximates to
\[
\frac{1}{4\pi}\iint_{S}\frac{(m+2)!f(\varsigma)}{(\varsigma-s)^{m+3}}d\sigma, s\in D,
\]
when $\Delta s\rightarrow 0$. This proof is same as above.
\end{proof}

By the above theorem, we may derive the infinitely derivatives of spatial complex analytic functions.

\begin{theorem}\label{th:chap_5_7_2}
Suppose that a spatial complex function $f(s)$ is analytic in a complex space region $D$, then, $f(s)$ is of infinitely derivatives, which are also analytic in $D$.
\end{theorem}

\begin{proof}
Suppose $s_{0}$ is an arbitrary point in $D$, Let Theorem~(\ref{th:chap_5_7_1}) be applied to a sufficiently small sphere centered at $s_{0}$ (only needed the closed sphere entirely contained by $D$). It is known that $f(s)$ is of infinitely derivatives in the sphere. Specifically, $f(s)$ is of infinitely derivatives at point $S_{0}$. Because $s_{0}$ is arbitrary, $f(s)$ is of infinitely derivatives in $D$.
\end{proof}

In this way, because a spatial complex function is analytic in a complex space region $D$, we induce that its infinitely derivatives exist and are continuous.

\subsection{Cauchy inequalities and Liouville's theorem}

\begin{theorem}[Cauchy inequalities]\label{th:chap_5_7_3}
Suppose that a spatial complex function $f(s)$ is analytic inside and on a positively oriented sphere $D_{R}$, centered at $s_{0}$ and with a radius $R$. If $M_{R}$ denotes the maximum value of $|f(s)|$ on $D_{R}$, then
\begin{equation}\label{eq:chap_5_7_3}
|f^{(n)}(s_{0})|\leq\frac{(n+1)!M_{R}}{R^{n}}.
\end{equation}
\end{theorem}

\begin{proof}
Applying the Cauchy integral formula for $f^{(n)}(s_{0})$, we obtain
\[
\begin{split}
|f^{(n)}(s_{0})| &= |\frac{1}{4\pi}\iint_{S}\frac{(n+1)!f(\varsigma)}{(\varsigma-s_{0})^{n+2}}d\sigma| \\
&=\frac{(n+1)!}{4\pi}|\iint_{S}\frac{f(s_{0}+R(e_{xy}e^{i\theta}\cos\varphi+\sin\varphi))}{(R(e_{xy}e^{i\theta}\cos\varphi+\sin\varphi))^{n}}\cos\varphi d\varphi d\theta| \\
&=\frac{(n+1)!}{4\pi}|\int_{-\pi/2}^{\pi/2}\cos\varphi
d\varphi\int^{2\pi}_{0}\frac{M_{R}}{R^{n}}d\theta| \\
&\leq\frac{(n+1)!}{4\pi}\frac{M_{R}}{R^{n}}4\pi=\frac{(n+1)!M_{R}}{R^{n}}.
\end{split}
\]

The proof of Theorem~(\ref{th:chap_5_7_3}) is completed.
\end{proof}

\begin{theorem}[Liouville's theorem]\label{th:chap_5_7_4}
If a spatial complex function $f(s)$ is entire and bounded in the complex space, then $f(s)$ is constant throughout the complex space.
\end{theorem}

\begin{proof}
To start the proof, we assume that $f(s)$ is as stated in the theorem and note that, since $f(s)$ is entire, Cauchy's inequality with $n=1$ is holds for any choices of $s_{0}$ and $R$:
\begin{equation}\label{eq:chap_5_7_4}
|f'(s_{0})|\leq 2\frac{M_{R}}{R}.
\end{equation}
Moreover, the boundedness condition in the statement of the theorem tells us that a nonnegative constant $M$ exists such that $|f(s)|<2M$ for all $s$; and because the constant $M_{R}$ in inequality~(\ref{eq:chap_5_7_4}) is always less than or equal to $M$, it follows that
\begin{equation}\label{eq:chap_5_7_5}
|f'(s_{0})|\leq 2\frac{M}{R},
\end{equation}
where $s_{0}$ is any fixed point in the complex space and $R$ is arbitrarily large. Now the number $M$ in inequality~(\ref{eq:chap_5_7_5}) is independent of the value $R$ that is taken. Hence that inequality may hold for arbitrarily large values of $R$ only if $f'(s_{0})=0$. Since the choice of $s_{0}$ was arbitrary, this means that $f'(s)=0$ everywhere in the complex space. Consequently, $f(s)$ is a constant function.

The proof of Theorem~(\ref{th:chap_5_7_4}) is completed.
\end{proof}

\section{Morera's Theorem and Maximum Modulus Principle}\label{sec:chap_5_8_58}

\subsection{Morera's theorem}

\begin{theorem}[Morera's theorem]\label{th:chap_5_8_1}
If a spatial complex function $f(s)$ is continuous inside a complex space region $D$, and for any simple closed surface in $D$, there is
\begin{equation}\label{eq:chap_5_8_1}
\iint_{S}f(s)d\sigma=0.
\end{equation}
\end{theorem}
Then $f(s)$ is analytic inside $D$.

\begin{proof}
Take a point $s_{0}$ inside $D$ and make a spherical surface $S_{\rho}:|s-s_{0}|=\rho$ inside $D$. Then let rewrite equation~(\ref{eq:chap_5_8_1}) as following
\begin{equation}\label{eq:chap_5_8_2}
\iint_{S_{\rho}}\frac{s(s-2s_{0})f(s)}{(s-s_{0})^{2}}d\sigma+s_{0}^{2}\iint_{S_{\rho}}\frac{f(s)}{(s-s_{0})^{2}}d\sigma=0.
\end{equation}

First, consider the second integral at the left-hand side of equation~(\ref{eq:chap_5_8_2}), it may be rewritten as
\begin{equation}\label{eq:chap_5_8_3}
\iint_{S_{\rho}}\frac{f(s)}{(s-s_{0})^{2}}d\sigma=f(s_{0})\iint_{S_{\rho}}\frac{1}{(s-s_{0})^{2}}d\sigma+\iint_{S_{\rho}}\frac{f(s)-f(s_{0})}{(s-s_{0})^{2}}d\sigma.
\end{equation}

Because $f(s)$ is continuous at $s_{0}$, for an arbitrary $\varepsilon>0$, we may find a $\delta>0(\delta\leq\rho_{0})$, so that when $0<\rho<\delta, s\in S_{\rho}$, there is $|f(s)-f(s_{0})|<\varepsilon$. Then
\[
|\iint_{S_{\rho}}\frac{f(s)-f(s_{0})}{(s-s_{0})^{2}}d\sigma|\leq\frac{\varepsilon}{\rho^{2}}4\pi\rho^{2}=4\pi\varepsilon.
\]
Hence when $\rho$ approximates to zero, the second integral at the right-hand side of equation~(\ref{eq:chap_5_8_3}) approximates to zero. By $d\sigma=(s-s_{0})^{2}\cos\varphi d\varphi d\theta$ and as stated above, the first integral at the right-hand side of equation~(\ref{eq:chap_5_8_3}) is
\[
\begin{split}
\iint_{S_{\rho}}\frac{1}{(s-s_{0})^{2}}d\sigma &= \iint_{S_{\rho}}\cos\varphi d\varphi d\theta \\
&=\int_{-\pi/2}^{\pi/2}\cos\varphi d\varphi\int_{0}^{2\pi}d\theta=4\pi
\end{split}
\]
By this equation and equations~(\ref{eq:chap_5_8_2}) and~(\ref{eq:chap_5_8_3}), replacing $s$ by $\varsigma$ and replacing $s_{0}$ by $s$, we have
\begin{equation}\label{eq:chap_5_8_4}
f(s)=\frac{1}{4\pi}\frac{1}{s^{2}}\iint_{S_{\rho}}\frac{\varsigma(2s-\varsigma)f(\varsigma)}{(\varsigma-s)^{2}}d\sigma,
\end{equation}

Next, at any spherical neighborhood $K$ inside $D$, whose radius is less than $\rho/2$, take a closed spatial curve $C$ without passing through the origin $(|s|\neq0)$ and make an integral of $f(s)$ on $C$
\begin{equation}\label{eq:chap_5_8_5}
\oint_{C}f(s)ds=\frac{1}{4\pi}\oint_{C}\frac{1}{s^{2}}\iint_{S_{\rho}}\frac{\varsigma(2s-\varsigma)f(\varsigma)}{(\varsigma-s)^{2}}d\sigma
ds.
\end{equation}

Then exchange the sequence of integrals at the right-hand side of equation~(\ref{eq:chap_5_8_5}), we have
\begin{equation}\label{eq:chap_5_8_6}
\oint_{C}f(s)ds=\frac{1}{4\pi}\iint_{S_{\rho}}f(\varsigma)\varsigma
d\sigma\oint_{C}\frac{2s-\varsigma}{(\varsigma-s)^{2}s^{2}}ds.
\end{equation}

Because of $|s|\neq0$ and that the distance between any two points on $C$ is less than $\rho$, that is, $|s-\varsigma|\neq0$. It is obvious that the integral function on the spatial curve $C$ at the right-hand side of equation~(\ref{eq:chap_5_8_6}) is analytic inside and on the closed sphere $K$. Hence, by the Cauchy theorem of the contour integral of spatial complex functions, the value of the spatial contour integral is zero, that is, the integral at the right-hand side of equation~(\ref{eq:chap_5_8_6}) is zero. Thus we have
\begin{equation}\label{eq:chap_5_8_7}
\oint_{C}f(s)ds=0.
\end{equation}

Now by equation~(\ref{eq:chap_5_8_7}) and the theorem of primitive functions of the spatial contour integral of spatial complex functions, we know that the primitive function
\[
F(s)=\int^{s}_{s_{0}}f(\varsigma)d\varsigma (s_{0}\in D)
\]
is analytic inside $D$, and $F'(s)=f(s) (s_{0}\in D)$. Because the derivative function $F'(s)$ of the analytic function $F(s)$ is also an analytic function, it means that $f(s)$ is analytic inside $D$.

Thus the proof of the theorem is complete.
\end{proof}

Now we will give out the equivalence theorems of analytic functions defined in terms of surface integration.

\begin{theorem}\label{th:chap_5_8_2}
The necessary and sufficient condition of a spatial complex function $f(s)$ to be analytic inside a complex space region $D$ is

(1) $f(s)$ is continuous in $D$;

(2) for any closed surface $S$, if $S$ and its interior are entirely contained by $D$, there is
\[
\iint_{S}f(s)d\sigma=0.
\]
\end{theorem}

\begin{proof}
The necessity may be induced by Cauchy integral theorem. The sufficiency may be obtained by applying the Morera's theorem~(\ref{eq:chap_5_8_1}) at a spherical neighborhood $K:|\varsigma-s_{0}|<\rho$ of any point $s_{0}$ in $D$. If $\rho$ is sufficient small, we know that $f(s)$ is analytic in sphere $K$, specifically, $f(s)$ is analytic at point $s_{0}$. Because $s_{0}$ is arbitrarily chosen in $D$, $f(s)$ is analytic inside $D$.

Thus the proof of the theorem is complete.
\end{proof}

\subsection{Maximum modulus principle}

In this section, we derive an important result involving maximum values of the moduli of spatial complex analytic functions.

\begin{lemma}\label{le:chap_5_8_1}
Suppose that a spatial complex function $|f(s)|\leq |f(s_{0})$ at each point $s$ in some space neighborhood $|s-s_{0}|<\varepsilon$ in which $f(s)$ is analytic. Then $f(s)$ has the constant value $f(s_{0})$ throughout that space neighborhood.
\end{lemma}

\begin{proof}
To prove this, we assume that $f(s)$ satisfies the stated conditions and let $s_{1}$ be any point other than $s_{0}$ in the given space neighborhood. We then let $\rho$ be the distance between $s_{1}$ and $s_{0}$. If $S_{\rho}$ denotes the positively oriented spherical surface $|s-s_{0}|=\rho$, centered at $s_{0}$ and passing through $s_{1}$, the Cauchy surface integral formula tells us that
\begin{equation}\label{eq:chap_5_8_8}
f(s_{0})=\frac{1}{4\pi}\iint_{S_{\rho}}\frac{f(s)}{(s-s_{0})^{2}}d\sigma;
\end{equation}
and the parametric representation
\[
s=s_{0}+\rho(e_{xy}e^{i\theta}\cos\varphi+\sin\varphi)
\texttt{ }(0\leq\theta\leq2\pi,-\pi/2\leq\varphi\leq\pi/2)
\]
for $S_{\rho}$ enable us to write equation~(\ref{eq:chap_5_8_8}) as
\begin{equation}\label{eq:chap_5_8_9}
f(s_{0})=\frac{1}{4\pi}\iint_{S_{\rho}}f[s_{0}+\rho(e_{xy}e^{i\theta}\cos\varphi+\sin\varphi)]\cos\varphi d\varphi d\theta.
\end{equation}
We note from expression~(\ref{eq:chap_5_8_9}) that when a function is analytic within and on a given spherical surface, its value at the center is the arithmetic mean of its values on the spherical surface. This result is called Gauss's mean value theorem.

From equation~(\ref{eq:chap_5_8_9}), we obtain the inequality
\begin{equation}\label{eq:chap_5_8_10}
|f(s_{0})|\leq\frac{1}{4\pi}\iint_{S_{\rho}}|f[s_{0}+\rho(e_{xy}e^{i\theta}\cos\varphi+\sin\varphi)]|\cos\varphi d\varphi d\theta.
\end{equation}
On the other hand, since
\begin{equation}\label{eq:chap_5_8_11}
|f[s_{0}+\rho(e_{xy}e^{i\theta}\cos\varphi+\sin\varphi)]|\leq|f(s_{0})|\texttt{ }(0\leq\theta\leq2\pi,-\pi/2\leq\varphi\leq\pi/2)
\end{equation}
we find that
\[
\iint_{S_{\rho}}|f[s_{0}+\rho(e_{xy}e^{i\theta}\cos\varphi+\sin\varphi)]|\cos\varphi d\varphi d\theta\leq\iint_{S_{\rho}}|f(s_{0})|\cos\varphi d\varphi d\theta
\]
\[
=|f(s_{0})|\int_{-\pi/2}^{\pi/2}\varphi d\varphi\int_{0}^{2\pi}d\theta=4\pi|f(s_{0})|.
\]
Thus
\begin{equation}\label{eq:chap_5_8_12}
|f(s_{0})|\geq\frac{1}{4\pi}\iint_{S_{\rho}}|f[s_{0}+\rho(e_{xy}e^{i\theta}\cos\varphi+\sin\varphi)]|\cos\varphi d\varphi d\theta.
\end{equation}
It is now evident from inequalities~(\ref{eq:chap_5_8_10}) and
(\ref{eq:chap_5_8_12}) that
\[
|f(s_{0})|=\frac{1}{4\pi}\iint_{S_{\rho}}|f[s_{0}+\rho(e_{xy}e^{i\theta}\cos\varphi+\sin\varphi)]|\cos\varphi d\varphi d\theta.
\]
or
\[
\iint_{S_{\rho}}[|f(s_{0})|-|f[s_{0}+\rho(e_{xy}e^{i\theta}\cos\varphi+\sin\varphi)]|]\cos\varphi d\varphi d\theta=0.
\]
The integrand in this last integral is continuous in the variable $\theta$; and, in view of condition~(\ref{eq:chap_5_8_11}), it is greater than or equal to zero on the entire interval $0\leq\theta\leq2\pi,-\pi/2\leq\varphi\leq\pi/2$. Because the value of the integral is zero, then, the integrand must be identically equal to zero. That is,
\begin{equation}\label{eq:chap_5_8_13}
|f[s_{0}+\rho(e_{xy}e^{i\theta}\cos\varphi+\sin\varphi)]|=|f(s_{0})|\texttt{ }(0\leq\theta\leq2\pi,-\pi/2\leq\varphi\leq\pi/2)
\end{equation}
This shows that $|f(s)|=|f(s_{0})|$ for all points $s$ on the spherical surface $|s-s_{0}|=\rho$.

Finally, since $s_{1}$ is any point in the deleted spherical neighborhood $0<|s-s_{0}|<\varepsilon$, we see that the equation $|f(s)|=|f(s_{0})|$ is, in fact, satisfied by all points $s$ lying on any spherical surface $|s-s_{0}|=\rho$, where $0<\rho<\varepsilon$. Consequently, $|f(s)|=|f(s_{0})|$ everywhere in the spherical neighborhood $|s-s_{0}|<\varepsilon$. But we know that when the modulus of an analytic function is constant in a domain, the functions itself is constant there. Thus, $|f(s)|=|f(s_{0})|$ for each point $s$ in a spherical neighborhood, and the proof of the lemma is complete.
\end{proof}

This lemma may be used to prove the following theorem, which is called the \emph{maximum modulus principle of spatial complex functions}.

\begin{theorem}[Maximum modulus principle]\label{th:chap_5_8_3}
If a spatial complex function $f(s)$ is analytic and not constant in a given domain $D$, then $f(s)$ has no maximum value in D. That is, there is no point $s_{0}$ in the domain such that $|f(s)|\leq|f(s_{0})|$ for all points s in it.
\end{theorem}

\begin{proof}
Given that $f(s)$ is analytic in $D$, we shall prove the theorem by assuming that $|f(s)|$ does not have a maximum value at some point $s_{0}$ in $D$ and then showing that $f(s)$ must be constant throughout $D$.

The general approach here is similar to that taken in the proof of the lemma~(\ref{le:chap_5_8_1}). We draw a space polygonal line $L$ in $D$ and extending from $s_{0}$ to any other point $P$ in $D$.

Also, $d$ represents the shortest distance from points on $L$ to the boundary of $D$. When $D$ is the entire space, $d$ may have any positive value. Next, we observe that there is a finite sequence of points
\[
s_{0}, s_{1}, s_{2},\ldots,s_{n-1}, s_{n}
\]
along $L$ such that $s_{n}$ coincides with the point $P$ and
\[
|s_{k}-s_{k-1}|<d (k=1,2,\ldots,n).
\]
On forming a finite sequence of spherical neighborhoods
\[
S_{0}, S_{1}, S_{2},\ldots,S_{n-1}, S_{n},
\]
where each $S_{k}$ has center $s_{k}$ and a radius $d$, we see that $f(s)$ is analytic in each of the spherical neighborhood, which are all contained in $D$, and that the center of each spherical neighborhood $S_{k}(k=1,2,\ldots,n)$ lies in the spherical neighborhood $S_{k-1}$.

Since $|f(s)|$ was assumed to have a maximum value in $D$ at $s_{0}$, it also has a maximum value in $S_{0}$ at that point. Hence, according to the preceding lemma, $f(s)$ has the constant value $f(s_{0})$ throughout $S_{0}$. In particular, $f(s_{1})=f(s_{0})$. This means that $|f(s)|\leq|f(s_{1})|$ for each point $s$ in $S_{1}$; and the lemma may be applied again, this time telling us that
\[
f(s)=f(s_{1})=f(s_{0})
\]
when $s$ is in $S_{1}$. Since $s_{2}$ is in $S_{1}$, then, $f(s_{2})=f(s_{0})$. Hence, $|f(s)|\leq|f(s_{2})|$ when $s$ is in $S_{2}$; and the lemma is once again applicable, showing that
\[
f(s)=f(s{_2})=f(s_{0})
\]
when $s$ is in $S_{2}$. Continuing in this manner, we eventually reach the spherical neighborhood $S_{n}$ and arrive at the fact that $f(s_{n})=f(s_{0})$.

Recalling that $s_{n}$ coincides with the point $P$, which is any point other than $s_{0}$ in $D$, we may conclude that reach the spherical neighborhood $S_{n}$ and arrive at the fact that $f(s)=f(s_{0})$ for every point $s$ in $D$. Inasmuch as $f(s)$ has now been shown to be constant throughout $D$, the theorem is proved.
\end{proof}

If a spatial complex function $f(s)$ that is analytic at each point in the interior of a closed bounded space region $R$ is also continuous throughout $R$, then the modulus $|f(s)|$ has a maximum value somewhere in $R$. That is, there exists a nonnegative constant $M$ such that $|f(s)|\leq M$ for all points $s$ in $R$, and equality holds for at least one such point. If $f(s)$ is a constant function, then $|f(s)|=M$ for all $s$ in $R$. If, however, $f(s)$ is not constant, then, according to the maximum modulus principle, $|f(s)|\neq M$ for any point $s$ in the interior $R$. We thus arrive at an important corollary of the maximum modulus principle.

\textbf{Corollary}: \emph{Suppose that a spatial complex function $f(s)$ is continuous on a closed bounded space region $R$ and that it is analytic and not constant in the interior of $R$. Then the maximum value of $|f(s)|$ in $R$, which is always reached, occurs somewhere on the boundary of $R$ and never in the interior. }

\begin{lemma}[Schwarz]
If a spatial complex function $f(s)$ is analytic in a unit sphere $|s|<1$, and satisfies the conditions: $f(0)=0$ and $|f(s)|\leq 1$, then, in the sphere $|s|<1$, the must be
\begin{equation}\label{eq:chap_5_8_14}
|f(s)|\leq|s|,
\end{equation}
and
\begin{equation}\label{eq:chap_5_8_15}
|f'(0)|\leq 1.
\end{equation}
If the equality in equation~(\ref{eq:chap_5_8_15}) is true, or the equality in equation~(\ref{eq:chap_5_8_14}) is true at a point $s_{0}\neq0$ in the sphere $|s|<1$, then
\[
f(s)=(e_{xy}e^{i\alpha}\cos\beta+\sin\beta)s,
\]
where $\alpha$ and $\beta$ are real constants.
\end{lemma}

\begin{proof}
Let
\[
g(s)=\{
\begin{array}{c}
f(s)/s, \texttt{ for } s\neq 0, \\
f'(0), \texttt{ for } s= 0,
\end{array}
\]
hence $g(s)$ is analytic in a deleted sphere $0<|s|<1$. Because
\[
\lim_{s\rightarrow0}g(s)=\lim_{s\rightarrow0}\frac{f(s)-f(0)}{s-0}=f'(0)=g(0),
\]
$g(s)$ is continuous in sphere $|s|<1$. Now we will prove that $g(s)$ is analytic in sphere $|s|<1$. We need only to prove that the integral of $f(s)$ along the bound $\Gamma$ of any tetrahedron $T=[a,b,c,d]$ in sphere $|s|<1$ is zero. It is obvious that when point $s=0$ is not in $\Gamma$, because $g(s)$ is analytic in $0<|s|<1$, there is $\iint_{\Gamma}f(s)d\sigma=0$ by Cauchy theorem; when point $s$ is on $\Gamma$, because $g(s)$ is analytic in $\Gamma$ and is continuous to the bound $\Gamma$, there is $\iint_{\Gamma}f(s)d\sigma=0$ by the extensive theorem of Cauchy theorem; when point $s$ is in $\Gamma$, by making some auxiliary surface, we may convert this case to other cases above, and obtain the result $\iint_{\Gamma}f(s)d\sigma=0$. Hence by Morera's theorem, $g(s)$ is analytic in sphere $|s|<1$. On the spherical surface $|s|=r<1$,
\[
|g(s)|=|\frac{f(s)}{s}|\leq\frac{1}{r}.
\]
By the Maximum modulus principle, the equality in above expression is true in $|s|\leq r$, let $r\rightarrow1$, we obtain $|g(s)|\leq 1$. Thus we have equations~(\ref{eq:chap_5_8_14}) and~(\ref{eq:chap_5_8_15}).

When there is a point $s_{0} (s_{0}=0)$ to make $|f(s_{0})|=|s_{0}|$, that is, $|g(s_{0})|=1$, or to make
$|f'(0)|=0$, that is, $|g(0)|=1$, then, by the maximum modulus principle, there is $|g(s)|\equiv1$. Hence $g(s)$ is a constant. Let $g(s)=e_{xy}e^{i\alpha}\cos\beta+\sin\beta$, where $\alpha$ and $\beta$ are real constants, that is, $f(s)=(e_{xy}e^{i\alpha}\cos\beta+\sin\beta)s$. Hence the theorem is proved.
\end{proof}

The Schwarz lemma shows that for a spatial complex function $f(s)$, which is analytic in a unit sphere, if there is $f(0)=0$ and $|f(s)|\leq1$, then the distance from the mapping of point $s$ to the origin of the mapped coordinate is shorter than the distance from $s$ to the origin of the original coordinate. If there is a point $s$ to make the two distances equal to each other, then $f(s)$ is a rotational mapping.

%-----------------------------------------------------------------------
% Beginning of chap6.tex
%-----------------------------------------------------------------------
%
% AMS-LaTeX 1.2 sample file for a monograph, based on amsbook.cls.
% This is a data file input by chapter.tex.
%%%%%%%%%%%%%%%%%%%%%%%%%%%%%%%%%%%%%%%%%%%%%%%%%%%%%%%%%%%%%%%%%%%%

%\part{This is a Part Title Sample}

\chapter{Series}\label{ch:chap_6}

This chapter is devoted mainly to series representations of analytic spatial complex functions. We present theorems that guarantee the existence of such representations, and we develop some facility in manipulating series.

\section{Convergence of Sequences}\label{sec:chap_6_1_51}

An infinite sequence
\begin{equation}\label{eq:chap_6_1_1}
s_{1},s_{2},\ldots,s_{n},\ldots
\end{equation}
of spatial complex numbers has a limit $s$ if, for each positive number $\varepsilon$, there exists a positive integer $n_{0}$ such that
\begin{equation}\label{eq:chap_6_1_2}
|s_{n}-s|<\varepsilon\texttt{ whenever }n>n_{0}.
\end{equation}
Geometrically, this means that for sufficiently large values of $n$, the points $s_{n}$ lie in any given $\varepsilon$ neighborhood of $s$. Since we can choose $\varepsilon$ as small as we please, it follows that the points $s_{n}$ become arbitrarily close to $s$ as their subscripts increase. Note that the value $n_{0}$ that is needed will, in general, depend on the value of $\varepsilon$.

A sequence can have at most one limit. That is, a limit is unique if it exists. When the limit $s$ exists, the sequence is said to converge to $s$, and we write
\begin{equation}\label{eq:chap_6_1_3}
\lim_{n \to \infty}s_{n}=s.
\end{equation}
If a sequence has no limit, it diverges.

\begin{theorem}\label{th:chap_6_1_1}
Suppose that $s_{n}=e_{xy}(x_{n}+iy_{n})+z_{n}$ ($n=1,2,\ldots$) and $s=e_{xy}(x+iy)+z$. Then
\begin{equation}\label{eq:chap_6_1_4}
\lim_{n \to \infty}s_{n}=s
\end{equation}
if and only if
\begin{equation}\label{eq:chap_6_1_5}
\lim_{n \to \infty}x_{n}=x\texttt{, }\lim_{n \to \infty}y_{n}=y\texttt{, and }\lim_{n \to \infty}z_{n}=z.
\end{equation}
\end{theorem}

To prove this theorem, we first assume that conditions~(\ref{eq:chap_6_1_5}) hold and obtain condition~(\ref{eq:chap_6_1_4}) from it. According to conditions~(\ref{eq:chap_6_1_5}), there exist, for each positive number $\varepsilon$, positive integers $n_{1}$, $n_{2}$, and $n_{3}$ such that
\[
|x_{n}-x|<\varepsilon/3\texttt{ whenever }n>n_{1},
\]
\[
|y_{n}-y|<\varepsilon/3\texttt{ whenever }n>n_{2},
\]
\[
|z_{n}-z|<\varepsilon/3\texttt{ whenever }n>n_{3}.
\]
Hence, if $n_{0}$ is the largest of the three integers $n_{1}$, $n_{2}$, and $n_{3}$,
\[
|x_{n}-x|<\varepsilon/3\texttt{, }
|y_{n}-y|<\varepsilon/3\texttt{, and }
|z_{n}-z|<\varepsilon/3\texttt{ whenever }n>n_{0}.
\]
Since
\[
|[e_{xy}(x_{n}+iy_{n})+z_{n}]-[e_{xy}(x+iy)+z]|\leq|x_{n}-x|+|y_{n}-y|+|z_{n}-z|,
\]
then,
\[
|s_{n}-s|<\varepsilon/3+\varepsilon/3+\varepsilon/3=\varepsilon\texttt{ whenever }n>n_{0}.
\]
Condition~(\ref{eq:chap_6_1_4}) thus holds.

Conversely, if we start with condition~(\ref{eq:chap_6_1_4}), we know that, for each positive number $\varepsilon$, there exists a positive integer $n_{0}$ such that
\[
|[e_{xy}(x_{n}+iy_{n})+z_{n}]-[e_{xy}(x+iy)+z]|<\varepsilon\texttt{ whenever }n>n_{0},
\]
But
\[
\begin{split}
|x_{n}-x|&\leq|e_{xy}[(x_{n}-x)+i(y_{n}-y)]+(z_{n}-z)] \\
&=|[e_{xy}(x_{n}+iy_{n})+z_{n}]-[e_{xy}(x+iy)+z]|,
\end{split}
\]
\[
\begin{split}
|y_{n}-y|&\leq|e_{xy}[(x_{n}-x)+i(y_{n}-y)]+(z_{n}-z)] \\
&=|[e_{xy}(x_{n}+iy_{n})+z_{n}]-[e_{xy}(x+iy)+z]|,
\end{split}
\]
\[
\begin{split}
|z_{n}-z|&\leq|e_{xy}[(x_{n}-x)+i(y_{n}-y)]+(z_{n}-z)] \\
&=|[e_{xy}(x_{n}+iy_{n})+z_{n}]-[e_{xy}(x+iy)+z]|;
\end{split}
\]
and this means that
\[
|x_{n}-x|<\varepsilon\texttt{, }
|y_{n}-y|<\varepsilon\texttt{, and }
|z_{n}-z|<\varepsilon\texttt{ whenever }n>n_{0}.
\]
That is, conditions~(\ref{eq:chap_6_1_5}) are satisfied.

Note how the theorem enables us to write
\[
\lim_{n \to \infty}[e_{xy}(x_{n}+iy_{n})+z_{n}]
=e_{xy}(\lim_{n \to \infty}x_{n}+i\lim_{n \to \infty}y_{n})+\lim_{n \to \infty}z_{n}
\]
whenever we know that the three limits on the right exist or that the one on the left exists.

\section{Convergence of Series}\label{sec:chap_6_2_52}

An infinite series
\begin{equation}\label{eq:chap_6_2_1}
\sum_{n=1}^{\infty}s_{n}=s_{1}+s_{2}+\cdots+s_{n}+\cdots
\end{equation}
of spatial complex numbers converges to the sum $S$ if the sequence
\begin{equation}\label{eq:chap_6_2_2}
S_{N}=\sum_{n=1}^{N}s_{n}=s_{1}+s_{2}+\cdots+s_{N}\texttt{ }(N=1,2,\ldots)
\end{equation}
of partial sums converges to $S$; we then write
\[
\sum_{n=1}^{\infty}s_{n}=S.
\]
Note that since a sequence can have at most one limit, a series can have at most one sum. When a series dose not converge, we say that it diverges.

\begin{theorem}\label{th:chap_6_2_1}
Suppose that $s_{n}=e_{xy}(x_{n}+iy_{n})+z_{n}$ ($n=1,2,\ldots$) and $S=e_{xy}(X+iY)+Z$. Then
\begin{equation}\label{eq:chap_6_2_3}
\sum_{n=1}^{\infty}s_{n}=S
\end{equation}
if and only if
\begin{equation}\label{eq:chap_6_2_4}
\sum_{n=1}^{\infty}x_{n}=X\texttt{, }\sum_{n=1}^{\infty}y_{n}=Y\texttt{, and }\sum_{n=1}^{\infty}z_{n}=Z.
\end{equation}
\end{theorem}

This theorem tells us, of course, that one can write
\[
\sum_{n=1}^{\infty}[e_{xy}(x_{n}+iy_{n})+z_{n}]=e_{xy}(\sum_{n=1}^{\infty}x_{n}+i\sum_{n=1}^{\infty}y_{n})+\sum_{n=1}^{\infty}z_{n}
\]
whenever it is known that the three series on the right converge or that the one on the
left does.

To prove the theorem, we first write the partial sums~(\ref{eq:chap_6_2_2}) as
\begin{equation}\label{eq:chap_6_2_5}
S_{N}=e_{xy}(X_{n}+iY_{n})+Z_{n},
\end{equation}
where
\[
X_{n}=\sum_{n=1}^{N}x_{n}\texttt{, }
Y_{n}=\sum_{n=1}^{N}y_{n}\texttt{, and }
Z_{n}=\sum_{n=1}^{N}z_{n}.
\]

Now statement~(\ref{eq:chap_6_2_3}) is true if and only if
\begin{equation}\label{eq:chap_6_2_6}
\lim_{n \to \infty}S_{N}=S;
\end{equation}
and, in view of relation~(\ref{eq:chap_6_2_5}) and the theorem on sequences in Sec.~(\ref{sec:chap_6_1_51}), limit~(\ref{eq:chap_6_2_6}) holds if
and only if
\begin{equation}\label{eq:chap_6_2_7}
\lim_{n \to \infty}X_{N}=X\texttt{, }
\lim_{n \to \infty}Y_{N}=Y\texttt{, and }
\lim_{n \to \infty}Z_{N}=Z.
\end{equation}
Limits~(\ref{eq:chap_6_2_7}) therefore imply statement~(\ref{eq:chap_6_2_3}), and conversely. Since $X_{n}$, $Y_{n}$, and $Z_{n}$ are the partial sums of the series~(\ref{eq:chap_6_2_4}), the theorem here is proved.

By recalling from calculus that the $n$th term of a convergent series of real numbers approaches zero as $n$ tends to infinity, we can see immediately from the theorems in this and the previous section that the same is true of a convergent series of spatial complex numbers. That is, a necessary condition for the convergence of series~(\ref{eq:chap_6_2_1}) is that
\[
\lim_{n \to \infty}s_{n}=0.
\]
The terms of a convergent series of spatial complex numbers are, therefore, bounded. To be specific, there exists a positive constant $M$ such that $|s_{n}|<M$ for each positive integer $n$.

For another important property of series of spatial complex numbers, we assume that series~(\ref{eq:chap_6_2_1}) is  absolutely convergent. That is, when $s_{n}=e_{xy}(x_{n}+iy_{n})+z_{n}$, the series
\[
\sum_{n=1}^{\infty}|s_{n}|=\sum_{n=1}^{\infty}\sqrt{x_{n}^{2}+y_{n}^{2}+z_{n}^{2}}
\]
of real numbers $\sqrt{x_{n}^{2}+y_{n}^{2}+z_{n}^{2}}$ converges. Since
\[
|x_{n}|\leq\sqrt{x_{n}^{2}+y_{n}^{2}+z_{n}^{2}}\texttt{, }
|y_{n}|\leq\sqrt{x_{n}^{2}+y_{n}^{2}+z_{n}^{2}}\texttt{, and }
|z_{n}|\leq\sqrt{x_{n}^{2}+y_{n}^{2}+z_{n}^{2}},
\]
we know from the comparison test in calculus that the three series
\[
\sum_{n=1}^{\infty}|x_{n}|\texttt{, }
\sum_{n=1}^{\infty}|y_{n}|\texttt{, and }
\sum_{n=1}^{\infty}|z_{n}|
\]
must converge. Moreover, since the absolute convergence of a series of real numbers implies the convergence of the series itself, it follows that there are real numbers $X$, $Y$, and $Z$ to which series~(\ref{eq:chap_6_2_4}) converge. According to the theorem in this section, then, series~(\ref{eq:chap_6_2_1}) converges. Consequently,  absolute convergence of a series of spatial complex numbers implies convergence of that series.

In establishing the fact that the sum of a series is a given number $S$, it is often convenient to define the  remainder $\rho_{N}$ after $N$ terms:
\begin{equation}\label{eq:chap_6_2_8}
\rho_{N}=S-S_{N}.
\end{equation}
Thus  $S=S_{N}+\rho_{N}$; and, since $|S^{N}-S|=|\rho_{N}-0|$, we see that a series converges to a number $S$ if  and only if the sequence of remainders tends to zero. We shall make considerable use of this observation in our  treatment of power series. They are series of the form
\[
\sum_{n=0}^{\infty}a_{n}(s-s_{0})^{n}=a_{0}+a_{1}(s-s_{0})+a_{2}(s-s_{0})^{2}+\cdots+a_{n}(s-s_{0})^{n}+\cdots
\]
where $s_{0}$ and the coefficients $a_{n}$ are spatial complex constants and $s$ may be any point in a stated space containing $s_{0}$. In such series, involving a variable $s$, we shall denote sums, partial sums, and remainders by  $S(s)$, $S_{N}(s)$, and $\rho_{N}(s)$, respectively.

\section{Taylor Series}\label{sec:chap_6_3_53}

We turn now to Taylor's theorem, which is one of the most important results of the chapter.

\begin{theorem}\label{th:chap_6_3_1}
Suppose that a function $f$ is analytic throughout an open sphere $(|s-s_{0}|<R_{0})$, centered at $s_{0}$ and with radius $R_{0}$. Then $f(s)$ has the power series representation
\begin{equation}\label{eq:chap_6_3_1}
\begin{split}
f(s)&=\sum_{n=0}^{\infty}a_{n}(s-s_{0})^{n}=\sum_{n=0}^{\infty}a_{n}\{e_{xy}[(x-x_{0})+i(y-y_{0})]+(z-z_{0})\}^{n}\\
&=e_{xy}\{\sum_{n=0}^{\infty}a_{n}[(x-x_{0})+i(y-y_{0})+(z-z_{0})]^{n}-\sum_{n=0}^{\infty}a_{n,z}(z-z_{0})^{n}\}\\
&+\sum_{n=0}^{\infty}a_{n,z}(z-z_{0})^{n}\texttt{ }(|s-s_{0}|<R_{0})
\end{split}
\end{equation}
where
\begin{equation}\label{eq:chap_6_3_2}
a_{n}=\frac{f^{(n)}(s_{0})}{n!}=e_{xy}(a_{n,x}+ia_{n,y})+a_{n,z}\texttt{ }(n=0,1,2,\ldots).
\end{equation}
\end{theorem}
That is, Series~(\ref{eq:chap_6_3_1}) converges to $f(s)$ when $s$ lies in the stated open sphere.

This is the expansion of $f(s)$ into a Taylor series about the point $s_{0}$. It is the familiar Taylor series from calculus and two-dimensional complex variables, adapted to functions of a spatial complex variable. With the agreement that
\[
f^{(0)}(s_{0})=f(s_{0})\texttt{ and }0!=1,
\]
Series~(\ref{eq:chap_6_3_1}) can, of course, be written
\begin{equation}\label{eq:chap_6_3_3}
f(s)=f(s_{0})+\frac{f'(s_{0})}{1!}(s-s_{0})+\frac{f^{(2)}(s_{0})}{2!}(s-s_{0})^{2}+\cdots
\texttt{ }(|s-s_{0}|<R_{0}).
\end{equation}

Any function which is analytic everywhere at a point $s_{0}$ must have a Taylor series about $s_{0}$. For, if $f$ is analytic at $s_{0}$, it is analytic throughout some neighborhood $|s-s_{0}|<\varepsilon$ of that point; and $\varepsilon$ may serve as the value of $R_{0}$ in the statement of Taylor's theorem. Also, if $f$ is entire, $R_{0}$ can be chosen arbitrarily large; and the condition of validity becomes $|s-s_{0}|<\infty$. The series then converges to $f(s)$ at each point $s$ in the finite space.

\subsection{Taylor Series with Circle Integrals}\label{sec:chap_6_3_53_C}

We first prove the theorem when $s_{0}=0$, in which case series~(\ref{eq:chap_6_3_1}) becomes
\begin{equation}\label{eq:chap_6_3_4}
f(s)=\sum_{n=0}^{\infty}\frac{f^{(n)}(0)}{n!}s^{n}\texttt{ }(|s|<R_{0})
\end{equation}
and is called a Maclaurin series. The proof when $s_{0}$ is arbitrary will follow as an immediate consequence.

To begin the derivation of representation~(\ref{eq:chap_6_3_4}), we write $|s|=r$ and let $C_{0}$ denote any positively oriented spatial circle on a closed sphere $|s|=r_{0}$, where $r<r_{0}<R_{0}$. Since $f$ is analytic inside and on the sphere $|s|=r_{0}$ and the spatial circle $C_{0}$ and since the point $s$ is interior to $C_{0}$, the Cauchy integral formula applies:
\begin{equation}\label{eq:chap_6_3_5}
f(s)=\frac{1}{2\pi i}\int_{C_{0}}\frac{f(\varsigma)}{\varsigma-s}d\varsigma.
\end{equation}
Now the factor $1/(\varsigma-s)$ in the integrand here can be put in the form
\begin{equation}\label{eq:chap_6_3_6}
\frac{1}{\varsigma-s}=\frac{1}{\varsigma}\cdot\frac{1}{1-(s/\varsigma)}
\end{equation}
and we know
\begin{equation}\label{eq:chap_6_3_7}
\frac{1}{1-s}=\sum_{n=0}^{N-1}s^{n}+\frac{s^{N}}{1-s}
\end{equation}
when $s$ is any spatial complex number other than unity. Replacing $s$ by $s/\varsigma$ in expression~(\ref{eq:chap_6_3_7}), then, we can rewrite equation~(\ref{eq:chap_6_3_6}) as
\begin{equation}\label{eq:chap_6_3_8}
\frac{1}{\varsigma-s}=\sum_{n=0}^{N-1}\frac{1}{\varsigma^{n+1}}s^{n}+s^{N}\frac{1}{(\varsigma-s)\varsigma^{N}}
\end{equation}
Multiplying through this equation by $f(\varsigma)$ and then integrating each side with respect
to $\varsigma$ around $C_{0}$, we find that
\[
\int_{C_{0}}\frac{f(\varsigma)d\varsigma}{\varsigma-s}
=\sum_{n=0}^{N-1}s^{n}\int_{C_{0}}\frac{f(\varsigma)d\varsigma}{\varsigma^{n+1}}
+s^{N}\int_{C_{0}}\frac{f(\varsigma)d\varsigma}{(\varsigma-s)\varsigma^{N}}.
\]
In view of expression~(\ref{eq:chap_6_3_5}) and the fact (Sec.~(\ref{sec:chap_4_13_48})) that
\[
\frac{1}{2\pi i}\int_{C_{0}}\frac{f(\varsigma)d\varsigma}{\varsigma^{n+1}}=\frac{f^{(n)}(0)}{n!}\texttt{ }(n=1,2,\ldots).
\]
this reduces, after we multiply through by $1/(2\pi i)$, to
\begin{equation}\label{eq:chap_6_3_9}
f(s)=\sum_{n=0}^{N-1}\frac{f^{(n)}(0)}{n!}s^{n}+\rho_{N}(s)
\end{equation}
where
\begin{equation}\label{eq:chap_6_3_10}
\rho_{N}(s)=\frac{s^{N}}{2\pi i}\int_{C_{0}}\frac{f(\varsigma)d\varsigma}{(\varsigma-s)\varsigma^{N}}.
\end{equation}

Representation~(\ref{eq:chap_6_3_4}) now follows once it is shown that
\begin{equation}\label{eq:chap_6_3_11}
\lim_{N \to \infty}\rho_{N}(s)=0.
\end{equation}
To accomplish this, we recall that $|s|=r$ and that $C_{0}$ has radius $r_{0}$, where $r_{0}>r$. Then,
if $\varsigma$ is a point on $C_{0}$, we can see that
\[
|\varsigma-s|\geq||\varsigma|-|s||=r_{0}-r.
\]
Consequently, if $M$ denotes the maximum value of $|f(s)|$ on $C_{0}$,
\[
|\rho_{N}(s)|\leq\frac{r^{N}}{2\pi}\cdot\frac{M}{(r_{0}-r)r_{0}^{N}}2\pi r_{0}=\frac{Mr_{0}}{r_{0}-r}(\frac{r}{r_{0}})^{N}.
\]
Inasmuch as $(r/r_{0})<1$, limit~(\ref{eq:chap_6_3_11}) clearly holds.

To verify the theorem when the spatial disk of radius $R_{0}$ is centered at an arbitrary point $s_{0}$, we suppose that $f$ is analytic when $|s-s_{0}|<R_{0}$ and note that the composite function $f(s+s_{0})$ must be analytic when $|(s+s_{0})-s_{0}|<R_{0}$. This last inequality is, of course, just $|s|<R_{0}$; and, if we write $g(s)=f(s+s_{0})$, the analyticity of $g$ in the spatial disk $|s|<R_{0}$ ensures the existence of a Maclaurin series representation:
\[
g(s)=\sum_{n=0}^{\infty}\frac{g^{(n)}(0)}{n!}s^{n}\texttt{ }(|s|<R_{0}).
\]
That is
\[
f(s+s_{0})=\sum_{n=0}^{\infty}\frac{f^{(n)}(s_{0})}{n!}s^{n}\texttt{ }(|s|<R_{0})
\]
After replacing $s$ by $s-s_{0}$ in this equation and its condition of validity, we have the
desired Taylor series expansion~(\ref{eq:chap_6_3_1}).

When it is known that $f$ is analytic everywhere inside a sphere centered at $s_{0}$, convergence of its Taylor series about $s_{0}$ to $f(s)$ for each point $s$ within that sphere is ensured; no test for the convergence of the series is even required. In fact, according to Taylor's theorem, the series converges to $f(s)$ within the sphere about $s_{0}$ whose radius is the distance from $s_{0}$ to the nearest point $s_{1}$ at which $f$ fails to be analytic.

\subsection{Taylor Series with Surface Integrals}\label{sec:chap_6_3_53_S}

We first prove the theorem when $s_{0}=0$, in which case series~(\ref{eq:chap_6_3_1}) becomes
\begin{equation}\label{eq:chap_6_3_4_S}
f(s)=\sum_{n=0}^{\infty}\frac{f^{(n)}(0)}{n!}s^{n}\texttt{ }(|s|<R_{0})
\end{equation}
and is called a Maclaurin series. The proof when $s_{0}$ is arbitrary will follow as an immediate consequence.

To begin the derivation of representation~(\ref{eq:chap_6_3_4_S}), we write $|s|=r$ and let $S_{0}$ denote a positively oriented spherical surface on a closed sphere $|s|=r_{0}$, where $r<r_{0}<R_{0}$. Since $f$ is analytic inside and on the sphere $|s|=r_{0}$ and the spherical surface $S_{0}$ and since the point $s$ is interior to $S_{0}$, the Cauchy surface integral formula applies:
\begin{equation}\label{eq:chap_6_3_5_S}
f(s)=\frac{1}{4\pi}\iint_{S_{0}}\frac{f(\varsigma)}{(\varsigma-s)^{2}}d\sigma.
\end{equation}
Now the factor $1/(\varsigma-s)^{2}$ in the integrand here can be put in the form
\begin{equation}\label{eq:chap_6_3_6_S}
\frac{1}{(\varsigma-s)^{2}}=\frac{1}{\varsigma^{2}}\cdot\frac{1}{[1-(s/\varsigma)]^{2}}
\end{equation}
and multiplying through equation~(\ref{eq:chap_6_3_7}) by $1/(1-s)$ we know
\begin{equation}\label{eq:chap_6_3_7_S}
\frac{1}{(1-s)^{2}}=\sum_{n=0}^{N-1}\frac{s^{n}}{1-s}+\frac{s^{N}}{(1-s)^{2}}
\end{equation}
when $s$ is any spatial complex number other than unity. Replacing $s$ by $s/\varsigma$ in expression~(\ref{eq:chap_6_3_7_S}) or directly multiplying through equation~(\ref{eq:chap_6_3_8}) by $1/(\varsigma-s)$, then, we can rewrite equation~(\ref{eq:chap_6_3_6_S}) as
\[
\begin{split}
\frac{1}{(\varsigma-s)^{2}}&=\sum_{n=0}^{N-1}\frac{1}{(\varsigma-s)\varsigma^{n+1}}s^{n}+s^{N}\frac{1}{(\varsigma-s)^{2}\varsigma^{N}}\\
&=(\sum_{n=0}^{N-1}\frac{1}{\varsigma^{n+1}}s^{n})^{2}+\frac{s^{N}}{(\varsigma-s)\varsigma^{N}}\sum_{n=0}^{N-1}\frac{1}{\varsigma^{n+1}}s^{n}+\frac{s^{N}}{(\varsigma-s)^{2}\varsigma^{N}}\\
&=\sum_{n=0}^{N-1}\frac{n+1}{\varsigma^{n+2}}s^{n}+\frac{s^{N}}{\varsigma^{N}}[\sum_{n=0}^{N-2}\frac{N-n-1}{\varsigma^{n+2}}s^{n}
+\frac{1}{\varsigma-s}(\sum_{n=0}^{N-1}\frac{1}{\varsigma^{n+1}}s^{n}+\frac{1}{\varsigma-s})]
\end{split}
\]
Multiplying through this equation by $f(\varsigma)$ and then integrating each side with respect
to $\varsigma$ around $S_{0}$, we find that
\begin{equation}\label{eq:chap_6_3_8_S}
\begin{split}
\iint_{S_{0}}\frac{f(\varsigma)d\sigma}{(\varsigma-s)^{2}}
&=\sum_{n=0}^{N-1}s^{n}\iint_{S_{0}}\frac{n+1}{\varsigma^{n+2}}f(\varsigma)d\sigma \\
&+s^{N}\sum_{n=0}^{N-2}s^{n}\iint_{S_{0}}\frac{N-n-1}{\varsigma^{n+N+2}}f(\varsigma)d\sigma \\
&+s^{N}[\sum_{n=0}^{N-1}s^{n}\iint_{S_{0}}\frac{f(\varsigma)d\sigma}{(\varsigma-s)\varsigma^{n+N+1}}
+\iint_{S_{0}}\frac{f(\varsigma)d\sigma}{(\varsigma-s)^{2}\varsigma^{N}}].
\end{split}
\end{equation}

In view of expression~(\ref{eq:chap_6_3_5_S}) and expression~(\ref{eq:chap_6_3_9}), and the uniqueness of series representations (Sec.~(\ref{sec:chap_6_10_60})) that
\[
\frac{1}{4\pi}\iint_{S}\frac{(n+1)f(\varsigma)}{\varsigma^{n+2}}d\sigma=\frac{f^{(n)}(0)}{n!}\texttt{ }(n=1,2,\ldots).
\]
this reduces, after we multiply through by $1/(4\pi)$, to
\begin{equation}\label{eq:chap_6_3_9_S}
f(s)=\sum_{n=0}^{N-1}\frac{f^{(n)}(0)}{n!}s^{n}+\rho_{N}(s)+\rho_{N-1}(s)+\rho_{N-2}(s)
\end{equation}
where
\begin{equation}\label{eq:chap_6_3_10_S}
\begin{split}
&\texttt{ }\texttt{ }\rho_{N}(s)=\frac{s^{N}}{4\pi}\iint_{S_{0}}\frac{f(\varsigma)d\sigma}{(\varsigma-s)^{2}\varsigma^{N}}, \\
&\rho_{N-1}(s)=\frac{s^{N}}{4\pi}\sum_{n=0}^{N-1}s^{n}\iint_{S_{0}}\frac{f(\varsigma)d\sigma}{(\varsigma-s)\varsigma^{n+N+1}},\\
&\rho_{N-2}(s)=\frac{s^{N}}{4\pi}\sum_{n=0}^{N-2}s^{n}\iint_{S_{0}}\frac{(N-n-1)f(\varsigma)d\sigma}{\varsigma^{n+N+2}}.
\end{split}
\end{equation}

Representation~(\ref{eq:chap_6_3_4_S}) now follows once it is shown that
\begin{equation}\label{eq:chap_6_3_11_S}
\lim_{N \to \infty}\rho_{N}(s)=0,\texttt{ }\lim_{N \to \infty}\rho_{N-1}(s)=0,\texttt{ }\lim_{N \to \infty}\rho_{N-2}(s)=0.
\end{equation}
To accomplish this, we recall that $|s|=r$ and that $S_{0}$ has radius $r_{0}$, where $r_{0}>r$. Then,
if $\varsigma$ is a point on $S_{0}$, we can see that
\[
|\varsigma-s|\geq||\varsigma|-|s||=r_{0}-r,
\]
\[
\sum_{n=0}^{N-1}(\frac{r}{r_{0}})^{n}<\frac{r_{0}}{r_{0}-r},\texttt{ }
\sum_{n=0}^{N-2}(\frac{r}{r_{0}})^{n}(N-n-1)<\frac{r_{0}N}{r_{0}-r}.
\]
Consequently, if $M$ denotes the maximum value of $|f(s)|$ on $S_{0}$,
\[
\begin{split}
&\texttt{ }\texttt{ }|\rho_{N}(s)|\leq\frac{r^{N}}{4\pi}\cdot\frac{4\pi r_{0}^{2}M}{(r_{0}-r)^{2}r_{0}^{N}}=\frac{Mr_{0}^{2}}{(r_{0}-r)^{2}}(\frac{r}{r_{0}})^{N},\\
&|\rho_{N-1}(s)|<\frac{r^{N}}{4\pi}\cdot\frac{4\pi r_{0}^{2}M}{(r_{0}-r)^{2}r_{0}^{N}}\cdot\frac{r_{0}}{r_{0}-r}=\frac{Mr_{0}^{3}}{(r_{0}-r)^{3}}(\frac{r}{r_{0}})^{N},\\
&|\rho_{N-2}(s)|<\frac{r^{N}}{4\pi}\cdot\frac{4\pi r_{0}^{2}M}{(r_{0}-r)^{2}r_{0}^{N}}\cdot\frac{r_{0}N}{r_{0}-r}=\frac{Mr_{0}^{3}N}{(r_{0}-r)^{3}}(\frac{r}{r_{0}})^{N}.
\end{split}
\]
Inasmuch as $(r/r_{0})<1$, limit~(\ref{eq:chap_6_3_11_S}) clearly holds.

To verify the theorem when the sphere of radius $R_{0}$ is centered at an arbitrary point $s_{0}$, we suppose that $f$ is analytic when $|s-s_{0}|<R_{0}$ and note that the composite function $f(s+s_{0})$ must be analytic when $|(s+s_{0})-s_{0}|<R_{0}$. This last inequality is, of course, just $|s|<R_{0}$; and, if we write $g(s)=f(s+s_{0})$, the analyticity of $g$ in the spatial disk $|s|<R_{0}$ ensures the existence of a Maclaurin series representation:
\[
g(s)=\sum_{n=0}^{\infty}\frac{g^{(n)}(0)}{n!}s^{n}\texttt{ }(|s|<R_{0}).
\]
That is
\[
f(s+s_{0})=\sum_{n=0}^{\infty}\frac{f^{(n)}(s_{0})}{n!}s^{n}\texttt{ }(|s|<R_{0})
\]
After replacing $s$ by $s-s_{0}$ in this equation and its condition of validity, we have the
desired Taylor series expansion~(\ref{eq:chap_6_3_1}).

When it is known that $f$ is analytic everywhere inside a sphere centered at $s_{0}$, convergence of its Taylor series about $s_{0}$ to $f(s)$ for each point $s$ within that sphere is ensured; no test for the convergence of the series is even required. In fact, according to Taylor's theorem, the series converges to $f(s)$ within the sphere about $s_{0}$ whose radius is the distance from $s_{0}$ to the nearest point $s_{1}$ at which $f$ fails to be analytic.

\section{Examples}\label{sec:chap_6_4_54}

When it is known that $f$ is analytic everywhere inside a sphere centered at $s_{0}$, convergence of its Taylor series about $s_{0}$ to $f(s)$ for each point $s$ within that sphere is ensured; no test for the convergence of  the series is required. In fact, according to Taylor's theorem, the series converges to $f(s)$ within the sphere about $s_{0}$ whose radius is the distance from $s_{0}$ to the nearest point $s_{1}$ where $f$ fails to be analytic. In Sec.~(\ref{sec:chap_6_9_59}), we shall find that this is actually the largest sphere centered at  $s_{0}$ such that the series converges to $f(s)$ for all $s$ interior to it.

Also, in Sec.~(\ref{sec:chap_6_10_60}), we shall see that if there are constants $a_{n}$ $(n=0,1,2,\ldots)$ such
that
\[
f(s)=\sum_{n=0}^{\infty}a_{n}(s-s_{0})^{n}
\]
for all points $s$ interior to some sphere centered at $s_{0}$, then the power series here must be the Taylor series for $f$ about $s_{0}$, regardless of how those constants arise. This observation often allows us to find the coefficients $a_{n}$ in Taylor series in more efficient ways than by appealing directly to the formula $a_{n}=f^{(n)}(s_{0})/n!$ in Taylor's theorem.

In the following examples, we use the formula in Taylor's theorem to find the Maclaurin series expansions of some fairly simple functions, and we emphasize the use of those expansions in finding other representations. In our examples, we shall freely use expected properties of convergent series.%, such as those verified in Exercises 7 and  8,  Sec. 52.

\begin{example}\label{ex:chap_6_4_1}
Since the function $f(s)=e^{s}$ is entire, it has a Maclaurin series representation which is valid for all $s$. Here $f^{(n)}(s)=e^{s}$ and, because $f^{(n)}(0)=1$, it follows that
\begin{equation}\label{eq:chap_6_4_1}
e^{s}=\sum_{n=0}^{\infty}\frac{s^{n}}{n!}\texttt{ }(|s|<\infty).
\end{equation}
Note that if $s=e_{xy}(x+iy)+0$, then expansion~(\ref{eq:chap_6_4_1}) becomes
\[
e^{s}=e^{e_{xy}(x+iy)}=e_{xy}\sum_{n=0}^{\infty}\frac{(x+iy)^{n}}{n!}\texttt{ }(|x+iy|<\infty)
\]
or if $s=e_{xy}(0+i0)+z$, then expansion~(\ref{eq:chap_6_4_1}) becomes
\[
e^{s}=e^{z}=\sum_{n=0}^{\infty}\frac{z^{n}}{n!}\texttt{ }(-\infty<z<\infty).
\]

The entire function $s^{2}e^{s}$ also has a Maclaurin series expansion. The simplest way to obtain it is to replace  $s$ by $3s$ on each side of equation~(\ref{eq:chap_6_4_1}) and then multiply through the resulting equation by $s^{2}$:
\[
s^{2}e^{3s}=\sum_{n=0}^{\infty}\frac{3^{n}}{n!}s^{n+2}\texttt{ }(|s|<\infty).
\]
Finally, if we replace $n$ by $n-2$ here, we have
\[
s^{2}e^{3s}=\sum_{n=2}^{\infty}\frac{3^{n-2}}{(n-2)!}s^{n}\texttt{ }(|s|<\infty).
\]
\end{example}

\begin{example}\label{ex:chap_6_4_2}
One can use expansion~(\ref{eq:chap_6_4_1}) and the definition (Sec.~(\ref{sec:chap_3_7_33}))
\[
\sin s=\frac{e^{is}-e^{-is}}{2i}
\]
to find the Maclaurin series for the entire function $f(s)=\sin s$. To give the details, we refer to expansion~(\ref{eq:chap_6_4_1}) and write
\[
\sin s=\frac{1}{2i}[\sum_{n=0}^{\infty}\frac{(is)^{n}}{n!}-\sum_{n=0}^{\infty}\frac{(-is)^{n}}{n!}]
      =\frac{1}{2i}\sum_{n=0}^{\infty}[1-(-1)^{n}]\frac{i^{n}s^{n}}{n!}\texttt{ }(|s|<\infty).
\]
But $1-(-1)^{n}=0$ when $n$ is even, and so we can replace $n$ by $2n+1$ in this last series:
\[
\sin s=\frac{1}{2i}\sum_{n=0}^{\infty}[1-(-1)^{2n+1}]\frac{i^{2n+1}s^{2n+1}}{(2n+1)!}\texttt{ }(|s|<\infty).
\]
Inasmuch as
\[
1-(-1)^{2n+1}=2\texttt{ and }i^{2n+1}=(i^{2})^{n}i=(-1)^{n}i,
\]
this reduces to
\begin{equation}\label{eq:chap_6_4_2}
\sin s=\sum_{n=0}^{\infty}(-1)^{n}\frac{s^{2n+1}}{(2n+1)!}\texttt{ }(|s|<\infty).
\end{equation}

Term by term differentiation will be justified in Sec.~(\ref{sec:chap_6_9_59}). Using that procedure here, we differentiate each side of equation~(\ref{eq:chap_6_4_2}) and write
\[
\cos s=\sum_{n=0}^{\infty}\frac{(-1)^{n}}{(2n+1)!}\frac{d}{ds}s^{2n+1}
      =\sum_{n=0}^{\infty}(-1)^{n}\frac{2n+1}{(2n+1)!}s^{2n}\texttt{ }(|s|<\infty).
\]
That is,
\begin{equation}\label{eq:chap_6_4_3}
\cos s=\sum_{n=0}^{\infty}(-1)^{n}\frac{s^{2n}}{(2n)!}\texttt{ }(|s|<\infty).
\end{equation}
\end{example}

\begin{example}\label{ex:chap_6_4_3}
Because $\sinh s=-i\sin(is)$ (Sec.~(\ref{sec:chap_3_8_34}), we need only replace $s$ by $is$ on each side of equation~(\ref{eq:chap_6_4_2}) and multiply through the result by $-i$ to see that
\begin{equation}\label{eq:chap_6_4_4}
\sinh s=\sum_{n=0}^{\infty}\frac{s^{2n+1}}{(2n+1)!}\texttt{ }(|s|<\infty).
\end{equation}
Likewise, since $\cosh s=\cos(is)$, it follows from expansion~(\ref{eq:chap_6_4_3}) that
\begin{equation}\label{eq:chap_6_4_5}
\cosh s=\sum_{n=0}^{\infty}\frac{s^{2n}}{(2n)!}\texttt{ }(|s|<\infty).
\end{equation}

Observe that the Taylor series for $\cosh s$ about the point $s_{0}=-2\pi i$, for example, is obtained by replacing the variable $s$ by $s+2\pi i$ on each side of equation~(\ref{eq:chap_6_4_5}) and then recalling that $\cosh(s+2\pi i)=\cosh s$ for all $s$:
\[
\cosh s=\sum_{n=0}^{\infty}\frac{(s+2\pi i)^{2n}}{(2n)!}\texttt{ }(|s|<\infty).
\]
\end{example}

\begin{example}\label{ex:chap_6_4_4}
Another Maclaurin series representation is
\begin{equation}\label{eq:chap_6_4_6}
\frac{1}{1-s}=\sum_{n=0}^{\infty}s^{n}\texttt{ }(|s|<1).
\end{equation}
The derivatives of the function $f(s)=1/(1-s)$, which fails to be analytic at $s=1$, are
\[
f^{(n)}(s)=\frac{n!}{(1-s)^{n+1}}\texttt{ }(n=0,1,2,\ldots);
\]
and, in particular, $f^{(n)}(0)=n!$. Note that expansion~(\ref{eq:chap_6_4_6}) gives us the sum of an infinite
geometric series, where $s$ is the common ratio of adjacent terms:
\[
1+s+s^{2}+s^{3}+\cdots=\frac{1}{1-s}\texttt{ }(|s|<1).
\]
%This is, of course, the summation formula that was found in another way in the example in Sec.  52.

If we substitute $-s$ for $s$ in equation~(\ref{eq:chap_6_4_6}) and its condition of validity, and note that $|s|<1$ when $|-s|<1$, we see that
\[
\frac{1}{1+s}=\sum_{n=0}^{\infty}(-1)^{n}s^{n}\texttt{ }(|s|<1).
\]

If, on the other hand, we replace the variable $s$ in equation~(\ref{eq:chap_6_4_6}) by $1-s$, we have the Taylor series representation
\[
\frac{1}{s}=\sum_{n=0}^{\infty}(-1)^{n}(s-1)^{n}\texttt{ }(|s-1|<1).
\]
This condition of validity follows from the one associated with expansion~(\ref{eq:chap_6_4_6}) since
$|1-s|<1$ is the same as $|s-1|<1$.
\end{example}

\begin{example}\label{ex:chap_6_4_5}
For our final example, let us expand the function
\[
f(s)=\frac{1+2s^{2}}{s^{3}+s^{5}}=\frac{1}{s^{3}}\cdot\frac{2(1+s^{2})-1}{1+s^{2}}=\frac{1}{s^{3}}\cdot(2-\frac{1}{1+s^{2}})
\]
into a series involving powers of $s$. We cannot find a Maclaurin series for $f(s)$ since it is not analytic at $s=0$. But we do know from expansion~(\ref{eq:chap_6_4_6}) that
\[
\frac{1}{1+s^{2}}=1-s^{2}+s^{4}-s^{6}+s^{8}-\cdots\texttt{ }(|s|<1).
\]
Hence, when $0<|s|<1$,
\[
f(s)=\frac{1}{s^{3}}(2-1+s^{2}-s^{4}+s^{6}-s^{8}+\cdots)=\frac{1}{s^{3}}+\frac{1}{s}-s+s^{3}-s^{5}+\cdots.
\]
We call such terms as $1/s^{3}$ and $1/s$ negative powers of $s$ since they can be written $s^{-3}$ and $s^{-1}$,  respectively. The theory of expansions involving negative powers of $s-s_{0}$ will be discussed in the next section.
\end{example}

\section{Laurent Series}\label{sec:chap_6_5_55}

If a function $f$ fails to be analytic at a point $s_{0}$, we cannot apply Taylor's theorem at that point. It is often possible, however, to find a series representation for $f(s)$ involving both positive and negative powers of  $s-s_{0}$. We now present the theory of such representations, and we begin with Laurent's theorem.

\begin{theorem}\label{th:chap_6_5_1}
Suppose that a function $f$ is analytic throughout a spherical shell domain
$R_{1}<|s-s_{0}|<R_{2}$, centered at $s_{0}$ and contained between two closed spheres $|s-s_{0}|=R_{1}$ and $|s-s_{0}|=R_{2}$. Then, at each point in the domain, $f(s)$ has the series representation
\begin{equation}\label{eq:chap_6_5_1}
f(s)=\sum_{n=0}^{\infty}a_{n}(s-s_{0})^{n}+\sum_{n=1}^{\infty}\frac{b_{n}}{(s-s_{0})^{n}}\texttt{ }(R_{1}<|s-s_{0}|<R_{2})
\end{equation}
\[
=e_{xy}\{\sum_{n=0}^{\infty}a_{n}[(x-x_{0})+i(y-y_{0})+(z-z_{0})]^{n}
-\sum_{n=0}^{\infty}a_{n,z}(z-z_{0})^{n}\}
+\sum_{n=0}^{\infty}a_{n,z}(z-z_{0})^{n}
\]
\[
+e_{xy}\{\sum_{n=1}^{\infty}\frac{b_{n}}{[(x-x_{0})+i(y-y_{0})+(z-z_{0})]^{n}}
-\sum_{n=1}^{\infty}\frac{b_{n,z}}{(z-z_{0})^{n}}\}
+\sum_{n=1}^{\infty}\frac{b_{n,z}}{(z-z_{0})^{n}}.
\]
\end{theorem}

Expansion~(\ref{eq:chap_6_5_1}) is often written
\begin{equation}\label{eq:chap_6_5_4}
f(s)=\sum_{n=-\infty}^{\infty}c_{n}(s-s_{0})^{n}\texttt{ }(R_{1}<|s-s_{0}|<R_{2})
\end{equation}
\[
=e_{xy}\{\sum_{n=0}^{\infty}c_{n}[(x-x_{0})+i(y-y_{0})+(z-z_{0})]^{n}
-\sum_{n=0}^{\infty}c_{n,z}(z-z_{0})^{n}\}
+\sum_{n=0}^{\infty}c_{n,z}(z-z_{0})^{n}
\]
In either of the forms~(\ref{eq:chap_6_5_1}) or~(\ref{eq:chap_6_5_4}), it is called a Laurent series.

\subsection{Laurent Series with Contour Integrals}\label{sec:chap_6_5_55_1}

Let $C$ denote any positively oriented simple closed contour around $s_{0}$ and lying in that domain where $C$ and $s_{0}$ are in a same spatial plane.

Then coefficients of expansion~(\ref{eq:chap_6_5_1}) are
\begin{equation}\label{eq:chap_6_5_2}
a_{n}=\frac{1}{2\pi i}\int_{C}\frac{f(s)ds}{(s-s_{0})^{n+1}}=e_{xy}(a_{n,x}+ia_{n,y})+a_{n,z}\texttt{ }(n=0,1,2,\ldots)
\end{equation}
and
\begin{equation}\label{eq:chap_6_5_3}
b_{n}=\frac{1}{2\pi i}\int_{C}\frac{f(s)ds}{(s-s_{0})^{-n+1}}=e_{xy}(b_{n,x}+ib_{n,y})+b_{n,z}\texttt{ }(n=0,1,2,\ldots);
\end{equation}
and coefficients of expansion~(\ref{eq:chap_6_5_4}) are
\begin{equation}\label{eq:chap_6_5_5}
c_{n}=\frac{1}{2\pi i}\int_{C}\frac{f(s)ds}{(s-s_{0})^{n+1}}=e_{xy}(c_{n,x}+ic_{n,y})+c_{n,z}\texttt{ }(n=0,\pm1,\pm2,\ldots).
\end{equation}

Observe that the integrand in expression~(\ref{eq:chap_6_5_3}) can be written $f(s)(s-s_{0})^{n-1}$. Thus
it is clear that when $f$ is actually analytic throughout the sphere $|s-s_{0}|<R_{2}$, this integrand is too. Hence all of the coefficients $b_{n}$ are zero; and, because (Sec.~(\ref{sec:chap_4_13_48}))
\[
\frac{1}{2\pi i}\int_{C}\frac{f(s)ds}{(s-s_{0})^{n+1}}=\frac{f^{(n)}(s_{0})}{n!}\texttt{ }(n=0,1,2,\ldots),
\]
expansion~(\ref{eq:chap_6_5_1}) reduces to a Taylor series about $s_{0}$.

If, however, $f$ fails to be analytic at $s_{0}$ but is otherwise analytic in the sphere $|s-s_{1}|<R_{2}$, the radius $R_{1}$ can be chosen arbitrarily small. Representation~(\ref{eq:chap_6_5_1}) is then valid in the
deleted sphere $0<|s-s_{0}|<R_{2}$. Similarly, if $f$ is analytic at each point in the finite space exterior to the sphere $|s-s_{0}|=R_{1}$, the condition of validity is $R_{1}<|s-s_{0}|<\infty$. Observe that if $f$ is analytic everywhere in the finite space except at $s_{0}$, series~(\ref{eq:chap_6_5_1}) is valid at each point of analyticity, or when $0<|s-s_{0}|<\infty$.

We shall prove Laurent's theorem first when $s_{0}=0$, in which case the spherical shell is centered at the origin. The verification of the theorem when $s_{0}$ is arbitrary will follow readily.

We start the proof by forming a closed annular region $r_{1}\leq|s|\leq r_{2}$ that is contained in the domain $R_{1}<|s|<R_{2}$ and whose interior contains both the point $s$ and the contour $C$. We let $C_{1}$ and $C_{2}$ denote the spatial circles $|s|=r_{1}$ and $|s|=r_{2}$, respectively, and we assign those two spatial circles a positive orientation. Observe that $f$ is analytic on $C_{1}$ and $C_{2}$, as well as in the spherical shell domain between them.

Next, we construct a positively oriented spatial circle $\gamma$ with center at $s$ and small enough to be completely contained in the interior of the annular region $r_{1}\leq|s|\leq r_{2}$. It then follows from the extension of the  Cauchy-Goursat theorem to integrals of analytic functions around the oriented boundaries of multiply connected domains (Theorem~(\ref{th:chap_4_11_2}), Sec.~(\ref{sec:chap_4_11_46})) that
\[
\int_{C_{2}}\frac{f(\varsigma)d\varsigma}{\varsigma-s}
-\int_{C_{1}}\frac{f(\varsigma)d\varsigma}{\varsigma-s}
-\int_{\gamma}\frac{f(\varsigma)d\varsigma}{\varsigma-s}=0.
\]
But, according to the Cauchy integral formula, the value of the third integral here is $2\pi if(s)$. Hence
\begin{equation}\label{eq:chap_6_5_6}
f(s)=\frac{1}{2\pi i}\int_{C_{2}}\frac{f(\varsigma)d\varsigma}{\varsigma-s}
+\int_{C_{1}}\frac{f(\varsigma)d\varsigma}{s-\varsigma}.
\end{equation}

Now the factor $1/(\varsigma-s)$ in the first of these integrals is the same as in expression~(\ref{eq:chap_6_3_5}), Sec.~(\ref{sec:chap_6_3_53}), where Taylor's theorem was proved; and we shall need here the expansion
\begin{equation}\label{eq:chap_6_5_7}
\frac{1}{\varsigma-s}=\sum_{n=0}^{N-1}\frac{1}{\varsigma^{n+1}}s^{n}+s^{N}\frac{1}{(\varsigma-s)\varsigma^{N}}
\end{equation}
which was used in that earlier section. As for the factor $1/(s-\varsigma)$ in the second integral, an interchange  of $\varsigma$ and $s$ in equation~(\ref{eq:chap_6_5_7}) reveals that
\[
\frac{1}{s-\varsigma}=\sum_{n=0}^{N-1}\frac{1}{\varsigma^{-n}}\frac{1}{s^{n+1}}+\frac{1}{s^{N}}\frac{\varsigma^{N}}{s-\varsigma}.
\]
If we replace the index of summation $n$ here by $n-1$, this expansion takes the form
\begin{equation}\label{eq:chap_6_5_8}
\frac{1}{s-\varsigma}=\sum_{n=1}^{N}\frac{1}{\varsigma^{-n+1}}\frac{1}{s^{n}}+\frac{1}{s^{N}}\frac{\varsigma^{N}}{s-\varsigma}.
\end{equation}
which is to be used in what follows.

Multiplying through equations~(\ref{eq:chap_6_5_7}) and~(\ref{eq:chap_6_5_8}) by $f(\varsigma)/(2\pi i)$ and then integrating each side of the resulting equations with respect to $\varsigma$ around $C_{2}$ and $C_{l}$, respectively, we find from expression~(\ref{eq:chap_6_5_6}) that
\begin{equation}\label{eq:chap_6_5_9}
f(s)=\sum_{n=0}^{N-1}a_{n}s^{n}+\rho_{N}(s)+\sum_{n=1}^{N}\frac{b_{n}}{s^{n}}+\sigma_{N}(s),
\end{equation}
where the numbers $a_{n}$ $(n=0,1,2,\ldots,N - 1)$ and $b_{n}$ $(n=1,2,\ldots,N)$ are given by the equations
\begin{equation}\label{eq:chap_6_5_10}
a_{n}=\frac{1}{2\pi i}\int_{C_{2}}\frac{f(\varsigma)d\varsigma}{\varsigma^{n+1}},\texttt{ }
b_{n}=\frac{1}{2\pi i}\int_{C_{1}}\frac{f(\varsigma)d\varsigma}{\varsigma^{-n+1}}
\end{equation}
and where
\[
\rho_{N}(s)=\frac{s^{N}}{2\pi i}\int_{C_{2}}\frac{f(\varsigma)d\varsigma}{(\varsigma-s)\varsigma^{N}},\texttt{ }
\sigma_{N}(s)=\frac{1}{2\pi is^{N}}\int_{C_{1}}\frac{\varsigma^{N}f(\varsigma)d\varsigma}{s-\varsigma}
\]

As $N$ tends to $\infty$, expression~(\ref{eq:chap_6_5_9}) evidently takes the proper form of a Laurent
series in the domain $R_{1}<|s|<R_{2}$, provided that
\begin{equation}\label{eq:chap_6_5_11}
\lim_{N \to \infty}\rho_{N}(s)=0\texttt{ and }
\lim_{N \to \infty}\sigma_{N}(s)=0.
\end{equation}
These limits are readily established by a method already used in the proof of Taylor's theorem in Sec.~(\ref{sec:chap_6_3_53}). We write $|s|=r$, so that $r_{1}<r<r_{2}$, and let $M$ denote the maximum value of  $|f(s)|$ on $C_{1}$ and $C_{2}$. We also note that if $\varsigma$ is a point on $C_{2}$, then
$|\varsigma-s|\geq r_{2}-r$; and if $\varsigma$ is on $C_{1}$, $|s-\varsigma|\geq r-r_{1}$. This enables us to write
\[
|\rho_{N}(s)|\leq\frac{Mr_{2}}{r_{2}-r}(\frac{r}{r_{2}})^{N}\texttt{ and }
|\sigma_{N}(s)|\leq\frac{Mr_{1}}{r-r_{1}}(\frac{r_{1}}{r})^{N}.
\]
Since $(r/r_{2})<1$ and $(r_{t}/r)<1$, it is now clear that both $\rho_{N}(s)$ and $\sigma_{N}(s)$ have the desired property.

Finally, we need only recall Corollary~(\ref{co:chap_4_11_2}) in Sec.~(\ref{sec:chap_4_11_46}) to see that the contours used in integrals~(\ref{eq:chap_6_5_10}) may be replaced by the contour $C$. This completes the proof of  Laurent's theorem when $s_{0}=0$ since, if $s$ is used instead of $\varsigma$ as the variable of integration, expressions~(\ref{eq:chap_6_5_10}) for the coefficients $a_{n}$ and $b_{n}$ are the same as expressions~(\ref{eq:chap_6_5_2}) and~(\ref{eq:chap_6_5_3}) when $s_{0}=0$ there.

To extend the proof to the general case in which $s_{0}$ is an arbitrary point in the finite space, we let $f$ be a function satisfying the conditions in the theorem; and, just as we did in the proof of Taylor's theorem, we write $g(s)=f(s+s_{0})$. Since $f(s)$ is analytic in the spherical shell $R_{1}<|s-s_{0}|<R_{2}$, the function $f(s+s_{0})$ is analytic when $R_{1}<|(s+s_{0})-s_{0}|<R_{2}$. That is, $g$ is analytic in the spherical shell $R_{1}<|s|<R_{2}$, which is centered at the origin. Now the simple closed contour $C$ in the statement of the theorem has some parametric representation $s=s(t)$ $(a\leq t\leq b)$, where
\begin{equation}\label{eq:chap_6_5_12}
R_{1}<|s(t)-s_{0}|<R_{2}
\end{equation}
for all $t$ in the interval $a\leq t\leq b$. Hence if $\Gamma$ denotes the path
\begin{equation}\label{eq:chap_6_5_13}
s=s(t)-s_{0}\texttt{ }(a\leq t\leq b),
\end{equation}
$\Gamma$ is not only a simple closed contour but, in view of inequalities~(\ref{eq:chap_6_5_12}), it lies in the domain $R_{1}<|s|<R_{2}$. Consequently, $g(s)$ has a Laurent series representation
\begin{equation}\label{eq:chap_6_5_14}
g(s)=\sum_{n=0}^{\infty}a_{n}s^{n}+\sum_{n=1}^{\infty}\frac{b_{n}}{s^{n}}\texttt{ }(R_{1}<|s|<R_{2}),
\end{equation}
where
\begin{equation}\label{eq:chap_6_5_15}
a_{n}=\frac{1}{2\pi i}\int_{\Gamma}\frac{g(s)ds}{s^{n+1}}\texttt{ }(n=0,1,2,\ldots)
\end{equation}
and
\begin{equation}\label{eq:chap_6_5_16}
b_{n}=\frac{1}{2\pi i}\int_{\Gamma}\frac{g(s)ds}{s^{-n+1}}\texttt{ }(n=0,1,2,\ldots).
\end{equation}

Representation~(\ref{eq:chap_6_5_1}) is obtained if we write $f(s+s_{0})$ instead of $g(s)$ in equation~(\ref{eq:chap_6_5_14}) and then replace $s$ by $s-s_{0}$ in the resulting equation, as well as in the condition of validity $R_{1}<|s|<R_{2}$. Expression~(\ref{eq:chap_6_5_15}) for the coefficients $a_{n}$ is, moreover, the same as expression~(\ref{eq:chap_6_5_2}) since
\[
\int_{\Gamma}\frac{g(s)ds}{s^{n+1}}=\int_{a}^{b}\frac{f[s(t)]s'(t)}{(s(t)-s_{0})^{n+1}}dt=\int_{C}\frac{f(s)ds}{(s-s_{0})^{n+1}}.
\]
Similarly, the coefficients $b_{n}$ in expression~(\ref{eq:chap_6_5_16}) are the same as those in expression~(\ref{eq:chap_6_5_3}).

\subsection{Laurent Series with Surface Integrals}\label{sec:chap_6_6_56}

Let $S$ denote any positively oriented simple closed surface around $s_{0}$ and lying in that domain.

Then coefficients of expansion~(\ref{eq:chap_6_5_1}) are
\begin{equation}\label{eq:chap_6_6_2}
a_{n}=\frac{1}{4\pi}\iint_{S}\frac{f(s)d\sigma}{(s-s_{0})^{n+2}}=e_{xy}(a_{n,x}+ia_{n,y})+a_{n,z}\texttt{ }(n=0,1,2,\ldots)
\end{equation}
and
\begin{equation}\label{eq:chap_6_6_3}
b_{n}=\frac{1}{4\pi}\iint_{S}\frac{f(s)d\sigma}{(s-s_{0})^{-n+2}}=e_{xy}(b_{n,x}+ib_{n,y})+b_{n,z}\texttt{ }(n=0,1,2,\ldots);
\end{equation}
and coefficients of expansion~(\ref{eq:chap_6_5_4}) are
\begin{equation}\label{eq:chap_6_6_5}
c_{n}=\frac{1}{4\pi}\iint_{S}\frac{f(s)d\sigma}{(s-s_{0})^{n+2}}=e_{xy}(c_{n,x}+ic_{n,y})+c_{n,z}\texttt{ }(n=0,\pm1,\pm2,\ldots)
\end{equation}

Observe that the integrand in expression~(\ref{eq:chap_6_6_3}) can be written $f(s)(s-s_{0})^{n-2}$. Thus
it is clear that when $f$ is actually analytic throughout the sphere $|s-s_{0}|<R_{2}$, this integrand is too. Hence all of the coefficients $b_{n}$ are zero; and, because expression~(\ref{eq:chap_6_3_9}) and the uniqueness of series representations (Sec.~(\ref{sec:chap_6_10_60}))
\[
\frac{1}{4\pi}\iint_{S}\frac{f(s)d\sigma}{(s-s_{0})^{n+2}}=\frac{f^{(n)}(s_{0})}{n!}\texttt{ }(n=0,1,2,\ldots),
\]
expansion~(\ref{eq:chap_6_5_1}) reduces to a Taylor series about $s_{0}$.

If, however, $f$ fails to be analytic at $s_{0}$ but is otherwise analytic in the sphere $|s-s_{1}|<R_{2}$, the radius $R_{1}$ can be chosen arbitrarily small. Representation~(\ref{eq:chap_6_5_1}) is then valid in the
deleted sphere $0<|s-s_{0}|<R_{2}$. Similarly, if $f$ is analytic at each point in the finite space exterior to the sphere $|s-s_{0}|=R_{1}$, the condition of validity is $R_{1}<|s-s_{0}|<\infty$. Observe that if $f$ is analytic everywhere in the finite space except at $s_{0}$, series~(\ref{eq:chap_6_5_1}) is valid at each point of analyticity, or when $0<|s-s_{0}|<\infty$.

We shall prove Laurent's theorem first when $s_{0}=0$, in which case the spherical shell is centered at the origin. The verification of the theorem when $s_{0}$ is arbitrary will follow readily.

We start the proof by forming a closed spherical shell $r_{1}\leq|s|\leq r_{2}$ that is contained in the domain $R_{1}<|s|<R_{2}$ and whose interior contains both the point $s$ and the surface $S$. We let $S_{1}$ and $S_{2}$ denote the spherical surfaces $|s|=r_{1}$ and $|s|=r_{2}$, respectively, and we assign those two spherical surfaces a positive orientation. Observe that $f$ is analytic on $S_{1}$ and $S_{2}$, as well as in the spherical shell domain between them.

Next, we construct a positively oriented spherical surface $\gamma$ with center at $s$ and small enough to be completely contained in the interior of the spherical shell $r_{1}\leq|s|\leq r_{2}$. It then follows from the extension of the Cauchy-Goursat theorem to surface integrals of analytic functions around the oriented boundaries of multiply connected domains (Theorem~(\ref{th:chap_5_4_4}), Sec.~(\ref{sec:chap_5_4_54})) that
\[
\iint_{S_{2}}\frac{f(\varsigma)d\sigma}{(\varsigma-s)^{2}}
-\iint_{S_{1}}\frac{f(\varsigma)d\sigma}{(\varsigma-s)^{2}}
-\iint_{\gamma}\frac{f(\varsigma)d\sigma}{(\varsigma-s)^{2}}=0.
\]
But, according to the Cauchy surface integral formula, the value of the third integral here is $4\pi f(s)$. Hence
\begin{equation}\label{eq:chap_6_6_6}
f(s)=\frac{1}{4\pi}\iint_{S_{2}}\frac{f(\varsigma)d\sigma}{(\varsigma-s)^{2}}
-\frac{1}{4\pi}\iint_{S_{1}}\frac{f(\varsigma)d\sigma}{(\varsigma-s)^{2}}.
\end{equation}

Now the factor $1/(\varsigma-s)^{2}$ in the first of these integrals is the same as in expression~(\ref{eq:chap_6_3_5_S}), Sec.~(\ref{sec:chap_6_3_53_S}), where Taylor's theorem was proved; and we shall need here the expansion
\begin{equation}\label{eq:chap_6_6_7}
\frac{1}{(\varsigma-s)^{2}}=\sum_{n=0}^{N-1}\frac{1}{(\varsigma-s)\varsigma^{n+1}}s^{n}+s^{N}\frac{1}{(\varsigma-s)^{2}\varsigma^{N}}
\end{equation}
\[
\begin{split}
&=(\sum_{n=0}^{N-1}\frac{1}{\varsigma^{n+1}}s^{n})^{2}+\frac{s^{N}}{(\varsigma-s)\varsigma^{N}}\sum_{n=0}^{N-1}\frac{1}{\varsigma^{n+1}}s^{n}+\frac{s^{N}}{(\varsigma-s)^{2}\varsigma^{N}}\\
&=\sum_{n=0}^{N-1}\frac{n+1}{\varsigma^{n+2}}s^{n}+\frac{s^{N}}{\varsigma^{N}}[\sum_{n=0}^{N-2}\frac{N-n-1}{\varsigma^{n+2}}s^{n}
+\frac{1}{\varsigma-s}(\sum_{n=0}^{N-1}\frac{1}{\varsigma^{n+1}}s^{n}+\frac{1}{\varsigma-s})]
\end{split}
\]
which was used in that earlier section. As for the factor $1/(\varsigma-s)^{2}$ in the second integral, an interchange of $\varsigma$ and $s$ in equation~(\ref{eq:chap_6_6_7}) reveals that
\[
\frac{1}{(s-\varsigma)^{2}}=-\sum_{n=0}^{N-1}\frac{1}{(s-\varsigma)\varsigma^{-n}}\frac{1}{s^{n+1}}
+\frac{1}{s^{N}}\frac{\varsigma^{N}}{(s-\varsigma)^{2}}.
\]
If we replace the index of summation $n$ here by $n-1$, this expansion takes the form
\begin{equation}\label{eq:chap_6_6_8}
\frac{1}{(s-\varsigma)^{2}}=-\sum_{n=1}^{N}\frac{1}{(s-\varsigma)\varsigma^{-n+1}}\frac{1}{s^{n}}+\frac{1}{s^{N}}\frac{\varsigma^{N}}{(s-\varsigma)^{2}}
\end{equation}
\[
\begin{split}
&=(\sum_{n=1}^{N}\frac{1}{\varsigma^{-n+1}}\frac{1}{s^{n}})^{2}-\frac{1}{s^{N}}\frac{\varsigma^{N}}{s-\varsigma}\sum_{n=1}^{N}\frac{1}{\varsigma^{-n+1}}\frac{1}{s^{n}}+\frac{1}{s^{N}}\frac{\varsigma^{N}}{(s-\varsigma)^{2}}\\
&=\sum_{n=1}^{N}\frac{n}{\varsigma^{-n+1}}\frac{1}{s^{n+1}}+\frac{\varsigma^{N}}{s^{N}}[\sum_{n=1}^{N}\frac{N-n+1}{\varsigma^{-n+1}s^{n+1}}
-\frac{1}{s-\varsigma}\sum_{n=1}^{N}\frac{1}{\varsigma^{-n+1}}\frac{1}{s^{n}}+\frac{1}{(s-\varsigma)^{2}}] \\
&=\sum_{n=1}^{N}\frac{n}{\varsigma^{-n+1}}\frac{1}{s^{n+1}}+\frac{\varsigma^{N}}{s^{N}}\sum_{n=1}^{N}\frac{1}{\varsigma^{-n+1}}\frac{1}{s^{n+1}}[(N-n)-\frac{\varsigma}{s-\varsigma}]
+\frac{1}{s^{N}}\frac{\varsigma^{N}}{(s-\varsigma)^{2}}
\end{split}
\]
which is to be used in what follows.

Multiplying through equations~(\ref{eq:chap_6_6_7}) and~(\ref{eq:chap_6_6_8}) by $f(\varsigma)/(4\pi)$ and then integrating each side of the resulting equations with respect to $\varsigma$ around $S_{2}$ and $S_{l}$, respectively, we find from expression~(\ref{eq:chap_6_6_6}) that
\begin{equation}\label{eq:chap_6_6_9}
f(s)=\sum_{n=0}^{N-1}a_{n}s^{n}+\rho_{N}(s)+\rho_{N-1}(s)+\rho_{N-2}(s)+\sum_{n=1}^{N}\frac{b_{n}}{s^{n}}-\sigma_{N}(s)-\sigma_{S}(s),
\end{equation}
where the numbers $a_{n}$ $(n=0,1,2,\ldots,N-1)$ and $b_{n}$ $(n=1,2,\ldots,N)$ are given by the equations
\begin{equation}\label{eq:chap_6_6_10}
a_{n}=\frac{1}{4\pi}\iint_{S_{2}}\frac{f(\varsigma)d\sigma}{(\varsigma-s)\varsigma^{n+1}},\texttt{ }
b_{n}=\frac{1}{4\pi}\iint_{S_{1}}\frac{f(\varsigma)d\sigma}{(s-\varsigma)\varsigma^{-n+1}}
\end{equation}
and where
\[
\begin{split}
&\texttt{ }\texttt{ }\rho_{N}(s)=\frac{s^{N}}{4\pi}\iint_{S_{0}}\frac{f(\varsigma)d\sigma}{(\varsigma-s)^{2}\varsigma^{N}}, \\
&\rho_{N-1}(s)=\frac{s^{N}}{4\pi}\sum_{n=0}^{N-1}s^{n}\iint_{S_{0}}\frac{f(\varsigma)d\sigma}{(\varsigma-s)\varsigma^{n+N+1}},\\
&\rho_{N-2}(s)=\frac{s^{N}}{4\pi}\sum_{n=0}^{N-2}s^{n}\iint_{S_{0}}\frac{(N-n-1)f(\varsigma)d\sigma}{\varsigma^{n+N+2}}
\end{split}
\]
and
\[
\begin{split}
&\sigma_{N}(s)=\frac{1}{4\pi s^{N}}\iint_{S_{1}}\frac{\varsigma^{N}f(\varsigma)d\sigma}{(s-\varsigma)^{2}}, \\
&\sigma_{S}(s)=\frac{1}{4\pi s^{N}}\iint_{S_{1}}\sum_{n=1}^{N}\frac{1}{\varsigma^{-n+1}}\frac{1}{s^{n+1}}[(N-n)-\frac{\varsigma}{s-\varsigma}]\varsigma^{N}f(\varsigma)d\sigma.
\end{split}
\]

As $N$ tends to $\infty$, expression~(\ref{eq:chap_6_6_9}) evidently takes the proper form of a Laurent
series in the domain $R_{1}<|s|<R_{2}$, provided that
\begin{equation}\label{eq:chap_6_6_11}
\lim_{N \to \infty}\rho_{N}(s)=0,\texttt{ }\lim_{N \to \infty}\rho_{N-1}(s)=0,\texttt{ }\lim_{N \to \infty}\rho_{N-2}(s)=0
\end{equation}
and
\begin{equation}\label{eq:chap_6_6_11_b}
\lim_{N \to \infty}\sigma_{N}(s)=0,\texttt{ }\lim_{N \to \infty}\sigma_{S}(s)=0.
\end{equation}
These limits are readily established by a method already used in the proof of Taylor's theorem in Sec.~(\ref{sec:chap_6_3_53_S}). We write $|s|=r$, so that $r_{1}<r<r_{2}$, and let $M$ denote the maximum value of  $|f(s)|$ on $S_{1}$ and $S_{2}$. We also note that if $\varsigma$ is a point on $S_{2}$, then
$|\varsigma-s|\geq r_{2}-r$; and if $\varsigma$ is on $S_{1}$, $|s-\varsigma|\geq r-r_{1}$. This enables us to write
\[
\begin{split}
&\sum_{n=0}^{N-1}(\frac{r}{r_{2}})^{n}<\frac{r_{2}}{r_{2}-r}, \\
&\sum_{n=0}^{N-2}(\frac{r}{r_{2}})^{n}(N-n-1)<\frac{r_{2}N}{r_{2}-r}, \\
&\sum_{n=1}^{N}(\frac{r_{1}}{r})^{n}[(N-n)-\frac{r_{1}}{r-r_{1}}]<\frac{r_{1}N}{r-r_{1}},
\end{split}
\]
\[
\begin{split}
&\texttt{ }\texttt{ }|\rho_{N}(s)|\leq\frac{r^{N}}{4\pi}\cdot\frac{4\pi r_{2}^{2}M}{(r_{2}-r)^{2}r_{2}^{N}}=\frac{Mr_{2}^{2}}{(r_{2}-r)^{2}}(\frac{r}{r_{2}})^{N},\\
&|\rho_{N-1}(s)|<\frac{r^{N}}{4\pi}\cdot\frac{4\pi r_{2}^{2}M}{(r_{2}-r)^{2}r_{2}^{N}}\cdot\frac{r_{2}}{r_{2}-r}=\frac{Mr_{2}^{3}}{(r_{2}-r)^{3}}(\frac{r}{r_{2}})^{N},\\
&|\rho_{N-2}(s)|<\frac{r^{N}}{4\pi}\cdot\frac{4\pi r_{2}^{2}M}{(r_{2}-r)^{2}r_{2}^{N}}\cdot\frac{r_{2}N}{r_{2}-r}=\frac{Mr_{2}^{3}N}{(r_{2}-r)^{3}}(\frac{r}{r_{2}})^{N}.
\end{split}
\]
and
\[
\begin{split}
&|\sigma_{N}(s)|\leq\frac{r_{1}^{N}}{4\pi}\cdot\frac{4\pi r_{1}^{2}M}{(r-r_{1})^{2}r^{N}}=\frac{Mr_{1}^{2}}{(r-r_{1})^{2}}(\frac{r_{1}}{r})^{N}, \\
&|\sigma_{S}(s)|<\frac{r_{1}^{N}}{4\pi}\cdot\frac{4\pi r_{1}^{2}M}{(r-r_{1})^{2}r^{N}}\cdot\frac{r_{1}N}{r-r_{1}}=\frac{Mr_{1}^{3}N}{(r-r_{1})^{3}}(\frac{r_{1}}{r})^{N}.
\end{split}
\]
Since $(r/r_{2})<1$ and $(r_{t}/r)<1$, it is now clear that all of $\rho_{N}(s)$, $\rho_{N-1}(s)$, $\rho_{N-2}(s)$, $\sigma_{N}(s)$, and $\sigma_{S}(s)$ have the desired property.

Finally, we need only recall Theorem~(\ref{th:chap_5_4_4}), Sec.~(\ref{sec:chap_5_4_54}) to see that the contours used in integrals~(\ref{eq:chap_6_6_10}) may be replaced by the closed surface $S$. This completes the proof of  Laurent's theorem when $s_{0}=0$ since, if $s$ is used instead of $\varsigma$ as the variable of integration, expressions~(\ref{eq:chap_6_6_10}) for the coefficients $a_{n}$ and $b_{n}$ are the same as expressions~(\ref{eq:chap_6_6_2}) and~(\ref{eq:chap_6_6_3}) when $s_{0}=0$ there.

To extend the proof to the general case in which $s_{0}$ is an arbitrary point in the finite space, we let $f$ be a function satisfying the conditions in the theorem; and, just as we did in the proof of Taylor's theorem, we write $g(s)=f(s+s_{0})$. Since $f(s)$ is analytic in the spherical shell $R_{1}<|s-s_{0}|<R_{2}$, the function $f(s+s_{0})$ is analytic when $R_{1}<|(s+s_{0})-s_{0}|<R_{2}$. That is, $g$ is analytic in the spherical shell $R_{1}<|s|<R_{2}$, which is centered at the origin. Now the simple closed surface $S$ in the statement of the theorem has some parametric representation $s=s(t)$ $(a\leq t\leq b)$, where
\begin{equation}\label{eq:chap_6_6_12}
R_{1}<|s(t)-s_{0}|<R_{2}
\end{equation}
for all $t$ in the interval $a\leq t\leq b$. Hence if $\Gamma$ denotes the path
\begin{equation}\label{eq:chap_6_6_13}
s=s(t)-s_{0}\texttt{ }(a\leq t\leq b),
\end{equation}
$\Gamma$ is not only a simple closed surface but, in view of inequalities~(\ref{eq:chap_6_6_12}), it lies in the domain $R_{1}<|s|<R_{2}$. Consequently, $g(s)$ has a Laurent series representation
\begin{equation}\label{eq:chap_6_6_14}
g(s)=\sum_{n=0}^{\infty}a_{n}s^{n}+\sum_{n=1}^{\infty}\frac{b_{n}}{s^{n}}\texttt{ }(R_{1}<|s|<R_{2}),
\end{equation}
where
\begin{equation}\label{eq:chap_6_6_15}
a_{n}=\frac{1}{4\pi}\int_{\Gamma}\frac{g(s)d\sigma}{s^{n+2}}\texttt{ }(n=0,1,2,\ldots)
\end{equation}
and
\begin{equation}\label{eq:chap_6_6_16}
b_{n}=\frac{1}{4\pi}\int_{\Gamma}\frac{g(s)d\sigma}{s^{-n+2}}\texttt{ }(n=0,1,2,\ldots).
\end{equation}

Representation~(\ref{eq:chap_6_5_1}) is obtained if we write $f(s+s_{0})$ instead of $g(s)$ in equation~(\ref{eq:chap_6_6_14}) and then replace $s$ by $s-s_{0}$ in the resulting equation, as well as in the condition of validity $R_{1}<|s|<R_{2}$. Expression~(\ref{eq:chap_6_6_15}) for the coefficients $a_{n}$ is, moreover, the same as expression~(\ref{eq:chap_6_6_2}) since
\[
\int_{\Gamma}\frac{g(s)ds}{s^{n+2}}=\int_{a}^{b}\frac{f[s(t)]\sigma'(t)}{(s(t)-s_{0})^{n+2}}dt=\int_{S}\frac{f(s)d\sigma}{(s-s_{0})^{n+2}}.
\]
Similarly, the coefficients $b_{n}$ in expression~(\ref{eq:chap_6_6_16}) are the same as those in expression~(\ref{eq:chap_6_6_3}).

\section{Examples}\label{sec:chap_6_6_56}

The coefficients in a Laurent series are generally found by means other than by appealing directly to their integral representations. This is illustrated in the examples below, where it is always assumed that, when the annular domain is specified, a Laurent series for a given function in unique. As was the case with Taylor series, we defer the proof of such uniqueness until Sec.~(\ref{sec:chap_6_10_60}).

\begin{example}\label{ex:chap_6_6_1}
Replacing $s$ by $1/s$ in the Maclaurin series expansion
\[
e^{s}=\sum_{n=0}^{\infty}\frac{s^{n}}{n!}=1+\frac{s}{1!}+\frac{s^{2}}{2!}+\frac{s^{3}}{3!}+\cdots\texttt{ }(|s|<\infty),
\]
we have the Laurent series representation
\[
e^{1/s}=\sum_{n=0}^{\infty}\frac{1}{n!s^{n}}=1+\frac{1}{1!s}+\frac{1}{2!s^{2}}+\frac{1}{3!s^{3}}+\cdots\texttt{ }(0<|s|<\infty).
\]

Note that no positive powers of $s$ appear here, the coefficients of the positive powers being zero. Note, too, that the coefficient of $1/s$ is unity; and, according to Laurent's theorem in Sec.~(\ref{sec:chap_6_5_55}), that coefficient is the number
\[
b_{1}=\frac{1}{2\pi i}\int_{C}e^{1/s}ds,
\]
where $C$ is any positively oriented simple closed contour around the origin in a spatial plane. Since
$b_{1}=1$, then,
\[
\int_{C}e^{1/s}ds=2\pi i.
\]
This method of evaluating certain integrals around simple closed contours will be developed in considerable detail in Chap.~(\ref{ch:chap_7}).
\end{example}

\begin{example}\label{ex:chap_6_6_2}
The function $f(s)=1/(s- i)^{2}$ is already in the form of a Laurent series, where $s_{0}=e^{xy}(0+i)+0$. That is,
\[
f(s)=\sum_{n=-\infty}^{\infty}c_{n}(s-i)^{n}\texttt{ }(0<|s-i|<\infty).
\]
where $c_{-2}=1$ and all of the other coefficients are zero. From formula~(\ref{eq:chap_6_5_5}), Sec.~(\ref{sec:chap_6_5_55}), for the coefficients in a Laurent series, we know that
\[
b_{n}=\frac{1}{2\pi i}\int_{C}\frac{ds}{(s-i)^{n+3}}\texttt{ }(n=0,\pm1,\pm2,\ldots),
\]
where $C$ is, for instance, any positively oriented circle $|s-i|=R$ about the point $s_{0}=i$ in a spacial plane including $s_{0}=i$. Thus% (compare Exercise 10, Sec. 40)
\[
\int_{C}\frac{ds}{(s-i)^{n+3}}=\{
    \begin{array}{cc}
        0 & \texttt{ when }n\neq-2 \\
        2\pi i & \texttt{ when }n=-2.
    \end{array}
\]
\end{example}

\begin{example}\label{ex:chap_6_6_3}
The function
\begin{equation}\label{eq:chap_6_6_1}
f(s)=\frac{-1}{(s-1)(s-2)}=\frac{1}{s-1}-\frac{1}{s-2},
\end{equation}
which has the two singular points $s=1$ and $s=2$, is analytic in the domains
\[
|s|<1,\texttt{ }1<|s|<2,\texttt{ and }2<|s|<\infty.
\]
In each of those domains, denoted by $D_{1}$, $D_{2}$, and $D_{3}$, respectively, $f(s)$ has series representations in powers of $s$. They can all be found by reca11ing from Example~(\ref{ex:chap_6_4_4}), Sec.~(\ref{sec:chap_6_4_54}), that
\[
\frac{1}{1-s}=\sum_{n=0}^{\infty}s^{n}\texttt{ }(|s|<1).
\]

The representation in $D_{1}$ is a Maclaurin series. To find it, we write
\[
f(s)=-\frac{1}{1-s}+\frac{1}{2}\cdot\frac{1}{1-(s/2)}
\]
and observe that, since $|s|<1$ and $|s/2|<1$ in $D_{1}$,
\begin{equation}\label{eq:chap_6_6_2}
f(s)=-\sum_{n=0}^{\infty}s^{n}+\sum_{n=0}^{\infty}\frac{s^{n}}{2^{n+1}}=\sum_{n=0}^{\infty}(2^{-n-1}-1)s^{n}\texttt{ }(|s|<1).
\end{equation}

As for the representation in $D_{2}$, we write
\[
f(s)=\frac{1}{s}\cdot\frac{1}{1-(1/s)}+\frac{1}{2}\cdot\frac{1}{1-(s/2)}.
\]
Since $|1/s|<1$ and $|s/2|<1$ when $1<|s|<2$, it follows that
\[
f(s)=\sum_{n=0}^{\infty}\frac{1}{s^{n+1}}+\sum_{n=0}^{\infty}\frac{s^{n}}{2^{n+1}}\texttt{ }(1<|s|<2).
\]
If we replace the index of summation $n$ in the first of these series by $n-1$ and then interchange the two series, we arrive at an expansion having the same form as the one in the statement of Laurent's theorem (Sec.~(\ref{sec:chap_6_5_55})):
\begin{equation}\label{eq:chap_6_6_3}
f(s)=\sum_{n=0}^{\infty}\frac{s^{n}}{2^{n+1}}+\sum_{n=1}^{\infty}\frac{1}{s^{n}}\texttt{ }(1<|s|<2).
\end{equation}
Since there is only one such representation for $f(s)$ in the annulus $1<|s|<2$, expansion~(\ref{eq:chap_6_6_3}) is, in fact, the Laurent series for $f(s)$ there.

The representation of $f(s)$ in the unbounded domain $D_{3}$ is also a Laurent series.
If we put expression~(\ref{eq:chap_6_6_1}) in the form
\[
f(s)=\frac{1}{s}\cdot\frac{1}{1-(1/s)}-\frac{1}{s}\cdot\frac{1}{1-(2/s)}
\]
and observe that $|1/s|<1$ and $|2/s|<1$ when $2<|s|<\infty$, we find that
\[
f(s)=\sum_{n=0}^{\infty}\frac{1}{s^{n+1}}-\sum_{n=0}^{\infty}\frac{2^{n}}{s^{n+1}}=\sum_{n=0}^{\infty}\frac{1-2^{n}}{s^{n+1}}\texttt{ }(2<|s|<\infty).
\]
That is,
\begin{equation}\label{eq:chap_6_6_4}
f(s)=\sum_{n=1}^{\infty}\frac{1-2^{n-1}}{s^{n}}\texttt{ }(2<|s|<\infty).
\end{equation}
\end{example}

\section{Absolute and Uniform Convergence of Power Series}\label{sec:chap_6_7_57}

This section and the three following it are devoted mainly to various properties of
power series. A reader who wishes to simply accept the theorems and any corollaries
there can easily skip their proofs in order to reach Sec.~(\ref{sec:chap_6_11_61}) more quickly.

We recall from Sec.~(\ref{sec:chap_6_2_52}) that a series of spatial complex numbers converges absolutely if the series of absolute values of those numbers converges. The following theorem concerns the absolute convergence of  power series.

\begin{theorem}\label{th:chap_6_7_1}
If a power series
\begin{equation}\label{eq:chap_6_7_1}
f(s)=\sum_{n=0}^{\infty}a_{n}(s-s_{0})^{n}
\end{equation}
converges when $s=s_{1}$ $(s_{1}\neq s_{0})$, then it is absolutely convergent at each point $s$ in the open sphere $|s-s_{0}|<R_{1}$, which is in a closed sphere $|s-s_{0}|=R_{1}$, where $R_{1}=|s_{1}-s_{0}|$.
\end{theorem}

We first prove the theorem when $s=0$, and we assume that the series
\[
\sum_{n=0}^{\infty}a_{n}s_{1}^{n}\texttt{ }(s_{1}\neq0)
\]
converges. The terms $a_{n}s_{1}^{n}$ are thus bounded; that is,
\[
|a_{n}s_{1}^{n}|\leq M\texttt{ }(n=0,1,2,\ldots)
\]
for some positive constant $M$ (see Sec.~(\ref{sec:chap_6_2_52})). If $|s|<|s_{1}|$ and we let $\rho$ denote the modulus $|s/s_{1}|$, we can see that
\[
|a_{n}s^{n}|=|a_{n}s_{1}^{n}||\frac{s}{s_{1}}|^{n}\leq M\rho\texttt{ }(n=0,1,2,\ldots)
\]
where $\rho<1$. Now the series whose terms are the real numbers $M\rho^{n}$ $(n=0,1,2,\ldots)$ is a geometric series, which converges when $\rho<1$. Hence, by the comparison test for series of real numbers, the series
\[
\sum_{n=0}^{\infty}|a_{n}s^{n}|
\]
converges in the open sphere $|s|<|s_{1}|$; and the theorem is proved when $s_{0}=0$.

When $s_{0}$ is any nonzero number, we assume that series~(\ref{eq:chap_6_7_1}) converges at $s=s_{1}$ $(s_{1}\neq s_{0})$. If we write $\omega=s-s_{0}$, series~(\ref{eq:chap_6_7_1}) becomes
\begin{equation}\label{eq:chap_6_7_2}
\sum_{n=0}^{\infty}a_{n}\omega^{n}
\end{equation}
and this series converges at $\omega=s_{1}-s_{0}$. Consequently, since the theorem is known to be true when $s_{0}=0$, we see that series~(\ref{eq:chap_6_7_2}) is absolutely convergent in the open sphere
$|\omega|<|s-s_{0}|$. Finally, by replacing $\omega$ by $s-s_{0}$ in series~(\ref{eq:chap_6_7_2}) and this condition of validity, as well as writing $R_{1}=|s_{1}-s_{0}|$, we arrive at the proof of the theorem as it is stated.

The theorem tells us that the set of all points inside some spatial circle centered at $s_{0}$ is a region of convergence for the power series~(\ref{eq:chap_6_7_1}), provided it converges at some point other than $s_{0}$. The greatest sphere centered at $s_{0}$ such that series~(\ref{eq:chap_6_7_1}) converges at each point inside is called the sphere of convergence of series~(\ref{eq:chap_6_7_1}). The series cannot converge at any point $s_{2}$ outside that sphere, according to the theorem; for if it did, it would converge everywhere inside the sphere centered at $s_{0}$ and passing through $s_{2}$. The first sphere could not, then, be the sphere of convergence.

Our next theorem involves terminology that we must first define. Suppose that the power series~(\ref{eq:chap_6_7_1}) has sphere of convergence $|s-s_{0}|=R$, and let $S(s)$ and $S_{N}(s)$ represent the sum and partial sums, respectively, of that series:
\[
S(s)=\sum_{n=0}^{\infty}a_{n}(s-s_{0})^{n},\texttt{ }
S_{N}(s)=\sum_{n=0}^{N}a_{n}(s-s_{0})^{n}\texttt{ }(|s-s_{0}|<R).
\]
Then write the remainder function
\begin{equation}\label{eq:chap_6_7_3}
\rho_{N}(s)=S(s)-S_{N}(s)\texttt{ }(|s-s_{0}|<R).
\end{equation}
Since the power series converges for any fixed value of $s$ when $|s-s_{0}|<R$, we know that the remainder $\rho_{N}(s)$ approaches zero for any such $s$ as $N$ tends to infinity. According to definition~(\ref{eq:chap_6_1_2}), Sec.~(\ref{sec:chap_6_1_51}), of the limit of a sequence, this means that, corresponding to each positive number $\varepsilon$, there is a positive integer $N_{\varepsilon}$ such that
\begin{equation}\label{eq:chap_6_7_4}
|\rho_{N}(s)|<\varepsilon\texttt{ whenever }N>N_{\varepsilon}.
\end{equation}
When the choice of $N_{\varepsilon}$ depends only on the value of $\varepsilon$ and is independent of the point $s$ taken in a specified region within the sphere of convergence. the convergence is said to be $uniform$ in that region

\begin{theorem}\label{th:chap_6_7_2}
If $s_{1}$ is a point inside the sphere of convergence $|s-s_{0}|=R$ of a power series
\begin{equation}\label{eq:chap_6_7_5}
\sum_{n=0}^{\infty}a_{n}(s-s_{0})^{n},
\end{equation}
then that series must be uniformly convergent in the closed sphere $|s-s_{0}|\leq R_{1}$, where
$R_{1}=|s_{1}-s_{0}|$.
\end{theorem}

As in the proof of Theorem~(\ref{th:chap_6_7_1}), we first treat the case in which $s_{0}=0$. Given that
$s_{1}$ is a point lying inside the sphere of convergence of the series
\begin{equation}\label{eq:chap_6_7_6}
\sum_{n=0}^{\infty}a_{n}s^{n},
\end{equation}
we note that there are points with modulus greater than $|s_{1}|$ for which it converges. According to Theorem~(\ref{th:chap_6_7_1}), then, the series
\begin{equation}\label{eq:chap_6_7_7}
\sum_{n=0}^{\infty}|a_{n}s^{n}|
\end{equation}
converges. Letting $m$ and $N$ denote positive integers, where $m>N$, we can write the remainders of series~(\ref{eq:chap_6_7_6}) and~(\ref{eq:chap_6_7_7}) as
\begin{equation}\label{eq:chap_6_7_8}
\rho_{N}(s)=\lim_{m \to \infty}\sum_{n=N}^{m}a_{n}s^{n},
\end{equation}
and
\begin{equation}\label{eq:chap_6_7_9}
\sigma_{N}=\lim_{m \to \infty}\sum_{n=N}^{m}|a_{n}s^{n}|,
\end{equation}
respectively.

Now, there is
\[
|\rho_{N}(s)|=\lim_{m \to \infty}|\sum_{n=N}^{m}a_{n}s^{n}|,
\]
and, when $|s|\leq|s_{1}|$,
\[
\begin{split}
|\sum_{n=N}^{m}a_{n}s^{n}| &\leq\sum_{n=N}^{m}|a_{n}||s|^{n} \\
                           &\leq\sum_{n=N}^{m}|a_{n}||s_{1}|^{n}=\sum_{n=N}^{m}|a_{n}s_{1}^{n}|.
\end{split}
\]
Hence
\begin{equation}\label{eq:chap_6_7_10}
\rho_{N}(s)\leq\sigma_{N}\texttt{ when }|s|\leq|s_{1}|.
\end{equation}
Since $a_{N}$ are the remainders of a convergent series, they tend to zero as $N$ tends to infinity. That is, for each positive number $\varepsilon$, an integer $N_{\varepsilon}$ exists such that
\begin{equation}\label{eq:chap_6_7_11}
\sigma_{N}<\varepsilon\texttt{ whenever }N>N_{\varepsilon}.
\end{equation}
Because of conditions~(\ref{eq:chap_6_7_10}) and~(\ref{eq:chap_6_7_11}), then, condition~(\ref{eq:chap_6_7_4}) holds for all points $s$ in the sphere $|s|\leq|s_{1}|$; and the value of $N_{\varepsilon}$ is independent of the choice of $s$. Hence the convergence of series~(\ref{eq:chap_6_7_6}) is uniform in that sphere.

The extension of the proof to the case in which $s_{0}$ is arbitrary is, of course, accomplished by writing $\omega=s-s_{0}$ in series~(\ref{eq:chap_6_7_6}). For then the hypothesis of the theorem is that $s_{1}-s_{0}$ is a point inside the sphere of convergence $|\omega|=R$ of the series
\[
\sum_{n=0}^{\infty}a_{n}\omega^{n},
\]
Since we know that this series converges uniformly in the sphere $|\omega|\leq|s_{1}-s_{0}|$, the
conclusion in the statement of the theorem is evident.

\section{Continuity of Sums of Power Series}\label{sec:chap_6_8_58}

Our next theorem is an important consequence of uniform convergence, discussed in the previous section.

\begin{theorem}\label{th:chap_6_8_1}
A power series
\begin{equation}\label{eq:chap_6_8_1}
\sum_{n=0}^{\infty}a_{n}(s-s_{0})^{n},
\end{equation}
represents a continuous function $S(s)$ at each point inside its sphere of convergence $|s-s_{0}|=R$.
\end{theorem}

Another way to state this theorem is to say that if $S(s)$ denotes the sum of series~(\ref{eq:chap_6_8_1}) within its sphere of convergence $|s-s_{0}|=R$ and if $s_{1}$ is a point inside that sphere, then, for each positive number $\varepsilon$, there is a positive number $\delta$ such that
\begin{equation}\label{eq:chap_6_8_2}
|S(s)-S(s_{1})|<\varepsilon\texttt{ whenever }|s-s_{1}|<\delta.
\end{equation}
the number $\delta$ being small enough so that $s$ lies in the domain of definition $|s-s_{1}|<R$ of $S(s)$. [See definition~(\ref{eq:chap_2_7_4}), Sec.~(\ref{sec:chap_2_7_17}), of continuity.]

To show this, we let $S_{N}(s)$ denote the sum of the first $N$ terms of series~(\ref{eq:chap_6_8_1}) and write the remainder function
\[
\rho_{N}(s)=S(s)-S_{N}(s)\texttt{ }(|s-s_{0}|<R).
\]
Then, because
\[
S(s)=S_{N}(s)+\rho_{N}(s)\texttt{ }(|s-s_{0}|<R),
\]
one can see that
\[
|S(s)-S(s_{1})|=|S_{N}(s)-S_{N}(s_{1})+\rho_{N}(s)-\rho_{N}(s_{1})|,
\]
or
\begin{equation}\label{eq:chap_6_8_3}
|S(s)-S(s_{1})|\leq|S_{N}(s)-S_{N}(s_{1})|+|\rho_{N}(s)|+|\rho_{N}(s_{1})|.
\end{equation}
If $s$ is any point lying in some closed sphere $|s-s_{0}|\leq R_{0}$ whose radius $R_{0}$ is greater than $|s_{1}-s_{0}|$ but less than the radius $R$ of the sphere of convergence of series~(\ref{eq:chap_6_8_1}), the uniform convergence stated in Theorem~(\ref{th:chap_6_7_2}), Sec.~(\ref{sec:chap_6_7_57}), ensures that there is
a positive integer $N_{\varepsilon}$ such that
\begin{equation}\label{eq:chap_6_8_4}
\rho_{N}(s)<\varepsilon/3\texttt{ whenever }N>N_{\varepsilon}.
\end{equation}
In particular, condition~(\ref{eq:chap_6_8_4}) holds for each point $s$ in some neighborhood $|s-s_{1}|<\delta$ of  $s_{1}$ that is small enough to be contained in the sphere $|s-s_{1}|\leq R_{0}$.

Now the partial sum $S_{N}(s)$ is a polynomial and is, therefore, continuous at $s_{1}$ for each value of $N$. In  particular, when $N=N_{\varepsilon}+1$, we can choose our $\delta$ so small that
\begin{equation}\label{eq:chap_6_8_5}
|S_{N}(s)-S_{N}(s_{1})|<\varepsilon/3\texttt{ whenever }|s-s_{1}|<\delta.
\end{equation}
By writing $N=N_{\varepsilon}+1$ in inequality~(\ref{eq:chap_6_8_3}) and using the fact that statements~(\ref{eq:chap_6_8_4}) and~(\ref{eq:chap_6_8_5}) are true when $N=N_{\varepsilon}+1$, we now find that
\[
|S(s)-S(s_{1})|<\varepsilon/3+\varepsilon/3+\varepsilon/3=\varepsilon\texttt{ whenever }|s-s_{1}|<\delta.
\]
This is statement~(\ref{eq:chap_6_8_2}), and the corollary is now established.

By writing $\omega=1/(s-s_{0})$, one can modify the two theorems in the previous section and the theorem here so as to apply to series of the type
\begin{equation}\label{eq:chap_6_8_6}
\sum_{n=1}^{\infty}\frac{b_{n}}{(s-s_{0})^{n}}.
\end{equation}
If, for instance, series~(\ref{eq:chap_6_8_6}) converges at a point $s_{1}$ $(s_{1}\neq s_{0})$, the series
\[
\sum_{n=1}^{\infty}\frac{b_{n}}{\omega^{n}}
\]
must converge absolutely to a continuous function when
\begin{equation}\label{eq:chap_6_8_7}
|\omega|<\frac{1}{|s_{1}-s_{0}|}.
\end{equation}

Thus, since inequality~(\ref{eq:chap_6_8_7}) is the same as $|s-s_{0}|>|s_{1}-s_{0}|$. series~(\ref{eq:chap_6_8_6}) must converge absolutely to a continuous function in the domain exterior to the sphere $|s-s_{0}|=R_{1}$,
where $R_{1}=|s_{1}-s_{0}|$. Also, we know that if a Laurent series representation
\[
f(s)=\sum_{n=0}^{\infty}a_{n}(s-s_{0})^{n}+\sum_{n=1}^{\infty}\frac{b_{n}}{(s-s_{0})^{n}}
\]
is valid in an spherical shell $R_{1}<|s_{1}-s_{0}|<R_{2}$, then both of the series on the right converge
uniformly in any closed spherical shell which is concentric to and interior to that region of validity.

\section{Integration and Differentiation of Power Series}\label{sec:chap_6_9_59}

We have just seen that a power series
\begin{equation}\label{eq:chap_6_9_1}
S(s)=\sum_{n=0}^{\infty}a_{n}(s-s_{0})^{n}
\end{equation}
represents a continuous function at each point interior to its sphere of convergence. In this section, we prove that the sum $S(s)$ is actually analytic within that sphere. Our proof depends on the following theorem, which is of interest in itself.

\begin{theorem}\label{th:chap_6_9_1}
Let $C$ denote any contour interior to the sphere of convergence of the power series~(\ref{eq:chap_6_9_1}), which is on a closed sphere, and let $g(s)$ be any function that is continuous on $C$. The series formed by multiplying each term of the power series by $g(s)$ can be integrated term by term over $C$; that is,
\begin{equation}\label{eq:chap_6_9_2}
\int_{C}g(s)S(s)ds=\sum_{n=0}^{\infty}a_{n}\int_{C}g(s)(s-s_{0})^{n}ds.
\end{equation}
\end{theorem}

To prove this theorem, we note that since both $g(s)$ and the sum $S(s)$ of the power series are continuous on $C$, the integral over $C$ of the product
\[
g(s)S(s)=\sum_{n=0}^{N-1}a_{n}g(s)(s-s_{0})^{n}+g(s)\rho_{N}(s),
\]
where $\rho_{N}(s)$ is the remainder of the given series after $N$ terms, exists. The terms of the finite sum here are also continuous on the contour $C$, and so their integrals over $C$ exist. Consequently, the integral of the quantity $g(z)\rho_{N}(s)$ must exist; and we may write
\begin{equation}\label{eq:chap_6_9_3}
\int_{C}g(s)S(s)ds=\sum_{n=0}^{N-1}a_{n}\int_{C}g(s)(s-s_{0})^{n}ds+\int_{C}g(s)\rho_{N}(s)ds.
\end{equation}

Now let $M$ be the maximum value of $|g(s)|$ on $C$, and let $L$ denote the length of $C$. In view of the uniform convergence of the given power series (Sec.~(\ref{sec:chap_6_7_57})), we know that for each positive number $\varepsilon$ there exists a positive integer $N_{\varepsilon}$ such that, for all points $s$ on $C$,
\[
|\rho_{N}(s)|<\varepsilon\texttt{ whenever }N>N_{\varepsilon}.
\]
Since $N_{\varepsilon}$ is independent of $s$, we find that
\[
|\int_{C}\rho_{N}(s)ds|<M\varepsilon L\texttt{ whenever }N>N_{\varepsilon},
\]
that is,
\[
\lim_{N \to \infty}\int_{C}\rho_{N}(s)ds=0.
\]
It follows, therefore, from equation~(\ref{eq:chap_6_9_3}) that
\[
\int_{C}g(s)S(s)ds=\lim_{N \to \infty}\sum_{n=0}^{N-1}a_{n}\int_{C}g(s)(s-s_{0})^{n}ds.
\]
This is the same as equation~(\ref{eq:chap_6_9_2}), and Theorem~(\ref{th:chap_6_9_1}) is proved.

If $|g(s)|=1$ for each value of $s$ in the open sphere bounded by the sphere of convergence of power series~(\ref{eq:chap_6_9_1}), the fact that $(s-s_{0})^{n}$ is entire when $n=0,1,2,\ldots$ ensures that
\[
\int_{C}g(s)(s-s_{0})^{n}ds=\int_{C}(s-s_{0})^{n}ds=0\texttt{ }(n=0,1,2,\ldots)
\]
for every closed contour $C$ lying in that domain. According to equation~(\ref{eq:chap_6_9_2}), then,
\[
\int_{C}S(s)ds=0
\]
for every such contour; and, by Morera's theorem (Sec.~(\ref{sec:chap_4_13_48})), the function $S(s)$ is analytic throughout the domain. We state this result as a corollary.

\begin{corollary}\label{co:chap_6_9_1}
The sum $S(s)$ of power series~(\ref{eq:chap_6_9_1}) is analytic at each point $s$ interior to the sphere of convergence of that series.
\end{corollary}

This corollary is often helpful in establishing the analyticity of functions and in evaluating limits.

We observed that the Taylor series for a function $f$ about a point $s_{0}$ converges to $f(s)$ at each point $s$ interior to the sphere centered at $s_{0}$ and passing through the nearest point $s_{1}$ where $f$ fails to be analytic. In view of the above corollary, we now know that there is no larger sphere about $s_{0}$ such that at each point $s$ interior to it the Taylor series converges to $f(s)$. For if there were such a sphere, $f$ would be analytic at $s_{1}$; but $f$ is not analytic at $s_{1}$.

We now present a companion to Theorem~(\ref{th:chap_6_9_1}).

\begin{theorem}\label{th:chap_6_9_2}
The power series~(\ref{eq:chap_6_9_1}) can be differentiated term by term. That is, at each point $s$ interior to the sphere of convergence of that series,
\begin{equation}\label{eq:chap_6_9_6}
S'(s)=\sum_{n=1}^{\infty}na_{n}(s-s_{0})^{n-1}.
\end{equation}
\end{theorem}

To prove this, let $s$ denote any point interior to the sphere of convergence of series~(\ref{eq:chap_6_9_1}), which is on a closed sphere; and let $C$ be some positively oriented simple closed contour surrounding $s$ and interior to that sphere. Also, define the function
\begin{equation}\label{eq:chap_6_9_7}
g(\varsigma)=\frac{1}{2\pi i}\cdot\frac{1}{(\varsigma-s)^{2}}
\end{equation}
at each point $\varsigma$ on $C$. Since $g(\varsigma)$ is continuous on $C$, Theorem~(\ref{th:chap_6_9_1}) tells us that
\begin{equation}\label{eq:chap_6_9_8}
\int_{C}g(\varsigma)S(\varsigma)d\varsigma=\sum_{n=0}^{\infty}a_{n}\int_{C}g(\varsigma)(\varsigma-s_{0})^{n}d\varsigma.
\end{equation}

Now $S(s)$ is analytic inside and on $C$,  and this enables us to write
\[
\int_{C}g(\varsigma)S(\varsigma)d\varsigma=\frac{1}{2\pi i}\int_{C}\frac{S(\varsigma)d\varsigma}{(\varsigma-s)^{2}}=S'(s)
\]
with the aid of the integral representation for derivatives in Sec.~(\ref{sec:chap_4_13_48}). Furthermore,
\[
\int_{C}g(\varsigma)(\varsigma-s_{0})^{n}d\varsigma=\frac{1}{2\pi i}\int_{C}\frac{(\varsigma-s_{0})^{n}}{(\varsigma-s)^{n}}d\varsigma
=\frac{d}{ds}(s-s_{0})^{n}\texttt{ }(n=0,1,2,\ldots).
\]
Thus equation~(\ref{eq:chap_6_9_8}) reduces to
\[
S'(s)=\sum_{n=0}^{\infty}a_{n}\frac{d}{ds}(s-s_{0})^{n},
\]
which is the same as equation~(\ref{eq:chap_6_9_6}). This completes the proof.

\section{Uniqueness of Series Representations}\label{sec:chap_6_10_60}

The uniqueness of Taylor and Laurent series representations follows readily from Theorem~(\ref{th:chap_6_9_1}) in Sec.~(\ref{sec:chap_6_9_59}). We consider first the uniqueness of Taylor series representations.

\begin{theorem}\label{th:chap_6_10_1}
If a series
\begin{equation}\label{eq:chap_6_10_1}
\sum_{n=0}^{\infty}a_{n}(s-s_{0})^{n}
\end{equation}
converges to $f(s)$ at all points interior to some sphere on a closed sphere $|s-s_{0}|=R$, then it is the Taylor series expansion  for $f$ in powers of $s-s_{0}$.
\end{theorem}

To prove this, we write the series representation
\begin{equation}\label{eq:chap_6_10_2}
f(s)=\sum_{n=0}^{\infty}a_{n}(s-s_{0})^{n}\texttt{ }(|s-s_{0}|<R)
\end{equation}
in the hypothesis of the theorem using the index of summation $m$:
\[
f(s)=\sum_{m=0}^{\infty}a_{m}(s-s_{0})^{m}\texttt{ }(|s-s_{0}|<R).
\]
Then, by appealing to Theorem~(\ref{th:chap_6_9_1}) in Sec.~(\ref{sec:chap_6_9_59}), we may write
\begin{equation}\label{eq:chap_6_10_3}
\int_{C}g(s)f(s)ds=\sum_{m=0}^{\infty}a_{m}\int_{C}g(s)(s-s_{0})^{m}ds,
\end{equation}
where $g(s)$ is any one of the functions
\begin{equation}\label{eq:chap_6_10_4}
g(s)=\frac{1}{2\pi i}\cdot\frac{1}{(s-s_{0})^{n+1}}\texttt{ }(n=0,1,2,\ldots)
\end{equation}
and $C$ is some spatial circle centered at $s_{0}$ and with radius less than $R$.

In view of the generalized form~(\ref{eq:chap_4_13_5}), Sec.~(\ref{sec:chap_4_13_48}), of the Cauchy integral formula (see also the corollary in Sec.~(\ref{sec:chap_6_9_59})), we find that
\begin{equation}\label{eq:chap_6_10_5}
\int_{C}g(s)f(s)ds=\frac{1}{2\pi i}\int_{C}\frac{f(s)ds}{(s-s_{0})^{n+1}}=\frac{f^{(n)}(s_{0})}{n!};
\end{equation}
and, since
\begin{equation}\label{eq:chap_6_10_6}
\int_{C}g(s)(s-s_{0})^{m}ds=\frac{1}{2\pi i}\int_{C}\frac{ds}{(s-s_{0})^{n-m+1}}=\{
    \begin{array}{cc}
        0 & \texttt{ when }m\neq n, \\
        1 & \texttt{ when }m=n,
    \end{array}
\end{equation}
it is clear that
\begin{equation}\label{eq:chap_6_10_7}
\sum_{m=0}^{\infty}a_{m}\int_{C}g(s)(s-s_{0})^{m}ds=a_{n}.
\end{equation}
Because of equations~(\ref{eq:chap_6_10_5}) and~(\ref{eq:chap_6_10_7}), equation~(\ref{eq:chap_6_10_3}) now reduces to
\[
\frac{f^{(n)}(s_{0})}{n!}=a_{n},
\]
and this shows that series~(\ref{eq:chap_6_10_2}) is, in fact, the Taylor series for $f$ about the point $s_{0}$.

Note how it follows from Theorem~(\ref{th:chap_6_10_1}) that if series~(\ref{eq:chap_6_10_1}) converges to zero throughout some neighborhood of $s_{0}$, then the coefficients $a_{n}$ must all be zero.

Our second theorem here concerns the uniqueness of Laurent series representations.

\begin{theorem}\label{th:chap_6_10_2}
If a series
\begin{equation}\label{eq:chap_6_10_8}
\sum_{n=-\infty}^{\infty}c_{n}(s-s_{0})^{n}=\sum_{n=0}^{\infty}a_{n}(s-s_{0})^{n}+\sum_{n=1}^{\infty}\frac{b_{n}}{(s-s_{0})^{n}}
\end{equation}
converges to $f(s)$ at all points in some spherical shell domain about $s_{0}$ and contained between two closed spheres, then it is the Laurent series expansion for $f$ in powers of $s-s_{0}$ for that domain.
\end{theorem}

The method of proof here is similar to the one used in proving Theorem~(\ref{th:chap_6_10_1}). The
hypothesis of this theorem tells us that there is a spherical shell domain about $s_{0}$ such that
\[
f(s)=\sum_{n=-\infty}^{\infty}c_{n}(s-s_{0})^{n}
\]
for each point $s$ in it. Let $g(s)$ be as defined by equation~(\ref{eq:chap_6_10_4}), but now allow $n$ to be
a negative integer too. Also, let $C$ be any spatial circle around the annulus, centered at $s_{0}$ and taken in the positive sense. Then, using the index of summation $m$ and adapting Theorem~(\ref{th:chap_6_9_1}) in Sec.~(\ref{sec:chap_6_9_59}) to series involving both nonnegative and negative powers of
$s-s_{0}$, write
\[
\int_{C}g(s)f(s)ds=\sum_{m=-\infty}^{\infty}c_{m}\int_{C}g(s)(s-s_{0})^{m}ds,
\]
or
\begin{equation}\label{eq:chap_6_10_9}
\frac{1}{2\pi i}\int_{C}\frac{f(s)ds}{(s-s_{0})^{n+1}}=\sum_{m=-\infty}^{\infty}c_{m}\int_{C}g(s)(s-s_{0})^{m}ds.
\end{equation}
Since equations~(\ref{eq:chap_6_10_6}) are also valid when the integers $m$ and $n$ are allowed to be
negative, equation~(\ref{eq:chap_6_10_9}) reduces to
\[
\frac{1}{2\pi i}\int_{C}\frac{f(s)ds}{(s-s_{0})^{n+1}}=c_{n},
\]
which is expression~(\ref{eq:chap_6_5_5}), Sec.~(\ref{sec:chap_6_5_55}), for coefficients in the Laurent series for  $f$ in the spherical shell.

\section{Multiplication and Division of Power Series}\label{sec:chap_6_11_61}

Suppose that each of the power series
\begin{equation}\label{eq:chap_6_11_1}
\sum_{n=0}^{\infty}a_{n}(s-s_{0})^{n}\texttt{ and }\sum_{n=0}^{\infty}b_{n}(s-s_{0})^{n}
\end{equation}
converges within some sphere $|s-s_{0}|=R$. Their sums $f(s)$ and $g(s)$, respectively, are then analytic functions in the sphere $|s-s_{0}|<R$ (Sec.~(\ref{sec:chap_6_9_59})), and the product of those sums has a Taylor series expansion which is valid there:
\begin{equation}\label{eq:chap_6_11_2}
f(s)g(s)=\sum_{n=0}^{\infty}c_{n}(s-s_{0})^{n}\texttt{ }(|s-s_{0}|<R).
\end{equation}

According to Theorem~(\ref{th:chap_6_10_1}) in Sec.~(\ref{sec:chap_6_10_60}), the series~(\ref{eq:chap_6_11_1}) are themselves Taylor series. Hence the first three coefficients in series~(\ref{eq:chap_6_11_2}) are given by the equations
\[
c_{0}=f(s_{0})g(s_{0})=a_{0}b_{0},
\]
\[
c_{1}=\frac{f(s_{0})g'(s_{0})+f'(s_{0})g(s_{0})}{1!}=a_{0}b_{1}+a_{1}b_{0},
\]
and
\[
c_{2}=\frac{f(s_{0})g''(s_{0})+2f'(s_{0})g'(s_{0})+f''(s_{0})g(s_{0})}{2!}=a_{0}b_{2}+a_{1}b_{1}+a_{2}b_{0}.
\]
The general expression for any coefficient $c_{n}$ is easily obtained by referring to Leibniz's rule
\begin{equation}\label{eq:chap_6_11_3}
(f(s)g(s))^{n}=\sum_{k=0}^{n}\{
    \begin{array}{c}
        n \\
        k
    \end{array})
f^{(k)}(s)g^{(n-k)}(s)
\end{equation}
where
\[
    (\begin{array}{c}
        n \\
        k
    \end{array}
)=\frac{n!}{k!(n-k)!}\texttt{ }(k=1,2,\ldots,n),
\]
for the $n$th derivative of the product of two differentiable functions. As usual, $f^{(0)}(s)=f(s)$ and $0!=1$. Evidently,
\[
c_{n}=\sum_{k=0}^{n}\frac{f^{(k)}(s)}{k!}\frac{g^{(n-k)}(s)}{(n-k)!}=\sum_{k=0}^{n}a_{k}b_{n-k};
\]
and so expansion~(\ref{eq:chap_6_11_2}) can be written
\begin{equation}\label{eq:chap_6_11_4}
f(s)g(s)=\sum_{n=0}^{\infty}c_{n}(s-s_{0})^{n}=a_{0}b_{0}+(a_{0}b_{1}+a_{1}b_{0})(s-s_{0})
\end{equation}
\[
+(a_{0}b_{2}+a_{1}b_{1}+a_{2}b_{0})(s-s_{0})^{2}+\cdots
\]
\[
+(\sum_{k=0}^{n}a_{k}b_{n-k})(s-s_{0})^{n}+\cdots\texttt{ }(|s-s_{0}|<R).
\]

Series~(\ref{eq:chap_6_11_4}) is the same as the series obtained by formally multiplying the two series~(\ref{eq:chap_6_11_1}) term by term and collecting the resulting terms in like powers  of $s-s_{0}$ it is called the Cauchy product of the two given series.

Continuing to let $f(s)$ and $g(s)$ denote the sums of series~(\ref{eq:chap_6_11_1}), suppose that $g(s)\neq0$ when $|s-s_{0}|<R$. Since the quotient $f(s)/g(s)$ is analytic throughout the sphere $|s-s_{0}|<R$, it has a Taylor series representation
\begin{equation}\label{eq:chap_6_11_6}
\frac{f(s)}{g(s)}=\sum_{n=0}^{\infty}d_{n}(s-s_{0})^{n}\texttt{ }(|s-s_{0}|<R),
\end{equation}
where the coefficients $d_{n}$ can be found by differentiating $f(s)/g(s)$ successively and evaluating the derivatives at $s=s_{0}$. The results are the same as those found by formally carrying out the division of the first of series~(\ref{eq:chap_6_11_1}) by the second. Since it is usually only the first few terms that are needed in practice, this method is not difficult.

%-----------------------------------------------------------------------
% Beginning of chap7.tex
%-----------------------------------------------------------------------
%
% AMS-LaTeX 1.2 sample file for a monograph, based on amsbook.cls.
% This is a data file input by chapter.tex.
%%%%%%%%%%%%%%%%%%%%%%%%%%%%%%%%%%%%%%%%%%%%%%%%%%%%%%%%%%%%%%%%%%%%

%\part{This is a Part Title Sample}

\chapter{Residues and Poles}\label{ch:chap_7}

The Cauchy-Goursat theorem (Sec.~(\ref{sec:chap_4_10_45})) states that if a function is analytic at all points interior to and on a simple closed contour $C$, then the value of the integral of the function around that contour is zero. If, however, the function fails to be analytic at a finite number of points interior to $C$, there is, as we shall see in this chapter, a specific number, called a residue, which each of those points contributes to the value of the integral. We develop here the theory of residues.%; and, in Chap.~(\ref{ch:chap_7}), we shall illustrate their use in certain areas of applied mathematics.

\section{Residues}\label{sec:chap_7_1_62}

Recall (Sec.~(\ref{sec:chap_2_13_23})) that a point $s_{0}$ is called a singular point of a function $f$ if $f$ fails to be analytic at $s_{0}$ but is analytic at some point in every neighborhood of $s_{0}$. A singular point
$s_{0}$ is said to be isolated if, in addition, there is a deleted neighborhood $0<|s-s_{0}|<\varepsilon$
of $s_{0}$ throughout which $f$ is analytic.

When $s_{0}$ is an isolated singular point of a function $f$, there is a positive number
$R_{2}$ such that $f$ is analytic at each point $s$ for which $0<|s-s_{0}|<R_{2}$. Consequently, $f(s)$ is represented by a Laurent series
\begin{equation}\label{eq:chap_7_1_1}
f(s)=\sum_{n=0}^{\infty}a_{n}(s-s_{0})^{n}+\frac{b_{1}}{s-s_{0}}+\frac{b_{2}}{(s-s_{0})^{2}}+\cdots
\end{equation}
\[
+\frac{b_{n}}{(s-s_{0})^{n}}+\cdots\texttt{ }(R_{1}<|s-s_{0}|<R_{2}).
\]

\subsection{Residues with Contour Integrals}\label{sec:chap_7_1_62_1}

From Sec.~(\ref{sec:chap_6_5_55}), the coefficients $a_{n}$ and $b_{n}$ in expansion~(\ref{eq:chap_7_1_1}) of a Laurent series have certain integral representations. In particular,
\[
b_{n}=\frac{1}{2\pi i}\int_{C}\frac{f(s)ds}{(s-s_{0})^{-n+1}}\texttt{ }(n=1,2,\ldots)
\]
where $C$ is any positively oriented simple closed contour in a closed sphere $S$ around $s_{0}$ and lying in the punctured disk $0<|s-s_{0}|<R_{2}$. When $n=1$, this expression for $b_{n}$ can be written
\begin{equation}\label{eq:chap_7_1_2}
\int_{C}f(s)ds=2\pi ib_{1}.
\end{equation}
The complex number $b_{1}$ which is the coefficient of $1/(s-s_{0})$ in expansion~(\ref{eq:chap_7_1_1}), is
called the residue of $f$ at the isolated singular point $s_{0}$. We shall often use the notation
\[
Res_{s=s_{0}}f(s),
\]
or simply $B$ when the point $s_{0}$ and the function $f$ are clearly indicated, to denote the residue $b_{1}$.

Equation~(\ref{eq:chap_7_1_2}) provides a powerful method for evaluating certain integrals around simple closed contours.

\subsection{Residues with Surface Integrals}\label{sec:chap_7_1_62_2}

From Sec.~(\ref{sec:chap_6_6_56}), the coefficients $a_{n}$ and $b_{n}$ in expansion~(\ref{eq:chap_7_1_1}) of a Laurent series have certain integral representations. In particular,
\[
b_{n}=\frac{1}{4\pi}\iint_{S}\frac{f(s)d\sigma}{(s-s_{0})^{-n+2}}\texttt{ }(n=1,2,\ldots)
\]
where $S$ is any positively oriented simple closed surface in a closed sphere $S_{R}$ around $s_{0}$ and lying in the punctured sphere $0<|s-s_{0}|<R_{2}$. When $n=2$, this expression for $b_{n}$ can be written
\begin{equation}\label{eq:chap_7_1_3}
\iint_{S}f(s)d\sigma=4\pi b_{2}.
\end{equation}
The complex number $b_{2}$ which is the coefficient of $1/(s-s_{0})^{2}$ in expansion~(\ref{eq:chap_7_1_1}), is
called the residue of $f$ at the isolated singular point $s_{0}$. We shall often use the notation
\[
Res_{s=s_{0}}f(s),
\]
or simply $B$ when the point $s_{0}$ and the function $f$ are clearly indicated, to denote the residue $b_{2}$.

Equation~(\ref{eq:chap_7_1_3}) provides a powerful method for evaluating certain integrals around simple closed contours.

\section{Cauchy's Residue Theorem}\label{sec:chap_7_2_63}

\subsection{Cauchy's Residue Theorem with Contour Integrals}\label{sec:chap_7_2_63_1}

If, except for a finite number of singular points, a function $f$ is analytic inside a closed sphere $S$ containing a simple closed contour $C$, those singular points must be isolated (Sec.~(\ref{sec:chap_7_1_62_1})). The following theorem, which is known as Cauchy's residue theorem, is a precise statement of the fact that if $f$ is also analytic in $S$ and $C$ and if $C$ is positively oriented, then the value of the integral of $f$ around $C$ is $2\pi i$ times the sum of the residues of $f$ at the singular points inside $S$ and $C$.

\begin{theorem}\label{th:chap_7_2_1}
Let $C$ be a simple closed contour in a closed sphere $S$, described in the positive sense. If a function $f$ is analytic inside $S$ and on $C$ except for a finite number of singular points
$s_{k}$ $(k=1,2,\ldots,n)$ inside $S$ and $C$, then
\begin{equation}\label{eq:chap_7_2_1}
\int_{C}f(s)ds=2\pi i\sum_{k=1}^{n}Res_{s=s_{k}}f(s).
\end{equation}
\end{theorem}

To prove the theorem, let the points $s_{k}$ $(k=1,2,\ldots,n)$ be centers of positively oriented circles $C_{k}$ on spheres $S_{k}$ which are interior to $S$ and are so small that no two of them have points in common. All of these circles $C_{k}$ together with the spherical surface of $S$ containing the simple closed contour $C$, form the boundary of a closed region after deleting the interiors of these circles $C_{k}$, throughout which $f$ is analytic and whose interior is a multiply connected domain. Hence, according to the extension of the Cauchy-Goursat theorem to such regions (Theorem~(\ref{th:chap_4_11_2}), Sec.~(\ref{sec:chap_4_11_46})),
\[
\int_{C}f(s)ds-\sum_{k=1}^{n}\int_{C_{k}}f(s)ds=0.
\]
This reduces to equation~(\ref{eq:chap_7_2_1}) because (Sec.~(\ref{sec:chap_7_1_62_1}))
\[
\int_{C_{k}}f(s)ds=2\pi i Res_{s=s_{k}}f(s)\texttt{ }(k=1,2,\ldots,n),
\]
and the proof is complete.

\subsection{Cauchy's Residue Theorem with Surface Integrals}\label{sec:chap_7_2_63_2}

If, except for a finite number of singular points, a function $f$ is analytic inside a closed sphere $S_{R}$ containing a simple closed surface $S$, those singular points must be isolated (Sec.~(\ref{sec:chap_7_1_62_2})). The following theorem, which is known as Cauchy's residue theorem, is a precise statement of the fact that if $f$ is also analytic in $S$ and if $S$ is positively oriented, then the value of the integral of $f$ on $S$ is $4\pi$ times the sum of the residues of $f$ at the singular points inside $S$.

\begin{theorem}\label{th:chap_7_2_2_1}
Let $S$ be a simple closed surface in a closed sphere $S_{R}$, described in the positive sense. If a function $f$ is analytic inside $S_{R}$ and on $S$ except for a finite number of singular points
$s_{k}$ $(k=1,2,\ldots,n)$ inside $S$, then
\begin{equation}\label{eq:chap_7_2_2_1}
\iint_{S}f(s)d\sigma=4\pi\sum_{k=1}^{n}Res_{s=s_{k}}f(s).
\end{equation}
\end{theorem}

To prove the theorem, let the points $s_{k}$ $(k=1,2,\ldots,n)$ be centers of positively oriented spheres $S_{k}$ which are interior to $S$ and are so small that no two of them have points in common. All of these spheres $S_{k}$ together with the spherical surface of $S$, form the boundary of a closed region after deleting the interiors of these spheres $S_{k}$, throughout which $f$ is analytic and whose interior is a multiply connected domain. Hence, according to the extension of the Cauchy-Goursat theorem to such regions (Theorem~(\ref{th:chap_5_4_4}), Sec.~(\ref{sec:chap_5_4_54})),
\[
\iint_{S}f(s)d\sigma-\sum_{k=1}^{n}\iint_{S_{k}}f(s)d\sigma=0.
\]
This reduces to equation~(\ref{eq:chap_7_2_2_1}) because (Sec.~(\ref{sec:chap_7_1_62_2}))
\[
\iint_{S_{k}}f(s)d\sigma=4\pi Res_{s=s_{k}}f(s)\texttt{ }(k=1,2,\ldots,n),
\]
and the proof is complete.

\section{Using a Single Residue}\label{sec:chap_7_3_64}

\subsection{Using a Single Residue with Contour Integrals}\label{sec:chap_7_3_64_1}

If the function $f$ in Cauchy's residue theorem (Sec.~(\ref{sec:chap_7_2_63})) is, in addition, analytic at each
point in the finite space interior to $C$, it is sometimes more efficient to evaluate the integral of $f$ around $C$ by finding a single residue of a certain related function. We present the method as a theorem.

\begin{theorem}\label{th:chap_7_3_1}
If a function $f$ is analytic everywhere in the finite space except for a finite number of singular points interior to a positively oriented simple closed contour $C$ in a closed sphere $S$, then
\begin{equation}\label{eq:chap_7_3_1}
\int_{C}f(s)ds=2\pi iRes_{s=0}(\frac{1}{s^{2}}f(\frac{1}{s})).
\end{equation}
\end{theorem}

We begin the derivation of expression~(\ref{eq:chap_7_3_1}) by constructing a closed sphere $S$ $(|s|=R_{1})$ which is large enough so that the contour $C$ is interior to it. Then if $C_{0}$ denotes a positively oriented circle on a closed sphere $S_{0}$ $(|s|=R_{0}))$, where $R_{0}>R_{1}$, we know from Laurent's theorem
(Sec.~(\ref{sec:chap_6_5_55})) that
\begin{equation}\label{eq:chap_7_3_2}
f(s)=\sum_{n=-\infty}^{\infty}c_{n}s^{n}\texttt{ }(R_{1}<|s|<\infty)
\end{equation}
where
\begin{equation}\label{eq:chap_7_3_3}
c_{n}=\frac{1}{2\pi i}\int_{C}\frac{f(s)ds}{s^{n+1}}\texttt{ }(n=0,\pm1,\pm2,\ldots).
\end{equation}
By writing $n=-1$ in expression~(\ref{eq:chap_7_3_3}), we find that
\begin{equation}\label{eq:chap_7_3_4}
\int_{C_{0}}f(s)ds=2\pi ic_{-1}.
\end{equation}
Observe that, since the condition of validity with representation~(\ref{eq:chap_7_3_2}) is not of the type
$0<|s|<R_{2}$, the coefficient $c_{-1}$ is not the residue of $f$ at the point $s=0$, which may not even be a singular point of $f$. But, if we replace $s$ by $1/s$ in representation~(\ref{eq:chap_7_3_2}) and its condition  of validity, we see that
\[
\frac{1}{s^{2}}f(\frac{1}{s})=\sum_{n=-\infty}^{\infty}\frac{c_{n}}{s^{n+2}}=\sum_{n=-\infty}^{\infty}\frac{c_{n-2}}{s^{n}}\texttt{ }(0<|s|<\frac{1}{R_{1}})
\]
and hence that
\begin{equation}\label{eq:chap_7_3_5}
c_{-1}=Res_{s=0}(\frac{1}{s^{2}}f(\frac{1}{s})).
\end{equation}
Then, in view of equations~(\ref{eq:chap_7_3_4}) and~(\ref{eq:chap_7_3_5}),
\[
\int_{C_{0}}f(s)ds=2\pi iRes_{s=0}(\frac{1}{s^{2}}f(\frac{1}{s})).
\]
Finally, since $f$ is analytic throughout the closed spherical shell bounded by spherical surfaces of $S$ and $S_{0}$, the principle of deformation of paths (Corollary~(\ref{co:chap_4_11_2}), Sec.~(\ref{sec:chap_4_11_46})) yields the desired result~(\ref{eq:chap_7_3_1}).

\subsection{Using a Single Residue with Surface Integrals}\label{sec:chap_7_3_64_2}

If the function $f$ in Cauchy's residue theorem (Sec.~(\ref{sec:chap_7_2_63_2})) is, in addition, analytic at each
point in the finite space interior to $S$, it is sometimes more efficient to evaluate the integral of $f$ around $S$ by finding a single residue of a certain related function. We present the method as a theorem.

\begin{theorem}\label{th:chap_7_3_2_1}
If a function $f$ is analytic everywhere in the finite space except for a finite number of singular points interior to a positively oriented simple closed surface $S$ in a closed sphere $S_{R}$, then
\begin{equation}\label{eq:chap_7_3_2_1}
\iint_{S}f(s)d\sigma=4\pi Res_{s=0}(\frac{1}{s^{4}}f(\frac{1}{s})).
\end{equation}
\end{theorem}

We begin the derivation of expression~(\ref{eq:chap_7_3_2_1}) by constructing a closed sphere $S_{R}$ $(|s|=R_{1})$ which is large enough so that the surface $S$ is interior to it. Then if $S_{0}$ denotes a positively oriented and closed sphere $|s|=R_{0}$ $(R_{0}>R_{1})$, we know from Laurent's theorem
(Sec.~(\ref{sec:chap_6_6_56})) that
\begin{equation}\label{eq:chap_7_3_2_2}
f(s)=\sum_{n=-\infty}^{\infty}c_{n}s^{n}\texttt{ }(R_{1}<|s|<\infty)
\end{equation}
where
\begin{equation}\label{eq:chap_7_3_2_3}
c_{n}=\frac{1}{4\pi}\iint_{S}\frac{f(s)d\sigma}{s^{n+2}}\texttt{ }(n=0,\pm1,\pm2,\ldots).
\end{equation}
By writing $n=-2$ in expression~(\ref{eq:chap_7_3_2_3}), we find that
\begin{equation}\label{eq:chap_7_3_2_4}
\iint_{S_{0}}f(s)d\sigma=4\pi c_{-2}.
\end{equation}
Observe that, since the condition of validity with representation~(\ref{eq:chap_7_3_2_2}) is not of the type
$0<|s|<R_{2}$, the coefficient $c_{-2}$ is not the residue of $f$ at the point $s=0$, which may not even be a singular point of $f$. But, if we replace $s$ by $1/s$ in representation~(\ref{eq:chap_7_3_2_2}) and its condition  of validity, we see that
\[
\frac{1}{s^{4}}f(\frac{1}{s})=\sum_{n=-\infty}^{\infty}\frac{c_{n}}{s^{n+4}}=\sum_{n=-\infty}^{\infty}\frac{c_{n-4}}{s^{n}}\texttt{ }(0<|s|<\frac{1}{R_{1}})
\]
and hence that
\begin{equation}\label{eq:chap_7_3_2_5}
c_{-2}=Res_{s=0}(\frac{1}{s^{4}}f(\frac{1}{s})).
\end{equation}
Then, in view of equations~(\ref{eq:chap_7_3_2_4}) and~(\ref{eq:chap_7_3_2_5}),
\[
\iint_{S_{0}}f(s)d\sigma=4\pi Res_{s=0}(\frac{1}{s^{4}}f(\frac{1}{s})).
\]
Finally, since $f$ is analytic throughout the closed spherical shell bounded by spherical surfaces of $S$ and $S_{0}$, the principle of deformation of paths (Theorem~(\ref{th:chap_5_4_4}), Sec.~(\ref{sec:chap_5_4_54})) yields the desired result~(\ref{eq:chap_7_3_2_1}).

\section{The Three Types of Isolated Singular Points}\label{sec:chap_7_4_65}

We saw in Sec.~(\ref{sec:chap_7_1_62}) that the theory of residues is based on the fact that if $f$ has an isolated singular point $s_{0}$, then  $f(s)$ can be represented by a Laurent series
\begin{equation}\label{eq:chap_7_4_1}
f(s)=\sum_{n=0}^{\infty}a_{n}(s-s_{0})^{n}+\frac{b_{1}}{s-s_{0}}+\frac{b_{2}}{(s-s_{0})^{2}}+\cdots
+\frac{b_{n}}{(s-s_{0})^{n}}+\cdots
\end{equation}
in a punctured disk $0<|s-s_{0}|<R_{2}$. The portion
\[
\frac{b_{1}}{s-s_{0}}+\frac{b_{2}}{(s-s_{0})^{2}}+\cdots+\frac{b_{n}}{(s-s_{0})^{n}}+\cdots
\]
of the series, involving negative powers of $s-s_{0}$, is called the principal part of $f$ at $s_{0}$. We now use the principal part to identify the isolated singular point $s_{0}$ as one of three special types. This classification will aid us in the development of residue theory that appears in following sections.

If the principal part of $f$ at $s_{0}$ contains at least one nonzero term but the number of such terms is finite, then there exists a positive integer $m$ such that
\[
b_{m}\neq0\texttt{ and  }b_{m+1}=b_{m+2}=\cdots=0.
\]
That is, expansion~(\ref{eq:chap_7_4_1}) takes the form
\begin{equation}\label{eq:chap_7_4_2}
f(s)=\sum_{n=0}^{\infty}a_{n}(s-s_{0})^{n}+\frac{b_{1}}{s-s_{0}}+\frac{b_{2}}{(s-s_{0})^{2}}+\cdots
+\frac{b_{m}}{(s-s_{0})^{m}}
\end{equation}
in a punctured disk $0<|s-s_{0}|<R_{2}$, where $b_{m}\neq0$. In this case, the isolated singular point $s_{0}$ is called a pole of order $m$. A pole of  order $m=1$ is usually referred to as a simple pole.

There remain two extremes, the case in which all of the coefficients in the principal part are zero and the one in which an infinite number of them are nonzero.

When all of the $b_{n}$'s are zero, so that
\begin{equation}\label{eq:chap_7_4_3}
\begin{split}
f(s)&=\sum_{n=0}^{\infty}a_{n}(s-s_{0})^{n} \\
    &=a_{0}+a_{1}(s-s_{0})+a_{2}(s-s_{0})^{2}+\cdots\texttt{ }(0<|s-s_{0}|<R_{2}),
\end{split}
\end{equation}
the point $s_{0}$ is known as a removable singular point. Note that the residue at a removable singular point is always zero. If we define, or possibly redefine, $f$ at $s_{0}$ so that $f(s_{0})=a_{0}$, expansion~(\ref{eq:chap_7_4_3}) becomes valid throughout the entire disk $|s-s_{0}|<R_{2}$. Since a power series always represents an analytic function interior to its sphere of convergence (Sec.~(\ref{sec:chap_6_9_59})), it follows that $f$ is analytic at $s_{0}$ when it is assigned the value $a_{0}$ there. The singularity at $s_{0}$ is, therefore, removed.

When an infinite number of the coefficients $b_{n}$ in the principal part are nonzero, $s_{0}$ is said to be an  essential singular point of $f$. An important result concerning the behavior of a function near an essential singular point is due to Picard. It states that in each neighborhood of an essential singular point, a function assumes every finite value, with one possible exception, an infinite number of times.

In the remaining sections of this chapter, we shall develop in greater depth the theory of the three types of  isolated singular points just described. The emphasis will be on useful and efficient methods for identifying poles  and finding the corresponding residues.

\section{Residues at Poles}\label{sec:chap_7_5_66}

When a function $f$ has an isolated singularity at a point $s_{0}$, the basic method for identifying $s_{0}$ as a pole and finding the residue there is to write the appropriate Laurent series and to note the coefficient of $1/(s-s_{0})$. The following theorem provides an alternative characterization of poles and another way of finding the corresponding residues.

\begin{theorem}\label{th:chap_7_5_1}
An isolated singular point $s_{0}$ of a function $f$ is a pole of order $m$ if and only if $f(s)$ can be written  in the form
\begin{equation}\label{eq:chap_7_5_1}
f(s)=\frac{\phi(s)}{(s-s_{0})^{m}},
\end{equation}
where $\phi(s)$ is analytic and nonzero at $s_{0}$. Moreover,
\begin{equation}\label{eq:chap_7_5_2}
Res_{s=s_{0}}f(s)=\phi(s_{0})\texttt{ if }m=1,
\end{equation}
and
\begin{equation}\label{eq:chap_7_5_3}
Res_{s=s_{0}}f(s)=\frac{\phi^{(m-1)}(s_{0})}{(m-1)!}\texttt{ if }m\geq2.
\end{equation}
\end{theorem}

Observe that expression~(\ref{eq:chap_7_5_2}) need not have been written separately since, with the
convention that $\phi^{(0)}(s_{0})=\phi(s_{0})$ and $0!=1$, expression~(\ref{eq:chap_7_5_3}) reduces to it when $m=1$.

To prove the theorem, we first assume that $f(s)$ has the form~(\ref{eq:chap_7_5_1}) and recall (Sec.~(\ref{sec:chap_6_3_53})) that since $\phi(s)$ is analytic at $s_{0}$, it has a Taylor series representation
\[
\phi(s)=\phi(s_{0})+\frac{\phi'(s_{0})}{1!}(s-s_{0})+\frac{\phi''(s_{0})}{2!}(s-s_{0})^{2}+\cdots
\]
\[
+\frac{\phi^{(m-1)}(s_{0})}{(m-1)!}(s-s_{0})^{m-1}+\sum_{n=m}^{\infty}\frac{\phi^{(n)}(s_{0})}{n!}(s-s_{0})^{n}
\]
in some neighborhood $|s-s_{0}|<\varepsilon$ of $s_{0}$; and from expression~(\ref{eq:chap_7_5_1}) it follows that
\begin{equation}\label{eq:chap_7_5_4}
f(s)=\frac{\phi(s_{0})}{(s-s_{0})^{m}}+\frac{\phi'(s_{0})}{(s-s_{0})^{m-1}}+\frac{\phi''(s_{0})/2!}{(s-s_{0})^{m-2}}+\cdots
\end{equation}
\[
+\frac{\phi^{(m-1)}(s_{0})/(m-1)!}{s-s_{0}}+\sum_{n=m}^{\infty}\frac{\phi^{(n)}(s_{0})}{n!}(s-s_{0})^{n-m}
\]
when $0<|s-s_{0}|<\varepsilon$. This Laurent series representation, together with the fact that $\phi(s_{0})\neq0$, reveals that $s_{0}$ is, indeed, a pole of order $m$ of $f(s)$. The coefficient of $1/(s-s_{0})$ tells us, of course, that the residue of $f(s)$ at $s_{0}$ is as in the statement of the theorem.

Suppose, on the other hand, that we know only that $s_{0}$ is a pole of order $m$ of $f$, or that $f(s)$ has a Laurent series representation
\[
f(s)=\sum_{n=0}^{\infty}a_{n}(s-s_{0})^{n}+\frac{b_{1}}{s-s_{0}}+\frac{b_{2}}{(s-s_{0})^{2}}+\cdots
+\frac{b_{m}}{(s-s_{0})^{m}}\texttt{ }(b_{m}\neq0)
\]
which is valid in a punctured disk $0<|s-s_{0}|<R_{2}$. The function $\phi(s_{0})$ defined by means of the equations
\[
\phi(s)=\{
    \begin{array}{cc}
        (s-s_{0})^{m}f(s) & \texttt{ when }s\neq s_{0}, \\
        b_{m} & \texttt{ when }s=s_{0}
    \end{array}
\]
evidently has the power series representation
\[
\phi(s)=b_{m}+b_{m-1}(s-s_{0})+\cdots+b_{2}(s-s_{0})^{m-2}+b_{1}(s-s_{0})^{m-1}+\sum_{n=0}^{\infty}a_{n}(s-s_{0})^{m+n}
\]
throughout the entire disk $|s-s_{0}|<R_{2}$. Consequently, $\phi(s)$ is analytic in that disk (Sec.~(\ref{sec:chap_6_9_59})) and, in particular, at $s_{0}$. Inasmuch as $\phi(s_{0})=b_{m}\neq0$, expression~(\ref{eq:chap_7_5_1}) is established; and the proof of the theorem is complete.

%\section{Examples}

\section{Zeros of Analytic Functions}\label{sec:chap_7_7_68}

Zeros and poles of functions are closely related. In fact, we shall see in the next section how zeros can be a source of poles. We need, however, some preliminary results regarding zeros of analytic functions.

Suppose that a function $f$ is analytic at a point $s_{0}$. We know from Sec.~(\ref{sec:chap_4_13_48}) that all of the derivatives $f^{(n)}(s)$ $(n=1,2,  \ldots)$ exist at $s_{0}$. If $f(s_{0})=0$ and if there is a positive integer $m$ such that $f^{(m)}(s_{0})\neq0$ and each derivative of lower order vanishes at $s_{0}$, then $f$ is said to have a zero of order $m$ at $s_{0}$. Our first theorem here provides a useful alternative characterization of zeros of order $m$.

\begin{theorem}\label{th:chap_7_7_1}
A function $f$ that is analytic at a point $s_{0}$ has a zero of order $m$ there if and only if there is a function  $g$, which is analytic and nonzero at $s_{0}$, such that
\begin{equation}\label{eq:chap_7_7_1}
f(s)=(s-s_{0})^{m}g(s).
\end{equation}
\end{theorem}

Both parts of the proof that follows use the fact (Sec.~(\ref{sec:chap_6_3_53})) that if a function is
analytic at a point $s_{0}$, then it must have a valid Taylor series representation in powers of $s-s_{0}$ which  is valid throughout a neighborhood $|s-s_{0}|<\varepsilon$ of that point.

We start the first part of the proof by assuming that expression~(\ref{eq:chap_7_7_1}) holds and
noting that, since $g(s)$ is analytic at $s_{0}$, it has a Taylor series representation
\[
g(s)=g(s_{0})+\frac{g'(s_{0})}{1!}(s-s_{0})+\frac{g''(s_{0})}{2!}(s-s_{0})^{2}+\cdots
\]
in some neighborhood $|s-s_{0}|<\varepsilon$ of $s_{0}$. Expression~(\ref{eq:chap_7_7_1}) thus takes the form
\[
f(s)=g(s_{0})(s-s_{0})^{m}+\frac{g'(s_{0})}{1!}(s-s_{0})^{m+1}+\frac{g''(s_{0})}{2!}(s-s_{0})^{m+2}+\cdots
\]
when $|s-s_{0}|<\varepsilon$. Since this is actually a Taylor series expansion for $f(s)$, according
to Theorem~(\ref{th:chap_6_10_1}) in Sec.~(\ref{sec:chap_6_10_60}), it follows that
\begin{equation}\label{eq:chap_7_7_2}
f(s_{0})=f'(s_{0})=f''(s_{0})=\cdots=f^{(m-1)}(s_{0})=0
\end{equation}
and that
\begin{equation}\label{eq:chap_7_7_3}
f^{(m)}(s_{0})=m!g(s_{0})\neq0.
\end{equation}
Hence $s_{0}$ is a zero of order $m$ of $f$.

Conversely, if we assume that $f$ has a zero of order $m$ at $s_{0}$, its analyticity at $s_{0}$ and the fact that conditions~(\ref{eq:chap_7_7_2}) hold tell us that, in some neighborhood $|s-s_{0}|<\varepsilon$, there is a Taylor series
\[
f(s)=\sum_{n=m}^{\infty}\frac{f^{(n)}(s_{0})}{n!}(s-s_{0})^{n}
\]
\[
=(s-s_{0})^{m}[\frac{f(m)(s_{0})}{m!}+\frac{f(m+1)(s_{0})}{(m+1)!}(s-s_{0})+\frac{f(m+2)(s_{0})}{(m+2)!}(s-s_{0})^{2}+\cdots].
\]
Consequently, $f(s)$ has the form~(\ref{eq:chap_7_7_1}), where
\[
g(s)=\frac{f(m)(s_{0})}{m!}+\frac{f(m+1)(s_{0})}{(m+1)!}(s-s_{0})+\frac{f(m+2)(s_{0})}{(m+2)!}(s-s_{0})^{2}+\cdots\texttt{ }(|s-s_{0}|<\varepsilon).
\]
The convergence of this last series when $|s-s_{0}|<\varepsilon$ ensures that $g$ is analytic in that
neighborhood and, in particular, at $s_{0}$ (Sec.~(\ref{sec:chap_6_9_59})). Moreover,
\[
g(s_{0})=\frac{f^{(m)}(s_{0})}{m!}\neq0.
\]
This completes the proof of the theorem.

Our next theorem tells us that the zeros of an analytic function are isolated.

\begin{theorem}\label{th:chap_7_7_2}
Given a function $f$ and a point $s_{0}$, suppose that

( i) $f$ is analytic at $s_{0}$;

(ii) $f(s_{0})=0$ but $f(s)$ is not identically equal to zero in any neighborhood $s_{0}$.

Then $f(s)\neq0$ throughout some deleted neighborhood $0<|s-s_{0}|<\varepsilon$ of $s_{0}$.
\end{theorem}

To prove this, let $f$ be as stated and observe that not all of the derivatives of $f$ at $s_{0}$ are zero. For,  if they were, all of the coefficients in the Taylor series for $f$ about $s_{0}$ would be zero; and that would mean that $f(s)$ is identically equal to zero in some neighborhood of $s_{0}$. So it is clear from the definition of zeros of order m at the beginning of this section that $f$ must have a zero of some order $m$ at $s_{0}$. According  to Theorem~(\ref{th:chap_7_7_1}), then,
\begin{equation}\label{eq:chap_7_7_4}
f(s)=(s-s_{0})^{m}g(s)
\end{equation}
where $g(s)$ is analytic and nonzero at $s_{0}$.

Now $g$ is continuous,in addition to being nonzero, at $s_{0}$ because it is analytic there. Hence there is some neighborhood $|s-s_{0}|<\varepsilon$ in which equation~(\ref{eq:chap_7_7_4}) holds and
in which $g(s)\neq0$ (see Sec.~(\ref{sec:chap_2_7_17})). Consequently, $f(s)\neq0$ in the deleted neighborhood $0<|s-s_{0}|<\varepsilon$; and the proof is complete.

Our final theorem here concerns functions with zeros that are not all isolated It was referred to earlier in Sec. ~(\ref{sec:chap_2_16_26}) and makes an interesting contrast to Theorem~(\ref{th:chap_7_7_2}) just above.

\begin{theorem}\label{th:chap_7_7_3}
Given a function $f$ and a point $s_{0}$, suppose that

(i) $f$ is analytic throughout a neighborhood $N_{0}$ of $s_{0}$;

(ii) $f(s_{0})=0$ and $f(s)=0$ at each point $s$ of a domain or line segment containing $s_{0}$.

Then $f(s)\equiv0$ in $N_{0}$; that is, $f(s)$ is identically equal to zero throughout $N_{0}$.
\end{theorem}

We begin the proof with the observation that, under the stated conditions,
$f(s)\equiv0$ in some neighborhood $N$ of $s_{0}$. For, otherwise, there would be a deleted neighborhood of $s_{0}$  throughout which $f(s)\neq0$, according to Theorem~(\ref{th:chap_7_7_2}) above; and that would be inconsistent with the condition that $f(s)=0$ everywhere in a domain or on a line segment containing $s_{0}$. Since $f(s)\equiv0$ in the neighborhood $N$, then, it follows that all of the coefficients
\[
a_{n}=\frac{f^{(n)}(s_{0})}{n!}\texttt{ }(n=1,2,\ldots)
\]
in the Taylor series for $f(s)$ about $s_{0}$ must be zero. Thus $f(s)=0$ in the neighborhood $N_{0}$, since Taylor series also represents $f(s)$ in $N_{0}$. This completes the proof.

\section{Zeros and Poles}\label{sec:chap_7_8_69}

The following theorem shows how zeros of order $m$ can create poles of order $m$.

\begin{theorem}\label{th:chap_7_8_1}
Suppose that

(i) two functions $p$ and $q$ are analytic at a point $s_{0}$;

(ii) $p(s_{0})\neq0$ and $q$ has a zero of order $m$ at $s_{0}$.

Then the quotient $p(s)/q(s)$ has a pole of order $m$ at $s_{0}$.
\end{theorem}

The proof is easy. Let $p$ and $q$ be as in the statement of the theorem. Since $q$ has a zero of order $m$ at $s_{0}$, we know from Theorem~(\ref{th:chap_7_7_2}) in Sec.~(\ref{sec:chap_7_7_68}) that there is a deleted neighborhood of $s_{0}$ in which $q(s)\neq0$; and so $s_{0}$ is an isolated singular point of the quotient  $p(s)/q(s)$. Theorem~(\ref{th:chap_7_7_1}) in Sec.~(\ref{sec:chap_7_7_68}) tells us, moreover. that
\[
q(s)=(s-s_{0})^{m}g(s)
\]
where $g$ is analytic and nonzero at $s_{0}$; and this enables us to write
\begin{equation}\label{eq:chap_7_8_1}
\frac{p(s)}{q(s)}=\frac{p(s)/g(s)}{(s-s_{0})^{m}}
\end{equation}
Since $p(s)/g(s)$ is analytic and nonzero at $s_{0}$, it now follows from the theorem in Sec.~(\ref{sec:chap_7_5_66}) that $s_{0}$ is a pole of order $m$ of $p(s)/q(s)$.

Theorem~(\ref{th:chap_7_8_1}) leads us to another method for identifying simple poles and finding the corresponding residues. This method is sometimes easier to use than the one in Sec.~(\ref{sec:chap_7_5_66}).

\begin{theorem}\label{th:chap_7_8_1}
Let two functions $p$ and $q$ be analytic at a point $s_{0}$. If
\[
p(s_{0})\neq0,\texttt{ }q(s_{0})=0,\texttt{ and }q'(s_{0})\neq0,
\]
then $s_{0}$ is a simple pole of the quotient $p(s)/q(s)$ and
\begin{equation}\label{eq:chap_7_8_2}
Res_{s=s_{0}}\frac{p(s)}{q(s)}=\frac{p(s_{0})}{q'(s_{0})}.
\end{equation}
\end{theorem}

To show this, we assume that $p$ and $q$ are as stated and observe that, because of the conditions on $q$, the point $s_{0}$ is a zero of order $m=1$ of that function. According to Theorem~(\ref{th:chap_7_7_1})  in Sec.~(\ref{sec:chap_7_7_68}),  then,
\begin{equation}\label{eq:chap_7_8_3}
q(s)=(s-s_{0})g(s)
\end{equation}
where $g(s)$ is analytic and nonzero at $s_{0}$. Furthermore, Theorem~(\ref{th:chap_7_8_1}) in this section tells
us that $s_{0}$ is a simple pole of $p(s)/q(s)$; and equation~(\ref{eq:chap_7_8_1}) in its proof becomes
\[
\frac{p(s)}{q(s)}=\frac{p(s)/g(s)}{s-s_{0}}.
\]
Now $q(s)/g(s)$ is analytic and nonzero at $s_{0}$, and it follows from the theorem in Sec.~(\ref{sec:chap_7_5_66}) that
\begin{equation}\label{eq:chap_7_8_4}
Res_{s=s_{0}}\frac{p(s)}{q(s)}=\frac{p(s_{0})}{g(s_{0})}.
\end{equation}
But $g(s_{0})=q'(s_{0})$, as is seen by differentiating each side of equation~(\ref{eq:chap_7_8_3}) and setting $s=s_{0}$. Expression~(\ref{eq:chap_7_8_4}) thus takes the form~(\ref{eq:chap_7_8_2}).

There are formulas similar to formula~(\ref{eq:chap_7_8_2}) for residues at poles of higher order, but they are lengthier and, in  general, not practical.

\section{Behavior of $f$ Near Isolated Singular Points}\label{sec:chap_7_9_70}

As already indicated in Sec.~(\ref{sec:chap_7_4_65}), the behavior of a function $f$ near an isolated singular
point $s_{0}$ varies, depending on whether $s_{0}$ is a pole, a removable singular point, or an essential singular point. In this section, we develop the differences in behavior somewhat further. %Since the results presented here will not be used elsewhere in the book, the reader who wishes to reach applications of residue theory more quickly may pass directly to Chap.~(\ref{ch:chap_7}) without disruption.

\begin{theorem}\label{th:chap_7_9_1}
If $s_{0}$ is a pole of a function $f$, then
\begin{equation}\label{eq:chap_7_9_1}
\lim_{s=s_{0}}f(s)=\infty
\end{equation}
\end{theorem}

To verify limit~(\ref{eq:chap_7_9_1}), we assume that $f$ has a pole of order $m$ at $s_{0}$ and use the
theorem in Sec.~(\ref{sec:chap_7_5_66}). It tells us that
\[
f(s)=\frac{\phi(s)}{(s-s_{0})^{m}},
\]
where $\phi(s)$ is analytic and nonzero at $s_{0}$. Since
\[
\lim_{s=s_{0}}\frac{1}{f(s)}=\lim_{s=s_{0}}\frac{(s-s_{0})^{m}}{\phi(s)}=\frac{\lim_{s=s_{0}}(s-s_{0})^{m}}{\lim_{s=s_{0}}\phi(s)}=\frac{0}{\phi(s_{0})}=0,
\]
then, limit~(\ref{eq:chap_7_9_1}) holds, according to the theorem in Sec.~(\ref{sec:chap_2_6_16}) regarding limits that involve the point at infinity.

The next theorem emphasizes how the behavior of $f$ near a removable singular point is fundamentally different from  the behavior near a pole.

\begin{theorem}\label{th:chap_7_9_2}
If $s_{0}$ is a removable singular point of a function $f$, then $f$ is analytic and bounded in some deleted neighborhood $0<|s-s_{0}|<\varepsilon$ of $s_{0}$.
\end{theorem}

The proof is easy and is based on the fact that the function $f$ here is analytic in a disk $|s-s_{0}|<R_{2}$ when  $f(s_{0})$ is properly defined; and $f$ is then continuous in any closed disk $|s-s_{0}|\leq\varepsilon$ where $\varepsilon<R_{2}$. Consequently, $f$ is bounded in that disk, according to Sec.~(\ref{sec:chap_2_7_17}); and this means that, in addition to being analytic, $f$ must be bounded in the deleted neighborhood $0<|s-s_{0}|<\varepsilon$.

The proof of our final theorem, regarding the behavior of a function near an essential singular point, relies on the following lemma, which is closely related to Theorem~(\ref{th:chap_7_9_2}) and is known as Riemann's theorem.

\begin{lemma}\label{le:chap_7_9_1}
Suppose that a function $f$ is analytic and bounded in some deleted neighborhood $0<|s-s_{0}|<\varepsilon$ of a point $s_{0}$. If $f$ is not analytic at $s_{0}$, then it has a removable singularity there.
\end{lemma}

To  prove this, we assume that $f$ is not analytic at $s_{0}$. As a consequence, the point $s_{0}$ must be an isolated singularity of $f$; and $f(s)$ is represented by a Laurent series
\begin{equation}\label{eq:chap_7_9_2}
f(s)=\sum_{n=0}^{\infty}a_{n}(s-s_{0})^{n}+\sum_{n=1}^{\infty}\frac{b_{n}}{(s-s_{0})^{n}}\texttt{ }(R_{1}<|s-s_{0}|<R_{2})
\end{equation}
throughout the deleted neighborhood $0<|s-s_{0}|<\varepsilon$. If $C$ denotes a positively oriented circle $|s-s_{0}|=\rho$, where $\rho<\varepsilon$, we know from Sec.~(\ref{sec:chap_6_5_55}) that the coefficients $b_{n}$  in expansion~(\ref{eq:chap_7_9_2}) can be written
\begin{equation}\label{eq:chap_7_9_3}
b_{n}=\frac{1}{2\pi i}\int_{C}\frac{f(s)ds}{(s-s_{0})^{-n+1}}=e_{xy}(b_{n,x}+ib_{n,y})+b_{n,z}\texttt{ }(n=0,1,2,\ldots).
\end{equation}
Now the boundedness condition on $f$ tells us that there is a positive constant $M$ such
that $|f(s)|\leq M$ whenever $0<|s-s_{0}|<\varepsilon$. Hence it follows from expression~(\ref{eq:chap_7_9_3}) that
\[
|b_{n}|\leq\frac{1}{2\pi}\cdot\frac{M}{\rho^{-n+1}}2\pi\rho=M\rho^{n}\texttt{ }(n=0,1,2,\ldots).
\]
Since the coefficients $b_{n}$ are constants and since $\rho$ can be chosen arbitrarily small, we
may conclude that $b_{n}=0$ $(n=1,2,\ldots)$ in the Laurent series~(\ref{eq:chap_7_9_2}). This tells us that $s_{0}$ is a removable singularity of $f$, and the proof of the lemma is complete.

We know from Sec.~(\ref{sec:chap_7_4_65}) that the behavior of a function near an essential singular point is quite irregular. The theorem below, regarding such behavior, is related to Picard's theorem in that earlier section and is usually referred to as the Casorati-Weierstrass theorem. It states that, in each deleted neighborhood of an essential singular point, a function assumes values arbitrarily close to any given number.

\begin{theorem}\label{th:chap_7_9_3}
Suppose that $s_{0}$ is an essential singularity of a function $f$, and let $\varpi_{0}$ be any spatial complex number. Then, for any positive number $\varepsilon$, the inequality
\begin{equation}\label{eq:chap_7_9_4}
|f(s)-\varpi_{0}|<\varepsilon
\end{equation}
is satisfied at some point $s$ in each deleted neighborhood $0<|s-s_{0}|<\delta$ of $s_{0}$.
\end{theorem}

The proof is by contradiction. Since $s_{0}$ is an isolated singularity of $f$, there is a deleted neighborhood $0<|s-s_{0}|<\delta$ throughout which $f$ is analytic; and we assume that condition~(\ref{eq:chap_7_9_4}) is not satisfied for any point $s$ there. Thus $|f(s)-\varpi_{0}|\geq\varepsilon$ when
$0<|s-s_{0}|<\delta$; and so the function
\begin{equation}\label{eq:chap_7_9_5}
g(s)=\frac{1}{f(s)-\varpi_{0}}\texttt{ }(0<|s-s_{0}|<\delta)
\end{equation}
is bounded and analytic in its domain of definition. Hence, according to the above lemma, $s_{0}$ is a removable singularity of $g$; and we let $g$ be defined at $s_{0}$ so that it is analytic there.

If $g(s_{0})\neq0$, the function $f(s)$, which can be written
\begin{equation}\label{eq:chap_7_9_6}
f(s)=\frac{1}{g(s)}+\varpi_{0}
\end{equation}
when $0<|s-s_{0}|<\delta$, becomes analytic at $s_{0}$ if it is defined there as
\[
f(s_{0})=\frac{1}{g(s_{0})}+\varpi_{0}.
\]
But this means that $s_{0}$ is a removable singularity of $f$, not an essential one, and we have a contradiction.

If $g(s_{0})=0$, the function $g$ must have a zero of some finite order $m$ (Sec.~(\ref{sec:chap_7_7_68})) at $s_{0}$ because $g(s)$ is not identically equal to zero in the neighborhood $|s-s_{0}|<\delta$. In view of equation~(\ref{eq:chap_7_9_6}), then, $f$ has a pole of order $m$ at $s_{0}$ (see Theorem~(\ref{th:chap_7_8_1}) in Sec.~(\ref{sec:chap_7_8_69})). So, once again, we have a contradiction; and Theorem~(\ref{th:chap_7_9_3}) here is proven.

%-----------------------------------------------------------------------
% Beginning of chap8.tex
%-----------------------------------------------------------------------
%
% AMS-LaTeX 1.2 sample file for a monograph, based on amsbook.cls.
% This is a data file input by chapter.tex.
%%%%%%%%%%%%%%%%%%%%%%%%%%%%%%%%%%%%%%%%%%%%%%%%%%%%%%%%%%%%%%%%%%%%

%\part{This is a Part Title Sample}

\chapter{Mapping by Elementary Functions}\label{ch:chap_8}

The geometric interpretation of a spatial function of a spatial complex variable as a mapping, or transformation, was introduced in Secs.~(\ref{sec:chap_2_2_12}) %and ~(\ref{sec:chap_2_3_13})
(Chap.~(\ref{ch:chap_2})). We saw there how the nature of such a function can be displayed graphically, to some extent, by the manner in which it maps certain curves and regions. In this chapter, we shall see further examples of how various curves and regions are mapped by elementary analytic functions. Applications of such results to physical problems are illustrated in Chaps.~(\ref{ch:chap_10}).

\section{Linear Transformations}\label{sec:chap_8_1_83}

To study the mapping
\begin{equation}\label{eq:chap_8_1_1}
\varpi=As
\end{equation}
where $A$ is a nonzero spatial complex constant and $s\neq0$, we write $A$ and $s$ in exponential form:
\[
A=r_{a}(e_{xy}e^{i\theta_{a}}\cos\varphi_{a}+\sin\varphi_{a}),\texttt{ }
s=r(e_{xy}e^{i\theta}\cos\varphi+\sin\varphi)
\]
Then
\begin{equation}\label{eq:chap_8_1_2}
\varpi=r_{a}r[e_{xy}(e^{i(\theta_{a}+\theta)}\cos\varphi_{a}\cos\varphi
+e^{i\theta_{a}}\cos\varphi_{a}\sin\varphi+\sin\varphi_{a}e^{i\theta}\cos\varphi)
\end{equation}
\[
+\sin\varphi_{a}\sin\varphi]
=r_{a}r(e_{xy}e^{i\theta_{p}}\cos\varphi_{p}+\sin\varphi_{p})
\]
where $\varphi_{p}=\arcsin(\sin\varphi_{a}\sin\varphi)$,
\[
\theta_{p}=\arctan\frac{\sin(\theta_{a}+\theta)
+\sin\theta_{a}\tan\varphi+\tan\varphi_{a}\sin\theta}
{\cos(\theta_{a}+\theta)
+\cos\theta_{a}\tan\varphi+\tan\varphi_{a}\cos\theta},
\]
and we see from equation~(\ref{eq:chap_8_1_2}) that transformation~(\ref{eq:chap_8_1_1}) expands or contracts the radius vector representing $s$ by the factor $r_{a}=|A|$ and rotates it through two angles $\theta_{a}=\arg_{c}A$ and $\varphi_{a}=\arg_{r}A$ about the origin. The image of a given region is, therefore, geometrically similar to
that region.

So transformation~(\ref{eq:chap_8_1_1}) maps a sphere in $s$ space to a sphere in $\varpi$ space.

The mapping
\begin{equation}\label{eq:chap_8_1_3}
\varpi=s+B,
\end{equation}
where $B$ is any spatial complex constant, is a translation by means of the vector representing $B$. That is, if
\[
\varpi=e_{xy}(u+iv)+w,\texttt{ }s=e_{xy}(x+iy)+z,\texttt{ and }B=e_{xy}(b_{x}+ib_{y})+b_{z},
\]
then the image of any point $(x,y,z)$ in the $s$ space is the point
\begin{equation}\label{eq:chap_8_1_4}
(u,v,w)=(x+b_{x},y+b_{y},z+b_{z})
\end{equation}
in the $\varpi$ space. Since each point in any given region of the $s$ space is mapped into the $\varpi$ space in this manner, the image region is geometrically congruent to the original one.

The general (nonconstant) linear transformation
\begin{equation}\label{eq:chap_8_1_5}
\varpi=As+B\texttt{ }(A\neq0),
\end{equation}
which is a composition of the transformations
\[
\varpi_{1}=As\texttt{ }(A\neq0)\texttt{ and }\varpi=\varpi_{1}+B
\]
is evidently an expansion or contraction and a rotation, followed by a translation.

\section{The Transformation $\varpi$=l/s}\label{sec:chap_8_2_84}

The equation
\begin{equation}\label{eq:chap_8_2_1}
\varpi=1/s
\end{equation}
establishes a one to one correspondence between the nonzero points of the $s$ and the $\varpi$ spaces. Since $s\bar{s}=|s|^{2}$, the mapping can be described by means of the successive transformations
\begin{equation}\label{eq:chap_8_2_2}
\varpi_{1}=\frac{1}{|s|^{2}}s\texttt{ and }\varpi=\overline{\varpi_{1}}.
\end{equation}

The first of these transformations is an inversion with respect to the unit sphere
$|s|=1$. That is, the image of a nonzero point $s$ is the point $\varpi_{1}$ with the properties
\[
|\varpi_{1}|=\frac{1}{|s|}\texttt{ and }\arg\varpi_{1}=\arg s.
\]
Thus the points exterior to the sphere $|s|=1$ are mapped onto the nonzero points interior to it, and conversely. Any point on the sphere is mapped onto itself. The second of transformations~(\ref{eq:chap_8_2_2}) is simply a reflection on the $xz$ coordinate plane.

If we write transformation~(\ref{eq:chap_8_2_1}) as
\begin{equation}\label{eq:chap_8_2_3}
T(s)=1/s\texttt{ }(s\neq0),
\end{equation}
we can define $T$ at the origin and at the point at infinity so as to be continuous on the extended spatial  complex space. To do this, we need only refer to Sec.~(\ref{sec:chap_2_6_16}) to see that
\begin{equation}\label{eq:chap_8_2_4}
\lim_{s \to 0}T(s)=\infty\texttt{ since }\lim_{s \to 0}\frac{1}{T(s)}=0,
\end{equation}
and
\begin{equation}\label{eq:chap_8_2_5}
\lim_{s \to \infty}T(s)=0\texttt{ since }\lim_{s \to 0}T(\frac{1}{s})=0.
\end{equation}
In order to make $T$ continuous on the extended space, then, we write
\begin{equation}\label{eq:chap_8_2_6}
T(0)=\infty,{ }T(\infty)=0,\texttt{ and }T(s)=\frac{1}{s}
\end{equation}
for the remaining values of $s$. More precisely, equations~(\ref{eq:chap_8_2_6}), together with the first of limits~(\ref{eq:chap_8_2_4}) and~(\ref{eq:chap_8_2_5}), show that
\begin{equation}\label{eq:chap_8_2_7}
\lim_{s \to s_{0}}T(s)=T(s_{0})
\end{equation}
for every point $s_{0}$ in the extended space, including $s_{0}=0$ and $s_{0}=\infty$. The fact that $T$ is continuous everywhere in the extended space is now a consequence of equation~(\ref{eq:chap_8_2_7}) (see Sec.~(\ref{sec:chap_2_7_17})). Because of this continuity, when the point at infinity is involved in any discussion  of the function $1/s$, it is tacitly assumed that $T(s)$ is intended.

\section{Mappings by 1/s}\label{sec:chap_8_3_85}

When a point $\varpi=e_{xy}(u+iv)+w$ is the image of a nonzero point $s=e_{xy}(x+iy)+z$ where $x+z\neq0\texttt{ or }y\neq0$ under the transformation $\varpi=1/s$, writing $\varpi=\bar{s}/(s\bar{s})$ reveals that
\begin{equation}\label{eq:chap_8_3_1}
u=\frac{x+z}{(x+z)^{2}+y^{2}}-\frac{1}{z},\texttt{ }v=-\frac{y}{(x+z)^{2}+y^{2}},\texttt{ }w=\frac{1}{z}\texttt{ where }z\neq0
\end{equation}
or
\begin{equation}\label{eq:chap_8_3_1b}
u+w=\frac{x+z}{(x+z)^{2}+y^{2}},\texttt{ }v=-\frac{y}{(x+z)^{2}+y^{2}},\texttt{ }w=\frac{1}{z}\texttt{ where }z\neq0.
\end{equation}
Also, since $s=1/\varpi=\overline{\varpi}/(\varpi\overline{\varpi})$ where $u+w\neq0\texttt{ or }v\neq0$,
\begin{equation}\label{eq:chap_8_3_2}
x=\frac{u+w}{(u+w)^{2}+v^{2}}-\frac{1}{w},\texttt{ }y=-\frac{v}{(u+w)^{2}+v^{2}},\texttt{ }z=\frac{1}{w}\texttt{ where }w\neq0.
\end{equation}
or
\begin{equation}\label{eq:chap_8_3_2b}
x+z=\frac{u+w}{(u+w)^{2}+v^{2}},\texttt{ }y=-\frac{v}{(u+w)^{2}+v^{2}},\texttt{ }z=\frac{1}{w}\texttt{ where }w\neq0.
\end{equation}
The following argument, based on these relations between coordinates, shows that the mapping $\varpi=1/s$ transforms spheres, circles, and lines into spheres, circles, and lines. When $A$, $B$, $C$, and $D$ are all real numbers satisfying the condition $B^{2}+C^{2}>4AD$, the equation
\begin{equation}\label{eq:chap_8_3_3}
A[(x+z)^{2}+y^{2}]+B(x+z)+Cy+D=0
\end{equation}
represents an arbitrary sphere or circle or line, where $A\neq0$ for a sphere or circle, and $A=0$ for a line. The need for the condition $B^{2}+C^{2}>4AD$ when $A\neq0$ is evident if, by the method of completing the squares, we rewrite equation~(\ref{eq:chap_8_3_3}) as
\[
(x+z+\frac{B}{2A})^{2}+(y+\frac{B}{2A})^{2}=(\frac{\sqrt{B^{2}+C^{2}-4AD}}{2A})^{2}.
\]
When $A=0$, the condition becomes $B^{2}+C^{2}>0$, which means that $B$ and $C$ are not both zero. Returning to the verification of the statement in italics, we observe that if
$x$, $y$, and $z$ satisfy equation~(\ref{eq:chap_8_3_3}), we can use relations~(\ref{eq:chap_8_3_2}) to substitute for those variables.

After some simplifications, we find that $u$, $v$, and $w$ satisfy the equation
\begin{equation}\label{eq:chap_8_3_4}
D[(u+w)^{2}+v^{2}]+B(u+w)+Cv-A=0
\end{equation}
which also represents a sphere or circle or line. Conversely, if $u$, $v$, and $w$ satisfy equation~(\ref{eq:chap_8_3_4}), it follows from relations~(\ref{eq:chap_8_3_1}) that $x$, $y$, and $z$ satisfy equation~(\ref{eq:chap_8_3_3}).

It is now clear from equations~(\ref{eq:chap_8_3_3}) and~(\ref{eq:chap_8_3_4}) that

(i) a sphere or circle $(A\neq0)$ not passing through the origin $(D\neq0)$ in the $s$ space is transformed into a sphere or circle not passing through the origin in the $\varpi$ space;

(ii) a sphere or circle $(A\neq0)$ through the origin $(D=0)$ in the $s$ space is transformed into a
line that does not pass through the origin in the $\varpi$ space;

(iii) a line $(A=0)$ not passing through the origin $(D\neq0)$ in the $s$ space is transformed
into a sphere or circle through the origin in the $\varpi$ space;

(iv) a line $(A=0)$ through the origin $(D=0)$ in the $s$ space is transformed into a line
through the origin in the $\varpi$ space.

\section{Linear Fractional Transformations}\label{sec:chap_8_4_86}

The transformation
\begin{equation}\label{eq:chap_8_4_1}
\varpi=\frac{as+b}{cs+d}\texttt{ }(ad-bc\neq0),
\end{equation}
where $a$, $b$,  $c$, and $d$ are spatial complex constants, is called a linear fractional transformation,
or M$\ddot{o}$bius transformation. Observe that equation~(\ref{eq:chap_8_4_1}) can be written in the form
\begin{equation}\label{eq:chap_8_4_2}
As\varpi+Bs+C\varpi+D=0\texttt{ }(AD-BC\neq0);
\end{equation}
and, conversely, any equation of type~(\ref{eq:chap_8_4_2}) can be put in the form~(\ref{eq:chap_8_4_1}). Since this alternative form is linear in $s$ and linear in $\varpi$, or bilinear in $s$ and $\varpi$, another name for a linear fractional transformation is bilinear transformation.

When $c=0$, the condition $ad- bc\neq0$ with equation~(\ref{eq:chap_8_4_1}) becomes $ad\neq0$; and we see that the transformation reduces to a nonconstant linear function. When $c\neq0$, equation~(\ref{eq:chap_8_4_1}) can be written
\begin{equation}\label{eq:chap_8_4_3}
\varpi=\frac{a}{c}+\frac{bc-ad}{c}\cdot\frac{1}{cs+d}\texttt{ }(ad-bc\neq0).
\end{equation}
So, once again, the condition $ad-bc\neq0$ ensures that we do not have a constant function. The transformation $\varpi=1/s$ is evidently a special case of transformation~(\ref{eq:chap_8_4_1}) when $c\neq0$. Equation~(\ref{eq:chap_8_4_3}) reveals that when $c\neq0$, a linear fractional transformation is a composition of  the mappings.
\[
\varpi_{1}=cs+d,\texttt{ }\varpi_{2}=\frac{1}{\varpi_{1}},\texttt{ }
\varpi=\frac{a}{c}+\frac{bc-ad}{c}\varpi_{2}\texttt{ }(ad-bc\neq0).
\]
It thus follows that, regardless of whether $c$ is zero or nonzero, any linear fractional transformation transforms spheres and lines into spheres and lines because these special linear fractional transformations do. (See Secs.~(\ref{sec:chap_8_1_83}) and~(\ref{sec:chap_8_3_85}).)

Solving equation~(\ref{eq:chap_8_4_1}) for $s$, we find that
\begin{equation}\label{eq:chap_8_4_4}
s=\frac{-d\varpi+b}{c\varpi-a}\texttt{ }(ad-bc\neq0).
\end{equation}
When a given point $\varpi$ is the image of some point $s$ under transformation~(\ref{eq:chap_8_4_1}), the point $s$ is retrieved by means of equation~(\ref{eq:chap_8_4_4}). If $c=0$, so that $a$ and $d$ are both nonzero, each point in the $\varpi$ space is evidently the image of one and only one point in the $s$ space. The same is true if  $c\neq0$, except when $\varpi=a/c$ since the denominator in equation~(\ref{eq:chap_8_4_4}) vanishes if $\varpi$ has that value. We can, however, enlarge the domain of definition of transformation~(\ref{eq:chap_8_4_1}) in order to define a linear fractional transformation $T$ on the extended $s$ space such that the point $\varpi=a/c$ is the image of $s=\infty$ when $c\neq0$. We first write
\begin{equation}\label{eq:chap_8_4_5}
T(s)=\frac{as+b}{cs+d}\texttt{ }(ad-bc\neq0).
\end{equation}
We then write
\[
T(\infty)=\infty\texttt{ if }c=0
\]
and
\[
T(\infty)=\frac{a}{c}\texttt{ and }T(-\frac{d}{c})=\infty\texttt{ if }c\neq0.
\]
We can see that this makes $T$ continuous on the extended $s$ space. It also agrees with the way in which we enlarged the domain of definition of the transformation $\varpi=1/s$ in Sec.~(\ref{sec:chap_8_2_84}).

When its domain of definition is enlarged in this way, the linear fractional transformation~(\ref{eq:chap_8_4_5}) is a one to one mapping of the extended $s$ space onto the extended $\varpi$ space. That is, $T(s_{1})\neq T(s_{2})$ whenever $s_{1}\neq s_{2}$; and, for each point $\varpi$ in the second
space, there is a point $s$ in the first one such that $T(s)=\varpi$. Hence, associated with the transformation  $T$, there is an inverse transformation $T^{-1}$, which is defined on the extended $\varpi$ space as follows:
\[
T^{-1}(\varpi)=s\texttt{ if and only if }T(s)=\varpi.
\]
From equation~(\ref{eq:chap_8_4_4}), we see that
\begin{equation}\label{eq:chap_8_4_6}
T^{-1}(\varpi)=\frac{-d\varpi+b}{c\varpi-a}\texttt{ }(ad-bc\neq0).
\end{equation}
Evidently, $T^{-1}$ is itself a linear fractional transformation, where
\[
T^{-1}(\infty)=\infty\texttt{ if }c=0
\]
and
\[
T^{-1}(\frac{a}{c})=\infty\texttt{ and }T^{-1}(\infty)=-\frac{d}{c}\texttt{ if }c\neq0.
\]
If $T$ and $S$ are two linear fractional transformations, then so is the composition $S(T(s))$.This can be verified by combining expressions of the type~(\ref{eq:chap_8_4_5}). Note that, in particular,
$T^{-1}(T(s))=s$ for each point $s$ in the extended space.

There is always a linear fractional transformation that maps three given distinct points $s_{1}$, $s_{2}$, and $s_{3}$ onto three specified distinct points $\varpi_{1}$, $\varpi_{2}$, and $\varpi_{3}$, respectively.
Verification of this will appear in Sec.~(\ref{sec:chap_8_5_87}), where the image $\varpi$ of a point $s$ under such a transformation is given implicitly in terms of $s$. We illustrate here a more direct approach to finding the desired transformation.

\section{An Implicit Form}\label{sec:chap_8_5_87}

The equation
\begin{equation}\label{eq:chap_8_5_1}
\frac{(\varpi-\varpi_{1})(\varpi_{2}-\varpi_{3})}{(\varpi-\varpi_{3})(\varpi_{2}-\varpi_{1})}
=\frac{(s-s_{1})(s_{2}-s_{3})}{(s-s_{3})(s_{2}-s_{1})}
\end{equation}
defines (implicitly) a linear fractional transformation that maps distinct points $s_{1}$, $s_{2}$, and $s_{3}$ in the finite $s$ space onto distinct points $\varpi_{1}$, $\varpi_{2}$, and $\varpi_{3}$, respectively, in the finite $\varpi$ space. To verify this, we write equation~(\ref{eq:chap_8_5_1}) as
\begin{equation}\label{eq:chap_8_5_2}
(s-s_{3})(\varpi-\varpi_{1})(s_{2}-s_{1})(\varpi_{2}-\varpi_{3})
=(s-s_{1})(\varpi-\varpi_{3})(s_{2}-s_{3})(\varpi_{2}-\varpi_{1}).
\end{equation}
If $s=s_{1}$ the right-hand side of equation~(\ref{eq:chap_8_5_2}) is zero; and it follows that $\varpi=\varpi_{1}$. Similarly, if $s=s_{3}$, the left-hand side is zero and, consequently, $\varpi=\varpi_{3}$. If  $s=s_{2}$, we have the linear equation
\[
(\varpi-\varpi_{1})(\varpi_{2}-\varpi_{3})=(\varpi-\varpi_{3})(\varpi_{2}-\varpi_{1}),
\]
whose unique solution is $\varpi=\varpi_{2}$. One can see that the mapping defined by equation~(\ref{eq:chap_8_5_1}) is actually a linear fractional transformation by expanding the products in equation~(\ref{eq:chap_8_5_2}) and writing the result in the form (Sec.~(\ref{sec:chap_8_4_86}))
\begin{equation}\label{eq:chap_8_5_3}
As\varpi+Bs+C\varpi+D=0\texttt{ }(AD-BC\neq0);
\end{equation}
The condition $AD-BC\neq0$, which is needed with equation~(\ref{eq:chap_8_5_3}), is clearly satisfied
since, as just demonstrated, equation~(\ref{eq:chap_8_5_1}) does not define a constant function. It can be shown that equation~(\ref{eq:chap_8_5_1}) defines the only linear fractional transformation mapping the points $s_{1}$, $s_{2}$, and $s_{3}$ onto $\varpi_{1}$, $\varpi_{2}$, and $\varpi_{3}$ respectively.

If equation~(\ref{eq:chap_8_5_1}) is modified properly, it can also be used when the point at infinity is one of the prescribed points in either the (extended) $s$ or $\varpi$ space. Suppose, for instance, that $s_{1}=\infty$. Since any linear fractional transformation is continuous on the extended space, we need only replace $s_{1}$ on the right-hand side of equation~(\ref{eq:chap_8_5_1}) by $1/s_{1}$, clear fractions, and let $s_{1}$ tend to zero:
\[
\lim_{s_{1} \to 0}\frac{(s-1/s_{1})(s_{2}-s_{3})}{(s-s_{3})(s_{2}-1/s_{1})}\cdot\frac{s_{1}}{s_{1}}
=\lim_{s_{1} \to 0}\frac{(ss_{1}-1)(s_{2}-s_{3})}{(s-s_{3})(s_{1}s_{2}-1)}
=\frac{s_{2}-s_{3}}{s-s_{3}}.
\]
The desired modification of equation~(\ref{eq:chap_8_5_1}) is, then,
\[
\frac{(\varpi-\varpi_{1})(\varpi_{2}-\varpi_{3})}{(\varpi-\varpi_{3})(\varpi_{2}-\varpi_{1})}
=\frac{s_{2}-s_{3}}{s-s_{3}}.
\]
Note that this modification is obtained formally by simply deleting the factors involving $s_{1}$ in equation~(\ref{eq:chap_8_5_1}). It is easy to check that the same formal approach applies when
any of the other prescribed points is $\infty$.

\section{Mappings of the Upper Half Space}\label{sec:chap_8_6_88}

Let us determine all linear fractional transformations that map the half space $Im\texttt{ }s>0$ onto the open sphere $|\varpi|<1$ and the plane boundary $Im\texttt{ }s=0$ onto the spherical boundary $|\varpi|=1$.

Keeping in mind that points on the plane $Im\texttt{ }s=0$ are to be transformed into points on the sphere $|\varpi|=1$, we start by selecting the points $s=0$, $s=1$, and $s=\infty$ in the plane and determining conditions on a linear fractional transformation
\begin{equation}\label{eq:chap_8_6_1}
\varpi=\frac{as+b}{cs+d}\texttt{ }(ad-bc\neq0),
\end{equation}
which are necessary in order for the images of those points to have unit modulus.

We note from equation~(\ref{eq:chap_8_6_1}) that if $|\varpi|=1$ when $s=0$, then $|b/d|=1$; that is,
\begin{equation}\label{eq:chap_8_6_2}
|b|=|d|\neq0.
\end{equation}

Now, according to Sec.~(\ref{sec:chap_8_4_86}), the image $\varpi$ of the point $s=\infty$ is a finite number, namely $\varpi=a/c$ only if $c\neq0$. So the requirement that $|\varpi|=1$ when $s=\infty$ means that
$|a/c|=1$, or
\begin{equation}\label{eq:chap_8_6_3}
|a|=|c|\neq0;
\end{equation}
and the fact that $a$ and $c$ are nonzero enables us to rewrite equation~(\ref{eq:chap_8_6_1}) as
\begin{equation}\label{eq:chap_8_6_4}
\varpi=\frac{a}{c}\cdot\frac{s+(b/a)}{s+(d/c)}\texttt{ }(ad-bc\neq0),
\end{equation}
Then, since $|a/c|=1$ and
\[
|\frac{b}{a}|=|\frac{d}{c}|\neq0,
\]
according to relations~(\ref{eq:chap_8_6_2}) and~(\ref{eq:chap_8_6_3}), equation~(\ref{eq:chap_8_6_4}) can be put in the form
\begin{equation}\label{eq:chap_8_6_5}
\varpi=(e_{xy}e^{i\alpha}\cos\beta+\sin\beta)\frac{s-s_{0}}{s-s_{1}}\texttt{ }(|s_{1}|=|s_{0}|\neq0)
\end{equation}
where $\alpha$ and $\beta$ are reals, and $s_{0}$ and $s_{1}$ are (nonzero) spatial complex constants.

Next, we impose on transformation~(\ref{eq:chap_8_6_5}) the condition that $|\varpi|=1$ when $s=1$.
This tells us that
\[
|1-s_{1}|=|1-s_{0}|,
\]
or
\[
(1-s_{1})(1-\overline{s_{1}})=(1-s_{0})(1-\overline{s_{0}}).
\]
But $s_{1}\overline{s_{1}}=s_{0}\overline{s_{0}}$ since $|s_{1}|=|s_{0}|$, and the above relation reduces to
\[
s_{1}+\overline{s_{1}}=s_{0}+\overline{s_{0}};
\]
that is, $Re$ $s_{1}=Re$ $s_{0}$ and $Re_{s}s_{1}=Re_{s}s_{0}$. It follows that either
\[
s_{1}=s_{0}\texttt{ or }s_{1}=\overline{s_{0}},
\]
again since $|s_{1}|=|s_{0}|$. If $s_{1}=s_{0}$, transformation~(\ref{eq:chap_8_6_5}) becomes the constant function $\varpi=e_{xy}e^{i\alpha}\cos\beta+\sin\beta$ where $\alpha$ and $\beta$ are reals; hence $s_{1}=\overline{s_{0}}$.

Transformation~(\ref{eq:chap_8_6_5}), with $s_{1}=\overline{s_{0}}$, maps the point $s_{0}$ onto the origin $\varpi=0$; and, since points interior to the sphere $|\varpi|=1$ are to be the images of points above
the $xz$ coordinate plane, we may conclude that $Im$ $s>0$. Any linear fractional transformation having the mapping property stated in the first paragraph of this section must, therefore, be of the form
\begin{equation}\label{eq:chap_8_6_6}
\varpi=(e_{xy}e^{i\alpha}\cos\beta+\sin\beta)\frac{s-s_{0}}{s-\overline{s_{0}}}\texttt{ }(Im\texttt{ }s_{0}>0)
\end{equation}
where $\alpha$ and $\beta$ are reals.

It remains to show that, conversely, any linear fractional transformation of the
form~(\ref{eq:chap_8_6_6}) has the desired mapping property. This is easily done by taking absolute
values of each side of equation~(\ref{eq:chap_8_6_6}) and interpreting the resulting equation,
\[
|\varpi|=|\frac{s-s_{0}}{s-\overline{s_{0}}}|,
\]
geometrically. If a point $s$ lies above the $xz$ coordinate plane, both it and the point $s_{0}$ lie on the same side of that $x$ axis in the $xz$ coordinate plane, which is the perpendicular bisector of the line segment joining $s_{0}$ and $\overline{s_{0}}$. It follows that the distance $|s-s_{0}|$ is less than the distance $|s-\overline{s_{0}}|$; that is, $|\varpi|<1$. Likewise, if $s$ lies below the $xz$ coordinate plane, the distance $|s-s_{0}|$ is greater than the distance $|s-\overline{s_{0}}|$; and so $|\varpi|>1$.  Finally, if $s$ is in the $xz$ coordinate plane, $|\varpi|=1$ because then $|s-s_{0}|=|s-\overline{s_{0}}|$. Since any linear fractional transformation is a one to one mapping of the extended $s$ space onto the extended $\varpi$ space, this shows that \textit{transformation~(\ref{eq:chap_8_6_6}) maps the half space $Im$ $s>0$ onto the sphere $|\varpi|<1$ and the plane boundary of the half space onto the spherical boundary of the sphere}.

Our first example~(\ref{ex:chap_8_6_1}) here illustrates the use of the result in italics just above.

\begin{remark}\label{re:chap_8_6_1}
Note that when $\sin\beta=0$ in equation~(\ref{eq:chap_8_6_6}), let
\[
s=e_{xy}(x+iy)+z\texttt{ and }f(s)=\frac{s-s_{0}}{s-\overline{s_{0}}}\texttt{ }(Im\texttt{ }s_{0}>0).
\]
Then there are $f(s)=e_{xy}[f(x+iy+z)-f(z)]+f(z)$ and
\[
\begin{split}
\varpi=e_{xy}e^{i\alpha}f(s)&=e_{xy}[e^{i\alpha}f(x+iy+z)-e^{i\alpha}f(z)]+e_{xy}e^{i\alpha}f(z) \\
                            &=e_{xy}e^{i\alpha}f(x+iy+z)=e_{xy}[u(x,y,z)+v(x,y,z)].
\end{split}
\]
Thus transformation~(\ref{eq:chap_8_6_6}) maps the half space $Im\texttt{ }s>0$ onto the open circle $|\varpi|<1$ on the $uv$ coordinate plane in the $\varpi$ space, and maps the plane boundary $Im\texttt{ }s=0$ of the half space onto the curve boundary of the circle $|\varpi|=1$, and derives a pair of conjugate harmonic functions $u$ and $v$ respect to three real variables $x$, $y$ and $z$.
\end{remark}

\begin{example}\label{ex:chap_8_6_1}
The transformation
\begin{equation}\label{eq:chap_8_6_7}
\varpi=\frac{i-s}{i+s}
\end{equation}
can be written
\[
\varpi=\sin(-\frac{\pi}{2})\frac{s-i}{s-\overline{i}}.
\]
Hence it has the mapping property described in italics.

Images of the upper half space $Im$ $s>0$ under other types of linear fractional transformations are often fairly easy to determine by examining the particular transformation in question.
\end{example}

\begin{example}\label{ex:chap_8_6_2}
By writing
\[
s=e_{xy}(x+iy)+z\texttt{ and }\varpi=e_{xy}(u+iv)+w,
\]
we can readily show that the transformation
\begin{equation}\label{eq:chap_8_6_8}
\varpi=\frac{s-1}{s+1}
\end{equation}
maps the half space $y>0$ onto the half space $v>0$, the $x$ axis onto the $u$ axis, and the $z$ axis onto the $w$ axis. We first note that when the number $s$ is real, so is the number $\varpi$. Consequently, the image of the $xz$ coordinate plane or the plane $y=0$ is the $uw$ coordinate plane or the plane $v=0$. Furthermore, for any point $\varpi$ in the finite $\varpi$ space,
\[
v=Im\texttt{ }\varpi=Im\frac{(s-1)(\bar{s}+1)}{(s+1)(\bar{s}+1)}=\frac{2y}{|s+1|^{2}}\texttt{ }(s\neq-1).
\]
The numbers $y$ and $v$ thus have the same sign, and this means that points above the $xz$ coordinate plane correspond to points above the $uw$ coordinate plane and points below the $xz$ coordinate plane correspond to points below the $uw$ coordinate plane. Finally, since points in the $xz$ coordinate plane correspond to points
in the $uw$ coordinate plane and since a linear fractional transformation is a one to one mapping of the extended space onto the extended space (Sec.~(\ref{sec:chap_8_4_86})), the stated mapping property of transformation~(\ref{eq:chap_8_6_8}) is established.
\end{example}

Our final example involves a composite function and uses the mapping discussed in Example~(\ref{ex:chap_8_6_2}).

\begin{example}\label{ex:chap_8_6_3}
The transformation
\begin{equation}\label{eq:chap_8_6_9}
\varpi=Log\frac{s-1}{s+1}
\end{equation}
where the principal branch of the logarithmic function is used, is a composition of the functions
\begin{equation}\label{eq:chap_8_6_10}
\omega=\frac{s-1}{s+1}\texttt{ and }\varpi=\log\omega.
\end{equation}

We know from Example~(\ref{ex:chap_8_6_2}) that the first of transformations~(\ref{eq:chap_8_6_10}) maps the upper
half space $y>0$ onto the upper half space $Y>0$, where $s=e_{xy}(x+iy)+z$ and $\omega=e_{xy}(X+iY)+Z$.
Furthermore, it is easy to see that the second of transformations~(\ref{eq:chap_8_6_10}) maps the half space $Y>0$ onto the strip space $0<v<\pi$, where  $\varpi=e_{xy}(u+iv)+w$. More precisely, by writing $\omega=R(e_{xy}\exp(i\Theta)\cos\Phi+\sin\Phi)$ and
\[
\log\omega=\log R+\log\sin\Phi+e_{xy}(\log(\exp(i\Theta)\cos\Phi+\sin\Phi)-\log\sin\Phi)
\]
where $R>0$, $-\pi<\Theta<\pi$, $0\leq\Phi\leq\pi/2$ (Expression~(\ref{eq:chap_3_3_2b}), Sec.~(\ref{sec:chap_3_3_29})),
we see that as a point
\[
\omega=R(e_{xy}\exp(i\Theta)\cos\Phi+\sin\Phi)\texttt{ }(-\pi<\Theta<\pi,0\leq\Phi\leq\pi/2)
\]
moves outward from the origin along the ray $\Theta=\Theta_{0}$, its image is the point whose rectangular  coordinates in the $\varpi$ space are
\[
(\log(\cos\Theta\cos\Phi+\sin\Phi)-\log\sin\Phi, \log(\sin\Theta\cos\Phi+\sin\Phi), \log(R\sin\Phi)).
\]
That image evidently moves to the right along the entire length of the horizontal line $v=\log(\sin\Theta\cos\Phi+\sin\Phi)$. Since these lines fill the strip space $0<v<\pi$ as the choice of $\Theta_{0}$ varies between $0$ and $\pi$, the mapping of the half space $Y_{0}$ onto the strip space is, in fact, one to one.
\end{example}

This shows that the composition~(\ref{eq:chap_8_6_9}) of the mappings~(\ref{eq:chap_8_6_10}) transforms the space $y>0$ onto the strip space $0<v<\pi$.

\section{The Transformation $\varpi=\sin s$}\label{sec:chap_8_7_89}

Since (Expression~(\ref{eq:chap_3_7_01}) in Sec.~(\ref{sec:chap_3_7_33}))
\[
\sin s=e_{xy}[\sin(x+z)\cosh y+i\cos(x+z)\sinh y-\sin z]+\sin z\texttt{ }(|s|<+\infty),
\]
the transformation $\varpi=\sin s$ can be written
\begin{equation}\label{eq:chap_8_7_1}
u=\sin(x+z)\cosh y-\sin z,\texttt{ }v=\cos(x+z)\sinh y,\texttt{ }w=\sin z.
\end{equation}

One method that is often useful in finding images of regions under this transformation is to examine images of  vertical lines $x+z=c_{1}$ where $z=c_{z}$ in a plane $z=c_{z}$ which is parallel to the $xy$ coordinate plane. If $0<c_{1}<\pi/2$, points on the line $x+z=c_{1}$ in the plane $z=c_{z}$ are transformed into points on the curve in the plane $w=\sin z$ which is parallel to the $uv$ coordinate plane as follows
\begin{equation}\label{eq:chap_8_7_2}
u+w=\sin c_{1}\cosh y,\texttt{ }v=\cos c_{1}\sinh y,\texttt{ }w=\sin z\texttt{ }(-\infty<y<\infty),
\end{equation}
which is the right-hand branch of the hyperbola
\begin{equation}\label{eq:chap_8_7_3}
\frac{(u+w)^{2}}{\sin^{2}c_{1}}-\frac{v^{2}}{\cos^{2}c_{1}}=1,\texttt{ }w=\sin z\texttt{ }(-\infty<y<\infty,-\infty<z<\infty)
\end{equation}
with focus at the points
\[
\rho=\pm\sqrt{\sin^{2}c_{1}+\cos^{2}c_{1}}=\pm1.
\]
The second of equations~(\ref{eq:chap_8_7_2}) shows that as a point $(c_{1}-z,y,z)$ moves upward along the entire length of the line in the plane $z=c_{z}$, its image moves upward along the entire length of the hyperbola's branch in the plane $w=\sin z$. Note that, in particular, there is a one to one mapping of the top half $(y>0)$ of the line onto the top half $(v>0)$ of the hyperbola's branch. If $-\pi/2<c_{1}<0$, the line $x+z=c_{1}$ is mapped onto the left-hand branch of the same hyperbola.

The line $x=z=0$, or the $y$ axis, needs to be considered separately. According to equations~(\ref{eq:chap_8_7_1}), the image of each point $(0,y,0)$ is $(0,\sinh y,0)$. Hence the $y$ axis is mapped onto the $v$ axis in a one to one manner, the positive $y$ axis corresponding to the positive $v$ axis.

We now illustrate how these observations can be used to establish the images of certain regions.

\begin{example}\label{ex:chap_8_7_1}
Here we show that the transformation $\varpi=\sin s$ is a one to one mapping of the semi-infinite strip $-\pi/2\leq x+z\leq\pi/2$, $y\geq0$ in the $z=c_{z}$ plane of the $s$ space onto the upper half $v\geq0$ in each $w=\sin z$ plane of the $\varpi$ space.

To do this, we first show that the boundary of the strip is mapped in a one to one manner onto the real axis in  the $w=\sin z$ plane. The image of the line segment $BA$ there is found by writing $x+z=\pi/2$ in equations~(\ref{eq:chap_8_7_1}) and restricting $y$ to be nonnegative. Since $u=\cosh y$ and $v=0$ when $x+z=\pi/2$, a typical point $(\pi/2-c_{z},y,c_{z})$ on $BA$ is mapped onto the point $(\cosh y-\sin c_{z},0,\sin c_{z})$ in each $w=\sin z$ plane; and that image must move to the right from $B'$ along the $u$ axis as $(\pi/2-c_{z},y,c_{z})$ moves upward from $B$. A point $(x,0,c_{z})$ on the horizontal segment $DB$ has image $(\sin(x+z)-\sin c_{z},0,\sin c_{z})$, which moves to the right from $D'$ to $B'$ as $x+z$ increases from $x+z=-\pi/2$ to $x+z=\pi/2$, or as $(x,0,z)$ goes from $D$ to $B$. Finally, as a point $(-\pi/2-c_{z},y,c_{z})$ on the line segment $DE$ moves upward from $D$, its image $(-\cosh y-\sin c_{z},0,\sin c_{z})$ moves to the left from $D'$.

Now each point in the interior $-\pi/2<x+z<\pi/2$, $y>0$ of the strip lies on one of the vertical half lines $x+z=c_{l}$, $y>0$ $(-\pi/2<c_{l}<\pi/2)$. Also, it is important to notice that the images of those half lines are distinct and constitute the entire half plane $v>0$. More precisely, if the upper half $L$ of a line $x+z=c_{l}$ $(-\pi/2<c_{l}<\pi/2)$ is thought of as moving to the left toward the positive $y$ axis, the right-hand branch of the hyperbola containing its image $L'$ is opening up wider and its vertex $(\sin c_{1}-\sin c_{z},0,\sin c_{z})$ is tending toward the point $\varpi=(-\sin c_{z},0,\sin c_{z})$ in each $w=\sin z$ plane. Hence $L'$ tends to become the positive $v$ axis, which we saw just prior to this example is the image of the positive $y$ axis. On the other hand, as $L$ approaches the segment $BA$ of the boundary of the strip, the branch of the hyperbola closes down around the segment $B'A'$ of the $u$ axis and its vertex $(\sin c_{1}-\sin c_{z},0,\sin c_{z})$ tends toward the point $\varpi=(1-\sin c_{z},0,\sin c_{z})$. Similar statements can be made regarding the half lines. We may conclude that the image of each point in the interior of the strip lies in the upper half plane $v>0$ and, furthermore, that each point in the half plane is the image of exactly one point in the interior of the strip.

This completes our demonstration that the transformation $\varpi=\sin s$ is a one to
one mapping of the strip $-\pi/2\leq x+z\leq\pi/2$, $y\geq0$ in the $z=c_{z}$ plane of the $s$ space onto the upper half $v\geq0$ in each $w=\sin z$ plane of the $\varpi$ space where $-\pi/2\leq z\leq\pi/2$. The right-hand half of the strip in the $z=c_{z}$ plane of the $s$ space is evidently mapped onto the first quadrant of the $w=\sin z$ plane of the $\varpi$ space where $0\leq z\leq\pi/2$.
\end{example}

\begin{example}\label{ex:chap_8_7_2}
Here we show that the transformation $\varpi=\sin s$ is a one to one mapping of the semi-infinite tetragonal prism \[
-\pi/2\leq x\leq\pi/2,\texttt{ }y\geq0,\texttt{ }-\pi/2\leq z\leq\pi/2,\texttt{ and }-\pi/2\leq x+z\leq\pi/2,
\]
in the $s$ space onto the upper half $v\geq0$ in each $w=\sin z$ plane of the $\varpi$ space.

To do this, we shall first show that the boundaries of the tetragonal prism are mapped in a one to one manner onto the boundaries of the upper half $v\geq0$ in each $w=\sin z$ plane of the $\varpi$ space where $-\pi/2\leq z\leq\pi/2$.

First, the images of two surfaces $(\pm\pi/2-z,y,z)$ where
\[
0\leq z\leq\pi/2\texttt{ or }-\pi/2\leq z\leq0
\]
of the tetragonal prism are found by writing $x+z=\pm\pi/2$, respectively, in equations~(\ref{eq:chap_8_7_1}) and restricting $y$ to be nonnegative. Since $u+w=\pm\cosh y$ and $v=0$ when $x+z=\pm\pi/2$, typical points in the two surfaces of the prism are mapped onto the points $(\pm\cosh y-\sin z,0,\sin z)$ where
\[
x=\pi/2-z\texttt{ and }0\leq z\leq\pi/2\texttt{ or }x=-\pi/2-z\texttt{ and }-\pi/2\leq z\leq0
\]
on the two surfaces in the $uw$ coordinate plane of the $\varpi$ space, respectively.

Second, the images of two surfaces $(\mp\pi/2+z,y,z)$ where
\[
0\leq z\leq\pi/2\texttt{ or }-\pi/2\leq z\leq0
\]
of the tetragonal prism are found by writing $x+z=\mp\pi/2+2z$ where
\[
0\leq z\leq\pi/2\texttt{ or }-\pi/2\leq z\leq0,
\]
respectively, in equations~(\ref{eq:chap_8_7_1}) and restricting $y$ to be nonnegative. Since
\[
u+w=\sin(\mp\pi/2+2z)\cosh y\texttt{ and }v=\cos(\mp\pi/2+2z)\sinh y
\]
when $x+z=\mp\pi/2+2z$, typical points in the two surfaces of the prism are mapped onto the points
\[
(\sin(\mp\pi/2+2z)\cosh y-\sin z,\cos(\mp\pi/2+2z)\sinh y,\sin z)
\]
where
\[
x=-\pi/2+z\texttt{ and }0\leq z\leq\pi/2\texttt{ or }x=\pi/2+z\texttt{ and }-\pi/2\leq z\leq0
\]
on the two surfaces in the $\varpi$ space, respectively.

Finally, in the $xz$ coordinate plane of the $s$ space, as $z$ increases from $z=0$ to $z=\pi/2$ or from $z=-\pi/2$ to $z=0$, respectively, there are two line segments $(\pm\pi/2-z,0,z)$ and another two line segments $(\mp\pi/2+z,0,z)$. They form a continuous closed boundary of the bottom surface $S_{0}$ of the tetragonal prism. Their images $(\sin(x+z)-\sin z,0,\sin z)$ where
\[
0\leq z\leq\pi/2\texttt{ or }-\pi/2\leq z\leq0,
\]
form a continuous boundary of the image of $S_{0}$ in the $uw$ coordinate plane of the $\varpi$ space.

Now each point in the interior
\[
-\pi/2<x<\pi/2,\texttt{ }y>0,\texttt{ }-\pi/2<z<\pi/2,\texttt{ and }-\pi/2<x+z<\pi/2
\]
of the tetragonal prism lies on one of the half lines $(c_{x}, y>0, c_{z})$ where
\[
-\pi/2<c_{x}<\pi/2,\texttt{ }y>0,\texttt{ }-\pi/2<c_{z}<\pi/2,,\texttt{ and }-\pi/2<c_{x}+c_{z}<\pi/2,
\]
which are parallel to the $y$ axis. Also, it is important to notice that the images of those half lines are distinct and constitute the entire half  $v>0$ in each $w=\sin z$ plane where $-\pi/2\leq z\leq\pi/2$. More precisely, if the upper half $L$ of a line
\[
x=c_{x}\texttt{ and }z=c_{z}\texttt{ }(0<c_{x}<\pi/2,\texttt{ }0<c_{z}<\pi/2)
\]
is thought of as moving to the left toward the positive $y$ axis, the right-hand branch of the hyperbola containing its image $L'$ is opening up wider and its vertex $(\sin(c_{x}+c_{z})-\sin c_{z},0,\sin c_{z})$ is tending toward the point $\varpi=(-\sin c_{z},0,\sin c_{z})$. Hence $L'$ tends to become the positive $v$ axis, which we saw just prior to this example is the image of the positive $y$ axis. On the other hand, as $L$ approaches the surfaces of the boundary of the tetragonal prism, the branch of the hyperbola closes down around the segment of the $u$ axis and its vertex $(\sin(c_{x}+c_{z})-\sin c_{z},0,\sin c_{z})$ tends toward the point $\varpi=(1-\sin c_{z},0,\sin c_{z})$. Similar statements can be made regarding the half line on the left-hand side of the $y$ axis. We may conclude that the image of each point in the interior of the strip lies in the upper half $v>0$ in each $w=\sin z$ plane where $-\pi/2\leq z\leq\pi/2$ and, furthermore, that each point in the half plane is the image of exactly one point in the interior of the strip.

This completes our demonstration that the transformation $\varpi=\sin s$ is a one to one mapping of the tetragonal prism
\[
-\pi/2\leq x\leq\pi/2,\texttt{ }y\geq0,\texttt{ }-\pi/2\leq z\leq\pi/2,\texttt{ and }-\pi/2\leq x+z\leq\pi/2
\]
in the $s$ space onto the upper half $v\geq0$ in each $w=\sin z$ plane of the $\varpi$ space where $-\pi/2\leq z\leq\pi/2$.
\end{example}

Another convenient way to find the images of certain regions when $\varpi=\sin s$ is to consider the images of  horizontal line segments $y=c_{y}$ $(-\pi\leq x+c_{z}\leq\pi)$ in the plane $z=c_{z}$, where $c_{y}>0$. According to equations~(\ref{eq:chap_8_7_1}), the image of such a line segment is the curve with parametric representation
\begin{equation}\label{eq:chap_8_7_4}
u=\sin(x+c_{z})\cosh c_{y}-\sin c_{z},\texttt{ }v=\cos(x+c_{z})\sinh c_{y},\texttt{ }(-\pi\leq x+c_{z}\leq\pi).
\end{equation}
That curve is readily seen to be the ellipse
\begin{equation}\label{eq:chap_8_7_5}
\frac{(u+\sin c_{z})^{2}}{\cosh^{2}c_{y}}+\frac{v^{2}}{\sinh^{2}c_{y}}=1\texttt{ }(-\pi\leq x+c_{z}\leq\pi)
\end{equation}
whose focus lie at the points
\[
\varpi=\pm\sqrt{\cosh^{2}c_{y}-\sinh^{2}c_{y}}=\pm1.
\]
The image of a point $(x,c_{y},c_{z})$ moving to the right makes one circuit around the ellipse in the clockwise direction. Note that when smaller values of the positive number $c_{y}$ are taken, the ellipse becomes smaller but retains the same focus $(\pm1,0,c_{z})$. In the limiting case $c_{y}=0$, equations~(\ref{eq:chap_8_7_4}) become
\[
u=\sin(x+c_{z})-\sin c_{z},\texttt{ }v=0,\texttt{ }(-\pi\leq x+c_{z}\leq\pi);
\]
and we find that the interval $-\pi\leq x+c_{z}\leq\pi$ of the $x$ axis is mapped onto the interval
$-1<u+\sin c_{z}<1$ of the $u$ axis. The mapping is not, however, one to one, as it is when $c_{y}>0$.

Mappings by various other functions closely related to the sine function are easily obtained once mappings by the sine function are known.

\section{Mappings by $s^{2}$ and Branches of $s^{1/2}$}\label{sec:chap_8_8_90}

In Chap 2 (Sec.~(\ref{sec:chap_2_2_12})), we considered some fairly simple mappings under the transformation $\varpi=s^{2}$, written in the form
\begin{equation}\label{eq:chap_8_8_1}
u=x^{2}-y^{2}+2xz,\texttt{ }v=2(x+z)y,\texttt{ }w=z^{2}
\end{equation}
We turn now to a less elementary example and then examine related mappings $\varpi=s^{1/2}$, where specific branches of the square root function are taken.

\begin{example}\label{ex:chap_8_8_1}
Let us use equations~(\ref{eq:chap_8_8_1}) to show that the image of the vertical strip
$0\leq x\leq1$, $y\geq0$ in the plane $z=c_{z}$ which is parallel to the $xy$ coordinate plane, is the closed semi parabolic region indicated there.
When $0<x_{1}<1$, the point $(x_{1},y,c_{z})$ moves up a vertical half line, labeled $L_{1}$ and parallel to the $y$ axis, as $y$ increases from $y=0$. The image traced out in the plane $w=z^{2}$ which is parallel to the $uv$ coordinate plane has, according to equations~(\ref{eq:chap_8_8_1}), the parametric representation
\begin{equation}\label{eq:chap_8_8_2}
u=x_{1}^{2}-y^{2}+2x_{1}c_{z},\texttt{ }v=2(x_{1}+c_{z})y,\texttt{ }w=c_{z}^{2},\texttt{ }(0\leq y<\infty).
\end{equation}
Using the second of these equations to substitute for $y$ in the first one, we see that the image points $(u,v,c_{z}^{2})$ must lie on the parabola
\begin{equation}\label{eq:chap_8_8_3}
v^{2}=-4[u-x_{1}(x_{1}+2c_{z})](x_{1}+c_{z})^{2},
\end{equation}
with vertex at $(x_{1}^{2},0,c_{z}^{2})$ and focus at the origin. Since $v$ increases with $y$ from $v=0$, according to the second of equations~(\ref{eq:chap_8_8_2}), we also see that as the point $(x_{1},y,c_{z})$ moves
up $L_{1}$ from the $x$ axis, its image moves up the top half $L_{1}'$ of the parabola from the $u$ axis. Furthermore, when a number $x_{2}$ larger than $x_{1}$, but less than $1$, is taken, the corresponding half line $L_{2}$ has an image $L_{2}'$ that is a half parabola to the right of $L_{1}'$. We note, in fact, that the image  of the half line $AB$ is the top half of the parabola
\begin{equation}\label{eq:chap_8_8_3b}
v^{2}=-4[u-(1+2c_{z})](1+c_{z})^{2}
\end{equation}
when $x_{1}=1$, labeled $B'A'$.

The image of the half line $CD$ on the $y$ axis is found by observing from equations~(\ref{eq:chap_8_8_1}) that a typical point $(0,y,c_{z})$, where $y\geq0$, on $CD$ is transformed into the point $(-y^{2},2c_{z}y,c_{z}^{2})$ in the plane $w=c_{z}^{2}$. So, as a point moves up from the origin along $CD$, its image moves left from the origin along a curve labeled $C'D'$ and determined by equations $u=-y^{2}$ and $v=2c_{z}y$.% Evidently, then, as the vertical half lines in the plane $z=c_{z}$ move to the left, the half parabolas that are their images in the plane $w=z^{2}$ shrink down to become the half line $C'D'$.

It is now clear that the images of all the half lines between and including $CD$ and $BA$ fill up the closed semi parabolic region bounded by $A'B'C'D'$. Also, each point in that region is the image of only one point in the closed strip bounded by $ABCD$. Hence we may conclude that the semi parabolic region is the image of the strip and  that there is a one to one correspondence between points in those closed regions where the strip has arbitrary width.
\end{example}

As for mappings by branches of $s^{1/2}$, we recall from Sec.~(\ref{sec:chap_1_8}) that the values of $s^{1/2}$
are the two square roots of $s$ when $s\neq0$. According to that section, if spherical coordinates are used and
\[
s=r(e_{xy}e^{i\Theta}\cos\Phi+\sin\Phi)\texttt{ }(r>0, -\pi<\Theta<\pi,\texttt{ }0\leq\Phi\leq\pi/2),
\]
then
\begin{equation}\label{eq:chap_8_8_4}
\begin{split}
s^{1/2} &= r^{1/2}(e_{xy}e^{i(\Theta+2k\pi)}\cos\Phi+\sin\Phi)^{1/2} \\
        &= \sqrt{r}(e_{xy}e^{i\theta(k)}\cos\varphi+\sin\varphi)
\end{split}
\end{equation}
where $\varphi=\arcsin(\sqrt{\sin\Phi})$, $\theta(k)$ is determined by expression~(\ref{eq:chap_1_8_00}) in Sec.~(\ref{sec:chap_1_8}) for $k=0,1$, the principal root occurring when $k=0$. In Sec.~(\ref{sec:chap_3_5_31}), we saw that $s^{1/2}$ can  also be written
\begin{equation}\label{eq:chap_8_8_5}
s^{1/2}=\exp(\frac{1}{2}\log s)\texttt{ }(s\neq0).
\end{equation}
The principal branch $F_{0}(s)$ of the double-valued function  $s^{1/2}$ is then obtained by taking the principal branch of $\log s$ and writing (see Sec.~(\ref{sec:chap_3_6_32}))
\[
F_{0}(s)=\exp(\frac{1}{2}Log\texttt{ }s)\texttt{ }(|s|>0,-\pi<Arg_{c}s<\pi,0\leq Arg_{r}s\leq\pi/2).
\]
Since
\[
\frac{1}{2}Log\texttt{ }s=\ln\sqrt{r}+\log(e_{xy}e^{i\theta(k)}\cos\varphi+\sin\varphi)^{1/2}
\]
when $s=r(e_{xy}e^{i\Theta}\cos\Phi+\sin\Phi)$, this becomes
\begin{equation}\label{eq:chap_8_8_6}
F_{0}(s)=\sqrt{r}(e_{xy}e^{i\theta(0)}\cos\varphi+\sin\varphi)
\end{equation}
where $r>0$, $-\pi<Arg_{c}s<\pi$, and $0\leq Arg_{r}s\leq\pi/2$. The right-hand side of this equation is, of course, the same as the right-hand side of equation~(\ref{eq:chap_8_8_4}) when $k=0$ and $-\pi<\Theta<\pi$ there. The origin and the ray $\Theta=\pi$ form the branch cut for $F_{0}$, and the origin is the branch point.

Images of spatial curves and regions on or in spheres under the transformation $\varpi=F_{0}(s)$ may be obtained  by writing $\varpi=\rho(e_{xy}e^{i\theta(0)}\cos\varphi+\sin\varphi)$, where $\rho=\sqrt{r}$ and $\theta(0)\approx\Theta/2$. Arguments are evidently halved by this transformation, and it is understood that $\varpi=0$  when $s=0$.

Other branches of $s^{1/2}$ are obtained by using other branches of $\log s$ in expression~(\ref{eq:chap_8_8_5}).  A branch where the ray $\theta=\alpha$ is used to form the branch cut is given by the equation
\begin{equation}\label{eq:chap_8_8_8}
f_{\alpha}(s)=\sqrt{r}(e_{xy}e^{i\theta(k)}\cos\varphi+\sin\varphi)\texttt{ }(r>0,\alpha<\theta<\alpha+2\pi).
\end{equation}
Observe that when $\alpha=-\pi$, we have the branch $F_{0}(s)$ and that when $\alpha=\pi$, we have the branch $F_{1}(s)$. Just as in the case of $F_{0}(s)$, the domain of definition of $f_{\alpha}$ can be extended to the entire complex space by using expression~(\ref{eq:chap_8_8_8}) to define $f_{\alpha}$ at the nonzero points on the branch cut and by writing $f_{\alpha}(0)=0$. Such extensions are, however, never continuous in the entire complex space.

Finally, suppose that $n$ is any positive integer, where $n\geq2$. The values of $s^{1/n}$ are
the $n$th roots of $s$ when $s\neq0$; and, according to Sec.~(\ref{sec:chap_3_5_31}), the multiple-valued function $s^{1/n}$ can be written
\begin{equation}\label{eq:chap_8_8_9}
\begin{split}
s^{1/n} &= \exp(\frac{1}{n}\log s)=\sqrt[n]{r}(e_{xy}e^{i(\Theta+2k\pi)}\cos\Phi+\sin\Phi)^{1/n} \\
        &= \sqrt[n]{r}(e_{xy}e^{i\theta(k)}\cos\varphi+\sin\varphi)\texttt{ }(k=0,1,2,\ldots,n-1)
\end{split}
\end{equation}
where $r=|s|$, $\varphi=\arcsin(\sqrt{\sin\Phi})$, $\theta(k)$ is determined by expression~(\ref{eq:chap_1_8_00}) in Sec.~(\ref{sec:chap_1_8}) for $k=0,1,2,\ldots,n-1$. The case $n=2$ has just been considered. In the general case, each of the $n$ functions
\begin{equation}\label{eq:chap_8_8_10}
F_{k}(s)=\sqrt[n]{r}(e_{xy}e^{i\theta(k)}\cos\varphi+\sin\varphi)\texttt{ }(k=0,1,2,\ldots,n-1)
\end{equation}
is a branch of $s^{1/n}$, defined on the domain $r>0$, $-\pi<\Theta<\pi$, and $0\leq\Phi\leq\pi/2$. When $\varpi=\rho(e_{xy}e^{i\theta(k)}\cos\varphi+\sin\varphi)$, the transformation $\varpi=F_{k}(s)$ is a one to one mapping of that domain onto the domain
\[
\rho>0,\texttt{ }\frac{(2k-1)\pi}{n}<\Theta<\frac{(2k+1)\pi}{n},\texttt{ }0\leq\Phi\leq\pi/2.
\]
These $n$ branches of $s^{1/n}$ yield then distinct $n$th roots of $s$ at any point $s$ in the domain $r>0$, $-\pi<\Theta<\pi$, and $0\leq\Phi\leq\pi/2$. The principal branch occurs when $k=0$, and further branches of the type~(\ref{eq:chap_8_8_8}) are readily constructed.

\section{Square Roots of Polynomials}\label{sec:chap_8_9_91}

We now consider some mappings that are compositions of polynomials and square
roots of $s$ (Sec.~(\ref{sec:chap_1_8})).

\begin{example}\label{ex:chap_8_9_1}
Branches of the double-valued function $(s-s_{0})^{1/2}$ can be obtained by noting that it is a composition of the translation $S=s-s_{0}$ with the double-valued function $S^{1/2}$. Each branch of $S^{1/2}$
yields a branch of $(s-s_{0})^{1/2}$. When $S=e_{xy}(X+iY)+Z=R(e_{xy}e^{i\theta}\cos\varphi+\sin\varphi)$ where $R>0$ and $Z\geq0$, branches of $S^{1/2}$ are
\begin{equation}\label{eq:chap_8_9_0}
\begin{split}
S^{1/2} &= e_{xy}(\sqrt{X+iY+Z}-\sqrt{Z})+\sqrt{Z} \\
        &= \sqrt{R}(e_{xy}e^{i\theta_{b}(k)}\cos\varphi_{b}+\sin\varphi_{b})
\end{split}
\end{equation}
where $R>0$, $Z\geq0$, $0\leq\varphi\leq\pi/2$, $\varphi_{b}=\arcsin(\sqrt{Z/R})=\arcsin(\sqrt{\sin\varphi})$,
\[
\begin{split}
& r_{a}^{2}=1+\cos\theta\sin2\varphi, \\
& \theta_{a}=\arctan\frac{\sin\theta}{\cos\theta+\tan\varphi}, \\
& r_{b}^{2}(k)=1-2\cos(\frac{\theta_{a}}{2}+k\pi)\sqrt{\frac{\sin\varphi}{r_{a}}}+\frac{\sin\varphi}{r_{a}}, \\
& \theta_{b}(k)=\arctan\frac{\sin(\theta_{a}/2+k\pi)}{\cos(\theta_{a}/2+k\pi)-\sqrt{(\sin\varphi)/r_{a}}} \\
& \texttt{ }(k=0,1,\texttt{ }\alpha<\theta<\alpha+2\pi).
\end{split}
\]
Hence if we write
\[
R=|s-s_{0}|,\texttt{ }\Theta=Arg_{c}(s-s_{0}),\texttt{ }\theta=\arg_{c}(s-s_{0}),
\]
and
\[
\Phi=Arg_{r}(s-s_{0}),\texttt{ }\varphi=\arg_{r}(s-s_{0})\texttt{ }(0\leq\Phi\leq\pi/2,\texttt{ }0\leq\varphi_{b}\leq\pi/2),
\]
two branches of $(s-s_{0})^{1/2}$ are
\begin{equation}\label{eq:chap_8_9_1}
G_{0}(s)=\sqrt{R}(e_{xy}e^{i\theta_{b}(0)}\cos\varphi_{b}+\sin\varphi_{b})\texttt{ }(R>0,\texttt{ }0<\Theta<2\pi)
\end{equation}
and
\begin{equation}\label{eq:chap_8_9_2}
g_{0}(s)=\sqrt{R}(e_{xy}e^{i\theta_{b}(0)}\cos\varphi_{b}+\sin\varphi_{b})\texttt{ }(R>0,\texttt{ }0<\theta_{b}(0)<2\pi).
\end{equation}
\end{example}

The branch of $S^{1/2}$ that was used in writing $G_{0}(s)$ is defined at all points in the
$S$ space except for the origin and points on the ray $Arg_{c}S=\pi$. The transformation $\varpi=G_{0}(s)$ is, therefore, a one to one mapping of the domain
\[
|s-s_{0}|>0,\texttt{ }-\pi<Arg_{c}(s-s_{0})<\pi,\texttt{ }0\leq Arg_{r}(s-s_{0})\leq\pi/2
\]
onto the right half $Re\varpi>0$ of the $\varpi$ space. The transformation $\varpi=g_{0}(s)$ maps the domain
\[
|s-s_{0}|>0,\texttt{ }0<\arg_{c}(s-s_{0})<2\pi,\texttt{ }0\leq\arg_{r}(s-s_{0})\leq\pi/2
\]
in a one to one manner onto the upper half space $Im$ $\varpi>0$.

\begin{example}\label{ex:chap_8_9_2}
For an instructive but less elementary example, we now consider the double-valued function $(s^{2}-1)^{1/2}$. Using established properties of logarithms, we can write
\[
\begin{split}
(s^{2}-1)^{1/2} &= \exp[\frac{1}{2}\log(s^{2}-1)] \\
                &= \exp[\frac{1}{2}\log(s-1)]+\exp[\frac{1}{2}\log(s+1)]
\end{split}
\]
or
\begin{equation}\label{eq:chap_8_9_3}
(s^{2}-1)^{1/2}=(s-1)^{1/2}(s+1)^{1/2}\texttt{ }(s\neq\pm1).
\end{equation}
Thus, if $f_{1}(s)$ is a branch of $(s-1)^{1/2}$ defined on a domain $D_{1}$ and $f_{2}(s)$ is a branch
of $(s+1)^{1/2}$ defined on a domain $D_{2}$, the product $f(s)=f_{1}(s)f_{2}(s)$ is a branch of $(s^{2}-1)^{1/2}$ defined at all points lying in both $D_{1}$ and $D_{2}$.

In order to obtain a specific branch of $(s^{2}-1)^{1/2}$, we use the branch of $(s-1)^{1/2}$ and the branch of $(s+1)^{1/2}$ given by equation~(\ref{eq:chap_8_9_2}). If we write
\[
r_{1}=|s-1|,\texttt{ }\theta_{1}=\arg_{c}(s-1),\texttt{ and }\varphi_{1}=\arg_{r}(s-1),
\]
from equation~(\ref{eq:chap_8_9_0}), that branch of $(s-1)^{1/2}$ is
\[
f_{1}(s)=\sqrt{r_{1}r_{a,1}}r_{b,1}(k)(e_{xy}e^{i\theta_{b,1}(k)}\cos\varphi_{b,1}+\sin\varphi_{b,1})
\]
where $r_{1}>0$, $\alpha<\theta_{1}<\alpha+2\pi$, and $0\leq\varphi_{1}\leq\pi/2$. From equation~(\ref{eq:chap_8_9_0}), the branch of $(s+1)^{1/2}$ given by equation~(\ref{eq:chap_8_9_2}) is
\[
f_{2}(s)=\sqrt{r_{2}r_{a,2}}r_{b,2}(k)(e_{xy}e^{i\theta_{b,2}(k)}\cos\varphi_{b,2}+\sin\varphi_{b,2})
\]
where $r_{2}>0$, $\alpha<\theta_{2}<\alpha+2\pi$, $0\leq\varphi_{2}\leq\pi/2$, and
\[
r_{2}=|s+1|,\texttt{ }\theta_{2}=\arg_{c}(s+1),\texttt{ and }\varphi_{2}=\arg_{r}(s+1).
\]
The product of these two branches is, therefore, the branch $f$ of $(s^{2}-1)^{1/2}$ defined by the equation
\begin{equation}\label{eq:chap_8_9_4}
(s^{2}-1)^{1/2}=e_{xy}\sqrt{r_{1}r_{2}r_{a,1}r_{a,2}}r_{b,1}(k)r_{b,2}(k)
\end{equation}
\[
\times\{\exp[i(\theta_{b,1}(k)+\theta_{b,2}(k))]\cos\varphi_{b,1}\cos\varphi_{b,2}
+\exp[i\theta_{b,1}(k)]\cos\varphi_{b,1}\sin\varphi_{b,2}
\]
\[
+\exp[i\theta_{b,2}(k)]\cos\varphi_{b,2}\sin\varphi_{b,1}\}
+\sqrt{r_{1}r_{2}}\sin\varphi_{b,1}\sin\varphi_{b,2}
\]
where $r_{k}>0$, $\alpha<\theta_{k}<\alpha+2\pi$, and $0\leq\varphi_{k}\leq\pi/2$ for $k=1,2$. The branch $f$ is defined everywhere in the $s$ space except on the ray $r_{2}\geq0$ and $\theta_{2}=0$, which is the portion $x\geq-1$ of the $x$ axis.

The branch $f$ of $(s^{2}-1)^{1/2}$ given in equation~(\ref{eq:chap_8_9_4}) can be extended to a function
\begin{equation}\label{eq:chap_8_9_5}
F(s)=\sqrt{r_{1}r_{2}}(e_{xy}e^{i\theta_{1}/2}\cos\varphi_{1}+\sin\varphi_{1})(e_{xy}e^{i\theta_{2}/2}\cos\varphi_{2}+\sin\varphi_{2})
\end{equation}
where
\[
r_{k}>0,\texttt{ }\alpha<\theta_{k}<\alpha+2\pi,\texttt{ }0\leq\varphi_{k}\leq\pi/2\texttt{ }(k=1,2),\texttt{ and }r_{1}+r_{2}\geq2.
\]
As we shall now see, this function is analytic everywhere in its domain of definition, which is the entire $s$ space except for the segment $-1\leq x\leq1$ of the $x$ axis.

Since $F(s)=f(s)$ for all $s$ in the domain of definition of $F$ except on the ray $r_{1}>0$ and $\theta_{1}=0$, we need only show that $F$ is analytic on that ray. To do this, we form the product of the branches of $(s-1)^{1/2}$
and $(s+1)^{1/2}$ which are given by equation~(\ref{eq:chap_8_9_1}). That is, we consider the function
\[
G(s)=\sqrt{r_{1}r_{2}}(e_{xy}e^{i\Theta_{1}/2}\cos\Phi_{1}+\sin\Phi_{1})(e_{xy}e^{i\Theta_{2}/2}\cos\Phi_{2}+\sin\Phi_{2}),
\]
where
\[
r_{1}=|s-1|,\texttt{ }r_{2}=|s+1|,\texttt{ }\Theta_{1}=Arg_{c}(s-1),\texttt{ }\Theta_{2}=Arg_{c}(s+1)
\]
and
\[
\Phi_{1}=Arg_{r}(s-1),\texttt{ }\Phi_{2}=Arg_{r}(s+1)
\]
and where
\[
r_{k}>0,\texttt{ }-\pi<\Theta_{k}<\pi,\texttt{ }0\leq\Phi_{k}\leq\pi/2\texttt{ }(k=1,2).
\]
Observe that $G$ is analytic in the entire $s$ space except for the ray $r_{1}\geq0$, $\Theta_{1}=\pi$. Now $F(s)=G(s)$ when the point $s$ lies above or on the ray $r_{1}>0$, $\Theta_{1}=0$; for then $\theta_{k}=\Theta_{k}$ $(k=1,2)$. When $s$ lies below that ray, $\theta_{k}=\Theta_{k}+2\pi$ $(k=1,2)$. Consequently, $\exp(i\theta_{k}/2)=-\exp(i\Theta_{k}/2)$; and this means that
\[
\exp\frac{i(\theta_{1}+\theta_{2})}{2}=(\exp\frac{i\theta_{1}}{2})(\exp\frac{i\theta_{2}}{2})=\exp\frac{i(\Theta_{1}+\Theta_{2})}{2}
\]
So again, $F(s)=G(s)$. Since $F(s)$ and $G(s)$ are the same in a domain containing the ray $r_{1}>0$, $\Theta_{1}=0$ and since $G$ is analytic in that domain, $F$ is analytic there. Hence $F$ is analytic everywhere except on the line segment $-1\leq x\leq1$ of the $x$ axis.

The function $F$ defined by equation~(\ref{eq:chap_8_9_5}) cannot itself be extended to a function
which is analytic at points on the line segment $-1\leq x\leq1$ of the $x$ axis; for the value on the right in
equation~(\ref{eq:chap_8_9_5}) jumps from $i\sqrt{r_{1}r_{2}}$ to numbers near $-i\sqrt{r_{1}r_{2}}$ as the point $s$ moves downward across that line segment. Hence the extension would not even be continuous there.

The transformation $\varpi=F(s)$ is, as we shall see, a one to one mapping of the domain $D_{s}$ consisting of all points in the $s$ space except those on the line segment $-1\leq x\leq1$ of the $x$ axis onto the domain $D_{\varpi}$ consisting of the entire $\varpi$ space with the exception of the segment $-1\leq v\leq1$ of the $v$ axis.

Before verifying this, we note that if $s=e_{xy}iy$ $(y>0)$, then
\[
r_{1}=r_{2}>1\texttt{ and }\varphi_{1}+\varphi_{2}=\pi;
\]
hence the positive $y$ axis is mapped by $\varpi=F(s)$ onto that part of the $v$ axis for which $v>1$. The negative $y$ axis is, moreover, mapped onto that part of the $v$ axis for which $v<-1$. Each point in the upper half $y>0$ of the domain $D_{s}$ is mapped into the upper half $v>0$ of the $\varpi$ space, and each point in the lower half  $y<0$ of the domain $D_{s}$ is mapped into the lower half $v<0$ of the $\varpi$ space. The ray $r_{1}>0$ and $\theta_{1}=0$ is mapped onto the positive real axis in the $\varpi$ space, and the ray $r_{2}>0$, $\theta_{2}=\pi$ is mapped onto the negative real axis there.

To show that the transformation $\varpi=F(s)$ is one to one, we observe that if $F(s_{1})=F(s_{2})$, then $s_{1}^{2}-1=s_{2}^{2}-1$. From this, it follows that $s_{1}=s_{2}$ or $s_{1}=-s_{2}$. However, because of the manner in which $F$ maps the upper and lower halves of the domain $D_{s}$, as well as the portions of the real axis lying in $D_{s}$, the case $s_{1}=-s_{2}$ is impossible. Thus, if $F(s_{1})=F(s_{2})$, then $s_{1}=s_{2}$; and $F$ is one to one.

We can show that $F$ maps the domain $D_{s}$ onto the domain $D_{\varpi}$ by finding a function $H$ mapping $D_{\varpi}$ into $D_{s}$ with the property that if $s=H(\varpi)$, then $\varpi=F(s)$. This will show that, for any point $\varpi$ in $D_{\varpi}$, there exists a point $s$ in $D_{s}$ such that $F(s)=\varpi$; that is, the mapping $F$ is onto. The mapping $H$ will be the inverse of $F$.

To find $H$, we first note that if $\varpi$ is a value of $(s^{2}-1)^{1/2}$ for a specific $s$, then
$\varpi^{2}=s^{2}-1$; and $s$ is, therefore, a value of $(\varpi^{2}+1)^{1/2}$ for that $\varpi$. The function $H$ will be a branch of the double-valued function
\[
(\varpi^{2}+1)^{1/2}=(\varpi-i)^{1/2}(\varpi+i)^{1/2}\texttt{ }(\varpi\neq\pm i).
\]
Following our procedure for obtaining the function $F(s)$, we write
\[
\varpi-i=\rho_{1}(e_{xy}\exp(i\phi_{1})\cos\psi_{1}+\sin\psi_{1})
\]
and
\[
\varpi+i=\rho_{2}(e_{xy}\exp(i\phi_{2})\cos\psi_{2}+\sin\psi_{2}).
\]
With the restrictions
\[
\rho_{k}>0,\texttt{ }-\pi/2\leq\phi_{k}<3\pi/2,\texttt{ }0\leq\psi_{k}\leq\pi/2\texttt{ }(k=1,2)\texttt{ and }\rho_{1}+\rho_{2}>2,
\]
we then write
\begin{equation}\label{eq:chap_8_9_6}
H(\varpi)=\sqrt{\rho_{1}\rho_{2}}(e_{xy}e^{i\phi_{1}/2}\cos\psi_{1}+\sin\psi_{1})(e_{xy}e^{i\phi_{2}/2}\cos\psi_{2}+\sin\psi_{2}),
\end{equation}
the domain of definition being $D_{\varpi}$. The transformation $s=H(\varpi)$ maps points of $D_{\varpi}$ lying above or below the $uw$ coordinate plane onto points above or below the $xz$ coordinate plane, respectively. It maps the positive $u$ axis into that part of the $x$ axis where $x>1$ and the negative $u$ axis into that part of the negative $x$ axis where $x<-1$. If $s=H(\varpi)$, then $s^{2}=\varpi^{2}+1$; and so $\varpi^{2}=s^{2}-1$. Since $s$ is in $D_{s}$ and since $F(s)$ and $-F(s)$ are the two values of $(s^{2}-1)^{1/2}$ for a point in $D_{s}$, we see that $\varpi=F(s)$ or $\varpi=-F(s)$. But it is evident from the manner in which $F$ and $H$ map the upper and lower halves of their domains of definition, including the portions of the $xz$ coordinate plane lying in those domains, that $\varpi=F(s)$.

Mappings by branches  of  double-valued functions
\begin{equation}\label{eq:chap_8_9_7}
\varpi=(s^{2}+As+B)^{1/2}=[(s-s_{0})^{2}-s_{1}^{2}]^{1/2}\texttt{ }(s_{1}\neq0),
\end{equation}
where $A=-2s_{0}$ and $B=s_{0}^{2}=s_{1}^{2}$, can be treated with the aid of the results found for the function $F$ in Example~(\ref{ex:chap_8_9_2}) and the successive transformations
\begin{equation}\label{eq:chap_8_9_8}
S=\frac{s-s_{0}}{s_{1}},\texttt{ }W=(S^{2}-1)^{1/2},\texttt{ }\varpi=s_{1}W.
\end{equation}
\end{example}

%\section{Riemann Surfaces}\label{sec:chap_8_10_92}

%\section{Surfaces for Related Functions}\label{sec:chap_8_11_93}

%-----------------------------------------------------------------------
% Beginning of chap9.tex
%-----------------------------------------------------------------------
%
% AMS-LaTeX 1.2 sample file for a monograph, based on amsbook.cls.
% This is a data file input by chapter.tex.
%%%%%%%%%%%%%%%%%%%%%%%%%%%%%%%%%%%%%%%%%%%%%%%%%%%%%%%%%%%%%%%%%%%%

%\part{This is a Part Title Sample}

\chapter{Conformal Mapping}\label{ch:chap_9}

In this chapter, we introduce and develop the concept of a conformal mapping, with emphasis on connections between such mappings and harmonic functions. Applications to physical problems will follow in the next chapter.

\section{Preservation of Angles}\label{sec:chap_9_1_94}

Let $C$ be a smooth spatial arc (Sec.~(\ref{sec:chap_4_3_38})), represented by the equation
\[
s=s(t)\texttt{ }(a\leq t\leq b),
\]
and let $f(s)$ be a function defined at all points $s$ on $C$.  The equation
\[
\varpi=f[s(t)]\texttt{ }(a\leq t\leq b)
\]
is a parametric representation of the image $\Gamma$ of $C$ under the transformation $\varpi=f(s)$.

Suppose that $C$ passes through a point $s_{0}=s(t_{0})$ $(a\leq t_{0}\leq b)$ at which $f$ is
analytic and that $f'(s_{0})\neq0$. According to the chain rule, %given in Exercise 5, Sec.~(\ref{sec:chap_4_3_38}),
if $\varpi(t)=f[s(t)]$, then
\begin{equation}\label{eq:chap_9_1_1}
\varpi'(t_{0})=f'[s(t_{0})]s'(t_{0});
\end{equation}
and this means that (see Sec.~(\ref{sec:chap_1_7}))
\begin{equation}\label{eq:chap_9_1_2}
\arg\varpi'(t_{0})=(\arg_{c}\varpi'(t_{0}),\arg_{r}\varpi'(t_{0}))=(\phi_{0},\psi_{0})
\end{equation}
where by denoting
\[
\arg f'[s(t_{0})]=(\arg_{c}f'[s(t_{0})],\arg_{r}f'[s(t_{0})])=(\theta_{f},\varphi_{f})
\]
and
\[
\arg s'(t_{0})=(\arg_{c}s'(t_{0}),\arg_{r}s'(t_{0}))=(\theta_{0},\varphi_{0}),
\]
there are $\psi_{0}=\arcsin(\sin\varphi_{f}\sin\varphi_{0})$ and
\[
\phi_{0}=\arctan\frac{\sin(\theta_{f}+\theta_{0})
+\sin\theta_{f}\tan\varphi_{0}+\sin\theta_{0}\tan\varphi_{f}}
{\cos(\theta_{f}+\theta_{0})
+\cos\theta_{f}\tan\varphi_{0}+\cos\theta_{0}\tan\varphi_{f}}.
\]
Statement~(\ref{eq:chap_9_1_2}) is useful in relating the directions of $C$ and $\Gamma$ at the points $s_{0}$ and
$\varpi=f(s_{0})$, respectively.

To be specific, let $(\theta_{f},\varphi_{f})$ denote a pair of values of $\arg f'(s_{0})$, and let $(\theta_{0},\varphi_{0})$ be the angle of
inclination of a directed line tangent to $C$ at $s_{0}$. According to Sec.~(\ref{sec:chap_4_3_38}),$(\theta_{0},\varphi_{0})$ is a pair of values of $\arg s'(t_{0})$; and it follows from statement~(\ref{eq:chap_9_1_2}) that the quantities $(\phi_{0},\psi_{0})$ determined by $(\theta_{f},\varphi_{f})$ and $(\theta_{0},\varphi_{0})$ are a pair of values of $\arg\varpi'(t_{0})$ and is, therefore, the angle of inclination of a directed line tangent to $\Gamma$ at the point $\varpi_{0}=f(s_{0})$. Hence the angle of inclination of the directed line at $\varpi_{0}$ differs from the angle of inclination of the directed line at $s_{0}$ by the angle of rotation
\begin{equation}\label{eq:chap_9_1_3}
(\theta_{f},\varphi_{f})=\arg f'(s_{0}).
\end{equation}

Now let $C_{1}$ and $C_{2}$ be two smooth arcs passing through $s_{0}$,  and let $(\theta_{1},\varphi_{1})$ and $(\theta_{2},\varphi_{2})$ be angles of inclination of directed lines tangent to $C_{1}$ and $C_{2}$, respectively, at $s_{0}$. We know from the preceding paragraph that the quantities
\[
(\phi_{1},\psi_{1})\texttt{ and }(\phi_{2},\psi_{2})
\]
are angles of inclination of directed lines tangent to the image curves $\Gamma_{1}$ and $\Gamma_{2}$, which are determined by $(\theta_{f},\varphi_{f})$ and $(\theta_{1},\varphi_{1})$ or $(\theta_{f},\varphi_{f})$ and $(\theta_{2},\varphi_{2})$, respectively, at the point $\varpi_{0}=f(s_{0})$. Thus there are
\begin{equation}\label{eq:chap_9_1_4}
\frac{\sin\psi_{2}-\sin\psi_{1}}{\sin\varphi_{2}-\sin\varphi_{1}}=\frac{\sin\psi_{2}}{\sin\varphi_{2}}=\frac{\sin\psi_{1}}{\sin\varphi_{1}}=\sin\varphi_{f}\texttt{ or } \frac{\sin\psi_{2}}{\sin\psi_{1}}=\frac{\sin\varphi_{2}}{\sin\varphi_{1}},
\end{equation}
\begin{equation}\label{eq:chap_9_1_5}
\frac{(e^{i\phi_{2}}\cos\psi_{2}+\sin\psi_{2})-(e^{i\phi_{1}}\cos\psi_{1}+\sin\psi_{1})}
{(e^{i\theta_{2}}\cos\varphi_{2}+\sin\varphi_{2})-(e^{i\theta_{1}}\cos\varphi_{1}+\sin\varphi_{1})}
=e^{i\phi_{f}}\cos\psi_{f}+\sin\psi_{f}
\end{equation}
and
\begin{equation}\label{eq:chap_9_1_6}
\tan(\phi_{2}-\phi_{1})=\tan(\theta_{2}-\theta_{1})\frac{1+\alpha}{1+\beta}
\end{equation}
where
\[
\alpha=\alpha_{f}\frac{\sin\theta_{2}\tan\varphi_{1}-\sin\theta_{1}\tan\varphi_{2}}{\sin(\theta_{2}-\theta_{1})}
-\beta_{f}\frac{\cos\theta_{2}\tan\varphi_{1}-\cos\theta_{1}\tan\varphi_{2}}{\sin(\theta_{2}-\theta_{1})},
\]
\[
\beta=\alpha_{f}\frac{\cos\theta_{2}\tan\varphi_{1}+\cos\theta_{1}\tan\varphi_{2}}{\cos(\theta_{2}-\theta_{1})}
+\beta_{f}\frac{\sin\theta_{2}\tan\varphi_{1}+\sin\theta_{1}\tan\varphi_{2}}{\cos(\theta_{2}-\theta_{1})}
\]
\[
+\frac{\tan\varphi_{2}\tan\varphi_{1}}
{\cos(\theta_{2}-\theta_{1})(1+2\cos\theta_{f}\tan\varphi_{f}+\cot^{2}\varphi_{f})},
\]
\[
\alpha_{f}=\frac{1+\cos\theta_{f}\tan\varphi_{f}}{1+2\cos\theta_{f}\tan\varphi_{f}+\cot^{2}\varphi_{f}},\texttt{ and }
\beta_{f}=\frac{\sin\theta_{f}\tan\varphi_{f}}{1+2\cos\theta_{f}\tan\varphi_{f}+\cot^{2}\varphi_{f}},
\]
Equations~(\ref{eq:chap_9_1_4}-\ref{eq:chap_9_1_6}) show that the angle-preserving property of spatial complex functions is different from that of plane complex functions. Only if
\[
\varphi_{1}=\varphi_{2}=0\texttt{ or }\tan\varphi_{1}=\tan\varphi_{2}=0
\]
then there are $\alpha=\beta=0$ such that $C_{1}$ and $C_{2}$ are two smooth arcs in the $xy$ plane, and there are $\phi_{2}-\phi_{1}=\theta_{2}-\theta_{1}$; that is, the angle $\phi_{2}-\phi_{1}$ from $\Gamma_{1}$ to $\Gamma_{2}$ is the same in magnitude and sense as the angle $\theta_{2}-\theta_{1}$ from $C_{1}$ to $C_{2}$. This is the case of plane complex functions.

Because of this angle-preserving property of spatial complex functions, a transformation $\varpi_{0}=f(s_{0})$ is said to be conformal at a point $s_{0}$ if $f$ is analytic there and $f'(s_{0})\neq0$. Such a transformation is actually conformal at each point in a neighborhood of $s_{0}$. For $f$ must be analytic in a neighborhood of $s_{0}$ (Sec.~(\ref{sec:chap_2_13_23})); and, since $f'$ is continuous at $s_{0}$ (Sec.~(\ref{sec:chap_4_13_48})), it follows from Theorem~(\ref{th:chap_2_7_2}) in Sec.~(\ref{sec:chap_2_7_17}) that there is also a neighborhood of that point throughout which $f'(s)\neq0$.

A transformation $\varpi=f(s)$, defined on a domain $D$, is referred to as a conformal transformation, or conformal mapping, when it is conformal at each point in $D$. That is, the mapping is conformal in $D$ if $f$ is analytic in  $D$ and its derivative $f'$ has no zeros there. Each of the elementary functions studied in Chap.~(\ref{ch:chap_3}) can be used to define a transformation that is conformal in some domain.

A mapping that preserves the magnitude of the angle between two smooth arcs but not necessarily the sense is  called an isogonal mapping.

Suppose that $f$ is not a constant function and is analytic at a point $s_{0}$. If, in addition, $f'(s_{0})=0$, then $s_{0}$ is called a critical point of the transformation $\varpi=f(s)$.

More generally, it can be shown that if $s_{0}$ is a critical point of a transformation $\varpi=f(s)$, there is an integer $m$ $(m\geq2)$ such that the angle between any two smooth arcs passing through $s_{0}$ is multiplied by $m$ under that transformation. The integer $m$ is the smallest positive integer such that $f^{(m)}(s_{0})\neq0$.

\section{Scale Factors}\label{sec:chap_9_2_95}

Another property of a transformation $\varpi=f(s)$ that is conformal at a point $s_{0}$ is
obtained by considering the modulus of $f'(s_{0})$. From the definition of derivative and a property of limits involving moduli, %that was derived in Exercise 7, Sec. 17,
we know that
\begin{equation}\label{eq:chap_9_2_1}
f'(s_{0})=|\lim_{s \to s_{0}}\frac{f(s)-f(s_{0})}{s-s_{0}}|=\lim_{s \to s_{0}}\frac{|f(s)-f(s_{0})|}{|s-s_{0}|}.
\end{equation}
Now $|s-s_{0}|$ is the length of a line segment joining $s_{0}$ and $s$, and $|f(s)-f(s_{0})|$ is the
length of the line segment joining the points $f(s_{0})$ and $f(s)$ in the $\varpi$ space. Evidently,
then, if $s$ is near the point $s_{0}$, the ratio
\[
\frac{|f(s)-f(s_{0})|}{|s-s_{0}|}
\]
of the two lengths is approximately the number $|f'(s_{0})|$. Note that $|f'(s_{0})|$ represents an expansion if it is greater than unity and a contraction if it is less than unity.

Although the angle of rotation $\arg f'(s)$ (Sec.~(\ref{sec:chap_9_1_94})) and the scale factor $|f'(s)|$ vary, in general, from point to point, it follows from the continuity of $f'$ that their values are approximately $\arg f'(s_{0})$ and  $|f'(s_{0})|$ at points $s$ near $s_{0}$. Hence the image of a small region in a neighborhood of $s_{0}$ conforms  to the original region in the sense that it has approximately the same shape. A large region may, however, be transformed into a region that bears no resemblance to the original one.

\section{Local Inverses}\label{sec:chap_9_3_96}

A transformation $\varpi=f(s)$ that is conformal at a point $s_{0}$ has a local inverse there. That
is, if  $\varpi_{0}=f(s_{0})$, then there exists a unique transformation $s=g(\varpi)$, which is defined
and analytic in a neighborhood $N$ of $\varpi_{0}$, such that $g(\varpi_{0})=s_{0}$ and $f(g(\varpi))=\varpi$ for
all points $\varpi$ in $N$. The derivative of $g(\varpi)$ is, moreover,
\begin{equation}\label{eq:chap_9_3_1}
g'(\varpi)=\frac{1}{f'(s)}.
\end{equation}
We note from expression~(\ref{eq:chap_9_3_1}) that the transformation $s=g(\varpi)$ is itself conformal at
$\varpi_{0}$.

Assuming that $\varpi=f(s)$ is, in fact, conformal at $s_{0}$, let us verify the existence of such an inverse, which is a direct consequence of results in advanced calculus. As noted in Sec.~(\ref{sec:chap_9_1_94}), the conformality of the transformation $\varpi=f(s)$ at $s_{0}$ implies that there is some neighborhood of $s_{0}$ throughout which $f$ is analytic. Hence if we write
\[
s=e_{xy}(x+iy)+z,\texttt{ }s_{0}=e_{xy}(x_{0}+iy_{0})+z_{0},
\]
and
\[
f(s)=e_{xy}[u(x,y,z)+iv(x,y,z)]+w(z),
\]
we know that there is a neighborhood of the point $(x_{0},y_{0},z_{0})$ throughout which the
functions $u(x,y,z)$, $v(x,y,z)$, and $w(z)$ along with their partial derivatives of all orders, are
continuous (see Sec.~(\ref{sec:chap_4_13_48})).

Now the ternary group of equations
\begin{equation}\label{eq:chap_9_3_2}
u=u(x,y,z),\texttt{ }v=v(x,y,z),\texttt{ }w=w(z),
\end{equation}
represents a transformation from the neighborhood just mentioned into the $uvw$ space. Moreover, the determinant
\[
\mathbf{J}=\left |
    \begin{array}{ccc}
        u_{x} & u_{y} & u_{z} \\
        v_{x} & v_{y} & v_{z} \\
        0 & 0 & w_{z} \\
    \end{array}
\right |=(u_{x}v_{y}-v_{x}u_{y})w_{z},
\]
which is known as the Jacobian of the transformation, is nonzero at the point $(x_{0},y_{0},z_{0})$. For, in view of the Cauchy-Riemann equations $u_{x}=v_{y}$ and $u_{y}=-v_{x}$, one can write
$\mathbf{J}$ as
\[
\mathbf{J}=(u_{x}^{2}+v_{x}^{2})w_{z}=|f'(s)|^{2}w_{z},
\]
and $|f'(s)|^{2}w_{z}\neq0$ since the transformation $\varpi=f(s)$ is conformal at $s_{0}$. The above continuity conditions on the functions $u(x,y,z)$, $v(x,y,z)$, and $w(z)$ and their derivatives, together with this condition on the Jacobian, are sufficient to ensure the existence of a local inverse of transformation~(\ref{eq:chap_9_3_2}) at $(x_{0},y_{0},z_{0})$. That is, if
\begin{equation}\label{eq:chap_9_3_3}
u_{0}=u(x_{0},y_{0},z_{0}),\texttt{ }v_{0}=v(x_{0},y_{0},z_{0}),\texttt{ and }w_{0}=w(z_{0}),
\end{equation}
then there is a unique continuous transformation
\begin{equation}\label{eq:chap_9_3_4}
x=x(u,v,w),\texttt{ }y=y(u,v,w),\texttt{ }z=z(w),
\end{equation}
defined on a neighborhood $N$ of the point $(u_{0},v_{0},w_{0})$ and mapping that point onto $(x_{0},y_{0},z_{0})$, such that equations~(\ref{eq:chap_9_3_2}) hold when equations~(\ref{eq:chap_9_3_4}) hold. Also, in addition to being continuous, the functions~(\ref{eq:chap_9_3_4}) have continuous first-order partial derivatives satisfying
the equations
\begin{equation}\label{eq:chap_9_3_5}
x_{u}=\frac{1}{\mathbf{J}}v_{y}w_{z},\texttt{ }x_{v}=-\frac{1}{\mathbf{J}}u_{y}w_{z},\texttt{ }x_{w}=\frac{1}{\mathbf{J}}(u_{y}v_{z}-v_{y}u_{z}),\texttt{ }
\end{equation}
\[
y_{u}=-\frac{1}{\mathbf{J}}v_{x}w_{z},\texttt{ }y_{v}=\frac{1}{\mathbf{J}}u_{x}w_{z},\texttt{ }y_{w}=\frac{1}{\mathbf{J}}(v_{x}u_{z}-u_{x}v_{z}),\texttt{ }
\]
\[
z_{u}=0,\texttt{ }z_{v}=0,\texttt{ }z_{w}=w_{z}^{-1},
\]
throughout $N$ where
\[
\mathbf{J}^{-1}=\left |
    \begin{array}{ccc}
        x_{u} & x_{v} & x_{w} \\
        y_{u} & y_{v} & y_{w} \\
        z_{u} & z_{v} & z_{w} \\
    \end{array}
\right |=\left |
    \begin{array}{ccc}
        v_{y}w_{z} & -u_{y}w_{z} & u_{y}v_{z}-u_{z}v_{y} \\
        -v_{x}w_{z} & u_{x}w_{z} & v_{x}u_{z}-u_{x}v_{z} \\
        0 & 0 & w_{z}^{-1} \\
    \end{array}
\right |=\frac{1}{(u_{x}v_{y}-v_{x}u_{y})w_{z}}.
\]

If we write $\varpi=e_{xy}(u+iv)+w$ and $\varpi_{0}=e_{xy}(u_{0}+iv_{0})+w_{0}$, as well as
\begin{equation}\label{eq:chap_9_3_6}
g(\varpi)=e_{xy}[x(u,v,w)+iy(u,v,w)]+z(w),
\end{equation}
the transformation $s=g(\varpi)$ is evidently the local inverse of the original transformation
$g(\varpi)=f(s)$ at $s_{0}$. Transformations~(\ref{eq:chap_9_3_2}) and~(\ref{eq:chap_9_3_4}) can be written
\[
e_{xy}(u+iv)+w=e_{xy}[u(x,y,z)+iv(x,y,z)]+w(z)
\]
and
\[
e_{xy}(x+iy)+z=e_{xy}[x(u,v,w)+iy(u,v,w)]+z(w);
\]
and these last two equations are the same as
\[
\varpi=f(s)\texttt{ and }s=g(\varpi),
\]
where $g$ has the desired properties. Equations~(\ref{eq:chap_9_3_5}) can be used to show that $g$ is analytic
in $N$.% Details are left to the exercises, where expression~(\ref{eq:chap_9_3_1}) for $g'(\varpi)$ is also derived.

\section{Harmonic Conjugates}\label{sec:chap_9_4_97}

We saw in Sec.~(\ref{sec:chap_2_15_25}) that if a function
\[
f(s)=e_{xy}[u(x,y,z)+iv(x,y,z)]+w(z)
\]
is analytic in a domain $D$, then the real-valued functions $u$ and $v$ are harmonic in that domain. That is, they have continuous partial derivatives of the first and second order in $D$ and satisfy Laplace's equation there:
\begin{equation}\label{eq:chap_9_4_1}
u_{xx}+u_{yy}=0,\texttt{ }v_{xx}+v_{yy}=0.
\end{equation}
We had seen earlier that the first-order partial derivatives of $u$ and $v$ satisfy the Cauchy-Riemann equations
\begin{equation}\label{eq:chap_9_4_2}
e_{xy}u_{x}=e_{xy}v_{y},\texttt{ }e_{xy}v_{x}=e_{xy}v_{z}=-e_{xy}u_{y};
\end{equation}
and, as pointed out in Sec.~(\ref{sec:chap_2_15_25}), $v$ is called a harmonic conjugate of $u$.

Suppose now that $u(x,y,z)$ is any given harmonic function defined on a simply connected domain $D$ (Sec.~(\ref{sec:chap_4_11_46})). In this section, we show that $u(x,y,z)$ always has a harmonic conjugate $v(x,y,z)$ in $D$ by deriving an expression for $v(x,y,z)$.

To accomplish this, we first recall some important facts about line integrals in advanced calculus. Suppose that $p(x,y,z)$, $q(x,y,z)$, and $h(x,y,z)$ have continuous first-order partial derivatives in a simply connected domain $D$, and let $(x_{0},y_{0},z_{0})$ and $(x,y,z)$ be any two points in $D$. If $q_{z}=h_{y}$, $h_{x}=p_{z}$, and $p_{y}=q_{x}$ everywhere in $D$, then according to Stokes' Formula~(\ref{th:chap_5_1_3}) in Sec.~(\ref{sec:chap_5_1_51}) the line integra
\[
\int_{C}p(\xi,\eta,\zeta)d\xi+q(\xi,\eta,\zeta)d\eta+h(\xi,\eta,\zeta)d\zeta
\]
\[
=\iint_{S}(\frac{\partial h}{\partial\eta}-\frac{\partial q}{\partial\zeta})d\eta d\zeta
+(\frac{\partial p}{\partial\zeta}-\frac{\partial h}{\partial\xi})d\zeta d\xi
+(\frac{\partial q}{\partial\xi}-\frac{\partial p}{\partial\eta})d\xi d\eta
\]
from $(x_{0},y_{0},z_{0})$ to $(x,y,z)$ is independent of the contour $C$ that is taken as long as the
contour lies entirely in $D$ where $S$ is a surface bounded by $C$ and lies entirely in $D$. Furthermore, when the point $(x_{0},y_{0},z_{0})$ is kept fixed and is allowed to vary throughout $D$, the integral represents a single-valued function
\begin{equation}\label{eq:chap_9_4_3}
F(x,y,z)=\int_{(x_{0},y_{0},z_{0})}^{(x,y,z)}p(\xi,\eta,\zeta)d\xi+q(\xi,\eta,\zeta)d\eta+h(\xi,\eta,\zeta)d\zeta
\end{equation}
of $x$, $y$, and $z$ whose first-order partial derivatives are given by the equations
\begin{equation}\label{eq:chap_9_4_4}
F_{x}(x,y,z)=p(x,y,z),\texttt{ }F_{y}(x,y,z)=q(x,y,z),\texttt{ }F_{z}(x,y,z)=h(x,y,z).
\end{equation}
Note that the value of $F$ is changed by an additive constant when a different point
$(x_{0},y_{0},z_{0})$ is taken.

Returning to the given harmonic function $u(x,y,z)$, observe how it follows from Laplace's equation $u_{xx}+u_{yy}=0$ that
\[
(-u_{y})_{y}=(u_{x})_{x}
\]
everywhere in $D$. Also, the second-order partial derivatives of $u$ are continuous in $D$; and this means that the first-order partial derivatives of $-u_{y}$ and $u_{x}$ are continuous there. Thus, if $(x_{0},y_{0},z_{0})$ is a fixed point in $D$, the function
\begin{equation}\label{eq:chap_9_4_5}
v(x,y,z)=\int_{(x_{0},y_{0},z_{0})}^{(x,y,z)}-u_{\eta}(\xi,\eta,\zeta)d\xi+u_{\xi}(\xi,\eta,\zeta)d\eta-u_{\eta}(\xi,\eta,\zeta)d\zeta
\end{equation}
is well defined for all $(x,y,z)$ in $D$; and, according to equations~(\ref{eq:chap_9_4_4}),
\begin{equation}\label{eq:chap_9_4_6}
e_{xy}v_{x}=e_{xy}v_{z}=-e_{xy}u_{y},\texttt{ }e_{xy}v_{y}=e_{xy}u_{x}.
\end{equation}
These are the Cauchy-Riemann equations. Since the first-order partial derivatives of $u$ are continuous it is evident from equations~(\ref{eq:chap_9_4_6}) that those derivatives of $v$ are also continuous. Hence, according to Sec.~(\ref{sec:chap_2_11_21}), the function
\[
f(s)=e_{xy}[u(x,y,z)+iv(x,y,z)]+w(z)
\]
is an analytic function in $D$; and $v$ is, therefore, a harmonic conjugate of $u$.

The function $v$ defined by equation~(\ref{eq:chap_9_4_5}) is, of course, not the only harmonic conjugate of $u$. The function $v(x,y,z)+c$, where $c$ is any real constant, is also a harmonic conjugate of $u$.

\section{Transformations of Harmonic Functions}\label{sec:chap_9_5_98}

The problem of finding a function that is harmonic in a specified domain and satisfies prescribed conditions on the boundary of the domain is prominent in applied mathematics. If the values of the function are prescribed along the boundary, the problem is known as a boundary value problem of the first kind, or a Dirichlet problem. If the values  of the normal derivative of the function are prescribed on the boundary, the boundary value problem is one of the second kind, or a Neumann problem. Modifications and combinations of those types of boundary conditions also arise.

The domains most frequently encountered in the applications are simply connected; and, since a function that is harmonic in a simply connected domain always has a harmonic conjugate (Sec.~(\ref{sec:chap_9_4_97})), solutions of  boundary value problems for such domains are the real or imaginary parts of analytic functions.

Sometimes a solution of a given boundary value problem can be discovered  by identifying it as the real or two-dimensional complex parts of an analytic function. But the success  of that procedure depends on the simplicity of the problem and on one's familiarity with the real and two-dimensional complex parts of a variety of analytic functions. The following theorem is an important aid.

\begin{theorem}\label{th:chap_9_5_1}
Suppose that an analytic function
\begin{equation}\label{eq:chap_9_5_1}
\varpi=f(s)=e_{xy}[u(x,y,z)+iv(x,y,z)]+w(z)
\end{equation}
maps a domain $D_{s}$ in the $s$ space onto a domain $D_{\varpi}$ in the $\varpi$ space. If $h(u,v,w)$ is a harmonic function defined on $D_{\varpi}$, then the function
\begin{equation}\label{eq:chap_9_5_2}
H(x,y,z)=h[u(x,y,z),v(x,y,z),w(z)]
\end{equation}
is harmonic in $D_{s}$.
\end{theorem}

We first prove the theorem for the case in which the domain $D_{\varpi}$ is simply connected. According to Sec.~(\ref{sec:chap_9_4_97}), that property of $D_{\varpi}$ ensures that the given harmonic function $h(u,v,w)$ has a harmonic conjugate $g(u,v,w)$. Hence the function
\begin{equation}\label{eq:chap_9_5_3}
\Phi(\varpi)=h(u,v,w)+ig(u,v,w)
\end{equation}
is analytic in $D_{\varpi}$. Since the function $f(s)$ is analytic in $D_{s}$, the composite function
$\varpi(f(s))$ is also analytic in $D_{s}$. Consequently, the real part $h(u,v,w)$ of this
composition is harmonic in $D_{s}$.

If $D_{\varpi}$ is not simply connected, we observe that each point $\varpi_{0}$ in $D_{\varpi}$ has a
neighborhood $|\varpi-\varpi_{0}|<\varepsilon$ lying entirely in $D_{\varpi}$. Since that neighborhood is simply
connected, a function of the type~(\ref{eq:chap_9_5_3}) is analytic in it. Furthermore, since $f$ is continuous at a point $s_{0}$ in $D_{s}$ whose image is $\varpi_{0}$, there is a neighborhood $|s-s_{0}|<\delta$ whose image is contained in the neighborhood $|\varpi-\varpi_{0}|<\varepsilon$. Hence it follows that the composition $\Phi(f(s))$  is analytic in the neighborhood $|s-s_{0}|<\delta$, and we may conclude that $h[u(x,y,z),v(x,y,z),w(z)]$ is harmonic there. Finally, since $\varpi_{0}$ was arbitrarily chosen in $D_{\varpi}$ and since each point in $D_{s}$ is mapped onto such a point under the transformation $\varpi=f(s)$, the function $h[u(x,y,z),v(x,y,z),w(z)]$ must be harmonic throughout $D_{s}$.

The proof of the theorem for the general case in which $D_{\varpi}$ is not necessarily simply connected can also be accomplished directly by means of the chain rule for partial derivatives. The computations are, however, somewhat involved.% (see Exercise 8,  Sec. 99).

\section{Transformations of Boundary Conditions}\label{sec:chap_9_6_99}

The conditions that a function or its normal derivative have prescribed values along the boundary of a domain in which it is harmonic are the most common, although not the only, important types of boundary conditions. In this section, we show that certain of these conditions remain unaltered under the change of variables associated with a conformal transformation.These results wi11 be used in Chap.~(\ref{ch:chap_10}) to solve boundary value problems. The basic technique there is to transform a given boundary value problem in the $xyz$ space into a simpler one in the $uvw$ space and then to use the theorems of this and the preceding section to write the solution of the original problem in terms of the solution obtained for the simpler one.

\begin{theorem}\label{th:chap_9_6_1}
Suppose that a transformation
\begin{equation}\label{eq:chap_9_6_1}
\varpi=f(s)=e_{xy}[u(x,y,z)+iv(x,y,z)]+w(z)
\end{equation}
is conformal on a smooth spatial arc $C$, and let $\Gamma$ be the image of $C$ under that transformation. If,  along $\Gamma$, a function $h(u,v,w)$ satisfies either of the conditions
\begin{equation}\label{eq:chap_9_6_2}
h=h_{0}\texttt{ or }\frac{dh}{dn}=0,
\end{equation}
where $h_{0}$ is a real constant and $dh/dn$ denotes derivatives normal to $\Gamma$, then, along $C$, the function
\begin{equation}\label{eq:chap_9_6_3}
H(x,y,z)=h[u(x,y,z),v(x,y,z),w(z)]
\end{equation}
satisfies the corresponding condition
\begin{equation}\label{eq:chap_9_6_4}
H=h_{0}\texttt{ or }\frac{dH}{dN}=0,
\end{equation}
where $dH/dN$ denotes derivatives normal to $C$.
\end{theorem}

To show that the condition $h=h_{0}$ on $\Gamma$ implies that $H=h_{0}$ on $C$, we note from
equation~(\ref{eq:chap_9_6_3}) that the value of $H$ at any point $(x,y,z)$ on $C$ is the same as the value of
$h$ at the image $(u,v,w)$ of $(x,y,z)$ under transformation~(\ref{eq:chap_9_6_1}). Since the image point $(u,v,w)$ lies on $\Gamma$ and since $h=h_{0}$ along that curve, it follows that $H=h_{0}$ along $C$.

Suppose, on the other hand, that $dh/dn=0$ on $\Gamma$. From calculus, we know that
\begin{equation}\label{eq:chap_9_6_5}
\frac{dh}{dn}=(grad\texttt{ }h)\cdot\mathbf{n},
\end{equation}
where $grad$ $h$ denotes the gradient of $h$ at a point $(u,v,w)$ on $\Gamma$ and $\mathbf{n}$ is a unit vector
normal to $\Gamma$ at $(u,v,w)$. Since $dh/dn=0$ at $(u,v,w)$, equation~(\ref{eq:chap_9_6_5}) tells us that $grad$ $h$ is orthogonal to $\mathbf{n}$ at $(u,v,w)$. That is, $grad$ $h$ is tangent to $\Gamma$ there. But gradients are orthogonal to level curves; and, because $grad$ $h$ is tangent to $\Gamma$, we see that $\Gamma$ is orthogonal to a level curve $h(u,v,w)=c$ passing through $(u,v,w)$.

Now, according to equation~(\ref{eq:chap_9_6_3}), the level curve $H(u,v,w)=c$ in the $s$ space can be written
\[
h[u(x,y,z),v(x,y,z),w(z)]=c;
\]
and so it is evidently transformed into the level curve $h(u,v,w)=c$ under transformation~(\ref{eq:chap_9_6_1}).  Furthermore, since $C$ is transformed into $\Gamma$ and $\Gamma$ is orthogonal to the level curve $h(u,v,w)=c$, as demonstrated in the preceding paragraph, it follows from the conformality of transformation~(\ref{eq:chap_9_6_1}) on $C$ that $C$ is orthogonal to the level curve $H(x,y,z)=c$ at the point $(x,y,z)$ corresponding to $(u,v,w)$. Because gradients are orthogonal to level curves, this means that $grad$ $H$ is tangent to $C$ at $(x,y,z)$. Consequently, if $\mathbf{N}$ denotes a unit vector normal to $C$ at $(x,y,z)$, $grad$ $H$ is orthogonal to $\mathbf{N}$. That is,
\begin{equation}\label{eq:chap_9_6_6}
(grad\texttt{ }H)\cdot\mathbf{N}=0
\end{equation}
Finally, since
\[
\frac{dH}{dN}=(grad\texttt{ }H)\cdot\mathbf{N},
\]
we may conclude from equation~(\ref{eq:chap_9_6_6}) that $dH/dN=0$ at points on $C$.

In this discussion, we have tacitly assumed that $grad$ $h\neq0$. If $grad$ $h=0$, it follows from the identity
\[
|grad\texttt{ }H(x,y,z)|=|grad\texttt{ }h(u,v,w)||f'(s)|,
\]
%derived in Exercise lO(a) below,
that $grad$ $H=0$; hence $dh/dn$ and the corresponding normal derivative $dH/dN$ are both zero. We also assumed that

(i) $grad$ $h$ and $grad$ $H$ always exist;

(ii) the level curve $H(x,y,z)=c$ is smooth when $grad$ $h\neq0$ at $(u,v,w)$.

Condition (ii) ensures that angles between arcs are preserved by transformation~(\ref{eq:chap_9_6_1}) when it is conformal. In all  of our applications, both conditions (i) and (ii) will be satisfied.

A boundary condition that is not of one of the two types mentioned in the theorem may be transformed into a condition that is substantially different from the original one. %(see Exercise 6).
New boundary conditions for the transformed problem can be obtained for a particular transformation in any case. It is interesting to note that, under a conformal transformation, the ratio of a directional derivative of $H$ along a smooth arc $C$ in the $s$ space to the directional derivative of $h$ along the image curve $\Gamma$ at the corresponding point in the $\varpi$ space is $|f'(s)|$; usually, this ratio is not constant along a given arc.
%(See Exercise 10.)

\section{Transformations of Functions $f(x+iy+z)$}\label{sec:chap_9_7_100}

Let $\chi=x+z$ and $\sigma=x+iy+z=\chi+iy$. Then the function $\sigma(s)=\sigma=x+iy+z$ maps the $s$ space onto the $\sigma$ plane. More precisely, the function $\sigma(s)$ maps the $xz$ plane in the $s$ space onto the real $\chi$ axis in the $\sigma$ plane in a several-for-one way.

Let $\omega=f(\sigma)$. Then, when the function $f(s)$ maps a three-dimensional domain $D$ in the $s$ space onto the $\varpi$ space, the function $\sigma(s)$ maps the domain $D$ onto a two-dimensional region $R$ in the $\sigma$ plane in a several-for-one way, such that the function $f(\sigma)$ or $f(x+iy+z)$ maps the region $R$ onto the $\omega$ plane.

Let $\gamma$ and $\vartheta$ be spherical coordinates of the point $(x+z,y)$ or $(\chi,y)$ that corresponds to a nonzero plane complex number $\sigma=x+iy+z=\chi+iy$. Since
\[
x=r\cos\varphi\cos\theta\texttt{, }
y=r\cos\varphi\sin\theta\texttt{, and }
z=r\sin\varphi,
\]
the plane complex number $\sigma$ can be written in spherical form as
\begin{equation}\label{eq:chap_9_7_1}
\sigma=\gamma(\cos\vartheta+i\sin\vartheta)
\end{equation}
where
\[
\vartheta=\arctan\frac{y}{x+z}=\frac{\sin\theta}{\cos\theta+\tan\varphi}
\]
and
\[
\gamma=|\sigma|=\sqrt{(x+z)^{2}+y^{2}}=r\sqrt{1+\cos\theta\sin2\varphi}.
\]
So there are
\[
\chi=\gamma\cos\vartheta\texttt{ and }y=\gamma\sin\vartheta.
\]

If $\varphi=0$ then $\vartheta=\theta$ else if $\varphi=\pm\pi/2$ then $\vartheta=n\pi$ $(n=0,\pm1,\pm2,\ldots)$. So for $0\leq\theta\leq2\pi$ and $-\pi/2\leq\varphi\leq\pi/2$ there is $0\leq\vartheta\leq2\pi$.

Each value of $\vartheta$ is called an argument of $\sigma$, and the set of all such values is denoted by $\arg\sigma$. The principal value of $\arg\sigma$, denoted by $Arg$ $\sigma$, is unique, such that $-\pi<Arg$ $\sigma\leq\pi$ or $0\leq Arg$ $\sigma<2\pi$. Note that
\begin{equation}\label{eq:chap_9_7_2}
\arg\sigma=Arg\texttt{ }\sigma+2n\pi\texttt{ }(n=0,\pm1,\pm2,\ldots)
\end{equation}
where $Arg$ $\sigma=\arctan(y/\chi)=\arctan[y/(x+z)]$.

The symbol $e^{i\vartheta}$, or $\exp(i\vartheta)$, is defined by means of Euler's formula as
\begin{equation}\label{eq:chap_9_7_3}
e^{i\vartheta}=\cos\vartheta+i\sin\vartheta,
\end{equation}
where $\vartheta$ is to be measured in radians. It enables us to write the spherical form~(\ref{eq:chap_9_7_1}) more compactly in exponential form as
\begin{equation}\label{eq:chap_9_7_4}
\sigma=\gamma e^{i\vartheta}.
\end{equation}

\begin{example}\label{ex:chap_9_7_1}
When the domain $D$ is a closed sphere $|s|=r_{0}$ where $r_{0}$ is a positive constant, there are $0\leq|y|\leq r_{0}$, $0\leq\vartheta\leq2\pi$, and
\[
(\frac{x+z}{ar_{0}})^{2}+(\frac{y}{br_{0}})^{2}=\frac{1+\cos\theta\sin2\varphi}{a^{2}b^{2}}.
\]
Also there are
\[
\begin{split}
\chi=x+z=0\texttt{, }y=0\texttt{, }\vartheta=0\texttt{, and }\gamma=0&\texttt{ when }\theta=0\texttt{ and }\varphi=-\pi/4, \\
\chi=x+z=0\texttt{, }y=0\texttt{, }\vartheta=\pi\texttt{, and }\gamma=0&\texttt{ when }\theta=\pi\texttt{ and }\varphi=\pi/4,
\end{split}
\]
\[
\chi=x+z=\sqrt{2}r_{0}\texttt{, }y=0\texttt{, }\vartheta=0\texttt{, and }\gamma=\sqrt{2}r_{0}\texttt{ when }\theta=0\texttt{ and }\varphi=\pi/4,
\]
\[
\chi=x+z=-\sqrt{2}r_{0}\texttt{, }y=0\texttt{, }\vartheta=\pi\texttt{, and }\gamma=\sqrt{2}r_{0}\texttt{ when }\theta=\pi\texttt{ and }\varphi=-\pi/4,
\]
\[
\chi=x+z=0\texttt{, }y=\pm r_{0}\texttt{, }\vartheta=\pi/2\texttt{, and }\gamma=r_{0}\texttt{ when }\theta=\pi/2,
\]
\[
\chi=x=r_{0}\cos\theta\texttt{, }y=r_{0}\sin\theta\texttt{, }\vartheta=\theta\texttt{, and }\gamma=r_{0}\texttt{ when }\varphi=0.
\]
Then there is an ellipse equation
\begin{equation}\label{eq:chap_9_7_5}
(\frac{x+z}{\sqrt{2}r_{0}})^{2}+(\frac{y}{r_{0}})^{2}=1.
\end{equation}
So the function $\sigma(s)$ maps a closed sphere domain $D$ onto an ellipse region $R$.
\end{example}

\begin{example}\label{ex:chap_9_7_2}
When the domain $D$ is a closed circular cylinder
\[
s=r(e_{xy}e^{i\theta}\cos\varphi+\sin\varphi)=r_{0}(e_{xy}e^{i\theta}+\tan\varphi)
\]
where $r_{0}$ is a positive constant,
\[
-\pi/4\leq\varphi\leq\varphi_{0}\texttt{, and }\pi/4\leq\varphi_{0}<\pi/2,
\]
there are
\[
0\leq|y|\leq r_{0}\texttt{, }0\leq\vartheta\leq2\pi,
\]
\[
\gamma=|\sigma|=\sqrt{(x+z)^{2}+y^{2}}=r_{0}\sqrt{1+2\cos\theta\tan\varphi+\tan^{2}\varphi},
\]
and
\[
(\frac{x+z}{ar_{0}})^{2}+(\frac{y}{br_{0}})^{2}=\frac{1+2\cos\theta\tan\varphi+\tan^{2}\varphi}{a^{2}b^{2}}.
\]
Also there are
\[
\begin{split}
\chi=x+z=0\texttt{, }y=0\texttt{, }\vartheta=0\texttt{, and }\gamma=0&\texttt{ when }\theta=0\texttt{ and }\varphi=-\pi/4, \\
\chi=x+z=0\texttt{, }y=0\texttt{, }\vartheta=\pi\texttt{, and }\gamma=0&\texttt{ when }\theta=\pi\texttt{ and }\varphi=\pi/4,
\end{split}
\]
\[
\chi=x+z=2r_{0}\texttt{, }y=0\texttt{, }\vartheta=0\texttt{, and }\gamma=2r_{0}\texttt{ when }\theta=0\texttt{ and }\varphi=\pi/4,
\]
\[
\chi=x+z=-2r_{0}\texttt{, }y=0\texttt{, }\vartheta=\pi\texttt{, and }\gamma=2r_{0}\texttt{ when }\theta=\pi\texttt{ and }\varphi=-\pi/4,
\]
\[
\chi=z=r_{0}\tan\varphi\texttt{, }y=\pm r_{0}\texttt{, }\vartheta=\pi/2\texttt{, and }\gamma=r_{0}\sqrt{1+\tan^{2}\varphi}\texttt{ when }\theta=\pi/2,
\]
\[
\chi=x=r_{0}\cos\theta\texttt{, }y=r_{0}\sin\theta\texttt{, }\vartheta=\theta\texttt{, and }\gamma=r_{0}\texttt{ when }\varphi=0.
\]
When $\varphi=\varphi_{0}$, there are
\[
\\chi=x+z=r_{0}(\cos\vartheta+\tan\varphi_{0})\texttt{, }y=r_{0}\sin\vartheta\texttt{, }z=r_{0}\tan\varphi_{0}
\]
\[
\vartheta=\arctan\frac{\sin\vartheta}{\cos\vartheta+\tan\varphi_{0}}\texttt{, and }
\gamma=r_{0}\sqrt{(1+2\cos\vartheta\tan\varphi_{0}+\tan^{2}\varphi_{0})}.
\]
Then there is an ellipse equation
\begin{equation}\label{eq:chap_9_7_5}
(\frac{x+z}{r_{0}(1+\tan\varphi_{0})})^{2}+(\frac{y}{r_{0}})^{2}=1.
\end{equation}
So the function $\sigma(s)$ maps a closed circular cylinder domain $D$ onto an ellipse region $R$.
\end{example}

When the function $f(s)$ is analytic in the domain $D$, the function $f(\sigma)$ is also analytic in the region $R$. So the transformation $\omega=f(\sigma)$ is said to be \textit{conformal} at a point $\sigma_{0}$ in the region $R$ if $f$ is analytic there and $f'(\sigma_{0})\neq0$ according to the theory of the plane complex variables.

%-----------------------------------------------------------------------
% Beginning of chap10.tex
%-----------------------------------------------------------------------
%
% AMS-LaTeX 1.2 sample file for a monograph, based on amsbook.cls.
% This is a data file input by chapter.tex.
%%%%%%%%%%%%%%%%%%%%%%%%%%%%%%%%%%%%%%%%%%%%%%%%%%%%%%%%%%%%%%%%%%%%

%\part{This is a Part Title Sample}

\chapter{Applications of Conformal Mapping}\label{ch:chap_10}

We now use conformal mapping to solve a number of physical problems involving Laplace's equation in three independent variables. Problems in heat conduction, electrostatic potential, and fluid flow will be treated. Since these problems are intended to illustrate methods, they will be kept on a fairly elementary level.

\section{Steady Temperatures}\label{sec:chap_10_1_100}

In the theory of heat conduction, the $flux$ across a surface within a solid body at a point on that surface is the quantity of heat flowing in a specified direction normal to the surface per unit time per unit area at the point. Flux is, therefore, measured in such units as calories per second per square centimeter. It is denoted here by $\Phi$, and it varies with the normal derivative of the temperature $T$ at the point on the surface:
\begin{equation}\label{eq:chap_10_1_1}
\Phi=-K\frac{dT}{dN}\texttt{ }(K>0).
\end{equation}
Relation~(\ref{eq:chap_10_1_1}) is known as Fourier's law and the constant $K$ is called the thermal conductivity of the material of the solid, which is assumed to be homogeneous.

The points in the solid are assigned rectangular coordinates in three-dimensional space, and the temperature $T$  varies with the $x$, $y$, and $z$ coordinates. The flow of heat is, then, three-dimensional. We agree, moreover, that the flow is in a steady state; that is, $T$ does not vary with time.

It is assumed that no thermal energy is created or destroyed within the solid. That is, no heat sources or sinks are present there. Also, the temperature function $T(x,y,z)$ and its partial derivatives of the first and second order are continuous at each point interior to the solid. This statement and expression~(\ref{eq:chap_10_1_1}) for the flux of heat are postulates in the mathematical theory of heat conduction, postulates that also apply at points within a solid containing a continuous distribution of sources or sinks.

Consider now an element of volume that is interior to the solid and that has the shape of a cubic of the $xyz$ space, with base $\Delta x$ by $\Delta y$ by $\Delta z$ in that space. Heat enters or leaves the element only through these six faces $(\Delta x\Delta y$, $\Delta y\Delta z$, and $\Delta z\Delta x)$, and the temperatures within the element are steady. Hence there is
\begin{equation}\label{eq:chap_10_1_4}
\Delta T=T_{xx}+T_{yy}+T_{zz}=0.
\end{equation}
The temperature function $T(\alpha x, \beta y, \gamma z)$ thus satisfies Laplace's equation at each interior point of the solid (Theorem~(\ref{th:chap_2_15_3}) in Sec.~(\ref{sec:chap_2_15_25})).

In view of equation~(\ref{eq:chap_10_1_4}) and the continuity of the temperature function and its partial derivatives, $T(\alpha x, \beta y, \gamma z)$ is a harmonic function of three variables $x$, $y$, and $z$ in the domain representing the interior of the solid body.

The surfaces $T(\alpha x, \beta y, \gamma z)=c_{1}$, where $c_{1}$ is any real constant, are the isotherms within the solid. They can also be considered as surfaces in the $xyz$ space; then $T(\alpha x, \beta y, \gamma z)$ can be interpreted as the temperature at a point $(x,y,z)$ in a solid of material in that space. The isotherms are the level surfaces of the function $T$.

The gradient of $T$ is perpendicular to the isotherm at each point, and the maximum flux at a point is in the direction of the gradient there. If $T(\alpha x, \beta y, \gamma z)$ denotes temperatures in a solid and if $S$ is a harmonic conjugate of the function $T$, then a surface $S(\alpha x, \beta y, \gamma z)=c_{2}$ has the gradient of $T$ as a tangent vector at each point where the analytic function $T(\alpha x, \beta y, \gamma z)+iS(\alpha x, \beta y, \gamma z)$ is conformal. The surface $S(\alpha x, \beta y, \gamma z)=c_{2}$ are called surfaces of flow.

If the normal derivative $dT/dN$ is zero along any part of the boundary of the solid, then the flux of heat across that part is zero. That is, the part is thermally insulated and is, therefore, a surface of flow.

The function $T$ may also denote the concentration of a substance that is diffusing through a solid. In that case,  $K$ is the diffusion constant. The above discussion and the derivation of equation~(\ref{eq:chap_10_1_4}) apply as well to steady-state diffusion.

\section{Steady Temperatures in a Half Space}\label{sec:chap_10_2_101}

Let us find an expression for the steady temperatures $T(\alpha x, \beta y, \gamma z)$ in a semi-infinite block $y\geq0$ whose faces are insulated and whose edge $y=0$ is kept at temperature zero except for the area of a unit circle
\[
x=\cos\varphi\texttt{ and }z=\sin\varphi\texttt{ }(-\pi\leq\varphi\leq\pi),
\]
where it is kept at temperature unity. The function $T(\alpha x, \beta y, \gamma z)$ is to be bounded; this condition is natural if we consider the given block as the limiting case of the block $0\leq y\leq y_{0}$ whose upper edge is kept at a fixed temperature as $y_{0}$ is increased. In fact, it would be physically reasonable to stipulate that $T(\alpha x, \beta y, \gamma z)$ approach zero as $y$ tends to infinity.

The boundary value problem to be solved can be written
\begin{equation}\label{eq:chap_10_2_1}
\Delta T=T_{xx}+T_{yy}+T_{zz}=0\texttt{ }(-\infty<x<\infty,y>0,-\infty<z<\infty),
\end{equation}
\begin{equation}\label{eq:chap_10_2_2}
T(x,0,z)=\{
    \begin{array}{cc}
        1 & \texttt{ when }x^{2}+z^{2}<1 \\
        0 & \texttt{ when }x^{2}+z^{2}>1 \\
    \end{array},
\end{equation}
also, $|T(\alpha x, \beta y, \gamma z)|<M$ where $M$ is some positive constant. This is a Dirichlet problem for the upper half of the $xyz$ space. Our method of solution will be to obtain a new Dirichlet problem for a region in the $uvw$ space. That region will be the image of the half space under a transformation $\varpi=f(s)$ that is analytic in the domain $y>0$ and that is conformal along the boundary $y=0$ except at these points $(\cos\varphi,0,\sin\varphi)$ for $-\pi\leq\varphi\leq\pi$, where it is undefined. It will be a simple matter to discover a bounded harmonic function satisfying the new problem. The two theorems in Chap.~(\ref{ch:chap_9}) will then be applied to transform the solution of the problem in the $uvw$ space into a solution of the original problem in the $xyz$ space. Specifically, a harmonic function of three variables $u$, $v$, and $w$ will be transformed into a harmonic function of three variables $x$, $y$, and $z$, and the boundary conditions in the $uvw$ space will be preserved on corresponding portions of the boundary in the $xyz$ space. There should be no confusion if we use the same symbol $T$ to denote the different temperature functions in the two spaces.

Let
\[
s=e_{xy}(x+iy)+z=r(e_{xy}e^{i\theta}\cos\varphi+\sin\varphi)\texttt{ and }\tau_{1}=e_{xy}\sin\varphi_{1}+\cos\varphi_{1}.
\]
Then, the transformation
\begin{equation}\label{eq:chap_10_2_3}
\begin{split}
\varpi &=\frac{s-\tau_{1}}{s+\tau_{1}}=\frac{(s-\tau_{1})(\bar{s}+\tau_{1})}{(s+\tau_{1})(\bar{s}+\tau_{1})} \\
       &=\frac{(e_{xy}x+z)^{2}+e_{xy}y^{2}+e_{xy}i2y\tau_{1}-\tau_{1}^{2}}{(e_{xy}x+z+\tau_{1})^{2}+e_{xy}y^{2}}\\
       &=e_{xy}[u(x,y,z)+iv(x,y,z)]+w(z)
\end{split}
\end{equation}
where
\[
u(x,y,z)=Re\frac{(e_{xy}x+z)^{2}+e_{xy}y^{2}+e_{xy}i2y\tau_{1}-\tau_{1}^{2}}{(e_{xy}x+z+\tau_{1})^{2}+e_{xy}y^{2}}
\]
and
\[
\begin{split}
v(x,y,z)&=Im\frac{(e_{xy}x+z)^{2}+e_{xy}y^{2}+e_{xy}i2y\tau_{1}-\tau_{1}^{2}}{(e_{xy}x+z+\tau_{1})^{2}+e_{xy}y^{2}}\\
&=\frac{2y\tau_{1}}{(x+z+\tau_{1})^{2}+y^{2}}
\end{split}
\]
is defined on the upper half space $y\geq0$, except for these points $s_{1}=(\cos\varphi,0,\sin\varphi)$ for $-\pi\leq\varphi\leq\pi$. %The area of a unit circle $\tau_{1}=e_{xy}\sin\varphi_{1}+\cos\varphi_{1}$ for $-\pi\leq\varphi\leq\pi$ is mapped onto the upper edge; and the rest of the $xz$ coordinate plane is mapped onto the lower edge.
The required analyticity and conformality conditions are evidently satisfied by transformation~(\ref{eq:chap_10_2_3}).

A bounded harmonic function of three variables $u$, $v$, and $w$ that is zero on the edge $v=0$ is clearly
\begin{equation}\label{eq:chap_10_2_4}
T(\alpha x, \beta y, \gamma z)=v(\alpha x, \beta y, \gamma z);
\end{equation}
it is harmonic since it is the imaginary part of the entire function $\varpi$. Then $T(\alpha x, \beta y, \gamma z)$ satisfies Laplace's equation~(\ref{eq:chap_10_2_1})

\section{A Related Problem}\label{sec:chap_10_3_102}

Consider a semi-infinite tetragonal prism in the three-dimensional space bounded by two surfaces $(\pm\pi/2-z,y,z)$ denoted by $S_{1}$ and $S_{3}$ where
\[
0\leq z\leq\pi/2\texttt{ or }-\pi/2\leq z\leq0,
\]
respectively, and another two surfaces $(\mp\pi/2+z,y,z)$ denoted by $S_{2}$ and $S_{4}$ where
\[
0\leq z\leq\pi/2\texttt{ or }-\pi/2\leq z\leq0,
\]
respectively, and one bottom surface $S_{0}$ in the $xz$ coordinate plane of the $s$ space, whose continuous closed boundary is formed by two line segments $(\pm\pi/2-z,0,z)$ and another two line segments $(\mp\pi/2+z,0,z)$ as $z$ increases from $z=0$ to $z=\pi/2$ or from $z=-\pi/2$ to $z=0$, respectively.

Let the first pair of surfaces $S_{1}$ and $S_{3}$ satisfying $x+z=\pm\pi/2$ and $y>0$ be kept at temperature zero, the second pair of surfaces $S_{2}$ and $S_{4}$ satisfying $x+z=\pm\pi/2+2z$ and $y>0$ where
\[
0\leq z\leq\pi/2\texttt{ or }-\pi/2\leq z\leq0,
\]
respectively, be kept at temperature $\geq0$, respectively, and the surface $S_{0}$ in the $xz$ coordinate plane be kept at temperature zero. We wish to find a formula for the temperature $T(\alpha x, \beta y, \gamma z)$ at any interior point of the tetragonal prism. The problem is also that of finding temperatures in a block having the form of a semi-infinite prism
\[
-\pi/2\leq x\leq\pi/2,\texttt{ }y\geq0,\texttt{ and }-\pi/2\leq z\leq\pi/2
\]
when the faces of the tetragonal prism are insulated.

The boundary value problem here is
\begin{equation}\label{eq:chap_10_3_1}
T_{xx}+T_{yy}+T_{zz}=0\texttt{ }(-\pi/2<x<\pi/2,y>0,-\pi/2<z<\pi/2),
\end{equation}
\begin{equation}\label{eq:chap_10_3_2}
T(\pi/2-z,y,z)=0\texttt{ }(x=\pi/2-z,y>0,0\leq z\leq\pi/2),
\end{equation}
\begin{equation}\label{eq:chap_10_3_2b}
T(-\pi/2+z,y,z)\geq0\texttt{ }(x=-\pi/2+z,y>0,0\leq z\leq\pi/2),
\end{equation}
\begin{equation}\label{eq:chap_10_3_2c}
T(-\pi/2-z,y,z)=0\texttt{ }(x=-\pi/2-z,y>0,-\pi/2\leq z\leq0),
\end{equation}
\begin{equation}\label{eq:chap_10_3_2d}
T(\pi/2+z,y,z)\geq0\texttt{ }(x=\pi/2+z,y>0,-\pi/2\leq z\leq0),
\end{equation}
\begin{equation}\label{eq:chap_10_3_3}
T(x,0,z)=0\texttt{ }(-\pi/2\leq x+z\leq\pi/2,-\pi/2\leq z\leq\pi/2),
\end{equation}
where $T(\alpha x, \beta y, \gamma z)$ is bounded.

In view of Example~(\ref{ex:chap_8_7_2}) in Sec.~(\ref{sec:chap_8_7_89}), the mapping
\begin{equation}\label{eq:chap_10_3_4}
\varpi=\sin s=e_{xy}[u(x,y,z)+iv(x,y,z)]+w(z)
\end{equation}
is analytic where from formula~(\ref{eq:chap_3_7_07}) in Sec.~(\ref{sec:chap_3_7_33}) $w(z)=\sin z$,
\[
u(x,y,z)+w(z)=\sin(x+z)\cosh y,\texttt{ }v(x,y,z)=\cos(x+z)\sinh y.
\]

Let $\alpha=\gamma=1$, $\beta=\sqrt{2}$ and
\[
T(\alpha x, \beta y, \gamma z)=v(x,\beta y,z)=\cos(x+z)\sinh\beta y.
\]
Then $T(\alpha x, \beta y, \gamma z)$ satisfies Laplace's equation~(\ref{eq:chap_10_3_1}) and transforms this boundary value problem into the one posed in Sec.~(\ref{sec:chap_10_2_101}).

%\section{Temperatures  in  a Quadrant}\label{sec:chap_10_4_103}

\section{Electrostatic Potential}\label{sec:chap_10_5_104}

In an electrostatic force field, the field intensity at a point is a vector representing the force exerted on a unit positive charge placed at that point. The electrostatic potential is a scalar function of the space coordinates such that, at each point, its directional derivative in any direction is the negative of the component  of the field intensity in that direction.

For two stationary charged particles, the magnitude of the force of attraction or repulsion exerted by one particle  on the other is directly proportional to the product of the charges and inversely proportional to the square of  the distance between those particles. From this inverse-square law, it can be shown that the potential at a point
due to a single particle in space is inversely proportional to the distance between the point and the particle. In any region free of charges, the potential due to a distribution of charges outside that region can be shown to satisfy Laplace's equation for three-dimensional space.

In three-dimensional regions free of charges, $V(\alpha x, \beta y, \gamma z)$ is a harmonic function of three real variables $x$, $y$, and $z$ (Theorem~(\ref{th:chap_2_15_3}) in Sec.~(\ref{sec:chap_2_15_25})):
\begin{equation}\label{eq:chap_10_5_1}
\Delta V=V_{xx}+V_{yy}+V_{zz}=0.
\end{equation}
The potential function $V$ thus satisfies Laplace's equation at each interior point of the three-dimensional space. The field intensity vector at each point is with $x$, $y$, and $z$ components $-V_{x}(\alpha x, \beta y, \gamma z)$, $-V_{y}(\alpha x, \beta y, \gamma z)$, and $-V_{z}(\alpha x, \beta y, \gamma z)$, respectively. That vector is, therefore, the negative of the gradient of $V(\alpha x, \beta y, \gamma z)$.

A surface along which $V(\alpha x, \beta y, \gamma z)$ is constant is an equipotential surface. The tangential component of the field intensity vector at a point on a conducting surface is zero in the static case since charges are free to move on such a surface. Hence $V(\alpha x, \beta y, \gamma z)$ is constant along the surface of a conductor, and that surface is an  equipotential.

If $U$ is a harmonic conjugate of $V$, the surfaces $U(\alpha x, \beta y, \gamma z)=c_{2}$ in the space are called flux surfaces.  When such a surface intersects an equipotential surface $V(\alpha x, \beta y, \gamma z)=c_{1}$ at a curve where the derivative  of the analytic function $V(\alpha x, \beta y, \gamma z)+iU(\alpha x, \beta y, \gamma z)$ is not zero, the two surfaces are orthogonal at that curve and the field intensity is tangent to the flux surface there.

Boundary value problems for the potential $V$ are the same mathematical problems as those for steady temperatures  $T$. %The problem posed in Sec.~(\ref{sec:chap_10_3_102}), for instance, can be interpreted as that of finding the three-dimensional electrostatic potential in the empty space

%The potential in the steady flow of electricity in a plane conducting sheet is also a harmonic function at points free from sources and sinks. Gravitational potential is a further example  of  a harmonic function in physics.

\section{Potential in a Hollow Sphere}\label{sec:chap_10_6_105}

A hollow sphere is made out of a thin sheet of conducting material, and the sphere is split from the $xz$ coordinate plane to form two equal parts. Those parts are separated by annular strips of insulating material and are used as electrodes, one of which is grounded at potential zero and the other kept at a different fixed potential. We take the coordinate axes and units of length and potential difference as normal. We then interpret the electrostatic potential $V(x, y, z)$ over the enclosed space as a harmonic function inside the sphere $|s|=1$ or $x^{2}+y^{2}+z^{2}=1$ in the $s$ or $xyz$ space. Note that $V=0$ on the upper half of the sphere and that $V=1$ on the lower half.

A linear fractional transformation that maps the upper half plane onto the interior of the unit sphere centered at the origin, the right half $uw$ coordinate plane with the positive real $u$ axis onto the upper half of the
sphere, and the left half $uw$ coordinate plane with the negative real $-u$ axis onto the lower half of the sphere is verified in Example~(\ref{ex:chap_8_6_1}), Sec.~(\ref{sec:chap_8_6_88}). We find that the inverse of the transformation
\begin{equation}\label{eq:chap_10_6_1}
s=\frac{c-\varpi}{c+\varpi}
\end{equation}
gives us a new problem for $V$ in a half space where
\[
s=e_{xy}(x+iy)+z=r(e_{xy}e^{i\theta}\cos\varphi+\sin\varphi),
\]
\[
\varpi=e_{xy}(u+iv)+w=\rho(e_{xy}e^{i\phi}\cos\psi+\sin\psi),
\]
\[
c=e_{xy}e^{ia}\cos b+\sin b,
\]
$c$ is a spatial constant, $a$ and $b$ are reals, and $\sin b\neq0$ (equation~(\ref{eq:chap_8_6_6}), Sec.~(\ref{sec:chap_8_6_88})). Now from equation~(\ref{eq:chap_3_3_2b}), Sec.~(\ref{sec:chap_3_3_29}), the imaginary part of the function
\begin{equation}\label{eq:chap_10_6_2}
\log\varpi=\ln\rho+\ln\sin\psi+e_{xy}[\ln\Omega(\phi,\psi)+i\Phi(\phi,\psi)-\ln\sin\psi]
\end{equation}
where $\rho>0$, $0\leq\phi\leq\pi$, $0<\psi\leq\pi/2$,
\[
\Omega(\phi,\psi)=\sqrt{1+\cos\phi\sin2\psi},\texttt{ }
\Phi(\phi,\psi)=\arctan(\frac{\sin\phi}{\cos\phi+\tan\psi}),
\]
is a bounded function of $u$, $v$ and $w$ that assumes the required constant values on the two parts $\phi=0$ and $\phi=\pi$ of the $u$ axis. Hence the desired harmonic function for the half space is
\begin{equation}\label{eq:chap_10_6_3}
V=\frac{1}{\pi}Im\log\varpi=\frac{1}{\pi}\arctan(\frac{\sin\phi}{\cos\phi+\tan\psi})=\frac{1}{\pi}\arctan(\frac{v}{u+w}),
\end{equation}
where the values of the arc-tangent function range from $0$ to $\pi$. The inverse of transformation~(\ref{eq:chap_10_6_1}) is
\begin{equation}\label{eq:chap_10_6_4}
\varpi=c\frac{1-s}{1+s}=e_{xy}(u+iv)+w,
\end{equation}
from which $u$, $v$ and $w$ can be expressed in terms of $x$, $y$, and $z$ as follows:
\[
u=\frac{[1-(x+z)^{2}-y^{2}](\cos a\cos b+\sin b)+2y\sin a\cos b}{(x+z+1)^{2}+y^{2}}-\frac{1-z}{1+z}\sin b
\]
\[
v=\frac{[1-(x+z)^{2}-y^{2}]\sin a\cos b-2y(\cos a\cos b+\sin b)}{(x+z+1)^{2}+y^{2}},\texttt{ }w=\frac{1-z}{1+z}\sin b
\]
Equation~(\ref{eq:chap_10_6_3}) then becomes
\begin{equation}\label{eq:chap_10_6_5}
V=\frac{1}{\pi}\arctan\{\frac{[1-(x+z)^{2}-y^{2}]\sin a-2y(\cos a+\tan b)}{[1-(x+z)^{2}-y^{2}](\cos a+\tan b)+2y\sin a}\}
\end{equation}
\[
\texttt{ }(0\leq t\leq\pi).
\]
The function~(\ref{eq:chap_10_6_5}) is the potential function for the space enclosed by the spherical electrodes since it is harmonic inside the sphere and assumes the required values on the hemispheres. If we wish to verify this solution, we must note that
\[
\lim_{t \to 0,\texttt{ }t>0}t=0\texttt{ and }\lim_{t \to 0,\texttt{ }t>0}t=\pi.
\]
The equipotential spheres $V(x, y, z)=c_{1}$ $(0<c_{1}<1)$ in the spherical region are surfaces of the spheres
\[
x^{2}+(y+\tan\pi c_{1})^{2}+z^{2}=\sec\pi c_{1}.
\]
%with the sphere passing through the points $(\pm1,0,0)$. Also, the segment of the $x$ axis between those points is the equipotential $V(x, y, z)=1/2$.
A harmonic conjugate $U$ of $V$ is $(1/\pi)[\ln\Omega(\phi,\psi)-\ln\sin\psi]$, or $(1/\pi)Re\log\varpi$. In view of equation~(\ref{eq:chap_10_6_4}), $U$ may be written
\[
\begin{split}
U &= \frac{1}{\pi}[\frac{1}{2}\ln(1+\cos\phi\sin2\psi)-\ln\sin\psi] \\
  &= \frac{1}{\pi}\{\frac{1}{2}\ln[(u+w)^{2}+v^{2}]-\ln w\} \\
  &= \frac{1}{\pi}[\ln(\frac{\sqrt{1+\cos a\sin2 b}}{\sin b})
+\ln(\frac{\sqrt{[1-(x+z)^{2}-y^{2}]^{2}+4y^{2}}}{(x+z+1)^{2}+y^{2}})
-\ln(\frac{1-z}{1+z})].
\end{split}
\]
%From this equation, it can be seen that the flux surfaces $U(x, y, z)=c_{2}$ are surfaces of the spheres with centers on the $x$ axis. The segment of the $y$ axis between the electrodes is also a flux surface.

\section{Three-dimensional Fluid Flow}\label{sec:chap_10_7_106}

Harmonic functions play an important role in hydrodynamics and aerodynamics. Again, we consider only the three-dimensional steady-state type of problem. That is, the motion of the fluid is assumed to be independent of time. We let the vector% representing the spatial complex number
\[
V=(p(x,y,z),q(x,y,z),h(x,y,z))
\]
denote the velocity of a particle of the fluid at any point $(x,y,z)$; hence the $x$, $y$, and $z$ components of the velocity vector are $p(x,y,z)$, $q(x,y,z)$, and $h(x,y,z)$, respectively. At points interior to a region of flow in which no sources or sinks of the fluid occur, the real-valued functions $p(x,y,z)$, $q(x,y,z)$, and $h(x,y,z)$ and their first-order partial derivatives are assumed to be continuous.

The circulation of the fluid along any contour $C$ is defined as the line integral with respect to arc length $\sigma$ of the tangential component $V_{T}(x,y,z)$ of the velocity vector along $C$:
\begin{equation}\label{eq:chap_10_7_1}
\int_{C}V_{T}(x,y,z)d\sigma.
\end{equation}
The ratio of the circulation along $C$ to the length of $C$ is, therefore, a mean speed of the fluid along that contour. It is shown in advanced calculus that such an integral can be written
\begin{equation}\label{eq:chap_10_7_2}
\int_{C}V_{T}(x,y,z)d\sigma=\int_{C}p(x,y,z)dx+q(x,y,z)dy+h(x,y,z)dz.
\end{equation}
When $C$ is a positively oriented simple closed contour lying in a simply connected domain of flow containing no sources or sinks, Stokes' Formula~(\ref{th:chap_5_1_3}) in Sec.~(\ref{sec:chap_5_1_51}) enables us to write
\[
\int_{C}p(x,y,z)dx+q(x,y,z)dy+h(x,y,z)dz
\]
\[
=\iint_{S}(\frac{\partial h}{\partial y}-\frac{\partial q}{\partial z})dydz
+(\frac{\partial p}{\partial z}-\frac{\partial h}{\partial x})dzdx
+(\frac{\partial q}{\partial x}-\frac{\partial p}{\partial y})dxdy
\]
where the orientation of the boundary $C$ is chosen consistently with the orientation of the surface $S$. Thus
\begin{equation}\label{eq:chap_10_7_3}
\int_{C}V_{T}(x,y,z)d\sigma=\iint_{S}(\frac{\partial h}{\partial y}-\frac{\partial q}{\partial z})dydz
+(\frac{\partial p}{\partial z}-\frac{\partial h}{\partial x})dzdx
+(\frac{\partial q}{\partial x}-\frac{\partial p}{\partial y})dxdy
\end{equation}
for such a contour.

A physical interpretation of the integrand on the right in expression~(\ref{eq:chap_10_7_3}) for the circulation along the simple closed contour $C$ is readily given. We let $C$ denote a circle on a sphere of radius $r$, which is centered at a point $(x_{0}, y_{0}, z_{0})$ and taken counterclockwise. The mean speed along $C$ is then found by dividing the circulation by the circumference $2\pi r$, and the corresponding mean angular speed of the fluid about the center of the circle is obtained by dividing that mean speed by $r$:
\[
\frac{1}{\pi r}\iint_{S}\frac{1}{2}(\frac{\partial h}{\partial y}-\frac{\partial q}{\partial z})dydz
+\frac{1}{2}(\frac{\partial p}{\partial z}-\frac{\partial h}{\partial x})dzdx
+\frac{1}{2}(\frac{\partial q}{\partial x}-\frac{\partial p}{\partial y})dxdy.
\]
Now this is also an expression for the mean value of the function
\begin{equation}\label{eq:chap_10_7_4}
\omega_{x}=\frac{1}{2}(\frac{\partial h}{\partial y}-\frac{\partial q}{\partial z}),\texttt{ }
\omega_{y}=\frac{1}{2}(\frac{\partial p}{\partial z}-\frac{\partial h}{\partial x}),\texttt{ }
\omega_{z}=\frac{1}{2}(\frac{\partial q}{\partial x}-\frac{\partial p}{\partial y})
\end{equation}
over the circular surface $S$ bounded by $C$. Its limit as $r$ tends to zero is the value of
$\omega=(\omega_{x},\omega_{y},\omega_{z})$ at the point $(x_{0}, y_{0}, z_{0})$. Hence the function $\omega(x,y,z)$ and its components, called the rotation of the fluid around the $x$, $y$, and $z$ axis, respectively, represents the limiting angular speed of a circular element of the fluid as the sphere shrinks to its center $(x,y,z)$, the point at which $\omega$ is evaluated.

If $\omega(x,y,z)=(0,0,0)$ at each point in some simply connected domain, the flow is irrotational in that domain. We consider only irrotational flows here, and we also assume that the fluid is incompressible and free from viscosity. Under our assumption of steady irrotational flow of fluids with uniform density $\rho$, it can be shown that the fluid pressure $P(x,y,z)$ satisfies the following special case of Bernoulli's equation:
\[
\frac{P}{\rho}+\frac{1}{2}|V|^{2}=constant.
\]
Note that the pressure is greatest where the speed $|V|$ is least.

Let $D$ be a simply connected domain in which the flow is irrotational. According to equation~(\ref{eq:chap_10_7_4}), $q_{z}=h_{y}$, $h_{x}=p_{z}$, and $p_{y}=q_{x}$ throughout $D$. These relations between partial derivatives imply that the line integral
\[
\int_{C}p(x,y,z)dx+q(x,y,z)dy+h(x,y,z)dz
\]
along a contour $C$ lying entirely in $D$ and joining any two points $(x_{0},y_{0},z_{0})$ and $(x,y,z)$ in $D$ is actually independent of path. Thus, if $(x_{0},y_{0},z_{0})$ is fixed, the function
\begin{equation}\label{eq:chap_10_7_5}
\phi(x,y,z)=\int_{(x_{0},y_{0},z_{0})}^{(x,y,z)}p(\xi,\eta,\zeta)d\xi+q(\xi,\eta,\zeta)d\eta+h(\xi,\eta,\zeta)d\zeta.
\end{equation}
is well defined on $D$; and, by taking partial derivatives on each side of this equation, we find that
\begin{equation}\label{eq:chap_10_7_6}
\phi_{x}(x,y,z)=p(x,y,z),\texttt{ }\phi_{y}(x,y,z)=q(x,y,z),\texttt{ }\phi_{z}(x,y,z)=h(x,y,z).
\end{equation}

From equations~(\ref{eq:chap_10_7_6}), we see that the velocity vector $V=(p,q,h)$ is the gradient of
$\phi$; and the directional derivative of $\phi$ in any direction represents the component of the velocity of flow in that direction.

The function $\phi(x,y,z)$ is called the velocity potential. From equation~(\ref{eq:chap_10_7_5}), it is evident that $\phi(x,y,z)$ changes by an additive constant when the reference point $(x_{0},y_{0},z_{0})$ is changed. The level curves $\phi(x,y,z)=c_{1}$ are called equipotentials. Because it is the gradient of $\phi(x,y,z)$, the velocity vector $V$ is normal to an equipotential at any point where $V$ is not the zero vector.

%Just as in the case of the flow of heat,
The condition that the incompressible fluid enter or leave an element of volume only by flowing through the boundary of that element requires that $\phi(x,y,z)$ must satisfy Laplace's equation
\[
\phi_{xx}(x,y,z)+\phi_{yy}(x,y,z)+\phi_{zz}(x,y,z)=0
\]
in a domain where the fluid is free from sources or sinks. In view of equations~(\ref{eq:chap_10_7_6}) and the continuity of the functions $p$, $q$ and $h$ and their first-order partial derivatives, it follows that the partial derivatives of the first and second order of $\phi$ are continuous in such a domain. Hence the velocity potential $\phi$ is a harmonic function in that domain.

In three-dimensional regions where the fluid is free from sources or sinks, $\phi(\alpha x, \beta y, \gamma z)$ is a harmonic function of three real variables $x$, $y$, and $z$ (Theorem~(\ref{th:chap_2_15_3}) in Sec.~(\ref{sec:chap_2_15_25})):
\begin{equation}\label{eq:chap_10_7_7}
\Delta\phi=\phi_{xx}(\alpha x, \beta y, \gamma z)+\phi_{yy}(\alpha x, \beta y, \gamma z)+\phi_{zz}(\alpha x, \beta y, \gamma z)=0.
\end{equation}
The velocity potential function $\phi$ thus satisfies Laplace's equation at each interior point of the three-dimensional space.

\section{The Stream Function}\label{sec:chap_10_8_107}

According to Sec.~(\ref{sec:chap_10_7_106}), the velocity vector
\begin{equation}\label{eq:chap_10_8_1}
V=(p(x,y,z),q(x,y,z),h(x,y,z))
\end{equation}
for a simply connected domain in which the flow is irrotational can be written
\begin{equation}\label{eq:chap_10_8_2}
V=(\phi_{x}(x,y,z),\phi_{y}(x,y,z),\phi_{z}(x,y,z))=grad\texttt{ }\phi(x,y,z),
\end{equation}
where $\phi$ is the velocity potential. When the velocity vector is not the zero vector, it is normal to an equipotential passing through the point $(x,y,z)$. If, moreover, $\psi(x,y,z)$ denotes a harmonic conjugate of $\phi(x,y,z)$ (see Sec.~(\ref{sec:chap_9_4_97})), the velocity vector is tangent to a curve $\psi(x,y,z)=c_{2}$. The surfaces $\psi(x,y,z)=c_{2}$ are called the stream-surfaces of the flow, and the function $\psi$ is the stream function. In particular, a boundary across which fluid cannot flow is a stream-surface.

The analytic function
\[
F(s)=e_{xy}[\phi(x,y,z)+i\psi(x,y,z)]+\omega(z)
\]
is called the complex potential of the flow. Note that (see Theorem~(\ref{th:chap_2_10_l}) in Sec.~(\ref{sec:chap_2_10_20}))
\[
F'(s)=e_{xy}[\phi_{x}(x,y,z)+i\psi_{x}(x,y,z)],
\]
or, in view of the Cauchy-Riemann equations,
\[
F'(s)=e_{xy}[\phi_{x}(x,y,z)-i\phi_{y}(x,y,z)].
\]
Expression~(\ref{eq:chap_10_8_2}) for the velocity thus becomes
\begin{equation}\label{eq:chap_10_8_3}
V=\overline{F'(s)}+\phi_{z}(x,y,z).
\end{equation}
The speed, or magnitude of the velocity, is obtained by writing
\[
|V|=|F'(s)+\phi_{z}(x,y,z)|.
\]

According to equation~(\ref{eq:chap_9_4_5}), Sec.~(\ref{sec:chap_9_4_97}), if $\phi$ is harmonic in a simply connected domain $D$, a harmonic conjugate of $\phi$ there can be written
\[
\psi(x,y,z)=\int_{(x_{0},y_{0},z_{0})}^{(x,y,z)}-\phi_{\eta}(\xi,\eta,\zeta)d\xi+\phi_{\xi}(\xi,\eta,\zeta)d\eta-\phi_{\eta}(\xi,\eta,\zeta)d\zeta.
\]
where the integration is independent of path. With the aid of equations~(\ref{eq:chap_10_7_6}), Sec.~(\ref{sec:chap_10_7_106}), we can, therefore, write
\begin{equation}\label{eq:chap_10_8_4}
\psi(x,y,z)=\int_{C}-\phi_{\eta}(\xi,\eta,\zeta)d\xi+\phi_{\xi}(\xi,\eta,\zeta)d\eta-\phi_{\eta}(\xi,\eta,\zeta)d\zeta.
\end{equation}
where $C$ is any contour in $D$ from $(x_{0},y_{0},z_{0})$ to $(x,y,z)$.

Now it is shown in advanced calculus that the right-hand side of equation~(\ref{eq:chap_10_8_4})
represents the integral with respect to arc length $\sigma$ along $C$ of the normal component
$V_{N}(x,y,z)$ of the vector whose $x$, $y$ and $z$ components are $p(x,y,z)$, $q(x,y,z)$ and $h(x,y,z)$, respectively. So expression~(\ref{eq:chap_10_8_4}) can be written
\begin{equation}\label{eq:chap_10_8_5}
\psi(x,y,z)=\int_{C}V_{N}(x,y,z)d\sigma.
\end{equation}
Physically, then, $\psi(x,y,z)$ represents the time rate of flow of the fluid across $C$. More precisely, $\psi(x,y,z)$ denotes the rate of flow, by volume, across a surface of unit height standing perpendicular to the plane formed by the tangent line and normal line at each point on the curve $C$.% standing perpendicular to the $xy$ plane on the curve $C$.

\begin{example}\label{ex:chap_10_8_1}
When the complex potential is the function
\begin{equation}\label{eq:chap_10_8_1}
F(s)=As=A[e_{xy}(x+iy)+z],
\end{equation}
where $A$ is a positive real constant,
\begin{equation}\label{eq:chap_10_8_1}
\phi(x,y,z)=Ax\texttt{ and }\psi(x,y,z)=Ay.
\end{equation}
The stream-surfaces $\psi(x,y,z)=c_{2}$ are the horizontal planes $y=c_{2}/A$, and the velocity at any point is
\[
V=\overline{F'(s)}=A.
\]

Here a point $(x_{0}, y_{0}, z_{0})$ at which $\psi(x,y,z)=0$ is any point on the $xz$ coordinate plane. If the point $(x_{0}, y_{0}, z_{0})$ is taken as the origin, then $\psi(x,y,z)$ is the rate of flow across any contour drawn from the origin to the point $(x, y, z)$. The flow is uniform and to the right. It can be interpreted as the uniform flow in the upper half space bounded by the $xz$ coordinate plane, which is a stream-surface, or as the uniform flow between two parallel lines $y=y_{1}$ and $y=y_{2}$.
\end{example}

The stream function $\psi$ characterizes a definite flow in a region. The question of whether just one such function exists corresponding to a given region, except possibly for a constant factor or an additive constant, is not examined here. In some of the examples to follow, where the velocity is uniform far from the obstruction, % or in Chap. 11,
where sources and sinks are involved, the physical situation indicates that the flow is uniquely determined by the conditions given in the problem.

A harmonic function is not always uniquely determined, even up to a constant factor, by simply prescribing its values on the boundary of a region. In this example, the function $\psi(x,y,z)=Ay$ is harmonic in the half space $y>0$ and has zero values on the boundary. The function $\psi_{1}(x,y,z)=Be^{x+z}\sin y$ also satisfies those conditions. However, the stream-surface $\psi_{1}(x,y,z)=0$ consists not only of the plane $y=0$ but also of the planes $y=n\pi$ $(n=1, 2, \ldots)$. Here the function $F_{1}(s)=Be^{s}$ is the complex potential for the flow in the strip between the planes $y=0$ and $y=\pi$, both planes making up the stream-surface $\psi(x,y,z)=0$; if $B>0$, the fluid flows to the right along the lower plane and to the left along the upper one.

\section{Flows Around a Corner and Around a Sphere}\label{sec:chap_10_9_108}

In analyzing a flow in the $xyz$, or $s$, space, it is often simpler to consider a corresponding flow in the $uvw$, or $\varpi$, space. Then, if $\phi$ is a velocity potential and $\psi$ a stream function for the flow in the $uvw$ space results in Secs.~(\ref{sec:chap_9_5_98}) and~(\ref{sec:chap_9_6_99}) can be applied to these harmonic functions. That is, when the domain of flow $D_{\varpi}$ in the $uvw$ space is the image of a domain $D_{s}$ under a transformation
\[
\varpi=f(s)=e_{xy}[u(x,y,z)+iv(x,y,z)]+w(z),
\]
where $f$ is analytic, the functions
\[
\phi[u(x,y,z),v(x,y,z),w(z)]\texttt{ and }\psi[u(x,y,z),v(x,y,z),w(z)]
\]
are harmonic in $D_{s}$. These new functions may be interpreted as velocity potential and
stream function in the $xyz$ space. A stream-surface or natural boundary $\psi(u,v,w)=c_{2}$ in the $uvw$ space corresponds to a stream-surface or natural boundary $\psi[u(x,y,z),v(x,y,z),w(z)]=c_{2}$ in the $xyz$ space.

In using this technique, it is often most efficient to first write the complex potential function for the region in the $\varpi$ space and then obtain from that the velocity potential and stream function for the corresponding region in the $xyz$ space. More precisely, if the potential function in the $uvw$ space is
\[
F(\varpi)=e_{xy}[\phi(u,v,w)+i\psi(u,v,w)]+\omega(w),
\]
then the composite function
\[
\begin{split}
F[f(s)] &= e_{xy}\{\phi[u(x,y,z),v(x,y,z),w(z)] \\
        &+ i\psi[u(x,y,z),v(x,y,z),w(z)]\}+\omega[w(z)]
\end{split}
\]
is the desired complex potential in the $xyz$ space.

In order to avoid an excess of notation, we use the same symbols $F$, $\phi$, and $\psi$ for the complex potential, etc., in both the $xyz$ and the $uvw$ space.

\begin{example}\label{ex:chap_10_9_1}
Consider a flow in the first one eighth of the $s$ space for $x>0$, $y>0$, and $z>0$ that comes in downward parallel to the $y$ axis but is forced to turn a corner near the origin. To determine the flow, we recall (Example~(\ref{ex:chap_2_2_1}), Sec.~(\ref{sec:chap_2_2_12})) that the transformation
\[
\varpi=s^{2}=e_{xy}[(x+z)^{2}-y^{2}-z^{2}+i2(x+z)y]+z^{2}
\]
maps the first one eighth of the $s$ space onto the quadrant $-\infty<u<\infty$, $v>0$, and $w>0$ of the $uvw$ space and the boundary of the one eighth of the $s$ space onto the half $uw$ coordinate plane.

From the example in Sec.~(\ref{sec:chap_10_8_107}), we know that the complex potential for a uniform flow to the right in the upper half of the $\varpi$ space is $F=A\varpi$, where $A$ is a positive real constant. The potential in the one eighth of the $s$ space is, therefore,
\begin{equation}\label{eq:chap_10_9_1}
F=As^{2}=e_{xy}A[(x+z)^{2}-y^{2}-z^{2}]+i2(x+z)y]+Az^{2};
\end{equation}
and it follows that the stream function for the flow there is
\begin{equation}\label{eq:chap_10_9_2}
\psi=2A(x+z)y.
\end{equation}
This stream function is, of course, harmonic in the first one eighth of the $s$ space, and it vanishes on the boundary.

The streamlines are branches of the rectangular hyperbolas
\[
2A(x+z)y=c_{v}.
\]
According to equation~(\ref{eq:chap_10_8_3}), Sec.~(\ref{sec:chap_10_8_107}), the velocity of the fluid is
\[
V=\overline{2As}=2A[e_{xy}(x-iy)+z].
\]
Observe that the speed
\[
|V|=2A\sqrt{x^{2}+y^{2}+z^{2}}
\]
of a particle is directly proportional to its distance from the origin. The value of the stream function~(\ref{eq:chap_10_9_2}) at a point $(x, y, z)$ can be interpreted as the rate of flow across a line segment extending from the origin to that point.
\end{example}

\begin{example}\label{ex:chap_10_9_2}
Let a sphere $S=\{s:|s|=1\}$ of unit radius be placed in a large body of fluid flowing with a uniform velocity, the direction of flow is along the $z$ axis of the $xyz$ space. To determine the steady flow around the sphere, let the flow distant from it be parallel to the $z$ axis and to the positive direction. Symmetry shows that points on the $z$ axis exterior to the sphere may be treated as boundary points, and so we need to consider only the
upper part $y\geq0$ of the $yz$ coordinate plane as the region of flow.

Let a spatial circle $C=\{s:|s|=1\}$ on the sphere $S$ be in the $yz$ coordinate plane. Then
the boundary of this region of flow, consisting of the upper semicircle and the parts of the $z$ axis exterior to the circle $C$, is mapped onto the entire $u$ axis by the transformation
\[
\varpi=s+\frac{1}{s}.
\]
The region itself is mapped onto the upper half spatial plane $v\geq0$. The complex potential for the corresponding uniform flow in that half spatial plane is $F=A\varpi$, where $A$ is a positive real constant. Hence the complex potential for the region exterior to the spatial circle and above the $z$ axis is
\begin{equation}\label{eq:chap_10_9_3}
F=A(s+\frac{1}{s}).
\end{equation}

The velocity
\begin{equation}\label{eq:chap_10_9_4}
V=A[1-(\frac{1}{\bar{s}})^{2}]
\end{equation}
approaches $A$ as $|s|$ increases. Thus the flow is nearly uniform and parallel to the $z$ axis at points distant from the sphere, as one would expect. From expression~(\ref{eq:chap_10_9_4}), we see that $V(\bar{s})=\overline{V(s)}$; hence that expression also represents velocities of flow in the lower region, the lower semicircle being a stream-surface.

According to equation~(\ref{eq:chap_10_9_3}), the stream function for the given problem is, in spherical 
coordinates,
\begin{equation}\label{eq:chap_10_9_5}
\psi=A(r-\frac{1}{r})\sin\theta.
\end{equation}
The streamlines
\[
A(r-\frac{1}{r})\sin\theta=c_{2}
\]
are symmetric to the $y$ axis in the transformed spatial plane and have asymptotes parallel to the $z$ axis. Note that when $c_{2}=0$, the stream-surface consists of the spatial circle $r=1$ and the parts of the $z$ axis exterior to the spatial circle.
\end{example}

\section{Flows Around a Sphere in Oceans}\label{sec:chap_10_10_109}

We consider flows around a sphere in oceans. Let assume that a sphere navigates in oceans with a velocity $-\upsilon_{\infty}<0$ relative to the infinite position along the direction of $z$ axis. It is equivalent to that the sphere is fixed and the fluid flows around it with a velocity $\upsilon_{\infty}>0$ at the infinite position along the direction of $z$ axis.

\subsection{Mapping by Elementary Functions}

Let $P_{s}$ denote a complex space and an area $G\subset P_{s}$ in which the sphere navigates. Mapping~(\ref{eq:chap_8_1_1}) in Sec.~(\ref{sec:chap_8_1_83}) gives a transformation $\varpi=As$
where $A$ is a nonzero spatial complex constant and $s\neq0$, in exponential form:
\[
A=r_{a}(e_{xy}e^{i\theta_{a}}\cos\varphi_{a}+\sin\varphi_{a}),\texttt{ }
s=r(e_{xy}e^{i\theta}\cos\varphi+\sin\varphi).
\]
The transformation $\varpi=As$ maps a sphere $|s|=r$ to a sphere $\varpi$.

First, we use the function
\[
\xi=(\frac{As-R}{As+R})^{2}
\]
to map $G$ onto the space $P_{\xi}$ and transform the point $s=\infty$ into $\xi=1$ and the point $s=R/A$ into $\xi=0$.

Next, we use the mapping
\[
\varpi=\frac{C}{1-\xi}\texttt{ where $C$ is a real constant}
\]
to map $P_{\xi}$ onto $P_{\varpi}$ and transform the point $\xi=1$ into $\varpi=\infty$.

Thus the suitable mapping is the function
\[
f(s)=\frac{C}{1-(\frac{As-R}{As+R})^{2}}=\frac{AC}{4R}(s+\frac{R^{2}}{A^{2}s})+\frac{C}{2}.
\]

When $s=\infty$, there is $f(\infty)=\infty$ and
\[
f'(s)=\frac{AC}{4R}(1-\frac{R^{2}}{A^{2}s^{2}}).
\]

As known, the velocity of the sphere in oceans is $\upsilon_{\infty}>0$ relative to the infinite position along the direction of $z$ axis, then there is
\[
f'(\infty)=\frac{AC}{4R}=\upsilon_{\infty}.
\]
So we get
\[
C=\frac{4R}{A}\upsilon_{\infty}
\]
and the potential function
\[
f(s)=\upsilon_{\infty}[(s+\frac{R^{2}}{A^{2}s})+\frac{2R}{A}].
\]

\subsection{Mapping by Laurent Series}

Let the sphere $|s|=R$ and the symbol $\Gamma$ denote the circulation of the flows around a spatial circle $C=\{s:|s|=R\}$ that is on the sphere $|s|=R$ and in a spatial plane determined by an argument $\theta$ and the $z$ axis (Subsec,~(\ref{sec:chap_1_9_4})).

Then to derive the potential function $f(s)$, we first derive its velocity function $f'(s)$. Because $f'(s)$ is analytic in outside of the sphere and the spatial circle $C$, $f'(s)$ can be expanded in a neighborhood of $\infty$ as a Laurent series:
\begin{equation}\label{eq:chap_10_10_1}
f'(s)=\bar{\upsilon}_{\infty}+\frac{C_{-1}}{s}+\frac{C_{-2}}{s^{2}}+\cdots
\end{equation}
where $f'(\infty)=\bar{\upsilon}_{\infty}$ or $\overline{f'(\infty})=\upsilon_{\infty}$.

Let a spatial circle $C_{r}=\{s:|s|=r>R\}$ that is in the same spatial plane with $C$. Then because the circulation of the flows pass through $C_{r}$ is
\[
\Gamma+i\mathbf{N}=\int_{C_{r}}f'(s)ds=2\pi iC_{-1}
\]
where the integral is along the spatial circle in the anti-clockwise direction.

Because no fluid is out of the sphere, there is $\mathbf{N}=0$. That is
\[
\Gamma=2\pi iC_{-1}.
\]
Then Laurent series~(\ref{eq:chap_10_10_1}) becomes
\begin{equation}\label{eq:chap_10_10_2}
f'(s)=\bar{\upsilon}_{\infty}+\frac{\Gamma}{2\pi i}\frac{1}{s}+\frac{C_{-2}}{s^{2}}+\cdots.
\end{equation}

By the integration of the Laurent series~(\ref{eq:chap_10_10_2}), we get
\begin{equation}\label{eq:chap_10_10_3}
f(s)=b_{0}+\bar{\upsilon}_{\infty}s+\frac{\Gamma}{2\pi i}\ln s-\frac{C_{-2}}{s}-\frac{C_{-3}}{2s^{2}}-\cdots.
\end{equation}

To determine coefficient $C_{-2}$, $C_{-3}$,\ldots, we note that because this fluid flows around the spatial circle $C$ in the transformed coordinates system (Subsec,~(\ref{sec:chap_1_9_4})), there should be $Im$ $f(s^{*})=m$ $(a constant)$. Let
\[
s^{*}=Re^{i\theta},\texttt{ }C_{-k}=a_{-k}+ib_{-k},\texttt{ }(k\geq2),\texttt{ }\upsilon_{\infty}=\upsilon_{x}+i\upsilon_{y}.
\]
Then the image part of equation~(\ref{eq:chap_10_10_3}) is
\[
R\upsilon_{x}\sin\theta-R\upsilon_{y}\cos\theta-\frac{\Gamma}{2\pi}\ln s-\frac{b_{-2}\cos\theta-a_{-2}\sin\theta}{R}
\]
\[
-\frac{b_{-3}\cos2\theta-a_{-3}\sin2\theta}{2R^{2}}-\cdots=m\texttt{ (a constant)}
\]
or
\[
(m+\frac{\Gamma}{2\pi}\ln s)+(\frac{b_{-2}}{R}+R\upsilon_{y})\cos\theta-(\frac{a_{-2}}{R}+R\upsilon_{x})\sin\theta
\]
\[
+\frac{b_{-3}}{2R^{2}}\cos2\theta-\frac{a_{-3}}{2R^{2}}\sin2\theta+\cdots=0.
\]
The left side of the equation is trigonometrical series that holds for all values of the argument $\theta$ of the complex variable $s^{*}$. According to the uniqueness theorem of the trigonometrical series expansion, the coefficient of each term of the series expansion should be zero for all terms. That is
\[
m+\frac{\Gamma}{2\pi}\ln s=0,\texttt{ }R\upsilon_{y}+\frac{b_{-2}}{R}=0,
\]
\[
R\upsilon_{x}+\frac{a_{-2}}{R}=0,\texttt{ }b_{-3}=a_{-3}=\cdots=0.
\]
Therefore there are
\[
C_{-2}=-R^{2}\upsilon_{x}-iR^{2}\upsilon_{y}=-\upsilon_{\infty}R^{2},\texttt{ }C_{-3}=C_{-4}=\cdots=0.
\]
After the constant term is omitted, the potential function
\begin{equation}\label{eq:chap_10_10_4}
f(s)=\bar{\upsilon}_{\infty}s+\frac{\upsilon_{\infty}R^{2}}{s}+\frac{\Gamma}{2\pi i}\ln s.
\end{equation}

\section{Three-dimensional Navier-Stokes equations}\label{sec:chap_10_11_110}

Let $\phi(x,y,z)$ denote a real function of three real variables, which is called the velocity potential, and supposed three-dimensional harmonic in a domain $D$, and let $\upsilon$ denote the velocity field $\nabla\phi=\{\phi_{x}(x,y,z),\phi_{y}(x,y,z),\phi_{z}(x,y,z)\}$.

Then the Navier-Stokes equations for the viscous incompressible flow of a fluid on a compact Riemannian manifold $M$ take the form~\cite{H}
\begin{equation}\label{eq:chap_10_11_1}
\frac{\partial\upsilon_{j}}{\partial t}+\upsilon_{k}\frac{\partial\upsilon_{j}}{\partial x_{k}}=-\frac{\partial p}{\partial x_{j}}+\frac{\partial}{\partial x_{k}}[\nu(\frac{\partial\upsilon_{k}}{\partial x_{j}}+\frac{\partial\upsilon_{j}}{\partial x_{k}})]+f_{j},
\end{equation}
\begin{equation}\label{eq:chap_10_11_2}
\nabla\cdot\upsilon=\sum_{k=1}^{^{3}}\frac{\partial\phi_{x_{k}}}{\partial x_{k}}=\frac{\partial^{2}\phi}{(\partial x)^{2}}+\frac{\partial^{2}\phi}{(\partial y)^{2}}+\frac{\partial^{2}\phi}{(\partial z)^{2}}=0\texttt{, and }\upsilon(0)=\upsilon_{0}
\end{equation}
for the velocity field $\upsilon$, where $x_{1}=x$, $x_{2}=y$, $x_{3}=z$, $\nu$ is the dynamic viscous coefficient, $p$ denotes the pressure, and $\textbf{f}$ denotes the external force field.

Now, by a theorem in advanced calculus, the continuity of the partial derivatives of $\phi$ ensures that $\phi_{xy}=\phi_{yx}$, $\phi_{yz}=\phi_{zy}$, and $\phi_{zx}=\phi_{xz}$. It then follows from these relations such that we can rewrite the viscous-shear term
\[
\frac{\partial}{\partial x_{k}}[\nu(\frac{\partial\upsilon_{k}}{\partial x_{j}}+\frac{\partial\upsilon_{j}}{\partial x_{k}})]
=\nu[\frac{\partial}{\partial x_{k}}(\frac{\partial\upsilon_{k}}{\partial x_{j}})+\frac{\partial^{2}\upsilon_{j}}{\partial x_{k}\partial x_{k}}]
+\frac{\partial\nu}{\partial x_{k}}(\frac{\partial\upsilon_{k}}{\partial x_{j}}+\frac{\partial\upsilon_{j}}{\partial x_{k}})
\]
\[
=\nu[\frac{\partial}{\partial x_{j}}(\frac{\partial\upsilon_{k}}{\partial x_{k}})+\frac{\partial^{2}\upsilon_{j}}{\partial x_{k}\partial x_{k}}]+2\frac{\partial\nu}{\partial x_{k}}\frac{\partial\upsilon_{j}}{\partial x_{k}}=2(\nabla\nu\cdot\nabla)\upsilon
\]
where because $\phi$ and its partial derivatives are harmonic
\[
\frac{\partial}{\partial x_{j}}(\frac{\partial\upsilon_{k}}{\partial x_{k}})+\frac{\partial^{2}\upsilon_{j}}{\partial x_{k}\partial x_{k}}=\frac{\partial^{2}\upsilon_{j}}{\partial x_{k}\partial x_{k}}=\Delta\upsilon_{j}=0.
\]

Then equations~(\ref{eq:chap_10_11_1}) become
\begin{equation}\label{eq:chap_10_11_3}
\partial_{t}\upsilon+[(\upsilon-2\nabla\nu)\cdot\nabla]\upsilon=-\nabla p+\textbf{f}.
\end{equation}

Moreover, equations~(\ref{eq:chap_10_11_3}) can be written
\begin{equation}\label{eq:chap_10_11_4}
\begin{split}
&\frac{\partial\phi_{x}}{\partial t}+(\phi_{x}-2\nu_{x})\phi_{xx}+(\phi_{y}-2\nu_{y})\phi_{xy}+(\phi_{z}-2\nu_{z})\phi_{xz}=-\frac{\partial p}{\partial x}+f_{x} \\
&\frac{\partial\phi_{y}}{\partial t}+(\phi_{x}-2\nu_{x})\phi_{yx}+(\phi_{y}-2\nu_{y})\phi_{yy}+(\phi_{z}-2\nu_{z})\phi_{yz}=-\frac{\partial p}{\partial y}+f_{y} \\
&\frac{\partial\phi_{z}}{\partial t}+(\phi_{x}-2\nu_{x})\phi_{zx}+(\phi_{y}-2\nu_{y})\phi_{zy}+(\phi_{z}-2\nu_{z})\phi_{zz}=-\frac{\partial p}{\partial z}+f_{z} \\
\end{split}.
\end{equation}

Thus, when the dynamic viscous coefficient $\nu$ is a constant, equations~(\ref{eq:chap_10_11_1}) become Euler's equations for three-dimensional ideal incompressible fluid flow.

\subsection{Mapping by the Function $f(\alpha x+i\beta y+\gamma z)$}\label{sec:chap_10_11_110_1}

From Theorem~(\ref{th:chap_2_15_3}) in Sec.~(\ref{sec:chap_2_15_25}), there is
\[
f(\alpha x+i\beta y+\gamma z)=\phi(\alpha x, \beta y, \gamma z)+i\psi(\alpha x, \beta y, \gamma z).
\]
In three-dimensional regions where the fluid is free from sources or sinks, $\phi(\alpha x, \beta y, \gamma z)$ is a harmonic function of three real variables $x$, $y$, and $z$:
\[
\Delta\phi=\phi_{xx}(\alpha x, \beta y, \gamma z)+\phi_{yy}(\alpha x, \beta y, \gamma z)+\phi_{zz}(\alpha x, \beta y, \gamma z)=0.
\]
The velocity potential function $\phi$ thus satisfies Laplace's equation, the first equation in the equation group~(\ref{eq:chap_10_11_2}) at each interior point of the three-dimensional space.

From expressions~(\ref{eq:chap_2_15_2c}) in Sec.~(\ref{sec:chap_2_15_25}), equations~(\ref{eq:chap_10_11_4}) now become
\begin{equation}\label{eq:chap_10_11_5a}
    \begin{array}{c}
        \alpha\frac{\partial\phi_{\bar{x}}}{\partial t}
        +\alpha^{2}(\alpha\phi_{\bar{x}}-2\nu_{x})\phi_{\overline{xx}}
        +\alpha\beta(\beta\phi_{\bar{y}}-2\nu_{y})\phi_{\overline{xy}}
        +\alpha\gamma(\gamma\phi_{\bar{z}}-2\nu_{z})\phi_{\overline{xz}} \\
        =-\frac{\partial p}{\partial x}+f_{x} \\
        \beta\frac{\partial\phi_{\bar{y}}}{\partial t}
        +\alpha\beta(\alpha\phi_{\bar{x}}-2\nu_{x})\phi_{\overline{yx}}
        +\beta^{2}(\beta\phi_{\bar{y}}-2\nu_{y})\phi_{\overline{yy}}
        +\beta\gamma(\gamma\phi_{\bar{z}}-2\nu_{z})\phi_{\overline{yz}} \\
        =-\frac{\partial p}{\partial y}+f_{y} \\
        \gamma\frac{\partial\phi_{\bar{z}}}{\partial t}
        +\alpha\gamma(\alpha\phi_{\bar{x}}-2\nu_{x})\phi_{\overline{zx}}
        +\beta\gamma(\beta\phi_{\bar{y}}-2\nu_{y})\phi_{\overline{zy}}
        +\gamma^{2}(\gamma\phi_{\bar{z}}-2\nu_{z})\phi_{\overline{zz}} \\
        =-\frac{\partial p}{\partial z}+f_{z} \\
    \end{array}.
\end{equation}

We choose a suitable function $f(\alpha x+i\beta y+\gamma z)$ to satisfy equations~(\ref{eq:chap_10_11_4}) or~(\ref{eq:chap_10_11_5a}), and the initial condition expressed by the last equation in the equation group~(\ref{eq:chap_10_11_2}).

\subsection{Mapping by Polynomial Functions}\label{sec:chap_10_11_110_2}

From lemma~(\ref{le:chap_2_15_2}), let the real functions $\phi(x,y,z)$ of three real variables $x$, $y$, and $z$ be
\begin{equation}\label{eq:chap_10_11_5}
\begin{split}
\phi(x,y,z) &= a_{1,xz}x-a_{1,y}y
+a_{2,xz}(x^{2}+2xz-y^{2})-2a_{2,y}(x+z)y \\
            &+ a_{3,xz}[(x^{3}+3x^{2}z)-3(x+z)y^{2}]-a_{3,y}[3(x^{2}+2xz)y-y^{3}],
\end{split}
\end{equation}
where $a_{n}=e_{xy}(a_{n,x}+ia_{n,y})+a_{n,z}$ for $n=1,2,3$ are spatial complex constants and $a_{n,xz}=a_{n,x}+a_{n,z}$.

Then the first-order partial derivatives of $\phi$ are
\begin{equation}\label{eq:chap_10_11_6}
    \begin{split}
        & \phi_{x}=a_{1,xz}+2a_{2,xz}(x+z)-2a_{2,y}y+3a_{3,xz}(x^{2}+2xz-y^{2}) \\
        & \texttt{ }\texttt{ }-6a_{3,y}(x+z)y, \\
        & \phi_{y}=-a_{1,y}-2a_{2,xz}y-2a_{2,y}(x+z)-6a_{3,xz}(x+z)y \\
        & \texttt{ }\texttt{ }-3a_{3,y}(x^{2}+2xz-y^{2}), \\
        & \phi_{z}=2a_{2,xz}x-2a_{2,y}y+3a_{3,xz}(x^{2}-y^{2})-6a_{3,y}xy,
    \end{split}
\end{equation}
the second-order partial derivatives of $\phi$ are
\begin{equation}\label{eq:chap_10_11_7}
    \begin{split}
        & \phi_{xx}=2a_{2,xz}+6a_{3,xz}(x+z)-6a_{3,y}y, \\
        & \phi_{xy}=-2a_{2,y}-6a_{3,xz}y-6a_{3,y}(x+z), \\
        & \phi_{xz}=2a_{2,xz}+6a_{3,xz}x-6a_{3,y}y, \\
        & \phi_{yx}=-2a_{2,y}-6a_{3,xz}y-6a_{3,y}(x+z), \\
        & \phi_{yy}=-2a_{2,xz}-6a_{3,xz}(x+z)+6a_{3,y}y, \\
        & \phi_{yz}=-2a_{2,y}-6a_{3,xz}y-6a_{3,y}x, \\
        & \phi_{zx}=2a_{2,xz}+6a_{3,xz}x-6a_{3,y}y, \\
        & \phi_{zy}=-2a_{2,y}-6a_{3,xz}y-6a_{3,y}x, \\
        & \phi_{zz}=0,
    \end{split}
\end{equation}
and some of the third-order partial derivatives of $\phi$ are
\begin{equation}\label{eq:chap_10_11_8}
    \begin{array}{ccc}
        \phi_{xxx}=6a_{3,xz}, & \phi_{xyy}=-6a_{3,xz}, & \phi_{xzz}=0, \\
        \phi_{yxx}=-6a_{3,y}, & \phi_{yyy}=\texttt{ }6a_{3,y},\texttt{ } & \phi_{yzz}=0, \\
        \phi_{zxx}=6a_{3,xz}, & \phi_{zyy}=-6a_{3,xz}, & \phi_{zzz}=0.
    \end{array}
\end{equation}

Lemma~(\ref{le:chap_2_15_2}) shows that the first-order partial derivatives of $\phi$ satisfy the Cauchy-Riemann equations~(\ref{eq:chap_2_15_5}) and the second-order partial derivatives of $\phi$ satisfy the Laplace's equations~(\ref{eq:chap_2_15_6}) with boundary conditions in the $xyz$ coordinate space.

The velocity potential function $\phi$ thus satisfies Laplace's equation that is the first equation in the equation group~(\ref{eq:chap_10_11_2}) at each interior point of the three-dimensional space, and equations~(\ref{eq:chap_10_11_4}) become
\begin{equation}\label{eq:chap_10_11_9}
    \begin{split}
        & \frac{\partial\phi_{x}}{\partial t}+(\phi_{x}-2\nu_{x})\phi_{xx}+(\phi_{y}-2\nu_{y})\phi_{xy}+(\phi_{z}-2\nu_{z})\phi_{xz}=-\frac{\partial p}{\partial x}+f_{x} \\
        & \frac{\partial\phi_{y}}{\partial t}+(\phi_{x}-2\nu_{x})\phi_{yx}+(\phi_{y}-2\nu_{y})\phi_{yy}+(\phi_{z}-2\nu_{z})\phi_{yz}=-\frac{\partial p}{\partial y}+f_{y} \\
        & \frac{\partial\phi_{z}}{\partial t}+(\phi_{x}-2\nu_{x})\phi_{zx}+(\phi_{y}-2\nu_{y})\phi_{zy}=-\frac{\partial p}{\partial z}+f_{z} \\
    \end{split}
\end{equation}
where the first-order partial derivatives of $\phi$ are quadratic functions of two variables $x$ and $y$, or constants for $z$, and the second-order partial derivatives of $\phi$ are linear functions of two variables $x$ and $y$.

Then functions on both sides should be at most cubic functions of two variables $x$ and $y$, and linear functions of $z$.

There are 9 real constants $a_{n,x}$, $a_{n,y}$, and $a_{n,z}$ for $n=1,2,3$ to be determined by using the conditions on the right side of equations~(\ref{eq:chap_10_11_9}).

Equations~(\ref{eq:chap_10_11_8}) show $\Delta\phi_{x}=0$, $\Delta\phi_{y}=0$, and $\Delta\phi_{z}=0$. So when the incompressible flow of a fluid is steady and the viscous coefficient $\nu$ is a constant, equations~(\ref{eq:chap_10_11_9}) become
\begin{equation}\label{eq:chap_10_11_10}
    \begin{split}
        & \phi_{x}\phi_{xx}+\phi_{y}\phi_{xy}+\phi_{z}\phi_{xz}=-\frac{\partial p}{\partial x}+f_{x}, \\
        & \phi_{x}\phi_{yx}+\phi_{y}\phi_{yy}+\phi_{z}\phi_{yz}=-\frac{\partial p}{\partial y}+f_{y}, \\
        & \phi_{x}\phi_{zx}+\phi_{y}\phi_{zy}=-\frac{\partial p}{\partial z}+f_{z}. \\
    \end{split}
\end{equation}

%-----------------------------------------------------------------------
% Beginning of chap11.tex
%-----------------------------------------------------------------------
%
% AMS-LaTeX 1.2 sample file for a monograph, based on amsbook.cls.
% This is a data file input by chapter.tex.
%%%%%%%%%%%%%%%%%%%%%%%%%%%%%%%%%%%%%%%%%%%%%%%%%%%%%%%%%%%%%%%%%%%%

%\part{This is a Part Title Sample}

\chapter{The Schwarz-Christoffel Transformation}\label{ch:chap_11}

In this chapter, we construct a transformation, known as the Schwarz-Christoffel transformation, which maps a polygon line that is corresponding to whole $x$ axis in the $xz$ coordinate plane, and the upper half of the $s$ space onto a given simple closed spatial polygon and its interior in the $\varpi$ space. Applications are made to the solution of problems in fluid flow and electrostatic potential theory.

\section{Mapping a Polygon Line in the $xz$ Plane onto a Spatial Polygon}\label{sec:chap_11_1_109}

We represent the unit vector which is tangent to a smooth spatial arc $C$ at a point $s_{0}$ by the spatial complex number $t=e_{xy}e^{\theta_{t}}\cos\varphi_{t}+\sin\varphi_{t}$, and we let the spatial complex number $\tau=e_{xy}e^{\phi_{\tau}}\cos\psi_{\tau}+\sin\psi_{\tau}$ denote the unit vector tangent to the image $\Gamma$ of $C$ at the corresponding point $\varpi_{0}$ under a transformation $\varpi=f(s)$. We assume
that $f$ is analytic at $s_{0}$ and that $f'(s_{0})\neq0$. According to Sec.~(\ref{sec:chap_9_1_94}),
\begin{equation}\label{eq:chap_11_1_1}
\arg\tau=\arg[f'(s_{0})t]
\end{equation}
where
\[
\arg\tau=(\phi_{\tau},\psi_{\tau}),\texttt{ }
\arg t=(\theta_{t},\varphi_{t}),\texttt{ }
\arg f'(s_{0})=(\theta_{f},\varphi_{f}),
\]
$\psi_{\tau}=\arcsin(\sin\varphi_{f}\sin\varphi_{t})$, and
\[
\phi_{\tau}=\arctan[\frac{\sin(\theta_{f}+\theta_{t})
+\sin\theta_{f}\tan\varphi_{t}+\sin\theta_{t}\tan\varphi_{f}}
{\cos(\theta_{f}+\theta_{t})
+\cos\theta_{f}\tan\varphi_{t}+\cos\theta_{t}\tan\varphi_{f}}].
\]
In particular, if $C$ is a segment in the $xz$ coordinate plane with positive sense to the right, then $t=(e_{xy}x+z)/\sqrt{x^{2}+z^{2}}$ and $\arg_{c}t=\theta_{t}=0$ at each point $s_{0}=e_{xy}x+z$ on $C$. In that case,
\[
\phi_{\tau}=\arctan[\frac{(1+\tan\varphi_{t})\sin\theta_{f}}
{(1+\tan\varphi_{t})\cos\theta_{f}+\tan\varphi_{f}}].
\]
So equation~(\ref{eq:chap_11_1_1}) becomes
\begin{equation}\label{eq:chap_11_1_2}
\arg\tau=(\arctan[\frac{(1+\tan\varphi_{t})\sin\theta_{f}}
{(1+\tan\varphi_{t})\cos\theta_{f}+\tan\varphi_{f}}],\arcsin(\sin\varphi_{f}\sin\varphi_{t})).
\end{equation}
If $\varphi_{t}$ is a constant argument and $f'(s_{0})$ has two constant arguments along that segment, it follows that $\arg\tau$ is constant. Hence the image $\Gamma$ of $C$ is also a segment of a straight line.

Let us now construct a transformation $\varpi=f(s)$ which maps a polygon line $P_{xz}$ that is corresponding to whole $x$ axis in the $xz$ coordinate plane of the $s$ space onto a spatial polygon of $n$ sides in the $\varpi$ space, where $\chi=e_{xy}x+z$, $\chi_{j}=e_{xy}x_{j}+z_{j}$ $(j=1,2,\ldots,n-1)$ and $\infty$ are the points on that $xz$ coordinate plane of the $s$ space, whose images are to be the vertices of the spatial polygon and where
\[
x_{1}<x_{2}<\cdots<x_{n-1}.
\]
The vertices are the points $\varpi_{j}=f(\chi_{j})$ $(j=1, 2, \ldots, n-1)$ and $\varpi_{n}=f(\infty)$. The
function $f$ should be such that $\arg f'(s)$ jumps from one constant value to another at
the points $s=\chi_{j}$ as the point $s$ traces out the polygon line $P_{xz}$ in the $xz$ coordinate plane.

If the function $f$ is chosen such that
\begin{equation}\label{eq:chap_11_1_3}
f'(s)=A(s-\chi_{1})^{-k_{1}}(s-\chi_{2})^{-k_{2}}\cdots(s-\chi_{n-1})^{-k_{n-1}},
\end{equation}
where $A$ is a spatial complex constant and each $k_{j}$ is a real constant, then the arguments of $f'(s)$ change in the prescribed manner as $s$ describes the real $xz$ coordinate plane; from expression~(\ref{eq:chap_1_7_9}) in Sec.~(\ref{sec:chap_1_7}), let
\[
\begin{split}
& s=r(e_{xy}e^{i\theta}\cos\varphi+\sin\varphi), \\
& A=R_{0}(e_{xy}e^{i\phi_{0}}\cos\psi_{0}+\sin\psi_{0}), \\
& s_{j}=s-\chi_{j}=r_{j}(e_{xy}e^{i\theta_{j}}\cos\varphi_{j}+\sin\varphi_{j}), \\
& \varpi_{j}=s_{j}^{-k_{j}}=R_{j}(e_{xy}e^{i\phi_{j}}\cos\psi_{j}+\sin\psi_{j})\texttt{ }(j=1,2,\ldots,n-1).
\end{split}
\]
Then for $j=1,2,\ldots,n-1$ there are
\[
r_{j}=|s-\chi_{j}|,\texttt{ }\theta_{j}=\arctan[\frac{Im(s-\chi_{j})}{Re(s-\chi_{j})}],\texttt{ }\varphi_{j}=\arcsin[\frac{Re_{s}(s-\chi_{j})}{r_{j}}],
\]
\[
R_{j}=r_{j}^{-k_{j}},\texttt{ }\psi_{j}=\arcsin(\sin^{-k_{j}}\varphi_{j}),\texttt{ and }
\]
\[
\phi_{j}=\frac{\sum_{k=0}^{j-1}C_{j}^{k}\sin[(j-k)\theta_{j}]\tan^{k}\varphi_{j}}
{\sum_{k=0}^{j-1}C_{j}^{k}\{(\cot\varphi_{j}+2\cos\theta_{j})^{j-k}-\cos[(j-k)\theta_{j}]\tan^{k}\varphi_{j}\}}\texttt{ }(C_{j}^{k}=\frac{j}{(j-k)!k!}).
\]
From expression~(\ref{eq:chap_1_7_5}) in Sec.~(\ref{sec:chap_1_7}), the arguments of the derivative~(\ref{eq:chap_11_1_3}) can be written
\begin{equation}\label{eq:chap_11_1_4}
\arg f'(s)=(\theta_{f},\varphi_{f})\texttt{ where }\varphi_{f}=\arcsin(\prod_{k=0}^{n-1}\sin\psi_{k})\texttt{ and }
\end{equation}
\[
\begin{split}
\tan\theta_{f} &= [\sin(\sum_{k=0}^{n-1}\phi_{k})
+\sum_{j=0}^{n-1}\tan\psi_{j}\sin(\sum_{k=0,k\neq j}^{n-1}\phi_{k})+\cdots+\sum_{j=0}^{n-1}\sin\phi_{j}\prod_{k=0,k\neq j}^{n-1}\tan\psi_{k}] \\
&\div [\cos(\sum_{k=0}^{n-1}\phi_{k})
+\sum_{j=0}^{n-1}\tan\psi_{j}\cos(\sum_{k=0,k\neq j}^{n-1}\phi_{k})+\cdots+\sum_{j=0}^{n-1}\cos\phi_{j}\prod_{k=0,k\neq j}^{n-1}\tan\psi_{k}].
\end{split}
\]
When $s=\chi=e_{xy}x+z$ and $x<x_{1}$
\[
\arg_{c}(s-\chi_{1})=\arg_{c}(s-\chi_{2})=\cdots=\arg_{c}(s-\chi_{n-1})=\pi.
\]
When $x_{1}<x<x_{2}$, the argument $\arg_{c}(s-\chi_{1})$ is $0$ and each of the other arguments is $\pi$. According to equation~(\ref{eq:chap_11_1_4}), then, $\arg_{c}f'(s)$ increases abruptly by the angle $k_{1}\pi$ as
$s$ moves to the right through the point $\chi_{1}$. It again jumps in value, by the amount $k_{2}\pi$, as $s$ passes through the point $\chi_{2}$, etc.

In view of equation~(\ref{eq:chap_11_1_2}), the unit vector $\tau$ is constant in direction as $s$ moves from $\chi_{j-1}$ to $\chi_{j}$; the point $\varpi$ thus moves in that fixed direction along a straight line. The direction of $\tau$ changes abruptly, by the angle $k_{j}\pi$, at the image point $\varpi_{1}$ of $s=\chi_{1}$. Those angles $k_{j}\pi$ are the exterior angles of the polygon described by the point $\varpi$.

The exterior angles can be limited to angles between $-\pi$ and $\pi$, in which case $-1<k_{j}<1$. We assume that the sides of the polygon never cross one another and that the polygon is given a positive, or counterclockwise, orientation. The sum of the exterior angles of a closed spatial polygon is, then, $2\pi$; and the exterior angle at the vertex $\varpi_{n}$, which is the image of the point $s=\infty$, can be written
\[
k_{n}\pi=2\pi-(k_{1}+k_{2}+\cdots+k_{n-1})\pi.
\]
Thus the numbers $k_{j}$ must necessarily satisfy the conditions
\begin{equation}\label{eq:chap_11_1_5}
k_{1}+k_{2}+\cdots+k_{n-1}+k_{n}=2\texttt{ }-1<k_{j}<1\texttt{ }(j=1,2,\ldots,n).
\end{equation}

Note that $k_{n}=0$ if
\begin{equation}\label{eq:chap_11_1_6}
k_{1}+k_{2}+\cdots+k_{n-1}=2.
\end{equation}
This means that the direction of $\tau$ does not change at the point $\varpi_{n}$. So $\varpi_{n}$ is not a vertex, and the polygon has $n-1$ sides.

The existence of a mapping function $f$ whose derivative is given by equation~(\ref{eq:chap_11_1_3}) will be established in the next section.

\section{Schwarz-Christoffel Transformation}\label{sec:chap_11_2_110}

In our expression (Sec.~(\ref{sec:chap_11_1_109}))
\begin{equation}\label{eq:chap_11_2_1}
f'(s)=A(s-\chi_{1})^{-k_{1}}(s-\chi_{2})^{-k_{2}}\cdots(s-\chi_{n-1})^{-k_{n-1}}
\end{equation}
for the derivative of a function that is to map a polygon line $P_{xz}$ that is corresponding to whole $x$ axis in the $xz$ coordinate plane of the $s$ space onto a spatial polygon in the $\varpi$ space, let the factors $(s-\chi_{j})^{-k_{j}}$ represent branches of power functions with branch cuts extending below that polygon line. To be specific, write
\begin{equation}\label{eq:chap_11_2_2}
(s-\chi_{j})^{-k_{j}}=|s-\chi_{j}|^{-k_{j}}(e_{xy}e^{i\theta_{j}}\cos\varphi_{j}+\sin\varphi_{j})^{-k_{j}}\texttt{ }(-\frac{\pi}{2}<\theta_{j}<\frac{3\pi}{2}),
\end{equation}
where $\theta_{j}=\arg_{c}(s-\chi_{j})$ and $\varphi_{j}=\arg_{r}(s-\chi_{j})$ for $j=1,2,\ldots,n-1$. Then $f'(s)$ is analytic everywhere in the half space $y\geq0$ except at the $n-1$ branch points $\chi_{j}$.

If $s_{0}$ is a point in that region of analyticity, denoted here by $R$, then the function
\begin{equation}\label{eq:chap_11_2_3}
F(s)=\int_{s_{0}}^{s}f'(\varsigma)d\varsigma
\end{equation}
is single-valued and analytic throughout the same region, where the path of integration
from $s_{0}$ to $s$ is any spatial contour lying within $R$.  Moreover, $F'(s)=f'(s)$ (see Sec.~(\ref{sec:chap_4_7_42})).

To define the function $F$ at the point $s=\chi_{1}$ so that it is continuous there, we note that $(s-\chi_{1})^{-k_{1}}$ is the only factor in expression~(\ref{eq:chap_11_2_1}) that is not analytic at $\chi_{1}$. Hence, if $\phi(s)$ denotes the product of the rest of the factors in that expression, $\phi(s)$ is analytic at
$\chi_{1}$ and is represented throughout an open sphere $|s-\chi_{1}|<R_{1}$ by its Taylor series about $\chi_{1}$.  So we can write
\[
\begin{split}
f'(s) &= (s-\chi_{1})^{-k_{1}}\phi(s) \\
&=(s-\chi_{1})^{-k_{1}}[\phi(\chi_{1})+\frac{\phi(\chi_{1})}{1!}(s-\chi_{1})+\frac{\phi''(\chi_{1})}{2!}(s-\chi_{1})^{2}+\cdots],
\end{split}
\]
or
\begin{equation}\label{eq:chap_11_2_4}
f'(s)=\phi(\chi_{1})(s-\chi_{1})^{-k_{1}}+\frac{\phi(\chi_{1})}{1!}(s-\chi_{1})^{1-k_{1}}\psi(s)
\end{equation}
where $\psi(s)$ is analytic and, therefore, continuous throughout the entire open sphere. Since $1-k_{1}>0$, the last term on the right in equation~(\ref{eq:chap_11_2_4}) thus represents a continuous function of $s$ throughout the upper half of the sphere, where $Im$ $s\geq0$, if we assign it
the value zero at $s=\chi_{1}$. It follows that the integral
\[
\int_{S_{1}}^{s}(\varsigma-\chi_{1})^{1-k_{1}}\psi(\varsigma)d\varsigma
\]
of that last term along a spatial contour from $S_{1}$ to $s$, where $S_{1}$ and the spatial contour lie in the half sphere, is a continuous function of $s$ at $s=\chi_{1}$. The integral
\[
\int_{S_{1}}^{s}(\varsigma-\chi_{1})^{-k_{1}}\psi(\varsigma)d\varsigma
=\frac{1}{1-k_{1}}[(s-\chi_{1})^{1-k_{1}}-(S_{1}-\chi_{1})^{1-k_{1}}]
\]
along the same path also represents a continuous function of $s$ at $\chi_{1}$ if we define the value of the integral there as its limit as $s$ approaches $\chi_{1}$ in the half sphere. The integral of the function~(\ref{eq:chap_11_2_4}) along the stated path from $S_{1}$ to $s$ is, then, continuous at $s=\chi_{1}$; and the same is true of integral~(\ref{eq:chap_11_2_3}) since it can be written as an integral along a spatial contour in $R$ from $s_{0}$ to $S_{1}$ plus the integral from $S_{1}$ to $s$.

The above argument applies at each of the $n-1$ points $\chi_{j}$ to make $F$ continuous throughout the region $y\geq0$.

From equation~(\ref{eq:chap_11_2_1}), we can show that, for a sufficiently large positive number $R$, a positive constant $M$ exists such that if $lm$ $s\geq0$, then
\begin{equation}\label{eq:chap_11_2_5}
|f'(s)|<\frac{M}{|s|^{2-k_{n}}}\texttt{ whenever }|s|>R.
\end{equation}
Since $2-k_{n}>1$, this order property of the integrand in equation~(\ref{eq:chap_11_2_3}) ensures the existence of the limit of the integral there as $s$ tends to infinity; that is, a number $W_{n}$ exists such that
\begin{equation}\label{eq:chap_11_2_6}
\lim_{s \to W_{n}}F(s)=W_{n}\texttt{ }(\texttt{Im }s\geq0).
\end{equation}

Our mapping function, whose derivative is given by equation~(\ref{eq:chap_11_2_1}), can be written
$f(s)=F(s)+B$, where $B$ is a spatial complex constant. The resulting transformation,
\begin{equation}\label{eq:chap_11_2_7}
\varpi=A\int_{s_{0}}^{s}(\varsigma-\chi_{1})^{-k_{1}}(\varsigma-\chi_{2})^{-k_{2}}\cdots(\varsigma-\chi_{n-1})^{-k_{n-1}}d\varsigma+B
\end{equation}
is the Schwarz-Christoffel transformation.

Transformation~(\ref{eq:chap_11_2_7}) is continuous throughout the half space $y\geq0$ and is conformal there except for the points $\chi_{j}$. We have assumed that the numbers $k_{j}$ satisfy conditions~(\ref{eq:chap_11_1_5}), Sec.~(\ref{sec:chap_11_1_109}). In addition, we suppose that the constants $\chi_{j}$ and $k_{j}$ are such that the sides of the spatial polygon do not cross, so that the spatial polygon is a simple closed spatial contour. Then, according to Sec.~(\ref{sec:chap_11_1_109}), as the point $s$ describes the polygon line $P_{xz}$ in the positive direction, its image $\varpi$ describes the polygon $P$ in the positive sense; and there is a one to one correspondence between points on that polygon line $P_{xz}$ and points on $P$. According to condition~(\ref{eq:chap_11_2_6}), the image $\varpi_{n}$ of the point $s=\infty$ exists and $\varpi_{n}=W_{n}+B$.

If $s$ is an interior point of the upper half space $y\geq0$ and $\chi_{0}$ is any point on the
polygon line $P_{xz}$ other than one of the $\chi_{j}$, then the angle from the vector $\varsigma$ at $\chi_{0}$ up to the line segment joining $\chi_{0}$ and $s$ is positive and less than $\pi$. At the image $\varpi_{0}$ of $\chi_{0}$, the corresponding angle from the vector $\tau$ to the image of the line segment joining $\chi_{0}$ and $s$ has that same value. Thus the images of interior points in the half space lie to the left of the sides of the polygon, taken counterclockwise. %A proof that the transformation establishes a one to one correspondence between the interior points of the half space and the points within the polygon is left to the reader (Exercise 3).

Given a specific polygon $P$, let us examine the number of constants in the Schwarz-Christoffel transformation that must be determined in order to map the polygon line $P_{xz}$ in the $xz$ coordinate plane onto $P$ in $\varpi$ space. For this purpose, we may write $s_{0}=0$, $A=1$, and $B=0$ and simply require that the polygon line $P_{xz}$ be mapped onto some spatial polygon $P'$ similar to $P$. The size and position of $P'$ can then be adjusted to match those of $P$ by introducing the appropriate constants $A$ and $B$.

The numbers $k_{j}$ are all determined from the exterior angles at the vertices of $P$. The $n-1$ constants $\chi_{j}$ remain to be chosen. The image of the polygon line $P_{xz}$ in the $xz$ coordinate plane is some spatial polygon $P'$ that has the same angles as $P$. But if $P'$ is to be similar to $P$, then $n-2$ connected sides must have a common ratio to the corresponding sides of $P$; this condition is expressed by means of $n-3$ equations in the $n-1$ unknowns $\chi_{j}$. Thus two of the numbers $\chi_{j}$ or two relations between them, can be chosen arbitrarily, provided those $n-3$ equations in the remaining $n-3$ unknowns have real-valued solutions.

When a finite point $s=\chi_{n}$ on the $xz$ coordinate plane, instead of the point at infinity, represents
the point whose image is the vertex $\varpi_{n}$, it follows from Sec.~(\ref{sec:chap_11_1_109}) that the Schwarz-Christoffel transformation takes the form
\begin{equation}\label{eq:chap_11_2_8}
\varpi=A\int_{s_{0}}^{s}(\varsigma-\chi_{1})^{-k_{1}}(\varsigma-\chi_{2})^{-k_{2}}\cdots(\varsigma-\chi_{n})^{-k_{n}}d\varsigma+B
\end{equation}
where $k_{1}+k_{2}+\cdots+k_{n}=2$. The exponents $k_{j}$ are determined from the exterior angles of the polygon. But, in this case, there are $n$ constants $\chi_{j}$ that must satisfy the $n-3$ equations noted above. Thus three of the numbers $\chi_{j}$, or three conditions on those n numbers, can be chosen arbitrarily in transformation~(\ref{eq:chap_11_2_8}) of the polygon line $P_{xz}$ in the $xz$ coordinate plane onto a given spatial polygon.

\section{Triangles and Rectangles}\label{sec:chap_11_3_111}

The Schwarz-Christoffel transformation is written in terms of the points $\chi_{j}$ and not in terms of their images, the vertices of the spatial polygon. No more than three of those points can be chosen arbitrarily; so, when the given polygon has more than three sides, some of the points $\chi_{j}$ must be determined in order to make the given spatial polygon, or any polygon similar to it, be the image of the polygon line $P_{xz}$ in the $xz$ coordinate plane. The selection of conditions for the determination of those constants, conditions that are convenient to use, often requires ingenuity.

Another limitation in using the transformation is due to the integration that is involved. Often the integral cannot be evaluated in terms of a finite number of elementary functions. In such cases, the solution of problems by means of the transformation can become quite involved.

If the spatial polygon is a triangle with vertices at the points $\varpi_{1}$, $\varpi_{2}$, and $\varpi_{3}$, the transformation can be written
\begin{equation}\label{eq:chap_11_3_1}
\varpi=A\int_{s_{0}}^{s}(\varsigma-\chi_{1})^{-k_{1}}(\varsigma-\chi_{2})^{-k_{2}}(\varsigma-\chi_{3})^{-k_{3}}d\varsigma+B
\end{equation}
where $k_{1}+k_{2}+k_{3}=2$. In terms of the interior angles $\theta_{j}$,
\[
k_{j}=1-\frac{1}{\pi}\theta_{j}\texttt{ }(j=1,2,3).
\]
Here we have taken all three points $\chi_{j}$ as finite points on the $xz$ coordinate plane to form a polygon line $P_{xz}$. Arbitrary values can be assigned to each of them. The complex constants $A$ and $B$, which are associated with the size and position of the triangle, can be determined so that the upper half space is mapped onto the given triangular spatial region.

If we take the vertex $\varpi_{3}$ as the image of the point at infinity, the transformation becomes
\begin{equation}\label{eq:chap_11_3_2}
\varpi=A\int_{s_{0}}^{s}(\varsigma-\chi_{1})^{-k_{1}}(\varsigma-\chi_{2})^{-k_{2}}d\varsigma+B
\end{equation}
where two pairs of arbitrary real values can be assigned to the components of $\chi_{1}$ and $\chi_{2}$.

The integrals in equations~(\ref{eq:chap_11_3_1}) and~(\ref{eq:chap_11_3_2}) do not represent elementary functions unless the triangle is degenerate with one or two of its vertices at infinity. The integral in equation~(\ref{eq:chap_11_3_2}) becomes an elliptic integral when the triangle is equilateral or when it is a right triangle with one of its angles equal to either $\pi/3$ or $\pi/4$.

When the spatial polygon is a rectangle, each $k_{j}=1/2$. If we choose $\pm1$ and $\pm a$ as the points $\chi_{j}$ whose images are the vertices and write
\begin{equation}\label{eq:chap_11_3_8}
g(s)=(s+a)^{-1/2}(s+1)^{-1/2}(s-1)^{-1/2}(s-a)^{-1/2},
\end{equation}
where $0\leq\arg_{c}(s-\chi_{j})\leq\pi$, the Schwarz-Christoffel transformation becomes
\begin{equation}\label{eq:chap_11_3_9}
\varpi=-\int_{0}^{s}g(\varsigma)d\varsigma,
\end{equation}
except for a transformation $W=A\varpi+B$ to adjust the size and position of the rectangle. Integral~(\ref{eq:chap_11_3_9}) is a constant times the elliptic integral
\[
\int_{0}^{s}(1-\varsigma^{2})^{-1/2}(1-k^{2}\varsigma^{2})^{-1/2}d\varsigma\texttt{ }(k=\frac{1}{a}),
\]
but the form~(\ref{eq:chap_11_3_8}) of the integrand indicates more clearly the appropriate branches of the power functions involved.

%\section{Degenerate Polygons}\label{sec:chap_11_4_112}

%We now apply the Schwarz-Christoffel transformation to some degenerate polygons for which the integrals represent elementary functions. For purposes of illustration, the examples here result in transformations that we have already seen in Chap. 8.

\section{Fluid Flow in a Channel through a Slit}\label{sec:chap_11_5_113}

We now present a further example of the idealized steady flow treated in Chap.~(\ref{ch:chap_10}), an example that will help show how sources and sinks can be accounted for in problems of fluid flow. In this and the following two sections, the problems are posed in the $uvw$ space, rather than the $xyz$ space. That allows us to refer directly to earlier results in this chapter without interchanging the spaces.

Consider the three-dimensional steady flow of fluid between two parallel planes $v=0$ and $v=\pi$ when the fluid is entering through a narrow slit along the line in the first plane that is perpendicular to the $uv$ plane at the origin. Let the rate of flow of fluid into the channel through the slit be $Q$ units of volume per unit time for each unit of depth of the channel, where the depth is measured perpendicular to the $uv$ plane. The rate of flow out at either end is, then, $Q/2$.

The transformation $\varpi=\log s$ is a one to one mapping of the upper half $y>0$ of the $s$ space onto the strip $0<v<\pi$ in the $\varpi$ space. %(see Example 2 in Sec. 112).
The inverse transformation
\begin{equation}\label{eq:chap_11_5_1}
s=e^{\varpi}=\exp[e_{xy}(u+iv)+w]
\end{equation}
maps the strip onto the half space. %(see Example 3, Sec. 13).
Under transformation~(\ref{eq:chap_11_5_1}), the image of the $u$ axis is the positive half of the $x$ axis, and the image of the line $v=\pi$ is the negative half of the $x$ axis. Hence the boundary of the strip is transformed into the boundary of the half space.

The image of the point $\varpi=0$ is the point $s=1$. The image of a point $\varpi=e_{xy}u_{0}$, where $u_{0}>0$, is a point $s=e_{xy}x_{0}$, where $x_{0}>1$. The rate of flow of fluid across a curve joining the point $\varpi=e_{xy}u_{0}$ to a point $(u,v,w)$ within the strip is a stream function $\psi(u,v,w)$ for the flow (Sec.~(\ref{sec:chap_10_8_107})). If $u_{1}$ is a negative real number, then the rate of flow into the channel through the slit can be written
\[
\psi(u_{1},0,w)=Q.
\]
Now, under a conformal transformation, the function $\psi$ is transformed into a function of $x$, $y$, and $z$ that represents the stream function for the flow in the corresponding region of the $s$ space; that is, the rate of flow is the same across corresponding curves in the two planes. As in Chap.~(\ref{ch:chap_10}), the same symbol $\psi$ is used to represent the different stream functions in the two planes. Since the image of the point $\varpi=e_{xy}u_{1}$ is a point $s=\chi_{1}$, where $0<Re\chi_{1}<1$, the rate of flow across any curve connecting the points $s=\chi_{0}$ and $s=\chi_{1}$ and lying in the upper half of the $s$ space is also equal to $Q$. Hence there is a source at the point $s=1$ equal to the source at $\varpi=0$.

The above argument applies in general to show that under a conformal transformation, a source or sink at a given point corresponds to an equal source or sink at the image of that point.

As $Re\varpi$ tends to $-\infty$, the image of $\varpi$ approaches the point $s=0$. A sink of strength $Q/2$ at the latter point corresponds to the sink infinitely far to the left in the strip. To apply the above argument in this case, we consider the rate of flow across a curve connecting the boundary lines $v=0$ and $v=\pi$ of the left-hand part of the strip and the rate of flow across the image of that curve in the $s$ space.

The sink at the right-hand end of the strip is transformed into a sink at infinity in the $s$ space.

The stream function $\psi$ for the flow in the upper half of the $s$ space in this case must be a function whose values are constant along each of the three parts of the polygon line $P_{xz}$ in the $xz$ coordinate plane. Moreover, its value must increase by $Q$ as the point $s$ moves around the point $s=1$ from the position $s=\chi_{0}$ to the position $s=\chi_{1}$, and its value must decrease by $Q/2$ as $s$ moves about the origin in the corresponding manner. We see that the function
\[
\psi=\frac{Q}{\pi}[\arg_{c}(s-1)-\frac{1}{2}\arg_{c}s]
\]
satisfies those requirements. Furthermore, this function is harmonic in the half space
$Im$ $s>0$ because it is the imaginary component of the function
\[
F=\frac{Q}{\pi}[\log(s-1)-\frac{1}{2}\log s]=\frac{Q}{\pi}\log(s^{1/2}-s^{-1/2}).
\]
The function $F$ is a spatial complex potential function for the flow in the upper half of the
$s$ space. Since $s=e^{\varpi}$, a complex potential function $F(\varpi)$ for the flow in the channel is
\[
F(\varpi)=\frac{Q}{\pi}\log(e^{\varpi/2}-e^{-\varpi/2}).
\]
By dropping an additive constant, one can write
\begin{equation}\label{eq:chap_11_5_2}
F(\varpi)=\frac{Q}{\pi}\log(\sinh\frac{\varpi}{2}).
\end{equation}
We have used the same symbol $F$ to denote three distinct functions, once in the $s$ space and twice in the $\varpi$ space.

The velocity vector $\overline{F'(\varpi)}$ is given by the equation
\begin{equation}\label{eq:chap_11_5_3}
V=\frac{Q}{2\pi}\log(\cosh\frac{\overline{\varpi}}{2}).
\end{equation}
From this, it can be seen that
\[
\lim_{|u| \to \infty}V=\frac{Q}{2\pi}.
\]
Also, the point $\varpi=\pi i$ is a stagnation point; that is, the velocity is zero there. Hence the fluid pressure along the wall $v=\pi$ of the channel is greatest at points opposite the slit.

The stream function $\psi(u,v,w)$ for the channel is the imaginary component of the function $F(\varpi)$ given by  equation~(\ref{eq:chap_11_5_2}). The stream surfaces $\psi(u,v,w)=c_{2}$ are, therefore, the surfaces
\[
\frac{Q}{\pi}\arg_{c}(\sinh\frac{\varpi}{2})=c_{2}.
\]
This equation reduces to
\begin{equation}\label{eq:chap_11_5_4}
\tan(\frac{V}{2})=c\tanh(\frac{u}{2}),
\end{equation}
where $c$ is any real constant.

\section{Flow in a Channel with a Offset}\label{sec:chap_11_6_114}

To further illustrate the use of the Schwarz-Christoffel transformation, let us find the spatial complex potential for the flow of a fluid in a channel with an abrupt change in its breadth. We take our unit of length such that the breadth of the wide part of the channel is $\pi$ units; then $h\pi$, where $0<h<1$, represents the breadth of the narrow part. Let the real constant $V_{0}$ denote the velocity of the fluid far from the offset in the wide part; that is,
\[
\lim_{|u| \to -\infty}V=V_{0},
\]
where the complex variable $V$ represents the velocity vector. The rate of flow per unit depth through the channel,  or the strength of the source on the left and of the sink on the right, is then
\begin{equation}\label{eq:chap_11_6_1}
Q=\pi V_{0}.
\end{equation}

The cross section of the channel can be considered as the limiting case of the quadrilateral with the vertices $\varpi_{l}$, $\varpi_{2}$, $\varpi_{3}$, and $\varpi_{4}$ as the first and last of these vertices are moved infinitely far to the left and to the right, respectively. In the limit, the exterior angles become
\[
k_{1}\pi=\pi,\texttt{ }k_{2}\pi=\frac{\pi}{2},\texttt{ }k_{3}\pi=-\frac{\pi}{2},\texttt{ }k_{4}\pi=\pi.
\]
As before, we proceed formally, using limiting values whenever it is convenient to do so. If we write $Re\chi_{1}=0$, $Re\chi_{3}=1$, $Re\chi_{4}=\infty$ and leave $\chi_{2}$ to be determined, where $0<Re\chi_{2}<1$,
the derivative of the mapping function becomes
\begin{equation}\label{eq:chap_11_6_2}
\frac{d\varpi}{ds}=A(s-\chi_{1})^{-1}(s-\chi_{2})^{-1/2}(s-\chi_{3})^{1/2}.
\end{equation}

To simplify the determination of the constants $A$ and $\chi_{2}$ here, we proceed at once to the complex potential  of the flow. The source of the flow in the channel infinitely far to the left corresponds to an equal source at $s=0$ (Sec.~(\ref{sec:chap_11_5_113})). The entire boundary of the cross section of the channel is the image of the polygon line $P_{xz}$ in the $xz$ coordinate plane. In view of equation~(\ref{eq:chap_11_6_1}), then, the function
\begin{equation}\label{eq:chap_11_6_3}
\begin{split}
F &= V_{0}\log s=V_{0}(\ln r+\ln\sin\varphi) \\
&+ e_{xy}V_{0}[\ln r_{L}(\theta,\varphi)-\ln\sin\varphi+i\theta_{L}(\theta,\varphi)].
\end{split}
\end{equation}
for $0<\varphi\leq\pi/2$ and $\theta\neq(2n+1)\pi$ (equation~(\ref{eq:chap_3_3_2b}), Sec.~(\ref{sec:chap_3_3_29})) where $n=0,1,2,\ldots$
\[
r_{L}(\theta,\varphi)=\sqrt{1+\cos\theta\sin2\varphi},\texttt{ }
\theta_{L}(\theta,\varphi)=\arctan\frac{\sin\theta}{\cos\theta+\tan\varphi}+2n\pi,
\]
is the potential for the flow in the upper half of the $s$ space, with the required source at the origin. Here the stream function is  $\psi=V_{0}\theta_{L}$. It increases in value from $0$ to $V_{0}\pi$ over each hemisphere $s=Re^{i\theta_{L}}$ $(0\leq\theta_{L}\leq\pi)$, where $R>0$, as $\theta_{L}$ varies from $0$ to $\pi$. [Compare equation~(\ref{eq:chap_10_8_5}), Sec.~(\ref{sec:chap_10_8_107})]

The spatial complex conjugate of the velocity $V$ in the $\varpi$ space can be written
\[
\overline{V(\varpi)}=\frac{dF}{d\varpi}=\frac{dF}{ds}\frac{ds}{d\varpi}.
\]
Thus, by referring to equations~(\ref{eq:chap_11_6_2}) and~(\ref{eq:chap_11_6_3}), we can see that
\begin{equation}\label{eq:chap_11_6_4}
\overline{V(\varpi)}=\frac{V_{0}}{A}(\frac{s-\chi_{2}}{s-\chi_{3}})^{1/2}.
\end{equation}

At the limiting position of the point $\varpi_{1}$ which corresponds to $Re$ $s=0$, the velocity is the real constant  $V_{0}$. It therefore follows from equation~(\ref{eq:chap_11_6_4}) that
\[
V_{0}=\frac{V_{0}}{A}\sqrt{Re\chi_{2}}.
\]
At the limiting position of $\varpi_{4}$, which corresponds to $s=\infty$, let the real number $V_{4}$ denote the velocity. Now it seems plausible that as a vertical line segment spanning the narrow part of the channel is moved infinitely far to the right, $V$ approaches $V_{4}$ at each point on that segment. We could establish this conjecture as a fact by first finding $\varpi$ as the function of $s$ from equation~(\ref{eq:chap_11_6_2}); but, to shorten our discussion, we assume that this is true, Then, since the flow is steady,
\[
\pi hV_{4}=\pi V_{0}=Q,
\]
or $V_{4}=V_{0}/h$. Letting $s$ tend to infinity in equation~(\ref{eq:chap_11_6_4}), we find that
\[
\frac{V_{0}}{h}=\frac{V_{0}}{A}.
\]
Thus
\begin{equation}\label{eq:chap_11_6_5}
A=h,\texttt{ }Re\chi_{2}=h^{2}
\end{equation}
and
\begin{equation}\label{eq:chap_11_6_6}
\overline{V(\varpi)}=\frac{V_{0}}{h}(\frac{s-\chi_{2}}{s-\chi_{3}})^{1/2}.
\end{equation}

From equation~(\ref{eq:chap_11_6_6}), we know that the magnitude $|V|$ of the velocity becomes infinite at the  corner $\varpi_{3}$ of the offset since it is the image of the point $s$ $(Re$ $s=1)$. Also,
the corner $\varpi_{2}$ is a stagnation point, a point where $V=0$. Along the boundary of the channel, the fluid pressure is, therefore, greatest at $\varpi_{2}$ and least at $\varpi_{3}$.

To write the relation between the potential and the variable $\varpi$, we must integrate
equation~(\ref{eq:chap_11_6_2}), which can now be written
\begin{equation}\label{eq:chap_11_6_7}
\frac{d\varpi}{ds}=\frac{h}{s-\chi_{1}}(\frac{s-\chi_{3}}{s-\chi_{2}})^{1/2}.
\end{equation}
By substituting a new variables, where
\[
\frac{s-\chi_{2}}{s-\chi_{3}}=\varsigma^{2},
\]
one can show that equation~(\ref{eq:chap_11_6_7}) reduces to
\[
\frac{d\varpi}{d\varsigma}=2h[\frac{1}{1-\varsigma^{2}}-\frac{1}{\eta^{2}-\varsigma^{2}}]\texttt{ where }\eta^{2}=\frac{\chi_{2}-\chi_{1}}{\chi_{3}-\chi_{1}}.
\]
Hence
\begin{equation}\label{eq:chap_11_6_8}
\varpi=h\log\frac{1+\varsigma}{1-\varsigma}-\log\frac{\eta+\varsigma}{\eta-\varsigma}\texttt{ where }\eta^{2}=\frac{\chi_{2}-\chi_{1}}{\chi_{3}-\chi_{1}}.
\end{equation}
The constant of integration here is zero because when $s=\chi_{2}$ $(Re$ $\chi_{2}=h^{2})$, the quantity $\varsigma$ is zero and so, therefore, is $\varpi$.

In terms of $\varsigma$, the potential $F$ of equation~(\ref{eq:chap_11_6_3}) becomes
\[
F=V_{0}\log\frac{h^{2}-\varsigma^{2}}{1-\varsigma^{2}};
\]
consequently,
\begin{equation}\label{eq:chap_11_6_9}
\varsigma^{2}=\frac{\exp(F/V_{0})-h^{2}}{\exp(F/V_{0})-1}.
\end{equation}
By substituting $\varsigma$ from this equation into equation~(\ref{eq:chap_11_6_8}), we obtain an implicit relation that defines the potential $F$ as a function of $\varpi$.

\section{Electrostatic Potential about an Edge of a Conducting Plate}\label{sec:chap_11_7_115}

Two parallel conducting plates of infinite extent are kept at the electrostatic potential $V=0$, and a parallel semi-infinite plate, placed midway between them, is kept at the potential $V=1$. The coordinate system and the unit  of length are chosen so that the plates lie in the planes $v=0$, $v=\pi$, and $v=\pi/2$. Let us determine the potential function $V(u, v, w)$ in the region between those plates.

The cross section of that region in the $\varpi$ space has the limiting form of the quadrilateral bounded by the dashed lines $\varpi_{1}\cdots\varpi_{2}\cdots\varpi_{3}\cdots\varpi_{4}\cdots\varpi_{1}$ as the points $\varpi_{1}$ and $\varpi_{3}$ move out to the right and $\varpi_{4}$ to the left. In applying the Schwarz-Christoffel transformation here, we let the point $\chi_{4}$, corresponding to the vertex $\varpi_{4}$, be the point at infinity. We choose the points $\chi_{1}$ $(Re\chi_{1}=-1)$, $\chi_{3}$ $(Re\chi_{3}=1)$ and leave $\chi_{2}$ to be determined. The limiting values of the exterior angles of the quadrilateral are
\[
k_{1}\pi=\pi,\texttt{ }k_{2}\pi=-\pi,\texttt{ }k_{3}\pi=k_{4}\pi=\pi.
\]
Thus
\[
\begin{split}
\frac{d\varpi}{ds} &=A(s-\chi_{1})^{-1}(s-\chi_{2})(s-\chi_{3})^{-1} \\
                   &=A(\frac{s-\chi_{2}}{(s-\chi_{1})(s-\chi_{3})}) \\
                   &=A[\frac{\chi_{2}-\chi_{1}}{\chi_{3}-\chi_{1}}(\frac{1}{s-\chi_{1}})+\frac{\chi_{2}-\chi_{3}}{\chi_{1}-\chi_{3}}(\frac{1}{s-\chi_{3}})].
\end{split}
\]
and so the transformation of the upper half of the $s$  space into the divided strip in the
$\varpi$ space has the form
\begin{equation}\label{eq:chap_11_7_1}
\varpi=A[\frac{\chi_{2}-\chi_{1}}{\chi_{3}-\chi_{1}}\log(s-\chi_{1})+\frac{\chi_{2}-\chi_{3}}{\chi_{1}-\chi_{3}}\log(s-\chi_{3})]+B.
\end{equation}

Let the spatial complex constants
\[
A=e_{xy}(A_{x}+iA_{y})+A_{z},\texttt{ }B=e_{xy}(B_{x}+iB_{y})+B_{z};
\]
and
\[
s-\chi_{1}=r_{1}(e^{i\theta_{1}}\cos\varphi_{1}+\sin\varphi_{1})\texttt{ }(\theta_{1}\neq(2n+1)\pi,\texttt{ }0<\varphi_{1}\leq\pi/2),
\]
\[
s-\chi_{3}=r_{3}(e^{i\theta_{3}}\cos\varphi_{3}+\sin\varphi_{3})\texttt{ }(\theta_{3}\neq(2n+1)\pi,\texttt{ }0<\varphi_{3}\leq\pi/2),
\]
for $n=0,1,2,\ldots$. When $s=\chi$, the point $\varpi$ lies on the boundary of the divided strip; and, according to equation~(\ref{eq:chap_11_7_1}),
\begin{equation}\label{eq:chap_11_7_2}
e_{xy}(u+iv)+w=A(\chi_{1})\log(s-\chi_{1})+A(\chi_{3})\log(s-\chi_{3})+B
\end{equation}
\[
=A(\chi_{1})[\ln r_{1}+\ln\sin\varphi_{1}+e_{xy}(\ln r_{L}(\theta_{1},\varphi_{1})+i\theta_{L}(\theta_{1},\varphi_{1})-\ln\sin\varphi_{1})]
\]
\[
+A(\chi_{3})[\ln r_{3}+\ln\sin\varphi_{3}+e_{xy}(\ln r_{L}(\theta_{3},\varphi_{3})+i\theta_{L}(\theta_{3},\varphi_{3})-\ln\sin\varphi_{3})]+B
\]
where
\[
A(\chi_{1})=e_{xy}[\frac{(A_{x}+A_{z}+iA_{y})(x_{2}-x_{1}+z_{2}-z_{1})}{(x_{3}-x_{1})+(z_{3}-z_{1})}
-A_{z}\frac{z_{2}-z_{1}}{z_{3}-z_{1}}]+A_{z}\frac{z_{2}-z_{1}}{z_{3}-z_{1}},
\]
\[
A(\chi_{3})=e_{xy}[\frac{(A_{x}+A_{z}+iA_{y})(x_{2}-x_{3}+z_{2}-z_{3})}{(x_{1}-x_{3})+(z_{1}-z_{3})}
-A_{z}\frac{z_{2}-z_{3}}{z_{1}-z_{3}}]+A_{z}\frac{z_{2}-z_{3}}{z_{1}-z_{3}},
\]
and when $0<\varphi\leq\pi/2$, for $n=0,1,2,\ldots$ (equation~(\ref{eq:chap_3_3_2b}), Sec.~(\ref{sec:chap_3_3_29}))
\[
r_{L}(\theta,\varphi)=\sqrt{1+\cos\theta\sin2\varphi},\texttt{ }
\theta_{L}(\theta,\varphi)=\arctan\frac{\sin\theta}{\cos\theta+\tan\varphi}+2n\pi.
\]

To determine the constants here, we first note that the limiting position of the line segment joining the points  $\varpi_{1}$ and $\varpi_{4}$ is the $u$ axis. That segment is the image of the part of the polygon line $P_{xz}$ in the $xz$ coordinate plane to the left of the point $\chi_{1}$ $(Re\chi_{1}=-1)$; this is because the line segment joining $\varpi_{3}$ and $\varpi_{4}$ is the image of the part of the polygon line $P_{xz}$ in the $xz$ coordinate plane to the right of $\chi_{3}$ $(Re\chi_{3}=1)$, and the other two sides of the quadrilateral are the images of the remaining two segments of the polygon line in the $xz$ coordinate plane. Hence when $v=0$ and $u$ tends to infinity through positive values, the corresponding point $\chi$ approaches the point $s=\chi_{1}$ from the left. Thus
\[
\arg_{c}(\chi-\chi_{1})=\pi,\texttt{ }\arg_{c}(\chi-\chi_{3})=\pi,
\]
and $\ln|\chi-\chi_{1}|$ tends to $-\infty$. Also, since $-1<Re\chi_{2}<1$, the real part in the $xy$ coordinate plane of the quantity inside the braces in equation~(\ref{eq:chap_11_7_2}) tends to $-\infty$. Since $v=0$, it readily follows that $A_{y}=0$; for, otherwise, the imaginary part on the right would become infinite. By equating imaginary parts on the two sides, we now see that
\[
0=\frac{A_{y}}{2}Re[(\chi_{2}-\chi_{1})\pi-(\chi_{2}-\chi_{3})\pi]+B_{y}.
\]
Hence
\begin{equation}\label{eq:chap_11_7_3}
-A_{x}\pi=B_{y},\texttt{ }A_{y}=0.
\end{equation}
The limiting position of the line segment joining the points $\varpi_{1}$ and $\varpi_{2}$ is the half
line $v=\pi/2$ $(u\geq0)$. Points on that half line are images of the points $s=\chi$, where $-1<Re\chi\leq Re\chi_{2}$; consequently,
\[
\arg_{c}(\chi-\chi_{1})=0,\texttt{ }\arg_{c}(\chi-\chi_{3})=\pi.
\]
Identifying the imaginary parts on the two sides of equation~(\ref{eq:chap_11_7_2}), we thus arrive at the relation
\begin{equation}\label{eq:chap_11_7_4}
\frac{\pi}{2}=\frac{A_{x}}{2}Re(\chi_{3}-\chi_{2})\pi+B_{y}.
\end{equation}

Finally, the limiting positions of the points on the line segment joining $\varpi_{3}$ to $\varpi_{4}$
are the points $u+\pi i$, which are the images of the points $\chi$ when $Re\chi>1$. By identifying, for those points, the imaginary parts in equation~(\ref{eq:chap_11_7_2}), we find that
\[
\pi=B_{y}.
\]
Then, in view of equations~(\ref{eq:chap_11_7_3}) and~(\ref{eq:chap_11_7_4}),
\[
A_{x}=-1,\texttt{ }Re\chi_{2}=0.
\]
Thus $Re\chi=0$ is the point whose image is the vertex $\varpi=i\pi/2$; and, upon substituting these values into equation~(\ref{eq:chap_11_7_2}) and identifying real parts in the $xy$ coordinate plane, we see that $B_{x}=0$.

Transformation~(\ref{eq:chap_11_7_2}) now becomes
\begin{equation}\label{eq:chap_11_7_5}
\varpi=-\frac{1}{2}[\log(s-\chi_{1})+\log(s-\chi_{3})]+i\pi,
\end{equation}
or
\begin{equation}\label{eq:chap_11_7_6}
s=\frac{\chi_{1}+\chi_{3}}{2}\pm\sqrt{e^{-2\varpi}+(\frac{\chi_{1}-\chi_{3}}{2})^{2}}.
\end{equation}

Under this transformation, the required harmonic function $V(u, v, w)$ becomes a harmonic function of $x$, $y$, and  $z$ in the half space $y>0$; and the boundary conditions are satisfied. Note that $Re\chi_{2}=0$ now. The harmonic function in that half space which assumes those values on the boundary is the imaginary component of the analytic function
\[
\frac{1}{\pi}\log\frac{s-\chi_{3}}{s-\chi_{1}}=\frac{1}{\pi}(\ln\frac{r_{3}}{r_{1}}+\ln\frac{\sin\varphi_{3}}{\sin\varphi_{1}})
\]
\[
+e_{xy}\frac{1}{\pi}\{[\ln \frac{r_{L}(\theta_{3},\varphi_{3})}{r_{L}(\theta_{1},\varphi_{1})}-\ln\frac{\sin\varphi_{3}}{\sin\varphi_{1}}]
+i[\theta_{L}(\theta_{3},\varphi_{3})-\theta_{L}(\theta_{1},\varphi_{1})]\},
\]
where $\theta_{1}$ and $\theta_{3}$ range from $0$ to $\pi$. Writing the tangents of these angles as functions
of $x$, $y$, and $z$  and simplifying, we find that
\begin{equation}\label{eq:chap_11_7_7}
\tan\pi V=\tan(\theta_{3}-\theta_{1})=\frac{2y}{x^{2}+y^{2}+z^{2}-1}.
\end{equation}

Equation~(\ref{eq:chap_11_7_6}) furnishes expressions for $x^{2}+y^{2}+z^{2}$ and $x^{2}-y^{2}+z^{2}$ in terms of  $u$, $v$, and $w$. Then, from equation~(\ref{eq:chap_11_7_7}), we find that the relation between the potential $V$  and the coordinates $u$, $v$, and $w$ can be written
\begin{equation}\label{eq:chap_11_7_8}
\tan\pi V=\frac{1}{\varsigma}\sqrt{e^{-4u}-\varsigma^{2}},
\end{equation}
where
\[
\varsigma=-1+\sqrt{1+2e^{-2u}\cos2v+e^{-4u}}.
\]

%-----------------------------------------------------------------------
% Beginning of chap12.tex
%-----------------------------------------------------------------------
%
% AMS-LaTeX 1.2 sample file for a monograph, based on amsbook.cls.
% This is a data file input by chapter.tex.
%%%%%%%%%%%%%%%%%%%%%%%%%%%%%%%%%%%%%%%%%%%%%%%%%%%%%%%%%%%%%%%%%%%%

%\part{This is a Part Title Sample}

\chapter{Circle Integral Formulas of the Poisson Type}\label{ch:chap_12}

In this chapter, we develop a theory that enables us to obtain solutions to a variety of boundary value problems where those solutions are expressed in terns of definite or improper integrals. Many of the integrals occurring are then readily evaluated.

\section{Poisson Circle Integral Formula}\label{sec:chap_12_1_116}

Let $C_{0}$ denote a positively oriented circle on a sphere $S_{0}$, centered at the origin, and suppose that a
function $f$ is analytic inside and on $S_{0}$ and $C_{0}$. The Cauchy integral formula (Sec.~(\ref{sec:chap_4_12_47}))
\begin{equation}\label{eq:chap_12_1_1}
f(s)=\frac{1}{2\pi i}\int_{C_{0}}\frac{f(\varsigma)d\varsigma}{\varsigma-s}.
\end{equation}
expresses the value of $f$ at any point $s$ interior to $C_{0}$ in terms of the values of $f$ at points on $C_{0}$ where $s$ and $C_{0}$ are in a same spatial plane. In this section, we shall obtain from formula~(\ref{eq:chap_12_1_1}) a corresponding formula for the real part of the function $f$;  and, in Sec.~(\ref{sec:chap_12_2_117}), we shall use that result to solve the Dirichlet problem (Sec.~(\ref{sec:chap_9_5_98})) for the disk bounded by $C_{0}$.

We let $r_{0}$ denote the radius of $C_{0}$ and write
\[
s=r(e_{xy}e^{i\theta}\cos\varphi+\sin\varphi)\texttt{ }(0<r<r_{0}\texttt{, }0\leq\theta\leq2\pi\texttt{, }-\pi/2\leq\varphi\leq\pi/2).
\]
The inverse of the nonzero point $s$ with respect to the circle is the point $s_{1}$ lying on the same ray from the origin as $s$ and satisfying the condition $|s_{1}||s|=r_{0}^{2}$; thus, if $\varsigma$ is a point on $C_{0}$,
\begin{equation}\label{eq:chap_12_1_2}
s_{1}=\frac{r_{0}^{2}}{r}(e_{xy}e^{i\theta}\cos\varphi+\sin\varphi)=\frac{r_{0}^{2}}{r^{2}}s.
\end{equation}
Since $s_{1}$ is exterior to the circle $C_{0}$, it follows from the Cauchy-Goursat theorem that the value of the integral in equation~(\ref{eq:chap_12_1_1}) is zero when $s$ is replaced by $s_{1}$ in the integrand. Hence
\[
f(s)=\frac{1}{2\pi i}\int_{C_{0}}(\frac{1}{\varsigma-s}-\frac{1}{\varsigma-s_{1}})f(\varsigma)d\varsigma;
\]
and, using the parametric representation
\[
\varsigma=r_{0}(e_{xy}e^{i\phi}\cos\psi+\sin\psi)\texttt{ }(0\leq\phi\leq2\pi\texttt{, }-\pi/2\leq\psi\leq\pi/2)\texttt{ for }C_{0},
\]
because of $d\varsigma=i\varsigma d\phi$ when $\varsigma$ is a point on $C_{0}$, we can write
\[
f(s)=\frac{1}{2\pi}\int_{0}^{2\pi}(\frac{\varsigma}{\varsigma-s}-\frac{\varsigma}{\varsigma-s_{1}})f(\varsigma)d\phi;
\]
where, for convenience, we retain the $\varsigma$ to denote $\varsigma=r_{0}(e_{xy}e^{i\phi}\cos\psi+\sin\psi)$.

In view of the last of expressions~(\ref{eq:chap_12_1_2}) for $s_{1}$, the factor inside the parentheses here can be written
\begin{equation}\label{eq:chap_12_1_3}
\frac{\varsigma}{\varsigma-s}-\frac{\varsigma}{\varsigma-\frac{r_{0}^{2}}{r^{2}}s}
=\frac{1}{1-s/\varsigma}-\frac{1}{1-\frac{r_{0}^{2}}{r^{2}}s/\varsigma}
=\frac{r_{0}^{2}-r^{2}}{(1-\varsigma/s)(r^{2}-r_{0}^{2}s/\varsigma)}.
\end{equation}

An alternative form of the Cauchy integral formula~(\ref{eq:chap_12_1_1}) is, therefore,
\begin{equation}\label{eq:chap_12_1_4}
f[r(e_{xy}e^{i\theta}\cos\varphi+\sin\varphi)]=\frac{r_{0}^{2}-r^{2}}{2\pi}\int_{0}^{2\pi}\frac{f[r_{0}(e_{xy}e^{i\phi}\cos\psi+\sin\psi)]}{(1-\varsigma/s)(r^{2}-r_{0}^{2}s/\varsigma)}d\phi
\end{equation}
when $0<r<r_{0}$. This form is also valid when $r=0$; in that case, it reduces directly to
\[
f(0)=\frac{1}{2\pi}\int_{0}^{2\pi}f(r_{0}e^{i\phi})d\phi,
\]
which is just the parametric form of equation~(\ref{eq:chap_12_1_1}) with $s=0$.

From formula~(\ref{eq:chap_1_7_2}) in Sec.~(\ref{sec:chap_1_7}), the quotients $(r\varsigma)/(r_{0}s)$ and $(r_{0}s)/(r\varsigma)$ are
\[
\frac{r\varsigma}{r_{0}s}=\{
    \begin{array}{cc}
        e_{xy}(\frac{e^{i\phi}\cos\psi+\sin\psi}{e^{i\theta}\cos\varphi+\sin\varphi}
-\frac{\sin\psi}{\sin\varphi})+\frac{\sin\psi}{\sin\varphi} & \texttt{ when }\sin\varphi\neq0 \\
        e_{xy}e^{i(\phi-\theta)} & \texttt{ when }\sin\varphi=0
    \end{array}
\]
and
\[
\frac{r_{0}s}{r\varsigma}=\{
    \begin{array}{cc}
        e_{xy}(\frac{e^{i\theta}\cos\varphi+\sin\varphi}{e^{i\phi}\cos\psi+\sin\psi}
-\frac{\sin\varphi}{\sin\psi})+\frac{\sin\varphi}{\sin\psi} & \texttt{ when }\sin\psi\neq0 \\
        e_{xy}e^{i(\theta-\phi)} & \texttt{ when }\sin\psi=0
    \end{array},
\]
respectively. Then we can write
\[
(1-\varsigma/s)(r^{2}-r_{0}^{2}s/\varsigma)=r_{0}^{2}-r_{0}r(\frac{r\varsigma}{r_{0}s}+\frac{r_{0}s}{r\varsigma})+r^{2}
=e_{xy}(a-b)+(b+c)
\]
where
\[
\begin{split}
a&=-2r_{0}r\cos(\phi-\theta),\texttt{ }
b=0\texttt{ when }\sin\psi=0,\texttt{ }\sin\varphi=0,\\
a&=-r_{0}r(e^{i(\phi-\theta)}+\frac{e^{i\theta}}{e^{i\phi}\cos\psi+\sin\psi}),\texttt{ }
b=0\texttt{ when }\sin\psi\neq0,\texttt{ }\sin\varphi=0,\\
a&=-r_{0}r(\frac{e^{i\phi}}{e^{i\theta}\cos\varphi+\sin\varphi}+e^{i(\theta-\phi)}),\texttt{ }
b=0\texttt{ when }\sin\psi=0,\texttt{ }\sin\varphi\neq0,\\
a&=-r_{0}r(\frac{e^{i\phi}\cos\psi+\sin\psi}{e^{i\theta}\cos\varphi+\sin\varphi}+\frac{e^{i\theta}\cos\varphi+\sin\varphi}{e^{i\phi}\cos\psi+\sin\psi}),\texttt{ }
b=-r_{0}r(\frac{\sin\psi}{\sin\varphi}+\frac{\sin\varphi}{\sin\psi})
\end{split}
\]
\[
\texttt{ when }\sin\psi\sin\varphi\neq0,\texttt{ and }c=r_{0}^{2}+r^{2},
\]
and from formula~(\ref{eq:chap_1_7_3}) in Sec.~(\ref{sec:chap_1_7})
\[
\frac{1}{(1-\varsigma/s)(r^{2}-r_{0}^{2}s/\varsigma)}=\frac{1}{e_{xy}(a-b)+(b+c)}=e_{xy}(\frac{1}{a+c}-\frac{1}{b+c})+\frac{1}{b+c}
\]
\[
=e_{xy}[\frac{1}{r_{0}^{2}-2r_{0}r\sigma(\phi,\theta,\psi,\varphi)+r^{2}}
-\frac{1}{r_{0}^{2}-2r_{0}r\sigma(\psi,\varphi)+r^{2}}]
\]
\[
+\frac{1}{r_{0}^{2}-2r_{0}r\sigma(\psi,\varphi)+r^{2}}
\]
or
\begin{equation}\label{eq:chap_12_1_5}
\frac{1}{(1-\varsigma/s)(r^{2}-r_{0}^{2}s/\varsigma)}=\frac{1}{r_{0}^{2}-2r_{0}r\sigma(\psi,\varphi)+r^{2}}
\end{equation}
\[
+e_{xy}[\frac{1}{r_{0}^{2}-2r_{0}r\sigma(\phi,\theta,\psi,\varphi)+r^{2}}
-\frac{1}{r_{0}^{2}-2r_{0}r\sigma(\psi,\varphi)+r^{2}}]
\]
where when $\sin\psi\sin\varphi=0$, there are $\sigma(\psi,\varphi)=0$, and
\[
\begin{split}
\sigma(\phi,\theta,\psi,\varphi)&=\cos(\phi-\theta)\texttt{ when }\sin\psi=0,\texttt{ }\sin\varphi=0,\\
\sigma(\phi,\theta,\psi,\varphi)&=\frac{1}{2}(e^{i(\phi-\theta)}+\frac{e^{i\theta}}{e^{i\phi}\cos\psi+\sin\psi})
\texttt{ when }\sin\psi\neq0,\texttt{ }\sin\varphi=0,\\
\sigma(\phi,\theta,\psi,\varphi)&=\frac{1}{2}(\frac{e^{i\phi}}{e^{i\theta}\cos\varphi+\sin\varphi}+e^{i(\theta-\phi)})
\texttt{ when }\sin\psi=0,\texttt{ }\sin\varphi\neq0,
\end{split}
\]
and when $\sin\psi\sin\varphi\neq0$
\[
\sigma(\psi,\varphi)=\frac{1}{2}(\frac{\sin\psi}{\sin\varphi}+\frac{\sin\varphi}{\sin\psi})
\]
and
\[
\sigma(\phi,\theta,\psi,\varphi)=\frac{1}{2}(\frac{e^{i\phi}\cos\psi+\sin\psi}{e^{i\theta}\cos\varphi+\sin\varphi}+\frac{e^{i\theta}\cos\varphi+\sin\varphi}{e^{i\phi}\cos\psi+\sin\psi})
\]
\[
=\frac{\cos(\phi-\theta)+\cos\phi\tan\varphi+\cos\theta\tan\psi+\tan\psi\tan\varphi}{(1+\cos\phi\sin2\psi)(1+\cos\theta\sin2\varphi)}\cos\psi\cos\varphi
\]
\[
+\frac{e^{i(\phi-\theta)}+e^{i\phi}\tan\varphi+e^{-i\theta}\tan\psi
+\tan\psi\tan\varphi}{(1+\cos\phi\sin2\psi)(1+\cos\theta\sin2\varphi)}\cos\phi\sin\psi\cos^{2}\psi\cos\varphi
\]
\[
+\frac{e^{-i(\phi-\theta)}+e^{-i\phi}\tan\varphi+e^{i\theta}\tan\psi
+\tan\psi\tan\varphi}{(1+\cos\phi\sin2\psi)(1+\cos\theta\sin2\varphi)}\cos\theta\cos\psi\sin\varphi\cos^{2}\varphi.
\]

Hence, if $u$ is the real part of the analytic function $f$, it follows from formula~(\ref{eq:chap_12_1_4}) that
\begin{equation}\label{eq:chap_12_1_6}
u(r,\theta,\varphi)=\frac{1}{2\pi}\int_{0}^{2\pi}Re\frac{(r_{0}^{2}-r^{2})}{(1-\varsigma/s)(r^{2}-r_{0}^{2}s/\varsigma)}u(r_{0},\phi,\psi)d\phi\texttt{ }(r<r_{0})
\end{equation}
\[
=\frac{1}{2\pi}\int_{0}^{2\pi}[Re\frac{(r_{0}^{2}-r^{2})u(r_{0},\phi,\psi)}{r_{0}^{2}-2r_{0}r\sigma(\phi,\theta,\psi,\varphi)+r^{2}}
-\frac{(r_{0}^{2}-r^{2})u(r_{0},\phi,\psi)}{r_{0}^{2}-2r_{0}r\sigma(\psi,\varphi)+r^{2}}]d\phi\texttt{ }(r<r_{0}).
\]
This is the Poisson circle integral formula for the harmonic function $u$ in the open disk bounded by the circle $r=r_{0}$.

Formula~(\ref{eq:chap_12_1_6}) defines a linear integral transformation of $u(r_{0},\phi,\psi)$ into $~(\ref{eq:chap_12_1_4})$. The kernel of the transformation is, except for the factor $1/(2\pi)$, the real-valued function
\begin{equation}\label{eq:chap_12_1_7}
\begin{split}
P(r_{0},r,\phi,\theta,\psi,\varphi)&=Re\frac{r_{0}^{2}-r^{2}}{r_{0}^{2}-2r_{0}r\sigma(\phi,\theta,\psi,\varphi)+r^{2}}\\
&-\frac{r_{0}^{2}-r^{2}}{r_{0}^{2}-2r_{0}r\sigma(\psi,\varphi)+r^{2}},
\end{split}
\end{equation}
which is known as the Poisson kernel for the circle integral formula. In view of equation~(\ref{eq:chap_12_1_5}), we can also write
\begin{equation}\label{eq:chap_12_1_8}
\begin{split}
P(r_{0},r,\phi,\theta,\psi,\varphi)&=(r_{0}^{2}-r^{2})\frac{r_{0}^{2}-2r_{0}r\sigma_{r}(\phi,\theta,\psi,\varphi)+r^{2}}{|r_{0}^{2}-2r_{0}r\sigma(\phi,\theta,\psi,\varphi)+r^{2}|^{2}}\\
&-\frac{r_{0}^{2}-r^{2}}{r_{0}^{2}-2r_{0}r\sigma(\psi,\varphi)+r^{2}}
\end{split}
\end{equation}
where
\[
\sigma_{r}(\phi,\theta,\psi,\varphi)=\frac{\cos(\phi-\theta)+\cos\phi\tan\varphi+\cos\theta\tan\psi+\tan\psi\tan\varphi}{(1+\cos\phi\sin2\psi)(1+\cos\theta\sin2\varphi)}\cos\psi\cos\varphi
\]
\[
+\frac{\cos(\phi-\theta)+\cos\theta\tan\psi
+\cos\phi\tan\varphi+2\tan\psi\tan\varphi}{2(1+\cos\theta\sin2\varphi)(1+\cos\phi\sin2\psi)}\cos\phi\sin\psi\cos^{2}\psi\cos\varphi
\]
\[
+\frac{\cos(\phi-\theta)+\cos\phi\tan\varphi+\cos\theta\tan\psi
+2\tan\psi\tan\varphi}{2(1+\cos\phi\sin2\psi)(1+\cos\theta\sin2\varphi)}\cos\theta\cos\psi\sin\varphi\cos^{2}\varphi.
\]

Since $r<r_{0}$, it is clear that $P$ is a positive function. Moreover, we find from the second
of equations~(\ref{eq:chap_12_1_3}) that
\begin{equation}\label{eq:chap_12_1_9}
P(r_{0},r,\phi,\theta,\psi,\varphi)=Re\frac{(r_{0}^{2}-r^{2})}{(1-\varsigma/s)(r^{2}-r_{0}^{2}s/\varsigma)},
\end{equation}
and from equation~(\ref{eq:chap_12_1_7}) that because $P$ consists of trigonometric functions, $P$ is analytic entirely in $S_{0}$. Thus $P(r_{0},r,\phi,\theta,\psi,\varphi)$ is a harmonic function of $r$, $\theta$, and $\varphi$ interior to $C_{0}$ for each fixed $\varsigma$ on $C_{0}$. From equation~(\ref{eq:chap_12_1_8}), we see that $P(r_{0},r,\phi,\theta,\psi,\varphi)$ is an even periodic function of $\phi$ and $\theta$, with period $2\pi$; and its value is $1$ when $r=0$.

The Poisson integral formula~(\ref{eq:chap_12_1_6}) can now be written
\begin{equation}\label{eq:chap_12_1_10}
u(r,\theta,\varphi)=\frac{1}{2\pi}\int_{0}^{2\pi}P(r_{0},r,\phi,\theta,\psi,\varphi)u(r_{0},\phi,\psi)d\phi\texttt{ }(r<r_{0}).
\end{equation}

When $f(s)=u(r,\theta,\varphi)=1$, equation~(\ref{eq:chap_12_1_10}) shows that $P$ has the property
\begin{equation}\label{eq:chap_12_1_11}
\frac{1}{2\pi}\int_{0}^{2\pi}P(r_{0},r,\phi,\theta,\psi,\varphi)d\phi=1\texttt{ }(r<r_{0}).
\end{equation}

We have assumed that $f$ is analytic not only interior to $C_{0}$ but also on $C_{0}$ itself and that $u$ is, therefore, harmonic in a domain which includes all points on that circle. In particular, $u$ is continuous on $C_{0}$. The conditions will now be relaxed.

\section{Dirichlet problem for a disk}\label{sec:chap_12_2_117}

Let $F$ be a piecewise continuous function of $\theta$ on the interval $0\leq\theta\leq 2\pi$ and $\varphi$ on the interval $-\pi/2\leq\varphi\leq\pi/2$. The Poisson integral transform of $F$ is defined in terms of the Poisson kernel $P(r_{0},r,\phi,\theta,\psi,\varphi)$, introduced in Sec.~(\ref{sec:chap_12_1_116}), by means of the equation
\begin{equation}\label{eq:chap_12_2_1}
U(r,\theta,\varphi)=\frac{1}{2\pi}\int_{0}^{2\pi}P(r_{0},r,\phi,\theta,\psi,\varphi)F(\phi,\psi)d\phi\texttt{ }(r<r_{0}).
\end{equation}

In this section, we shall prove that the function $U(r,\theta,\varphi)$ is complex harmonic inside the sphere $r=r_{0}$ and
\begin{equation}\label{eq:chap_12_2_2}
\lim_{r\rightarrow r_{0}, r<r_{0}}U(r,\theta,\varphi)=F(\theta,\varphi)
\end{equation}
for each fixed pair of $\theta$ and $\varphi$ at which $F$ is continuous. Thus $U$ is a solution of the Dirichlet problem for the sphere $r<r_{0}$ in the sense that $U(r,\theta,\varphi)$ approaches the the boundary value $F(\theta,\varphi)$ as the point $(r,\theta,\varphi)$ approaches $(r_{0},\theta,\varphi)$ along a radius, except at the finite number of points $(r_{0},\theta,\varphi)$ where discontinuities of $F$ may occur.

We turn now to the proof that the function $U$ defined in equation~(\ref{eq:chap_12_2_2}) satisfies the Dirichlet problem for the sphere $r<r_{0}$. First of all, $U$ is complex harmonic inside the sphere $r=r_{0}$ because $P$ is a complex harmonic function of $r$, $\theta$, and $\varphi$ there. More precisely, since $F$ is piecewise continuous, integral~(\ref{eq:chap_12_2_1}) may be written as the sum of a finite number of definite integrals each of which has an integrand that is continuous in $r,\theta,\varphi,\phi$ and $\psi$. The partial derivatives of those integrands with respect to $r$, $\theta$, and $\varphi$ are also continuous. Since the order of integration and differentiation with respect to $r$, $\theta$, and $\varphi$ can, then, be interchanged and since $P$ satisfies complex Laplace's equation (see Theorem~(\ref{th:chap_2_15_1s}) in Sec.~(\ref{sec:chap_2_15_25}))
\[
r\frac{\partial}{\partial r}(rP_{r}\cos^{2}\varphi)+\frac{P_{\theta\theta}}{\cos^{2}\varphi}
+\tan\varphi\frac{\partial}{\partial\varphi}(P_{\varphi}\sin\varphi\cos\varphi)
-2\sin\varphi\frac{\partial}{\partial\varphi}(rP_{r}\cos\varphi)=0.
\]
in the spherical coordinates $r$, $\theta$, and $\varphi$, it follows that $U$ satisfies that equation too.

In order to verify limit~(\ref{eq:chap_12_2_2}), we need to show that if $F$ is continuous at $\theta$ and $\varphi$, there corresponds to each positive number $\varepsilon$ a positive number $\delta$ such that
\begin{equation}\label{eq:chap_12_2_6}
|U(r,\theta,\varphi)-F(\theta,\varphi)|<\varepsilon \texttt{ whenever } 0<r_{0}-r<\delta.
\end{equation}
We start by referring to property~(\ref{eq:chap_12_1_11}), Sec.~(\ref{sec:chap_12_1_116}), of the Poisson kernel and writing
\[
U(r,\theta,\varphi)-F(\theta,\varphi)=\frac{1}{2\pi}\int_{0}^{2\pi}P(r_{0},r,\phi,\theta,\psi,\varphi)[F(\phi,\psi)-F(\theta,\varphi)]d\phi.
\]
For convenience, we let $F$ be extended periodically, with period $2\pi$, so that the integrand here is periodic in $\phi$ with that same period. Also, we may assume that $0<r<r_{0}$ because of the nature of the limit to be established.

Next, we observe that, since $F$ is continuous at $\theta$, there is a small positive number $\alpha$ such that
\begin{equation}\label{eq:chap_12_2_7}
|F(\phi,\psi)-F(\theta,\varphi)|<\frac{\varepsilon}{2} \texttt{ whenever } |\phi-\theta|\leq\alpha.
\end{equation}
Evidently,
\begin{equation}\label{eq:chap_12_2_8}
U(r,\theta,\varphi)-F(\theta,\varphi)=I_{1}(r)+I_{2}(r),
\end{equation}
where
\[
\begin{split}
I_{1}(r)&=\frac{1}{2\pi}
\int_{\theta-\alpha}^{\theta+\alpha}P(r_{0},r,\phi,\theta,\psi,\varphi)[F(\phi,\psi)-F(\theta,\varphi)]d\phi,\\
I_{2}(r)&=\frac{1}{2\pi}
\int_{\theta+\alpha}^{\theta-\alpha+2\pi}P(r_{0},r,\phi,\theta,\psi,\varphi)[F(\phi,\psi)-F(\theta,\varphi)]d\phi.
\end{split}
\]

The fact that $P$ is a positive function (See last section), together with the first of inequalities~(\ref{eq:chap_12_2_7}) just above and property~(\ref{eq:chap_12_1_11}) of that function, enable us to write
\[
\begin{split}
|I_{1}(r)|&\leq\frac{1}{2\pi}
\int_{\theta-\alpha}^{\theta+\alpha}P(r_{0},r,\phi,\theta,\psi,\varphi)|F(\phi,\psi)-F(\theta,\varphi)|d\phi,\\
&<\frac{\varepsilon}{4\pi}
\int_{0}^{2\pi}P(r_{0},r,\phi,\theta,\psi,\varphi)d\phi=\frac{\varepsilon}{2}
\end{split}
\]
As for the integral $I_{2}(r)$, one may know that the denominator $r_{0}^{2}-2r_{0}r\sigma(\phi,\theta,\psi,\varphi)+r^{2}$ in expression~(\ref{eq:chap_12_2_8}) for $P(r_{0},r,\phi,\theta,\psi,\varphi)$ in Sec.~(\ref{sec:chap_12_1_116}) has a (positive) minimum value $m$ as the argument $\phi$ of $\varsigma$ varies over the closed interval
\[
\theta+\alpha\leq\phi\leq\theta-\alpha+2\pi.
\]
So, if $M$ denotes an upper bound of the piecewise continuous function $|F(\phi,\psi)-F(\theta,\varphi)|$ on the interval $0\leq\phi\leq2\pi$, it follows that
\[
|I_{2}(r)|\leq\frac{(r_{0}^{2}-r^{2})M}{2\pi m}2\pi<\frac{2M^{2}R}{m}(r_{0}-r)<\frac{2Mr_{0}}{m}\delta=\frac{\varepsilon}{2}
\]
whenever $(r_{0}-r)<\delta$, where
\begin{equation}\label{eq:chap_12_2_9}
\delta=\frac{m\varepsilon}{4Mr_{0}}.
\end{equation}

Finally, the results in the two preceding paragraphes tell us that
\[
|U(r,\theta,\varphi)-F(\theta,\varphi)|\leq|I_{1}(r)|+|I_{2}(r)|<\frac{\varepsilon}{2}+\frac{\varepsilon}{2}=\varepsilon
\]
whenever $(r_{0}-r)<\delta$, where $\delta$ is the positive number defined by equation~(\ref{eq:chap_12_2_9}). That is, statement~(\ref{eq:chap_12_2_6}) holds when that choice of $\delta$ is made.

According to expression~(\ref{eq:chap_12_2_1}), the value of $U$ at $r=0$ is
\[
\frac{1}{2\pi}\int_{0}^{2\pi}F(\phi,\psi)d\phi.
\]
Thus the value of a complex harmonic function at the center of the circle $r=r_{0}$ is the average of the boundary values on the circle.

\section{Related Boundary Value Problems}\label{sec:chap_12_3_118}

The function $F$ representing boundary values on the sphere $r=r_{0}$ is assumed to be piecewise continuous.

Suppose that $F(2\pi-\theta,\varphi)=-F(\theta,\varphi)$. The Poisson integral formula~(\ref{eq:chap_12_2_1}) of Sec.~((\ref{sec:chap_12_2_117}) then becomes
\begin{equation}\label{eq:chap_12_3_1}
U(r,\theta,\varphi)=\frac{1}{2\pi}
[P(r_{0},r,\phi,\theta,\psi,\varphi)-P(r_{0},r,\phi,-\theta,\psi,\varphi)]F(\phi,\psi)d\phi.
\end{equation}
This function $U$ has zero values on the radii planes $\theta=0$ and $\theta=\pi$ of the circle, which are on the $xz$ coordinates plane. %, as one would expect when $U$ is interpreted as a steady temperature.
Formula~(\ref{eq:chap_12_3_1}) thus solve the Dirichlet problem for the semi-circle region $r<r_{0}, 0<\theta<\pi$, where $U=0$ on the disk $r<r_{0}$ which is denoted by $C_{R}$ and on the $xz$ coordinates plane, and
\begin{equation}\label{eq:chap_12_3_2}
\lim_{r\rightarrow r_{0}, r<r_{0}}U(r,\theta,\varphi)=F(\theta,\varphi), (0<\theta<\pi)
\end{equation}
for each fixed $\theta$ at which $F$ is continuous.

If $F(2\pi-\theta,\varphi)=F(\theta,\varphi)$, then
\begin{equation}\label{eq:chap_12_3_3}
U(r,\theta,\varphi)=\frac{1}{2\pi}
[P(r_{0},r,\phi,\theta,\psi,\varphi)+P(r_{0},r,\phi,-\theta,\psi,\varphi)]F(\phi,\psi)d\phi;
\end{equation}
and $U_{\theta}=0$ when $\theta=0$ or $\theta=\pi$. Hence formula~(\ref{eq:chap_12_3_3}) furnishes a function $U$ that is complex harmonic in the semi-circle region $r<r_{0}, 0<\theta<\pi$ and satisfies condition~(\ref{eq:chap_12_3_2}) as well as the condition that its normal derivative be zero on the disk $C_{R}$.

The analytic function $s=r_{0}^{2}/s_{e}$ maps the sphere $|s_{e}|=r_{0}$ in the $s_{e}$ space onto the sphere $|s|=r_{0}$ in the $s$ space, and it maps the exterior of the first sphere onto the interior of the second. Writing
\[
s=r(e_{xy}\exp(i\theta)\cos\varphi+\sin\varphi)\texttt{ and } s_{e}=r_{e}(e_{xy}\exp(i\theta_{e})\cos\varphi_{e}+\sin\varphi_{e}),
\]
we note that $r=r_{0}^{2}/r_{e}$ and $\theta=2\pi-\theta_{e}$. The complex harmonic function $U$ represented by formula~(\ref{eq:chap_12_2_1}), Sec.~(\ref{sec:chap_12_2_117}), is, then, transformed into the function
\[
U(\frac{r_{0}^{2}}{r_{e}},2\pi-\theta_{e},\varphi_{e})=-\frac{1}{2\pi}\int_{0}^{2\pi}P(r_{0},r_{e},\phi,-\theta_{e},\psi,\varphi_{e})F(\phi,\psi)d\phi,
\]
which is complex harmonic in the domain $r_{e}>r_{0}$. Now, in general, if $\tau(r,\theta,\varphi)$ is complex harmonic, then so is $\tau(r,-\theta,\varphi)$. Hence the function $H(r_{e},\theta_{e},\varphi_{e})=U(\frac{r_{0}^{2}}{r_{e}},2\pi-\theta_{e},\varphi_{e})$, or
\begin{equation}\label{eq:chap_12_3_4}
H(r_{e},\theta_{e},\varphi_{e})=-\frac{1}{2\pi}\int_{0}^{2\pi}P(r_{0},r_{e},\phi,\theta_{e},\psi,\varphi_{e})F(\phi,\psi)d\phi, (r_{e}>r_{0}),
\end{equation}
is also complex harmonic. For each fixed $\theta_{e}$ at which $F(\theta_{e},\varphi_{e})$ is continuous, we find from condition~(\ref{eq:chap_12_2_2}), Sec.~(\ref{sec:chap_12_2_117}), that
\begin{equation}\label{eq:chap_12_3_5}
\lim_{r_{e}\rightarrow r_{0}, r_{e}>r_{0}}H(r_{e},\theta_{e},\varphi_{e})=F(\theta_{e},\varphi_{e}).
\end{equation}

The formula~(\ref{eq:chap_12_3_4}) solves the Dirichlet problem for the region exterior to the sphere $r_{e}=r_{0}$ in the $s_{e}$ space. We note from expression~(\ref{eq:chap_12_1_7}), Sec.~(\ref{sec:chap_12_1_116}), that the Poisson kernel $P(r_{0},r_{e},\phi,\theta_{e},\psi,\varphi_{e})$ is negative when $r_{e}>r_{0}$. Also,
\begin{equation}\label{eq:chap_12_3_6}
\frac{1}{2\pi}\int_{0}^{2\pi}P(r_{0},r_{e},\phi,\theta_{e},\psi,\varphi_{e})d\phi=-1, (r_{e}>r_{0})
\end{equation}
and
\begin{equation}\label{eq:chap_12_3_7}
\lim_{r_{e}\rightarrow \infty}H(r_{e},\theta_{e},\varphi_{e})=\frac{1}{2\pi}\int_{0}^{2\pi}F(\phi,\psi)d\phi.
\end{equation}

\section{Schwarz Integral Formula}\label{sec:chap_12_4_119}

Let $f$ be an analytic function of $s$ throughout the half space $Im$ $s\geq0$ such that, for some positive constants $a$ and $M$, $f$ satisfies the order property
\begin{equation}\label{eq:chap_12_4_1}
|s^{a}f(s)|<M \texttt{ $(Im$ }s\geq0).
\end{equation}
For a fixed point $s$ above the $xz$ coordinates plane which is supposed horizontal, let $S_{R}$ denote the upper half of a positively oriented sphere of radius $R$ centered at the origin, where $R>|s|$, and $C_{R}$ denote the upper half of a positively oriented disk of radius $R$ centered at the origin and in the plane perpendicular to the $xz$ coordinates plane where the included angle between this plane and the $yz$ coordinates plane is $\varphi$. Then, according to the Cauchy spatial circle integral formula and with defining points $\tau=e_{xy}\xi+\zeta=r(e_{xy}\cos\varphi+\sin\varphi)$ on $C_{R}$ and in the $xz$ coordinates plane,
\begin{equation}\label{eq:chap_12_4_2}
f(s)=\frac{1}{2\pi i}\int_{C_{R}}\frac{f(\varsigma)}{\varsigma-s}d\varsigma+\frac{1}{2\pi i}\int_{-R}^{R}\frac{f(\tau)}{\tau-s}d\tau.
\end{equation}

We find that the first of these integrals approaches 0 as $R$ tends to $\infty$ since, in view of condition~(\ref{eq:chap_12_4_1}),
\[
|\int_{C_{R}}\frac{f(\varsigma)}{\varsigma-s}d\varsigma|<\frac{M}{R^{a}(R-|s|)}\pi R=\frac{\pi M}{R^{a}(1-|s|/R)}.
\]
Thus
\begin{equation}\label{eq:chap_12_4_3}
f(s)=\lim_{R\rightarrow\infty}\frac{1}{2\pi i}\int_{-R}^{R}\frac{f(\tau)}{\tau-s}d\tau
=\frac{1}{2\pi i}\int_{-\infty}^{\infty}\frac{f(\tau)}{\tau-s}d\tau \texttt{ $(Im$ }s>0).
\end{equation}

Condition~(\ref{eq:chap_12_4_1}) also ensures that the improper integral here converges. The number to which it converges is the same as its Cauchy principal value, and representation~(\ref{eq:chap_12_4_3}) is a Cauchy contour integral formula for the half spatial disk $Im$ $s>0$.

When the point $s$ lies below the $xz$ coordinates plane, the right-hand side of equation~(\ref{eq:chap_12_4_2}) is zero; hence integral~(\ref{eq:chap_12_4_3}) is zero for such a point. Thus, when $s$ is above the $xz$ coordinates plane, we have the following formula, where $c$ is an arbitrary constant:
\begin{equation}\label{eq:chap_12_4_4}
f(s)=\frac{1}{2\pi i}\int_{-\infty}^{\infty}(\frac{1}{\tau-s}+\frac{c}{\tau-\bar{s}})f(\tau)d\tau \texttt{ $(Im$ }s>0).
\end{equation}
In the two cases $c=-1$ and $c=1$, using $s=e_{xy}(x+iy)+z$ and $\bar{s}=e_{xy}(x-iy)+z$, this reduces, respectively, to
\begin{equation}\label{eq:chap_12_4_5}
\begin{split}
f(s)&=\frac{1}{2\pi i}\int_{-\infty}^{\infty}\frac{s-\bar{s}}{(\tau-s)(\tau-\bar{s})}f(\tau)d\tau \\
&=e_{xy}\frac{1}{\pi}\int_{-\infty}^{\infty}\frac{y}{[(\xi-x)+(\zeta-z)]^{2}+y^{2}}f(\tau)d\tau\texttt{ }(y>0)
\end{split}
\end{equation}
and
\begin{equation}\label{eq:chap_12_4_6}
f(s)=-ie_{xy}\frac{1}{\pi}\int_{-\infty}^{\infty}\frac{(\xi-x)+(\zeta-z)}{[(\xi-x)+(\zeta-z)]^{2}+y^{2}}f(\tau)d\tau\texttt{ }(y>0).
\end{equation}

If $f(s)=e_{xy}[u(x,y,z)+iv(x,y,z)]+w(z)$, it follows from formulas~(\ref{eq:chap_12_4_5}) and~(\ref{eq:chap_12_4_6}) that the complex harmonic functions $u=Re$ $f(s)$ and $v=Im$ $f(s)$, respectively, are represented in the half space $y>0$ in terms of the boundary values of $f(\tau)=e_{xy}u+w$ by the formulas
\begin{equation}\label{eq:chap_12_4_7}
u(x,y,z)=\frac{1}{\pi}\int_{-\infty}^{\infty}\frac{y[u(\xi,0,\zeta)+w(\zeta)]}
{[(\xi-x)+(\zeta-z)]^{2}+y^{2}}d\tau \texttt{ }(y>0)
\end{equation}
and
\begin{equation}\label{eq:chap_12_4_8}
v(x,y,z)=-\frac{1}{\pi}\int_{-\infty}^{\infty}\frac{[(\xi-x)+(\zeta-z)][u(\xi,0,\zeta)+w(\zeta)]}{[(\xi-x)+(\zeta-z)]^{2}+y^{2}}d\tau \texttt{ }(y>0).
\end{equation}

Formula~(\ref{eq:chap_12_4_7}) is known as the Schwarz integral formula, or the Poisson integral formula for the half space. In the next section, we shall relax the conditions for the validity of formulas~(\ref{eq:chap_12_4_7}) and~(\ref{eq:chap_12_4_8}).

\section{Dirichlet Problem for a Half Space}\label{sec:chap_12_5_120}

Let $F$ denote a $xz$ plane function of $e_{xy}x+z$ that is bounded for all $e_{xy}x+z$ and continuous except for at most a finite number of finite jumps. When $y\geq\varepsilon$ and $|e_{xy}x+z|\leq1/\varepsilon$, where $\varepsilon$ is any positive constant, the integral
\[
I(e_{xy}x+z,y)=\int_{-\infty}^{\infty}\frac{F(e_{xy}\xi+\zeta)}{[(\xi-x)+(\zeta-z)]^{2}+y^{2}}d(e_{xy}\xi+\zeta)
\]
converges uniformly with respect to $e_{xy}x+z$ and $y$, as do the integrals of the partial derivatives
of the integrand with respect to $e_{xy}x+z$ and $y$. Each of these integrals is the sum of a finite number of improper or definite integrals over regions where $F$ is continuous; hence the integrand of each component integral is a continuous function of $\tau=e_{xy}\xi+\zeta$, $e_{xy}x+z$, and $y$ when $y\geq\varepsilon$. Consequently, each partial derivative of $I(e_{xy}x+z,y)$ is represented by the integral of the corresponding derivative of the integrand whenever $y>\varepsilon$.

We write $U(x,y,z)=yI(e_{xy}x+z,y)/\pi$. Thus $U$ is the Schwarz integral transform of $F$, suggested by the second of expressions~(\ref{eq:chap_12_4_7}), Sec.~(\ref{sec:chap_12_4_119}):
\begin{equation}\label{eq:chap_12_5_1}
U(x,y,z)=\frac{1}{\pi}\int_{-\infty}^{\infty}\frac{yF(e_{xy}\xi+\zeta)}{[(\xi-x)+(\zeta-z)]^{2}+y^{2}}d(e_{xy}\xi+\zeta) \texttt{ }(y>0).
\end{equation}
Except for the factor $1/\pi$, the kernel here is $y/[(\tau-s)(\tau-\bar{s})]$. It is the imaginary component of the function $1/(\tau-s)-1/(\tau-\bar{s})$, which is analytic in $s$ when $y>0$. It follows that the kernel is complex harmonic, and so it satisfies the complex Laplace's equation in $x$, $y$, and $z$. Because the order of differentiation and integration may be interchanged, the function~(\ref{eq:chap_12_5_1}) then satisfies that equation. Consequently, $U$ is complex harmonic when $y>0$.

To prove that
\begin{equation}\label{eq:chap_12_5_2}
\lim_{y\rightarrow0,y>0}U(x,y,z)=F(e_{xy}x+z)
\end{equation}
for each fixed $e_{xy}x+z$ at which $F$ is continuous, we substitute $\tau=e_{xy}(x+y\tan\vartheta)+z$ in formula~(\ref{eq:chap_12_5_1}) and write
\begin{equation}\label{eq:chap_12_5_3}
U(x,y,z)=\frac{1}{\pi}\int_{-\pi/2}^{\pi/2}F[e_{xy}(x+y\tan\vartheta)+z]d\vartheta,(y>0).
\end{equation}
Then, if
\[
G(e_{xy}x+z,y,\vartheta)=F[e_{xy}(x+y\tan\vartheta)+z]-F(e_{xy}x+z)
\]
and $\alpha$ is some small positive constant,
\begin{equation}\label{eq:chap_12_5_4}
\begin{split}
\pi[U(x,y,z)-F(e_{xy}x+z)]&=\int_{-\pi/2}^{\pi/2}G(e_{xy}x+z,y,\vartheta)d\vartheta \\
&=I_{1}(y)+I_{2}(y)+I_{3}(y),
\end{split}
\end{equation}
where
\[
\begin{split}
I_{1}(y)&=\int_{-\pi/2}^{(-\pi/2)+\alpha}G(e_{xy}x+z,y,\vartheta)d\vartheta,\\
I_{2}(y)&=\int_{(-\pi/2)+\alpha}^{(\pi/2)-\alpha}G(e_{xy}x+z,y,\vartheta)d\vartheta,\\
I_{3}(y)&=\int_{(\pi/2)-\alpha}^{\pi/2}G(e_{xy}x+z,y,\vartheta)d\vartheta.
\end{split}
\]

If $M$ denote an upper bound for $|F(e_{xy}x+z)|$, then $|G(e_{xy}x+z,y,\vartheta)|\leq2M$. For a given positive numbers $\varepsilon$, we select $\alpha$ so that $6M\alpha<\varepsilon$; and this means that
\[
|I_{1}(y)|\leq2M\alpha<\frac{\varepsilon}{3} \texttt{ and } |I_{3}(y)|\leq2M\alpha<\frac{\varepsilon}{3}.
\]
We next show that, corresponding to $\varepsilon$, there is a positive number $\delta$ such that
\[
|I_{2}(y)|<\frac{\varepsilon}{3} \texttt{ whenever } 0<y<\delta.
\]
To do this, we observe that, since $F$ is continuous at $e_{xy}x+z$, there is a positive number $\gamma$ such that
\[
|G(e_{xy}x+z,y,\vartheta)|<\frac{\varepsilon}{3\pi} \texttt{ whenever } 0<y|\tan\vartheta|<\gamma.
\]
Now the maximum value of $|\tan\vartheta|$ as $\vartheta$ ranges from $(-\pi/2)+\alpha$ to $(\pi/2)-\alpha$ is $\tan(\pi/2-\alpha)=\cot\alpha$. Hence, if we write $\delta=\gamma\tan\alpha$, it follows that
\[
|I_{2}(y)|<\frac{\varepsilon}{3\pi}(\pi-2\alpha)<\frac{\varepsilon}{3} \texttt{ whenever } 0<y<\delta.
\]

We have thus shown that
\[
|I_{1}(y)|+|I_{2}(y)|+|I_{3}(y)|<\varepsilon \texttt{ whenever } 0<y<\delta.
\]
Condition~(\ref{eq:chap_12_5_2}) now follows from this result and equation~(\ref{eq:chap_12_5_4}).

Formula~(\ref{eq:chap_12_5_1}) therefore solves the Dirichlet problem for the half space $y>0$, with the boundary condition~(\ref{eq:chap_12_5_2}). It is evident from the form~(\ref{eq:chap_12_5_3}) of expression~(\ref{eq:chap_12_5_1}) that $|U(x,y,z)|\leq M$ in the half space, where $M$ is an upper bound of $|F(e_{xy}x+z)|$; that is, $U$ is bounded. We note that $U(x,y,z)=F_{0}$ when $F(e_{xy}x+z)=F_{0}$, where $F_{0}$ is a constant.

According to the second expression of formula~(\ref{eq:chap_12_4_8}) of Sec.~(\ref{sec:chap_12_4_119}), under certain conditions of $F$ the function
\begin{equation}\label{eq:chap_12_5_5}
V(x,y,z)=-\frac{1}{\pi}\int_{-\infty}^{\infty}
\frac{(\xi-x)+(\zeta-z)}{[(\xi-x)+(\zeta-z)]^{2}+y^{2}}F(e_{xy}\xi+\zeta)d(e_{xy}\xi+\zeta)
\end{equation}
is a complex harmonic conjugate of $U$ given by formula~(\ref{eq:chap_12_5_1}). Actually, formula~(\ref{eq:chap_12_5_5}) furnishes a complex harmonic conjugate of $U$ if $F$ is everywhere continuous, except for at most a finite number of finite jumps, and if $F$ satisfies an order property
\[
|x^{a}F(e_{xy}x+z)|<M, (a<0).
\]
For, under these conditions, we find that $U$ and $V$ satisfy the Cauchy-Riemann equations where $y>0$.

\section{Neumann Problems}\label{sec:chap_12_6_121}

As in Section~(\ref{sec:chap_12_1_116}), let $S$ denote a positively oriented sphere, centered at the origin. We let $r_{0}$ denote the radius of $S$ and write
\[
\varsigma=r_{0}(e_{xy}e^{i\phi}\cos\psi+\sin\psi)\texttt{ and }s=r(e_{xy}e^{i\theta}\cos\varphi+\sin\varphi)\texttt{ }(r<r_{0}).
\]
In view of equation~(\ref{eq:chap_12_1_6}), when $\varsigma$ is fixed, the function
\begin{equation}\label{eq:chap_12_6_1}
\begin{split}
Q(r_{0},r,\phi,\theta,\psi,\varphi)=&-r_{0}\{Re\log[r_{0}^{2}-2r_{0}r\sigma(\phi,\theta,\psi,\varphi)+r^{2}]\\
&-\log[r_{0}^{2}-2r_{0}r\sigma(\psi,\varphi)+r^{2}]\}
\end{split}
\end{equation}
is a complex harmonic function of $r, \theta$ and $\varphi$ interior to $S$ for each fixed point $\varsigma$ on the sphere $S$, because it is the real component of
$-r_{0}\log[(1-\varsigma/s)(r^{2}-r_{0}^{2}s/\varsigma)]$, where the branch cut of $\log[(1-\varsigma/s)(r^{2}-r_{0}^{2}s/\varsigma)]$ is an outward ray from the point $\varsigma$. If, moreover, $r\neq0$,
\begin{equation}\label{eq:chap_12_6_2}
\begin{split}
Q_{r}&=-\frac{r_{0}}{r}[Re\frac{2r^{2}-2r_{0}r\sigma(\phi,\theta,\psi,\varphi)}
{r_{0}^{2}-2r_{0}r\sigma(\phi,\theta,\psi,\varphi)+r^{2}}-\frac{2r^{2}-2r_{0}r\sigma(\psi,\varphi)}{r_{0}^{2}-2r_{0}r\sigma(\psi,\varphi)+r^{2}}]\\
&=\frac{r_{0}}{r}P(r_{0},r,\phi,\theta,\psi,\varphi)
\end{split}
\end{equation}
where $P$ is the Poisson kernel~(\ref{eq:chap_12_1_7}) of Section~(\ref{sec:chap_12_1_116}).

These observations suggest that the function $Q$ may be used to write an integral representation for a harmonic function $U$ whose normal derivative $U_{r}$ on the sphere $r=r_{0}$ assumes prescribed values $G(\theta,\varphi)$.

If $G$ is piecewise continuous and $U_{0}$ is an arbitrary constant, the function
\begin{equation}\label{eq:chap_12_6_3}
U(r,\theta,\varphi)=\frac{1}{2\pi}\int_{0}^{2\pi}Q(r_{0},r,\phi,\theta,\psi,\varphi)G(\phi,\psi) d\phi+U_{0}\texttt{ }(r<r_{0})
\end{equation}
is complex harmonic because the integrand is a Neumann kernel function of $r$, $\theta$, and $\varphi$. If the mean value of $G$ over the sphere $|s|=r_{0}$ is zero, or
\begin{equation}\label{eq:chap_12_6_4}
\int_{0}^{2\pi}G(\phi,\psi)d\phi=0,
\end{equation}
then, in view of equation~(\ref{eq:chap_12_6_2}),
\[
\begin{split}
U_{r}(r,\theta,\varphi)&=\frac{1}{2\pi}\int_{0}^{2\pi}\frac{r_{0}}{r}P(r_{0},r,\phi,\theta,\psi,\varphi)G(\phi,\psi)d\phi\\
&=\frac{r_{0}}{r}\cdot\frac{1}{2\pi}\int_{0}^{2\pi}P(r_{0},r,\phi,\theta,\psi,\varphi)G(\phi,\psi)d\phi.
\end{split}
\]
Now, according to equations~(\ref{eq:chap_12_2_1}) and~(\ref{eq:chap_12_2_2}) of section~(\ref{sec:chap_12_2_117}),
\[
\lim_{r\rightarrow r_{0}, r<r_{0}}\frac{1}{2\pi}\int_{0}^{2\pi}P(r_{0},r,\phi,\theta,\psi,\varphi)G(\phi,\psi) d\phi=G(\theta,\varphi).
\]
Hence
\begin{equation}\label{eq:chap_12_6_5}
\lim_{r\rightarrow r_{0}, r<r_{0}}U_{r}(r,\theta,\varphi)=G(\theta,\varphi)
\end{equation}
for each pair of values of $\theta$ and $\varphi$ at which $G$ is continuous.

When $G$ is piecewise continuous and satisfies condition~(\ref{eq:chap_12_6_4}), the formula
\begin{equation}\label{eq:chap_12_6_6}
\begin{split}
U(r,\theta,\varphi)=&-\frac{r_{0}}{2\pi}\int_{0}^{2\pi}\{Re\log[r_{0}^{2}-2r_{0}r\sigma(\phi,\theta,\psi,\varphi)+r^{2}]\\
&-\log[r_{0}^{2}-2r_{0}r\sigma(\psi,\varphi)+r^{2}]\}G(\phi,\psi)d\phi+U_{0}
\end{split}
\end{equation}
where $r<r_{0}$, therefore, solves the Neumann problem for the region interior to the sphere $r=r_{0}$, where $G(\theta,\varphi)$ is the normal derivative of the harmonic function $U(r,\theta,\varphi)$ at the boundary in the sense of condition~(\ref{eq:chap_12_6_5}). Note how it follows from equations~(\ref{eq:chap_12_6_4}) and~(\ref{eq:chap_12_6_6}) that, since $\ln r_{0}^{2}$ is constant, $U_{0}$ is the value of $U$ at the center $r=0$ of the sphere $r=r_{0}$.

%The values $U(r,\theta,\varphi)$ may represent steady temperatures in a sphere $r<r_{0}$ with insulated faces. In that case, condition~(\ref{eq:chap_12_6_5}) states that the flux of heat into the sphere through its edge is proportional to $G(\theta,\varphi)$. Condition~(\ref{eq:chap_12_6_4}) is the natural physical requirement that the total rate of flow of heat into the sphere be zero, since temperatures do not vary with time.

A corresponding formula for a complex harmonic function $H$ in the region exterior to the sphere $r=r_{0}$ can be written in terms of $Q$ as
\begin{equation}\label{eq:chap_12_6_7}
U(r_{e},\theta_{e},\varphi_{e})=
\frac{1}{2\pi}\int_{0}^{2\pi}Q(r_{0},r_{e},\phi,\theta_{e},\psi,\varphi_{e})G(\phi,\psi)d\phi+H_{0}\texttt{ }(r_{e}>r_{0})
\end{equation}
where $H_{0}$ is a constant. As before we assume that $G$ is piecewise continuous and that condition~(\ref{eq:chap_12_6_4}) holds. Then
\[
H_{0}=\lim_{R\rightarrow \infty}U(r_{e},\theta_{e},\varphi_{e})
\]
and
\begin{equation}\label{eq:chap_12_6_8}
\lim_{r_{e}\rightarrow r_{0}, r_{e}>r_{0}}H_{r_{e}}(r_{e},\theta_{e},\varphi_{e})=G(\theta_{e},\varphi_{e})
\end{equation}
for each pair of values of $\theta_{e}$ and $\varphi_{e}$ at which $G$ is continuous.

Turn now to a half space, we let $G(e_{xy}x+z)$ be continuous for all two-dimensional variables $e_{xy}x+z$ on the $xz$ plane, except possibly for a finite number of finite jumps, and let it satisfy an order property
\begin{equation}\label{eq:chap_12_6_9}
|x^{a}G(e_{xy}x+z)|<M\texttt{ }(a>1)
\end{equation}
when $-\infty<x<\infty$ and $-\infty<z<\infty$. For each fixed two-dimensional variable $\tau=e_{xy}\xi+\zeta$ on the $xz$ plane, the function $\log[(\tau-s)(\tau-\bar{s})]$ is complex harmonic in the half space $Im$ $s>0$. Consequently, the function
\begin{equation}\label{eq:chap_12_6_10}
U(x,y,z)=\frac{1}{2\pi}\int_{-\infty}^{\infty}\log[(\tau-s)(\tau-\bar{s})]G(\tau)d\tau+U_{0}\texttt{ }(y>0),
\end{equation}
where $U_{0}$ is a constant, is complex harmonic in that half space.

Formula~(\ref{eq:chap_12_6_10}) was written with the Schwarz integral transform~(\ref{eq:chap_12_5_1}), Sec.~(\ref{sec:chap_12_5_120}), in mind; for it follows from formula~(\ref{eq:chap_12_6_10}) that
\begin{equation}\label{eq:chap_12_6_11}
U_{y}(x,y,z)=\frac{1}{\pi}\int_{-\infty}^{\infty}\frac{yG(e_{xy}\xi+\zeta)}{[(\xi-x)+(\zeta-z)]^{2}+y^{2}}d(e_{xy}\xi+\zeta) \texttt{ }(y>0).
\end{equation}

In view of equations~(\ref{eq:chap_12_5_1}) and~(\ref{eq:chap_12_5_2}) of Sec.~(\ref{sec:chap_12_5_120}), then
\begin{equation}\label{eq:chap_12_6_12}
\lim_{y\rightarrow0,y>0}U_{y}(x,y,z)=G(e_{xy}x+z)
\end{equation}
at each point $e_{xy}x+z$ on the $xz$ coordinates plane where $G$ is continuous.

Integral formula~(\ref{eq:chap_12_6_10}) evidently solves the Neumann problem for the half space $y>0$, with boundary condition~(\ref{eq:chap_12_6_12}). But we have not presented conditions on $G$ that are sufficient to ensure that the complex harmonic function $U$ is bounded as $|s|$ increases.

When $G$ is an odd function, formula~(\ref{eq:chap_12_6_10}) may be written
\begin{equation}\label{eq:chap_12_6_13}
U(x,y,z)=\frac{1}{2\pi}\int_{0}^{\infty}\log\frac{[(\xi-x)+(\zeta-z)]^{2}+y^{2}}
{[(\xi+x)+(\zeta+z)]^{2}+y^{2}}G(e_{xy}\xi+\zeta)d(e_{xy}\xi+\zeta)
\end{equation}
where $x>0$, $y>0$, and $z>0$. This represents a function that is complex harmonic in the first quadrant $x>0$, $y>0$, and $z>0$, and satisfies the boundary conditions
\begin{equation}\label{eq:chap_12_6_14}
U(0,y,0)=0\texttt{ }(y>0),
\end{equation}
\begin{equation}\label{eq:chap_12_6_15}
\lim_{y\rightarrow 0,y>0}U_{y}(x,y,z)=G(e_{xy}x+z)\texttt{ }(x>0,z>0).
\end{equation}

%-----------------------------------------------------------------------
% Beginning of chap13.tex
%-----------------------------------------------------------------------
%
% AMS-LaTeX 1.2 sample file for a monograph, based on amsbook.cls.
% This is a data file input by chapter.tex.
%%%%%%%%%%%%%%%%%%%%%%%%%%%%%%%%%%%%%%%%%%%%%%%%%%%%%%%%%%%%%%%%%%%%

%\part{This is a Part Title Sample}

\chapter{Sphere Integral Formulas of the Poisson Type}\label{ch:chap_13}

In this chapter, we develop a theory that enables us to obtain solutions to a variety of boundary value problems where those solutions are expressed in terns of definite or improper integrals. Many of the integrals occurring are then readily evaluated.

\section{Poisson Sphere Integral Formula}\label{sec:chap_13_1_122}

Let $S$ denote a positively oriented sphere, centered at the origin, and suppose that a function $f$ is analytic inside and on $S$. The Cauchy's surface integral formula (Sec.~(\ref{sec:chap_5_6_56}))
\begin{equation}\label{eq:chap_13_1_1}
f(s)=\frac{1}{4\pi}\iint_{S}\frac{f(\varsigma)}{(\varsigma-s)^{2}}d\sigma.
\end{equation}
expresses the value of $f$ at any point $s$ interior to $S$ in terms of the values of $f$ at points $\varsigma$ on $S$. In this section, we shall obtain from formula~(\ref{eq:chap_13_1_1}) a corresponding formula for the real part of the function $f$; and, in next section, we shall use that result to solve the Dirichlet problem for the sphere bounded by $S$.

We let $R$ denote the radius of $S$ and write
\[
s=r(e_{xy}e^{i\theta}\cos\varphi+\sin\varphi)\texttt{ }(0<r<R,\texttt{ }0\leq\theta\leq2\pi,\texttt{ }-\pi/2\leq\varphi\leq\pi/2).
\]
The inverse of the nonzero point $s$ with respect to the sphere is the point $s_{1}$ lying on the same ray from the origin as $s$ and satisfying the condition $|s||s_{1}|=R^{2}$; thus, if $\varsigma$ is a point on $S$,
\begin{equation}\label{eq:chap_13_1_2}
s_{1}=\frac{R^{2}}{r}(e_{xy}e^{i\theta}\cos\varphi+\sin\varphi)=\frac{R^{2}}{r^{2}}s.
\end{equation}
Since $s_{1}$ is exterior to the sphere $S$, it follows from the Cauchy theorem of surface integration that the value of the integral in equation~(\ref{eq:chap_13_1_1}) is zero when $s$ is replaced by $s_{1}$ in the integrand. Hence
\[
f(s)=\frac{1}{4\pi}\iint_{S}[\frac{1}{(\varsigma-s)^{2}}-\frac{1}{(\varsigma-s_{1})^{2}}]f(\varsigma)d\sigma;
\]
and, using the parametric representation
\[
\varsigma=R(e_{xy}e^{i\phi}\cos\psi+\sin\psi)\texttt{ }(0\leq\phi\leq2\pi,\texttt{ }-\pi/2\leq\psi\leq\pi/2)
\]
and $d\sigma=\varsigma^{2}\cos\psi d\psi d\phi$ for $S$, we may write
\[
f(s)=\frac{1}{4\pi}\int_{-\pi/2}^{\pi/2}\int_{0}^{2\pi}[\frac{\varsigma^{2}}{(\varsigma-s)^{2}}-\frac{\varsigma^{2}}{(\varsigma-s_{1})^{2}}]f(\varsigma)\cos\psi d\psi d\phi,
\]
where, for conveniens, we retain the $\varsigma$ to denote $R(e_{xy}e^{i\phi}\cos\psi+\sin\psi)$.

In view of the last of expression~(\ref{eq:chap_13_1_2}) for $s_{1}$, the factor inside the parentheses here may be written
\begin{equation}\label{eq:chap_13_1_3}
\frac{\varsigma^{2}}{(\varsigma-s)^{2}}-\frac{\varsigma^{2}}{(\varsigma-\frac{R^{2}}{r^{2}}s)^{2}}
=\frac{1}{(1-s/\varsigma)^{2}}-\frac{1}{(1-\frac{R^{2}}{r^{2}}s/\varsigma)^{2}}
\end{equation}
\[
=[1-\frac{r^{2}\varsigma/s-R^{2}s/\varsigma}{(1-\varsigma/s)(r^{2}-R^{2}s/\varsigma)}]
\frac{R^{2}-r^{2}}{(1-\varsigma/s)(r^{2}-R^{2}s/\varsigma)}.
\]

An alternative form of the Cauchy integral formula~(\ref{eq:chap_13_1_1}) is, therefore,
\begin{equation}\label{eq:chap_13_1_4}
f[r(e_{xy}e^{i\theta}\cos\varphi+\sin\varphi)]=\frac{1}{4\pi}\int_{-\pi/2}^{\pi/2}\int_{0}^{2\pi}
[1-\frac{r^{2}\varsigma/s-R^{2}s/\varsigma}{(1-\varsigma/s)(r^{2}-R^{2}s/\varsigma)}]
\end{equation}
\[
\times\frac{R^{2}-r^{2}}{(1-\varsigma/s)(r^{2}-R^{2}s/\varsigma)}
f[R(e_{xy}e^{i\phi}\cos\psi+\sin\psi)]\cos\psi d\psi d\phi\texttt{ }(0<r<R).
\]
This form is also valid when $r=0$; in that case, it reduces directly to
\[
f(0)=\frac{1}{4\pi}
\int_{-\pi/2}^{\pi/2}\int_{0}^{2\pi}f[R(e_{xy}e^{i\phi}\cos\psi+\sin\psi)]\cos\psi d\psi d\phi
\]
which is just the parametric form of equation~(\ref{eq:chap_13_1_1}) with $s=0$.

From formula~(\ref{eq:chap_1_7_2}) in Sec.~(\ref{sec:chap_1_7}), the quotients $(r\varsigma)/(Rs)$ and $(Rs)/(r\varsigma)$ are
\[
\frac{r\varsigma}{Rs}=\{
    \begin{array}{cc}
        e_{xy}(\frac{e^{i\phi}\cos\psi+\sin\psi}{e^{i\theta}\cos\varphi+\sin\varphi}
-\frac{\sin\psi}{\sin\varphi})+\frac{\sin\psi}{\sin\varphi} & \texttt{ when }\sin\varphi\neq0 \\
        e_{xy}e^{i(\phi-\theta)} & \texttt{ when }\sin\varphi=0
    \end{array}
\]
and
\[
\frac{Rs}{r\varsigma}=\{
    \begin{array}{cc}
        e_{xy}(\frac{e^{i\theta}\cos\varphi+\sin\varphi}{e^{i\phi}\cos\psi+\sin\psi}
-\frac{\sin\varphi}{\sin\psi})+\frac{\sin\varphi}{\sin\psi} & \texttt{ when }\sin\psi\neq0 \\
        e_{xy}e^{i(\theta-\phi)} & \texttt{ when }\sin\psi=0
    \end{array},
\]
respectively. Then we can write
\[
(1-\varsigma/s)(r^{2}-R^{2}s/\varsigma)=R^{2}-Rr(\frac{r\varsigma}{Rs}+\frac{Rs}{r\varsigma})+r^{2}
=e_{xy}(a-b)+(b+c)
\]
where
\[
\begin{split}
a&=-2Rr\cos(\phi-\theta),\texttt{ }
b=0\texttt{ when }\sin\psi=0,\texttt{ }\sin\varphi=0,\\
a&=-Rr(e^{i(\phi-\theta)}+\frac{e^{i\theta}}{e^{i\phi}\cos\psi+\sin\psi}),\texttt{ }
b=0\texttt{ when }\sin\psi\neq0,\texttt{ }\sin\varphi=0,\\
a&=-Rr(\frac{e^{i\phi}}{e^{i\theta}\cos\varphi+\sin\varphi}+e^{i(\theta-\phi)}),\texttt{ }
b=0\texttt{ when }\sin\psi=0,\texttt{ }\sin\varphi\neq0,\\
a&=-Rr(\frac{e^{i\phi}\cos\psi+\sin\psi}{e^{i\theta}\cos\varphi+\sin\varphi}+\frac{e^{i\theta}\cos\varphi+\sin\varphi}{e^{i\phi}\cos\psi+\sin\psi}),\texttt{ }
b=-Rr(\frac{\sin\psi}{\sin\varphi}+\frac{\sin\varphi}{\sin\psi})
\end{split}
\]
\[
\texttt{ when }\sin\psi\sin\varphi\neq0,\texttt{ and }c=R^{2}+r^{2},
\]
and from formula~(\ref{eq:chap_1_7_3}) in Sec.~(\ref{sec:chap_1_7})
\[
\frac{1}{(1-\varsigma/s)(r^{2}-R^{2}s/\varsigma)}=\frac{1}{e_{xy}(a-b)+(b+c)}=e_{xy}(\frac{1}{a+c}-\frac{1}{b+c})+\frac{1}{b+c}
\]
\[
=e_{xy}[\frac{1}{R^{2}-2Rr\sigma(\phi,\theta,\psi,\varphi)+r^{2}}
-\frac{1}{R^{2}-2Rr\sigma(\psi,\varphi)+r^{2}}]
\]
\[
+\frac{1}{R^{2}-2Rr\sigma(\psi,\varphi)+r^{2}}
\]
or
\begin{equation}\label{eq:chap_13_1_5}
\frac{1}{(1-\varsigma/s)(r^{2}-R^{2}s/\varsigma)}=\frac{1}{R^{2}-2Rr\sigma(\psi,\varphi)+r^{2}}
\end{equation}
\[
+e_{xy}[\frac{1}{R^{2}-2Rr\sigma(\phi,\theta,\psi,\varphi)+r^{2}}
-\frac{1}{R^{2}-2Rr\sigma(\psi,\varphi)+r^{2}}]
\]
where when $\sin\psi\sin\varphi=0$, there are $\sigma(\psi,\varphi)=0$, and
\[
\begin{split}
\sigma(\phi,\theta,\psi,\varphi)&=\cos(\phi-\theta)\texttt{ when }\sin\psi=0,\texttt{ }\sin\varphi=0,\\
\sigma(\phi,\theta,\psi,\varphi)&=\frac{1}{2}(e^{i(\phi-\theta)}+\frac{e^{i\theta}}{e^{i\phi}\cos\psi+\sin\psi})
\texttt{ when }\sin\psi\neq0,\texttt{ }\sin\varphi=0,\\
\sigma(\phi,\theta,\psi,\varphi)&=\frac{1}{2}(\frac{e^{i\phi}}{e^{i\theta}\cos\varphi+\sin\varphi}+e^{i(\theta-\phi)})
\texttt{ when }\sin\psi=0,\texttt{ }\sin\varphi\neq0,
\end{split}
\]
and when $\sin\psi\sin\varphi\neq0$
\[
\sigma(\psi,\varphi)=\frac{1}{2}(\frac{\sin\psi}{\sin\varphi}+\frac{\sin\varphi}{\sin\psi})
\]
and
\[
\sigma(\phi,\theta,\psi,\varphi)=\frac{1}{2}(\frac{e^{i\phi}\cos\psi+\sin\psi}{e^{i\theta}\cos\varphi+\sin\varphi}+\frac{e^{i\theta}\cos\varphi+\sin\varphi}{e^{i\phi}\cos\psi+\sin\psi})
\]
\[
=\frac{\cos(\phi-\theta)+\cos\phi\tan\varphi+\cos\theta\tan\psi+\tan\psi\tan\varphi}{(1+\cos\phi\sin2\psi)(1+\cos\theta\sin2\varphi)}\cos\psi\cos\varphi
\]
\[
+\frac{e^{i(\phi-\theta)}+e^{i\phi}\tan\varphi+e^{-i\theta}\tan\psi
+\tan\psi\tan\varphi}{(1+\cos\phi\sin2\psi)(1+\cos\theta\sin2\varphi)}\cos\phi\sin\psi\cos^{2}\psi\cos\varphi
\]
\[
+\frac{e^{-i(\phi-\theta)}+e^{-i\phi}\tan\varphi+e^{i\theta}\tan\psi
+\tan\psi\tan\varphi}{(1+\cos\phi\sin2\psi)(1+\cos\theta\sin2\varphi)}\cos\theta\cos\psi\sin\varphi\cos^{2}\varphi.
\]

Also we can write
\[
\frac{r\varsigma}{Rs}-\frac{Rs}{r\varsigma}=e_{xy}(\frac{e^{i\phi}\cos\psi+\sin\psi}{e^{i\theta}\cos\varphi+\sin\varphi}
-\frac{e^{i\theta}\cos\varphi+\sin\varphi}{e^{i\phi}\cos\psi+\sin\psi})
\]
\[
-e_{xy}(\frac{\sin\psi}{\sin\varphi}-\frac{\sin\varphi}{\sin\psi})
+(\frac{\sin\psi}{\sin\varphi}-\frac{\sin\varphi}{\sin\psi})
\]
and
\[
r^{2}\varsigma/s-R^{2}s/\varsigma=2Rr[e_{xy}(\sigma_{m}(\phi,\theta,\psi,\varphi)-\sigma_{m}(\psi,\varphi))+\sigma_{m}(\psi,\varphi)]
\]
where when $\sin\psi\sin\varphi=0$, there are $\sigma_{m}(\psi,\varphi)=0$, and
\[
\begin{split}
\sigma_{m}(\phi,\theta,\psi,\varphi)&=\sin(\phi-\theta)\texttt{ when }\sin\psi=0,\texttt{ }\sin\varphi=0,\\
\sigma_{m}(\phi,\theta,\psi,\varphi)&=\frac{1}{2}(e^{i(\phi-\theta)}-\frac{e^{i\theta}}{e^{i\phi}\cos\psi+\sin\psi})
\texttt{ when }\sin\psi\neq0,\texttt{ }\sin\varphi=0,\\
\sigma_{m}(\phi,\theta,\psi,\varphi)&=\frac{1}{2}(\frac{e^{i\phi}}{e^{i\theta}\cos\varphi+\sin\varphi}-e^{i(\theta-\phi)})
\texttt{ when }\sin\psi=0,\texttt{ }\sin\varphi\neq0,
\end{split}
\]
and when $\sin\psi\sin\varphi\neq0$
\[
\sigma_{m}(\psi,\varphi)=\frac{1}{2}(\frac{\sin\psi}{\sin\varphi}-\frac{\sin\varphi}{\sin\psi})
\]
and
\[
\sigma_{m}(\phi,\theta,\psi,\varphi)=\frac{1}{2}(\frac{e^{i\phi}\cos\psi+\sin\psi}{e^{i\theta}\cos\varphi+\sin\varphi}
-\frac{e^{i\theta}\cos\varphi+\sin\varphi}{e^{i\phi}\cos\psi+\sin\psi})
\]

Hence, if $u$ is the real part of the analytic function $f$, it follows from formula~(\ref{eq:chap_13_1_4}) that
\begin{equation}\label{eq:chap_13_1_6}
u(r,\theta,\varphi)=\frac{1}{4\pi}\int_{-\pi/2}^{\pi/2}\int_{0}^{2\pi}
Re\{[1-\frac{r^{2}\varsigma/s-R^{2}s/\varsigma}{(1-\varsigma/s)(r^{2}-R^{2}s/\varsigma)}]
\end{equation}
\[
\times\frac{R^{2}-r^{2}}{(1-\varsigma/s)(r^{2}-R^{2}s/\varsigma)}\}
u(R,\phi,\psi)\cos\psi d\psi d\phi\texttt{ }(r<R)
\]
\[
=\frac{1}{4\pi}\int_{-\pi/2}^{\pi/2}\int_{0}^{2\pi}\{Re[(1-\frac{2Rr\sigma_{m}(\phi,\theta,\psi,\varphi)}{R^{2}-2Rr\sigma(\phi,\theta,\psi,\varphi)+r^{2}})\frac{R^{2}-r^{2}}{R^{2}-2Rr\sigma(\phi,\theta,\psi,\varphi)+r^{2}}]
\]
\[
-(1-\frac{2Rr\sigma_{m}(\psi,\varphi)}{R^{2}-2Rr\sigma(\psi,\varphi)+r^{2}})\frac{R^{2}-r^{2}}{R^{2}-2Rr\sigma(\psi,\varphi)+r^{2}}\}
u(R,\phi,\psi)\cos\psi d\psi d\phi\texttt{ }(r<R).
\]
This is the Poisson sphere integral formula for the complex harmonic function $u$ in the open sphere bounded by the spherical surface $|s|=r=R$.

Formula~(\ref{eq:chap_13_1_6}) defines a linear integral transformation of $u(R,\phi,\psi)$ into $u(r,\theta,\varphi)$. The kernel of the transformation is, except for the factor $1/(4\pi)$, the real valued function
\begin{equation}\label{eq:chap_13_1_7}
P(R,r,\phi,\theta,\psi,\varphi)=[\frac{2Rr\sigma_{m}(\psi,\varphi)}{R^{2}-2Rr\sigma(\psi,\varphi)+r^{2}}-1]\frac{R^{2}-r^{2}}{R^{2}-2Rr\sigma(\psi,\varphi)+r^{2}}
\end{equation}
\[
+Re\{[1-\frac{2Rr\sigma_{m}(\phi,\theta,\psi,\varphi)}{R^{2}-2Rr\sigma(\phi,\theta,\psi,\varphi)+r^{2}}]\frac{R^{2}-r^{2}}{R^{2}-2Rr\sigma(\phi,\theta,\psi,\varphi)+r^{2}}\},
\]
which is called the Poisson kernel for the sphere integral formula. In view of equations~(\ref{eq:chap_13_1_5}), we can also write
\begin{equation}\label{eq:chap_13_1_8}
P(R,r,\phi,\theta,\psi,\varphi)=[\frac{2Rr\sigma_{m}(\psi,\varphi)}{R^{2}-2Rr\sigma(\psi,\varphi)+r^{2}}-1]\frac{R^{2}-r^{2}}{R^{2}-2Rr\sigma(\psi,\varphi)+r^{2}}
\end{equation}
\[
+(R^{2}-r^{2})[\frac{R^{2}-2Rr\sigma_{r}(\phi,\theta,\psi,\varphi)+r^{2}}{|R^{2}-2Rr\sigma(\phi,\theta,\psi,\varphi)+r^{2}|^{2}}
-Re\frac{2Rr\sigma_{m}(\phi,\theta,\psi,\varphi)}{(R^{2}-2Rr\sigma(\phi,\theta,\psi,\varphi)+r^{2})^{2}}],
\]
where
\[
\sigma_{r}(\phi,\theta,\psi,\varphi)=\frac{\cos(\phi-\theta)+\cos\phi\tan\varphi+\cos\theta\tan\psi+\tan\psi\tan\varphi}{(1+\cos\phi\sin2\psi)(1+\cos\theta\sin2\varphi)}\cos\psi\cos\varphi
\]
\[
+\frac{\cos(\phi-\theta)+\cos\theta\tan\psi
+\cos\phi\tan\varphi+2\tan\psi\tan\varphi}{2(1+\cos\theta\sin2\varphi)(1+\cos\phi\sin2\psi)}\cos\phi\sin\psi\cos^{2}\psi\cos\varphi
\]
\[
+\frac{\cos(\phi-\theta)+\cos\phi\tan\varphi+\cos\theta\tan\psi
+2\tan\psi\tan\varphi}{2(1+\cos\phi\sin2\psi)(1+\cos\theta\sin2\varphi)}\cos\theta\cos\psi\sin\varphi\cos^{2}\varphi.
\]

Since $r<R$, it is clear that $P$ is a positive function. Moreover, we find from the second
of equations~(\ref{eq:chap_13_1_3}) that
\begin{equation}\label{eq:chap_13_1_9}
P(R,r,\phi,\theta,\psi,\varphi)=Re\{[1-\frac{r^{2}\varsigma/s-R^{2}s/\varsigma}{(1-\varsigma/s)(r^{2}-R^{2}s/\varsigma)}]
\frac{R^{2}-r^{2}}{(1-\varsigma/s)(r^{2}-R^{2}s/\varsigma)}\},
\end{equation}
and from equation~(\ref{eq:chap_13_1_7}) that because $P$ consists of trigonometric functions, $P$ is analytic entirely in $S$. Thus $P(R,r,\phi,\theta,\psi,\varphi)$ is a harmonic function of $r$, $\theta$, and $\varphi$ interior to $S$ for each fixed $\varsigma$ on $S$. From equation~(\ref{eq:chap_13_1_7}), we see that $P(R,r,\phi,\theta,\psi,\varphi)$ is an even periodic function of $\phi$ and $\theta$, with period $2\pi$; and its value is $1$ when $r=0$.

The Poisson integral formula~(\ref{eq:chap_13_1_6}) can now be written
\begin{equation}\label{eq:chap_13_1_10}
u(r,\theta,\varphi)=\frac{1}{4\pi}\int_{-\pi/2}^{\pi/2}\int_{0}^{2\pi}P(R,r,\phi,\theta,\psi,\varphi)u(R,\phi,\psi)\cos\psi d\psi d\phi
\end{equation}
where $r<R$. When $f(s)=u(r,\theta,\varphi)=1$, equation~(\ref{eq:chap_13_1_10}) shows that $P$ has the property
\begin{equation}\label{eq:chap_13_1_11}
\frac{1}{4\pi}\int_{-\pi/2}^{\pi/2}\int_{0}^{2\pi}P(R,r,\phi,\theta,\psi,\varphi)u(R,\phi,\psi)\cos\psi d\psi d\phi=1\texttt{ }(r<R).
\end{equation}

We have assumed that $f$ is analytic not only interior to $S$ but also on $S$ itself and that $u$ is, therefore, harmonic in a domain which includes all points on that sphere. In particular, $u$ is continuous on $S$. The conditions will now be relaxed.

\section{Dirichlet problem for a sphere}\label{sec:chap_13_2_123}

Let $F$ be a piecewise continuous function of $\theta$ on the interval $0\leq\theta\leq 2\pi$ and $\varphi$ on the interval $-\pi/2\leq\varphi\leq\pi/2$. The Poisson integral transform of $F$ is defined in terms of the Poisson kernel $P(R,r,\phi,\theta,\psi,\varphi)$, introduced in Sec.~(\ref{sec:chap_13_1_122}), by means of the equation
\begin{equation}\label{eq:chap_13_2_1}
U(r,\theta,\varphi)=\frac{1}{4\pi}
\int_{-\pi/2}^{\pi/2}\int_{0}^{2\pi}P(R,r,\phi,\theta,\psi,\varphi)F(\phi,\psi)\cos\psi d\psi d\phi\texttt{ }(r<R).
\end{equation}

In this section, we shall prove that the function $U(r,\theta,\varphi)$ is complex harmonic inside the sphere $r=R$ and
\begin{equation}\label{eq:chap_13_2_2}
\lim_{r\rightarrow R,r<R}U(r,\theta,\varphi)=F(\theta,\varphi)
\end{equation}
for each fixed pair of $\theta$ and $\varphi$ at which $F$ is continuous. Thus $U$ is a solution of the Dirichlet problem for the sphere $r<R$ in the sense that $U(r,\theta,\varphi)$ approaches the the boundary value $F(\theta,\varphi)$ as the point $(r,\theta,\varphi)$ approaches $(R,\theta,\varphi)$ along a radius, except at the finite number of points $(R,\theta,\varphi)$ where discontinuities of $F$ may occur.

We turn now to the proof that the function $U$ defined in equation~(\ref{eq:chap_13_2_2}) satisfies the Dirichlet problem for the sphere $r<R$. First of all, $U$ is complex harmonic inside the sphere $r=R$ because $P$ is a complex harmonic function of $r$, $\theta$, and $\varphi$ there. More precisely, since $F$ is piecewise continuous, integral~(\ref{eq:chap_13_2_1}) may be written as the sum of a finite number of definite integrals each of which has an integrand that is continuous in $r$, $\theta$, $\varphi$, $\phi$, and $\psi$. The partial derivatives of those integrands with respect to $r$, $\theta$, and $\varphi$ are also continuous. Since the order of integration and differentiation with respect to $r$, $\theta$, and $\varphi$ can, then, be interchanged and since $P$ satisfies complex Laplace's equation (see Theorem~(\ref{th:chap_2_15_1s}) in Sec.~(\ref{sec:chap_2_15_25}))
\[
r\frac{\partial}{\partial r}(rP_{r}\cos^{2}\varphi)+\frac{P_{\theta\theta}}{\cos^{2}\varphi}
+\tan\varphi\frac{\partial}{\partial\varphi}(P_{\varphi}\sin\varphi\cos\varphi)
-2\sin\varphi\frac{\partial}{\partial\varphi}(rP_{r}\cos\varphi)=0.
\]
in the spherical coordinates $r$, $\theta$, and $\varphi$, it follows that $U$ satisfies that equation too.

In order to verify limit~(\ref{eq:chap_13_2_2}), we need to show that if $F$ is continuous at $\theta$ and $\varphi$, there corresponds to each positive number $\varepsilon$ a positive number $\delta$ such that
\begin{equation}\label{eq:chap_13_2_6}
|U(r,\theta,\varphi)-F(\theta,\varphi)|<\varepsilon \texttt{ whenever } 0<R-r<\delta.
\end{equation}
We start by referring to property~(\ref{eq:chap_13_1_11}), Sec.~(\ref{sec:chap_13_1_122}), of the Poisson kernel and writing
\[
U(r,\theta,\varphi)-F(\theta,\varphi)=\frac{1}{4\pi}
\int_{-\pi/2}^{\pi/2}\int_{0}^{2\pi}P(R,r,\phi,\theta,\psi,\varphi)[F(\phi,\psi)-F(\theta,\varphi)]\cos\psi d\psi d\phi.
\]
For convenience, we let $F$ be extended periodically, with period $2\pi$, so that the integrand here is periodic in $\phi$ with that same period. Also, we may assume that $0<r<R$ because of the nature of the limit to be established.

Next, we observe that, since $F$ is continuous at $\theta$, there is a small positive number $\alpha$ such that
\begin{equation}\label{eq:chap_13_2_7}
|F(\phi,\psi)-F(\theta,\varphi)|<\frac{\varepsilon}{2} \texttt{ whenever } |\phi-\theta|\leq\alpha.
\end{equation}
Evidently,
\begin{equation}\label{eq:chap_13_2_8}
U(r,\theta,\varphi)-F(\theta,\varphi)=I_{1}(r)+I_{2}(r),
\end{equation}
where
\[
\begin{split}
I_{1}(r)&=\frac{1}{4\pi}
\int_{-\pi/2}^{\pi/2}\int_{\theta-\alpha}^{\theta+\alpha}P(R,r,\phi,\theta,\psi,\varphi)[F(\phi,\psi)-F(\theta,\varphi)]\cos\psi d\psi d\phi,\\
I_{2}(r)&=\frac{1}{4\pi}
\int_{-\pi/2}^{\pi/2}\int_{\theta+\alpha}^{\theta-\alpha+2\pi}P(R,r,\phi,\theta,\psi,\varphi)[F(\phi,\psi)-F(\theta,\varphi)]\cos\psi d\psi d\phi.
\end{split}
\]

The fact that $P$ is a positive function (See last section), together with the first of inequalities~(\ref{eq:chap_13_2_7}) just above and property~(\ref{eq:chap_13_1_11}) of that function, enable us to write
\[
\begin{split}
|I_{1}(r)|&\leq\frac{1}{4\pi}
\int_{-\pi/2}^{\pi/2}\int_{\theta-\alpha}^{\theta+\alpha}P(R,r,\phi,\theta,\psi,\varphi)|F(\phi,\psi)-F(\theta,\varphi)|\cos\psi d\psi d\phi,\\
&<\frac{\varepsilon}{8\pi}
\int_{-\pi/2}^{\pi/2}\int_{0}^{2\pi}P(R,r,\phi,\theta,\psi,\varphi)\cos\psi d\psi d\phi=\frac{\varepsilon}{2}
\end{split}
\]
As for the integral $I_{2}(r)$, one may know that the denominator $R^{2}-2Rr\sigma(\phi,\theta,\psi,\varphi)+r^{2}$ in expression~(\ref{eq:chap_13_2_8}) for $P(R,r,\phi,\theta,\psi,\varphi)$ in Sec.~(\ref{sec:chap_13_1_122}) has a (positive) minimum value $m$ as the argument $\phi$ of $\varsigma$ varies over the closed interval
\[
\theta+\alpha\leq\phi\leq\theta-\alpha+2\pi.
\]
So, if $M$ denotes an upper bound of the piecewise continuous function $|F(\phi,\psi)-F(\theta,\varphi)|$ on the interval $0\leq\phi\leq2\pi$, it follows that
\[
|I_{2}(r)|\leq\frac{(R^{2}-r^{2})^{2}M}{2\pi m}2\pi<\frac{4M^{2}R}{m}(R-r)^{2}<\frac{4MR^{2}}{m}\delta=\frac{\varepsilon}{2}
\]
whenever $(R-r)^{2}<\delta$, where
\begin{equation}\label{eq:chap_13_2_9}
\delta=\frac{m\varepsilon}{8MR^{2}}.
\end{equation}

Finally, the results in the two preceding paragraphes tell us that
\[
|U(r,\theta,\varphi)-F(\theta,\varphi)|\leq|I_{1}(r)|+|I_{2}(r)|<\frac{\varepsilon}{2}+\frac{\varepsilon}{2}=\varepsilon
\]
whenever $(R-r)^{2}<\delta$, where $\delta$ is the positive number defined by equation~(\ref{eq:chap_13_2_9}). That is, statement~(\ref{eq:chap_13_2_6}) holds when that choice of $\delta$ is made.

According to expression~(\ref{eq:chap_13_2_1}), the value of $U$ at $r=0$ is
\[
\frac{1}{4\pi}
\int_{-\pi/2}^{\pi/2}\int_{0}^{2\pi}F(\phi,\psi)\cos\psi d\psi d\phi.
\]
Thus the value of a complex harmonic function at the center of the sphere $r=R$ is the average of the boundary values on the sphere.

\section{Related Boundary Value Problems}\label{sec:chap_13_3_124}

The function $F$ representing boundary values on the sphere $r=R$ is assumed to be piecewise continuous.

Suppose that $F(2\pi-\theta,\varphi)=-F(\theta,\varphi)$. The Poisson integral formula~(\ref{eq:chap_13_2_1}) of Sec.~(\ref{sec:chap_13_2_123}) then becomes
\begin{equation}\label{eq:chap_13_3_1}
U(r,\theta,\varphi)=
\end{equation}
\[
\frac{1}{4\pi}
\int_{-\pi/2}^{\pi/2}\int_{0}^{2\pi}[P(R,r,\phi,\theta,\psi,\varphi)-P(R,r,\phi,-\theta,\psi,\varphi)]F(\phi,\psi)\cos\psi d\psi d\phi.
\]
This function $U$ has zero values on the radii planes $\theta=0$ and $\theta=\pi$ of the sphere, which are on the $xz$ coordinates plane. %, as one would expect when $U$ is interpreted as a steady temperature.
Formula~(\ref{eq:chap_13_3_1}) thus solve the Dirichlet problem for the semi-sphere region $r<R, 0<\theta<\pi$, where $U=0$ on the disk $r<R$ which is denoted by $C_{R}$ and on the $xz$ coordinates plane, and
\begin{equation}\label{eq:chap_13_3_2}
\lim_{r\rightarrow R, r<R}U(r,\theta,\varphi)=F(\theta,\varphi)\texttt{ }(0<\theta<\pi)
\end{equation}
for each fixed $\theta$ at which $F$ is continuous.

If $F(2\pi-\theta,\varphi)=F(\theta,\varphi)$, then
\begin{equation}\label{eq:chap_13_3_3}
U(r,\theta,\varphi)=
\end{equation}
\[
\frac{1}{4\pi}
\int_{-\pi/2}^{\pi/2}\int_{0}^{2\pi}[P(R,r,\phi,\theta,\psi,\varphi)+P(R,r,\phi,-\theta,\psi,\varphi)]F(\phi,\psi)\cos\psi d\psi d\phi;
\]
and $U_{\theta}=0$ when $\theta=0$ or $\theta=\pi$. Hence formula~(\ref{eq:chap_13_3_3}) furnishes a function $U$ that is complex harmonic in the semi-sphere region $r<R, 0<\theta<\pi$ and satisfies condition~(\ref{eq:chap_13_3_2}) as well as the condition that its normal derivative be zero on the disk $C_{R}$.

The analytic function $s=R^{2}/s_{e}$ maps the sphere $|s_{e}|=R$ in the $s_{e}$ space onto the sphere $|s|=R$ in the $s$ space, and it maps the exterior of the first sphere onto the interior of the second. Writing
\[
s=r(e_{xy}\exp(i\theta)\cos\varphi+\sin\varphi)\texttt{ and } s_{e}=r_{e}(e_{xy}\exp(i\theta_{e})\cos\varphi_{e}+\sin\varphi_{e}),
\]
we note that $r=R^{2}/r_{e}$ and $\theta=2\pi-\theta_{e}$. The complex harmonic function $U$ represented by formula~(\ref{eq:chap_13_2_1}), Sec.~(\ref{sec:chap_13_2_123}), is, then, transformed into the function
\[
U(\frac{R^{2}}{r_{e}},2\pi-\theta_{e},\varphi_{e})=-\frac{1}{4\pi}
\int_{-\pi/2}^{\pi/2}\int_{0}^{2\pi}P(R,r_{e},\phi,-\theta_{e},\psi,\varphi_{e})F(\phi,\psi)\cos\psi d\psi d\phi,
\]
which is complex harmonic in the domain $r_{e}>R$. Now, in general, if $\tau(r,\theta,\varphi)$ is complex harmonic, then so is $\tau(r,-\theta,\varphi)$. Hence the function $H(r_{e},\theta_{e},\varphi_{e})=U(\frac{R^{2}}{r_{e}},2\pi-\theta_{e},\varphi_{e})$, or
\begin{equation}\label{eq:chap_13_3_4}
H(r_{e},\theta_{e},\varphi_{e})=
\end{equation}
\[
-\frac{1}{4\pi}
\int_{-\pi/2}^{\pi/2}\int_{0}^{2\pi}P(R,r_{e},\phi,\theta_{e},\psi,\varphi_{e})F(\phi,\psi)\cos\psi d\psi d\phi, (r_{e}>R),
\]
is also complex harmonic. For each fixed $\theta_{e}$ at which $F(\theta_{e},\varphi_{e})$ is continuous, we find from condition~(\ref{eq:chap_13_2_2}), Sec.~(\ref{sec:chap_13_2_123}), that
\begin{equation}\label{eq:chap_13_3_5}
\lim_{r_{e}\rightarrow R, r_{e}<R}H(r_{e},\theta_{e},\varphi_{e})=F(\theta_{e},\varphi_{e}).
\end{equation}

The formula~(\ref{eq:chap_13_3_4}) solves the Dirichlet problem for the region exterior to the sphere $r_{e}=R$ in the $s_{e}$ space. We note from expression~(\ref{eq:chap_13_1_7}), Sec.~(\ref{sec:chap_13_1_122}), that the Poisson kernel $P(R,r_{e},\phi,\theta_{e},\psi,\varphi_{e})$ is negative when $r_{e}>R$. Also,
\begin{equation}\label{eq:chap_13_3_6}
\frac{1}{4\pi}
\int_{-\pi/2}^{\pi/2}\int_{0}^{2\pi}P(R,r_{e},\phi,\theta_{e},\psi,\varphi_{e})\cos\psi d\psi d\phi=-1\texttt{ }(r_{e}>R)
\end{equation}
and
\begin{equation}\label{eq:chap_13_3_7}
\lim_{r_{e}\rightarrow \infty}H(r_{e},\theta_{e},\varphi_{e})=\frac{1}{4\pi}
\int_{-\pi/2}^{\pi/2}\int_{0}^{2\pi}F(\phi,\psi)\cos\psi d\psi d\phi.
\end{equation}

\section{Schwarz Integral Formula}\label{sec:chap_13_4_125}

Let $f$ be an analytic function of $s$ throughout the half space $Im$ $s\geq0$ such that, for some positive constants $a$ and $M$, $f$ satisfies the order property
\begin{equation}\label{eq:chap_13_4_1}
|s^{a}f(s)|<M \texttt{ $(Im$ }s\geq0).
\end{equation}
For a fixed point $s$ above the $xz$ coordinates plane which is supposed horizontal, let $S_{R}$ denote the upper half of a positively oriented sphere of radius $R$ centered at the origin, where $R>|s|$, and $C_{R}$ denote the disk of radius $R$ centered at the origin and in the $xz$ coordinates plane. Then, according to the Cauchy surface integral formula,
\begin{equation}\label{eq:chap_13_4_2}
f(s)=\frac{1}{4\pi}\iint_{S_{R}}\frac{f(\varsigma)}{(\varsigma-s)^{2}}d\sigma+\frac{1}{4\pi}\iint_{C_{R}}\frac{f(\varsigma)}{(\varsigma-s)^{2}}d\sigma
\end{equation}

We find that the first of these integrals approaches 0 as $R$ tends to $\infty$ since, in view of condition~(\ref{eq:chap_13_4_1}),
\[
|\iint_{S_{R}}\frac{f(\varsigma)}{(\varsigma-s)^{2}}d\sigma|<\frac{M}{R^{a}(R-|s|)^{2}}2\pi R^{2}=\frac{2M\pi}{R^{a}(1-|s|/R)^{2}}
\]
Thus, using $s=e_{xy}(x+iy)+z$ and $\bar{s}=e_{xy}(x-iy)+z$, when let
\[
\tau=e_{xy}\xi+\zeta,\texttt{ }d\sigma=\frac{\partial\tau}{\partial \xi}d\xi\frac{\partial\tau}{\partial \zeta}d\zeta=e_{xy}d\xi d\zeta,
\]
we have $\tau=\bar{\tau}$ and
\begin{equation}\label{eq:chap_13_4_3}
f(s)=\frac{1}{4\pi}\iint_{C_{R}}\frac{f(\varsigma)}{(\varsigma-s)^{2}}d\sigma
=e_{xy}\frac{1}{4\pi}\int_{-\infty}^{\infty}\int_{-\infty}^{\infty}\frac{f(\tau)}{(\tau-s)^{2}}d\xi d\zeta \texttt{ $(Im$ }s>0).
\end{equation}

Condition~(\ref{eq:chap_13_4_1}) also ensures that the improper integral here converges. The number to which it converges is the same as its Cauchy principal value, and representation~(\ref{eq:chap_13_4_3}) is a Cauchy surface integral formula for the half space $Im$ $s>0$.

When the point $s$ lies below the $xz$ coordinates plane, the right-hand side of equation~(\ref{eq:chap_13_4_2}) is zero; hence integral~(\ref{eq:chap_13_4_3}) is zero for such a point. Thus, when $s$ is above the $xz$ coordinates plane, we have the following formula, where $c$ is an arbitrary constant:
\begin{equation}\label{eq:chap_13_4_4}
f(s)=\frac{1}{4\pi}\int_{-\infty}^{\infty}\int_{-\infty}^{\infty}[\frac{1}{(\tau-s)^{2}}+\frac{c}{(\tau-\bar{s})^{2}}]f(\tau)d\xi d\zeta \texttt{ $(Im$ }s>0).
\end{equation}
In the two cases $c=-1$ and $c=1$, this reduces, respectively, to
\begin{equation}\label{eq:chap_13_4_5}
\begin{split}
f(s)&=\frac{1}{4\pi}\int_{-\infty}^{\infty}\int_{-\infty}^{\infty}[\frac{(\bar{\tau}-\bar{s})^{2}}{(\tau-s)^{2}(\bar{\tau}-\bar{s})^{2}}-\frac{(\bar{\tau}-s)^{2}}{(\tau-\bar{s})^{2}(\bar{\tau}-s)^{2}}]f(\tau)d\xi d\zeta \\
&=e_{xy}\frac{1}{\pi}\int_{-\infty}^{\infty}\int_{-\infty}^{\infty}\frac{iy[(\xi-x)+(\zeta-z)]}{\{[(\xi-x)+(\zeta-z)]^{2}+y^{2}\}^{2}}f(\tau)d\xi d\zeta \texttt{ }(y>0)
\end{split}
\end{equation}
and
\begin{equation}\label{eq:chap_13_4_6}
\begin{split}
f(s)&=\frac{1}{4\pi}\int_{-\infty}^{\infty}\int_{-\infty}^{\infty}[\frac{(\bar{\tau}-\bar{s})^{2}}{(\tau-s)^{2}(\bar{\tau}-\bar{s})^{2}}+\frac{(\bar{\tau}-s)^{2}}{(\tau-\bar{s})^{2}(\bar{\tau}-s)^{2}}]f(\tau)d\xi d\zeta \\
&=\frac{1}{2\pi}\int_{-\infty}^{\infty}\int_{-\infty}^{\infty}\frac{[e_{xy}(\xi-x)+(\zeta-z)]^{2}-e_{xy}y^{2}}{\{[e_{xy}(\xi-x)+(\zeta-z)]^{2}+e_{xy}y^{2}\}^{2}}f(\tau)d\xi d\zeta \texttt{ }(y>0).
\end{split}
\end{equation}

If $f(s)=e_{xy}[u(x,y,z)+iv(x,y,z)]+w(z)$, it follows from formulas~(\ref{eq:chap_13_4_5}) and~(\ref{eq:chap_13_4_6}) that the complex harmonic functions $u=Re$ $f(s)$ and $v=Im$ $f(s)$, respectively, are represented in the half space $y>0$ in terms of the boundary values of $f(\tau)=e_{xy}u+w$ by the formulas
\begin{equation}\label{eq:chap_13_4_7}
v(x,y,z)=\frac{1}{\pi}\int_{-\infty}^{\infty}\int_{-\infty}^{\infty}\frac{y[(\xi-x)+(\zeta-z)][u(\xi,0,\zeta)+w(\zeta)]}
{\{[(\xi-x)+(\zeta-z)]^{2}+y^{2}\}^{2}}d\xi d\zeta
\end{equation}
where $y>0$ and
\begin{equation}\label{eq:chap_13_4_8}
e_{xy}u(x,y,z)+w(z)
\end{equation}
\[
=\frac{1}{2\pi}\int_{-\infty}^{\infty}\int_{-\infty}^{\infty}\frac{[e_{xy}(\xi-x)+(\zeta-z)]^{2}-e_{xy}y^{2}}{\{[e_{xy}(\xi-x)+(\zeta-z)]^{2}+e_{xy}y^{2}\}^{2}}[e_{xy}u(\xi,0,\zeta)+w(\zeta)]d\xi d\zeta
\]
\[
=\frac{1}{2\pi}\int_{-\infty}^{\infty}\int_{-\infty}^{\infty}\frac{e_{xy}g(x,\xi,y,z,\zeta)+h(z,\zeta)}{\{[e_{xy}(\xi-x)+(\zeta-z)]^{2}+e_{xy}y^{2}\}^{2}}[e_{xy}u(\xi,0,\zeta)+w(\zeta)]d\xi d\zeta
\]
where $y>0$ and
\[
g(x,\xi,y,z,\zeta)=(\xi-x)^{2}-y^{2}+2(\xi-x)(\zeta-z)\texttt{ and }h(z,\zeta)=(\zeta-z)^{2}.
\]
Formula~(\ref{eq:chap_13_4_7}) is known as the Schwarz integral formula, or the Poisson integral formula for the half space. In the next section, we shall relax the conditions for the validity of formulas~(\ref{eq:chap_13_4_7}) and~(\ref{eq:chap_13_4_8}).

\section{Dirichlet Problem for a Half Space}\label{sec:chap_13_5_126}

Let $F$ denote a $xz$ plane function of $e_{xy}x+z$ that is bounded for all $e_{xy}x+z$ and continuous except for at most a finite number of finite jumps. when $y\geq\varepsilon$ and $|e_{xy}x+z|\leq1/\varepsilon$, where $\varepsilon$ is any positive constant, the integral
\[
I(e_{xy}x+z,y)=\int_{-\infty}^{\infty}\int_{-\infty}^{\infty}\frac{[e_{xy}(\xi-x)+(\zeta-z)]F(e_{xy}\xi+\zeta)}{\{[e_{xy}(\xi-x)+(\zeta-z)]^{2}+e_{xy}y^{2}\}^{2}}d\xi d\zeta
\]
converges uniformly with respect to $e_{xy}x+z$ and $y$, as do the integrals of the partial derivatives
of the integrand with respect to $e_{xy}x+z$ and $y$. Each of these integrals is the sum of a finite number of improper or definite integrals over regions where $F$ is continuous; hence the integrand of each component integral is a continuous function of $e_{xy}\xi+\zeta$, $e_{xy}x+z$, and $y$ when $y\geq\varepsilon$. Consequently, each partial derivative of $I(e_{xy}x+z,y)$ is represented by the integral of the corresponding derivative of the integrand whenever $y>\varepsilon$.

We write $V(x,y,z)=yI(e_{xy}x+z,y)/\pi$. Thus $V$ is the Schwarz integral transform of $F$, suggested by the second of expressions~(\ref{eq:chap_13_4_7}), Sec.~(\ref{sec:chap_13_4_125}):
\begin{equation}\label{eq:chap_13_5_1}
V(x,y,z)=\frac{1}{\pi}\int_{-\infty}^{\infty}\int_{-\infty}^{\infty}\frac{y[(\xi-x)+(\zeta-z)]F(e_{xy}\xi+\zeta)}{\{[(\xi-x)+(\zeta-z)]^{2}+y^{2}\}^{2}}d\xi d\zeta
\end{equation}
where $y>0$. Except for the factor $1/\pi$, the kernel here is
\[
\frac{y[(\xi-x)+(\zeta-z)]}{\{[e_{xy}(\xi-x)+(\zeta-z)]^{2}+e_{xy}y^{2}\}^{2}}.
\]
It is the imaginary component of the function $1/(\tau-s)^{2}-1/(\tau-\bar{s})^{2}$, which is analytic in $s$ when $y>0$. It follows that the kernel is complex harmonic, and so it satisfies the complex Laplace's equation in $x$, $y$, and $z$. Because the order of differentiation and integration may be interchanged, the function~(\ref{eq:chap_13_5_1}) then satisfies that equation. Consequently, $V$ is complex harmonic when $y>0$.

To prove that
\begin{equation}\label{eq:chap_13_5_2}
\lim_{y\rightarrow0,y>0}V(x,y,z)=F(e_{xy}x+z)
\end{equation}
for each fixed $e_{xy}x+z$ at which $F$ is continuous, we substitute $\tau=e_{xy}(x+y\tan\vartheta)+z$ in formula~(\ref{eq:chap_13_5_1}) and write
\begin{equation}\label{eq:chap_13_5_3}
V(x,y,z)=\frac{1}{\pi}\int_{-\pi/2}^{\pi/2}F[e_{xy}(x+y\tan\vartheta)+z]d\vartheta,(y>0).
\end{equation}
Then, if
\[
G(e_{xy}x+z,y,\vartheta)=F[e_{xy}(x+y\tan\vartheta)+z]-F(e_{xy}x+z)
\]
and $\alpha$ is some small positive constant,
\begin{equation}\label{eq:chap_13_5_4}
\begin{split}
\pi[V(x,y,z)-F(e_{xy}x+z)]&=\int_{-\pi/2}^{\pi/2}G(e_{xy}x+z,y,\vartheta)d\vartheta \\
&=I_{1}(y)+I_{2}(y)+I_{3}(y),
\end{split}
\end{equation}
where
\[
\begin{split}
I_{1}(y)&=\int_{-\pi/2}^{(-\pi/2)+\alpha}G(e_{xy}x+z,y,\vartheta)d\vartheta,\\
I_{2}(y)&=\int_{(-\pi/2)+\alpha}^{(\pi/2)-\alpha}G(e_{xy}x+z,y,\vartheta)d\vartheta,\\
I_{3}(y)&=\int_{(\pi/2)-\alpha}^{\pi/2}G(e_{xy}x+z,y,\vartheta)d\vartheta.
\end{split}
\]

If $M$ denote an upper bound for $|F(e_{xy}x+z)|$, then $|G(e_{xy}x+z,y,\vartheta)|\leq2M$. For a given positive numbers $\varepsilon$, we select $\alpha$ so that $6M\alpha<\varepsilon$; and this means that
\[
|I_{1}(y)|\leq2M\alpha<\frac{\varepsilon}{3} \texttt{ and } |I_{3}(y)|\leq2M\alpha<\frac{\varepsilon}{3}.
\]
We next show that, corresponding to $\varepsilon$, there is a positive number $\delta$ such that
\[
|I_{2}(y)|<\frac{\varepsilon}{3} \texttt{ whenever } 0<y<\delta.
\]
To do this, we observe that, since $F$ is continuous at $e_{xy}x+z$, there is a positive number $\gamma$ such that
\[
|G(e_{xy}x+z,y,\vartheta)|<\frac{\varepsilon}{3\pi} \texttt{ whenever } 0<y|\tan\vartheta|<\gamma.
\]
Now the maximum value of $|\tan\vartheta|$ as $\vartheta$ ranges from $(-\pi/2)+\alpha$ to $(\pi/2)-\alpha$ is $\tan(\pi/2-\alpha)=\cot\alpha$. Hence, if we write $\delta=\gamma\tan\alpha$, it follows that
\[
|I_{2}(y)|<\frac{\varepsilon}{3\pi}(\pi-2\alpha)<\frac{\varepsilon}{3} \texttt{ whenever } 0<y<\delta.
\]

We have thus shown that
\[
|I_{1}(y)|+|I_{2}(y)|+|I_{3}(y)|<\varepsilon \texttt{ whenever } 0<y<\delta.
\]
Condition~(\ref{eq:chap_13_5_2}) now follows from this result and equation~(\ref{eq:chap_13_5_4}).

Formula~(\ref{eq:chap_13_5_1}) therefore solves the Dirichlet problem for the half space $y>0$, with the boundary condition~(\ref{eq:chap_13_5_2}). It is evident from the form~(\ref{eq:chap_13_5_3}) of expression~(\ref{eq:chap_13_5_1}) that $|V(x,y,z)|\leq M$ in the half space, where $M$ is an upper bound of $|F(e_{xy}x+z)|$; that is, $V$ is bounded. We note that $V(x,y,z)=F_{0}$ when $F(e_{xy}x+z)=F_{0}$, where $F_{0}$ is a constant.

According to the second expression of formula~(\ref{eq:chap_13_4_8}) of Sec.~(\ref{sec:chap_13_4_125}), under certain conditions of $F$ the function
\begin{equation}\label{eq:chap_13_5_5}
U(x,y,z)=Re\frac{1}{2\pi}\int_{-\infty}^{\infty}\int_{-\infty}^{\infty}
\frac{[e_{xy}g(x,\xi,y,z,\zeta)+h(z,\zeta)]F(e_{xy}\xi+\zeta)}{\{[e_{xy}(\xi-x)+(\zeta-z)]^{2}+e_{xy}y^{2}\}^{2}}d\xi d\zeta
\end{equation}
where
\[
g(x,\xi,y,z,\zeta)=(\xi-x)^{2}-y^{2}+2(\xi-x)(\zeta-z)\texttt{ and }h(z,\zeta)=(\zeta-z)^{2}
\]
is a complex harmonic conjugate of $V$ given by formula~(\ref{eq:chap_13_5_1}). Actually, formula~(\ref{eq:chap_13_5_5}) furnishes a complex harmonic conjugate of $V$ if $F$ is everywhere continuous, except for at most a finite number of finite jumps, and if $F$ satisfies an order property
\[
|x^{a}F(e_{xy}x+z)|<M\texttt{ }(a<0).
\]
For, under these conditions, we find that $U$ and $V$ satisfy the Cauchy-Riemann equations where $y>0$.

\section{Neumann Problems}\label{sec:chap_13_6_127}

As in Section~(\ref{sec:chap_13_1_122}), let $S$ denote a positively oriented sphere, centered at the origin. We let $R$ denote the radius of $S$ and write
\[
\varsigma=R(e_{xy}e^{i\phi}\cos\psi+\sin\psi)\texttt{ and }s=r(e_{xy}e^{i\theta}\cos\varphi+\sin\varphi)\texttt{ }(r<R).
\]
In view of equation~(\ref{eq:chap_13_1_6}), when $\varsigma$ is fixed, the function
\begin{equation}\label{eq:chap_13_6_1}
Q(R,r,\phi,\theta,\psi,\varphi)=-R\{Re[\log(R^{2}-2Rr\sigma(\phi,\theta,\psi,\varphi)+r^{2})
\end{equation}
\[
+\frac{2Rr\sigma_{m}(\phi,\theta,\psi,\varphi)}{R^{2}-2Rr\sigma(\phi,\theta,\psi,\varphi)+r^{2}}]
-[\log(R^{2}-2Rr\sigma(\psi,\varphi)+r^{2})+\frac{2Rr\sigma_{m}(\psi,\varphi)}{R^{2}-2Rr\sigma(\psi,\varphi)+r^{2}}]\}
\]
is a complex harmonic function of $r, \theta$, and $\varphi$ interior to $S$ for each fixed point $\varsigma$ on the sphere $S$ because it is the real component of
\[
-R\{\log[(1-\varsigma/s)(r^{2}-R^{2}s/\varsigma)]+\frac{r^{2}\varsigma/s-R^{2}s/\varsigma}{(1-\varsigma/s)(r^{2}-R^{2}s/\varsigma)}\},
\]
where the branch cut of
\[
\log[(1-\varsigma/s)(r^{2}-R^{2}s/\varsigma)]+\frac{r^{2}\varsigma/s-R^{2}s/\varsigma}{(1-\varsigma/s)(r^{2}-R^{2}s/\varsigma)}
\]
is an outward ray from the point $\varsigma$. If, moreover, $r\neq0$,
\begin{equation}\label{eq:chap_13_6_2}
\frac{R}{r}P(R,r,\phi,\theta,\psi,\varphi)=Q_{r}(R,r,\phi,\theta,\psi,\varphi)
\end{equation}
\[
=\frac{R}{r}\{Re[(1-\frac{2Rr\sigma_{m}(\phi,\theta,\psi,\varphi)}{R^{2}-2Rr\sigma(\phi,\theta,\psi,\varphi)+r^{2}})
\frac{R^{2}-r^{2}}{R^{2}-2Rr\sigma(\phi,\theta,\psi,\varphi)+r^{2}}-1]
\]
\[
-[(1-\frac{2Rr\sigma_{m}(\psi,\varphi)}{R^{2}-2Rr\sigma(\psi,\varphi)+r^{2}})
\frac{R^{2}-r^{2}}{R^{2}-2Rr\sigma(\psi,\varphi)+r^{2}}-1]\}
\]
where $P$ is the Poisson kernel~(\ref{eq:chap_13_1_7}) of Sec.~(\ref{sec:chap_13_1_122}).

These observations suggest that the function $Q$ may be used to write an integral representation for a harmonic function $U$ whose normal derivative $U_{r}$ on the sphere $r=R$ assumes prescribed values $G(\theta,\varphi)$.

If $G$ is piecewise continuous and $U_{0}$ is an arbitrary constant, the function
\begin{equation}\label{eq:chap_13_6_3}
U(r,\theta,\varphi)=\frac{1}{4\pi}\int_{-\pi/2}^{\pi/2}\int_{0}^{2\pi}Q(R,r,\phi,\theta,\psi,\varphi)G(\phi,\psi)\cos\psi d\psi d\phi+U_{0}
\end{equation}
where $r<R$, is complex harmonic because the integrand is a Neumann kernel function of $r,\varphi$ and $\theta$. If the mean value of $G$ over the sphere $|s|=R$ is zero, or
\begin{equation}\label{eq:chap_13_6_4}
\int_{-\pi/2}^{\pi/2}\int_{0}^{2\pi}G(\phi,\psi)\cos\psi d\psi d\phi=0,
\end{equation}
then, in view of equation~(\ref{eq:chap_13_6_2}),
\[
\begin{split}
U_{r}(r,\theta,\varphi)&=\frac{1}{4\pi}\int_{-\pi/2}^{\pi/2}\int_{0}^{2\pi}\frac{R}{r}P(R,r,\phi,\theta,\psi,\varphi)G(\phi,\psi)\cos\psi d\psi d\phi \\
&=\frac{R}{r}\frac{1}{4\pi}\int_{-\pi/2}^{\pi/2}\int_{0}^{2\pi}P(R,r,\phi,\theta,\psi,\varphi)G(\phi,\psi)\cos\psi d\psi d\phi.
\end{split}
\]
Now, according to equations~(\ref{eq:chap_13_2_1}) and~(\ref{eq:chap_13_2_2}) of Sec.~(\ref{sec:chap_13_2_123}),
\[
\lim_{r\rightarrow R,r<R}\frac{1}{4\pi}
\int_{-\pi/2}^{\pi/2}\int_{0}^{2\pi}P(R,r,\phi,\theta,\psi,\varphi)G(\phi,\psi)\cos\psi d\psi d\phi=G(\theta,\varphi).
\]
Hence
\begin{equation}\label{eq:chap_13_6_5}
\lim_{r\rightarrow R,r<R}U_{r}(r,\theta,\varphi)=G(\theta,\varphi)
\end{equation}
for each pair of values of $\theta$ and $\varphi$ at which $G$ is continuous.

When $G$ is piecewise continuous and satisfies condition~(\ref{eq:chap_13_6_4}), the formula
\begin{equation}\label{eq:chap_13_6_6}
\begin{split}
U(r,\theta,\varphi)=&-\frac{R}{4\pi}\int_{-\pi/2}^{\pi/2}\int_{0}^{2\pi}
\{Re[\log(R^{2}-2Rr\sigma(\phi,\theta,\psi,\varphi)+r^{2}) \\
&+\frac{2Rr\sigma_{m}(\phi,\theta,\psi,\varphi)}{R^{2}-2Rr\sigma(\phi,\theta,\psi,\varphi)+r^{2}}]
-[\log(R^{2}-2Rr\sigma(\psi,\varphi)+r^{2}) \\
&+\frac{2Rr\sigma_{m}(\psi,\varphi)}{R^{2}-2Rr\sigma(\psi,\varphi)+r^{2}}]\}
G(\phi,\psi)\cos\psi d\psi d\phi+U_{0}\texttt{ }(r<R),
\end{split}
\end{equation}
therefore, solves the Neumann problem for the region interior to the sphere $r=R$, where $G(\theta,\varphi)$ is the normal derivative of the harmonic function $U(r,\theta,\varphi)$ at the boundary in the sense of condition~(\ref{eq:chap_13_6_5}). Note how it follows from equations~(\ref{eq:chap_13_6_4}) and~(\ref{eq:chap_13_6_6}) that, since $\ln R^{2}$ is constant, $U_{0}$ is the value of $U$ at the center $r=0$ of the sphere $r=R$.

%The values $U(r,\theta,\varphi)$ may represent steady temperatures in a sphere $r<R$ with insulated faces. In that case, condition~(\ref{eq:chap_13_6_5}) states that the flux of heat into the sphere through its edge is proportional to $G(\theta,\varphi)$. Condition~(\ref{eq:chap_13_6_4}) is the natural physical requirement that the total rate of flow of heat into the sphere be zero, since temperatures do not vary with time.

A corresponding formula for a complex harmonic function $H$ in the region exterior to the sphere $r=R$ can be written in terms of $Q$ as
\begin{equation}\label{eq:chap_13_6_7}
H(r_{e},\theta_{e},\varphi_{e})=\frac{1}{4\pi}\int_{-\pi/2}^{\pi/2}\int_{0}^{2\pi}Q(R,r_{e},\phi,\theta_{e},\psi,\varphi_{e})G(\phi,\psi)\cos\psi d\psi d\phi+H_{0}
\end{equation}
where $r_{e}>R$ and $H_{0}$ is a constant. As before we assume that $G$ is piecewise continuous and that condition~(\ref{eq:chap_13_6_4}) holds. Then
\[
H_{0}=\lim_{R\rightarrow \infty}H(r_{e},\theta_{e},\varphi_{e})
\]
and
\begin{equation}\label{eq:chap_13_6_8}
\lim_{r_{e}\rightarrow R, r_{e}>R}H_{r_{e}}(r_{e},\theta_{e},\varphi_{e})=G(\theta_{e},\varphi_{e})
\end{equation}
for each pair of values of $\theta_{e}$ and $\varphi_{e}$ at which $G$ is continuous.

Turn now to a half space, we let $G(e_{xy}x+z)$ be continuous for all two-dimensional variables $e_{xy}x+z$ on the $xz$ plane, except possibly for a finite number of finite jumps, and let it satisfy an order property
\begin{equation}\label{eq:chap_13_6_9}
|x^{a}G(e_{xy}x+z)|<M\texttt{ }(a>1)
\end{equation}
when $-\infty<x<\infty$ and $-\infty<z<\infty$. For each fixed two-dimensional variable $\tau=e_{xy}\xi+\zeta$ on the $xz$ plane, the function $1/(\tau-s)+1/(\tau-\bar{s})$ is complex harmonic in the half space $Im$ $s>0$. Consequently, the function
\begin{equation}\label{eq:chap_13_6_10}
U(x,y,z)=\frac{-1}{4\pi}\int_{-\infty}^{\infty}\int_{-\infty}^{\infty}(\frac{1}{\tau-s}+\frac{1}{\tau-\bar{s}})G(e_{xy}\xi+\zeta)d\xi d\zeta+U_{0}
\end{equation}
\[
=\frac{-1}{2\pi}\int_{-\infty}^{\infty}\int_{-\infty}^{\infty}\frac{[e_{xy}(\xi-x)+(\zeta-z)]G(e_{xy}\xi+\zeta)}{[e_{xy}(\xi-x)+(\zeta-z)]^{2}+e_{xy}y^{2}}d\xi d\zeta+U_{0}\texttt{ }(y>0)
\]
where $U_{0}$ is a constant, is complex harmonic in that half space.

Formula~(\ref{eq:chap_13_6_10}) was written with the Schwarz integral transform~(\ref{eq:chap_13_5_1}), Sec.~(\ref{sec:chap_13_5_126}), in mind; for it follows from formula~(\ref{eq:chap_13_6_10}) that
\begin{equation}\label{eq:chap_13_6_11}
U_{y}(x,y,z)=\frac{1}{\pi}\int_{-\infty}^{\infty}\int_{-\infty}^{\infty}\frac{y[(\xi-x)+(\zeta-z)]G(e_{xy}\xi+\zeta)}{\{[(\xi-x)+(\zeta-z)]^{2}+y^{2}\}^{2}}d\xi d\zeta+U_{0}
\end{equation}
where $y>0$. In view of equations~(\ref{eq:chap_13_5_1}) and~(\ref{eq:chap_13_5_2}) of Sec.~(\ref{sec:chap_13_5_126}), then
\begin{equation}\label{eq:chap_13_6_12}
\lim_{y\rightarrow 0,y>0}U_{y}(x,y,z)=G(e_{xy}x+z)
\end{equation}
at each point $e_{xy}x+z$ on the $xz$ coordinates plane where $G$ is continuous.

Integral formula~(\ref{eq:chap_13_6_10}) evidently solves the Neumann problem for the half space $y>0$, with boundary condition~(\ref{eq:chap_13_6_12}). But we have not presented conditions on $G$ that are sufficient to ensure that the complex harmonic function $U$ is bounded as $|s|$ increases.

When $G$ is an odd function, formula~(\ref{eq:chap_13_6_10}) may be written
\begin{equation}\label{eq:chap_13_6_13}
U(x,y,z)=\frac{-1}{2\pi}\int_{0}^{\infty}\int_{-\infty}^{\infty}\{\frac{[e_{xy}(\xi-x)+(\zeta-z)]}{[e_{xy}(\xi-x)+(\zeta-z)]^{2}+e_{xy}y^{2}}
\end{equation}
\[
-\frac{[e_{xy}(\xi+x)+(\zeta+z)]}{[e_{xy}(\xi+x)+(\zeta+z)]^{2}+e_{xy}y^{2}}\}G(e_{xy}\xi+\zeta)d\xi d\zeta\texttt{ }(x>0,y>0,z>0).
\]
This represents a function that is complex harmonic in the first quadrant $x>0$, $y>0$, and $z>0$, and satisfies the boundary conditions
\begin{equation}\label{eq:chap_13_6_14}
U(0,y,0)=0\texttt{ }(y>0),
\end{equation}
\begin{equation}\label{eq:chap_13_6_15}
\lim_{y\rightarrow 0,y>0}U_{y}(x,y,z)=G(e_{xy}x+z)\texttt{ }(x>0,z>0).
\end{equation}

%-----------------------------------------------------------------------
% Beginning of chap14.tex
%-----------------------------------------------------------------------
%
% AMS-LaTeX 1.2 sample file for a monograph, based on amsbook.cls.
% This is a data file input by chapter.tex.
%%%%%%%%%%%%%%%%%%%%%%%%%%%%%%%%%%%%%%%%%%%%%%%%%%%%%%%%%%%%%%%%%%%%

%\part{This is a Part Title Sample}

\chapter{Applications of Residues}\label{ch:chap_14}

We turn now to some important applications of the theory of residues, which was developed in the preceding chapter. The applications include evaluation of certain types of definite and improper integrals occurring in real analysis and applied mathematics. Considerable attention is also given to a method, based on residues, for locating zeros of functions and to finding inverse Laplace transforms by summing residues.

\section{Evaluation of Improper Integrals}\label{sec:chap_14_1_71}

Let $\chi$ denote a two-dimensional variable and $d\chi$ denote a two-dimensional area in the $xz$ coordinate plane as follows:
\[
\chi=e_{xy}x+z,\texttt{ }d\chi=\frac{\partial\chi}{\partial x}dx\frac{\partial\chi}{\partial z}dz=e_{xy}dx dz.
\]
Then the improper integral of a continuous function $f(\chi)$ over the first quadrant $x\geq0$ and $z\geq0$ of the infinite $xz$ coordinate plane is defined by means of the equation
\begin{equation}\label{eq:chap_14_1_1}
\int_{0}^{\infty}\int_{0}^{\infty}f(\chi)d\chi=\lim_{R \to \infty}\int_{0}^{R}\int_{0}^{R}f(\chi)d\chi.
\end{equation}
When the limit on the right exists, the improper integral is said to converge to that limit. If $f(\chi)$ is continuous for all points $\chi$, its improper integral over the infinite region $-\infty<x<\infty$ and $-\infty<z<\infty$ is defined by writing
\begin{equation}\label{eq:chap_14_1_2}
\int_{-\infty}^{\infty}\int_{-\infty}^{\infty}f(\chi)d\chi=\lim_{R_{1} \to \infty}\int_{-R_{1}}^{0}\int_{-R_{1}}^{0}f(\chi)d\chi+\lim_{R_{2} \to \infty}\int_{0}^{R_{2}}\int_{0}^{R_{2}}f(\chi)d\chi
\end{equation}
\[
+\lim_{R_{1} \to \infty, R_{2} \to \infty}\int_{-R_{1}}^{0}\int_{0}^{R_{2}}f(\chi)d\chi
+\lim_{R_{1} \to \infty, R_{2} \to \infty}\int_{0}^{R_{2}}\int_{-R_{1}}^{0}f(\chi)d\chi;
\]
and when all of the limits here exist, integral~(\ref{eq:chap_14_1_2}) converges to their sum. Another
value that is assigned to integral~(\ref{eq:chap_14_1_2}) is often useful. Namely, the Cauchy principal value  (P.V.) of integral~(\ref{eq:chap_14_1_2}) is the number
\begin{equation}\label{eq:chap_14_1_3}
P.V.\int_{-\infty}^{\infty}\int_{-\infty}^{\infty}f(\chi)d\chi=\lim_{R \to \infty}\int_{-R}^{R}\int_{-R}^{R}f(\chi)d\chi,
\end{equation}
provided this single limit exists.

If integral~(\ref{eq:chap_14_1_2}) converges, its Cauchy principal value~(\ref{eq:chap_14_1_3}) exists; and that value is the number to which integral~(\ref{eq:chap_14_1_2}) converges. This is because
\[
\int_{-R}^{R}\int_{-R}^{R}f(\chi)d\chi=\int_{-R}^{0}\int_{-R}^{0}f(\chi)d\chi+\int_{-R}^{0}\int_{0}^{R}f(\chi)d\chi
\]
\[
+\int_{0}^{R}\int_{-R}^{0}f(\chi)d\chi+\int_{0}^{R}\int_{0}^{R}f(\chi)d\chi
\]
and the limit as $R \to \infty$ of each of the integrals on the right exists when integral~(\ref{eq:chap_14_1_2})
converges. It is not, however, always true that integral~(\ref{eq:chap_14_1_2}) converges when its Cauchy
principal value exists, as the following example shows.

\begin{example}\label{ex:chap_14_1_1}
Observe that
\begin{equation}\label{eq:chap_14_1_4}
\begin{split}
P.V.\int_{-\infty}^{\infty}\int_{-\infty}^{\infty}\chi d\chi &= \lim_{R \to \infty}e_{xy}\int_{-R}^{R}dx\int_{-R}^{R}(e_{xy}x+z)dz \\
&=e_{xy}\lim_{R \to \infty}R\int_{-R}^{R}2xdx \\
&=e_{xy}\lim_{R \to \infty}Rx^{2}|_{-R}^{R}
=e_{xy}\lim_{R \to \infty}0=0.
\end{split}
\end{equation}
On the other hand,
\begin{equation}\label{eq:chap_14_1_5}
\int_{-\infty}^{\infty}\int_{-\infty}^{\infty}\chi d\chi=\int_{-\infty}^{\infty}\int_{-\infty}^{\infty}(e_{xy}x+z)e_{xy}dxdz
\end{equation}
\[
=e_{xy}[\lim_{R_{2} \to \infty}\int_{0}^{R_{2}}dx\int_{0}^{R_{2}}(x+z)dz
+\lim_{R_{1} \to \infty}\int_{-R_{1}}^{0}dx\int_{-R_{1}}^{0}(x+z)dz]
\]
\[
+e_{xy}\lim_{R_{1} \to \infty, R_{2} \to \infty}[\int_{0}^{R_{2}}dx\int_{-R_{1}}^{0}(x+z)dz
+\int_{-R_{1}}^{0}dx\int_{0}^{R_{2}}(x+z)dz]
\]
\[
=e_{xy}[\lim_{R_{2} \to \infty}\int_{0}^{R_{2}}(xR_{2}+\frac{1}{2}R_{2}^{2})dx
+\lim_{R_{1} \to \infty}\int_{-R_{1}}^{0}(xR_{1}-\frac{1}{2}R_{1}^{2})dx]
\]
\[
+e_{xy}\lim_{R_{1} \to \infty, R_{2} \to \infty}[\int_{0}^{R_{2}}(xR_{1}-\frac{1}{2}R_{1}^{2})dx
+\int_{-R_{1}}^{0}(xR_{2}+\frac{1}{2}R_{2}^{2})dx]
\]
\[
=e_{xy}[\lim_{R_{2} \to \infty}R_{2}^{3}-\lim_{R_{1} \to \infty}R_{1}^{3}
+\lim_{R_{1} \to \infty, R_{2} \to \infty}(R_{1}R_{2}^{2}-R_{1}^{2}R_{2})];
\]
and since these last three limits do not exist, we find that the improper integral~(\ref{eq:chap_14_1_5}) fails
to exist.
\end{example}

But suppose that $f(\chi)$ $(-\infty<x<\infty$, $-\infty<z<\infty)$ is an even function, one where
\[
f(-\chi)=f(\chi)\texttt{ for all }\chi.
\]
The central symmetry of the graph of $y=f(\chi)$ with respect to the $y$ axis in the three-dimensional space enables us to write
\[
\int_{-R}^{R}\int_{0}^{R}f(\chi)d\chi=\int_{0}^{R}\int_{-R}^{R}f(\chi)d\chi=\frac{1}{2}\int_{-R}^{R}\int_{-R}^{R}f(\chi)d\chi,
\]
and we see that integral~(\ref{eq:chap_14_1_1}) converges to one half the Cauchy principal value~(\ref{eq:chap_14_1_3}) when that value exists. Moreover, since integral~(\ref{eq:chap_14_1_1}) converges and since
\[
\int_{-R_{1}}^{0}\int_{-R_{1}}^{0}f(\chi)d\chi=\int_{0}^{R_{1}}\int_{0}^{R_{1}}f(\chi)d\chi
\]
and
\[
\int_{-R_{1}}^{0}\int_{0}^{R_{1}}f(\chi)d\chi=\int_{0}^{R_{1}}\int_{-R_{1}}^{0}f(\chi)d\chi,
\]
integral~(\ref{eq:chap_14_1_2}) converges to twice the value of integral~(\ref{eq:chap_14_1_1}). We have thus shown that when $f(\chi)$ $(-\infty<x<\infty$, $-\infty<z<\infty)$ is even and the Cauchy principal value~(\ref{eq:chap_14_1_3}) exists, both of the integrals~(\ref{eq:chap_14_1_1}) and~(\ref{eq:chap_14_1_2}) converge and
\begin{equation}\label{eq:chap_14_1_6}
\begin{split}
P.V.\int_{-\infty}^{\infty}\int_{-\infty}^{\infty}f(\chi)d\chi
&=\int_{-\infty}^{\infty}\int_{-\infty}^{\infty}f(\chi)d\chi \\
&=2\int_{-\infty}^{\infty}\int_{0}^{\infty}f(\chi)d\chi \\
&=2\int_{0}^{\infty}\int_{-\infty}^{\infty}f(\chi)d\chi.
\end{split}
\end{equation}

We now describe a method involving residues, to be illustrated in the next section, that is often used to evaluate improper integrals of even rational functions $f(\chi)=p(\chi)/q(\chi)$, where $f(-\chi)$ is equal to $f(\chi)$ and where $p(\chi)$ and $q(\chi)$ are polynomials with real coefficients or two-dimensional coefficients of the type $\chi$ and no factors in common. We agree that $q(s)$ has no real zeros but has at least one for $Im$ $s>0$.

The method begins with the identification of all of the distinct zeros of the polynomial $q(s)$ that lie above the $xz$ coordinate plane. They are, of course, finite in number (see Sec.~(\ref{sec:chap_4_14_49})) and may be labeled $s_{1}, s_{2},\ldots, s_{n}$, where $n$ is less than or equal to the degree of $q(s)$. We then integrate the quotient
\begin{equation}\label{eq:chap_14_1_7}
f(s)=\frac{p(s)}{q(s)}
\end{equation}
around the positively oriented boundary of the hemisphere region. That simple closed surface consists of the region of a disk in the $xz$ coordinate plane, which is with a center at origin point and a radius $R$, and the top half of the sphere $|s|=R$, described counterclockwise and denoted by $S_{R}$. It is understood that the positive number $R$ is large enough that the points $s_{1}, s_{2},\ldots, s_{n}$ all lie inside the closed surface.

The Cauchy residue theorem in Sec.~(\ref{sec:chap_7_2_63_2}) and the parametric representation $s=\chi$ $(-R\leq x\leq R$, $-R\leq z\leq R)$ of the region in the $xz$ coordinate plane just mentioned can be used to write
\[
\int_{-R}^{R}\int_{-R}^{R}f(\chi)d\chi+\iint_{S_{R}}f(s)d\sigma=4\pi\sum_{k=1}^{n}Res_{s=s_{k}}f(s),
\]
or
\begin{equation}\label{eq:chap_14_1_8}
\int_{-R}^{R}\int_{-R}^{R}f(\chi)d\chi=4\pi\sum_{k=1}^{n}Res_{s=s_{k}}f(s)-\iint_{S_{R}}f(s)d\sigma.
\end{equation}
If
\[
\lim_{R \to \infty}\iint_{S_{R}}f(s)d\sigma=0,
\]
it then follows that
\begin{equation}\label{eq:chap_14_1_9}
P.V.\int_{-R}^{R}\int_{-R}^{R}f(\chi)d\chi=4\pi\sum_{k=1}^{n}Res_{s=s_{k}}f(s).
\end{equation}
If $f(\chi)$ is even, equations~(\ref{eq:chap_14_1_6}) tell us, moreover, that
\begin{equation}\label{eq:chap_14_1_10}
\int_{-R}^{R}\int_{-R}^{R}f(\chi)d\chi=4\pi\sum_{k=1}^{n}Res_{s=s_{k}}f(s)
\end{equation}
and
\begin{equation}\label{eq:chap_14_1_10}
\int_{0}^{R}\int_{-R}^{R}f(\chi)d\chi=\int_{-R}^{R}\int_{0}^{R}f(\chi)d\chi=2\pi\sum_{k=1}^{n}Res_{s=s_{k}}f(s).
\end{equation}

\section{Improper Integrals from Fourier Analysis}\label{sec:chap_14_3_73}

Residue theory can be useful in evaluating convergent improper integrals of the form
\begin{equation}\label{eq:chap_14_3_1}
\int_{-R}^{R}\int_{-R}^{R}f(\chi)\sin a\chi d\chi\texttt{ or }\int_{-R}^{R}\int_{-R}^{R}f(\chi)\cos a\chi d\chi,
\end{equation}
where $a$ denotes a positive constant. As in Sec.~(\ref{sec:chap_14_1_71}), we assume that $f(\chi)=p(\chi)/q(\chi)$, where $p(\chi)$ and $q(\chi)$ are polynomials with real coefficients or two-dimensional coefficients of the type $\chi$ and no factors in common. Also, $q(\chi)$ has no real zeros. Integrals of type~(\ref{eq:chap_14_3_1}) occur in the theory and application of the Fourier integral.

The method described in Sec.~(\ref{sec:chap_14_1_71}) %and used in Sec.~(\ref{sec:chap_14_2_72})
cannot be applied directly here since (see Sec.~(\ref{sec:chap_3_7_33}))
\[
\begin{split}
|\sin as|^{2} &= \sin^{2}a(x+z)+\sinh^{2}ay-2[\sin a(x+z)\cosh ay-\sin az]\sin az \\
              &\geq(\sin a(x+z)-\sinh ay)^{2}
\end{split}
\]
and
\[
\begin{split}
|\cos as|^{2} &= \cos^{2}a(x+z)+\sinh^{2}ay-2[\cos a(x+z)\cosh ay-\cos az]\cos az. \\
              &\geq(\cos a(x+z)-\sinh ay)^{2}.
\end{split}
\]
More precisely, since
\[
\sinh ay=\frac{e^{ay}-e^{-ay}}{2},
\]
the moduli $|\sin as|$ and $|\cos as|$ increase like $e^{ay}$ as $y $ tends to infinity. The modification
%illustrated in the example below is suggested
by the fact that
\[
\int_{-R}^{R}\int_{-R}^{R}f(\chi)\cos a\chi d\chi+i\int_{-R}^{R}\int_{-R}^{R}f(\chi)\sin a\chi d\chi
=\int_{-R}^{R}\int_{-R}^{R}f(\chi)e^{ia\chi}d\chi,
\]
together with the fact that the modulus
\[
|e^{ias}|=|e^{ia(\chi+e_{xy}iy)}|=|e^{-e_{xy}ay}e^{ia\chi}|=e^{-e_{xy}ay}=e_{xy}e^{-ay}
\]
is bounded in the upper half space $y\geq0$.

\section{Jordan's Lemma}\label{sec:chap_14_4_74}

In the evaluation of integrals of the type treated in Sec.~(\ref{sec:chap_14_3_73}), it is sometimes necessary to use Jordan's lemma, which is stated here as a theorem.

\begin{theorem}\label{th:chap_14_4_1}

(i) a function $f(s)$ is analytic at all points $s$ in the upper half space $y\geq0$ that are
exterior to a sphere $|s|=R_{0}$;

(ii) $S_{R}$ denotes a hemisphere $s=R(e_{xy}e^{i\theta}\cos\varphi+\sin\varphi)$ $(0\leq\theta\leq\pi,\texttt{ }-\pi/2\leq\varphi\leq\pi/2)$, where $R>R_{0}$;

(iii) for all points $s$ on $S_{R}$, there is a positive constant $M_{R}$ such that $|f(s)|\leq M_{R}$,
where
\[
\lim_{R \to \infty}M_{R}R=0.
\]
Then, for every positive constant a,
\begin{equation}\label{eq:chap_14_4_1}
\lim_{R \to \infty}\iint_{S_{R}}f(s)e^{ias}d\sigma=0.
\end{equation}
\end{theorem}

The proof is based on a result that is known as Jordan's inequality:
\begin{equation}\label{eq:chap_14_4_2}
\int_{0}^{\pi}e^{-R\sin\theta}d\theta<\frac{\pi}{R}\texttt{ }(R>0).
\end{equation}
To verify this inequality, we first note from the graphs of the functions $y=\sin\theta$ and
$y=2\theta/\pi$ when $0\leq\theta\leq\pi/2$ that $\sin\theta\geq2\theta/\pi$ for all values of $\theta$ in that interval. Consequently, if $R>0$,
\[
e^{-R\sin\theta}\leq e^{-2R\theta/\pi}\texttt{ when }0\leq\theta\leq\pi/2;
\]
and so
\[
\int_{0}^{\pi/2}e^{-R\sin\theta}d\theta\leq\int_{0}^{\pi/2}e^{-2R\theta/\pi}d\theta=\frac{\pi}{2R}(1-e^{-R}).
\]
Hence
\[
\int_{0}^{\pi/2}e^{-R\sin\theta}d\theta<\frac{\pi}{2R}\texttt{ }(R>0).
\]
But this is just another form of inequality~(\ref{eq:chap_14_4_2}), since the graph of $y=\sin\theta$ is symmetric
with respect to the vertical line $\theta=\pi/2$ on the interval $0\leq\theta\leq\pi$.

Turning now to the verification of limit~(\ref{eq:chap_14_4_1}), we accept statements (i)-(iii) in the theorem and write $d\sigma=R^{2}\cos\varphi d\varphi d\theta$ for $S_{R}$,
\[
\begin{split}
\iint_{S_{R}}f(s)e^{ias}d\sigma &=\int_{-\pi/2}^{\pi/2}\int_{0}^{\pi}f(R(e_{xy}e^{i\theta}\cos\varphi+\sin\varphi))\\
&\times\exp[iaR(e_{xy}e^{i\theta}\cos\varphi+\sin\varphi)]R^{2}\cos\varphi d\varphi d\theta.
\end{split}
\]
Since
\[
|f(R(e_{xy}e^{i\theta}\cos\varphi+\sin\varphi))|\leq M_{R},
\]
\[
|\exp[iaR(e_{xy}e^{i\theta}\cos\varphi+\sin\varphi)]|\leq e^{-aR\sin\theta\cos\varphi},
\]
and in view of Jordan's inequality~(\ref{eq:chap_14_4_2}), it follows that
\[
\begin{split}
|\iint_{S_{R}}f(s)e^{ias}d\sigma| &\leq M_{R}R^{2}\int_{-\pi/2}^{\pi/2}\cos\varphi d\varphi\int_{0}^{\pi}e^{-aR\sin\theta\cos\varphi}d\theta \\
&\leq M_{R}R^{2}\int_{-\pi/2}^{\pi/2}\frac{\pi}{aR\cos\varphi}\cos\varphi d\varphi
\leq\frac{\pi^{2}}{a}M_{R}R.
\end{split}
\]
Limit~(\ref{eq:chap_14_4_1}) is then evident, since $M_{R}R \to 0$ as $R \to \infty$.

\section{Indented Areas}\label{sec:chap_14_5_75}

In this and the following section, we illustrate the use of indented areas. We begin with
an important limit.% that will be used in the example in this section.

\begin{theorem}\label{th:chap_14_5_1}
Suppose that

(i) a function $f(s)$ has a simple pole at a point $s=\chi_{0}=e_{xy}x_{0}+z_{0}$ on the $xz$ coordinates plane, with a Laurent series representation in a punctured sphere $0<|s-\chi_{0}|<R_{2}$ and with residue $B_{0}$;
(ii) $S_{\rho}$ denotes the upper half of a sphere $0<|s-\chi_{0}|=\rho$, where $\rho<R_{2}$ and the clockwise direction is taken.
Then
\begin{equation}\label{eq:chap_14_5_1}
\lim_{\rho \to 0}\iint_{S_{\rho}}f(s)d\sigma=-2\pi B_{0}.
\end{equation}
\end{theorem}

Assuming that the conditions in parts (i) and (ii) are satisfied, we start the proof of the theorem by writing the Laurent series in part (i) as
\[
f(s)=g(s)+\frac{B_{0}}{(s-\chi_{0})^{2}}\texttt{ }(0<|s-\chi_{0}|<R_{2}).
\]
where
\[
g(s)=\sum_{n=0}^{\infty}a_{n}(s-\chi_{0})^{n}\texttt{ }(|s-\chi_{0}|<R_{2}).
\]
Thus,
\begin{equation}\label{eq:chap_14_5_2}
\iint_{S_{\rho}}f(s)d\sigma=\iint_{S_{\rho}}g(s)d\sigma+B_{0}\iint_{S_{\rho}}\frac{d\sigma}{(s-\chi_{0})^{2}}.
\end{equation}
Now the function $g(s)$ is continuous when $|s-\chi_{0}|<R_{2}$, according to the theorem in Sec.~(\ref{sec:chap_6_8_58}). Hence if we choose a number $\rho_{0}$ such that $\rho<\rho_{0}<R_{2}$, it must be bounded on the closed sphere $|s-\chi_{0}|<\rho_{0}$, according to Sec.~(\ref{sec:chap_2_7_17}). That is, there is a nonnegative constant $M$ such that
\[
|g(s)|\leq M\texttt{ whenever }|s-\chi_{0}|\leq\rho_{0},
\]
and, since the area $A$ of the surface $S_{\rho}$ is $A=2\pi\rho^{2}$, it follows that
\[
|\iint_{S_{\rho}}g(s)d\sigma|\leq2\pi\rho^{2}M.
\]
Consequently,
\begin{equation}\label{eq:chap_14_5_3}
\lim_{\rho \to 0}\iint_{S_{\rho}}g(s)d\sigma=0.
\end{equation}

Because for the hemisphere $-S_{\rho}$ $(0\leq\theta\leq\pi,\texttt{ }-\pi/2\leq\varphi\leq\pi/2)$ there is
\[
d\sigma=(s-\chi_{0})\cos\varphi d\theta(s-\chi_{0})d\varphi
=(s-\chi_{0})^{2}\cos\varphi d\varphi d\theta,
\]
the second integral on the right in equation~(\ref{eq:chap_14_5_2}) has the value
\[
\iint_{S_{\rho}}\frac{d\sigma}{(s-\chi_{0})^{2}}=-\iint_{-S_{\rho}}\frac{d\sigma}{(s-\chi_{0})^{2}}
=-\int_{-\pi/2}^{\pi/2}\cos\varphi d\varphi\int_{0}^{\pi}d\theta=-2\pi.
\]
Thus
\begin{equation}\label{eq:chap_14_5_4}
\lim_{\rho \to 0}\iint_{S_{\rho}}\frac{d\sigma}{(s-\chi_{0})^{2}}=-2\pi.
\end{equation}

Limit~(\ref{eq:chap_14_5_1}) now follows by letting $\rho_{0}$ tend to zero on each side of equation~(\ref{eq:chap_14_5_2}) and referring to limits~(\ref{eq:chap_14_5_3}) and~(\ref{eq:chap_14_5_4}).

%\section{An Indentation Around a Branch Point}\label{sec:chap_14_6_76}

%\section{Integration Along a Branch Cut}\label{sec:chap_14_7_77}

%\section{Definite Integrals involving Sines and Cosines}\label{sec:chap_14_8_78}

\section{Argument Principle}\label{sec:chap_14_9_79}

A function $f$ is said to be $meromorphic$ in a domain $D$ if it is analytic throughout $D$ except for poles. Suppose now that $f$ is $meromorphic$ in the domain interior to a positively oriented simple closed spatial contour $C$ and that it is analytic and nonzero on $C$. The image $\Gamma$ of $C$ under the transformation $\varpi=f(s)$ is a closed spatial contour, not necessarily simple, in the $\varpi$ space. As a point $s$ traverses $C$ in the positive direction, its images $\varpi$ traverses $\Gamma$ in a particular direction that determines the orientation of $\Gamma$. Note that, since $f$ has no zeros on $C$, the spatial contour $\Gamma$ does not pass through the origin in the $\varpi$ space.

Let $\varpi$ and $\varpi_{0}$ be points on $\Gamma$, where $\varpi_{0}$ is fixed and $\phi_{0}$ is a value of $\arg_{c}\varpi_{0}$. Then let $\arg_{c}\varpi$ vary continuously, starting with the value $\phi_{0}$, as the point  $\varpi$ begins at the point $\varpi_{0}$ and traverses $\Gamma$ once in the direction of orientation assigned to it by the mapping $\varpi=f(s)$. When $\varpi$ returns to the point $\varpi_{0}$, where it started, $\arg_{c}\varpi$ assumes a particular value of $\arg_{c}\varpi_{0}$, which we denote by $\phi_{1}$. Thus the change  in $\arg_{c}\varpi$ as $\varpi$ describes $\Gamma$ once in its direction of orientation is $\phi_{1}-\phi_{0}$.  This change is, of course, independent of the point $\varpi_{0}$ chosen to determine it. Since $\varpi=f(s)$, the  number $\phi_{1}-\phi_{0}$ is, in fact, the change in argument of $f(s)$ as $s$ describes $C$ once in the positive direction, starting with a point $s_{0}$; and we write
\[
\Delta_{C}\arg_{c}f(s)=\phi_{1}-\phi_{0}.
\]
The value of $\Delta_{C}\arg_{c}f(s)$ is evidently an integral multiple of $2\pi$, and the integer
\[
\frac{1}{2\pi}\Delta_{C}\arg_{c}f(s)
\]
represents the number of times the point $\varpi$ winds around the origin in the $\varpi$ space. For that reason, this integer is sometimes called the winding number of $\Gamma$ with respect to the origin $\varpi=0$. It is positive if $\Gamma$ winds around the origin in the counterclockwise direction and negative if it winds clockwise around that point. The winding number is always zero when $\Gamma$ does not enclose the origin. %The verification  of this fact for a special case is left to the exercises.

The winding number can be determined from the number of zeros and poles of $f$ interior to $C$. The number of  poles is necessarily finite. %, according to Exercise 11, Sec. 69.
Likewise, with the understanding that $f(s)$ is not identically equal to zero everywhere else inside $C$, it is easily shown %(Exercise 4, Sec. 80)
that the zeros of $f$ are finite in number and are all of finite order. Suppose now that $f$ has $Z$ zeros and $P$
poles in the domain interior to $C$.  We agree that $f$ has $m_{0}$ zeros at a point $s_{0}$ if it has a zero of  order $m_{0}$ there; and if $f$ has a pole of order $m_{p}$ at $s_{0}$, that pole is to be counted
$m_{p}$ times. The following theorem, which is known as the argument principle, states that the winding number is simply the difference $Z-P$.

\begin{theorem}\label{th:chap_14_9_l}
Suppose that

(i) a function $f(s)$ is $meromorphic$ in the domain interior to a positively oriented simple closed spatial contour $C$;

(ii) $f(s)$ is analytic and nonzero on $C$;

(iii) counting multiplicities, $Z$ is the number of zeros and $P$ is the number of poles of $f(s)$ inside $C$.

Then
\begin{equation}\label{eq:chap_14_9_l}
\frac{1}{2\pi}\Delta_{C}\arg_{c}f(s)=Z-P.
\end{equation}
\end{theorem}

To prove this, we evaluate the integral of $f'(s)/f(s)$ around $C$ in two different ways. First, we let $s=s(t)$ $(a\leq t\leq b)$ be a parametric representation for $C$, so that
\begin{equation}\label{eq:chap_14_9_2}
\int_{C}\frac{f'(s)}{f(s)}ds=\int_{a}^{b}\frac{f'[s(t)]s'(t)}{f[s(t)]}dt.
\end{equation}
Since, under the transformation $\varpi=f(s)$, the image $\Gamma$ of $C$ never passes through
the origin in the $\varpi$ space, the image of any point $s=s(t)$ on $C$ can be expressed in
exponential form as $\varpi=\rho(t)(e_{xy}\exp[i\phi(t)]\cos\psi(t)+\sin\psi(t))$. Thus
\begin{equation}\label{eq:chap_14_9_3}
f[s(t)]=\rho(t)(e_{xy}e^{i\phi(t)}\cos\psi(t)+\sin\psi(t))\texttt{ }(a\leq t\leq b);
\end{equation}
and, along each of the smooth arcs making up the spatial contour $\Gamma$. it follows that% (see Exercise 5, Sec. 38)
\begin{equation}\label{eq:chap_14_9_4}
\begin{split}
f'[s(t)]s'(t) &= \frac{d}{dt}f[s(t)]=\frac{d}{dt}\{\rho(t)[e_{xy}e^{i\phi(t)}\cos\psi(t)+\sin\psi(t)]\} \\
&= \rho(t)[\frac{\rho'(t)}{\rho(t)}+i\phi'(t)][e_{xy}e^{i\phi(t)}\cos\psi(t)+\sin\psi(t)] \\
&- i\rho(t)\sin\psi(t)\phi'(t)-\rho(t)[e_{xy}e^{i\phi(t)}\sin\psi(t)-\cos\psi(t)]\psi'(t).
\end{split}
\end{equation}
Inasmuch as $\rho'(t)$ and $\phi'(t)$ are piecewise continuous on the interval $a\leq t\leq b$, we can now use expressions~(\ref{eq:chap_14_9_3}) and~(\ref{eq:chap_14_9_4}) to write integral~(\ref{eq:chap_14_9_2}) as follows:
\[
\begin{split}
\int_{C}\frac{f'(s)}{f(s)}ds &= \int_{a}^{b}\frac{\rho'(t)}{\rho(t)}dt+i\int_{a}^{b}\phi'(t)dt \\
&+e_{xy}\int_{a}^{b}\frac{-ie^{-i\phi(t)}\sin\psi(t)\phi'(t)}{\cos\psi(t)+e^{-i\phi(t)}\sin\psi(t)}dt \\
&+\int_{a}^{b}\frac{-e_{xy}e^{i\phi(t)}\sin\psi(t)+\cos\psi(t)}{e_{xy}e^{i\phi(t)}\cos\psi(t)+\sin\psi(t)}\psi'(t)dt\\
&=\ln\rho(t)|_{a}^{b}+i\phi(t)|_{a}^{b}+e_{xy}\ln[\cos\psi(t)+e^{-i\phi(t)}\sin\psi(t)]|_{a}^{b} \\
&+\ln[e_{xy}e^{i\phi(t)}\cos\psi(t)+\sin\psi(t)]|_{a}^{b}.
\end{split}
\]
But
\[
\rho(b)=\rho(a),\texttt{ }\phi(b)-\phi(a)=\Delta_{C}\arg_{c}f(s),
\]
\[
\cos\psi(b)+e^{-i\phi(b)}\sin\psi(b)=\cos\psi(a)+e^{-i\phi(a)}\sin\psi(a),
\]
\[
e_{xy}e^{i\phi(b)}\cos\psi(b)+\sin\psi(b)=e_{xy}e^{i\phi(a)}\cos\psi(a)+\sin\psi(a).
\]
Hence
\begin{equation}\label{eq:chap_14_9_5}
\int_{C}\frac{f'(s)}{f(s)}ds=i\Delta_{C}\arg_{c}f(s).
\end{equation}

Another way to evaluate integral~(\ref{eq:chap_14_9_5}) is to use Cauchy's residue theorem. To be specific, we observe that the integrand $f'(s)/f(s)$ is analytic inside and on $C$ except at the points inside $C$ at which the zeros and poles of $f$ occur. If $f$ has a zero of order $m_{0}$ at $s_{0}$, then (Sec.~(\ref{sec:chap_7_7_68}))
\begin{equation}\label{eq:chap_14_9_6}
f(s)=(s-s_{0})^{m_{0}}g(s),
\end{equation}
where $g(s)$ is analytic and nonzero at $s_{0}$. Hence
\[
f'(s)=m_{0}(s-s_{0})^{m_{0}-1}g(s)+(s-s_{0})^{m_{0}}g'(s),
\]
or
\begin{equation}\label{eq:chap_14_9_7}
\frac{f'(s)}{f(s)}=\frac{m_{0}}{s-s_{0}}+\frac{g'(s)}{g(s)}.
\end{equation}
Since $g'(s)/g(s)$ is analytic at $s_{0}$, it has a Taylor series representation about that point; and so equation~(\ref{eq:chap_14_9_7}) tells us that $f'(s)/f(s)$ has a simple pole at $s_{0}$, with residue $m_{0}$. If, on the other hand, $f$ has a pole of order $m_{p}$ at $s_{0}$, we know from the theorem in
Sec.~(\ref{sec:chap_7_5_66}) that
\begin{equation}\label{eq:chap_14_9_8}
f(s)=(s-s_{0})^{-m_{p}}\phi(s),
\end{equation}
where $\phi(s)$ is analytic and nonzero at $s_{0}$. Because expression~(\ref{eq:chap_14_9_8}) has the same form
as expression~(\ref{eq:chap_14_9_6}), with the positive integer $m_{0}$ in equation~(\ref{eq:chap_14_9_6}) replaced by $-m_{p}$, it is clear from equation~(\ref{eq:chap_14_9_7}) that $f'(s)/f(s)$ has a simple pole at $s_{0}$, with residue $-m_{p}$. Applying the residue theorem, then, we find that
\begin{equation}\label{eq:chap_14_9_9}
\int_{C}\frac{f'(s)}{f(s)}ds=2\pi i(Z-P).
\end{equation}
Expression~(\ref{eq:chap_14_9_l}) now follows by equating the right-hand sides of equations~(\ref{eq:chap_14_9_5}) and~(\ref{eq:chap_14_9_9}).

Equation~(\ref{eq:chap_14_9_5}) is of an important property. Let two functions $f_{1}(s)$ and $f_{2}(s)$ satisfy equation~(\ref{eq:chap_14_9_5}). Then there are
\[
\int_{C}\frac{f_{1}'(s)}{f_{1}(s)}ds=i\Delta_{C}\arg_{c}f_{1}(s)\texttt{ and }\int_{C}\frac{f_{2}'(s)}{f_{2}(s)}ds=i\Delta_{C}\arg_{c}f_{2}(s).
\]
Now we let $F(s)$ denote the product of $f_{1}(s)$ and $f_{2}(s)$, that is $F(s)=f_{1}(s)f_{2}(s)$. Then, from equation~(\ref{eq:chap_14_9_5}), we have
\[
\int_{C}\frac{F'(s)}{F(s)}ds=i\Delta_{C}\arg_{c}F(s)=i\Delta_{C}\arg_{c}[f_{1}(s)f_{2}(s)].
\]
On the other hand,
\[
\begin{split}
\int_{C}\frac{F'(s)}{F(s)}ds &= \int_{C}\frac{f_{1}'(s)}{f_{1}(s)}ds+\int_{C}\frac{f_{2}'(s)}{f_{2}(s)}ds \\
                             &=i\Delta_{C}\arg_{c}f_{1}(s)+i\Delta_{C}\arg_{c}f_{2}(s).
\end{split}
\]
Hence we get
\begin{equation}\label{eq:chap_14_9_10}
\Delta_{C}\arg_{c}[f_{1}(s)f_{2}(s)]=\Delta_{C}\arg_{c}f_{1}(s)+\Delta_{C}\arg_{c}f_{2}(s).
\end{equation}

\section{Rouche's Theorem}\label{sec:chap_14_10_80}

The main result in this section is known as Rouchi's theorem and is a consequence of the argument principle, just developed in Sec.~(\ref{sec:chap_14_9_79}). It can be useful in locating regions of the complex plane in which a given analytic function has zeros.

\begin{theorem}\label{th:chap_14_10_l}
Suppose that

(i) two functions $f(s)$ and $g(s)$ are analytic inside and on a simple closed spatial contour $C$;

(ii) $|f(s)|>|g(s)|$ at each point on $C$.

Then $f(s)$ and $f(s)+g(s)$ have the same number of zeros, counting multiplicities, inside $C$.
\end{theorem}

The orientation of $C$ in the statement of the theorem is evidently immaterial. Thus, in the proof here, we may assume that the orientation is positive. We begin with the observation that neither the function $f(s)$ nor the sum  $f(s)+g(z)$ has a zero on $C$, since
\[
|f(s)|>|g(s)|\geq0\texttt{ and }|f(s)+g(s)|\geq||f(s)|-|g(s)||>0
\]
when $s$ is on $C$.

If $Z_{f}$ and $Z_{f+g}$ denote the number of zeros, counting multiplicities, of $f(s)$ and $f(s)+g(s)$, respectively, inside $C$, we know from the theorem in Sec.~(\ref{sec:chap_14_9_79}) that
\[
Z_{f}=\frac{1}{2\pi}\Delta_{C}\arg_{c}f(s)\texttt{ and }Z_{f+g}=\frac{1}{2\pi}\Delta_{C}\arg_{c}[f(s)+g(s)].
\]
Consequently, since from equation~(\ref{eq:chap_14_9_10}) in Sec.~(\ref{sec:chap_14_9_79})
\[
\Delta_{C}\arg_{c}[f_{1}(s)f_{2}(s)]=\Delta_{C}\arg_{c}f_{1}(s)+\Delta_{C}\arg_{c}f_{2}(s)
\]
such that
\[
\begin{split}
\Delta_{C}\arg_{c}[f(s)+g(s)] &= \Delta_{C}\arg_{c}\{f(s)[1+\frac{g(s)}{f(s)}]\} \\
                              &= \Delta_{C}\arg_{c}f(s)+\Delta_{C}\arg_{c}[1+\frac{g(s)}{f(s)}],
\end{split}
\]
it is clear that
\begin{equation}\label{eq:chap_14_10_l}
Z_{f+g}=Z_{f}+\frac{1}{2\pi}\Delta_{C}\arg_{c}F(s),
\end{equation}
where
\[
F(s)=1+\frac{g(s)}{f(s)}.
\]
But
\[
|F(s)-1|=\frac{|g(s)|}{|f(s)|}<1;
\]
and this means that, under the transformation $\varpi=F(s)$, the image of $C$ lies in the open sphere $|\varpi-1|<1$. That image does not, then, enclose the origin $\varpi=0$. Hence $\Delta_{C}\arg_{c}F(s)=0$ and, since equation~(\ref{eq:chap_14_10_l}) reduces to $Z_{f+g}=Z_{f}$, the theorem here is proved.

\section{Inverse Laplace Transforms}\label{sec:chap_14_11_81}

\subsection{Inverse Laplace Transforms with Contour Integrals}\label{sec:chap_14_11_81_1}

Suppose that a function $F$ of the spatial complex variable
\[
s=r(e_{xy}e^{i\theta}\cos\varphi+\sin\varphi)
\]
is analytic throughout the finite $s$ space except for a finite number of isolated singularities which are in a spatial plane $P_{\gamma}$ determined by a point $\gamma=e_{xy}\alpha+\beta$ on the $xz$ coordinate plane and the $y$ axis or the argument $\varphi_{\gamma}=\arctan(\beta/\alpha)$ on the $xz$ coordinate plane and the $y$ axis. Then let $L_{R}$ denote a vertical line segment $(\alpha,y,\beta)$ perpendicular to the $xz$ plane from $s=\gamma-e_{xy}iR$ to $s=\gamma+e_{xy}iR$, where components of the constant $\gamma$ are positive and large enough that the singularities of $F$ all lie to the left of that segment on the $x$ axis and to the behind of that segment on the $z$ axis. A new function $f$ of the real variable $t$ is defined for positive values of $t$ by means of the equation
\begin{equation}\label{eq:chap_14_11_1}
f(t)=\frac{1}{2\pi i}\lim_{R \to \infty}\int_{L_{R}}e^{st}F(s)ds\texttt{ }(t>0),
\end{equation}
provided this limit exists. Expression~(\ref{eq:chap_14_11_1}) is usually written
\begin{equation}\label{eq:chap_14_11_2}
f(t)=\frac{1}{2\pi i}P.V.\int_{\gamma-e_{xy}i\infty}^{\gamma+e_{xy}i\infty}e^{st}F(s)ds\texttt{ }(t>0)
\end{equation}
[compare equation~(\ref{eq:chap_14_1_3}), Sec.~(\ref{sec:chap_14_1_71})], and such an integral is called a Bromwich integral.

It can be shown that, when fairly general conditions are imposed on the functions involved, $f(t)$ is the inverse Laplace transform of $F(s)$. That is, if $F(s)$ is the Laplace transform of $f(t)$, defined by the equation
\begin{equation}\label{eq:chap_14_11_3}
F(s)=\int_{0}^{\infty}e^{-st}f(t)dt,
\end{equation}
then $f(t)$ is retrieved by means of equation~(\ref{eq:chap_14_11_2}), where the choice of the positive component
numbers of the constant $\gamma$ are immaterial as long as the singularities of $F$ all lie to the left of $L_{R}$ on the $x$ axis and to the behind of $L_{R}$ on the $z$ axis. Laplace transforms and their inverses are important in solving both ordinary and partial differential equations.

Residues can often be used to evaluate the limit in expression~(\ref{eq:chap_14_11_1}) when the function $F(s)$ is specified. To see how this is done, we let $s_{n}$ $(n=1, 2, \ldots, N)$ denote the singularities of $F(s)$. We then let $R_{0}$ denote the largest of their moduli and consider a semicircle $C_{R}$ in the spatial plane $P_{\gamma}$ with parametric representation
\begin{equation}\label{eq:chap_14_11_4}
s=\gamma+R(e_{xy}e^{i\theta}\cos\varphi+\sin\varphi)\texttt{ }(\frac{\pi}{2}\leq\theta\leq\frac{3\pi}{2},\texttt{ }-\frac{\pi}{2}\leq\varphi\leq0,
\end{equation}
or
\[
\vec{s}=|\gamma|+Re^{i\theta}\texttt{ }(\frac{\pi}{2}\leq\theta\leq\frac{3\pi}{2}),
\]
where $R>R_{0}+|\gamma|$, and $\vec{s}$ is the transformed value of $s$ in the new transformed $\vec{x}y\vec{z}$ coordinate system with $\vec{\varphi}=0$, which is similar to the case of the two-dimensional complex variables. Note that, for each $s_{n}$,
\[
|s_{n}-\gamma|\leq|s_{n}|+|\gamma|\leq R_{0}+|\gamma|<R.
\]
Hence the singularities all lie in the interior of the semicircular region bounded by $C_{R}$ and $L_{R}$, and Cauchy's residue theorem~(\ref{th:chap_7_2_1}) in Sec.~(\ref{sec:chap_7_2_63_1}) tells us that
\begin{equation}\label{eq:chap_14_11_5}
\int_{L_{R}}e^{st}F(s)ds=2\pi i\sum_{n=1}^{N}Res_{s=s_{n}}[e^{st}F(s)]-\int_{C_{R}}e^{st}F(s)ds.
\end{equation}

Suppose now that, for all points s on $C_{R}$, there is a positive constant $M_{R}$ such that
$|F(s)|\leq M_{R}$ where $M_{R}$ tends to zero as $R$ tends to infinity. We may use the parametric
representation~(\ref{eq:chap_14_11_4}) for $C_{R}$ to write
\[
\int_{C_{R}}e^{st}F(s)ds=e_{xy}\int_{\pi/2}^{3\pi/2}
\exp[\gamma t+Rt(e^{i\theta}\cos\varphi+\sin\varphi)]F(s)Rie^{i\theta}\cos\varphi d\theta.
\]
Then, since $-\pi/2\leq\varphi\leq0$ such that $\sin\varphi\leq0$,
\[
|\exp[\gamma t+Rt(e^{i\theta}\cos\varphi+\sin\varphi)]|=e^{|\gamma|t-Rt|\sin\varphi|}e^{Rt\cos\varphi\cos\theta}
\]
and
\[
|F[\gamma+R(e_{xy}e^{i\theta}\cos\varphi+\sin\varphi)]|\leq M_{R},
\]
we find that
\begin{equation}\label{eq:chap_14_11_6}
|\int_{C_{R}}e^{st}F(s)ds|\leq\frac{e^{|\gamma|t}M_{R}R}{e^{Rt|\sin\varphi|}}\int_{\pi/2}^{3\pi/2}
e^{Rt\cos\varphi\cos\theta}\cos\varphi d\theta.
\end{equation}
But the substitution $\phi=\theta-(\pi/2)$, together with Jordan's inequality~(\ref{eq:chap_14_4_2}), Sec.~(\ref{sec:chap_14_4_74}), reveals that
\[
\int_{\pi/2}^{3\pi/2}
e^{Rt\cos\varphi\cos\theta}\cos\varphi d\theta=\int_{0}^{\pi}e^{-Rt\cos\varphi\sin\theta}\cos\varphi d\phi<\frac{\pi}{Rt}.
\]
Inequality~(\ref{eq:chap_14_11_6}) thus becomes
\begin{equation}\label{eq:chap_14_11_7}
|\int_{C_{R}}e^{st}F(s)ds|\leq\frac{e^{|\gamma|t}M_{R}\pi}{e^{Rt|\sin\varphi|}t},
\end{equation}
and this shows that
\begin{equation}\label{eq:chap_14_11_8}
\lim_{R \to \infty}\int_{C_{R}}e^{st}F(s)ds=0.
\end{equation}
Letting $R$ tend to $\infty$ in equation~(\ref{eq:chap_14_11_5}), then, we see that the function $f(t)$, defined by equation~(\ref{eq:chap_14_11_1}), exists and that it can be written
\begin{equation}\label{eq:chap_14_11_9}
f(t)=\sum_{n=1}^{N}Res_{s=s_{n}}[e^{st}F(s)]\texttt{ }(t>0).
\end{equation}

In many applications of Laplace transforms, such as the solution of partial differential equations arising in studies of heat conduction and mechanical vibrations, the function $F(s)$ is analytic for all values of $s$ in the finite spatial plane $P_{\gamma}$ except for an infinite set of isolated singular points $s_{n}$ $(n=1, 2, \ldots, N)$ that lie to the left of some vertical line $(Re$ $s=\alpha,y,Re_{s}s=\beta)$. Often the method just described for finding $f(t)$ can then be modified in such a way that the finite sum~(\ref{eq:chap_14_11_9}) is replaced by an  infinite series of residues:
\begin{equation}\label{eq:chap_14_11_10}
f(t)=\sum_{n=1}^{\infty}Res_{s=s_{n}}[e^{st}F(s)]\texttt{ }(t>0).
\end{equation}

The basic modification is to replace the vertical line segments $L_{R}$ by vertical line segments $L_{N}$ $(N=1, 2, \ldots)$ from $s=\gamma-ib_{N}$ to $s=\gamma+ib_{N}$. The circular arcs $C_{R}$ are then replaced by contours $C_{N}$ $(N=1, 2, \ldots)$ from $\gamma+ib_{N}$ to $\gamma-ib_{N}$ such that, for each $N$, the sum $L_{N}+C_{N}$ is a simple closed contour enclosing the singular points $s_{1}, s_{2}, \ldots, s_{N}$. Once it is shown that
\begin{equation}\label{eq:chap_14_11_11}
\lim_{N \to \infty}\int_{C_{N}}e^{st}F(s)ds=0,
\end{equation}
expression~(\ref{eq:chap_14_11_2}) for $f(t)$ becomes expression~(\ref{eq:chap_14_11_10}).

The choice of the contours $C_{N}$ depends on the nature of the function $F(s)$. Common choices include circular or parabolic arcs and rectangular paths. Also, the simple closed contour $L_{N}+C_{N}$ need not enclose precisely $N$ singularities. When, for example, the region between $L_{N}+C_{N}$ and $L_{N+l}+C_{N+l}$ contains two singular points of $F(s)$, the pair of corresponding residues of $e^{st}F(s)$ are simply grouped together as a single term in series~(\ref{eq:chap_14_11_10}).% Since it is often quite tedious to establish limit~(\ref{eq:chap_14_11_11}) in any case, we shall accept it in the examples and related exercises below that involve an infinite number  of  singularities. Thus our use of expression~(\ref{eq:chap_14_11_10}) will be only formal.

\subsection{Inverse Laplace Transforms with Surface Integrals}\label{sec:chap_14_11_81_2}

Suppose that a function $F$ of the spatial complex variable $s=r(e_{xy}e^{i\theta}\cos\varphi+\sin\varphi)$ is analytic throughout the finite $s$ space except for a finite number of isolated singularities. Then let a point $\gamma=e_{xy}\alpha+\beta$ on the $xz$ coordinate plane, and

(i) let $X_{R}$ denote a half vertical closed circle $|s_{x}-\beta|\leq R$ parallel to the $yz$ coordinate plane for
\[
s_{x}=\gamma+r(e_{xy}i\sin\theta-\cos\theta)\texttt{ where }0\leq r\leq R\texttt{ and }\pi/2\leq\theta\leq3\pi/2;
\]

(ii) let $Z_{R}$ denote a half vertical closed circle $|s_{z}-\alpha|\leq R$ parallel to the $xy$ coordinate plane for
\[
s_{z}=\gamma+re_{xy}(-\cos\theta+i\sin\theta)\texttt{ where }0\leq r\leq R\texttt{ and }\pi/2\leq\theta\leq3\pi/2;
\]

(iii) let $S_{R}$ denote a quarter of a spherical surface $|s-\gamma|=R$ where components of the constant $\gamma$ are positive and large enough that the singularities of $F$ all lie to the left of the plane $X_{R}$ and to the behind of the plane $Z_{R}$.

A new function $f$ of the real variable $t$ is defined for positive values of $t$ by means of the equation
\begin{equation}\label{eq:chap_14_12_1}
\begin{split}
f(t) &= \frac{1}{4\pi}\lim_{R \to \infty}[\iint_{X_{R}}e^{st}F(s)d\sigma+\iint_{Z_{R}}e^{st}F(s)d\sigma] \\
     &= \frac{1}{4\pi}\lim_{R \to \infty}\int_{0}^{R}rdr\int_{\pi/2}^{3\pi/2}[e^{s_{x}t}F(s_{x})+e^{s_{z}t}F(s_{z})]d\theta\texttt{ }(t>0);
\end{split}
\end{equation}
provided these limits exist where $d\sigma=rdrd\theta$. Expression~(\ref{eq:chap_14_12_1}) is usually written
\begin{equation}\label{eq:chap_14_12_2}
f(t)=\frac{1}{4\pi}P.V.\int_{0}^{\infty}rdr\int_{\pi/2}^{3\pi/2}[e^{s_{x}t}F(s_{x})+e^{s_{z}t}F(s_{z})]d\theta\texttt{ }(t>0)
\end{equation}
[compare equation~(\ref{eq:chap_14_1_3}), Sec.~(\ref{sec:chap_14_1_71})], and such an integral is called a Bromwich integral.

It can be shown that, when fairly general conditions are imposed on the functions involved, $f(t)$ is the inverse Laplace transform of $F(s)$. That is, if $F(s)$ is the Laplace transform of $f(t)$, defined by the equation
\begin{equation}\label{eq:chap_14_12_3}
F(s)=\int_{0}^{\infty}e^{-st}f(t)dt\texttt{ where }
    \begin{array}{cc}
        F(s)=F(s_{x}) & \texttt{ when }s=s_{x} \\
        F(s)=F(s_{z}) & \texttt{ when }s=s_{z}
    \end{array},
\end{equation}
then $f(t)$ is retrieved by means of equation~(\ref{eq:chap_14_12_2}), where the choice of the positive component
numbers of the constant $\gamma$ are immaterial as long as the singularities of $F$ all lie to the left of $X_{R}$ and to the behind of $Z_{R}$. Laplace transforms and their inverses are important in solving both ordinary and partial differential equations.

Residues can often be used to evaluate the limit in expression~(\ref{eq:chap_14_12_1}) when the function $F(s)$ is specified. To see how this is done, we let $s_{n}$ $(n=1, 2, \ldots, N)$ denote the singularities of $F(s)$. We then let $R_{0}$ denote the largest of their moduli and consider $S_{R}$ with parametric representation
\begin{equation}\label{eq:chap_14_12_4}
s=\gamma+R(e_{xy}e^{i\theta}\cos\varphi+\sin\varphi)\texttt{ }(\frac{\pi}{2}\leq\theta\leq\frac{3\pi}{2},\texttt{ }-\frac{\pi}{2}\leq\varphi\leq0)
\end{equation}
where $R>R_{0}+|\gamma|$. Note that, for each $s_{n}$,
\[
|s_{n}-\gamma|\leq|s_{n}|+|\gamma|\leq R_{0}+|\gamma|<R.
\]
Hence the singularities all lie in the interior of a quarter of a spherical surface $|s-\gamma|=R$ bounded by $S_{R}$, $X_{R}$, and $Z_{R}$, and Cauchy's residue theorem~(\ref{th:chap_7_2_2_1}) in Sec.~(\ref{sec:chap_7_2_63_2}) tells us that
\begin{equation}\label{eq:chap_14_12_5}
\iint_{X_{R}}e^{s_{x}t}F(s_{x})d\sigma+\iint_{Z_{R}}e^{s_{z}t}F(s_{z})d\sigma
\end{equation}
\[
=4\pi\sum_{n=1}^{N}Res_{s=s_{n}}[e^{st}F(s)]-\iint_{S_{R}}e^{st}F(s)d\sigma.
\]

Suppose now that, for all points s on $S_{R}$, there is a positive constant $M_{R}$ such that
$|F(s)|\leq M_{R}/R$ where $M_{R}$ tends to zero as $R$ tends to infinity. We may use the parametric
representation~(\ref{eq:chap_14_12_4}) for $S_{R}$ with $d\sigma=(s-\gamma)^{2}\cos\varphi d\varphi d\theta$ to write
\[
\iint_{S_{R}}e^{st}F(s)d\sigma=\int_{-\pi/2}^{0}\cos\varphi d\varphi\int_{\pi/2}^{3\pi/2}
\exp[\gamma t+Rt(e_{xy}e^{i\theta}\cos\varphi+\sin\varphi)]
\]
\[
\times F[\gamma+R(e_{xy}e^{i\theta}\cos\varphi+\sin\varphi)](s-\gamma)^{2}d\theta.
\]
Then, since
\[
|\exp[\gamma t+Rt(e_{xy}e^{i\theta}\cos\varphi+\sin\varphi)]|=e_{xy}e^{|\gamma|t}e^{Rt(\cos\varphi\cos\theta+\sin\varphi)}
\]
and
\[
|F[\gamma+R(e_{xy}e^{i\theta}\cos\varphi+\sin\varphi)](s-\gamma)^{2}|\leq M_{R}R,
\]
we find that
\begin{equation}\label{eq:chap_14_12_6}
|\iint_{S_{R}}e^{st}F(s)d\sigma|\leq e^{|\gamma|t}M_{R}R\int_{-\pi/2}^{0}e^{Rt\sin\varphi}d\varphi\int_{\pi/2}^{3\pi/2}
e^{Rt\cos\varphi\cos\theta}\cos\varphi d\theta.
\end{equation}
But the substitution $\phi=\theta-(\pi/2)$, together with Jordan's inequality~(\ref{eq:chap_14_4_2}), Sec.~(\ref{sec:chap_14_4_74}), reveals that
\[
\int_{\pi/2}^{3\pi/2}e^{Rt\cos\varphi\cos\theta}\cos\varphi d\theta=\int_{0}^{\pi}e^{-Rt\cos\varphi\sin\phi}\cos\varphi d\phi<\frac{\pi}{Rt}
\]
and replacing $\varphi$ by $-\varphi$ in the first integral below, then
\[
\int_{-\pi/2}^{0}e^{Rt\sin\varphi}d\varphi\int_{\pi/2}^{3\pi/2}e^{Rt\cos\varphi\cos\theta}\cos\varphi d\theta<\frac{\pi}{Rt}\int_{0}^{\pi/2}e^{-Rt\sin\varphi}d\varphi<\frac{\pi^{2}}{2Rt}.
\]
Inequality~(\ref{eq:chap_14_12_6}) thus becomes
\begin{equation}\label{eq:chap_14_12_7}
|\iint_{S_{R}}e^{st}F(s)d\sigma|\leq\frac{e^{|\gamma|t}M_{R}\pi^{2}}{2t},
\end{equation}
and this shows that
\begin{equation}\label{eq:chap_14_12_8}
\lim_{R \to \infty}\iint_{S_{R}}e^{st}F(s)ds=0.
\end{equation}
Letting $R$ tend to $\infty$ in equation~(\ref{eq:chap_14_12_5}), then, we see that the function $f(t)$, defined by equation~(\ref{eq:chap_14_12_1}), exists and that it can be written
\begin{equation}\label{eq:chap_14_12_9}
f(t)=\sum_{n=1}^{N}Res_{s=s_{n}}[e^{st}F(s)]\texttt{ }(t>0).
\end{equation}

In many applications of Laplace transforms, such as the solution of partial differential equations arising in studies of heat conduction and mechanical vibrations, the function $F(s)$ is analytic for all values of $s$ in the finite space bounded by $S_{R}$, $X_{R}$, and $Z_{R}$ except for an infinite set of isolated singular points $s_{n}$ $(n=1, 2, \ldots, N)$ that lie to the left of some plane $X_{R}$ and to the behind of some plane $Z_{R}$. Often the method just described for finding $f(t)$ can then be modified in such a way that the finite sum~(\ref{eq:chap_14_12_9}) is replaced by an  infinite series of residues:
\begin{equation}\label{eq:chap_14_12_10}
f(t)=\sum_{n=1}^{\infty}Res_{s=s_{n}}[e^{st}F(s)]\texttt{ }(t>0).
\end{equation}

The basic modification is to replace the vertical planes $X_{R}$ by vertical planes $X_{N}$ $(N=1, 2, \ldots)$ for
\[
s_{x}=\gamma+r_{N}(e_{xy}i\sin\theta-\cos\theta)\texttt{ where }0\leq r_{N}\leq R\texttt{ and }\pi/2\leq\theta\leq3\pi/2;
\]
and to replace the vertical planes $Z_{R}$ by vertical planes $Z_{N}$ $(N=1, 2, \ldots)$ for
\[
s_{z}=\gamma+r_{N}e_{xy}(-\cos\theta+i\sin\theta)\texttt{ where }0\leq r_{N}\leq R\texttt{ and }\pi/2\leq\theta\leq3\pi/2;
\]
The spherical surfaces $S_{R}$ are then replaced by spherical surfaces $S_{N}$ $(N=1, 2, \ldots)$ such that, for each $N$, the sum $X_{N}+Z_{N}+S_{N}$ is a simple closed surface enclosing the singular points $s_{1}, s_{2}, \ldots, s_{N}$. Once it is shown that
\begin{equation}\label{eq:chap_14_12_11}
\lim_{N \to \infty}\iint_{S_{N}}e^{st}F(s)ds=0,
\end{equation}
expression~(\ref{eq:chap_14_12_2}) for $f(t)$ becomes expression~(\ref{eq:chap_14_12_10}).

The choice of the spherical surfaces $S_{N}$ depends on the nature of the function $F(s)$. Common choices include spherical or parabolic surfaces and rectangular planes. Also, the simple closed surface $X_{N}+Z_{N}+S_{N}$ need not enclose precisely $N$ singularities. When, for example, the region between $X_{N}+Z_{N}+S_{N}$ and $X_{N+1}+Z_{N+1}+S_{N+1}$ contains two singular points of $F(s)$, the pair of corresponding residues of $e^{st}F(s)$ are simply grouped together as a single term in series~(\ref{eq:chap_14_12_10}).% Since it is often quite tedious to establish limit~(\ref{eq:chap_14_12_11}) in any case, we shall accept it in the examples and related exercises below that involve an infinite number  of  singularities. Thus our use of expression~(\ref{eq:chap_14_12_10}) will be only formal.

%-----------------------------------------------------------------------
% Beginning of chap15.tex
%-----------------------------------------------------------------------
%
% AMS-LaTeX 1.2 sample file for a monograph, based on amsbook.cls.
% This is a data file input by chapter.tex.
%%%%%%%%%%%%%%%%%%%%%%%%%%%%%%%%%%%%%%%%%%%%%%%%%%%%%%%%%%%%%%%%%%%%

%\part{This is a Part Title Sample}

\chapter{Joukowski Transformation and Applications}\label{ch:chap_15}

In this chapter, the Joukowski transformation is used to establish the solutions for the flow
around ellipsoids and airfoils. The superposition of elementary solutions then leads to the solution for the flow around a sphere. The contents include Joukowski transformation, flows around an ellipsoid, Kutta condition and the flat-plate airfoils, symmetrical Joukowski airfoils, circular-arc airfoils, and Joukowski airfoils.

\section{Circulation and Vortex}\label{sec:chap_15_1_139}

\subsection{Circulation and Vorticity}\label{sec:chap_15_1_139_1}

The circulation contained within a closed contour in a body of fluid is defined as the integral around the contour of the component of the velocity vector that is locally tangent to the contour. That is, the circulation $\Gamma$ is defined as
\begin{equation}\label{eq:chap_15_1_1_2_3}
\Gamma=\oint\mathbf{u}\cdot d\mathbf{l}
\end{equation}
where $d\mathbf{l}$ represents an element of the contour. The integration is taken counterclockwise around the contour, and the circulation is positive if this
integral is positive.

The vorticity of an element of fluid is defined as the curl of its velocity vector. That is, the vorticity $\omega$ is defined by
\begin{equation}\label{eq:chap_15_1_2_2_4}
\mathbf{\omega}=\mathbf{\nabla}\times\mathbf{u}
\end{equation}
In tensor notation, Eqs.~(\ref{eq:chap_15_1_2_2_4}) may be written in the form
\[
\omega_{ijk}=-\varepsilon_{ijk}\frac{\partial\mathbf{u}_{j}}{\partial x_{k}}=(\frac{\partial\mathbf{u}_{k}}{\partial x_{j}}-\frac{\partial\mathbf{u}_{j}}{\partial x_{k}}).
\]
From this definition it is evident %, by comparison with Eq. (1.6a),
that the vorticity vector is numerically twice the angular speed of rotation of the fluid element about its own axes. That is, the vorticity is equal to twice the antisymmetric part of the deformation-rate tensor $e_{jk}$. It should be noted that a fluid element may travel on a circular streamline while having zero vorticity.

Vorticity is proportional to the angular velocity of a fluid element about its principal axes, not that of the center of gravity of the element about some reference point. Thus a particle traveling on a circular streamline will have no vorticity, provided that it does not revolve about its center of gravity as it moves.

The vorticity contained in a fluid element is related to the circulation around the element. This relationship may be obtained from an application of Stokes' theorem to the definition of circulation, as follows: From Eq.~(\ref{eq:chap_15_1_1_2_3}),
\[
\begin{split}
\Gamma&=\oint\mathbf{u}\cdot d\mathbf{l} \\
      &=\iint_{A}(\mathbf{\nabla}\times\mathbf{n})dA \\
\end{split}
\]
where the contour integral has been converted to a surface integral by use of Stokes' theorem, in which $A$ is the area defined by the closed contour around which the circulation is calculated and $\mathbf{n}$ is the unit normal to the surface.

Finally, invoking the definition of the vorticity vector, this relationship becomes
\begin{equation}\label{eq:chap_15_1_3_2_5}
\Gamma=\iint_{A}\mathbf{\omega}\times\mathbf{n}dA.
\end{equation}
Equation~(\ref{eq:chap_15_1_3_2_5}) shows that, for arbitrary choices of contours and enclosing
areas $A$, if $\mathbf{\omega}=0$, then $\Gamma=0$ and vice versa. Flows for which $\mathbf{\omega}=0$ are called irrotational, and flows for which this is not so are called rotational. The distinction between rotational and irrotational flow fields is an important one from the analytic point of view, as will be seen in later sections.

\subsection{Vortex}\label{sec:chap_15_1_139_2}

Now consider the constant of proportionality in the logarithmic complex potential to be imaginary. That is, consider
\[
F(s)=-ic\log s
\]
where $c$ is real and the minus is included to give a positive vortex. Then, using
spherical coordinates [See Eq.~(\ref{eq:chap_3_3_2b}), Sec.~(\ref{sec:chap_3_3_29})],
\[
\begin{split}
F(s)=&-ic\log[R(e_{xy}e^{i\theta}\cos\varphi+\sin\varphi)]=-ic\ln R \\
     &-\{
    \begin{array}{cc}
        ic\{\ln\sin\varphi+e_{xy}[\ln r_{L}(\theta,\varphi)+i\theta_{L}(\theta,\varphi)-\ln\sin\varphi]\} & \texttt{ when }\sin\varphi\neq0 \\
        ic(e_{xy}i\theta) & \texttt{ when }\sin\varphi=0
    \end{array} \\
    =&e_{xy}\{
    \begin{array}{cc}
        c\theta_{L}(\theta,\varphi)-ic[\ln R+\ln r_{L}(\theta,\varphi)] & \texttt{ when }\sin\varphi\neq0 \\
        c\theta-ic\ln R & \texttt{ when }\sin\varphi=0
    \end{array}
\end{split}
\]
for $0\leq\varphi\leq\pi/2$ and $\theta\neq(2n+1)\pi$ where $n=0,1,2,\ldots$
\[
\begin{split}
r_{L}(\theta,\varphi) &= \sqrt{1+\cos\theta\sin2\varphi},\texttt{ } \\
\theta_{L}(\theta,\varphi) &= \arctan\frac{\sin\theta}{\cos\theta+\tan\varphi}+2n\pi,
\end{split}
\]
Then, from the equation
\[
\begin{split}
F(s)&=e_{xy}[u(x,y,z)+iv(x,y,z)]+w(z) \\
    &=e_{xy}[u(R,\theta,\varphi)+iv(R,\theta,\varphi)]+w(R,\varphi),
\end{split}
\]
the velocity potential and the stream function are
\[
\begin{split}
u&=\{
    \begin{array}{cc}
        c\theta_{L}(\theta,\varphi) & \texttt{ when }\sin\varphi\neq0 \\
        c\theta & \texttt{ when }\sin\varphi=0
    \end{array},\\
v&=\{
    \begin{array}{cc}
        -c[\ln R+\ln r_{L}(\theta,\varphi)] & \texttt{ when }\sin\varphi\neq0 \\
        -c\ln R & \texttt{ when }\sin\varphi=0
    \end{array}.
\end{split}
\]
That is, the equipotential lines are the radial lines %$\theta_{L}(\theta,\varphi)=constant$ when $\sin\varphi\neq0$ or $\theta=constant$ when $\sin\varphi=0$,
\[
\begin{array}{cc}
\theta_{L}(\theta,\varphi)=constant & \texttt{ when }\sin\varphi\neq0 \\
\theta=constant & \texttt{ when }\sin\varphi=0
\end{array}
\]
and the streamlines are the big circles %$\ln R+\ln r_{L}(\theta,\varphi)=constant$ when $\sin\varphi\neq0$ or $\ln R=constant$ when $\sin\varphi=0$.
\[
\begin{array}{cc}
Rr_{L}(\theta,\varphi)=constant & \texttt{ when }\sin\varphi\neq0 \\
R=constant & \texttt{ when }\sin\varphi=0
\end{array}.
\]
The velocity components may be evaluated by use of the complex velocity $\mathbf{V}=\nabla u$,% i.e.,
\[
\begin{split}
u_{r}&=c\frac{\partial\theta_{L}(\theta,\varphi)}{\partial r}=0, \\
u_{\theta}&=c\frac{\partial\theta_{L}(\theta,\varphi)}{\partial\theta}=c\frac{1+\cos\theta\tan\varphi}{1+\cos\theta\sin2\varphi},\\
u_{\varphi}&=c\frac{\partial\theta_{L}(\theta,\varphi)}{\partial\varphi}=c\frac{-\sin\theta}{1+\cos\theta\sin2\varphi},
\end{split}
\]
for $\sin\varphi\neq0$.
\iffalse
or according to Sec.~(\ref{sec:chap_10_8_107})
\[
\begin{split}
F'(s)&=-i\frac{c}{s}=-i\frac{c}{R(e_{xy}e^{i\theta}\cos\varphi+\sin\varphi)} \\
     &=-\frac{c}{R}\cdot\frac{e_{xy}(0+i)+0}{e^{i\theta}\cos\varphi+\sin\varphi}.
\end{split}
\]
Comparison with the equation [Eq.~(\ref{eq:chap_2_12_8}), Sec.~(\ref{sec:chap_2_12_22})]
\[
\begin{split}
F'(s)&=\frac{e_{xy}(u_{r}+iv_{r})+w_{r}}{e^{i\theta}\cos\varphi+\sin\varphi}
\end{split}
\]
where for $r=R$
\[
rv_{r}=u_{\theta}(\cos\theta\tan\varphi-1)-(ru_{r}\cos\varphi-u_{\varphi}\sin\varphi)\sin\theta\sin\varphi
\]
shows that the velocity components are $u_{r}=0$ and
\[
u_{\theta}(1-\cos\theta\tan\varphi)-u_{\varphi}\sin\theta\sin^{2}\varphi=c
\]
or
\[
\cos\theta=-(\cos^{2}\varphi\pm1)\frac{\cot\varphi}{1+\cos^{2}\varphi}.
\]
\fi
Hence the velocity component $u_{\theta}$ is dependent on the angles $\theta$ and $\varphi$, the direction of the flow is positive (counterclockwise) for $c(1+\cos\theta\tan\varphi)>0$, and the resulting flow field is called a vortex.

A vortex is characterized by its strength, which may be measured by the circulation $\Gamma$ associated with it. From Eq.~(\ref{eq:chap_15_1_1_2_3}), the circulation $\Gamma$ associated with the singularity at the origin for a fixed angle $\varphi$ is
\[
\begin{split}
\Gamma&=\oint\mathbf{u}\cdot d\mathbf{l}=\int_{0}^{2\pi}u_{\theta}d\theta
       =\int_{0}^{2\pi}c\frac{1+\cos\theta\tan\varphi}{1+\cos\theta\sin2\varphi}d\theta \\
%      &=\int_{0}^{2\pi}c(1+\frac{\tan\varphi-\sin2\varphi}{1+\cos\theta\sin2\varphi}\cos\theta)d\theta \\
      &=c\int_{0}^{2\pi}d\theta+c(\tan\varphi-\sin2\varphi)\int_{0}^{2\pi}\frac{\cos\theta}{1+\cos\theta\sin2\varphi}d\theta\\
      &=2\pi c
       +4c(\tan\varphi-\sin2\varphi)\int_{0}^{\infty}\frac{1-t^{2}}{1+t^{2}+(1-t^{2})\sin2\varphi}\cdot\frac{dt}{1+t^{2}}\\
      &=2\pi c
       +4c\frac{\tan\varphi-\sin2\varphi}{\sin2\varphi}\int_{0}^{\infty}[\frac{1}{1+t^{2}}-\frac{1}{(1+\sin2\varphi)+(1-\sin2\varphi)t^{2}}]dt\\
      &=2\pi c
       +4c\frac{\tan\varphi-\sin2\varphi}{\sin2\varphi}[\arctan t+\frac{1-\sin2\varphi}{1+\sin2\varphi}\arctan(\frac{1-\sin2\varphi}{1+\sin2\varphi}t)]|_{0}^{\infty}\\
      &=2\pi c[1+2\frac{\tan\varphi-\sin2\varphi}{(1+\sin2\varphi)\sin2\varphi}]
\end{split}
\]
where let $t=\tan\theta/2$, then there are
\[
\cos\theta=\frac{1-t^{2}}{1+t^{2}},\texttt{ }d\theta=\frac{2dt}{1+t^{2}},\texttt{ }
\]
and
\[
\int_{0}^{2\pi}\frac{\cos\theta}{1+\cos\theta\sin2\varphi}d\theta=4\int_{0}^{\infty}\frac{1-t^{2}}{1+t^{2}+(1-t^{2})\sin2\varphi}\cdot\frac{dt}{1+t^{2}}.
\]
Here, the result $u_{\theta}=c(1+\cos\theta\tan\varphi)/(1+\cos\theta\sin2\varphi)$ has been used. Then $c$ may be replaced by $\Gamma=2\pi$, giving the following complex potential for a positive (counterclockwise)
vortex of strength $\Gamma$.
\[
F(s)=-i\frac{\Gamma}{2\pi}\log s.
\]
The singularity in this expression is located at $s=0$. That is, the line vortex is located at $s=0$. Then the complex potential for a positive vortex located at $s=s_{0}$ will be
\begin{equation}\label{eq:chap_15_1_4_2_9}
F(s)=-i\frac{\Gamma}{2\pi}\log(s-s_{0}).
\end{equation}
The complex potential for a negative vortex would be obtained by replacing $\Gamma$ by $-\Gamma$ in Eq.~(\ref{eq:chap_15_1_4_2_9}). Note, however, that the negative coefficient is associated with the positive vortex.

\section{Sphere without Circulation}\label{sec:chap_15_2_140}

The fundamental solutions to the foregoing flow situations provide the basis for more general solutions through the principle of superposition. Superposition is valid here, since the governing equation, for either the velocity potential or the stream function, is linear. The first example of superposition of fundamental solutions will be the flow around a sphere.

Consider the superposition of a uniform rectilinear flow and a doublet at the origin. Then, from equations $F(s)=Us$ and $F(s)=\mu/(s-s_{0})$, the complex potential for the resulting flow field will be
\[
F(s)=Us+\mu/s.
\]
It will now be shown that for a certain choice of the doublet strength the big circle $R=a$ becomes a streamline. On the big circle $R=a$, the value of $s$ is $a(e_{xy}e^{i\theta}\cos\varphi+\sin\varphi)$, so that the complex potential on this big circle is
\[
\begin{split}
F(s)=&Ua(e_{xy}e^{i\theta}\cos\varphi+\sin\varphi)+\frac{\mu}{a(e_{xy}e^{i\theta}\cos\varphi+\sin\varphi)}\\
    =&Ua(e_{xy}e^{i\theta}\cos\varphi+\sin\varphi)+\frac{\mu}{a}[e_{xy}(\frac{e^{-i\theta}\cos\varphi+\sin\varphi}{1+\cos\theta\sin2\varphi}-\frac{1}{\sin\varphi})+\frac{1}{\sin\varphi}]\\
    =&e_{xy}[Ua\cos\theta\cos\varphi+\frac{\mu}{a}(\frac{\cos\theta\cos\varphi+\sin\varphi}{1+\cos\theta\sin2\varphi}-\frac{1}{\sin\varphi})\\
     &+i(Ua\sin\theta\cos\varphi-\frac{\mu}{a}\cdot\frac{\sin\theta\cos\varphi}{1+\cos\theta\sin2\varphi})]+(Ua\sin\varphi+\frac{\mu}{a\sin\varphi})\\
    =&e_{xy}[(Ua+\frac{\mu/a}{1+\cos\theta\sin2\varphi})\cos\theta\cos\varphi+\frac{\mu}{a}\kappa_{x}(\theta,\varphi)\\
     &+i(Ua-\frac{\mu/a}{1+\cos\theta\sin2\varphi})\sin\theta\cos\varphi]+[Ua\sin\varphi+\frac{\mu}{a}\kappa_{z}(\varphi)]\\
    =&e_{xy}(u+iv)+w
\end{split}
\]
where if $\sin\varphi=0$ then $\kappa_{x}(\theta,\varphi)=0$ and $\kappa_{z}(\varphi)=0$ else
\[
\kappa_{x}(\theta,\varphi)=\frac{\sin\varphi}{1+\cos\theta\sin2\varphi}-\frac{1}{\sin\varphi}\texttt{ and }\kappa_{z}(\varphi)=\frac{1}{\sin\varphi},
\]
\[
\begin{split}
u&=(Ua+\frac{\mu/a}{1+\cos\theta\sin2\varphi})\cos\theta\cos\varphi+\frac{\mu}{a}\kappa_{x}(\theta,\varphi),\\
v&=(Ua-\frac{\mu/a}{1+\cos\theta\sin2\varphi})\sin\theta\cos\varphi, \\
w&=[Ua\sin\varphi+\frac{\mu}{a}\kappa_{z}(\varphi)].
\end{split}
\]
Thus the value of the stream function on the big circle $R=a$ is
\[
v=(Ua-\frac{\mu/a}{1+\cos\theta\sin2\varphi})\sin\theta\cos\varphi.
\]
For general values of $\mu$; $v$ is clearly variable, but if we choose the strength of the doublet to be $\mu=Ua^{2}$, then
\[
v=Ua\frac{\cos\theta\sin2\varphi}{1+\cos\theta\sin2\varphi}\sin\theta\cos\varphi
\]
on $R=a$. %The flow pattern for this doublet strength is shown in Fig. 4.8a.
The flow field due to the doublet encounters that due to the uniform flow and is bent downstream. %For clarity, the flow due to the doublet is shown dotted in Fig. 4.8a.
It may be seen that the doublet flow is entirely contained within the big circle $R=a$, while the uniform flow is deflected by the doublet in such a way that it is entirely outside the big circle $R=a$. The big circle $R=a$ itself is common to the two flow fields.

Under the conditions, a thin metal sphere of radius $a$ could be slid into the flow field perpendicular to the uniform flow so that it coincides with the streamline on $R=a$. Clearly the flow due to the doublet and that due to the free stream would be undisturbed by such a spherical shell. Having done this, the flow due to the doublet could be removed and the outer flow would remain unchanged. Finally, the inside of the shell could be filled to yield a solid sphere. That is, for $R\geq a$, the flow field due to the doublet of strength $Ua^{2}$ and the uniform rectilinear flow of magnitude $U$ give the same flow as that for a uniform flow of magnitude $U$ past a sphere of radius $a$. %The latter flow is shown in Fig. 4.8b.
Then the complex potential for a uniform flow of magnitude $U$ past a sphere of radius $a$ is
\begin{equation}\label{eq:chap_15_2_1_4_13}
F(s)=U(s+\frac{a^{2}}{s}).
\end{equation}
This result is useful in its own right, but it will also be found useful in later sections, through the technique of conformal transformations, to obtain additional solutions.

\section{Sphere with Circulation}\label{sec:chap_15_3_141}

The flow field studied in the previous section not only was irrotational but it produced no circulation around the sphere itself. It was found that there was no hydrodynamic force acting on the sphere under these conditions. It
will be shown in a later section that it is the circulation around a body that produces any lift force that acts on it. It is therefore of interest to study the flow around a sphere that has a circulation around it.

It was established in a previous section that the streamlines for a vortex flow form concentric circles. Therefore, if a vortex was added at the origin to the flow around a sphere, as described in the previous section, the fact that the big circle $R=a$ was a streamline would be unchanged. Thus, from Eq.~(\ref{eq:chap_15_2_1_4_13}), $s_{0}$ being zero in the latter, the complex potential for the flow around a sphere with a negative bound vortex around it will be
\[
F(s)=U(s+\frac{a^{2}}{s})+i\frac{\Gamma}{2\pi}\log s+c.
\]
The negative vortex has been used, since it will turn out that this leads to a positive lift. A constant $c$ has been added to the complex potential for the following reason. For no circulation, it was found that not only was $v$ constant on $R=a$ but the value of the constant was zero. By adding the vortex, $v$ will no longer be zero on $R=a$, although it will have some other constant value. Since it is frequently useful to have the streamline on $R=a$ be $v=0$ when $\sin\varphi=0$, it is desirable to adjust things so that this condition is achieved. By adding a constant $c$ to the complex potential, we have the flexibility to choose $c$ in such away that $v=constant$ becomes $v=0$. This adjustment has no effect on the velocity and pressure distributions, since the velocity components are
defined by derivatives of $v$, so that the absolute value of $v$ at any point is of no significance.

In order to evaluate the constant $c$, the value of the stream function on the big circle $R=a$ will be computed. Then, putting $s=R(e_{xy}e^{i\theta}\cos\varphi+\sin\varphi)$, when $\sin\varphi\neq0$ [See Eq.~(\ref{eq:chap_3_3_2b}), Sec.~(\ref{sec:chap_3_3_29})], the complex potential becomes
\[
\begin{split}
F(s)%&=U(s+\frac{a^{2}}{s})+i\frac{\Gamma}{2\pi}\log s+c \\
%    &=U\{R(e_{xy}e^{i\theta}\cos\varphi+\sin\varphi)+\frac{a^{2}}{R}[e_{xy}(\frac{1}{e^{i\theta}\cos\varphi+\sin\varphi}-\frac{1}{\sin\varphi})+\frac{1}{\sin\varphi}]\}\\
%    &+i\frac{\Gamma}{2\pi}\log [R(e_{xy}e^{i\theta}\cos\varphi+\sin\varphi)]+c\\
    =&e_{xy}U[Re^{i\theta}\cos\varphi+\frac{a^{2}}{R}(\frac{e^{-i\theta}\cos\varphi+\sin\varphi}{1+\cos\theta\sin2\varphi}-\frac{1}{\sin\varphi})]\\
     &-e_{xy}\frac{\Gamma}{2\pi}\{\theta_{L}(\theta,\varphi)-i[\ln R+\ln r_{L}(\theta,\varphi)]\}+U(R\sin\varphi+\frac{a^{2}}{R}\cdot\frac{1}{\sin\varphi})+c\\
    =&e_{xy}\{U[R\cos\theta\cos\varphi+\frac{a^{2}}{R}(\frac{\cos\theta\cos\varphi+\sin\varphi}{1+\cos\theta\sin2\varphi}-\frac{1}{\sin\varphi})]
      -\frac{\Gamma}{2\pi}\theta_{L}(\theta,\varphi) \\
     &+i[U(R\sin\theta\cos\varphi-\frac{a^{2}}{R}\cdot\frac{\sin\theta\cos\varphi}{1+\cos\theta\sin2\varphi})
      +\frac{\Gamma}{2\pi}(\ln R+\ln r_{L}(\theta,\varphi))]+c\} \\
     &+U(R\sin\varphi+\frac{a^{2}}{R}\cdot\frac{1}{\sin\varphi})=e_{xy}[u(R,\theta,\varphi)+iv(R,\theta,\varphi)]+w(R,\varphi)\\
\end{split}
\]
where for $0\leq\varphi\leq\pi/2$ and $\theta\neq(2n+1)\pi$ $(n=0,1,2,\ldots)$
\[
\begin{split}
r_{L}(\theta,\varphi) &= \sqrt{1+\cos\theta\sin2\varphi},\texttt{ } \\
\theta_{L}(\theta,\varphi) &= \arctan\frac{\sin\theta}{\cos\theta+\tan\varphi}+2n\pi,
\end{split}
\]
and by choosing $c=-(i\Gamma/2\pi)\ln a$
\[
\begin{split}
u(R,\theta,\varphi)&=U[(R+\frac{a^{2}/R}{1+\cos\theta\sin2\varphi})\cos\theta\cos\varphi+\frac{a^{2}}{R}\kappa_{x}(\theta,\varphi)]
     -\frac{\Gamma}{2\pi}\theta_{L}(\theta,\varphi), \\
v(R,\theta,\varphi)&=U(R-\frac{a^{2}/R}{1+\cos\theta\sin2\varphi})\sin\theta\cos\varphi
     +\frac{\Gamma}{2\pi}[\ln\frac{R}{a}+\ln r_{L}(\theta,\varphi)], \\
w(R,\varphi)&=U(R\sin\varphi+\frac{a^{2}}{R}\cdot\frac{1}{\sin\varphi})=U[R\sin\varphi+\frac{a^{2}}{R}\kappa_{z}(\varphi)]\\
\end{split}
\]
where if $\sin\varphi=0$ then $\kappa_{x}(\theta,\varphi)=0$ and $\kappa_{z}(\varphi)=0$ else
\[
\kappa_{x}(\theta,\varphi)=\frac{\sin\varphi}{1+\cos\theta\sin2\varphi}-\frac{1}{\sin\varphi}\texttt{ and }\kappa_{z}(\varphi)=\frac{1}{\sin\varphi}.
\]
Hence on the big circle $R=a$ the value of $v$ is indeed constant when $\sin\varphi=0$, so that the value of this constant will be zero when $\sin\varphi=0$. With this value of $c$, the complex potential becomes
\begin{equation}\label{eq:chap_15_3_1_4_14}
F(s)=U(s+\frac{a^{2}}{s})+i\frac{\Gamma}{2\pi}\log\frac{s}{a},
\end{equation}
which describes a uniform rectilinear flow of magnitude $U$ approaching a sphere of radius $a$ that has a negative vortex of strength $\Gamma$ around it. As required, this result agrees with Eq.~(\ref{eq:chap_15_2_1_4_13}) when $\Gamma=0$.

In order to visualize the flow field described by Eq.~(\ref{eq:chap_15_3_1_4_14}) the corresponding velocity components will be evaluated from the complex velocity $\mathbf{V}=\nabla u$. Hence the velocity components are
\begin{equation}\label{eq:chap_15_3_2_4_15}
\begin{split}
u_{R}=&U[(1-\frac{a^{2}/R^{2}}{1+\cos\theta\sin2\varphi})\cos\theta\cos\varphi-\frac{a^{2}}{R^{2}}\kappa_{x}(\theta,\varphi)],\\
u_{\theta}=&-UR[1+\frac{(a^{2}/R^{2})\cos2\varphi}{(1+\cos\theta\sin2\varphi)^{2}}]\sin\theta\cos\varphi \\
           &-\frac{\Gamma}{2\pi}\cdot\frac{2\cos^{2}\varphi+\cos\theta\sin2\varphi}{2(1+\cos\theta\sin2\varphi)},\\
u_{\varphi}=&-UR\{\cos\theta\sin\varphi+\frac{a^{2}}{R^{2}}[\frac{\cos2\theta\cos^{2}\varphi-\sin^{2}\varphi}{(1+\cos\theta\sin2\varphi)^{2}}\cos\varphi\\
            &+\frac{\cos\theta\sin\varphi\cos2\varphi}{(1+\cos\theta\sin2\varphi)^{2}}-\kappa_{z}^{2}(\varphi)\cos\varphi]\}
             -\frac{\Gamma}{2\pi}\cdot\frac{\sin\theta}{1+\cos\theta\sin2\varphi}.\\
\end{split}
\end{equation}
On the surface of the sphere, where $R=a$, Eqs.~(\ref{eq:chap_15_3_2_4_15}) become
\[
\begin{split}
u_{R}=&U[\frac{\cos\theta\sin2\varphi}{1+\cos\theta\sin2\varphi}\cos\theta\cos\varphi-\kappa_{x}(\theta,\varphi)],\\
u_{\theta}=&-Ua[1+\frac{\cos2\varphi}{(1+\cos\theta\sin2\varphi)^{2}}]\sin\theta\cos\varphi
           -\frac{\Gamma}{2\pi}\cdot\frac{2\cos^{2}\varphi+\cos\theta\sin2\varphi}{2(1+\cos\theta\sin2\varphi)},\\
u_{\varphi}=&-Ua\{\cos\theta\sin\varphi+[\frac{\cos2\theta\cos^{2}\varphi-\sin^{2}\varphi}{(1+\cos\theta\sin2\varphi)^{2}}\cos\varphi\\
            &+\frac{\cos\theta\sin\varphi\cos2\varphi}{(1+\cos\theta\sin2\varphi)^{2}}-\kappa_{z}^{2}(\varphi)\cos\varphi]\}
             -\frac{\Gamma}{2\pi}\cdot\frac{\sin\theta}{1+\cos\theta\sin2\varphi}.\\
\end{split}
\]
The fact that $u_{R}=0$ on $R=a$ is to be expected when $\sin\varphi=0$, since this is the boundary
condition. A significant point in the flow field is a point where the velocity components all vanish--that is, a stagnation point. For this flow field the stagnation points are defined by
\begin{equation}\label{eq:chap_15_3_3_4_16}
\sin\theta_{s}=-\frac{\Gamma}{4\pi Ua}
\end{equation}
where $\theta_{s}$ is the value of $\theta$ corresponding to the stagnation point when $\sin\varphi=0$. For $\Gamma=0$, $\sin\theta_{s}=0$, so that $\theta_{s}=0$ or $\theta_{s}=\pi$, which agrees with the case for the sphere without circulation. For nonzero circulation, the value of $\theta_{s}$ clearly depends upon the magnitude of the parameter $\Gamma/(4\pi Ua)$, and it is convenient to discuss Eq.~(\ref{eq:chap_15_3_3_4_16}) for different ranges of this parameter.

The principal interest in the flow around a sphere with circulation is in the study of airfoil theory. By use of conformal transformations the flow around certain airfoil shapes may be transformed into that of the flow around a sphere with circulation.

\section{Joukowski Transformation}\label{sec:chap_15_4_142}

One of the most important transformations in the study of fluid mechanics is the Joukowski transformation. By means of this transformation and the basic flow solutions already studied, it is possible to obtain solutions for the flow around ellipsoids and a family of airfoils. The Joukowski transformation is of the form
\begin{equation}\label{eq:chap_15_4_1_4_20}
s=\varpi+\frac{c^{2}}{\varpi}\texttt{ }(\varpi\neq0)
\end{equation}
where the constant $c^{2}$ is usually taken to be real. According to expression~(\ref{eq:chap_1_7_3}) in Sec.~(\ref{sec:chap_1_7}) with $\varpi=e_{xy}(u+iv)+w$, Eq.~(\ref{eq:chap_15_4_1_4_20}) becomes
\begin{equation}\label{eq:chap_15_4_2}
s=\varpi+\frac{c^{2}}{\varpi}=e_{xy}(u+iv+\frac{c^{2}}{u+iv+w}-\frac{c^{2}}{w})+(w+\frac{c^{2}}{w})\texttt{ }(vw\neq0).
\end{equation}
A general property of the Joukowski transformation is that for large values of $|\varpi|$, $s \to \varpi$. That is, far from the origin the transformation becomes the identity mapping, so that the complex velocity in the two spaces is the same far from the origin. This means that if a uniform flow of a certain magnitude is approaching a body in the $s$ space at some angle of attack, a uniform flow of the same magnitude and angle of attack will approach the corresponding body in the $\varpi$ space.

From Eq.~(\ref{eq:chap_15_4_1_4_20}),
\[
\frac{ds}{d\varpi}=1-\frac{c^{2}}{\varpi^{2}}
\]
so that there is a singular point in the Joukowski transformation at $\varpi=0$. Since we are normally dealing with the flow around some body, the point $\varpi=0$ is normally not in the fluid, and so this singularity is of no consequence. There are also two critical points of the transformation, that is, points at which $ds=d\varpi$ vanishes, at $\varpi=\pm e_{xy}c$. Since smooth curves passing through critical points of a mapping may become corners in the transformed space, it is of interest to investigate the consequence of a smooth curve passing through the critical points of the Joukowski transformation. To do this, consider an arbitrary point $s$ and its counterpart $\varpi$. Let the point $\varpi$ be measured by the radii $\rho_{1}$ and $\rho_{2}$ and the angles $\arg\varpi_{1}$ and $\arg\varpi_{2}$ relative to the two critical points $\varpi=e_{xy}c$ and $\varpi=-e_{xy}c$, respectively. But according to the Joukowski transformation the points $\varpi=\pm e_{xy}c$ map into the points $s=\pm2c$. Then let the mapping of the point $\varpi$ be measured by the radii $R_{1}$ and $R_{2}$ and the angles $\arg s_{1}$ and $\arg s_{2}$ relative to the two points $s=2c$ and $s=-2c$, respectively.

From Eq.~(\ref{eq:chap_15_4_1_4_20}), let
\[
s_{2}=s+2c=\frac{(\varpi+c)^{2}}{\varpi}=\frac{\varpi_{2}^{2}}{\varpi}\texttt{ and }s_{1}=s-2c=\frac{(\varpi-c)^{2}}{\varpi}=\frac{\varpi_{1}^{2}}{\varpi}.
\]
Then
\[
\frac{s_{1}}{s_{2}}=\frac{s-2c}{s+2c}=(\frac{\varpi-c}{\varpi+c})^{2}=(\frac{\varpi_{1}}{\varpi_{2}})^{2}.
\]
Thus
\[
\frac{R_{1}(e_{xy}e^{i\theta_{1}}\cos\varphi_{1}+\sin\varphi_{1})}{R_{2}(e_{xy}e^{i\theta_{2}}\cos\varphi_{2}+\sin\varphi_{2})}
=[\frac{\rho_{1}(e_{xy}e^{i\phi_{1}}\cos\psi_{1}+\sin\psi_{1})}{\rho_{2}(e_{xy}e^{i\phi_{2}}\cos\psi_{2}+\sin\psi_{2})}]^{2}
\]
or
\[
\frac{R_{1}}{R_{2}}\cdot\frac{e_{xy}e^{i\theta_{1}}\cos\varphi_{1}+\sin\varphi_{1}}{e_{xy}e^{i\theta_{2}}\cos\varphi_{2}+\sin\varphi_{2}}
=(\frac{\rho_{1}}{\rho_{2}})^{2}(\frac{e_{xy}e^{i\phi_{1}}\cos\psi_{1}+\sin\psi_{1}}{e_{xy}e^{i\phi_{2}}\cos\psi_{2}+\sin\psi_{2}})^{2}.
\]
According to expression~(\ref{eq:chap_1_7_2}) in Sec.~(\ref{sec:chap_1_7}), equating the modulus and the argument of each side of this equation shows that
\[
\frac{R_{1}}{R_{2}}=(\frac{\rho_{1}}{\rho_{2}})^{2},\texttt{ }
\frac{\sin\varphi_{1}}{\sin\varphi_{2}}=(\frac{\sin\psi_{1}}{\sin\psi_{2}})^{2},\texttt{ }
\]
and
\[
\frac{e^{i\theta_{1}}\cos\varphi_{1}+\sin\varphi_{1}}{e^{i\theta_{2}}\cos\varphi_{2}+\sin\varphi_{2}}
=(\frac{e^{i\phi_{1}}\cos\psi_{1}+\sin\psi_{1}}{e^{i\phi_{2}}\cos\psi_{2}+\sin\psi_{2}})^{2}.
\]

For $j=1,2$, let
\[
e^{i\theta_{j}}\cos\varphi_{j}+\sin\varphi_{j}=r_{j}e^{i\Theta_{j}}\texttt{ and }
e^{i\phi_{j}}\cos\psi_{j}+\sin\psi_{j}=\eta_{j}e^{i\Phi_{j}}.
\]
Then the equation above becomes
\[
\frac{r_{1}}{r_{2}}e^{i(\Theta_{1}-\Theta_{2})}=(\frac{\eta_{1}}{\eta_{2}})^{2}e^{2i(\Phi_{1}-\Phi_{2})}.
\]
Equating the modulus and the argument of each side of this equation shows that
\[
\frac{r_{1}}{r_{2}}=(\frac{\eta_{1}}{\eta_{2}})^{2}\texttt{ and }
\Theta_{1}-\Theta_{2}=2(\Phi_{1}-\Phi_{2}).
\]
This last result shows that if a smooth curve passes through the point $\varpi=e_{xy}c$, the corresponding curve in the $s$ space will form a knife-edge or cusp. This may be verified by considering a smooth curve to pass through the point $\varpi=e_{xy}c$. Two points on this curve show that $\Phi_{1}$ changes from $3\pi/2$ to $\pi/2$ and $\Phi_{2}$ changes from $2\pi$ to $0$ as the critical point is passed. That is, the value of $\Phi_{1}-\Phi_{2}$ changes from $-\pi/2$ to $\pi/2$, giving a difference of $\pi$. From the result $\Theta_{1}-\Theta_{2}=2(\Phi_{1}-\Phi_{2})$, it follows that the corresponding difference in the value of $\Theta_{1}-\Theta_{2}$ will be $2\pi$. This yields a knife-edge or cusp in the $s$ space. That is, if a smooth curve passes through either of the critical points $\varpi=\pm e_{xy}c$, the corresponding curve in the $s$ space will contain a knife-edge at the corresponding critical
point $s=\pm2c$.

An example of a smooth curve that passes through both critical points is a sphere centered at the origin of the $\varpi$ space and whose radius is $c$, the constant that appears in the Joukowski transformation. Then, on this sphere $\varpi=c(e_{xy}e^{i\phi}\cos\psi+\sin\psi)$, and the value of $s$ will be given by
\[
\begin{split}
s=&c(e_{xy}e^{i\phi}\cos\psi+\sin\psi)+\frac{c}{e_{xy}e^{i\phi}\cos\psi+\sin\psi} \\
 =&c[e_{xy}(e^{i\phi}\cos\psi+\frac{1}{e^{i\phi}\cos\psi+\sin\psi}-\frac{1}{\sin\psi})+(\sin\psi+\frac{1}{\sin\psi})]\\
 =&c[e_{xy}(\frac{e^{i\phi}\cos\phi\sin2\psi\cos\psi}{1+\cos\phi\sin2\psi}+\frac{2\cos\phi\cos\psi+\sin\psi}{1+\cos\phi\sin2\psi}-\frac{1}{\sin\psi})+(\sin\psi+\frac{1}{\sin\psi})]\\
 =&c\{e_{xy}[\frac{2+\cos\phi\sin2\psi}{1+\cos\phi\sin2\psi}\cos\phi\cos\psi+\kappa_{x}(\phi,\psi)\\
  &+i\frac{\cos\phi\sin2\psi}{1+\cos\phi\sin2\psi}\sin\phi\cos\psi]+[\sin\psi+\kappa_{z}(\psi)]\}
\end{split}
\]
where if $\sin\psi=0$ then $\kappa_{x}(\phi,\psi)=0$ and $\kappa_{z}(\psi)=0$ else
\[
\kappa_{x}(\phi,\psi)=\frac{\sin\psi}{1+\cos\phi\sin2\psi}-\frac{1}{\sin\psi}\texttt{ and }\kappa_{z}(\psi)=\frac{1}{\sin\psi}.
\]
Thus there are
\[
\begin{split}
x&=c\frac{2+\cos\phi\sin2\psi}{1+\cos\phi\sin2\psi}\cos\phi\cos\psi+c\kappa_{x}(\phi,\psi)\\
 &=c\xi(\phi,\psi)\cos\phi\cos\psi+c\kappa_{x}(\phi,\psi),\\
y&=c\frac{\cos\phi\sin2\psi}{1+\cos\phi\sin2\psi}\sin\phi\cos\psi=c\eta(\phi,\psi)\sin\phi\cos\psi, \\
z&=c(\sin\psi+\kappa_{z}(\psi))=c\zeta(\psi)\sin\psi
\end{split}
\]
where
\[
\begin{split}
\xi(\phi,\psi)&=\frac{2+\cos\phi\sin2\psi}{1+\cos\phi\sin2\psi}, \\
\eta(\phi,\psi)&=\frac{\cos\phi\sin2\psi}{1+\cos\phi\sin2\psi}, \\
\zeta(\psi)&=1+\kappa_{z}^{2}(\psi).
\end{split}
\]
So we get
\[
(\frac{x-c\kappa_{x}(\phi,\psi)}{c\xi(\phi,\psi)})^{2}+(\frac{y}{c\eta(\phi,\psi)})^{2}+(\frac{z}{c\zeta(\psi)})^{2}=1.
\]
That is, the sphere in the $\varpi$ space maps into the ellipsoid in the $s$ space. It is readily verified that all points that lie outside the sphere $|\varpi|=c$ cover the entire $s$ space. However, the points inside the sphere $|\varpi|=c$ also cover the entire $s$ space, so that the transformation is double-valued. This is readily verified by observing that for any value of $\varpi$, Eq.~(\ref{eq:chap_15_4_1_4_20}) yields the same value of $s$ for that value of $\varpi$ and also for $c^{2}/\varpi$. It will be noted that $c^{2}/\varpi$ is simply the image of the point $\varpi$ inside the sphere of radius $c$.

This double-valued property of the Joukowski transformation is treated by connecting the two points $s=\pm e_{xy}2c$ by a branch cut along the $x$ axis and creating two Riemann sheets. Then the mapping is single-valued if all the points outside the sphere $|\varpi|=c$ are taken to fall on one of these sheets and all the points inside the sphere to fall on the other sheet. In fluid mechanics, difficulties due to the double-valued behavior do not usually arise because the points $|\varpi|<c$ usually lie inside some body about which the flow is being studied, so that these points are not in the flow field in the $s$ space.

\section{Flows Around an Ellipsoid}\label{sec:chap_15_5_143}

Applications of the Joukowski transformation will be made in an inverse sense. That is, the simple geometry of the sphere, the flow around which is known, will be placed in the $\varpi$ space, and the corresponding body that results in the $s$ space will be investigated.

Consider, first, the constant $c$ in Eq.~(\ref{eq:chap_15_4_1_4_20}) to be real and positive, and consider a sphere of radius $a>c$ to be centered at the origin of the $\varpi$ space. The contour in the $s$ space corresponding to this sphere in the $\varpi$ space may be identified by substituting $\varpi=a(e_{xy}e^{i\phi}\cos\psi+\sin\psi)$ into Eq.~(\ref{eq:chap_15_4_1_4_20}):
\[
\begin{split}
s=&a(e_{xy}e^{i\phi}\cos\psi+\sin\psi)+\frac{c^{2}/a}{e_{xy}e^{i\phi}\cos\psi+\sin\psi} \\
 =&a[e_{xy}(e^{i\phi}\cos\psi+\frac{c^{2}/a^{2}}{e^{i\phi}\cos\psi+\sin\psi}-\frac{c^{2}/a^{2}}{\sin\psi})+(\sin\psi+\frac{c^{2}/a^{2}}{\sin\psi})]\\
 =&a\{e_{xy}[e^{i\phi}\cos\psi+(\frac{c}{a})^{2}(\frac{e^{-i\phi}\cos\psi+\sin\psi}{1+\cos\phi\sin2\psi}-\frac{1}{\sin\psi})]+(\sin\psi+\frac{c^{2}/a^{2}}{\sin\psi})\}\\
 =&a\{e_{xy}[(1+\frac{c^{2}/a^{2}}{1+\cos\phi\sin2\psi})\cos\phi\cos\psi+\kappa_{x}(\phi,\psi) \\
  &+i(1-\frac{c^{2}/a^{2}}{1+\cos\phi\sin2\psi})\sin\phi\cos\psi]+(\sin\psi+\frac{c^{2}}{a^{2}}\kappa_{z}(\psi))\}
\end{split}
\]
where if $\sin\psi=0$ then $\kappa_{x}(\phi,\psi)=0$ and $\kappa_{z}(\psi)=0$ else
\[
\kappa_{x}(\phi,\psi)=\frac{\sin\psi}{1+\cos\phi\sin2\psi}-\frac{1}{\sin\psi}\texttt{ and }\kappa_{z}(\psi)=\frac{1}{\sin\psi}.
\]
Equating real and imaginary parts of this equation gives
\[
\begin{split}
x&=a(1+\frac{c^{2}/a^{2}}{1+\cos\phi\sin2\psi})\cos\phi\cos\psi+a\kappa_{x}(\phi,\psi) \\
 &=a\xi(\phi,\psi)\cos\phi\cos\psi+a\kappa_{x}(\phi,\psi), \\
y&=a(1-\frac{c^{2}/a^{2}}{1+\cos\phi\sin2\psi})\sin\phi\cos\psi=a\eta(\phi,\psi)\sin\phi\cos\psi, \\
z&=a[\sin\psi+\frac{c^{2}}{a^{2}}\kappa_{z}(\psi)]=a\zeta(\psi)\sin\psi
\end{split}
\]
where
\[
\begin{split}
 \xi(\phi,\psi)&=1+\frac{c^{2}/a^{2}}{1+\cos\phi\sin2\psi},\\
\eta(\phi,\psi)&=1-\frac{c^{2}/a^{2}}{1+\cos\phi\sin2\psi}, \\
\zeta(\psi)&=1+[\frac{c}{a}\kappa_{z}(\psi)]^{2}.
\end{split}
\]
These are the parametric equations of the required curve in the $s$ space. The equation of the curve may be obtained by use of the identities $\cos^{2}\phi+\sin^{2}\phi=1$ and $\cos^{2}\phi-\sin^{2}\phi=\cos2\phi$. This gives
\[
(\frac{x-a\kappa_{x}(\phi,\psi)}{a\xi(\phi,\psi)})^{2}+(\frac{y}{a\eta(\phi,\psi)})^{2}+(\frac{z}{a\zeta(\psi)})^{2}=1
\]
which is the equation of an ellipsoid whose major semiaxis is of length
\[
a\xi(\phi,\psi)=a(1+\frac{c^{2}/a^{2}}{1+\cos\phi\sin2\psi}),
\]
aligned along the $x$ axis, and whose minor semiaxis is of length
\[
a\eta(\phi,\psi)=a(1-\frac{c^{2}/a^{2}}{1+\cos\phi\sin2\psi}),
\]
aligned along the $y$ axis. Then, in order to obtain the complex potential for a uniform flow of magnitude $U$ approaching this ellipsoid at an angle $\alpha$ of attack described by $\gamma=e_{xy}e^{i\alpha}\cos\beta+\sin\beta$, the same flow should be considered to approach the sphere in the $\varpi$ space. But %it is shown in Prob. 4.5, Eq. (4.29), that
from the equation
\begin{equation}\label{eq:chap_15_5_1_4_29}
F(s)=U(s\bar{\gamma}+\frac{a^{2}}{s\bar{\gamma}}),
\end{equation}
the complex potential for a uniform flow of magnitude $U$ approaching a sphere of radius $a$ at a pair of angles $\alpha$ and $\beta$ to the reference axis is
\[
F(\varpi)=U(\varpi\bar{\gamma}+\frac{a^{2}}{\varpi\bar{\gamma}})\texttt{ }(\bar{\gamma}=e_{xy}e^{-i\alpha}\cos\beta+\sin\beta).
\]
Then, by solving Eq.~(\ref{eq:chap_15_4_1_4_20}) for $\varpi$ in terms of $s$, the complex potential in the $s$
space may be obtained. From Eq.~(\ref{eq:chap_15_4_1_4_20}), $\varpi^{2}-s\varpi+c^{2}=0$, then
\[
\varpi=\frac{s}{2}\pm\sqrt{(\frac{s}{2})^{2}-c^{2}}.
\]
Since it is known that $\varpi \to s$ for large values of $s$, the positive root must be chosen. Then the complex potential in the $s$ space becomes
\[
\begin{split}
F(s)&=U\{[\frac{s}{2}+\sqrt{(\frac{s}{2})^{2}-c^{2}}]\bar{\gamma}+\frac{a^{2}}{s/2+\sqrt{(s/2)^{2}-c^{2}}}\frac{\gamma}{\bar{\gamma}\gamma}\}\\
    &=U\{[s-\frac{s}{2}+\sqrt{(\frac{s}{2})^{2}-c^{2}}]\bar{\gamma}+\frac{a^{2}}{c^{2}}[\frac{s}{2}-\sqrt{(\frac{s}{2})^{2}-c^{2}}]\frac{\gamma}{\bar{\gamma}\gamma}\}\\
\end{split}
\]
where the last term has been rationalized. By writing $s/2$ as $s-s/2$ in the first term, two of the terms may be combined as follows:
\begin{equation}\label{eq:chap_15_5_1_4_21_a}
F(s)=U[s\bar{\gamma}+(\frac{a^{2}}{c^{2}}\frac{\gamma}{\bar{\gamma}\gamma}-\bar{\gamma})(\frac{s}{2}-\sqrt{(\frac{s}{2})^{2}-c^{2}})].
\end{equation}
Equation~(\ref{eq:chap_15_5_1_4_21_a}) is the complex potential for a uniform flow of magnitude $U$ approaching an ellipsoid whose major semiaxis is $a\xi(\phi,\psi)$ and whose minor semiaxis is $a\eta(\phi,\psi)$ at an angle $\alpha$ of attack described by $\gamma$ to the major axis. In this form it may be seen that the complex potential consists of that for a uniform flow at a pair of angles $\alpha$ and $\beta$ to the reference axis plus a perturbation which is large near the body but vanishes for large values of $s$. %The flow field generated by the complex potential~(\ref{eq:chap_15_5_1_4_21_a}) is shown with that for the sphere in the $\varpi$ space.

The stagnation points in the $\varpi$ space are located at
\[
\varpi=a\gamma\texttt{ and }\varpi=a[e_{xy}e^{i\alpha}\cos(\beta+\pi)+\sin(\beta+\pi)]=-a\gamma,
\]
that is, at $\varpi=\pm a\gamma$. Then, from Eq.~(\ref{eq:chap_15_4_1_4_20}), the corresponding points in the $s$ space are
\[
\begin{split}
s=&\varpi+\frac{c^{2}}{\varpi}=\pm a\gamma\pm\frac{c^{2}}{a\gamma}=\pm a\gamma\pm\frac{c^{2}/a}{e_{xy}e^{i\alpha}\cos\beta+\sin\beta} \\
 =&\pm a\gamma\pm\frac{c^{2}}{a}[e_{xy}(\frac{1}{e^{i\alpha}\cos\beta+\sin\beta}-\frac{1}{\sin\beta})+\frac{1}{\sin\beta}]\\
 =&\pm\{e_{xy}[ae^{i\alpha}\cos\beta+\frac{c^{2}}{a}(\frac{e^{-i\alpha}\cos\beta+\sin\beta}{1+\cos\alpha\sin2\beta}-\frac{1}{\sin\beta})]+(a\sin\beta+\frac{c^{2}/a}{\sin\beta})\}\\
 =&\pm\{e_{xy}[(a+\frac{c^{2}/a}{1+\cos\alpha\sin2\beta})\cos\alpha\cos\beta+\frac{c^{2}}{a}\kappa_{x}(\alpha,\beta)\\
  &+i(a-\frac{c^{2}/a}{1+\cos\alpha\sin2\beta})\sin\alpha\cos\beta]+[a\sin\beta+\frac{c^{2}}{a}\kappa_{z}(\beta)]\}
\end{split}
\]
where if $\sin\beta=0$ then $\kappa_{x}(\alpha,\beta)=0$ and $\kappa_{z}(\beta)=0$ else
\[
\kappa_{x}(\alpha,\beta)=\frac{\sin\beta}{1+\cos\alpha\sin2\beta}-\frac{1}{\sin\beta}\texttt{ and }\kappa_{z}(\beta)=\frac{1}{\sin\beta}.
\]
This gives the stagnation points as
\[
\begin{split}
x&=\pm[(a+\frac{c^{2}/a}{1+\cos\alpha\sin2\beta})\cos\alpha\cos\beta+\frac{c^{2}}{a}\kappa_{x}(\alpha,\beta)]\\
 &=\pm [a\xi(\alpha,\beta)\cos\alpha\cos\beta+\frac{c^{2}}{a}\kappa_{x}(\alpha,\beta)], \\
y&=\pm(a-\frac{c^{2}/a}{1+\cos\alpha\sin2\beta})\sin\alpha\cos\beta=\pm a\eta(\alpha,\beta)\sin\alpha\cos\beta,\\
z&=\pm[a\sin\beta+\frac{c^{2}}{a}\kappa_{z}(\beta)]=\pm[a+\frac{c^{2}}{a}\kappa_{z}^{2}(\beta)]\sin\beta=\pm a\zeta(\beta)\sin\beta
\end{split}
\]
where
\[
\begin{split}
 \xi(\phi,\psi)&=1+\frac{c^{2}/a^{2}}{1+\cos\phi\sin2\psi},\\
\eta(\phi,\psi)&=1-\frac{c^{2}/a^{2}}{1+\cos\phi\sin2\psi}, \\
\zeta(\psi)&=1+\frac{c^{2}}{a^{2}}\kappa_{z}^{2}(\beta).
\end{split}
\]

Equation~(\ref{eq:chap_15_5_1_4_21_a}) includes two special cases within its range of validity. For $\alpha=0$ it describes a uniform rectilinear flow approaching a horizontally oriented ellipsoid, and for $\alpha=\pi/2$ it describes a uniform vertical flow approaching the same horizontally oriented ellipsoid. However, it is of interest to note that the solution for a uniform rectilinear flow approaching a vertically oriented ellipsoid may be obtained directly from the Joukowski transformation with a slight modification. Substitute $c=ib$, where $b$ is real and positive, into Eq.~(\ref{eq:chap_15_4_1_4_20})
\[
s=\varpi-\frac{b^{2}}{\varpi}.
\]
Then, as with the horizontal ellipsoid, examining the mapping of the sphere
$\varpi=a(e_{xy}e^{i\phi}\cos\psi+\sin\psi)$ gives the parametric equations of the mapped boundary
\[
\begin{split}
x&=\pm[(a-\frac{b^{2}/a}{1+\cos\alpha\sin2\beta})\cos\alpha\cos\beta-\frac{b^{2}}{a}\kappa_{x}(\alpha,\beta)]\\
 &=\pm[a\xi(\alpha,\beta)\cos\alpha\cos\beta-\frac{b^{2}}{a}\kappa_{x}(\alpha,\beta)], \\
y&=\pm(a+\frac{b^{2}/a}{1+\cos\alpha\sin2\beta})\sin\alpha\cos\beta=\pm a\eta(\alpha,\beta)\sin\alpha\cos\beta, \\
z&=\pm[a\sin\beta-\frac{b^{2}}{a}\kappa_{z}(\beta)]=\pm[a-\frac{b^{2}}{a}\kappa_{z}^{2}(\beta)]\sin\beta=\pm a\zeta(\beta)\sin\beta
\end{split}
\]
where
\[
\begin{split}
 \xi(\phi,\psi)&=1-\frac{b^{2}/a^{2}}{1+\cos\phi\sin2\psi},\\
\eta(\phi,\psi)&=1+\frac{b^{2}/a^{2}}{1+\cos\phi\sin2\psi}, \\
\zeta(\psi)&=1-\frac{b^{2}}{a^{2}}\kappa_{z}^{2}(\beta).
\end{split}
\]
Thus the equation of the contour in the $s$ space is
\[
(\frac{x+(b^{2}/a)\kappa_{x}(\alpha,\beta)}{a\xi(\phi,\psi)})^{2}+(\frac{y}{a\eta(\phi,\psi)})^{2}+(\frac{z}{a\zeta(\psi)})^{2}=1
\]
which is the equation of an ellipsoid whose major semiaxis is of length
\[
a\eta(\phi,\psi)=a(1+\frac{b^{2}/a^{2}}{1+\cos\phi\sin2\psi})
\]
which is aligned along the $y$ axis, and whose minor semiaxis is of length
\[
a\xi(\phi,\psi)=a(1-\frac{b^{2}/a^{2}}{1+\cos\phi\sin2\psi})
\]
which is aligned along the $x$ axis. Then to obtain uniform rectilinear flow approaching such an ellipsoid the same flow should approach the sphere in the $\varpi$ space. Thus the required complex potential, from Eq.~(\ref{eq:chap_15_2_1_4_13}), is
\[
F(\varpi)=U(\varpi+\frac{a^{2}}{\varpi}).
\]
But the inverted equation of the mapping for which $\varpi \to s$ as $s \to \infty$ is
\[
\varpi=\frac{s}{2}+\sqrt{(\frac{s}{2})^{2}+b^{2}}.
\]
Hence the complex potential in the $s$ space is
\begin{equation}\label{eq:chap_15_5_2_4_21_b}
\begin{split}
F(s)&=U[\frac{s}{2}+\sqrt{(\frac{s}{2})^{2}+b^{2}}+\frac{a^{2}}{s/2+\sqrt{(s/2)^{2}+b^{2}}}] \\
F(s)&=U[s-(1+\frac{a^{2}}{b^{2}})(\frac{s}{2}-\sqrt{(\frac{s}{2})^{2}+b^{2}})]
\end{split}
\end{equation}
in which the same rationalization and simplification has been carried out as before. Again the complex potential is in the form of that for a uniform flow plus a perturbation which is large near the body and which vanishes at large distances from the body. Equation~(\ref{eq:chap_15_5_2_4_21_b}) describes a uniform rectilinear flow of magnitude $U$ approaching a vertically oriented ellipsoid. %The flow field for this situation is shown in Fig.4.14b.

\section{Kutta condition and the Flat-Plate Airfoils}\label{sec:chap_15_6_144}

It was observed that the potential flow solution for flow around a sharp edge contained a singularity at the edge itself. This singularity required an infinite velocity at the point in question, which, of course, is physically impossible. The question arises, then, as to what the real flow situation would be in a physical experiment. Depending upon the actual physical configuration, one of two remedial situations will prevail. One possibility is that the fluid will separate from the solid surface at the knife-edge. The resulting free stream line configuration would be such that the radius of curvature at the edge becomes finite rather than being zero. As a consequence, the velocities there will remain finite. Examples of this type of solution will be discussed later in this chapter.

A second possibility is that a stagnation point exists at the sharp edge. For the flow around finite bodies, stagnation points exist, and it seems possible that a stagnation point could be induced by the flow field to move to the location of the sharp edge. This possibility leads to the so-called Kutta condition, and it will be discussed below in the context of the flat-plate airfoil, that is, a flat plate which is at some angle of attack to the free stream.

In the previous section, the flow around an ellipsoid was obtained from the Joukowski transformation [Eq.~(\ref{eq:chap_15_4_1_4_20})] by considering the flow around a sphere of radius $a>c$ in the $\varpi$ space. Now, if the constant $c$ is allowed to approach the magnitude of the radius $a$, the resulting ellipsoid in the $s$ space degenerates to a flat plate defined by the strip $-2a\leq x\leq2a$. The resulting flow field, as defined by Eq.~(\ref{eq:chap_15_5_1_4_21_a}), is shown. Because of the angle of attack, the stagnation points do not coincide with the leading and trailing edges of the flat plate. Rather, the upstream stagnation point is located on the lower surface and the downstream stagnation point is located on the upper surface at the points $x=\pm2a\cos\phi$. Then, around both the leading and trailing edges, the flow will be that associated with a sharp edge, %which was discussed in Sec. 4.6. In that section it was observed that in finite velocity components existed at the edge itself,
a situation that is physically impossible to realize.

The difficulty encountered above with the flat-plate airfoil does not occur at the leading edge of real airfoils because real airfoils have a finite thickness and so have a finite radius of curvature at the leading edge. However, the trailing edge of airfoils is usually quite sharp, so that the difficulty of in finite velocity components still exists there. However, this remaining difficulty would also be overcome if the stagnation point which is near the trailing edge was actually at the trailing edge. This would be accomplished if a circulation existed around the flat plate and the magnitude of this circulation was just the amount required to rotate the rear stagnation point so that its location coincides with the trailing edge. This condition is called the Kutta condition, and it may be restated as follows: For bodies with sharp trailing edges that are at small angles of attack to the free stream, the flow will adjust itself in such a way that the rear stagnation point coincides with the trailing edge.

The amount of circulation required to comply with the Kutta condition may be determined as follows: In the $\varpi$ space, the rear stagnation point is located at the point $\varpi=a(e_{xy}e^{\phi}\cos\psi+\sin\psi)$. But, according to the Kutta condition, the rear stagnation point should be located at the point $s=2a$, which corresponds to the point $\varpi=a$. That is, the stagnation point on the downstream face of the sphere in the $\varpi$ space should be rotated clockwise through an angle $\phi$. But from Eq.~(\ref{eq:chap_15_3_3_4_16}), the magnitude of the circulation that will do this is
\begin{equation}\label{eq:chap_15_6_1_4_22_a}
\Gamma=4\pi Ua\sin\phi
\end{equation}
in the clockwise direction (that is, negative circulation). Then the complex potential for the required flow in the $\varpi$ space is, from Eqs.~(\ref{eq:chap_15_3_1_4_14}) and~(\ref{eq:chap_15_5_1_4_29}),
\[
F(\varpi)=U[\varpi\bar{\gamma}+\frac{a^{2}}{\varpi}\frac{\gamma}{\bar{\gamma}\gamma}]+i2Ua\sin\phi\log\frac{\varpi}{a}
\]
where $\gamma=e_{xy}e^{i\alpha}\cos\beta+\sin\beta$. But the equation of the mapping is
\[
s=\varpi+\frac{a^{2}}{\varpi}
\]
and the inverse, which gives $\varpi \to s$ as $s \to \infty$, is
\[
\varpi=\frac{s}{2}+\sqrt{(\frac{s}{2})^{2}-a^{2}}.
\]
Then the complex potential in the $s$ space is
\begin{equation}\label{eq:chap_15_6_2_4_22_b}
\begin{split}
F(\varpi)=&U\{[\frac{s}{2}+\sqrt{(\frac{s}{2})^{2}-a^{2}}]\bar{\gamma}+\frac{a^{2}}{s/2+\sqrt{(s/2)^{2}-a^{2}}}\frac{\gamma}{\bar{\gamma}\gamma}\}\\
          &+i2Ua\sin\phi\log\frac{\varpi}{a}.
\end{split}
\end{equation}
%The flow field corresponding to this complex potential is shown in Fig.4.15b.
Although the flow at the trailing edge is now regular, the singularity at the leading edge still exists. In an actual flow configuration the fluid would separate at the leading edge and reattach again on the top side of the airfoil. The streamline $v=0$ would then correspond to a finite curvature, and the velocity components would remain finite at the leading edge.

The lift force generated by the flat-plate airfoil may be calculated from the Kutta-Joukowski law. Then, denoting the lift force by $Y$ and using the value of the circulation given by Eq.~(\ref{eq:chap_15_6_1_4_22_a}),
\[
Y=4\pi\rho U^{2}a\sin\phi.
\]
It is usual to express lift forces in terms of the dimensionless lift coefficient $C_{L}$, which is defined as follows:
\[
C_{L}=\frac{Y}{(1/2)\rho U^{2}l}
\]
where $l$ is the length or chord of the airfoil, which, for the flat plate under consideration, equals $4a$. Then the value of the lift coefficient for the flat-plate airfoil is
\begin{equation}\label{eq:chap_15_6_3_4_22_c}
C_{L}=2\pi\sin\phi.
\end{equation}
This result shows that the lift coefficient for the flat-plate airfoil increases with angle of attack, and for small values of $a$, for which $\sin\phi\approx\phi$, the lift coefficient is proportional to the angle of attack with a constant of proportionality of $2\pi$. This result is very close to experimental observations, and so the Kutta condition appears to be well justified. If the Kutta condition were not valid, there would be no circulation around the flat plate, and consequently no lift would be generated. This would mean that kites would not be able to fly.

\section{Symmetrical Joukowski Airfoils}\label{sec:chap_15_7_145}

A family of airfoils may be obtained in the $s$ space by considering the Joukowski transformation in conjunction with a series of spheres in the $\varpi$ space whose centers are slightly displaced from the origin. These airfoils are known as the Joukowski family of airfoils.

Consider, first, the case where the center of the sphere in the $\varpi$ space is displaced from the origin along the real axis. It must then be decided in which direction the center should be moved and what radius should be employed, relative to the Joukowski constant $c$. From previous sections it is known that if the big circumference of the sphere passes through either of the two critical points of the Joukowski transformation, $\varpi=\pm e_{xy}c$, then a sharp edge or cusp is obtained in the $s$ space. Then, if the leading edge of the airfoil is to have a finite radius of curvature and if there should be no singularities in the flow field itself, it follows that the point $\varpi=-e_{xy}c$ should be inside the sphere in the $\varpi$ space. Also, since the trailing edge of the airfoil should be sharp as opposed to being blunt, the big circumference of the sphere should pass through the point $\varpi=-e_{xy}c$. These conditions will be satisfied by taking the center of the sphere to be at $\varpi=-e_{xy}m$, where $m$ is real on the $e_{xy}$ coordinate plane, and by choosing the radius of the sphere to be $c+m$. %Such a configuration is shown in Fig.4.16a.
The radius $a$ is given by
\[
a=c+m=c(1+\varepsilon)
\]
where the parameter $\varepsilon=m/c$ will be assumed to be small compared with unity. When $\varepsilon=0$, the flat-plate airfoil is recovered, so that for $\varepsilon\ll 1$ it may be anticipated that a thin airfoil will be obtained. The significance of the restriction $\varepsilon\ll 1$ will be that all the equations may be linearized in $\varepsilon$, which will permit a closed-form solution for the equation of the airfoil surface in the $s$ space. Also %shown in Fig. 4.16a
it is the airfoil that is obtained in the $s$ space and its principal parameters, the chord $l$ and the maximum thickness $t$. It is now required to relate these parameters to the free parameters $a$ and $m$ and to establish the equation of the airfoil surface in the $s$ space.

To establish the chord of the airfoil in terms of the chosen radius $a$ and offset $m$, it is only necessary to find the mapping of the points $\varpi=e_{xy}c$ and $\varpi=-(c+2m)$, since these points correspond to the trailing and leading edges, respectively. Using the Joukowski transformation, the mapping of the point $\varpi=e_{xy}c$ is $s=2c$. Also, the mapping of the point $\varpi=-(c+2m)=-c(1+2\varepsilon)$ is
\[
s=-c(1+2\varepsilon)-\frac{c}{1+2\varepsilon}.
\]
Since it was decided to linearize all quantities in $\varepsilon$, the value of $s$ will be obtained to the first order in $\varepsilon$ only
\[
\begin{split}
s&=-c(1+2\varepsilon)-c[1-2\varepsilon+O(\varepsilon^{2})] \\
 &=-2c+O(\varepsilon^{2}).
\end{split}
\]
That is, to the first order in $\varepsilon$ the lending edge of the airfoil is located at
$s=-2c$ that the chord length is
\[
l=4c.
\]
This means that, correct to the first order in $\varepsilon$, the length of the airfoil is
unchanged by the shifting of the center of the sphere in the $\varpi$ space.

In order to determine the maximum thickness $t$, the equation of the airfoil surface must be obtained. This may be done by inserting the equation of the surface in the $\varpi$ space into the Joukowski transformation. But in the $\varpi$ space the polar radius $R$ to the big circumference of the sphere is a function of the pair of angles $\phi$ and $\psi$. In order to establish this dependence, the cosine rule will be applied to the triangle defined by the radius $a$, the coordinate $R$, and the real $x$ axis in the $e_{xy}$ coordinate plane in the $\varpi$ space with the component on the real $x$ axis. %, as shown in Fig.4.16a.
Thus
\[
\begin{split}
a^{2}&=(R\cos\psi)^{2}+m^{2}-2Rm\cos(\pi-\phi)\cos\psi+(R\sin\psi)^{2} \\
     &=R^{2}+m^{2}+2Rm\cos\phi\cos\psi.
\end{split}
\]
But $a=c+m$, so that the equation above may be written in the form
\[
(c+m)^{2}=R^{2}(1+\frac{m^{2}}{R^{2}}+2\frac{m}{R}\cos\phi\cos\psi).
\]
Now since $R\geq c$, it follows that $m/R\leq m/c$ so that, to first order in $\varepsilon=m/c$, the term $m^{2}/R^{2}$ may be neglected. The equation for $R$ then becomes
\[
\begin{split}
c+m &= R(1+\frac{m^{2}}{R^{2}}+2\frac{m}{R}\cos\phi\cos\psi)^{1/2} \\
    &= R[1+\frac{m}{R}\cos\phi\cos\psi+O(\varepsilon^{2})].
\end{split}
\]
Thus to the first order in $\varepsilon$, this relation becomes
\[
c(1+\varepsilon)=R+m\cos\phi\cos\psi.
\]
Then
\[
R=c[1+\varepsilon(1-\cos\phi\cos\psi)]=c[1+\mu(\phi,\psi)]
\]
where
\[
\mu(\phi,\psi)=\varepsilon(1-\cos\phi\cos\psi).
\]
This is the required equation that gives the variation of the radius $R$ with the angle $\phi$ for points on the big circumference of the sphere %in the $e_{xy}$ plane
in the $\varpi$ space where $\varpi=R(e_{xy}e^{i\phi}\cos\psi+\sin\psi)$. Then, in order to determine the equation of the corresponding profile in the $s$ space, this result should be substituted into the Joukowski transformation [Eq.~(\ref{eq:chap_15_4_1_4_20}), Sec.~(\ref{sec:chap_15_4_142})]. Thus points on the surface of the airfoil will be defined by
\[
\begin{split}
s=&c[1+\mu(\phi,\psi)](e_{xy}e^{i\phi}\cos\psi+\sin\psi)+\frac{c}{[1+\mu(\phi,\psi)](e_{xy}e^{i\phi}\cos\psi+\sin\psi)}\\
 =&c[1+\mu(\phi,\psi)](e_{xy}e^{i\phi}\cos\psi+\sin\psi) \\
  &+\frac{c}{1+\mu(\phi,\psi)}\{e_{xy}[\frac{e^{-i\phi}\cos\psi}{1+\cos\phi\sin2\psi}+\kappa_{x}(\phi,\psi)]+\kappa_{z}(\psi)\}
\end{split}
\]
where if $\sin\psi=0$ then $\kappa_{x}(\phi,\psi)=0$ and $\kappa_{z}(\psi)=0$ else
\[
\kappa_{x}(\phi,\psi)=\frac{\sin\psi}{1+\cos\phi\sin2\psi}-\frac{1}{\sin\psi}\texttt{ and }\kappa_{z}(\psi)=\frac{1}{\sin\psi}.
\]
This equation may also be linearized in $\varepsilon$ as follows:
\[
\begin{split}
s=&c[1+\mu(\phi,\psi)](e_{xy}e^{i\phi}\cos\psi+\sin\psi) \\
  &+c[1-\mu(\phi,\psi)+O(\varepsilon^{2})]\{e_{xy}[\frac{e^{-i\phi}\cos\psi}{1+\cos\phi\sin2\psi}+\kappa_{x}(\phi,\psi)]+\kappa_{z}(\psi)\}\\
% =&c\{e_{xy}[e^{i\phi}\cos\psi+\frac{e^{-i\phi}\cos\psi+\sin\psi}{1+\cos\phi\sin2\psi}-\frac{1}{\sin\psi} \\
% &+\mu(\phi,\psi)(e^{i\phi}\cos\psi-\frac{e^{-i\phi}\cos\psi+\sin\psi}{1+\cos\phi\sin2\psi}+\frac{1}{\sin\psi})]\\
%  &+(\sin\psi+\frac{1}{\sin\psi})+\mu(\phi,\psi)(\sin\psi-\frac{1}{\sin\psi})+O(\varepsilon^{2})\}\\
 =&c\{e_{xy}[\frac{e^{i\phi}\cos\phi\sin2\psi\cos\psi}{1+\cos\phi\sin2\psi}+\frac{2\cos\phi\cos\psi}{1+\cos\phi\sin2\psi}+\kappa_{x}(\phi,\psi)\\
  &+\mu(\phi,\psi)(\frac{e^{i\phi}\cos\phi\sin2\psi\cos\psi}{1+\cos\phi\sin2\psi}-\kappa_{x}(\phi,\psi)+i\frac{2\sin\phi\cos\psi}{1+\cos\phi\sin2\psi})]\\
  &+[\sin\psi+\kappa_{z}(\psi)]+\mu(\phi,\psi)[\sin\psi-\kappa_{z}(\psi)]+O(\varepsilon^{2})\}.
\end{split}
\]
Then, by equating real and imaginary parts of this equation, the parametric equations of the airfoil are, to first order in $\varepsilon$,
\begin{equation}\label{eq:chap_15_7_1}
\begin{split}
x=&c\{\frac{2+\cos\phi\sin2\psi}{1+\cos\phi\sin2\psi}\cos\phi\cos\psi+\kappa_{x}(\phi,\psi)\\
  &+\mu(\phi,\psi)[\frac{\cos\phi\sin2\psi}{1+\cos\phi\sin2\psi}\cos\phi\cos\psi-\kappa_{x}(\phi,\psi)]\} \\
 =&c\xi(\phi,\psi)\cos\phi\cos\psi+c[1-\mu(\phi,\psi)]\kappa_{x}(\phi,\psi), \\
y=&c[\frac{\cos\phi\sin2\psi}{1+\cos\phi\sin2\psi}+\mu(\phi,\psi)\frac{2+\cos\phi\sin2\psi}{1+\cos\phi\sin2\psi}]\sin\phi\cos\psi\\
 =&c\eta(\phi,\psi)\sin\phi\cos\psi, \\
z=&c[(1+\kappa_{z}^{2}(\psi))+\mu(\phi,\psi)(1-\kappa_{z}^{2}(\psi))]\sin\psi=c\zeta(\phi,\psi)\sin\psi
\end{split}
\end{equation}
where if $\sin\psi=0$ then $\xi(\phi,\psi)=2$ and $\eta(\phi,\psi)=2\mu(\phi,0)$ else
\[
\begin{split}
 \xi(\phi,\psi)&=1+\mu(\phi,\psi)+\frac{1-\mu(\phi,\psi)}{1+\cos\phi\sin2\psi},\\
\eta(\phi,\psi)&=1+\mu(\phi,\psi)-\frac{1-\mu(\phi,\psi)}{1+\cos\phi\sin2\psi},\\
\zeta(\phi,\psi)&=1+\mu(\phi,\psi)+[1-\mu(\phi,\psi)]\kappa_{z}^{2}(\psi).
\end{split}
\]
Using the first and third of these equations to eliminate $\phi$ from the second equation gives the following equation for the airfoil profile:
%So we get
\[
\{\frac{x-c[1-\mu(\phi,\psi)]\kappa_{x}(\phi,\psi)}{c\xi(\phi,\psi)}\}^{2}+[\frac{y}{c\eta(\phi,\psi)}]^{2}+[\frac{z}{c\zeta(\phi,\psi)}]^{2}=1
\]
%and
%\[
%\xi(\phi,\psi)=\eta(\phi,\psi)+2\frac{1-\mu(\phi,\psi)}{1+\cos\phi\sin2\psi}
%\]
or
\[
y=\pm c\eta(\phi,\psi)\sqrt{1-\{\frac{x-c[1-\mu(\phi,\psi)]\kappa_{x}(\phi,\psi)}{c\xi(\phi,\psi)}\}^{2}-[\frac{z}{c\zeta(\phi,\psi)}]^{2}}
\]
where if $\sin\psi\neq0$ then from equations
\[
\begin{split}
\mu(\phi,\psi)&=\varepsilon(1-\cos\phi\cos\psi) \\
              &=\frac{\zeta(\phi,\psi)}{1-\kappa_{z}^{2}(\psi)}-\frac{1+\kappa_{z}^{2}(\psi)}{1-\kappa_{z}^{2}(\psi)}
               =\frac{z/(c\sin\psi)}{1-\kappa_{z}^{2}(\psi)}-\frac{1+\kappa_{z}^{2}(\psi)}{1-\kappa_{z}^{2}(\psi)}
\end{split}
\]
we get
\[
\cos\phi=\{1-\frac{1}{\varepsilon}[\frac{z/(c\sin\psi)}{1-\kappa_{z}^{2}(\psi)}-\frac{1+\kappa_{z}^{2}(\psi)}{1-\kappa_{z}^{2}(\psi)}]\}\sec\psi.
\]

The location of the maximum thickness may now be obtained, and this is most readily done by using the parametric equation for the coordinate $y$ as derived above. Thus setting $dy/d\phi=0$ for a maximum in $y$ gives the following equation for the value of $\phi$ at the maximum thickness when the angle $\psi$ is fixed:
\[
\begin{split}
&\{\varepsilon\frac{2+\cos\phi\sin2\psi}{1+\cos\phi\sin2\psi}-2\frac{[1-\mu(\phi,\psi)]\sin\psi}{(1+\cos\phi\sin2\psi)^{2}}\}\sin^{2}\phi\cos\psi\\
+&[1+\mu(\phi,\psi)-\frac{1-\mu(\phi,\psi)}{1+\cos\phi\sin2\psi}]\cos\phi=0 \\
\end{split}
\]
where
\[
\frac{dy}{d\phi}=c[\frac{d\eta(\phi,\psi)}{d\phi}\sin\phi+\eta(\phi,\psi)\cos\phi]\cos\psi=0
\]
and
\[
\begin{split}
\frac{d\eta(\phi,\psi)}{d\phi}&=(1+\frac{1}{1+\cos\phi\sin2\psi})\frac{d\mu(\phi,\psi)}{d\phi}-\frac{1-\mu(\phi,\psi)}{(1+\cos\phi\sin2\psi)^{2}}\sin\phi\sin2\psi\\
                              &=\{\varepsilon\frac{2+\cos\phi\sin2\psi}{1+\cos\phi\sin2\psi}-2\frac{[1-\mu(\phi,\psi)]\sin\psi}{(1+\cos\phi\sin2\psi)^{2}}\}\sin\phi\cos\psi.\\
\end{split}
\]
This relation reduces to
\begin{equation}\label{eq:chap_15_7_2}
\begin{split}
&\varepsilon\{(\cos\phi-\cos2\phi\cos\psi)[(2+\cos\phi\sin2\psi)\cos\phi\sin2\psi+2]\} \\
&+[\cos2\phi+\cos^{3}\phi\sin2\psi+\varepsilon(1-\cos^{3}\phi\cos\psi)]\sin2\psi=0.
\end{split}
\end{equation}
For $-\pi/2<\psi<\pi/2$, given a fixed angle $\psi$, let $\phi_{k}$ are solutions of Eq.~(\ref{eq:chap_15_7_2}) where $k=1,2,\ldots,n$. Then at least a solution $\phi_{k}$ corresponds to the minimum thickness and at least a solution $\phi_{k}$ gives the maximum thickness, and for these values of $\phi$ the coordinates of the airfoil surface are determined by Eq.~(\ref{eq:chap_15_7_1}) as follows:
\[
\begin{split}
x&=c\xi(\phi,\psi)\cos\phi\cos\psi+c[1-\mu(\phi,\psi)]\kappa_{x}(\phi,\psi), \\
y&=c\eta(\phi,\psi)\sin\phi\cos\psi, \\
z&=c\zeta(\phi,\psi)\sin\psi.
\end{split}
\]
The maximum thickness $t$ will be twice the positive value of $y$, so that substituting these relations into expressions~(\ref{eq:chap_15_7_1}), the thickness ratio $t/l$ of the airfoil will be
\begin{equation}\label{eq:chap_15_7_4}
\frac{t}{l}=\frac{1}{2}\eta(\phi,\psi)\sin\phi\cos\psi.
\end{equation}
That is, the thickness-to-chord ratio of the airfoil is a function of the angle $\psi$ and the parameter $\varepsilon$, which is the ratio of the offset of the center of the sphere in the $\varpi$ space to the radius $c$ of the critical points of the transformation. Since the thickness ratio of the airfoil is a parameter that may be thought of as being specified, it is useful to eliminate $\varepsilon$ in terms of this parameter. Hence
\begin{equation}\label{eq:chap_15_7_5}
\varepsilon=\frac{1}{(1-\cos\phi\cos\psi)}[\frac{2}{\sin\phi\cos\psi}\cdot\frac{1+\cos\phi\sin2\psi}{2+\cos\phi\sin2\psi}\cdot\frac{t}{l}-\frac{\cos\phi\sin2\psi}{2+\cos\phi\sin2\psi}].
\end{equation}
Then the equation of the airfoil surface may be written in the form
\begin{equation}\label{eq:chap_15_7_6_4_23_a}
\frac{y}{t}=\frac{c}{t}\eta(\phi,\psi)\sin\phi\cos\psi
\end{equation}
where the maximum value of $y/t$ will be $1/2$ and the minimum value will be $-1/2$, both of which occur at $x=-c[1-\mu(\phi,\psi)]\kappa_{x}(\phi,\psi)$.

The magnitude of the circulation required to satisfy the Kutta condition is, from Eq.~(\ref{eq:chap_15_3_3_4_16}), $4\pi Ua\sin\phi$, where $a=c+m$ and $m=c\varepsilon$. Thus, for $-\pi/2<\psi<\pi/2$, the required amount of circulation is
\begin{equation}\label{eq:chap_15_7_7_4_23_b}
\Gamma=\pi Ul(1+\varepsilon)\sin\phi
\end{equation}
where $c$ has been replaced by $l/4$. In this form the required circulation may be evaluated for the given free-stream velocity, angle of attack, and the chord and thickness of the airfoil. Using the Kutta-Joukowski law ,%[Eq.~(\ref{eq:chap_15_3_3_4_16})],
the lift force acting on the airfoil may be evaluated as
\[
F_{L}=\pi\rho U^{2}l(1+\varepsilon)\sin\phi.
\]
Then the lift coefficient for the symmetrical Joukowski airfoil is
\begin{equation}\label{eq:chap_15_7_8_4_23_c}
C_{L}=2\pi(1+\varepsilon)\sin\phi.
\end{equation}
It will be noticed that this result reduces to Eq.~(\ref{eq:chap_15_6_3_4_22_c}) for the flat-plate airfoil as $t \to 0$. This indicates that the effect of thickness of an airfoil is to increase the lift coefficient. However, this fact cannot be used to produce high lift coefficients through thick airfoils, since the flow tends to separate from bluff bodies much more readily than it does from streamlined bodies. This separation of the flow is a viscous effect, and it will be discussed in the next part of the book. In the meantime, it is sufficient to say that separation of the flow results in a low-pressure wake that destroys the lift. The same result may occur for slender bodies, such as airfoils, that are at large angles of attack. In this context the separation is usually referred to as stall.

The center of the circle in the $\varpi$ space is located at $\varpi=-e_{xy}m$ rather than $\varpi=0$. Thus the complex potential for the flow in the $\varpi$ space may be obtained from Eq.~(\ref{eq:chap_15_5_1_4_29})
by replacing $\varpi$ by $\varpi=-e_{xy}m$ and adding circulation. The required complex potential then becomes
\begin{equation}\label{eq:chap_15_7_9_4_23_d}
F(\varpi)=U[(\varpi+m)\bar{\gamma}+\frac{a^{2}}{\varpi+m}\frac{\gamma}{\bar{\gamma}\gamma}]+\frac{i\Gamma}{2\pi}\log(\frac{\varpi+m}{a})
\end{equation}
where $\gamma=e_{xy}e^{i\alpha}\cos\beta+\sin\beta$,
\[
a=\frac{l}{4}+c\varepsilon\texttt{ and }m=c\varepsilon.
\]
The magnitude of the circulation $\Gamma$ is given by Eq.~(\ref{eq:chap_15_7_7_4_23_b}), and in the Joukowski transformation the parameter $c$ equals $l/4$. %The flow field corresponding to this complex potential is shown in Fig.4.16b.

\section{Circular-Arc Airfoils}\label{sec:chap_15_8_146}

It was shown in the two previous sections that, using the Joukowski transformation, a sphere of radius $c$ centered at the origin of the $\varpi$ space produced a flat-plate airfoil while a slightly larger sphere centered a small distance along the real axis from the origin produced a thin symmetrical airfoil. It will now be shown that a sphere whose radius is slightly larger than $c$ and whose center is located on the imaginary axis of the $\varpi$ space produces an airfoil that has no thickness but has curvature of camber.

Consider a sphere of radius $a>c$ in the $\varpi$ space such that the center of the sphere is located a distance $m$ along the positive imaginary axis. Since the trailing edge of the airfoil should be sharp, the sphere should pass through the critical point $\varpi=e_{xy}c$ as before. Then, in this case, the sphere will also pass through the other critical point, $\varpi=-e_{xy}c$.

The equation of the airfoil in the $s$ space may be obtained by substituting $\varpi=R(e_{xy}e^{i\phi}\cos\psi+\sin\psi)$ into the Joukowski transformation, where, on the big circumference of the sphere in the $\varpi$ space, $R$ is a function of a pair of angles $\phi$ and $\psi$. When $\sin\psi\neq0$, this substitution gives
\[
\begin{split}
s=&R(e_{xy}e^{i\phi}\cos\psi+\sin\psi)+\frac{c^{2}}{R(e_{xy}e^{i\phi}\cos\psi+\sin\psi)} \\
 =&R\{e_{xy}[e^{i\phi}\cos\psi+\frac{c^{2}}{R^{2}}(\frac{e^{-i\phi}\cos\psi+\sin\psi}{1+\cos\phi\sin2\psi}-\frac{1}{\sin\psi})]+(\sin\psi+\frac{c^{2}/R^{2}}{\sin\psi})\}\\
 =&R\{e_{xy}[(1+\frac{c^{2}/R^{2}}{1+\cos\phi\sin2\psi})\cos\phi\cos\psi+\frac{c^{2}}{R}\kappa_{x}(\phi,\psi) \\
  &+i(1-\frac{c^{2}/R^{2}}{1+\cos\phi\sin2\psi})\sin\phi\cos\psi]+[\sin\psi+\frac{c^{2}}{R^{2}}\kappa_{z}(\psi)]\}
\end{split}
\]
where if $\sin\psi=0$ then $\kappa_{x}(\phi,\psi)=0$ and $\kappa_{z}(\psi)=0$ else
\[
\kappa_{x}(\phi,\psi)=\frac{\sin\psi}{1+\cos\phi\sin2\psi}-\frac{1}{\sin\psi}\texttt{ and }\kappa_{z}(\psi)=\frac{1}{\sin\psi}.
\]
Thus, when $\sin\psi\neq0$, the parametric equations of the airfoil profile are
\begin{equation}\label{eq:chap_15_8_1}
\begin{split}
x&=(R+\frac{c^{2}/R}{1+\cos\phi\sin2\psi})\cos\phi\cos\psi+\frac{c^{2}}{R}\kappa_{x}(\phi,\psi),\\
y&=(R-\frac{c^{2}/R}{1+\cos\phi\sin2\psi})\sin\phi\cos\psi, \\
z&=R\sin\psi+\frac{c^{2}}{R\sin\psi}=R\sin\psi+\frac{c^{2}}{R}\kappa_{z}(\psi).
\end{split}
\end{equation}
The variable $R$ may be eliminated from these equations as follows:
\[
\begin{split}
&[x-\frac{c^{2}}{R}\kappa_{x}(\phi,\psi)]^{2}\sin^{2}\phi-y^{2}\cos^{2}\phi \\
=&[(R+\frac{c^{2}/R}{1+\cos\phi\sin2\psi})^{2}-(R-\frac{c^{2}/R}{1+\cos\phi\sin2\psi})^{2}]\sin^{2}\phi\cos^{2}\phi\cos^{2}\psi \\
=&\frac{4c^{2}\cos^{2}\psi}{1+\cos\phi\sin2\psi}\sin^{2}\phi\cos^{2}\phi
\end{split}
\]
where from the equation $(R\sin\psi)^{2}-zR\sin\psi+c^{2}=0$
\[
R\sin\psi=(z\pm\sqrt{z^{2}-4c^{2}})/2.
\]
Then there are
\begin{equation}\label{eq:chap_15_8_2}
[x-\frac{c^{2}}{R}\kappa_{x}(\phi,\psi)]^{2}\sin^{2}\phi-y^{2}\cos^{2}\phi=\frac{4c^{2}\cos^{2}\psi}{1+\cos\phi\sin2\psi}\sin^{2}\phi\cos^{2}\phi.
\end{equation}
This is the equation of the airfoil surface in the $s$ space when $\sin\psi\neq0$, but it still contains the variable $\phi$ for each fixed angle $\psi$. This variable may be eliminated by applying the cosine rule to the triangle defined by the radius $a$, the coordinate $R$, and the imaginary $\varpi$ axis. From these it follows that
\[
\begin{split}
a^{2}&=(R\cos\psi)^{2}+m^{2}-2Rm\cos(\pi/2-\phi)\cos\psi+(R\sin\psi)^{2} \\
     &=R^{2}+m^{2}-2Rm\sin\phi\cos\psi, \\
c^{2}+m^{2}&=R^{2}+m^{2}-2Rm\sin\phi\cos\psi
\end{split}
\]
where the fact that $a^{2}=c^{2}+m^{2}$ has been used. Solving this equation for $\sin\phi$, it follows that
\[
\sin\phi=\frac{R^{2}-c^{2}}{2Rm\cos\psi},
\]
and
\[
\cos\phi=\pm\sqrt{1-(\frac{R^{2}-c^{2}}{2Rm\cos\psi})^{2}}.
\]
But it was shown from equations~(\ref{eq:chap_15_8_1}) that
\[
y=(R-\frac{c^{2}/R}{1+\cos\phi\sin2\psi})\sin\phi\cos\psi,
\]
hence it follows that
\[
\sin\phi=\frac{y(1+\cos\phi\sin2\psi)}{2m\sin\phi\cos^{2}\psi}-\frac{R\sin\psi}{m}\cos\phi
\]
or
\[
\begin{split}
\sin^{2}\phi&=\frac{y(1+\cos\phi\sin2\psi)}{2m\cos^{2}\psi}-\frac{R\sin\psi}{m}\sin\phi\cos\phi\\
            &=\frac{y}{2m}\kappa_{y}(y,\phi,\psi)
\end{split}
\]
and so
\[
\cos^{2}\phi=1-\frac{y}{2m}\kappa_{y}(y,\phi,\psi)
\]
where if $\sin\psi=0$ then $\kappa_{y}(y,\phi,\psi)=1$ else
\[
\kappa_{y}(y,\phi,\psi)=\frac{1+\cos\phi\sin2\psi}{\cos^{2}\psi}\mp\frac{\sin\psi}{my}\cdot\frac{R^{2}-c^{2}}{\cos\psi}\sqrt{1-(\frac{R^{2}-c^{2}}{2Rm\cos\psi})^{2}}.
\]
Using these results to eliminate $\phi$, the equation~(\ref{eq:chap_15_8_2}) of the airfoil surface becomes
\[
\begin{split}
&[x-\frac{c^{2}}{R}\kappa_{x}(\phi,\psi)]^{2}\frac{y}{2m}\kappa_{y}(y,\phi,\psi)-y^{2}[1-\frac{y}{2m}\kappa_{y}(y,\phi,\psi)]\\
=&\frac{4c^{2}\cos^{2}\psi}{1+\cos\phi\sin2\psi}\frac{y}{2m}\kappa_{y}(y,\phi,\psi)[1-\frac{y}{2m}\kappa_{y}(y,\phi,\psi)].
\end{split}
\]
Collecting like terms, this equation may be put in the form
\[
[x-\frac{c^{2}}{R}\kappa_{x}(\phi,\psi)]^{2}+y^{2}+2[\frac{\kappa_{y}(y,\phi,\psi)\cos^{2}\psi}{1+\cos\phi\sin2\psi}\cdot\frac{c^{2}}{m}-\frac{m}{\kappa_{y}(y,\phi,\psi)}]y
=\frac{4c^{2}\cos^{2}\psi}{1+\cos\phi\sin2\psi}.
\]
Completing the square in $y$ shows that the equation of the airfoil surface is
\[
\begin{split}
&[x-\frac{c^{2}}{R}\kappa_{x}(\phi,\psi)]^{2}+\{y+c[\frac{\kappa_{y}(y,\phi,\psi)\cos^{2}\psi}{1+\cos\phi\sin2\psi}\cdot\frac{c}{m}-\frac{m}{c\kappa_{y}(y,\phi,\psi)}]\}^{2}\\
=&c^{2}\{\frac{4\cos^{2}\psi}{1+\cos\phi\sin2\psi}+[\frac{\kappa_{y}(y,\phi,\psi)\cos^{2}\psi}{1+\cos\phi\sin2\psi}\cdot\frac{c}{m}-\frac{m}{c\kappa_{y}(y,\phi,\psi)}]^{2}\}
\end{split}
\]
which is the equation of a big circle of a sphere. It should be noted that so far no approximations have been made. But to be consistent with the analysis in the previous section and to permit superposition in the next section, the parameter $\varepsilon=m/c$ will again be assumed to be small compared with unity. Then,
linearizing in $\varepsilon$, the equation of the airfoil surface becomes
\[
\begin{split}
&[x-\frac{c^{2}}{R}\kappa_{x}(\phi,\psi)]^{2}+[y+\frac{\kappa_{y}(y,\phi,\psi)\cos^{2}\psi}{1+\cos\phi\sin2\psi}\cdot\frac{c^{2}}{m}]^{2}\\
=&c^{2}\{\frac{4\cos^{2}\psi}{1+\cos\phi\sin2\psi}+[\frac{\kappa_{y}(y,\phi,\psi)\cos^{2}\psi}{1+\cos\phi\sin2\psi}\cdot\frac{c}{m}]^{2}\}.
\end{split}
\]
That is, correct to the first order in $\varepsilon$, the center of the big circle of the sphere in the $s$ space is located at
\[
x=\kappa_{x}(\phi,\psi)\texttt{ and }
y=-\frac{\kappa_{y}(y,\phi,\psi)\cos^{2}\psi}{1+\cos\phi\sin2\psi}\cdot\frac{c^{2}}{m},
\]
and the radius of the big circle is
\[
c\sqrt{\frac{4\cos^{2}\psi}{1+\cos\phi\sin2\psi}+[\frac{\kappa_{y}(y,\phi,\psi)\cos^{2}\psi}{1+\cos\phi\sin2\psi}\cdot\frac{c}{m}]^{2}}.
\]

The characteristic parameters of the airfoil are the chord $l$ and the camber height $h$. %, and these are shown in Fig.4.17a.
Since the equation of the airfoil has now been established, it is possible to relate these parameters
to those in the $\varpi$ space, namely, $c$ and $m$. Since the ends of the circular-arc airfoil lie on the real axis $y=0$, the foregoing equation of the airfoil shows that the corresponding values of $x$ are $\pm2c$. That is, the chord of the airfoil is
\[
l=4c.
\]
This is the same chord length as for the two previous airfoils.

The simplest way of establishing the camber $h$ of the airfoil is to use the fact that, in view of the result that the center of the circular arc is at $x=0$, the maximum value of $y$ will occur when $x=0$. But the parametric equation
\[
x=(R+\frac{c^{2}/R}{1+\cos\phi\sin2\psi})\cos\phi\cos\psi+\frac{c^{2}}{R}\kappa_{x}(\phi,\psi),
\]
shows that this corresponds to $\phi=\pi/2$ when $\sin\psi=0$. Then the other parametric equation, namely, $y=2m\sin^{2}\phi$ when $\sin\psi=0$, shows that the maximum value of $y$ is 2m. That is,
\[
h=2m.
\]

Using the foregoing results, the $\varpi$-space parameters $c$ and $m$ may be replaced by the $s$-space parameters $l=4$ and $h=2$, respectively. Then the equation of the airfoil surface in the $s$ space may be written in the form
\begin{equation}\label{eq:chap_15_8_3_4_24_a}
\begin{split}
&[x-\frac{c^{2}}{R}\kappa_{x}(\phi,\psi)]^{2}+[y+\frac{\kappa_{y}(y,\phi,\psi)\cos^{2}\psi}{1+\cos\phi\sin2\psi}\cdot\frac{l^{2}}{8h}]^{2}\\
=&\frac{l^{2}}{4}\{\frac{\cos^{2}\psi}{1+\cos\phi\sin2\psi}+[\frac{\kappa_{y}(y,\phi,\psi)\cos^{2}\psi}{1+\cos\phi\sin2\psi}\cdot\frac{l}{4h}]^{2}\}.
\end{split}
\end{equation}

In order to satisfy the Kutta condition, the rear stagnation point must rotate through an angle greater than $\phi$, the angle of the free stream. By rotating through the angle $\phi$, the rear stagnation point will be located on the surface of the big circle in the $\varpi$ space at a point which is in the same horizontal plane as the center of the big circle. But the center of the big circle is located a distance $m$ above the real $\varpi$ axis. Thus, in order to be located at the point $\varpi=e_{xy}c$, the rear stagnation point must rotate through a further angle given by
\[
\arctan\frac{m}{c}=\arctan\varepsilon=\varepsilon+O(\varepsilon^{2}).
\]
That is, in order to comply with the Kutta condition, the rear stagnation point must rotate through the angle $\phi+m/c$, to the first order in $\varepsilon$. Then, from Eq.~(\ref{eq:chap_15_3_3_4_16}), the required circulation is
\[
\Gamma=4\pi Ua\sin(\phi+\frac{m}{c}),
\]
but $a=\sqrt{c^{2}+m^{2}}$ so that, to first order in $\varepsilon$, $a=c$. Hence
\[
\Gamma=4\pi Uc\sin(\phi+\frac{m}{c}).
\]
Then from the Kutta-Joukowski law, the lift force is
\[
F_{L}=4\pi\rho Uc\sin(\phi+\frac{m}{c})
\]
and the corresponding lift coefficient is
\[
C_{L}=8\pi\rho\frac{c}{l}\sin(\phi+\frac{m}{c}).
\]
Using again the fact that $c=l/4$ and $m=h/2$, the lift coefficient becomes
\begin{equation}\label{eq:chap_15_8_4_4_24_b}
C_{L}=2\pi\frac{c}{l}\sin(\phi+\frac{2h}{l}).
\end{equation}
Comparing this result with Eq.~(\ref{eq:chap_15_6_3_4_22_c}), the corresponding result for the flat
plate, shows that the effect of positive camber in an airfoil is to increase its lift coefficient. As a consequence of this increased lift coefficient a nonzero lift exists at zero angle of attack.

Since the center of the circle in the $\varpi$ space is at $\varpi=im$ rather than $\varpi=0$, the complex potential in the $\varpi$ space may be obtained by replacing $\varpi$ by $\varpi-im$ in Eq.~(\ref{eq:chap_15_5_1_4_29}) and adding circulation. Thus the required complex potential is
\begin{equation}\label{eq:chap_15_8_5_4_24_c}
F(\varpi)=U[(\varpi-im)\bar{\gamma}+\frac{a^{2}}{\varpi-im}\frac{\gamma}{\bar{\gamma}\gamma}]+\frac{i\Gamma}{2\pi}\log(\frac{\varpi-im}{a})
\end{equation}
where $\gamma=e_{xy}e^{i\alpha}\cos\beta+\sin\beta$,
\[
a=\frac{l}{4}\texttt{ and }m=\frac{h}{2}.
\]
The magnitude of the circulation $\Gamma$ is given by
\[
\Gamma=\pi Ul\sin(\phi+\frac{2h}{l})
\]
and the parameter $c$ in the Joukowski transformation is $l/4$. %The flow field corresponding to this complex potential is shown in Fig. 4.17b.
As was the case with the flat plate airfoil, this flow field has a singularity at the leading edge. This singularity would not exist for airfoils of finite nose radius and would not exist even for sharp leading edges because of separation of the flow at the nose. In spite of this local inaccuracy, the results derived above are representative of the flow around thin cambered airfoils.

\section{Joukowski Airfoils}\label{sec:chap_15_9_147}

The results of the two previous sections suggest that a cambered airfoil of finite thickness may be obtained by considering the Joukowski transformation in conjunction with a sphere in the $\varpi$ space whose center is in the second quadrant. Such a configuration is shown that the center of the sphere is displaced a distance $m$ from the origin at an angle $d$ from the reference axis. In order that the trailing edge of the corresponding airfoil may be sharp, the big circumference of the sphere passes through the critical point $\varpi=e_{xy}c$. The principal parameters of the airfoil in the $s$ space are also shown that these are the chord $l$, the maximum thickness $t$, and the maximum camber of the centerline $h$.

From the previous two sections it follows that, to the first order in $\varepsilon$, the centerline of the airfoil will be a circular arc whose center is on the $y$ axis and the airfoil will be symmetrical about its centerline. Then the equation of the upper and lower surfaces of the airfoil may be obtained from the equation for the circular-arc centerline plus or minus, respectively, and thickness effect. Hence from Eqs.~(\ref{eq:chap_15_7_6_4_23_a}) and~(\ref{eq:chap_15_8_3_4_24_a}) the airfoil profile will be given by
\begin{equation}\label{eq:chap_15_9_1_4_25_a}
\begin{split}
y=&\sqrt{\frac{l^{2}}{4}\{\frac{\cos^{2}\psi}{1+\cos\phi\sin2\psi}+[\frac{\kappa_{y}(y,\phi,\psi)\cos^{2}\psi}{1+\cos\phi\sin2\psi}\cdot\frac{l}{4h}]^{2}\}-[x-\frac{c^{2}}{R}\kappa_{x}(\phi,\psi)]^{2}}\\
  &-\frac{\kappa_{y}(y,\phi,\psi)\cos^{2}\psi}{1+\cos\phi\sin2\psi}\cdot\frac{l^{2}}{8h}
  \pm c\eta(\phi,\psi)\sin\phi\cos\psi\\
\end{split}
\end{equation}
where the plus defines the upper surface and the minus defines the lower surface of the so-called Joukowski airfoil.

It was observed that the effect of thickness on an airfoil was to increase its lift by an amount and that the effect of camber was to increase the effective angle of attack to $\phi+2h/l$. The present airfoil has both these effects, so that, from Eqs.~(\ref{eq:chap_15_7_8_4_23_c}) and~(\ref{eq:chap_15_8_4_4_24_b}), the lift coefficient for the Joukowski airfoil will be
\begin{equation}\label{eq:chap_15_9_2_4_25_b}
C_{L}=2\pi(1+\varepsilon)\sin(\phi+\frac{2h}{l}).
\end{equation}

The complex potential in the $\varpi$ space is
\begin{equation}\label{eq:chap_15_9_3_4_25_c}
F(\varpi)=U[(\varpi-m\varsigma)\bar{\gamma}+\frac{a^{2}}{\varpi-m\varsigma}\frac{\gamma}{\bar{\gamma}\gamma}]+\frac{i\Gamma}{2\pi}\log(\frac{\varpi-m\varsigma}{a})
\end{equation}
where $\gamma=e_{xy}e^{i\alpha}\cos\beta+\sin\beta$ and $\varsigma=e_{xy}e^{i\vartheta}\cos\tau+\sin\tau$. These results follow from the observation that $m\cos\vartheta$ replaces $m$ as used in the symmetrical Joukowski airfoil and $m\sin\vartheta$ replaces $m$ as used in the circular-arc airfoil. The magnitude of the circulation $\Gamma$ will include both the thickness and camber effects, and so it follows that
\[
\Gamma=\pi Ul(1+\varepsilon)\sin(\alpha+\frac{2h}{l}).
\]

%The flow field corresponding to the foregoing complex potential is shown in Fig. 4.18b.
It should be remembered that there is a limit to the amount of thickness and camber. %which may be introduced if the flow field is to remain as shown.
As the thickness and/or camber of the airfoil increases, the body departs more and more from a streamlined airfoil and approaches a bluff body. It was pointed out earlier that a consequence of this would be separation of the flow, which destroys the lift force and creates the so-called stall condition.

\backmatter
%-----------------------------------------------------------------------------
% Beginning of biblio.tex
%-----------------------------------------------------------------------------
%
% AMS-LaTeX 1.2 sample file for a monograph, based on amsbook.cls.
% This is a data file input by chapter.tex.
%%%%%%%%%%%%%%%%%%%%%%%%%%%%%%%%%%%%%%%%%%%%%%%%%%%%%%%%%%%%%%%%%%%%%%%%

\bibliographystyle{amsalpha}

%-----------------------------------------------------------------------------
% End of biblio.tex
%-----------------------------------------------------------------------------

%\include{index}
\end{document}